\documentclass{amsbook}

\usepackage{amssymb}
\usepackage{amsfonts}
\usepackage{amsmath}
\usepackage{lastpage}
\usepackage[legalpaper,bookmarks=true,colorlinks=true,linkcolor=blue,citecolor=blue]{hyperref}
\usepackage{graphicx}%
\usepackage{fancyhdr}
\usepackage{color}
\usepackage[mathlines]{lineno}
\usepackage{lscape}
\usepackage{epsfig}
\usepackage{natbib}

\newtheorem{theorem}{Theorem}
\theoremstyle{plain}

\newtheorem{claim}{Claim}

\newtheorem{definition}{Definition}

\numberwithin{equation}{section}

\begin{document}
\Large
\pagenumbering{roman}
\begin{center}

\huge \textbf{Gane Samb LO}\\
\vskip 6cm 
\bigskip \textbf{Measure Theory and Integration By and For the Learner}
\vskip 6cm

\huge \textit{\textbf{Statistics and Probability African Society (SPAS) Books Series}.\\
 \textbf{Saint-Louis, Senegal / Calgary, Alberta. 2017 (version 2017-01)}}.\\

\bigskip \Large  \textbf{DOI} : http://dx.doi.org/10.16929/sbs/2016.0005\\
\bigskip \textbf{ISBN} : 978-2-9559183-5-7
\end{center}

\newpage
\begin{center}
\huge \textbf{SPAS TEXTBOOKS SERIES}
\end{center}

\bigskip \bigskip

\Large

 \begin{center}
 \textbf{GENERAL EDITOR of SPAS EDITIONS}
 \end{center}

\bigskip
\noindent \textbf{Prof Gane Samb LO}\\
gane-samb.lo@ugb.edu.sn, gslo@ugb.edu.ng\\
Gaston Berger University (UGB), Saint-Louis, SENEGAL.\\
African University of Sciences and Technology, AUST, Abuja, Nigeria.\\

\bigskip

\begin{center}
\Large \textbf{ASSOCIATED EDITORS}
\end{center}

\bigskip
\noindent \textbf{KEhinde Dahud SHANGODOYIN}\\
shangodoyink@mopipi.ub.bw\\
University of Botswana (Botswana)\\

\noindent \textbf{Blaise SOME}\\
some@univ-ouaga.bf\\
Chairman of LANIBIO, UFR/SEA\\
Joseph Ki-Zerbo University (Ouaga I), Burkina-Faso.\\

\bigskip
\begin{center}
\Large \textbf{ADVISORS}
\end{center}

\bigskip

\noindent \textbf{Ahmadou Bamba SOW}\\
ahmadou-bamba.sow@ugb.edu.sn\\
Gaston Berger University, Senegal.\\

\noindent \textbf{Tchilabalo Abozou KPANZOU}\\
kpanzout@yahoo.fr\\
Kara University, Togo.\\

\newpage

\huge \textbf{SPAS Series Books}\\
\bigskip \bigskip

\Large \textbf{List of published or scheduled books in English}\\

\noindent $\square$ Weak Convergence (IA) - Sequences of Random Vectors. Gane Samb LO, Modou NGOM and Tchilabalo A. KPANZOU. 2016.\\
Doi : 10.16929/sbs/2016.0001. ISBN 978-2-9559183-1-9\\

\noindent $\square$ A Course on Elementary Probability Theory. Gane Samb LO. 2017.\\
Doi : 10.16929/sbs/2016.0003. ISBN 978-2-9559183-3-3\\

\newpage
\noindent \textbf{Library of Congress Cataloging-in-Publication Data}\\

\noindent Gane Samb LO, 1958-\\

\noindent Measure Theory and Integration By and For the Learner.\\

\noindent SPAS Books Series, 2017.\\

\noindent To be published SPAS Eds - Statistics and Probability African Society (SPAS).\\

\noindent \textit{DOI} : 10.16929/sbs/2016.0005\\

\noindent \textit{ISBN} : 978-2-9559183-5-7\\

\noindent Version 2017-01

\newpage

\noindent \textbf{Author : Gane Samb LO}\\
\bigskip

\bigskip
\noindent \textbf{Emails}:\\
\noindent gane-samb.lo@ugb.edu.sn, ganesamblo@ganesamblo.net.\\

\bigskip
\noindent \textbf{Url's}:\\
\noindent www.ganesamblo@ganesamblo.net\\
\noindent www.statpas.net/cva.php?email.ganesamblo@yahoo.com.\\

\bigskip \noindent \textbf{Affiliations}.\\
Main affiliation : University Gaston Berger, UGB, SENEGAL.\\
African University of Sciences and Technology, AUST, ABuja, Nigeria.\\
Affiliated as a researcher to : LSTA, Pierre et Marie Curie University, Paris VI, France.\\

\noindent \textbf{Teaches or has taught} at the graduate level in the following universities:\\
Saint-Louis, Senegal (UGB)\\
Banjul, Gambia (TUG)\\
Bamako, Mali (USTTB)\\
Ouagadougou - Burkina Faso (UJK)\\
African Institute of Mathematical Sciences, Mbour, SENEGAL, AIMS.\\
Franceville, Gabon\\

\bigskip \noindent \textbf{Dedication}.\\

\noindent \textbf{To all the students who followed my courses on Measure Theory and Integration in many universities}.\\

\bigskip \noindent \textbf{Acknowledgment of Funding}.\\

\noindent The author acknowledges continuous support of the World Bank Excellence Center in Mathematics, Computer Sciences and Intelligence Technology, CEA-MITIC. His research projects in 2014, 2015 and 2016 are funded by the University of Gaston Berger in different forms and by CEA-MITIC.

\title{Measure Theory and Integration By and For the Learner}

\begin{abstract} $ $\\

\noindent  Measure Theory and Integration is exposed with the clear aim to help beginning learners to perfectly master its essence. In opposition of a delivery of the contents in an academic and vertical course, the knowledge is broken into exercises which are left to the learners for solutions. Hints are present at any corner to help readers to achieve the solutions. In that way, the knowledge is constructed by the readers by summarizing the results of one or a group of exercises.\\

\noindent Each chapter is organized into Summary documents which contain the knowledge, Discovery documents which give the learner the opportunity to extract the knowledge himself through exercises and into Solution Documents which offer detailed answers for the exercises. Exceptionally, a few number of results (A key lemma related the justification of definition of the integral of a non-negative function, the Caratheodory's theorem and the Lebesgue-Stieljes measure on $\mathbb{R}^d$) are presented in appendix documents and given for reading in small groups.\\

\noindent The full theory is presented in the described way. We highly expect that any student who goes through the materials, alone or in a small group or under the supervision of an assistant will gain a very solid knowledge in the subject and by the way ensure a sound foundation for studying disciplines such as Probability Theory, Statistics, Functional Analysis, etc.\\

\noindent The materials have been successfully used as such in normal real analysis classes several times.\\

\noindent \textbf{Keywords} Sets, sets operations, measurable sets; measurable applications; Borel sets; Borel Applications; measures; Lebesgue-Stieljes Measures;
Counting measures; finite products measures; finite and sigma-finite measures; extension of measures; Caratheodory's Theorem; Convergence theorems; almost-sure convergence; convergence in measure; integration; Riemann-Stieljes integration; Lebesgue-Stieljes integration; series and integrals with respect to counting measures; Continuity and differentiability of integral depending on continuous parameter; Fubini's and Tonelli's Theorem; Hahn-Jordan decomposition Theorem; Lebesgue Decomposition Theorem; Radon-Nikodym's Theorem; High dimensional Lebesgue-Stieljes Integrals. $L^p$ spaces. 
\\

\noindent \textbf{AMS 2010 Classification Subjects :} 28A05; 28A10; 28A12; 28A20; 28A25; 28A33; 28A35
\end{abstract}

\maketitle

\frontmatter
\tableofcontents
%\listoftables
\mainmatter
\Large

\chapter*{General Preface}

\noindent \textbf{This textbook} is one of the elements of a series whose ambition is to cover a broad part of Probability Theory and Statistics.  These textbooks are intended to help learners and readers, of all levels, to train themselves.\\

\noindent As well, they may constitute helpful documents for professors and teachers for both courses and exercises.  For more ambitious  people, they are only starting points towards more advanced and personalized books. So, these textbooks are kindly put at the disposal of professors and learners.

\bigskip \noindent \textbf{Our textbooks are classified into categories}.\\

\noindent \textbf{A series of introductory  books for beginners}. Books of this series are usually destined to students of first year in universities and to any individual wishing to have an initiation on the subject. They do not require advanced mathematics.  Books on elementary probability theory (See \cite{ips-prob-ang}, for instance) and descriptive statistics are to be put in that category. Books of that kind are usually introductions to more advanced and mathematical versions of the same theory. Books of the first kind also prepare the applications of those of the second.\\

\noindent \textbf{A series of books oriented to applications}. Students or researchers in very related disciplines  such as Health studies, Hydrology, Finance, Economics, etc.  may be in need of Probability Theory or Statistics. They are not interested in these disciplines  by themselves.  Rather, they need to apply their findings as tools to solve their specific problems. So, adapted books on Probability Theory and Statistics may be composed to focus on the applications of such fields. A perfect example concerns the need of mathematical statistics for economists who do not necessarily have a good background in Measure Theory.\\

\noindent \textbf{A series of specialized books on Probability theory and Statistics of high level}. This series begins with a book on Measure Theory, a book on its probability theory version, and an introductory book on topology. On that basis, we will have, as much as possible,  a coherent presentation of branches of Probability theory and Statistics. We will try  to have a self-contained approach, as much as possible, so that anything we need will be in the series.\\

\noindent Finally, a series of \textbf{research monographs} closes this architecture. This architecture should be so diversified and deep that the readers of monograph booklets will find all needed theories and inputs in it.\\

\bigskip \noindent We conclude by saying that, with  only an undergraduate level, the reader will  open the door of anything in Probability theory and statistics with \textbf{Measure Theory and integration}. Once this course validated, eventually combined with two solid courses on topology and functional analysis, he will have all the means to get specialized in any branch in these disciplines.\\

\bigskip \noindent Our collaborators and former students are invited to make live this trend and to develop it so  that  the center of Saint-Louis becomes or continues to be a re-known mathematical school, especially in Probability Theory and Statistics.

\chapter*{Introduction}

\noindent \textbf{Measure Theory and Integration}.\\

\noindent Undoubtedly, Measure Theory and Integration is one of the most important part of Modern Analysis, with Topology and Functional Analysis for example. Indeed, Modern mathematics is based on functional analysis, which is a combination of the Theory of Measure and Integration, and Topology.\\

\noindent The application of mathematics is very pronounced in many fields, such as finance (through stochastic calculus), mathematical economics (through stochastic calculus), econometrics, which is a contextualization of statistical regression to economic problems, physic statistics. Probability Theory and Statistics has become an important tool for the analysis of biological phenomena and genetics.\\

\noindent In all these cases, a solid basis in probability and mathematical statistics is required. This can be achieved only with Measure theory and integration through its powerful theorems of convergence, differentiation and integrability.\\

\noindent \textbf{The pedagogical orientation and organization of the contents}.\\

\noindent In this book, the contents of Measure Theory and Integration are presented to the learners as solutions of exercises they are asked to solve. The opportunity of extracting the knowledge themselves is given to the learner. Undoubtedly, this makes the learner feel more confident in himself, makes him fully participate in the classes and leads to a great impact of the course on his level.\\

\noindent This method is used with the clear aim to help to beginning learners to perfectly master the essence of the theory. In opposition of a delivery of the contents in an academic and vertical course, the knowledge is broken into exercises which are left to the learners for solutions. Hints are present at any corner to help readers to achieve the solutions. In that way, the knowledge is constructed by the readers by summarizing the results of one or a group of exercises.\\

\noindent Let us first describe how the book is organized. We tried to present all the chapters in the same way.\\

\noindent (a) First, we give one, two or three documents in which the knowledge is summarized. The learner is recommended to know the statements in these summary as best as possible.\\

\noindent (b) Second, we put the knowledge in summary documents as conclusions of a series or series of exercises. This means that we allow the learner to discover himself the the elements of the courses.\\

\noindent These discovery exercises are accompanied with all the hints that allow the learner to proceed alone, or in a small group.\\

\noindent (c) If we face a very heavy material (in length or in complexity), we prefer to give it in an Appendix. The number of such documents does not exceed three in the whole courses.\\

\noindent An assistant may be appointed to accompany small groups of learners to read and repeat the proofs in these Appendix texts.\\

\noindent (d) The solutions of the discovery exercises are finally given. In these documents, the exercises are listed in the same order and the solution of each exercise is put just at its end. In general, the learner will see that all the reasoning has already been given as hints.\\

\noindent At the beginning of each chapter, a table will present the documents with the clear indications of the type : $S$ for Summary, $D$ for Discover Exercises, $SD$ for Solutions of the Discovery Exercises, $A$ for Appendix.\\

\noindent So, each chapter must have at least three documents.\\

\noindent These documents of this book my be used as class materials, allowing the teacher to let learners advancing at their own pace and next to help them to get the knowledge.\\

\noindent Strongly motivated students are free to go faster and to finish the course in a short time. But they are recommended to attend the classes and to participate in the exchanges.\\

\noindent It is important to notice that the textbook begins by Chapter 0. This means that all the content of that chapter is supposed to be known in the former courses of General Algebra in the first year of University. But, we need them as the basis of the whole theory. We exposed it seriously as the other chapters by supposing that the reader might have forgotten them.\\

\bigskip \noindent A parallel textbook entirely devoted to advanced exercises and using the knowledge developed in this book will be available in a recent future.\\

\noindent We hope that the learners will profit of this way of teaching measure Theory and Integration, and will have the strongest basis to advance to probability theory and Functional Analysis.\\

\noindent \textbf{References}.\\

\noindent Personally, my first contact with the subject happened in the Integration course in 1983 in Dakar University with professor Sakhir Thiam. The course was based on the Daniel integral 
(see Doc  09-04 and Doc 09-05, pages \pageref{doc09-04} and \pageref{doc09-05}, for an introduction to that approach) and  some elements of Modern measure Theorem from the books of \cite{genet} on Measure and Integration  \cite{metivier} on Probability Theory. That approach heavily depends on topological structures, mainly on locally compact spaces.\\

\noindent During my graduate studies in Paris (1983-1986) I discovered the nice book of \cite{loeve} and I have never stopped using it. Although the book is the result of more than two decades of
our teaching, it is remains greatly influenced by \cite{loeve}.\\

\noindent I use many other books for finding different proofs, statements, examples, etc. For example, the approach of the product measure is taken from \cite{revuz}. We used notations from \cite{kacimi} to propose a mode advanced way to make a comparison between the Lebesgue and the Riemann integrals based on semi-continuous functions. The readings of the seminal book of  \cite{halmosm} and the extremely rich books of \cite{bogachev1} and \cite{bogachev2} allowed to provide very useful examples and orientation. For example, the exercises on Daniel's integration approach in Doc 09-04 is inspired by \cite{bogachev2}.\\

\noindent But at the arrival, most of the results been stated as felt in classes by the students and myself.

%\cite{billinsgleypm} \cite{bogachev1} \cite{bogachev2} \cite{henri} \cite{foatafuchs} \cite{gramain} \cite{genet} \cite{gutt} \cite{halmosm}(1950) \cite{kacimi}
%\cite{loeve} \cite{marle}\cite{metivier} \cite{partha}  \cite{sacks} Stanislaw Sacks (1933) 

\LARGE
\noindent \textbf{Detailed Contents}.\\
\Large

%\noindent \textbf{Chapter \ref{00_sets} : Sets, Sets operations and applications}.\\

\begin{table}[htbp]
	\centering
		\caption{\textbf{Chapter \ref{00_sets} : Sets, Sets operations and applications}}
		\begin{tabular}{llll}
		\hline
		Type& Name & Title  & page\\
		S& Doc 00-01 &  Sets operations - A summary  & \pageref{doc00-01}\\
		D& Doc 00-02 & Exercises on sets operations   & \pageref{doc00-02} \\
		SD& Doc00-03  &Exercises on sets operations \textit{with solutions}   & \pageref{doc00-03} \\
		\hline
		\end{tabular}
\end{table} 

\vskip 0.5cm

%\noindent \textbf{Chapter \ref{01_setsmes} : Measurable Sets}.\\ 

\begin{table}[htbp]
	\centering
		\caption{\textbf{Chapter \ref{01_setsmes} : Measurable Sets}}
		\begin{tabular}{llll}
		\hline
		Type& Name & Title  & page\\
		\hline
		S & Doc 01-01 & Measurable sets - An outline  & \pageref{doc01-01}\\
		D &  Doc 01-02& Exercises on Measurable sets  & \pageref{doc01-02} \\
		SD& Doc 01-03 &    Exercises on Measurable sets with Solutions& \pageref{doc01-03}\\
		D& Doc 01-04 & Measurable sets. Exercises on $\lambda$ and $\pi$ systems& \pageref{doc01-04}\\
		SD& Doc 01-05 &   Measurable sets. Exercises on $\lambda$ and $\pi$ systems with solutions  & \pageref{doc01-05} \\
		D& Doc 01-05 & More  Exercises  & \pageref{doc01-06}\\
		SD& Doc 01-07 & More  Exercises - with solutions  & \pageref{doc01-07}\\
		\hline
		\end{tabular}
\end{table} 

\vskip 0.5cm

%\noindent \textbf{Chapter  \ref{02_applimess} : Measurable applications}.\\

\begin{table}[htbp]
	\centering
		\caption{\textbf{Chapter  \ref{02_applimess} : Measurable applications}}
		\begin{tabular}{llll}
		\hline
		Type& Name & Title  & page\\
		\hline
		S & Doc 02-01 & Measurable Applications - A summary& \pageref{doc02-01}\\
		S & Doc 02-02 & What you cannot ignore on Measurability  & \pageref{doc02-02}\\
		S& Doc 02-03 & Exercises on Measurable Applications   & \pageref{doc02-03}\\
		D& Doc 03-04 & Exercises on Measurable Applications with solutions& \pageref{doc02-04}\\
		\hline
		\end{tabular}
\end{table} 

\vskip 0.5cm

%\noindent \textbf{Chapter \ref{03_setsmes_applimes_cas_speciaux}: Measurability is Usual Spaces}.\\

\begin{table}[htbp]
	\centering
		\caption{\textbf{Chapter \ref{03_setsmes_applimes_cas_speciaux}: Measurability is Usual Spaces}}
		\begin{tabular}{llll}
		\hline
		Type& Name & Title  & page\\
		\hline
		S & Doc 03-01 & Borel sets  and $\mathbb{R}^d$-valued Measurable functions& \pageref{doc03-01}\\
		S & Doc 03-02 & Elementary functions  & \pageref{doc03-02} \\
		S& Doc 03-03 & Measurable applications - A summary  & \pageref{doc03-03}\\
		D& Doc 03-04 & Discover Exercises on Measurability 02 01& \\
		&   & Borel sets and Borel functions.   &  \pageref{doc03-04}\\
		D& Doc 03-05 & Discover Exercises on Measurability 02& \\
		&   & Elementary functions, Marginal functions and applications   &  \pageref{doc03-05}\\
		D& Doc 03-06 & Discover Exercises on Measurability 03& \\
		&   & Classical functions in Real Analysis  &  \pageref{doc03-06}\\
		SD& Doc 03-07 & Discover Exercises on Measurability 01 with solutions & \pageref{doc03-07}\\
		SD& Doc 03-08 & Discover Exercises on Measurability 02 with solutions & \pageref{doc03-08}\\
		SD& Doc 03-09 & Discover Exercises on Measurability 03 with solutions & \pageref{doc03-09}\\
		\hline
		\end{tabular}
\end{table} 

\vskip 0.5cm

%\noindent \textbf{Chapter  \ref{04_measures} : Measures}.\\

\begin{table}[htbp]
	\centering
		\caption{\textbf{Chapter  \ref{04_measures} : Measures}}
		\begin{tabular}{llll}
		\hline
		Type& Name & Title  & page\\
		\hline
		S&Doc 04-01  &Measures - A summary   &  \pageref{doc04-01}\\
		D&Doc 04-02  & Introduction to Measures - Exercises  & \pageref{doc04-02}\\
		S&Doc 04-03  & Carath\'{e}odory's Theorem - Lebesgue-Stieljes Measures   & \pageref{doc04-03}\\
		D&Doc 04-04  & Exercises on Lebesgue-Stieljes Measure  & \pageref{doc04-04}\\
		D&Doc 04-05  & General exercises on Measures  & \pageref{doc04-05}\\
		SD&Doc 04-06 & Introduction to Measures - Exercises with solutions & \pageref{doc04-06}\\
		SD&Doc 04-07 & Exercises on Lebesgue-Stieljes Measure with solutions & \pageref{doc04-07}\\
		SD&Doc 04-08 & General exercises on Measures with solutions  & \pageref{doc04-08}\\
		A&Doc 04-09  & The existence of the Lebesgue-Stieljes measure  & \pageref{doc04-09}\\
		A&Doc 04-10  & Proof the Carath\'{e}odory's Theorem  & \pageref{doc04-10}\\
		A&Doc 04-11  & Application of exterior measures to Lusin's theorem & \pageref{doc04-11}\\
		\hline
		\end{tabular}
\end{table} 

\vskip 0.5cm

%\noindent \textbf{Chapter \ref{05_integration} : Integration}.\\ 

\begin{table}[htbp]
	\centering
		\caption{\textbf{Chapter \ref{05_integration} : Integration}}
		\begin{tabular}{llll}
		\hline
		Type& Name & Title  & page\\
		\hline
		S & Doc 05-01 & Integration with respect to a measure - A summary   & \pageref{doc05-01}\\
		S & Doc 05-02 & Integrals of elementary functions - Exercises & \pageref{doc05-02}\\
		D & Doc 05-03 & Integration with the counting measure - Exercises  & \pageref{doc05-03}\\
		D & Doc 05-04 & Lebesgue/Riemann-Stieljes integrals on $\mathbb{R}$  - Exercises & \pageref{doc05-04}\\
		SD& Doc 05-05 & Integrals of elementary functions - Exercises with solutions & \pageref{doc05-05}\\
		SD& Doc 05-06 & Integration with the counting measure - Solutions & \pageref{doc05-06}\\
		SD & Doc 05-07 & Lebesgue/Riemann-Stieljes integrals on $\mathbb{R}$ - Solutions  & \pageref{doc05-07}\\
		A& Doc 05-08 &  Technical document on the definition of the integral & \\
		&           &  of a function of constant sign & \pageref{doc05-08}\\
		\hline
		\end{tabular}
\end{table} 

\vskip 0.5cm

%\noindent \textbf{Chapter \ref{06_convergence} : Introduction to Convergence Theorems and Applications}.\\ 

\begin{table}[htbp]
	\centering
		\caption{\textbf{Chapter \ref{06_convergence} : Introduction to Convergence Theorems and Applications}}
		\begin{tabular}{llll}
		\hline
		Type& Name & Title  & page\\
		S& Doc 06-01 &  The two Main Convergence Types - A summary & \pageref{doc06-01}\\
		S& Doc 06-02 & The two Main Convergence Theorems   & \pageref{doc06-02} \\
		S& Doc 06-03  & First Applications of The Convergence Theorems & \pageref{doc06-03} \\
		D& Doc 06-04  & Convergence Types - Exercises & \pageref{doc06-04} \\
		D& Doc 00-05  & Convergence Theorems - Exercises & \pageref{doc06-05} \\
		D& Doc 00-06  & Applications of Convergence Theorems - Exercises & \pageref{doc06-06} \\
		D& Doc 00-07  & Lebesgue-Stieljes/Riemann-Stieljes integrals - Exercises & \pageref{doc06-07} \\
		SD& Doc 06-08  & Convergence Types - Exercises with Solutions& \pageref{doc06-08} \\
		SD& Doc 06-09  & Convergence Theorems - Exercises with Solutions & \pageref{doc06-09} \\
		SD& Doc 06-10  & Applications of Convergence Theorems - Solutions & \pageref{doc06-10} \\
		SD& Doc 06-11  & Lebesgue-Stieljes/Riemann-Stieljes integrals - Solutions & \pageref{doc06-11} \\
		\hline
		\end{tabular}
\end{table} 

\vskip 0.5cm

%\noindent \textbf{Chapter \ref{07_finiteProductMeasure}  : Finite Product Measure Space and Fubini's Theorem}.\\

\begin{table}[htbp]
	\centering
		\caption{\textbf{Chapter \ref{07_finiteProductMeasure}  : Finite Product Measure Space and Fubini's Theorem}}
		\begin{tabular}{llll}
		\hline
		Type& Name & Title  & page\\
		S& 	Doc 07-01 & Finite Product Measures and Fubini's  Theorem - A summary  & \pageref{doc07-01}\\
		D& 	Doc 07-02 & Finite Product Measures - Exercises & \pageref{doc07-02} \\
		SD& Doc 07-03  & Finite Product Measures - Exercises with Solutions & \pageref{doc07-03} \\
				\hline
		\end{tabular}
\end{table} 

\vskip 0.5cm

%\noindent \textbf{Chapter \ref{08_nikodym} : $\sigma$-additive sets applications and Radon-Nikodym's Theorem}.\\

\begin{table}[htbp]
	\centering
		\caption{\textbf{Chapter \ref{08_nikodym} : $\sigma$-additive sets applications and Radon-Nikodym's Theorem}}
		\begin{tabular}{llll}
		\hline
		Type& Name & Title  & page\\
		S& Doc 08-01 &  Indefinite Integrals and Radon-Nikodym's Theorem & \pageref{doc08-01}\\
		D& Doc 08-02 &  $\sigma$-additive applications - Exercises & \pageref{doc08-02}\\
		D& Doc 08-03 &  Radon-Nikodym's Theorem - Exercises & \pageref{doc08-03}\\
		SD& Doc 08-04 &  $\sigma$-additive applications - Exercises/Solutions & \pageref{doc08-04}\\
		SD& Doc 08-05 &  Radon-Nikodym's Theorem - Exercises/Solutions & \pageref{doc08-05}\\
		SD& Doc 08-06 &  Exercise on a continuous family of measures & \pageref{doc08-06}\\
		%SD& Doc 08-05 &  Exercise on a continuous family of measures/Solutions & \pageref{doc08-07}\\
		\hline
		\end{tabular}
\end{table} 

\vskip 0.5cm

%\noindent \textbf{Chapter \ref{09_lp} : Introduction to spaces $L^p$, $p\in[1,+\infty]$}.\\

\begin{table}[htbp]
	\centering
		\caption{\textbf{Chapter \ref{09_lp} : Introduction to spaces $L^p$, $p\in[1,+\infty]$}}
		\begin{tabular}{llll}
		\hline
		Type& Name & Title  & page\\
		S& Doc 09-01 & Spaces $L^p$ - A summary & \pageref{doc09-01}\\
		D& Doc 09-02 & $L^p$ spaces by Exercises   & \pageref{doc09-02} \\
		SD& Doc 09-03  & $L^p$ spaces by Exercises with Exercises & \pageref{doc09-03} \\
		\hline
		\end{tabular}
\end{table} 

\vskip 0.5cm

%\noindent \textbf{Chapter  \ref{10_lsm} : Lebesgue-Stieljes measures on $\mathbb{R}^m$}.\\

\begin{table}[htbp]
	\centering
		\caption{\textbf{Chapter  \ref{10_lsm} : Lebesgue-Stieljes measures on $\mathbb{R}^m$}}
		\begin{tabular}{llll}
		\hline
		Type& Name & Title  & page\\
		A& Doc 10-01 &  Lebesgue-Stieljes measure and Probability Theory  & \pageref{doc10-01}\\
		A& 10-02 & Proof of the Existence of the Lebesgue-Stieljes measure on $\mathbb{R}^k$  & \pageref{doc10-02} \\
		A& Doc 10-03  & Lebesgue-Stieljes Integrals and Riemann-Stieljes Integrals & \pageref{doc10-03} \\
		\hline
		\end{tabular}
\end{table} 

\vskip 0.5cm

%\textbf{Chapter  \label{11_appendix} : General Appendix}\\

\begin{table}[htbp]
	\centering
		\caption{\textbf{Chapter  \label{11_appendix} : General Appendix}}
		\begin{tabular}{llll}
		\hline
		Type& Name & Title  & page\\
		SD& Doc 11-01 &  What should not be ignored on limits in $\overline{\mathbb{R}}$ & \pageref{doc11-01}\\
		\hline
		\end{tabular}
		\end{table} 
		
		\vskip 0.5cm

%\noindent \textbf{Chapter \ref{12_conbib} : Conclusion and Bibliography} 

\begin{table}[htbp]
	\centering
		\caption{\textbf{Chapter \ref{12_conbib} : Conclusion and Bibliography}}
		\begin{tabular}{llll}
		\hline
		Type& Name & Title  & page\\
		S& Doc 12-01 & Conclusion & \pageref{doc12-01}\\
		R& Doc 12-02 & Bibliography  & \pageref{doc12-02} \\
		\hline
		\end{tabular}
\end{table}

\part{Sets, Applications and Measurability}
%\chapter{Sets Operations}
\chapter{Sets, Sets operations and applications} \label{00_sets}

\noindent \textbf{Content of the Chapter}

\begin{table}[htbp]
	\centering
		\begin{tabular}{llll}
		\hline
		Type& Name & Title  & page\\
		S& Doc 00-01 &  Sets operations - A summary  & \pageref{doc00-01}\\
		D& Doc 00-02 & Exercises on sets operations   & \pageref{doc00-02} \\
		SD& Doc00-03  &Exercises on sets operations \textit{with solutions}   & \pageref{doc00-03} \\
		%&  &   & \\
		%&  &   & \\
		%&  &   & \\
		\hline
		\end{tabular}
\end{table}

\newpage
\noindent \LARGE \textbf{DOC 00-01 : Sets operations - A summary}. \label{doc00-01}\\

\Large

\bigskip

\noindent \textbf{NOTATION.}.\newline

\noindent $\Omega $ \ is a reference set. Its elements are denoted by $%
\omega $.\newline

\noindent $\emptyset$ is the empty set.\newline

\noindent $A$, $B$, $C$, .... are subsets of $\Omega $ (Capital normal letters).\newline

\noindent $\left( A_{n}\right) _{n\geq 0},$ $\left( B_{n}\right) _{n\geq 0}$, $\left(
C_{n}\right) _{n\geq 0}$ denote sequences of subsets of $\Omega .$\newline

\noindent $\mathcal{P}(\Omega $) is the power set of $\Omega $ : the collection of all
subsets of $\Omega .$\newline

\noindent $\mathcal{A}$, $\mathcal{B}$, $\mathcal{C}$, $\mathcal{D}$, $\mathcal{S}$,
..., denote collections of subsets of $\Omega $, that is subsets of $%
\mathcal{P}(\Omega )$ (calligraphic font).\newline

\bigskip

\noindent \textbf{I - BASIC DEFINITIONS}.\newline

\noindent The operations on subsets of $\Omega $ are the following.\newline

\bigskip

\noindent \textbf{(00.00)} Inclusion.\\

\noindent A set $A$ is contained in a set B or is equal to $B$, denoted $A \subseteq B$ or $A \subset B$ if and only if :

$$
\forall x \in A, \ x \in B.
$$

\bigskip \noindent A set $A$ is strictly contained in a set $B$, denoted $A \subsetneq B$ if and only if :

$$
\biggr(\forall x \in A, \ x \in B\biggr) \text{ and } \biggr(\exists x \in A, \ x \notin B\biggr)
$$

\bigskip \noindent or, equivalently,

$$
A \subseteq B  \text{ and } A \neq B.
$$

\bigskip \noindent \textbf{(00.01)} Complement of a subset $A$.

\begin{equation*}
A^{c}=\{x\in \Omega ,x\notin A\}
\end{equation*}

\bigskip \noindent \textbf{(00.02)} Finite union.\newline

\noindent (a) Union of two sets $A$ and $B$ :

\begin{equation*}
A\cup B=\{x\in \Omega ,x\in A\text{ or }x\in B\}.
\end{equation*}

\bigskip \noindent (b) Symmetry of the \textit{union} operator.

\begin{equation*}
A\cup B=B\cup A.
\end{equation*}

\bigskip \noindent (c) Associativity of the \textit{union} operator.

\begin{equation*}
(A\cup B)\cup C=A\cup (B\cup C)=A\cup B\cup C.
\end{equation*}

\bigskip \noindent (c) Union of a finite number $k$, $(k\geq 2)$ of sets%
\begin{equation*}
\bigcup\limits_{n=1}^{k}A_{n}=\{x\in \Omega ,(\exists n,1\leq n\leq k),x\in
A_{n}\}.
\end{equation*}

\bigskip \noindent \textbf{(00.03)} Finite intersection.\newline

\bigskip \noindent (a) Intersection of two sets $A$ and $B$.

\begin{equation*}
A\cap B=\{x\in \Omega ,x\in A\text{ and }x\in B\}.
\end{equation*}

\bigskip \noindent (b) Symmetry of the \textit{intersection} operator.%
\newline

\begin{equation*}
A\cap B=B\cap A.
\end{equation*}

\bigskip \noindent (c) Associativity of the \textit{intersection} operator.%
\newline

\begin{equation*}
(A\cap B)\cap C=A\cap (B\cap C)=A\cap B\cap C.
\end{equation*}

\bigskip \noindent (c) Intersection of a finite number $k$ $(k\geq 2)$ of sets

\begin{equation*}
\bigcap\limits_{n=1}^{k}A_{n}=\{x\in \Omega ,(\forall n, \ 1\leq n\leq k),x\in
A_{n}\}.
\end{equation*}

\bigskip \noindent \textbf{(00.04)} Sum of subsets : If $A$ and $B$ are
disjoint, that is $A\cap B=\varnothing,$ we denote

\begin{equation*}
A\cup B=A+B.
\end{equation*}

\bigskip \noindent \textbf{(00.05)} Difference $A$ from $B$

\begin{equation*}
B\text{ }\diagdown \text{ }A=\{x\in \Omega ,x\in B\text{ and }x\notin
A\}=BA^{c}.
\end{equation*}

\bigskip \noindent \textbf{(00.06)} Symmetrical difference of A and B :

\begin{equation*}
A\text{ }\Delta \text{ }B=(B\text{ }\diagdown \text{ }A)+(A\text{ }\diagdown 
\text{ }B).
\end{equation*}

\newpage

\bigskip \noindent \textbf{(00.07)} Countable union (Union of a class of
subsets indexed by a subset of the set of natural integers $\mathbb{N}$).

\begin{equation*}
\bigcup\limits_{n\geq 0}A_{n}=\{x\in \Omega ,(\exists n\geq 0),x\in A_{n}\}.
\end{equation*}

\bigskip \noindent \textbf{(00.08)} Countable sum of subsets : If the elements
of the sequence $(A)_{n\geq 0}$ are mutually disjoint, we write

\begin{equation*}
\bigcup\limits_{n\geq 0}A_{n}=\sum_{n\geq 0}A_{n}.
\end{equation*}

\bigskip \noindent \textbf{(00.09)} Countable intersection (Union of a class
of subsets indexed by a subset of the set of natural integers $\mathbb{N}$) :

\begin{equation*}
\bigcap\limits_{n\geq 0}A_{n}=\{x\in \Omega ,(\forall n\geq 0),x\in A_{n}\}.
\end{equation*}

\bigskip \bigskip

\noindent \textbf{II : MAIN PROPERTIES OF SETS OPERATIONS}.

\bigskip

\noindent \textbf{(00.10)} Complementation.

\begin{equation*}
\Omega ^{c}=\varnothing \text{, }\varnothing ^{c}=\Omega ,\text{ }%
(A^{c})^{c}=A;
\end{equation*}%
\begin{equation*}
(A\cup B)^{c}=A^{c}\cap B^{c},\text{ \ }(A\cap B)^{c}=A^{c}\cup B^{c};
\end{equation*}%
\begin{equation*}
\left( \bigcup\limits_{n\geq 0}A_{n}\right) ^{c}=\bigcap\limits_{n\geq
0}A_{n}^{c},\text{ \ }\left( \bigcap\limits_{n\geq 0}A_{n}\right)
^{c}=\bigcup\limits_{n\geq 0}A_{n}^{c}.
\end{equation*}

\bigskip \noindent \textbf{(00.11)} Distribution between union over intersection and vice-versa.\\

\begin{equation*}
A\cap (B\cup C)=(A\cap B)\cup (A\cap C),\text{ \ \ }A\cup (B\cap C)=(A\cup
B)\cap (A\cup C);
\end{equation*}

\begin{equation*}
A\bigcap \left( \bigcup\limits_{n\geq 0}B_{n}\right) =\bigcup\limits_{n\geq
0}\left( A\bigcap B_{n}\right);
\end{equation*}

\begin{equation*}
A\bigcup \left( \bigcap\limits_{n\geq 0}B_{n}\right) =\bigcap\limits_{n\geq
0}\left( A\bigcup B_{n}\right).
\end{equation*}

\bigskip \noindent \textbf{(00.12)} Difference between unions :
\begin{equation*}
\left( \bigcup\limits_{n\geq 0}A_{n}\right) \diagdown \left(
\bigcup\limits_{n\geq 0}B_{n}\right) \subseteq \bigcup\limits_{n\geq
0}(A_{n}\diagdown B_{n})
\end{equation*}%

\bigskip\noindent and

\begin{equation*}
\left( \bigcup\limits_{n\geq 0}A_{n}\right) \Delta \left(
\bigcup\limits_{n\geq 0}B_{n}\right) \subseteq \bigcup\limits_{n\geq
0}(A_{n}\Delta B_{n}).
\end{equation*}

\newpage \noindent \textbf{III : INDICATOR FUNCTIONS OF SUBSETS}.\newline

\bigskip

\noindent \textbf{(00.13)} The indicator function of the subset $A$ is a
function defined on $\Omega $ with values in $\mathbb{R}$ such that

\begin{equation*}
\forall \omega \in \Omega ,1_{A}(\omega )=\left\{ 
\begin{tabular}{lll}
$1$ & if & $\omega \in A$ \\ 
$0$ & if & $\omega \notin A$%
\end{tabular}%
\right. .
\end{equation*}

\bigskip \noindent \textbf{(00.14)} Properties :

\begin{equation*}
1_{\Omega }=1,\text{ }1_{\varnothing }=0,
\end{equation*}

\begin{equation*}
1_{A^{c}}=1-1_{A},
\end{equation*}

\begin{equation*}
1_{AB}=1_{A}\times 1_{B}, \ \ \  (PP)
\end{equation*}

\bigskip \noindent and

\begin{equation*}
1_{A\cup B}=1_{A}+1_{B}-1_{AB}.
\end{equation*}

\bigskip \noindent If  $A$ and $B$ are disjoint, we have

\begin{equation*}
1_{A\cup B}=1_{A}+1_{B}. \ \ (PS1)
\end{equation*}

\bigskip \noindent If the elements of the sequence $(A)_{n\geq 0}$ are mutually disjoint 
\begin{equation*}
1_{\bigcup_{n\geq 0}A_{n}}=1_{\sum_{n\geq 0}A_{n}}=\sum_{n\geq 0}1_{A_{n}}. \ \ (PS2)
\end{equation*}

\bigskip \noindent We have the Poincarr\'e's Formulas :

\begin{eqnarray*}
&&1_{\bigcup_{1\leq j\leq n}A_{j}} \notag \\
&=&\sum_{1\leq j\leq n} 1_{A_{j}}+\sum_{r=2}^{n}(-1)^{r+1} \sum_{1\leq i_{1}<...<t_{r}\leq
n}1_{A_{i_{1}}...A_{ir}} \ \ (HP)
\end{eqnarray*}

\noindent and

\begin{equation*}
Card(\bigcup_{1\leq j\leq n}A_n)=\sum_{1\leq j\leq n}Card(A_{j})+\sum_{r=2}^{n}(-1)^{r+1} \sum_{1\leq i_{1}<...<t_{r}\leq n} Card(A_{i_{1}}...A_{ir}). \ (HP)
\end{equation*}

\bigskip \noindent \textbf{NB} Because of Formula $(PP)$, we may write, if no confusion is possible, 

$$
A \bigcap B = AB.
$$

\bigskip 
\noindent \textbf{Warning}. For a countable intersection, we do not use the symbol $\Pi$ to avoid the confusion with the product of sets.\\

\bigskip \noindent \textbf{NB}. Because of Formulas $(PS1)$ and $(PS2)$, we may write, if $A$ and $B$ are disjoint, or if the $A_n$ are mutually disjoint.

$$
A \bigcup B = A+B.
$$

\bigskip 
\noindent and

$$
\bigcup_{n\geq 0} A_n = \sum_{n\geq 0} A_n.
$$

\newpage

\noindent \textbf{IV - LIMITS OF SEQUENCES OF SUBSETS}

\bigskip

\noindent \textbf{(00.15)} The inferior limit and superior limit of the
sequence of subsets $(A)_{n\geq 0}$ are subsets of $\Omega $ defined as
follows :%
\begin{equation*}
\underline{\lim }A_{n}=\liminf_{n\rightarrow +\infty}
=\bigcup\limits_{n\geq 0}\left( \bigcap\limits_{m\geq n}A_{m}\right)
\end{equation*}%
\begin{equation*}
\overline{\lim }\text{ }A_{n}=\limsup_{n\rightarrow +\infty }
=\bigcap\limits_{n\geq 0}\left( \bigcup\limits_{m\geq n}A_{m}\right)
\end{equation*}

\bigskip \noindent \textbf{(00.16)} Properties 
\begin{equation*}
\underline{\lim }A_{n}\subseteq \text{ }\overline{\lim }\text{ }A_{n}
\end{equation*}%
\begin{equation*}
\left( \underline{\lim }A_{n}\right) ^{c}=\text{ }\overline{\lim }\text{ }%
A_{n}^{c},\text{ \ \ \ }\left( \overline{\lim }\text{ }A_{n}\right) ^{c}=%
\underline{\lim }A_{n}^{c}\text{ }
\end{equation*}

\bigskip

\noindent \textbf{(00.17)} Criteria.

\begin{equation*}
\omega \in \text{ }\overline{\lim }\text{ }A_{n}\Leftrightarrow \text{ }%
\sum_{n\geq 0}1_{A_{n}}(\omega)=+\infty
\end{equation*}

\bigskip \noindent that is : $\ \omega \in $ $\overline{\lim }$ $A_{n}$ if
and only if $\omega $ belongs to an infinite number of $A_{n},$ ($\omega \in
A_{n}$ infinitely often), also abbreviated in $(\omega \in A_{n}$ $i.o).$

\begin{equation*}
\omega \in \text{ } \underline{\lim} A_{n}\Leftrightarrow \text{ }%
\sum_{n\geq 0}1_{A_{n}^{c}}(\omega)<+\infty
\end{equation*}

\bigskip \noindent that is : $\ \omega \in \underline{\lim} A_{n}$ if and
only if $\omega $ belongs, at most, to a finite number of $A_{n}^{c},$ ($%
\omega \in A_{n}^{c}$ finitely often), also abbreviated in $(\omega \in
A_{n}^{c}$ $f.o).$

\bigskip \noindent \textbf{(00.18)} Limits of monotone sequences.\newline

\noindent (a) If the sequence sequence $(A)_{n\geq 0}$ is non-decreasing, then

\begin{equation*}
\underline{\lim }A_{n}=\overline{\lim }A_{n}=\bigcup_{n\geq 0}A_{n}
\end{equation*}

\noindent and this common set is the limit of the sequence :

\begin{equation*}
\lim_{n\rightarrow +\infty} A_n = \bigcup\limits_{n\geq 0} A_{n}
\end{equation*}

\noindent (b) If the sequence $(A)_{n\geq 0}$ is non-increasing, then

\begin{equation*}
\underline{\lim} A_{n}=\overline{\lim} A_{n}=\bigcap_{n\geq 0} A_{n}
\end{equation*}

\noindent and this common set is the limit of the sequence :

\begin{equation*}
\lim_{n\rightarrow +\infty }A_{n}=\bigcap\limits_{n\geq 0}A_{n}
\end{equation*}

\bigskip \noindent \textbf{(00.19)} Limits of sequences of subsets.\\

\noindent (a) Limits of monotone subsets are special cases of the following.\\

\bigskip \noindent (b) A sequence sequence $(A)_{n\geq 0}$ of subsets of $%
\Omega $ has a limit if and only if its superior limit is equal to its
inferior limit and this common set is called the limit of the sequence :

\begin{equation*}
\lim_{n\rightarrow \infty }A_{n}=\underline{\lim }A_{n}=\overline{\lim }%
A_{n}=\bigcup_{n\geq 0}A_{n}
\end{equation*}

\newpage

\bigskip \noindent \textbf{V - product spaces}.\newline

\bigskip \noindent Let $\Omega _{1}$ and $\Omega _{2}$ be two non-empty
sets. The product space $\Omega =\Omega _{1}$ $\times \Omega _{2}$ is the
set of all couples $(\omega _{1},\omega _{2}),$ where $\omega _{1}\in \Omega
_{1}$ and $\omega _{2}\in \Omega _{2}$, that is%
\begin{equation*}
\Omega =\Omega _{1}\times \Omega _{2}=\{(\omega _{1},\omega _{2}),\omega
_{1}\in \Omega _{1}\text{ and }\omega _{2}\in \Omega _{2}\}.
\end{equation*}

\bigskip \noindent A subset of the product space is a rectangle if it is the
product of two subsets of $\Omega _{1}$ and $\Omega _{2},$ that is%
\begin{equation*}
A\times B\text{ \ with }A\subseteq \Omega _{1}\text{, }B\subseteq \Omega _{2}
\end{equation*}%
The rectangle $A\times B$ is defined as%
\begin{equation*}
A\times B=\{(\omega _{1},\omega _{2})\in \Omega ,\text{ }\ \omega _{2}\in A%
\text{ and }\omega _{2}\in B\}
\end{equation*}

\bigskip \noindent \textbf{(00.20)} Not every subset of $\Omega $ is a
rectangle (Unfortunately, taking all subsets of a product space as
rectangles is a very recurrent mistake). Here is a counterexample.\newline

\bigskip \noindent Let $\Omega =\Omega _{0}\times \Omega _{0}=\Omega _{0}^{2}
$ be the product of one space $\Omega _{0}$ by itself and consider its
diagonal%
\begin{equation*}
D=\{(\omega _{1},\omega _{2})\in \Omega ,\text{ }\ \omega _{1}=\omega _{2}\}.
\end{equation*}

\bigskip \noindent Suppose that $\Omega _{0}$ has at least two elements. So $%
D$ is not a rectangle. Hint : To see that, suppose that $D=A\times B.$ Show
that this is not possible by considering two cases : (1) Both $A$ and $B$
are singletons, (2) $A$ or $B$ has at least two elements.

\bigskip \noindent \textbf{(00.21) Properties of rectangles}.\newline

\begin{equation*}
A\times B=\emptyset \text{ if and only if }A=\emptyset \text{ or }%
B=\emptyset \text{ }
\end{equation*}%
\begin{equation*}
(A\times B)\cap (A^{\prime }\times B^{\prime })=\left( A\cap A^{\prime
}\right) \times \left( B\cap B^{\prime }\right)
\end{equation*}%
\begin{equation*}
\left( A\times B\right) ^{c}=A^{c}\times B+A\times B^{c}+A^{c}\times B^{c}
\end{equation*}

\newpage \bigskip \noindent \textbf{0.22 Product of }$k$\textbf{\ spaces}.\\

\noindent Let $\Omega _{1},...,\Omega _{k}$ be $k\geq 2$ non-empty sets. Their product
space is defined and denoted by 
\begin{equation*}
\Omega =\prod\limits_{j=1}^{k}\Omega _{j}=\{(\omega _{1},\omega
_{2},...,\omega _{k}),\forall (1\leq j\leq k),\omega _{j}\in \Omega _{j}\}.
\end{equation*}

\bigskip \noindent A cylinder (we may still call them rectangles) is a
subset de $\Omega $ of the form 
\begin{equation*}
\prod\limits_{j=1}^{k}A_{j}=\{(\omega _{1},\omega _{2},...,\omega _{k})\in
\Omega ,\forall (1\leq j\leq k),\omega _{j}\in A_{j}\}
\end{equation*}

\bigskip \noindent where each $A_{j}$ is a subset of $\Omega _{j}$ for $%
j=1,...,k.$ We have the same properties 

\begin{equation*}
\prod\limits_{j=1}^{k}A_{j}=\emptyset \text{ if and only if }\left( \exists
1\leq j\leq k\right) ,A_{j}=\emptyset.
\end{equation*}%

\begin{equation}
\left( \prod\limits_{j=1}^{k}A_{j}\right) \cap \left(
\prod\limits_{j=1}^{k}B_{j}\right) =\prod\limits_{j=1}^{k}(A_{j}\cap B_{j}).
\tag{IREC}  \label{mes_01_intesect_rect}
\end{equation}

\begin{equation*}
\left( \prod\limits_{j=1}^{k}A_{j}\right) ^{c}=\sum\limits_{\substack{ %
\varepsilon =(\varepsilon _{1},\varepsilon _{2},...,\varepsilon _{k})\in
\{0,1\}^{k} \\ \varepsilon _{1}+\varepsilon _{2}+...+\varepsilon _{k}>0}}%
\prod\limits_{j=1}^{k}A_{j}^{(\varepsilon _{j})} \tag{CPS}
\end{equation*}

\bigskip \noindent where $A_{j}^{(\varepsilon _{j})}=A_{j}$ if \ $\varepsilon _{j}=0$\ and $A_{j}^{(\varepsilon _{j})}=A_{j}^{c}$ if $\varepsilon _{j}=1$.\\

\newpage
\noindent \LARGE \textbf{DOC 00-02 : Exercises on sets operations}. \label{doc00-02}\\

\Large

\bigskip

\bigskip \noindent \textbf{Exercise 1}.  \label{exercise01_doc00-02} Prove the assertions in \textit{(00.12)}.%
\newline

\begin{equation}
\left( \bigcup\limits_{n\geq 0}A_{n}\right) \diagdown \left(
\bigcup\limits_{n\geq 0}B_{n}\right) \subseteq \bigcup\limits_{n\geq
0}(A_{n}\diagdown B_{n}).  \tag{D1}  \label{D1}
\end{equation}%
\begin{equation}
\left( \bigcup\limits_{n\geq 0}A_{n}\right) \Delta \left(
\bigcup\limits_{n\geq 0}B_{n}\right) \subseteq \bigcup\limits_{n\geq
0}(A_{n}\Delta B_{n}).  \tag{D2}
\end{equation}

\bigskip

\bigskip \noindent \textbf{Exercise 2}. \label{exercise02_doc00-02}  Prove the assertions in \textit{(00.16)}.%
\newline

\begin{equation}
\underline{\lim }A_{n}\subseteq \text{ }\overline{\lim }\text{ }A_{n}.  \tag{L1}
\end{equation}%
\begin{equation}
\left( \underline{\lim }A_{n}\right) ^{c}=\text{ }\overline{\lim }\text{ }%
A_{n}^{c},\text{ \ \ \ }\left( \overline{\lim }\text{ }A_{n}\right) ^{c}=%
\underline{\lim }A_{n}^{c}\text{ }  \tag{L2}.
\end{equation}

\bigskip

\bigskip \noindent \textbf{Exercise 3}. \label{exercise03_doc00-02} Prove the assertions in \textit{(00.17)} Prove the assertions in Point \textit{(00.17)} in \textit{Doc 00-01}, by
using \textit{(00.15)} and \textit{(00.16)} : 

\begin{equation}
\omega \in \text{ }\overline{\lim }\text{ }A_{n}\Leftrightarrow \text{ }%
\sum_{n\geq 0}1_{A_{n}}(\omega )=+\infty,  \tag{L3}
\end{equation}

\bigskip \noindent that is : $\omega \in $ $\overline{\lim }$ $A_{n}$ if and
only if $\omega $ belongs to an infinite number of $A_{n},$ ($\omega \in
A_{n}$ infinitely often), also abbreviated in $(\omega \in A_{n}$ $i.o).$

\begin{equation}
\omega \in \text{ }\underline{\lim }A_{n}\Leftrightarrow \text{ }\sum_{n\geq
0}1_{A_{n}^{c}}(\omega )<+\infty.  \tag{L4}
\end{equation}

\bigskip

\bigskip \noindent \textbf{Exercise 4}. \label{exercise4_doc00-02} Useful tool for transforming a union
of subsets into a sum of subsets. Let ${A_n,\ \ n\geq 0}$ be a sequence of
subsets of $\Omega$. Set 
\begin{equation*}
B_0=A_0, B_1=A_{0}^{c}A_{1}, \text{ } B_2= A_{0}^{c}A_{1}^{c}A_{2}, \text{ }%
...., \text{ } B_k= A_{0}^{c}A_{1}^{c}...A_{k-1}^{c}A_{k}, \text{ } k\geq 2. 
\end{equation*}

\bigskip \noindent \textbf{(a)} Show that

\begin{equation*}
\bigcup\limits_{n\geq 0} A_n = \sum_{k\geq 0} B_k. 
\end{equation*}

\bigskip \noindent \textbf{(b)} Interpret the meaning of $k$ for $\omega \in
B_k$, for any $\omega \in \bigcup\limits_{n\geq 0} A_n$.

\bigskip

\bigskip \noindent \textbf{Exercise 5}. \label{exercise05_doc00-02}\newline

\bigskip \noindent (a) Generalize the formula in \textit{(00.14)}

\begin{equation*}
1_{A\cup B}=1_{A}+1_{B}-1_{A\cap B}
\end{equation*}

\bigskip \noindent to any countable union. Precisely, by applying it to three sets and next to four sets, guess the general formula :

\begin{eqnarray*}
1_{\bigcup_{1\leq j\leq n}A_{j}} \notag \\
&=&\sum_{1\leq j\leq n} 1_{A_{j}}+\sum_{r=2}^{n}(-1)^{r+1} \sum_{1\leq i_{1}<...<t_{r}\leq
n}1_{A_{i_{1}}...A_{ir}}. \ \ (HP)
\end{eqnarray*}

\bigskip \noindent (b) From the formula

$$
Card( A \bigcup B ) =Card(A) + Card(B)- Card(A \bigcap B),
$$

\bigskip 
\noindent proceed similarly and justify the formula 

\begin{equation*}
Card(\bigcup_{1\leq j\leq n}A_j)=\sum_{1\leq j\leq n}Card(A_{j})+\sum_{r=2}^{n}(-1)^{r+1} \sum_{1\leq i_{1}<...<t_{r}\leq n} Card(A_{i_{1}}...A_{ir}). \ (HP)
\end{equation*}

\bigskip \noindent Both Formulas are versions of the Pointcarr\'e's Formula.\\

\noindent \textbf{NB} It is only asked to guess the formulas by applying the two-set formula for three sets and next for four sets and by using the associativity of the union operation on sets. But a proof of these formula will be given in the solutions.\\

\bigskip \noindent \textbf{Exercise 6}. \label{exercise06_doc00-02}\newline

\bigskip \noindent Show that :

\begin{equation*}
1_{\bigcup\limits_{n\geq 0} A_n} = \sup_{n\geq 0} 1_{A_n}. 
\end{equation*}

\bigskip 
\noindent and

\begin{equation*}
1_{\bigcap\limits_{n\geq 0}A_{n}}=\inf_{n\geq 0}1_{A_{n}}.
\end{equation*}

\bigskip \noindent \textbf{Exercise 7}. \label{exercise07_doc00-02} Let $\Omega =\Omega _{1}$ $\times \Omega _{2}$ be a non-empty
product defined by%
\begin{equation*}
\Omega =\Omega _{1}\times \Omega _{2}=\{(\omega _{1},\omega _{2}),\omega
_{1}\in \Omega _{1}\text{ and }\omega _{2}\in \Omega _{2}\}.
\end{equation*}

\bigskip \noindent Define a rectangle%
\begin{equation*}
A\times B=\{(\omega _{1},\omega _{2})\in \Omega ,\ \omega _{2}\in and\omega
_{2}\in B\},
\end{equation*}%

\bigskip 
\noindent where $A\subseteq \Omega _{1}$, $B\subseteq \Omega _{2}$.\\

\noindent (a) Show that the intersection of rectangles is a rectangle.\\

\noindent (b) Let $A\times B$ be a rectangle. Show that \\
\begin{equation*}
(A\times B)^{c}=A^{c}\times B+A\times B^{c}+A^{c}\times B^{c}.
\end{equation*}

\bigskip 
\noindent (c) Show that not all subsets of are rectangles. Use the following example.\\

\noindent Let $\Omega _{0}=\Omega _{1}=\Omega _{2}$ so that $\Omega =\Omega _{0}\times
\Omega _{0}=\Omega _{0}^{2}$ is the product of one space $\Omega _{0}$ by
itself and consider its diagonal
\begin{equation*}
D=\{(\omega _{1},\omega _{2})\in \Omega ,\text{ }\ \omega _{1}=\omega _{2}\}.
\end{equation*}

\bigskip \noindent Suppose that $\Omega _{0}$ has at least two elements.
Suppose that $D=A\times B.$ Show that this is not possible by considering
two cases : (1) Both $A$ and $B$ are singletons, (2) $A$ or $B$ has at least
two elements.

\newpage
\noindent \LARGE \textbf{DOC 00-03 : Exercises on sets operations \textit{with solutions}}.\label{doc00-03}\\

\Large

\bigskip

\bigskip \noindent \textbf{Exercise 1}.  \label{exercise01_sol_doc00-03} Prove the assertions in Point \textit{(00.12)} in Doc 00-01 : 
\begin{equation}
\left( \bigcup\limits_{n\geq 0}A_{n}\right) \diagdown \left(
\bigcup\limits_{n\geq 0}B_{n}\right) \subseteq \bigcup\limits_{n\geq
0}(A_{n}\diagdown B_{n}).  \tag{D1}  \label{D1}
\end{equation}%
\begin{equation}
\left( \bigcup\limits_{n\geq 0}A_{n}\right) \Delta \left(
\bigcup\limits_{n\geq 0}B_{n}\right) \subseteq \bigcup\limits_{n\geq
0}(A_{n}\Delta B_{n}).  \tag{D2}
\end{equation}

\bigskip \noindent \textbf{Solutions}. \\

\noindent \textbf{Proof of (D1)}.\newline

\bigskip \noindent \textbf{First method}. First remark that the notation of
index $n$ in $x\bigcup\limits_{n\geq 0}B_{n}$ does not matter and we may use
any other notation, like $\bigcup\limits_{p\geq 0}B_{p}.$ So 
\begin{eqnarray*}
\left( \bigcup\limits_{n\geq 0}A_{n}\right) \diagdown \left(
\bigcup\limits_{n\geq 0}B_{n}\right) &=&\left( \bigcup\limits_{n\geq
0}A_{n}\right) \bigcap \left( \bigcup\limits_{n\geq 0}B_{n}\right)
^{c}\\
&=&\left( \bigcup\limits_{n\geq 0}A_{n}\right) \bigcap \bigcap\limits_{p\geq
0}B_{p}^{c}.
\end{eqnarray*}

\bigskip \noindent Use the distributivity of the intersection over the union
to get%
\begin{equation}
\left( \bigcup\limits_{n\geq 0}A_{n}\right) \diagdown \left(
\bigcup\limits_{n\geq 0}B_{n}\right) =\bigcup\limits_{n\geq 0}\left\{
A_{n}\bigcap \left( \bigcap\limits_{p\geq 0}B_{p}^{c}\right) \right\}. 
\tag{D3}
\end{equation}

\bigskip \noindent Remark that, for $n$ fixed, we may use the inclusion $%
\bigcap\limits_{p\geq 0}B_{p}^{c}\subseteq B_{n}^{c}$ in $A_{n}\bigcap \left(
\bigcap\limits_{p\geq 0}B_{p}^{c}\right) $, to get%
\begin{equation*}
A_{n}\bigcap \left( \bigcap\limits_{p\geq 0}B_{p}^{c}\right) \subseteq 
A_{n}\bigcap B_{n}^{c}=A_{n}\setminus B_{n}.
\end{equation*}

\bigskip \noindent Plug-in this in (D3) to get%
\begin{equation*}
\left( \bigcup\limits_{n\geq 0}A_{n}\right) \diagdown \left(
\bigcup\limits_{n\geq 0}B_{n}\right) \subseteq \bigcup\limits_{n\geq 0}\left(
A_{n}\setminus B_{n}\right) .
\end{equation*}

\bigskip \noindent \textbf{Second method}. Use equivalences, that is 
\begin{equation*}
\omega \in \left( \bigcup\limits_{n\geq 0}A_{n}\right) \diagdown \left(
\bigcup\limits_{n\geq 0}B_{n}\right) \Longleftrightarrow \left( \omega \in
\bigcup\limits_{n\geq 0}A_{n}\text{ and }\omega \notin \bigcup\limits_{n\geq
0}B_{n}\right)
\end{equation*}
\begin{equation}
\Longleftrightarrow \left( \exists n\geq 0,\omega \in A_{n}\text{ and }%
\forall p\geq 0,\omega \notin B_{p}\right) .  \tag{D4}  \label{D}
\end{equation}

\bigskip \noindent So, if $\omega \in \left( \bigcup\limits_{n\geq
0}A_{n}\right) \diagdown \left( \bigcup\limits_{n\geq 0}B_{n}\right) ,$ we
fixed one $n$ such that $\omega \in A_{n}.$ Since $\omega \notin B_{p}$ for
all $p\geq 0$, then for the particular value of $n$, we also have $\omega
\notin B_{n}.$Thus (D4) implies

\begin{equation*}
\left( \exists n\geq 0),(\omega \in A_{n}\text{ and }\omega \notin
B_{n}\right) \Longleftrightarrow \omega \in \bigcup\limits_{n\geq
0}(A_{n}\setminus B_{n}).
\end{equation*}

\bigskip \noindent Hence we conclude that
\begin{equation*}
\left( \bigcup\limits_{n\geq 0}A_{n}\right) \diagdown \left(
\bigcup\limits_{n\geq 0}B_{n}\right) \subseteq \bigcup\limits_{n\geq
0}(A_{n}\setminus B_{n}).
\end{equation*}

\bigskip \noindent \textbf{Proof of (D2)}. By definition, we have

\begin{equation*}
\left( \bigcup\limits_{n\geq 0}A_{n}\right) \Delta \left(
\bigcup\limits_{n\geq 0}B_{n}\right) =\left\{ \left( \bigcup\limits_{n\geq
0}A_{n}\right) \diagdown \left( \bigcup\limits_{n\geq 0}B_{n}\right)
\right\} \bigcup \left\{ \left( \bigcup\limits_{n\geq 0}B_{n}\right)
\diagdown \left( \bigcup\limits_{n\geq 0}A_{n}\right) \right\} .
\end{equation*}

\bigskip \noindent By (D1)

\begin{eqnarray*}
\left( \bigcup\limits_{n\geq 0}A_{n}\right) \Delta \left(
\bigcup\limits_{n\geq 0}B_{n}\right) &\subseteq & \left\{ \bigcup\limits_{n\geq
0}(A_{n}\setminus B_{n})\right\} \bigcup \left\{ \bigcup\limits_{n\geq
0}(A_{n}\setminus B_{n})\right\}\\
& =&\bigcup\limits_{n\geq 0}\left\{(A_{n}\setminus B_{n})\bigcup (B_{n}\setminus A_{n})\right\}\\
&=&\bigcup\limits_{n\geq 0}(A_{n}\Delta \setminus B_{n}).
\end{eqnarray*}

\bigskip \noindent This completes the proof.

\bigskip

\bigskip \noindent \textbf{Exercise 2}. \label{exercise02_sol_doc00-03} Prove the assertions in Point \textit{(00.16)} in \textit{Doc 00-01} : 
\begin{equation}
\underline{\lim }A_{n}\subseteq \text{ }\overline{\lim }\text{ }A_{n}.  \tag{L1}
\end{equation}%
\begin{equation}
\left( \underline{\lim }A_{n}\right) ^{c}=\text{ }\overline{\lim }\text{ }%
A_{n}^{c},\text{ \ \ \ }\left( \overline{\lim }\text{ }A_{n}\right) ^{c}=%
\underline{\lim }A_{n}^{c}.\text{ }  \tag{L2}
\end{equation}

\bigskip \bigskip \textbf{Solutions}.\\

\noindent \textbf{Proof of L1}.\\

\bigskip \noindent We have, for each $n\geq 1$%
\begin{equation*}
\bigcap\limits_{n\geq m}A_{m}\subseteq \bigcup\limits_{m\geq n}A_{n}.
\end{equation*}

\bigskip \noindent \textbf{Proof of L2}. L2 has two expressions. They are
proved in the same way. So, we are going to prove the first expression only.
we have%
\begin{equation*}
\underline{\lim }A_{n}=\bigcup\limits_{n\geq 0}\bigcap\limits_{n\geq m}A_{m}.
\end{equation*}

\bigskip \noindent Thus, we have
\begin{equation*}
\left( \underline{\lim }A_{n}\right) ^{c}=\left( \bigcup\limits_{n\geq
0}\bigcap\limits_{n\geq m}A_{m}\right) ^{c}.
\end{equation*}

\bigskip \noindent Next, use the Morgan's law to get%
\begin{equation*}
\left( \underline{\lim }A_{n}\right) ^{c}=\bigcap\limits_{n\geq
0}\bigcup\limits_{m\geq n}A_{m}^{c}=\overline{\lim }\text{ }A_{n}^{c}.
\end{equation*}

\bigskip

\bigskip \noindent \textbf{Exercise 3}. \label{exercise03_sol_doc00-03} Prove the assertions in Point \textit{(00.17)} in \textit{Doc 00-01}, by
using \textit{(00.15)} and \textit{(00.16)} : 
\begin{equation}
\omega \in \text{ }\overline{\lim }\text{ }A_{n}\Leftrightarrow \ \
\sum_{n\geq 0}1_{A_{n}}(\omega )=+\infty.  \tag{L3}
\end{equation}

\bigskip \noindent that is : $\ \omega \in \overline{\lim } A_{n}$ if
and only if $\omega $ belongs to an infinite number of $A_{n},$ ($\omega \in
A_{n}$ infinitely often), also abbreviated in $(\omega \in A_{n}$ $i.o).$

\begin{equation}
\omega \in \text{ }\underline{\lim }A_{n}\Leftrightarrow \text{ }\sum_{n\geq
0}1_{A_{n}^{c}}(\omega )<+\infty  \tag{L4}
\end{equation}

\bigskip \bigskip \noindent \textbf{Solutions}.\\

\noindent \textbf{(A) Proof of (L3)}.\newline

\bigskip \noindent \textbf{(i) Prove the direct implication}

\begin{equation*}
\omega \in \overline{\lim }A_{n}\Longrightarrow \sum_{n\geq
0}1_{A_{n}}(\omega )=+\infty .
\end{equation*}

\bigskip \noindent Let $\omega \in \overline{\lim }A_{n}.$ We are going to
construct a sequence of indices $\left( n_{k}\right) _{k\geq 0}$ such that :
(i) $n_{0}<n_{1}<...<n_{k}<n_{k+1}<....$, \ and (ii) for any $k\geq 0,\omega
\in A_{n_{k}}$.

\noindent To do this, use the definition of $\overline{\lim }A_{n}$ : 
\begin{equation*}
\overline{\lim }A_{n}=\bigcap\limits_{n\geq 0}\bigcup\limits_{m\geq n}A_{m},
\end{equation*}

\bigskip \noindent to see that
\begin{equation*}
\text{For all }n\geq 0,\omega \in B_{n}=\bigcup\limits_{m\geq n}A_{m}.
\end{equation*}

\bigskip \noindent Begin by $n=0$, that is 
\begin{equation*}
\omega \in B_{n}=\bigcup\limits_{m\geq 0}A_{m}.
\end{equation*}

\bigskip \noindent Then there exists $n_{0}\geq 0$ such that%
\begin{equation*}
\omega \in A_{n_{0}}.
\end{equation*}

\bigskip \noindent Next, say that%
\begin{equation*}
\omega \in B_{n_{_{0}}+1}=\bigcup\limits_{m\geq n_{_{0}}+1}A_{m}.
\end{equation*}

\bigskip \noindent So there exists $n_{1}\geq n_{_{0}}+1$ such that%
\begin{equation*}
\omega \in A_{n_{1}}.
\end{equation*}

\bigskip \noindent Next use 
\begin{equation*}
\omega \in B_{n_{1}+1}=\bigcup\limits_{m\geq n_{1}+1}A_{m},
\end{equation*}

\bigskip \noindent to get $n_{2}\geq n_{_{1}}+1$ such that 
\begin{equation*}
\omega \in A_{n_{2}}
\end{equation*}

\bigskip \noindent And so forth. We get an infinite strictly increasing
sequence $\left( n_{k}\right) _{k\geq 0}$ such that : (1) for all $k\geq 0,$n%
$_{k+1}\geq n_{k}+1$ and (2) for all $k\geq 0,\omega \in A_{n_{k}}.$\newline

\bigskip \noindent Finally, since the series has non-negative terms, we get that
\begin{equation*}
\sum_{n\geq 0}1_{A_{n}}(\omega )\geq \sum_{k\geq 0}1_{A_{n_{k}}}(\omega
)=\sum_{k\geq 0}1=+\infty .
\end{equation*}

\bigskip \noindent \textbf{(ii) Prove the indirect implication}.\newline

\begin{equation*}
\sum_{n\geq 0}1_{A_{n}}(\omega )=+\infty \Longrightarrow \omega \in 
\overline{\lim }A_{n}.
\end{equation*}

\bigskip \noindent Let 
\begin{equation*}
\sum_{n\geq 0}1_{A_{n}}(\omega )=+\infty .
\end{equation*}

\bigskip \noindent Since the terms of the series, that is the $1_{A_{n}}(\omega ),$ are zero or one, if it is infinite, we necessarily that
an infinite number of the $1_{A_{n}}(\omega )$ are one. Denote these terms with the subscripts $n_{k}$ in a strictly increasing order. This implies
that the sequence $n_{k}$ is unbounded and
\begin{equation*}
\text{For all }k\geq 0,\omega \in A_{n_{k}}.
\end{equation*}

\bigskip \noindent Apply this to prove that $\omega \in \overline{\lim} A_{n}$, that is 
\begin{equation*}
\text{For all }n\geq 0,\omega \in B_{n}=\bigcup\limits_{m\geq n}A_{m}.
\end{equation*}

\bigskip \noindent Remind that the sequence $n_{k}$ is unbounded. So, for
all $n\geq 0,$ there exists $n_{k}$ such that 
\begin{equation*}
n\leq n_{k},
\end{equation*}

\bigskip \noindent and since $\omega \in A_{n_{k}},$ we get
\begin{equation*}
\omega \in A_{n_{k}}\subseteq \bigcup\limits_{m\geq n}A_{m},
\end{equation*}

\bigskip \noindent and then 
\begin{equation*}
\text{For all }n\geq 0,\omega \in B_{n}=\bigcup\limits_{m\geq n}A_{m},
\end{equation*}

\bigskip \noindent and hence 
\begin{equation*}
\omega \in \overline{\lim }A_{n}.
\end{equation*}

\bigskip \noindent We have proved (L3). And $(L3)$ means exactly that $%
\omega \in \overline{\lim }A_{n}$ if and only iff $\omega $ belongs to an
infinite number of $A_{n}$, abbreviated : $(\omega \in A_{n},i.o)$.\newline

\bigskip \noindent \textbf{(iii) Use (L3) to prove (L4)}.

\bigskip \noindent Let us use the complementation :%
\begin{equation*}
\omega \in \underline{\lim }A_{n}=\left( \underline{\lim }A_{n}^{c}\right)
^{c}\Longleftrightarrow \omega \notin \underline{\lim }A_{n}^{c}
\end{equation*}

\begin{equation*}
\Leftrightarrow Not\text{ }(\omega \in A_{n}^{c},i.o)\Longleftrightarrow
(\omega \in A_{n}^{c}\text{ finitely often)}
\end{equation*}%
\begin{equation*}
\Longleftrightarrow \sum_{n\geq 0}1_{A_{n}^{c}}(\omega )<+\infty .
\end{equation*}

\bigskip
\noindent \textbf{Exercise 4}. \label{exercise04_sol_doc00-03} Useful tool for transforming a union of
subsets into a sum of subsets. Let ${A_{n},\ \ n\geq 0}$ be a sequence of
subsets of $\Omega $. Set

\begin{equation*}
B_{0}=A_{0},B_{1}=A_{0}^{c}A_{1},B_{2}=A_{0}^{c}A_{1}^{c}A_{2},....,B_{k}=A_{0}^{c}A_{1}^{c}...A_{k-1}^{c}A_{k},k\geq 2.
\end{equation*}

\bigskip \noindent \textbf{(a)} Show that 
\begin{equation*}
\bigcup\limits_{n\geq 0} A_n = \sum_{k\geq 0} B_k.
\end{equation*}

\bigskip \noindent \textbf{(b)} Interpret the meaning of $k$ for $\omega \in B_{k}$, for any $\omega \in \bigcup\limits_{n\geq 0}A_{n}$.\\

\bigskip \noindent \textbf{Solution}.\newline

\noindent To prove 
\begin{equation*}
\bigcup\limits_{n\geq 0}A_{n}=\sum_{k\geq 0}B_{k}.
\end{equation*}

\bigskip \noindent We have to justify three things : the two unions $\bigcup\limits_{n\geq 0}B_{n}$ and $\bigcup\limits_{n\geq 0}A_{n}$ are equal
and $(B_{n})_{n\geq 0}$ is a sequence pairwise disjoint sets.\newline

\bigskip \noindent (a) We want to prove that : $\bigcup\limits_{n\geq
0}B_{n}\subseteq \bigcup\limits_{n\geq 0}A_{n}.$\newline

\bigskip \noindent Remark that for each $k\geq 0,$ 
\begin{equation*}
B_{k}=A_{0}^{c}A_{1}^{c}...A_{k-1}^{c}A_{k}\subseteq A_{k},
\end{equation*}

\bigskip since $B_{k}$\ is the intersection of $A_{k}$ and an other set. Then%
\begin{equation*}
(\forall k\geq 0,B_{k}\subseteq A_{k})\Longrightarrow \bigcup\limits_{k\geq
0}B_{k}\subseteq \bigcup\limits_{k\geq 0}A_{k}.
\end{equation*}

\bigskip \noindent Notice that in the expression $\bigcup\limits_{k\geq 0}B_{k}$, the notation of the subscript does not matter. So (a) is proved that is
\begin{equation}
\bigcup\limits_{n\geq 0}B_{n}\subseteq \bigcup\limits_{n\geq 0}A_{n}.  \tag{I-A}
\end{equation}

\bigskip \noindent (b) We want to prove that $\bigcup\limits_{n\geq 0}A_{n}\subseteq \bigcup\limits_{n\geq 0}B_{n}.$\newline

\bigskip \noindent Let $\omega \in \bigcup\limits_{n\geq 0}A_{n}$. Then there exists $n_{0}$ such that $\omega \in A_{n_{0}}.$ It equivalent to say
that $\omega \in \bigcup\limits_{n\geq 0}A_{n}$ if and only is the set
\begin{equation*}
N=\{n\in \mathbb{N},\omega \in A_{n}\}.
\end{equation*}

\bigskip \noindent is not empty. Denote by $k$ the minimum of the set $N$. This means that : (i) $k$ is member of the set, that is $\omega \in A_{k}$
and (ii) all integer less than $k$ are not member of this set, that is :  $(0\leq j\leq k-1)\Longrightarrow \omega \notin A_{j}$. Then, we have%
\begin{equation*}
\left\{ 
\begin{tabular}{l}
$\omega \in A_{k}$ \\ 
$\forall (0\leq j\leq k-1),\omega \in A_{j}^{c}$.
\end{tabular}%
\right.
\end{equation*}

\bigskip \noindent Thus 
\begin{equation*}
\omega \in A_{0}^{c}A_{1}^{c}....A_{k-1}^{c}A_{k}=B_{k} \subset \bigcup\limits_{p\geq 0}B_{p}.
\end{equation*}

\bigskip \noindent Hence 
\begin{equation*}
\omega \in \bigcup\limits_{n\geq 0}A_{n}\Longrightarrow \omega
\bigcup\limits_{n\geq 0}B_{n},
\end{equation*}

\bigskip \noindent that is 
\begin{equation}
\bigcup\limits_{n\geq 0}A_{n}\subseteq \bigcup\limits_{n\geq 0}B_{n}.
\label{I-B}
\end{equation}

\bigskip \noindent (c) We want to show that $(B_{n})_{n\geq 0}$ is a sequence pairwise disjoint sets.\newline

\bigskip \noindent Let \ $n$ and $m$ two different non-negative integers. Let $n$ be less that $m $ for example, that is $n<m,$ we have$%
\begin{tabular}{l}
\\ 
\end{tabular}%
$%
\begin{equation*}
\left\{ 
\begin{tabular}{lll}
$B_{n}=A_{0}^{c}A_{1}^{c}....A_{n-1}^{c}A_{n}$ & $\subseteq $ & $A_{n}$ \\ 
$%
B_{m}=A_{0}^{c}A_{1}^{c}....A_{n-1}^{c}A_{n}^{c}A_{n+1}^{c}...A_{m-1}^{c}A_{m} 
$ & $\subseteq $ & $A_{n}^{c}$%
\end{tabular}%
\right. .
\end{equation*}

\bigskip \noindent So%
\begin{equation*}
B_{n}\cap B_{m}\subseteq A_{n}\cap A_{n}^{c}=\emptyset.
\end{equation*}

\bigskip \noindent Thus get
\begin{equation*}
B_{n}\cap B_{m}=\emptyset. \ \ \square
\end{equation*}

\bigskip \noindent CONCLUSION : (a), (b) and (c) together prove that%
\begin{equation*}
\bigcup\limits_{n\geq 0}A_{n}=\bigcup\limits_{n\geq
0}B_{n}=\bigcup\limits_{n\geq 0}A_{n}\subseteq \sum_{n\geq 0}B_{n}.
\end{equation*}

\bigskip \noindent (d) Meaning of $k$ for $\omega \in B_{k}.$ For $\omega
\in \bigcup\limits_{n\geq 0}A_{n},$ $\omega \in B_{k}$ means that the first
index $j$ such that $\omega \in A_{j}$. In probability terminology, $k$ is
the first time $\omega $ belongs to one the the $A_{n}$.\\

\bigskip \noindent \textbf{Exercise 5}. \label{exercise05_sol_doc00-03}\\

\bigskip \noindent (a) Generalize the formula in \textit{(00.14)} :

\begin{equation*}
1_{A\cup B}=1_{A}+1_{B}-1_{A\cap B}
\end{equation*}

\bigskip \noindent to any countable union. Precisely, by applying it to three sets and next to four sets, guess the general formula :

\begin{eqnarray*}
1_{\bigcup_{1\leq j\leq n}A_{j}} \notag \\
&=&\sum_{1\leq j\leq n} 1_{A_{j}}+\sum_{r=2}^{n}(-1)^{r+1} \sum_{1\leq i_{1}<...<t_{r}\leq
n}1_{A_{i_{1}}...A_{ir}}. \ \ (HP)
\end{eqnarray*}

\bigskip \noindent (b) From the formula

$$
Card( A \bigcup B ) =Card(A) + Card(B)- Card(A \bigcap B),
$$

\bigskip \noindent proceed similarly and justify the formula 

\begin{equation*}
Card(\bigcup_{1\leq j\leq n}A_j)=\sum_{1\leq j\leq n}Card(A_{j})+\sum_{r=2}^{n}(-1)^{r+1} \sum_{1\leq i_{1}<...<t_{r}\leq n} Card(A_{i_{1}}...A_{ir}). \ (HP)
\end{equation*}

\bigskip \noindent Both Formulas versions of the Pointcarr\'e Formula.\\

\noindent \textbf{NB} It is only asked to guess the formulas by applying the two-set formula for three sets and next for four sets and by using the associativity of the union operation on sets. But a proof of these formula will be given in the solutions.\\

\bigskip \noindent \textbf{Solution}.\\

\noindent Question (a). First show that

\begin{equation*}
1_{A\cup B}=1_{A}+1_{B}-1_{A\cap B},\
\end{equation*}

\bigskip \noindent that is 
\begin{equation}
\forall \omega \in \Omega ,\text{ }1_{A\cup B}(\omega )=1_{A}(\omega
)+1_{B}(\omega )-1_{A\cap B}(\omega ).  \tag{II-1}
\end{equation}

\bigskip \noindent We have to consider $4$ cases based on the decomposition
: $\Omega =AB+A^{c}B+AB^{c}+A^{c}B^{{c}}$.\newline

\bigskip \noindent Case 1 : $\omega \in AB\Longleftrightarrow (\omega \in A$
and $\omega \in B)$.\newline

\bigskip \noindent In this case, check that the equality in (II-1), is 
\begin{equation*}
1=1+1-1.
\end{equation*}

\bigskip \noindent Case 2 : $\omega \in A^{c}B\Longleftrightarrow (\omega
\in A^{c}$ and $\omega \in B).$\newline

\bigskip \noindent In this case, check that the equality in (II-1), is%
\begin{equation*}
1=0+1-0.
\end{equation*}

\bigskip \noindent Case 2 : $\omega \in AB^{c}\Longleftrightarrow (\omega
\in A$ and $\omega \in B^{c}).$\newline

\bigskip \noindent In this case, check that the equality in (II-1), is%
\begin{equation*}
1=1+0-0.
\end{equation*}

\bigskip \noindent Case 3 : $\omega \in A^{c}B^{c}\Longleftrightarrow
(\omega \in A^{c}$ and $\omega \in B^{c}).$\newline

\bigskip \noindent In this case, check that the equality in (II-1), is%
\begin{equation*}
0=0+0-0.
\end{equation*}

\bigskip \noindent \textbf{Generalization}. Begin to derive this for the Formula of Exercise 1.

\bigskip

\bigskip \noindent case 1 of two sets : $A\cup B=A+A^{c}B.$ Then%
\begin{equation*}
1_{A\cup
B}=1_{A+A^{c}B}=1_{A}+1_{A^{c}B}=1_{A}+1_{A^{c}}1_{B}=1_{A}+(1-1)1_{B}=1_{A}+1_{B}-1_{A}1_{B}
\end{equation*}%
\begin{equation*}
=1_{A}+1_{B}-1_{A}1_{B}.
\end{equation*}

\bigskip \noindent Case 2 of three sets : $A\cup B\cup C$. Write $A\cup B\cup
C=A\cup (B\cup C). $ Apply the formula of Case 1 to the sets $A$ and $B\cup
C $ to get
\begin{equation*}
1_{A\cup (B\cup C)}=1_{A}+1_{B\cup C}-1_{A\cap (B\cup C)}=1_{A}+1_{B\cup
C}-1_{A}1_{(B\cup C)}.
\end{equation*}

\bigskip \noindent Apply again formula of Case 1 to $B\cup C$ to get%
\begin{equation*}
1_{A\cup (B\cup
C)}=1_{A}+1_{B}+1_{C+}1_{C}1_{B}-1_{A}(1_{B}+1_{C}-1_{C}1_{B})
\end{equation*}

\begin{equation*}
=\left\{ 1_{A}+1_{B}+1_{C}\right\} -\left\{
1_{A}1_{B}+1_{A}1_{C}+1_{C}1_{B}\right\} +1_{A}1_{B}1_{C}
\end{equation*}

\begin{equation*}
=\left\{ 1_{A}+1_{B}+1_{C}\right\} -\left\{ 1_{AB}+1_{AC}+1_{CB}\right\}.
-1_{ABC}
\end{equation*}

\bigskip \noindent Now, we may writing this formula for $A=A_{1}$, $B=A_{2}$ and $C=A_{3}$ in the form
\begin{equation*}
1_{\bigcup\limits_{1\leq j\leq 3}A_{j}}=\sum\limits_{1\leq j\leq
3}A_{j}+(-1)^{2+1}\sum\limits_{1\leq i<j\leq 3}1_{A_{i}A_{j}}+(-1)^{3+1}\sum\limits_{1\leq i<j<k\leq 3}1_{A_{i}A_{j}A_{k}}.
\end{equation*}

\bigskip \noindent This leads to the Poincarr\'{e} Formula :\\

\noindent \textbf{Inclusion-Exclusion Formula} : for all $n\geq 2$, 

\begin{eqnarray*}
&&1_{\bigcup\limits_{1\leq j\leq n}A_{j}} \notag \\
&=&\sum\limits_{1\leq j\leq n} 1_{A_{j}}+\sum_{r=2}^{n}(-1)^{r+1} \sum\limits_{1\leq i_{1}<...<t_{r}\leq
n}1_{A_{i_{1}}...A_{ir}}. \ \ \ (\mathcal{P}_n)
\end{eqnarray*}

\bigskip \noindent As promised, here is the proof of the formula.\\

\noindent \textbf{Proof of the Pointcarr\'e's Formula}.\\

\noindent To prove the general formula (say $\mathcal{P}_n$, for $n\geq 2$), we will proceed by induction. Since we know that it holds for $n=2$, we have to prove the induction principle, that is : $\mathcal{P}_n \Rightarrow \mathcal{P}_{n+1}$ for all $n\geq 2$.\\

\noindent Now, suppose that $\mathcal{P}_k$ holds for $k\leq n\geq 2$. We have to prove that $\mathcal{P}_{n+1}$ holds. We begin with

$$
\bigcup_{1\leq j\leq n+1}A_{j} = \bigcup_{1\leq j\leq n}A_{j} \bigcup A_{n+1}.
$$

\bigskip \noindent which implies that

\begin{eqnarray*}
&&1_{\bigcup\limits_{1\leq j\leq n+1}A_{j}} \notag \\
&=&\sum\limits_{1\leq j\leq n+1} 1_{A_{j}} \ (I_{n+1})\\
&+& \sum_{r=2}^{n+1}(-1)^{r+1} \sum\limits_{1\leq i_{1}<...<t_{r}\leq n+1}1_{A_{i_{1}}...A_{ir}}. \ (II_{n+1}) \ \ \ (\mathcal{P}_{n+1})
\end{eqnarray*}

\bigskip \noindent To be in a better position to compare $\mathcal{P}_n$ and $\mathcal{P}_{n+1}$, let us also write

\begin{eqnarray*}
&&1_{\bigcup\limits_{1\leq j\leq n}A_{j}} \notag \\
&=&\sum\limits_{1\leq j\leq n+1} 1_{A_{j}} \ (I_n)\\
&+& \sum_{r=2}^{n+1}(-1)^{r+1} \sum\limits_{1\leq i_{1}<...<t_{r}\leq n}1_{A_{i_{1}}...A_{ir}}. \ (II_{n}) \ \ \ (\mathcal{P}_{n})
\end{eqnarray*}

\bigskip 
\noindent Now, since we have

$$
\bigcup_{1\leq j\leq n+1}A_{j} = \bigcup_{1\leq j\leq n} A_{j} \bigcup A_{n+1}, 
$$

\bigskip 
\noindent we may apply the induction principle for two sets to have

$$
1_{\bigcup_{1\leq j\leq n+1}A_{j}}= 1_{\bigcup_{1\leq j\leq n}A_{j}} + 1_{A_{n+1}} - 1_{\bigcup_{1\leq j\leq n}A_{j} \bigcap A_{n+1}}. \ (I_2)
$$

\bigskip \noindent By developing the first term (that is  $\mathcal{P}_n$) and by adding $1_{A_{n+1}}$, we get

\begin{eqnarray*}
&&1_{\bigcup_{1\leq j\leq n+1}A_{j}} \notag \\
&=& \sum_{1\leq j\leq n+1} 1_{A_{j} + 1_{A_{n+1}}} \ \ (L1) \\
&+& \sum_{1\leq i_{1}<t_{2}\leq n} 1_{A_{i_1}A_{i_2}} \ \ (L2).\\
&+& \sum_{r=3}^{n} (-1)^{r+1} \sum_{1\leq i_{1}<...<t_{r}\leq n} 1_{A_{i_1}...A_{i_r}}. \  \ (L3).
\end{eqnarray*}

\bigskip 
\noindent Further, let us  develop $(I_2)$ by using the hypothesis that the formula holds for $n$ sets. We get

\begin{eqnarray*}
&&-1_{\bigcup_{1\leq j\leq n}A_{j} \bigcap A_{n+1}}\\
&=&-1_{\bigcup_{1\leq j\leq n} \left(A_{j}\bigcap A_{n+1}\right) }\\
&=& -\sum_{1\leq j\leq n+1} 1_{A_{j} \cap A_{n+1}}\ (L4)\\
&+& \sum_{r=2}^{n-1} (-1)^{r+1} \sum_{1\leq i_{1}<...<t_{r}\leq n} 1_{A_{i_1}...A_{i_r}A_{n+1}} \ (L5) \\
&+&  (-1)^{n+1} 1_{A_{i}...A_{n+1}}. \ (L6) \\
\end{eqnarray*}

\bigskip  \noindent Please, pay attention to the two last lines : we separate the two cases : $2 \leq r \leq n-1$ and $r=n-1$. For $2 \leq r \leq n-1$, we actually have the sum

$$
\sum_{r=3}^{n} (-1)^{r+1} \sum_{1\leq i_{1}<...<t_{r}\leq n} 1_{A_{i_1}...A_{i_r}},
$$

\bigskip \noindent where one of the terms $A_{i_j}$, $1\leq j \leq r$ is $A_{n+1}$. For $r=n$, and because of all the terms $A_{i_j}$ should be distinct between them and distinct from $A_{n+1}$, we get one term.\\

\noindent Let us conclude by comparing the two last blocks of formulas and Formula $\mathcal{P}_{n+1}$. We can draw the following facts :\\

\noindent 1) the line $I_{n+1}$ in $\mathcal{P}_{n+1}$ us exactly line $(L1)$ above,\\

\noindent 2) the line $II_{n+1}$ in $\mathcal{P}_{n+1}$ can be decomposed as follows :\\

\hskip 1cm 2a) for $r=2$, by the sums of the terms $1_{A_{i_1}A_{i_2}}$  when both $A_{i_1}$ and $A_{i_2}$ are different from $A_{n+1}$ and this corresponds to line $(L2)$, and by the sums of the terms 
$1_{A_{i_1}A_{i_2}}$ when one of $A_{i_1}$ and $A_{i_2}$ is $A_{n+1}$ and this corresponds to $(L4)$. By taking into account the common negative sign, the summation of Lines $(L2)$ and $(L3)$ gives the term in $II_{n+1}$ in $\mathcal{P}_{n+1}$ corresponding to $r=2$,\\

\hskip 1cm 2b) for $3\leq r\leq n$, by the sums of the terms when all the $A_{i_j}$, $1\leq j \leq r$, are different from $A_{n+1}$ and this corresponds to Line $(L3)$, and by the sums of such terms when one of the $A_{i_j}$, $1\leq j \leq r$, is $A_{n+1}$ and this corresponds to Line $(L5)$ by the remark below (L6),\\

\hskip 1cm 2c) for $r=n+1$, we have one case in which we take all the sets $A_1$, ..., $A_{n+1}$ corresponding to Line (L6). In total, the summation of Lines (L3), (L5) and (L6) give 
$II_{n+1}$ in $\mathcal{P}_{n+1}$ corresponding to $\geq 2$.\\

\noindent By putting together all these facts, we see that $\mathcal{P}_{n+1}$ is implied by $\mathcal{P}_{n}$ and $\mathcal{P}_{2}$. Hence the induction reasoning is complete.\\

\bigskip \noindent Question (b). The proof used in the solution of Question (b) may be reproduced exactly in the same way and lead to the same results.\\

\bigskip \noindent \textbf{Exercise 6} \label{exercise06_sol_doc00-03}\newline

\noindent Show that :

\begin{equation*}
1_{\bigcup\limits_{n\geq 0} A_n} = \sup_{n\geq 0} 1_{A_n}.
\end{equation*}

\bigskip 
\noindent and

\begin{equation*}
1_{\bigcap\limits_{n\geq 0}A_{n}}=\inf_{n\geq 0}1_{A_{n}}.
\end{equation*}

\bigskip \noindent \textbf{SOLUTIONS}.\newline

\bigskip \noindent (A) We want to show that%
\begin{equation*}
1_{\bigcup\limits_{n\geq 0}A_{n}}=\sup_{n\geq }1_{A_{n}},
\end{equation*}

\bigskip \noindent that is 
\begin{equation}
\forall \omega \in \Omega, \ 1_{\bigcup\limits_{n\geq 0}A_{n}}(\omega
)=\sup_{n\geq }1_{A_{n}}(\omega ).  \tag{III-A}
\end{equation}

\bigskip \noindent We are going to show the validity of (III-A) for the two
following possible cases.\newline

\bigskip \noindent (a) Case $\omega \in \bigcup\limits_{n\geq 0}A_{n}.$%
\newline

\bigskip \noindent In this case : 
\begin{eqnarray*}
&&\omega \in \bigcup\limits_{n\geq 0}A_{n}\text{ ( that is }%
1_{\bigcup\limits_{n\geq 0}A_{n}}(\omega )=1\text{ ) }\Longleftrightarrow
(\exists n_{0}\geq 0,\omega \in A_{n_{0}})\\
&&\Longleftrightarrow (\exists
n_{0}\geq 0,1_{A_{n_{0}}}(\omega )=1).
\end{eqnarray*}

\bigskip \noindent From this equivalence, we may say that :\newline

\bigskip \noindent (i) The left-hand member of the equality in (III-A) is $1$ since 
$\omega \in \bigcup\limits_{n\geq 0}A_{n}$.\newline

\bigskip \noindent (ii) The right-hand member of the equality in, is the supremum
of numbers that are either equal to $0$ or to $1,$ and one one them is $1$ \ 
$(for$ $n=n_{0}).$ Then the supremum is $1.$\newline

\bigskip \noindent So the equality in (III-A) holds for $\omega \in
\bigcup\limits_{n\geq 0}A_{n}$.\newline

\bigskip \noindent (b) Case $\omega \neq \bigcup\limits_{n\geq 0}A_{n}.$%
\newline

\bigskip \noindent In this case, $\omega $ does not belong to any of the $%
A_{n}.$ We have\newline

\bigskip

\begin{eqnarray*}
&&\omega \neq \bigcup\limits_{n\geq 0}A_{n}\text{ ( that is }%
1_{\bigcup\limits_{n\geq 0}A_{n}}(\omega )=0\text{ ) }\Longleftrightarrow
(\forall n\geq 0,\omega \notin A_{n})\\
&&\Longleftrightarrow (\forall n\geq
0,1_{A_{n}}(\omega )=0).
\end{eqnarray*}

\bigskip \noindent From this equivalence, we may say\newline

\noindent (i) The left side of the equality in (III-A) is $0$ since $\omega
\notin \bigcup\limits_{n\geq 0}A_{n}$.\newline

\bigskip \noindent (ii) The right side of the equality in (III-A) is the
supremum of numbers that are all equal to $0$. Then the supremum is $0$.%
\newline

\bigskip \noindent So the equality in (III-A) holds for $\omega \notin
\bigcup\limits_{n\geq 0}A_{n}.$.\newline

\bigskip \noindent This completes our answer.\newline

\bigskip \noindent (B) We want to show that%
\begin{equation*}
1_{\bigcap\limits_{n\geq 0}A_{n}}=\inf_{n\geq 0}1_{A_{n}},
\end{equation*}

\bigskip \noindent that is 
\begin{equation}
\forall \omega \in \Omega ,1_{\bigcap\limits_{n\geq 0}A_{n}}(\omega
)=\inf_{n\geq 0}1_{A_{n}}(\omega ).  \tag{III-B}
\end{equation}

\bigskip \noindent We are going to show the validity of (III-B) for the two
following possible cases.\newline

\bigskip \noindent (a) Case $\omega \in \bigcap\limits_{n\geq 0}A_{n}.$%
\newline

\bigskip \noindent In this case : 
\begin{eqnarray*}
&&\omega \in \bigcap\limits_{n\geq 0}A_{n}\text{ ( that is } 1_{\bigcap\limits_{n\geq 0}A_{n}}(\omega )=1\text{ ) }\Longleftrightarrow
(\forall n\geq 0,\omega \in A_{n})\\
&&\Longleftrightarrow (\forall n\geq 0,1_{A_{n}}(\omega )=1).
\end{eqnarray*}

\bigskip \noindent From this equivalence, we may say that :\newline

\bigskip \noindent (i) The left side of the equality in (III-B) is $1$ since $\omega \in \bigcap\limits_{n\geq 0}A_{n}.$\newline

\bigskip \noindent (ii) The right side of the equality (III-B) is the infinum of numbers that are all equal to $1$. Then the infinum is $1$.%
\newline

\bigskip \noindent So the equality in (III-B) holds for $\omega \in \bigcap\limits_{n\geq 0}A_{n}$.\newline

\bigskip \noindent \bigskip (b) Case $\omega \notin \bigcap\limits_{n\geq
0}A_{n}.$\newline

\bigskip \noindent In this case, $\omega$ does not belong, at least, to one
the $A_{n}.$ We have \newline

\begin{eqnarray*}
&&\omega \in \bigcap\limits_{n\geq 0}A_{n}\text{ ( that is }%
1_{\bigcap\limits_{n\geq 0}A_{n}}(\omega )=0\text{ ) }\Longleftrightarrow
(\exists n_{0}\geq 0,\omega \notin A_{n_{0}})\\
&&\Longleftrightarrow (\exists
n_{0}\geq 0,1_{A_{n_{0}}}(\omega )=0).
\end{eqnarray*}

\bigskip \noindent From this equivalence, we may say that : \newline

\bigskip \noindent (i) The left side of (III-B) is $0$ since $\omega \notin
\bigcap\limits_{n\geq 0}A_{n}$.\newline

\bigskip \noindent (ii) The right side is the infinum of numbers that are
either equal to $0$ or to $1,$ and one one them is $0$ (for $n=n_{0}$). Then
the infinum is $0$.\newline

\bigskip \noindent So the equality in (III-B) holds for $\omega \notin
\bigcap\limits_{n\geq 0}A_{n}$.\newline

\bigskip \noindent Other method : Use one of the formula to show the other.%
\newline

\bigskip \noindent Suppose that 
\begin{equation*}
1_{\bigcup\limits_{n\geq 0}A_{n}}=\sup_{n\geq }1_{A_{n}}.
\end{equation*}

\bigskip \noindent But 
\begin{equation*}
1_{\bigcap\limits_{n\geq 0}A_{n}}=1_{\left( \bigcup\limits_{n\geq
0}A_{n}^{c}\right) ^{c}}=1-1_{\bigcup\limits_{n\geq 0}A_{n}^{c}}.
\end{equation*}

\bigskip \noindent and this is

\begin{eqnarray*}
&=&1-\sup_{n\geq }1_{A_{n}^{c}}=1-\sup_{n\geq 0}(1-1_{A_{n}})=1-\left\{
1+\sup_{n\geq 0}(-1_{A_{n}})\right\} \\
&=&1-\left\{ 1-\inf_{n\geq 0}1_{A_{n}}\right\}\\
&=&\inf_{n\geq 0}1_{A_{n}}.
\end{eqnarray*}

\bigskip \bigskip  

\noindent \textbf{Exercise 7}.  \label{exercise07_sol_doc00-03}Let $\Omega =\Omega _{1}$ $\times \Omega _{2}$ be a
non-empty product defined by%
\begin{equation*}
\Omega =\Omega _{1}\times \Omega _{2}=\{(\omega _{1},\omega _{2}),\omega
_{1}\in \Omega _{1}\text{ and }\omega _{2}\in \Omega _{2}\}.
\end{equation*}

\bigskip \noindent Define a rectangle%
\begin{equation*}
A\times B=\{(\omega _{1},\omega _{2})\in \Omega ,\ \omega _{1}\in A\text{
and }\omega _{2}\in B\},
\end{equation*}

\bigskip 
\noindent where $A\subseteq \Omega _{1}$, $B\subseteq \Omega _{2}.$\\

\noindent (a) Show that the intersection of rectangles is a rectangle.\\

\noindent \textbf{Solution} : Let $A\times B$ and $C\times D$ be two rectangles. Let us form
the intersection :

\begin{eqnarray*}
(A\times B)\cap (C\times D) &=&\{(\omega _{1},\omega _{2})\in \Omega
,(\omega _{1},\omega _{2})\in A\times B\text{ and }(\omega _{1},\omega
_{2})\in C\times D\} \\
&=&\{(\omega _{1},\omega _{2})\in \Omega ,\omega _{1}\in A\text{ and }\omega
_{2}\in B\text{ and }\omega _{1}\in C\text{ and }\omega _{2}\in D\} \\
&=&\{(\omega _{1},\omega _{2})\in \Omega ,\omega _{1}\in A\cap C\text{ and }%
\omega _{2}\in B\cap D\} \\
&&A\cap C\times B\cap D.
\end{eqnarray*}

\bigskip 
\noindent (b) Let $A\times B$ be a rectangle. Show that%
\begin{equation*}
(A\times B)^{c}=A^{c}\times B+A\times B^{c}+A^{c}\times B^{c}.
\end{equation*}

\bigskip \noindent \textbf{Solution}. We may partition $\Omega $ in the following way : for
each $(\omega _{1},\omega _{2})\in \Omega ,$ the four exclusive cases my
occur : $\omega _{1}\in A$ and $\omega _{2}\in B$ or $\omega _{1}\notin A$
and $\omega _{2}\in B$ or $\omega _{1}\in A$ and $\omega _{2}\notin B$ or $%
\omega _{1}\notin A$ and $\omega _{2}\notin B$. This gives

\begin{equation*}
\Omega =A\times B+A^{c}\times B+A\times B^{c}+A^{c}\times B^{c}.
\end{equation*}

\bigskip 
\noindent Then the complement of $A\times B$ is $A^{c}\times B+A\times
B^{c}+A^{c}\times B^{c}.$\\

\noindent (c) Show that not all subsets of are rectangles. Use the following example.\\

\noindent Let $\Omega _{0}=\Omega _{1}=\Omega _{2}$ so that $\Omega =\Omega _{0}\times
\Omega _{0}=\Omega _{0}^{2}$ is the product of one space $\Omega _{0}$ by
itself and consider its diagonal%
\begin{equation*}
D=\{(\omega _{1},\omega _{2})\in \Omega ,\text{ }\ \omega _{1}=\omega _{2}\}.
\end{equation*}

\bigskip \noindent Suppose that $\Omega _{0}$ has at least two elements.
Suppose that $D=A\times B.$ Show that this is not possible by considering
two cases : (1) Both $A$ and $B$ are singletons, (2) $A$ or $B$ has at least
two elements.\\

\noindent \textbf{Solution}. Suppose that $\omega _{1},$ $\omega _{2}$ be two elements of $%
\Omega _{0}$. \ Suppose that $D=A\times B.$\\

\noindent (1) Suppose that $A$ and $B$ are singletons. Let $A=\{a\}$ and $B=\{b\}.$ Then 
\begin{equation*}
A\times B=\{(a,b)\}.
\end{equation*}

\bigskip \noindent If $D=A\times B=\{(a,b)\},$ then $a=b$ and $D=\{(a,a)\}$ is a singleton.
But, at least, $D$ has at least the two elements $(\omega _{1},\omega _{1})$
and $(\omega _{2},\omega _{2}).$ This is impossible.\\

\noindent (2) Suppose that $A$ has at least two elements $\omega _{1}$ and $\omega
_{2}.$ Let $\omega _{3}\in B.$ But  $D=A\times B$ implies that $(\omega
_{1},\omega _{3})\in D$ and $(\omega _{2},\omega _{3})\in D.$ This implies
that $\omega _{1}=\omega _{2}=\omega _{3}$. This is impossible since $\omega
_{1}\neq \omega _{2}$.

\newpage
%\include{mes_doc_00_04_ap}
%\newpage
%\include{mes_doc_00_05_ap}
%\chapter{Measurable Sets}

\chapter{Measurable Sets} \label{01_setsmes}

\bigskip
\noindent \textbf{Content of the Chapter}

\begin{table}[htbp]
	\centering
		\begin{tabular}{llll}
		\hline
		Type& Name & Title  & page\\
		\hline
		S & Doc 01-01 & Measurable sets - An outline  & \pageref{doc01-01}\\
		D &  Doc 01-02& Exercises on Measurable sets  & \pageref{doc01-02} \\
		SD& Doc 01-03 &    Exercises on Measurable sets with Solutions& \pageref{doc01-03}\\
		D& Doc 01-04 & Measurable sets. Exercises on $\lambda$ and $\pi$ systems& \pageref{doc01-04}\\
		SD& Doc 01-05 &   Measurable sets. Exercises on $\lambda$ and $\pi$ systems with solutions  & \pageref{doc01-05} \\
		D& Doc 01-05 & More  Exercises  & \pageref{doc01-06}\\
		SD& Doc 01-07 & More  Exercises - with solutions  & \pageref{doc01-07}\\
		\hline
		\end{tabular}
\end{table}

\newpage
\noindent \LARGE \textbf{DOC 01-01 :  Measurable sets - A summary}. \label{doc01-01}\\

 \Large

\bigskip \noindent Speaking about measurable sets on $\Omega $ requires
first one defines a field \ or a $\sigma $-algebra $\mathcal{A}$\ of subsets
on $\Omega $. The pair $(\Omega ,\mathcal{A})$ is called a measurable space.%
\newline

\bigskip \noindent \textbf{I - $\sigma $-algebras}.

\bigskip \noindent \textbf{(01.01)} (Definition) A collection $\mathcal{A}$\
of subsets of $\Omega ,$ that is :\ $\mathcal{A}$\ $\subset \mathcal{P}%
(\Omega ),$ is a field or a $\sigma $-algebra of subsets of $\Omega ,$ if
and only if (it contains $\Omega $ and $\emptyset ),$ and it is stable by
all finite and countable sets operations.\newline

\noindent \textbf{NB}. In this document and related ones, \textbf{countable} means \textbf{finite} or \textbf{infinite countable}.\\

\bigskip \noindent \textbf{(01.02)} If $\mathcal{A}$\ is $\sigma $-algebra of
subsets of $\Omega ,$ the couple $(\Omega ,\mathcal{A})$ is said to be a
measurable space. The elements of $\mathcal{A}$ are called measurable sets.%
\newline

\bigskip \noindent \textbf{(01.03)} Warning : Never use measurability,
without having in mind the underlying $\sigma $-algebra we are working with.%
\newline

\bigskip \noindent \textbf{(01.04) Criterion} : To prove that $\mathcal{A}$\
is $\sigma $-algebra of subsets of $\Omega ,$ we only need to check a few
number of conditions. For example, we may check.\newline

\bigskip \noindent \textbf{(i)} $\Omega \in \mathcal{A}$ \ or $\emptyset \in 
\mathcal{A}$.\newline

\noindent \textbf{(ii)} $\mathcal{A\ }$is stable under complementation : $%
A\in \mathcal{A}\Longrightarrow $ $A^{c}\in \mathcal{A}$.\newline

\noindent \textbf{(iii)} ($\mathcal{A}$ is stable under countable unions : $%
\{A_{n},n\geq 0\}\subset \mathcal{A}$ $\Longrightarrow \bigcup\limits_{n\geq 0}A_{n}\in \mathcal{A}$ ) or ($\mathcal{A}$ is stable
under countable intersections : $\{A_{n},n\geq 0\}\subset \mathcal{A}$ $%
\Longrightarrow $ $\bigcap\limits_{n\geq 0}A_{n}\in \mathcal{A}$ ).\newline

\bigskip \noindent \textbf{(01.05)} Examples of $\sigma$-algebras : (1) $\mathcal{A}$\ $=\{\emptyset ,\Omega \}$,
 (2) $\mathcal{A}=\mathcal{P}(\Omega)$, and
  
$$
\mathcal{A}=\{ A\subset \Omega,\text{ } A \text{ is countable or } A^{c} \text{ is countable }\}.
$$

\bigskip \noindent \textbf{(01.06)} In general, usual $\sigma $-algebras are
not explicitly known as in \textit{(01.05)}. A $\sigma $-algebra may be obtained by
giving a specific sub-collection $\mathcal{C}$ of it, in the following way :
There exists a smallest $\sigma $-algebra (in the sense of inclusion) that
contains $\mathcal{C}$. This $\sigma $-algebra is denoted $\sigma $($%
\mathcal{C}$) and is called : the $\sigma $-algebra generated by $\mathcal{C%
}$. In summary, we simultaneously have three points 
\begin{equation*}
\left\{ 
\begin{tabular}{l}
(i) $\sigma $($\mathcal{C}$) is a $\sigma $-algebra, \\ 
(ii) $\mathcal{C}\subset \sigma $($\mathcal{C}$), \\ 
(iii) For any other $\sigma $-algebra $\mathcal{D}$\ such that \ $\mathcal{C%
}\subset $\ $\mathcal{D}$, we have : $\sigma $($\mathcal{C}$) $\subset $ $%
\mathcal{D}$.
\end{tabular}%
\right.
\end{equation*}

\bigskip \noindent \textbf{(01.07)} Example : (1) The usual $\sigma $-algebra
on $\mathbb{R}$ is generated by the class $\mathcal{I}$ of all intervals of $%
\mathbb{R}$. For a topological space $(E,T)$ where $T$ is class of open
sets, the usual \ $\sigma $-algebra on $\mathcal{B}(T)$, called the Borel $%
\sigma $-algebra, is generated by open sets. It is also generated by closed
sets.\newline

\bigskip \noindent \textbf{(01.08)} Very important example : The product $\sigma $-algebra. Let $(\Omega _{i},\mathcal{A}_{i})$ two measurable spaces and let $\Omega =\Omega _{1}\times \Omega _{2}=\{(\omega _{1},\omega
_{2}),\omega _{1}\in \Omega _{1},\omega _{2}\in \Omega _{2}\}$ be the Cartesian product of $\Omega _{1}$ and $\Omega _{2}.$ The class of cylinders (or rectangles) is given by

\begin{equation*}
\mathcal{S}=\{A_{1}\times A_{2},\text{ }A_{1}\in \mathcal{A}_{1},A_{2}\in 
\mathcal{A}_{2}\}.
\end{equation*}

\bigskip \noindent The product $\sigma $-algebra on $\Omega =\Omega_{1}\times \Omega _{2},$ denoted as \ $\mathcal{A}_{1}\otimes \mathcal{A}_{2} $, is generated by $\mathcal{S}$ : 

\begin{equation*}
\mathcal{A}_{1}\otimes \mathcal{A}_{2}=\sigma (\mathcal{S}).
\end{equation*}

\bigskip \noindent $(\Omega _{1}\times \Omega _{2}, \mathcal{A}_{1}\otimes 
\mathcal{A}_{2})$ is called the product measurable space.\newline

\newpage

\noindent \textbf{II - Algebra}.\newline

\bigskip \noindent Usually, the class $\mathcal{C}$ that generates our
working \ $\sigma $-algebra, is an algebra. What is an algebra of subsets
of $\Omega$?\newline

\bigskip \noindent \textbf{(01.09)} (Definition) A collection $\mathcal{C}$\
of subsets of $\Omega ,$ that is \ $\mathcal{C}$\ $\subset \mathcal{P}%
(\Omega)$ is an algebra of subsets of $\Omega ,$ if and only if (it
contains $\Omega $ and $\emptyset ),$ and it is stable by all finite sets
operations.\newline

\bigskip \noindent \textbf{(01.10)} Criterion : To prove that $\mathcal{C}$\
is an algebra of subsets of $\Omega$, we only need to check a few number
of conditions. For example, we may check :

\bigskip \noindent \textbf{(i)} $\Omega \in \mathcal{C}$ \ or $\emptyset \in 
\mathcal{C}$\newline
\noindent \textbf{(ii)} $\mathcal{C}$ is stable under complementation : $%
A\in \mathcal{C}$ $\Longrightarrow $ $A^{c}\in \mathcal{A}$.\newline
\noindent \textbf{(iii)} ($\mathcal{C}$ is stable under finite unions : $%
\{A,B\}\subset \mathcal{C}\Longrightarrow A\cup B\in \mathcal{C}$ ) or ($%
\mathcal{C}$ is stable under finite intersections : $\{A,B\}\subset \mathcal{%
C}$ $\Longrightarrow $ $A\cap B\in \mathcal{C}$)

\bigskip \noindent \textbf{(01.11)} Examples : (1) $\mathcal{C}$\ $%
=\{\emptyset,\Omega \},$ $\mathcal{C}$\ $=\mathcal{P}(\Omega ),$ $\mathcal{C}
$\ $=\{A\subset \Omega ,A$ is finite or $A^{c}$ is finite$\}$.\newline

\bigskip \noindent \textbf{(01.12)} In general, usual algebras, are not
explicitly known as in (01.11). An algebra may be obtained by giving a
specific sub-collection $\mathcal{S}$ of it, in the following way : There
exists a smallest algebra (in the sens of inclusion) that contains $%
\mathcal{S}$. This algebra is denoted $a$($\mathcal{S}$) and is called :
the algebra generated by $\mathcal{S}$. In summary, we simultaneously have
three points :\newline

\begin{equation*}
\left\{ 
\begin{tabular}{l}
(i) $a$($\mathcal{S}$) is an algebra, \\ 
(ii) $\mathcal{S}\subset \sigma $($\mathcal{S}$), \ \  \\ 
(iii) For any other algebra $\mathcal{D}$\ such that \ $\mathcal{S}\subset $%
\ $\mathcal{D}$, we have : $\sigma $($\mathcal{C}$) $\subset $ $\mathcal{D}$%
\end{tabular}.
\right.
\end{equation*}

\newpage

\noindent \textbf{III - Semi-algebras}.\newline

\bigskip \noindent Usually, the class $\mathcal{S}$ that generates our
working \ algebra, is a semi-algebra. What is a semi-algebra of subsets of $%
\Omega .$\newline

\bigskip \noindent \textbf{(01.013)} (Definition) A collection $\mathcal{S}$\
of subsets of $\Omega ,$ that is \ $\mathcal{S}$\ $\subset \mathcal{P}%
(\Omega ),$ is a semi-algebra of subsets of $\Omega ,$ if and only if :

\begin{equation*}
\left\{ 
\begin{tabular}{l}
$\Omega \in \mathcal{S}$ \\ 
$\mathcal{S}$ is stable by finite intersections. \\ 
The complement of any element of $\mathcal{S}$ is a finite sum of elements of $\mathcal{S}$\ .%
\end{tabular}%
\right.
\end{equation*}

\bigskip \noindent \textbf{NB} : We take the elements of $\mathcal{S}$ as finite sums of elements of $\mathcal{S}$ with one term. Besides, the empty $\empty$ belongs to $\mathcal{S}$ and each element of 
$\mathcal{S}$ is a sum with the empty set.\\

\bigskip \noindent \textbf{(01.14)} : The algebra generated by a semi-algebra
\ $\mathcal{S}$ is the class of all finite sums of its elements 
\begin{equation*}
a(\mathcal{S})=\{A_{1}+A_{1}+...+A_{k},\text{ }A_{i}\in \mathcal{S},\text{ }%
k\geq 1\}.
\end{equation*}

\bigskip \noindent \textbf{(01.15)} Examples.\newline

\bigskip \noindent (a) The class $\mathcal{I}$ of all intervals of $\mathbb{R%
}$ is a semi-algebra (see (01.07)).\newline

\bigskip \noindent (b) The class $\mathcal{I}=\{]a,b], -\infty \leq a \leq b
\leq +\infty \}$ is a semi-algebra on $\mathbb{R}$.\newline

\bigskip \noindent NB. On $\mathbb{R}$, interpret $]a,+\infty] \equiv
]a,+\infty[$.\newline

\bigskip \noindent (c) The class of cylinders in a product space (see (01.08))
is a semi-algebra.\newline

\newpage \noindent \textbf{IV - Special subsets classes}.\newline

\bigskip \noindent \textbf{(01.16)} \textbf{Monotone class}. A nonempty class $%
\mathcal{M}$ of subsets of $\Omega$ is said to be monotone if and only if
its contains the limits of its monotone sequences. That is, if $(A_n)_{n\geq
0}$ is a non-decreasing sequence of subsets in $\mathcal{M}$, then $%
\cup_{n\geq 0} A_n$ is still in $\mathcal{M}$, and as well if $(A_n)_{n\geq
0}$ is a non-increasing sequence of subsets in $\mathcal{M}$, then $%
\cap_{n\geq 0} A_n$ is still in $\mathcal{M}$.\newline

\bigskip \noindent \textbf{(01.17)} (Dynkin system or $\lambda$-system). A
collection $\mathcal{D}$ of subsets of $\Omega$. is a Dynkin system if and
only if

\noindent \textbf{(i)} $\Omega \in \mathcal{D}$,\newline
\noindent \textbf{(ii)} For any $(A,B)\in \mathcal{D}^{2}$ with $A\subset B$%
, $B\backslash A\in \mathcal{D}$,\newline
\noindent \textbf{(iii)} For any sequence $(A_{n})_{n\geq 0})$ of pairwise
disjoint elements of $\mathcal{D}$, $\sum_{n\geq 1}A_{n}\in \mathcal{D}$.%
\newline

\bigskip 

\noindent \textbf{(01.18)} ($\pi $-system). A nonempty collection $\Pi $ of
subsets $\Omega $. is a $\pi $-system if and only if it is stable under
finite intersection.

\bigskip

\noindent NB : These classes are used for demonstration purposes only. For example, the $\lambda -\mu$ method is base on them (See Exercise 11 in Doc 04-02, Chapter \ref{04_measures}, page \pageref{exercise11_doc04-02}).

\newpage
\noindent \LARGE \textbf{DOC 01-02 : Exercises on Measurable sets}. \label{doc01-02}\\

\Large 

\bigskip \noindent \textbf{Exercise 1}. \label{exercise01_doc01-02} Show that the class $\mathcal{I}$ of all intervals of $\mathbb{R}$ is a semi-algebra.\\

\bigskip \bigskip \noindent \textbf{Exercise 2}. \label{exercise02_doc01-02}(01.08) Show that the class of rectangles $\mathcal{S}$ is a semi-algebra. Extend this result to a
product of $k$ measurable spaces, $k\geq 2$.\\

\bigskip \bigskip \noindent \textbf{Exercise 3}. \label{exercise03_doc01-02} Le $\Omega$ be a nonempty set and set

\begin{equation*}
\mathcal{A}=\{A\subset \Omega, \ \ A \ \ is \ \ finite \ \ or \newline
A^{c} \ \ is \ finite \}.
\end{equation*}

\bigskip \noindent 
\noindent (a) Show that $\mathcal{A}$ is an algebra.\newline

\noindent (b) Show that $\mathcal{A}$ is a $\sigma$-algebra if $\Omega$ is finite.\newline

\noindent (c) If $\Omega $ is infinite, consider an infinite subset $C=\{{\omega _{0},\omega _{1},\omega _{2},..\}}$ and set $A_{n}=\{{\ \omega _{2n}\} }$, $n\geq 0$. Show that the $A_{n}$ are in $\mathcal{A}$ and their union is
not.\newline

\noindent (d) Conclude that $\mathcal{A}$ is a $\sigma$-algebra \textbf{if and only if} $\Omega$ is finite.\newline

\bigskip

\bigskip \noindent \textbf{Exercise 4}. \label{exercise04_doc01-02} Let $\Omega$ be a nonempty set and
put

\begin{equation*}
\mathcal{A}=\{A\subset \Omega, \ \ A \ \ is \ \ countable \ \ or \newline
A^{c} \ \ is \ countable \}.
\end{equation*}

\bigskip \noindent 
\noindent (a) Show that $\mathcal{A}$ is a $\sigma$-algebra.\newline

\noindent (b) Let $\Omega=[0,1]$. Show that $[0,1/2]$ is not measurable in the measurable space $(\Omega,\mathcal{A})$.\newline

\bigskip

\bigskip \noindent \textbf{Exercise 5}. \label{exercise05_doc01-02} Prove 01.14 : The algebra generated by a semi-algebra \ $\mathcal{S}$ is the class of all finite sums of its elements 
\begin{equation*}
a(\mathcal{S})=\{A_{1}+A_{2}+...+A_{k},\text{ }A_{i}\in \mathcal{S},\text{ }%
k\geq 1\}.
\end{equation*}

\bigskip

\bigskip \noindent \textbf{Exercise 6}. \label{exercise06_doc01-02} Let $X$ be an application from $\Omega _{1}$ to $\Omega _{2}$ and suppose that $\Omega _{2}$ is endowed with an $\sigma $-algebra $\mathcal{A}_{2}$. Consider this class of subsets of $\Omega _{1}$ :

\begin{equation*}
\mathcal{A}_{X}=\{X^{-1}(B),B\in \mathcal{A}_{2}\},
\end{equation*}

\bigskip \noindent 
\noindent where $X^{-1}(B)$ denotes inverse image of $B$ by $X$ defined by

\begin{equation*}
X^{-1}(B)={\omega \in \Omega_1, X(\omega) \in B}.
\end{equation*}

\bigskip \noindent 
\noindent Show that $\mathcal{A}_{X}$ is a $\sigma $-algebra, called $\sigma$-algebra generated by $X$ (and by $\mathcal{A}_{2}$).\newline

\bigskip

\bigskip \noindent \textbf{Exercise 7}. \label{exercise07_doc01-02}  $(\Omega,\mathcal{A})$ be a measurable space and let $A$ a subset of $\Omega$. Consider this class of subsets of $\Omega$ :

\begin{equation*}
\mathcal{A}_{A}=\{B\cap A,B\in \mathcal{A}\}.
\end{equation*}

\noindent S (a) Show that $\mathcal{A}_{A}$ is a $\sigma $-algebra on the space $A$, called induced $\sigma $-algebra by $\mathcal{A}$ on $A$.\\

\noindent (b) show that is also the $\sigma $-algebra generated by the canonical injection

\begin{equation*}
\begin{tabular}{llll}
$i_{A}:$ & $A$ & $\longmapsto $ & $\Omega $ \\ 
& $x$ & $\hookrightarrow $ & $x$.
\end{tabular}
\end{equation*}

\bigskip \bigskip \noindent \textbf{Exercise 8} \label{exercise08_doc01-02}  Let $\Omega$ be a nonempty set and let $\mathcal{D}$ a collection of subsets of $\Omega$. We are going to justify the existence of $a(\mathcal{D})$.

\bigskip \noindent a) Let $(\mathcal{C})_{(t\in T)}$ a non-empty family of algebras on $\Omega $. Show that the intersection of these algebras :

\begin{equation*}
\mathcal{C}_{\infty} = \bigcap\limits_{t\in T} \mathcal{C}_t.
\end{equation*}

\noindent S is an algebra.\newline

\bigskip \noindent b) (Application). Consider the family $(\mathcal{C})_{(t\in T)}$ of algebras containing $\mathcal{D}$. Say why this class is not empty. Show that the intersection $\mathcal{C}_{\infty }$ is the smallest algebra containing $\mathcal{D}$.\newline

\bigskip \noindent c) Show that $a(\mathcal{D})=\mathcal{D}$ if and only if $%
\mathcal{D}$ is an algebra.\newline

\bigskip

\bigskip \noindent \textbf{Exercise 9} \label{exercise09_doc01-02} (01.06) Let $\Omega$ be a nonempty set and let $\mathcal{D}$ a collection of subsets of $\Omega$. We are going to justify the existence of $\sigma(\mathcal{D})$.\newline

\bigskip \noindent a) Let $(\mathcal{A})_{(t\in T)}$ a non-empty family of $\sigma $-algebras on $\Omega $. Show that the intersection of these $\sigma$-algebras :

\begin{equation*}
\mathcal{C}_{\infty} = \bigcap\limits_{t\in T} \mathcal{A}_t.
\end{equation*}

\noindent S is a $\sigma$-algebra.\newline

\bigskip \noindent b) (Application). Consider the family $(\mathcal{A})_{(t\in T)}$ of $\sigma $-algebras containing $\mathcal{D}$. Say why this class is not empty. Show that the intersection $(\mathcal{A})_{\infty }$ is the smallest $\sigma $-algebra containing $\mathcal{D}$.\\

\bigskip \noindent c) Show that $\sigma (\mathcal{D})=\mathcal{D}$ if and only if $\mathcal{D}$ is a $\sigma $-algebra.\newline

\bigskip

\bigskip \noindent \textbf{Exercise 10}. \label{exercise10_doc01-02} Let $\Omega$ be a nonempty set and let $\mathcal{D}$ a collection of subsets of $\Omega$. Show that

\begin{equation*}
\sigma(\mathcal{D})=\sigma(a(\mathcal{D})).
\end{equation*}

\bigskip

\bigskip \noindent \textbf{Exercise 11} (Extensions). \label{exercise11_doc01-02} In the spirit of Exercises 8 and 9, justify the existence of smallest monotone class
generated by some nonempty class $\mathcal{D}$ denoted $m(\mathcal{D})$, and that of smallest Dynkin system generated by some nonempty class $\mathcal{D}$
denoted $d(\mathcal{D})$.\newline

\bigskip

\bigskip \noindent \textbf{Exercise 12}. \label{exercise12_doc01-02} Let $\mathcal{D}$\ be a Dynkin System. Show that $\mathcal{D}$ contains the limits of its non-decreasing
sequences.\newline

\bigskip \noindent \textbf{Exercise 13}. \label{exercise13_doc01-02} Let $X$ be an application from $\Omega _{1}$ to $\Omega _{2}$ and suppose that $\Omega _{2}$ is endowed with
an $\sigma $-algebra $\mathcal{A}_{2}$. Consider this class of subsets of $\Omega _{1}$. Define the inverse image of $B$ in $\Omega _{2}$ by $X$ as follows

\begin{equation*}
X^{-1}(B)=\{{\omega \in \Omega _{1},X(\omega )\in B\}}.
\end{equation*}

\bigskip \noindent Show that the inverse image preserves all sets operations in the following sense :

\bigskip \noindent (i) $X^{-1}(\Omega _{2})=\Omega _{1},$ $X^{-1}(\emptyset)=\emptyset $.\newline

\bigskip \noindent (ii) $\forall B\subset \Omega_{2},X^{-1}(B^{c})=X^{-1}(B)^{c}$.\newline

\bigskip \noindent (iii) $\forall B_{1}\subset \Omega _{2},B_{2}\subset\Omega _{2}$, 
\begin{equation*}
X^{-1}(B_{1}\cup B_{2})=X^{-1}(B_{1})\cup X^{-1}(B_{2})
\end{equation*}

\bigskip \noindent and
\begin{equation*}
,X^{-1}(B_{1}\cap B_{2})=X^{-1}(B_{1})\cap X^{-1}(B_{2}).
\end{equation*}

\bigskip \noindent (iv) $\forall B_{1}\subset \Omega _{2},B_{2}\subset \Omega _{2}$,

\begin{equation*}
X^{-1}(B_{1}\backslash B_{2})=X^{-1}(B_{1})\backslash X^{-1}(B_{2})
\end{equation*}

\bigskip \noindent and
\begin{equation*}
X^{-1}(B_{1}\Delta B_{2})=X^{-1}(B_{1})\Delta X^{-1}(B_{2}).
\end{equation*}

\bigskip \noindent (v) $\forall B_{1}\subset \Omega _{2},B_{2}\subset \Omega_{2},$ $B_{1}$ and $B_{2}$ disjoint

\begin{equation*}
X^{-1}(B_{1}+B_{2})=X^{-1}(B_{1})+X^{-1}(B_{2}).
\end{equation*}

\bigskip \noindent (vi) $\forall (B_{n})_{n\geq 0}\subset \mathcal{P}(\Omega_{2}),$ 
\begin{equation*}
X^{-1}(\cup _{n\geq 0}B_{n})=\cup _{n\geq 0}X^{-1}(B_{n})
\end{equation*}

\bigskip \noindent and
\begin{equation*}
X^{-1}(\cap _{n\geq 0}B_{n})=\cap _{n\geq 0}X^{-1}(B_{n}).
\end{equation*}

\bigskip \noindent (vii) $\forall (B_{n})_{n\geq 0}\subset \mathcal{P}(\Omega _{2}),$ pairwise disjoint,

\begin{equation*}
X^{-1}\left( \sum\limits_{n\geq 0}{}B_{n}\right) =\sum\limits_{n\geq
0}{}X^{-1}\left( B_{n}\right). 
\end{equation*}

\bigskip \noindent \textbf{Exercise 14}. \label{exercise14_doc01-02}\\

\noindent Let $\mathcal{F}$ be an arbitrary collection of disjoint elements of an algebra $\mathcal{C}$ of subsets of $\Omega$. Let $\mathcal{S}$ be the collection finite unions (sums) of elements of $\mathcal{F}$.\\

\noindent Show that $\mathcal{S}$ is stable under finite unions.\\ 

\newpage
\noindent \LARGE \textbf{DOC 01-03 : Exercises on Measurable sets \textit{with solutions}}. \label{doc01-03}\\

 \Large

\bigskip

\bigskip \noindent \textbf{Exercise 1}. \label{exercise01_sol_doc01-03}Show that the class $%
\mathcal{I}$ of all intervals of $\mathbb{R}$ is a semi-algebra.\newline

\bigskip \noindent \textbf{Solution}.\newline

\bigskip \noindent Check the three points.\newline

\bigskip \noindent (a) $\Omega =\mathbb{R}$ is an interval, so $\mathbb{R}%
\in \mathcal{I}.$\newline

\bigskip \noindent (b) The intersection of two intervals of $\mathbb{R}$ is
an interval of $\mathbb{R}$. So $\mathcal{I}$ is stable under finite
intersections. You should remind this formula :

\begin{eqnarray*}
\lbrack x_{1},x_{2}]\cap \lbrack y_{1},y_{2}] &=&[\max (x_{1},y_{1}),\min
(x_{2},y_{2})] \\
]x_{1},x_{2}]\cap ]y_{1},y_{2}] &=&]\max (x_{1},y_{1}),\min (x_{2},y_{2})] \\
\lbrack x_{1},x_{2}[\cap \lbrack y_{1},y_{2}[ &=&[\max (x_{1},y_{1}),\min
(x_{2},y_{2})[ \\
]x_{1},x_{2}[\cap \lbrack y_{1},y_{2}[ &=&]\max (x_{1},y_{1}),\min
(x_{2},y_{2})[.
\end{eqnarray*}

\bigskip \noindent \textbf{Rule}. Intersection of intervals of the same type
is an interval of that type.\newline

\bigskip \noindent (c) To find the the complements of intervals of $\mathbb{R%
}$, consider all the cases in the following table where numbers $a$ and $b$
are finite such that $a<b$

\begin{equation*}
\begin{tabular}{lll}
case & Different types of intervals & Their complements \\ 
00 & $]-\infty ,+\infty \lbrack $ & $\emptyset =]0,0[$ \\ 
01 & $\emptyset =]0,0[$ & $]-\infty ,+\infty \lbrack $ \\ 
1 & $]-\infty ,b[$ & $[b,+\infty \lbrack $ \\ 
2 & $]-\infty ,b]$ & $]b,+\infty \lbrack $ \\ 
3 & $]a,+\infty \lbrack $ & $]-\infty ,a]$ \\ 
4 & $[a,+\infty \lbrack $ & $]-\infty ,a[$ \\ 
5 & $[a,b]$ & $]-\infty ,a[+]b,+\infty \lbrack $ \\ 
6 & $]a,b[$ & $]-\infty ,a]+[b,+\infty \lbrack $ \\ 
7 & $]a,b]$ & $]-\infty ,a]+]b,+\infty \lbrack $ \\ 
8 & $[a,b[$ & $]-\infty ,a[+[b,+\infty \lbrack $.
\end{tabular}%
\end{equation*}

\bigskip \noindent And you see that the complement of each interval is
either an interval or a sum of two intervals.\newline

\bigskip \noindent \textbf{Exercise 2}. \label{exercise02_sol_doc01-03} Show that the class of
rectangles $\mathcal{S}$ is a semi-algebra. Extend this result to a product
of $k$ measurable spaces, $k\geq 2$.\newline

\bigskip \noindent \textbf{Solution}.\newline

\bigskip \noindent Let $(\Omega _{i},\mathcal{A}_{i})$ two measurable spaces
and let $\Omega =\Omega _{1}\times \Omega _{2}=\{(\omega _{1},\omega
_{2}),\omega _{1}\in \Omega _{1},\omega _{2}\in \Omega _{2}\}$ be the
Cartesian product of $\Omega _{1}$ and $\Omega _{2}.$ The class of cylinders
(or rectangles) is given by 
\begin{equation*}
\mathcal{S}=\{A_{1}\times A_{2},\text{ }A_{1}\in \mathcal{A}_{1},A_{2}\in 
\mathcal{A}_{2}\}.
\end{equation*}

\bigskip \noindent We have to show that $\mathcal{S}$ is a
semi-algebra.Check the three conditions.\newline

\bigskip \noindent (a) It is obvious that $\Omega =\Omega _{1}\times \Omega
_{2}\in \mathcal{S}$.\newline

\bigskip \noindent (b) Let $A=A_{1}\times A_{2}\in \mathcal{S},$with $%
A_{1}\in \mathcal{A}_{1},A_{2}\in \mathcal{A}_{2}$ and $B=B_{1}\times
B_{2}\in \mathcal{S},$ with $B_{1}\in \mathcal{A}_{1},B_{2}\in \mathcal{A}%
_{2}$. We have 
\begin{eqnarray*}
A\cap B &=&\{(\omega _{1},\omega _{2})\in \Omega ,(\omega _{1},\omega
_{2})\in A\text{ and }(\omega _{1},\omega _{2})\in B\} \\
&=&\{(\omega _{1},\omega _{2})\in \Omega ,(\omega _{1},\omega _{2})\in
A_{1}\times A_{2}\text{ and }(\omega _{1},\omega _{2})\in A_{1}\times A_{2}\}
\\
&=&\{(\omega _{1},\omega _{2})\in \Omega ,\omega _{1}\in A_{1}\text{ and }%
\times \omega _{2}\in A_{2}\text{ and }\omega _{1}\in B_{1}\text{ and }%
\omega _{2}\in B_{2}\} \\
&=&\{(\omega _{1},\omega _{2})\in \Omega ,(\omega _{1}\in A_{1}\text{ and }%
\times \omega _{1}\in B_{1})\text{ and (}\omega _{2}\in A_{2}\text{ and }%
\omega _{2}\in B_{2})\} \\
&=&(A_{1}\cap B_{1})\times (A_{2}\cap B_{2}).
\end{eqnarray*}

\bigskip \noindent Since 
\begin{equation*}
A_{1}\cap B_{1})\in \mathcal{A}_{1}\text{ and }(A_{2}\cap B_{2})\in \mathcal{%
A}_{2},
\end{equation*}

\bigskip \noindent we get that $A\cap B\in \mathcal{S}$.\newline

\bigskip \noindent (c) Finally, let $A=A_{1}\times A_{2}\in \mathcal{S},$ with $A_{1}\in \mathcal{A}_{1},A_{2}\in \mathcal{A}_{2}.$ We easily
establish that

\begin{equation*}
A^{c}=(A_{1}\times A_{2})^{c}=A_{1}^{c}\times A_{2}+(A_{1}\times
A_{2}^{c}+A_{1}^{c}\times A_{2}^{c}.
\end{equation*}

\bigskip \noindent So, the complement of an element of $\mathcal{S}$ is a sum of elements of $\mathcal{S}$.

\bigskip \noindent \textbf{Exercise 3}. \label{exercise03_sol_doc01-03}Le $\Omega $ be a nonempty set and
set

\begin{equation*}
\mathcal{A}=\{A\subset \Omega, \ \ A \ \ is \ \ finite \ \ or \newline
A^{c} \ \ is \ finite \}.
\end{equation*}

\bigskip \noindent (a) Show that $\mathcal{A}$ is an algebra.\newline

\noindent (b) Show that $\mathcal{A}$ is a $\sigma$-algebra if $\Omega$ is
finite.\newline

\noindent (c) If $\Omega $ is infinite, consider an infinite subset $C=\{{%
\omega _{0},\omega _{1},\omega _{2},..\}}$ where the ${\omega _{i}}$ are
dsijoint and set $A_{n}=\{{\omega _{2n}\}}$, $n\geq 0$. Show that the $A_{n}$
are in $\mathcal{A}$ and their union is not.\newline

\noindent (d) Conclude that $\mathcal{A}$ is a $\sigma $-algebra \textbf{if
and only if} $\Omega $ is finite.\newline

\bigskip \noindent \textbf{Solution}.\newline

\noindent (a) Let us check that $\mathcal{A}$ is an algebra.\newline

\noindent (i) $\Omega $ belongs to $\mathcal{A}$ since its
complement is empty and then finite.\newline

\noindent (ii) The roles of $A$ and $A^{c}$ are symmetrical in the
definition of $\mathcal{A}.$ So $A$ belongs to $\mathcal{A}$ whenever $A^{c}$
belongs to $\mathcal{A}$ and vice-verse.\newline

\noindent (iii) Let us show that $\mathcal{A}$ is stable under
finite unions. Let $A$ and $B$ be two elements of $\mathcal{A}.$ we have two
case. Either $A$ and $B$ are finite sets and then $A\cup B$ is finite and
hence $A\cup B$ belongs to $\mathcal{A}$. Or one of them is not finite.
Suppose that $A$ is not finite. Since $A$ is in $\mathcal{A}$ and is not
finite, this implies that $A^{c}$ is finite (otherwise A would not be in $%
\mathcal{A}$). Then, we have

\begin{equation*}
(A\cup B)^{c}=A^{c}\cap B^{c}\subset A^{c}.
\end{equation*}

\bigskip \noindent Thus $(A\cup B)^{c}$ is finite as part of a finite set.
Hence $A\cup B\in \mathcal{A}$. Then in all possibles cases, $A\cup B$
belongs to $\mathcal{A}$ whenever $A$ and $B$ are in $\mathcal{A}.$\newline

\noindent (b) Since $\mathcal{A}$ is already an algebra, if
suffices to show that is is stable under countable unions or intersections.
But if $\Omega $ is a finite set, the collection of subsets is also finite
and so, a countable union of elements of $\mathcal{A}$ is merely a finite
union of element of $\mathcal{A}$ and then remains in $\mathcal{A}.$ Thus,
it is clear that $\mathcal{A}$ is an sigma-algebra if $\Omega $ is finite.\ 

\noindent (c) Clearly the sets $A_{n}$ are finite and then, they
belong to $\mathcal{A} $. But their union

\begin{equation*}
A=\cup _{n\geq 0}A_{n}=\{\omega _{0},\omega _{2},...\}
\end{equation*}

\bigskip \noindent is neither finite nor co-finite since%
\begin{equation*}
\{\omega _{1},\omega _{3},...\}\subset A^{c}.
\end{equation*}

\bigskip \noindent (d) The conclusion is : $\mathcal{A}$ is a sigma-algebra
if and only if $\Omega $ is finite.\ 

\bigskip \noindent \textbf{Exercise 4}. \label{exercise04_sol_doc01-03} Let $\Omega$ be a nonempty set and
put

\begin{equation*}
\mathcal{A}=\{A\subset \Omega ,\ \ A\ \ is\ \ countable\ \ or\newline
\text{ }A^{c}\ \ is\ countable\}.
\end{equation*}

\noindent (a) S Show that $\mathcal{A}$ is a $sigma$-algebra.\newline

\noindent (b) Let $\Omega =[0,1]$. Show that $[0,1/2]$ is not measurable in
the measurable space $(\Omega ,\mathcal{A})$.

\bigskip

\bigskip \noindent \textbf{Solution}.\newline

\noindent (a) Let us check that $\mathcal{A}$ is an algebra.\newline

\noindent (i) $\Omega $ belongs to $\mathcal{A}$ since its
complement is empty and then countable.\newline

\noindent (ii) The roles of $A$ and $A^{c}$ are symmetrical in the
definition of $\mathcal{A}.$ So $A$ belongs to $\mathcal{A}$ whenever $A^{c}$
belongs to $\mathcal{A}$ and vice-verse.\newline

\noindent (iii) Let us show that $\mathcal{A}$ is stable under
countable unions. Let ($A_{n})_{n\geq 0}$ a sequence of elements of $%
\mathcal{A}.$ We have two cases. Either all the $A_{n}$ are countable, then
the countable union $\cup _{n\geq 0}A_{n}$ is countable and hence $\cup
_{n\geq 0}A_{n}$ belongs to $\mathcal{A}$. Or one of them, say $A_{n_{0}}$
is not countable. Since $A_{n_{0}}$ is in $\mathcal{A}$ and is not
countable, this implies that $A_{n_{0}}^{c}$ is countable (otherwise $%
A_{n_{0}}$ would not be in $\mathcal{A}$). Then%
\begin{equation*}
(\cup _{n\geq 0}A_{n})^{c}=\cap _{n\geq 0}A_{n}^{c}\subset A_{n_{0}}^{c}.
\end{equation*}

\bigskip \noindent Thus $(\cup _{n\geq 0}A_{n})^{c}$ is countable as part of
a countable set. Hence $\cup _{n\geq 0}A_{n}\in \mathcal{A}$.\newline

\noindent (b) In this case $A=[0,1/2]$ is not countable and $%
A^{c}=]1/2,1]$ is not countable. Conclusion : In the measurable space $%
([0,1],\mathcal{A})$, $[0,1/2]$ is not measurable.\newline

\noindent But the set of irrational numbers is measurable, \textbf{%
in this measurable space}, since its complement, the set of rational numbers
is countable.\newline

\bigskip \noindent \textbf{Exercise 5}. \label{exercise05_sol_doc01-03} Prove 01.14 : The algebra generated
by a semi-algebra \ $\mathcal{S}$ is the class of all finite sums of its
elements 
\begin{equation*}
a(\mathcal{S})=\{A_{1}+A_{1}+...+A_{k},\text{ }A_{i}\in \mathcal{S},\text{ }%
k\geq 1\}.
\end{equation*}

\bigskip \noindent \textbf{Solution}.\newline

\noindent Denote
\begin{equation*}
\mathcal{C}_{0}=\{A_{1}+A_{1}+...+A_{k},\text{ }A_{i}\in \mathcal{S},\text{ }%
k\geq 1\}.
\end{equation*}

\bigskip \noindent It is clear that any algebra including $\mathcal{S}$
contains finite unions of $\mathcal{S}$\ and then contains $\mathcal{C}_{0}$%
. This means that $\mathcal{C}_{0}$\ will be the algebra generated by $%
\mathcal{S}$\ \ whenever $\mathcal{C}_{0}$\ is an algebra. Let us show that $%
\mathcal{C}_{0}$\ is an algebra by checking :\newline

\noindent (i) $\Omega \in \mathcal{C}_{0}.$ This is immediate
since, by definition, $\Omega \in $ $\mathcal{S}$ and $\mathcal{C}_{0}$ is
composed with the elements of $\mathcal{C}_{0}$\ and the finite sums of
elements of $\mathcal{S}$.\newline

\noindent (ii) $\mathcal{C}_{0}$ is stable under finite intersection. This obviously derives from the mutual distributivity of the union with
respect to the intersection and to the finite stability $\mathcal{S}$ under intersection.  Indeed, let $A=A_{1}+...+A_{k}\in \mathcal{C}_{0}$ with $A_{i}\in \mathcal{S},$ for $i=1,...,k$ 
and $B=BA_{1}+...+B_{h}\in \mathcal{C}_{0}$ with $B0_{j}\in \mathcal{S}$ for $j=1,...,h$. We have

$$
A \bigcap B = \sum_{i=1}^{k} \sum_{j=1}^{h} A_iB_j \in \mathcal{C}_{0}. \ \square
$$

\bigskip \noindent (iii) Let $A=A_{1}+...+A_{k}\in $ $\mathcal{C}_{0}$ with $%
A_{i}\in \mathcal{S},$ for $i=1,...,k.$ Check that $A^{c}\in \mathcal{C}_{0}.
$.\noindent But
\begin{equation*}
A^{c}=\bigcap\limits_{i=1}^{k}A_{i}^{c}.
\end{equation*}

\bigskip \noindent By definition of the semi-algebra $\mathcal{S}$, each \ $%
A_{i}^{c}$\ is a finite sum \ of elements of \ \ $\mathcal{S},$ say%
\begin{equation*}
A_{i}^{c}=A_{i,1}+...+A_{i,m(i)}.
\end{equation*}

\bigskip \noindent Be careful, there is no reason that the $A_{i}^{c}$ are
sums of the same number of non empty elements of $\mathcal{S}$. \ Hence, we
have : 
\begin{equation*}
A^{c}=\bigcap\limits_{i=1}^{k}\left( A_{i,1}+...+A_{i,m(i)}\right) .
\end{equation*}

\bigskip \noindent To see that $A^{c}$ is a finite sum of elements of $\mathcal{S}$, we may make
the expansion of the right-hand member of the formula above for simple cases.\\

\noindent For $k=2$, $m(1)=1$ and $m(2)=2,$ we have%
\begin{equation*}
A_{1,1}\cap (A_{2,1}+A_{2,1})=A_{1,1}A_{2,1}+A_{1,1}A_{2,1}.
\end{equation*}

\bigskip \noindent For $k=2$, $m(1)=2$ and $m(2)=2,$ we have%
\begin{equation*}
(A_{1,1}+A_{1,2})\cap (A_{2,1}+A_{2,1})=A_{1,1}A_{2,1}+A_{1,1}A_{2,1}.
\end{equation*}

\bigskip \noindent For $k=2$, $m(1)=1$ and $m(2)=2,$ we have%
\begin{equation*}
A_{1,1}\cap
(A_{2,1}+A_{2,1})=A_{1,1}A_{2,1}+A_{1,1}A_{2,1}+A_{1,2}A_{2,1}+A_{1,2}A_{2,1}.
\end{equation*}

\bigskip \noindent We easily see that in the general case,%
\begin{equation*}
A^{c}=\bigcap\limits_{i=1}^{k}\left( A_{i,1}+...+A_{i,m(i)}\right) 
\end{equation*}

\bigskip \noindent is a finite sum of elements of the form%
\begin{equation*}
\prod\limits_{h=1}^{k}A_{i_{h},j_{h}}\in \mathcal{S}.
\end{equation*}

\bigskip \noindent where $i_{h}\in \{1,2,...,k\}$, $j_{h}\in \{1,2,...,m(i_{h})\}.$ So $A^{c}\in a(\mathcal{S}).$

\bigskip \noindent We conclude that $A^{c}\in \mathcal{C}_{0}$.\newline

\bigskip \noindent Since all the three points are checked, we say that $%
\mathcal{C}_{0}$ is an algebra including any algebra including $\mathcal{S}$. So%
\begin{equation*}
\mathcal{C}_{0}=a(\mathcal{S}).
\end{equation*}

\bigskip \noindent \textbf{Exercise 6}. \label{exercise06_sol_doc01-03} Let $X$ be an application from $%
\Omega _{1}$ to $\Omega _{2}$ and suppose that $\Omega _{2}$ is endowed with
an $\sigma $-algebra $\mathcal{A}_{2}$. Consider this class of subsets of $%
\Omega _{1}$ :

\begin{equation*}
\mathcal{A}_{X}=\{X^{-1}(B),B\in \mathcal{A}_{2}\},
\end{equation*}

\bigskip \noindent where $X^{-1}(B)$ denotes inverse image of $B$ by $X$ defined by

\begin{equation*}
X^{-1}(B)={\omega \in \Omega_1, X(\omega) \in B}.
\end{equation*}

\bigskip \noindent Show that $\mathcal{A}_{X}$ is a $\sigma $-algebra, called $\sigma 
$-algebra generated by $X$ (and by $\mathcal{A}_{2}$).\\

\bigskip \noindent \textit{Solution}. We recall, as we will revise it in
Exercise 13 below, that the inverse image preserves all sets operations. In
particular, we have\newline

\bigskip \noindent (i) $X^{-1}(\Omega _{2})=\Omega _{1}$.\newline

\noindent (ii) $\forall B\subset \Omega_{2},X^{-1}(B^{c})=X^{-1}(B)^{c}$.\newline

\noindent (iii) $\forall (B_{n})_{n\geq 0}\subset \mathcal{P}(\Omega _{2}),$ $X^{-1}(\cup _{n\geq 0}B_{n})=\cup _{n\geq 0}X^{-1}(B_{n})$.%
\newline

\noindent etc.\newline

\bigskip \noindent Based on these three points, we have\newline

\bigskip \noindent (a) $\Omega _{1}\in $ \ $\mathcal{A}_{X}$ since $\Omega
_{1}=X^{-1}(\Omega _{2})$ and $\Omega _{2}\in \mathcal{A}_{2}.$\newline

\noindent (b) Let $A\in $ \ $\mathcal{A}_{X}$, that is there exists 
$B\in \mathcal{A}_{2}$ such that $A=X^{-1}(B).$ Then $%
A^{c}=X^{-1}(B)^{c}=X^{-1}(B^{c})\in $ $\mathcal{A}_{X}$ since $B^{c}\in 
\mathcal{A}_{2}$.\newline

\noindent (c) Let $\forall (A_{n})_{n\geq 0}\subset \mathcal{A}_{X}$. Then for each $n\geq 1,$ there exists $B_{n}\in \mathcal{A}_{2}$
such that $A_{n}=X^{-1}(B_{n}).$ Then $\cup _{n\geq 0}A_{n}=\cup _{n\geq
0}X^{-1}(B_{n})=X^{-1}(\cup _{n\geq 0}B_{n})\mathcal{A}_{X}$ since $\cup
_{n\geq 0}B_{n}\in \mathcal{A}_{2}$.\newline

\noindent Conclusion : $\mathcal{A}_{X}$ is a $\sigma $-algebra and:\\

\noindent (1) Clearly $X$ is measurable from ($\Omega _{1},\mathcal{A}_{X})$ to ($%
\Omega _{2},\mathcal{A}_{2})$.\newline

\noindent (2) $X$ is measurable from ($\Omega _{1},\mathcal{A}_{1})$
to ($\Omega _{2},\mathcal{A}_{2})$ if and only if $\mathcal{A}_{X}\subset 
\mathcal{A}_{1}$.\newline

\noindent (3) Final conclusion : $\mathcal{A}_{X}$ is the smallest
sigma-algebra on $\Omega _{1}$ such that $X:\Omega _{1}\longmapsto (\Omega
_{2},\mathcal{A}_{2})$ is measurable.\newline

\bigskip \noindent \textbf{Exercise 7}. \label{exercise07_sol_doc01-03} $(\Omega,\mathcal{A})$ be a
measurable space and let $A$ a subset of $\Omega$. Consider this class of
subsets of $\Omega$ :

\begin{equation*}
\mathcal{A}_{A}=\{B\cap A,B\in \mathcal{A}\}.
\end{equation*}

\bigskip \noindent (a) Show that $\mathcal{A}_{A}$ is a $\sigma $-algebra on the space $A$, called induced $\sigma $-algebra by $\mathcal{A}$ on $A$.\newline

\noindent (b) show that is also the $\sigma $-algebra generated by the canonical injection

\begin{equation*}
\begin{tabular}{llll}
$i_{A}:$ & $A$ & $\longmapsto $ & $\Omega $ \\ 
& $x$ & $\hookrightarrow $ & $x$%
\end{tabular}%
\end{equation*}%
\newline

\bigskip \noindent \textbf{Solution}.\newline

\noindent (a) Check the three points.\newline

\bigskip \noindent (i) $A=\Omega \cap A$ and $\Omega \in \mathcal{A}$. Then $%
A\in \mathcal{A}_{A}$.\newline

\noindent (ii) Let $C=B\cap A\in \mathcal{A}_{A}$ with $B\in $ $%
\mathcal{A}.$\ Then the complement of $C$ in $\ A$ is $C^{c}=B^{c}\cap A\in 
\mathcal{A}_{A}$. (Make a drawing to see this last formula).\newline

\noindent (iii) Let $C_{n}=B_{n}\cap A\in \mathcal{A}_{A}$ with $%
B_{n}\in $ $\mathcal{A}$ for each $n\geq 1.$ Then, by distributivity of
union over intersection, $\cup _{n\geq 1}=C_{n}=\cup _{n\geq 1}(B_{n}\cap
A)=\left( \cup _{n\geq 1}B_{n}\right) \cup A\in \mathcal{A}_{A}$ since $\cup
_{n\geq 1}B_{n}\in \mathcal{A}$.\newline

\noindent (b) For any $B\in \mathcal{A},$%
\begin{equation*}
i_{A}^{-1}(B)=\{\omega \in A,i_{A}(\omega )\in B\}=\{\omega \in A,\omega \in
B\}=A\cap B.
\end{equation*}

\bigskip \noindent Then 
\begin{equation*}
\mathcal{A}_{A}=\{B\cap A,B\in \mathcal{A}\}=\{i_{A}^{-1}(B),B\in \mathcal{A}%
\}.
\end{equation*}

\bigskip \noindent So $\mathcal{A}_{A}$ is exactly the sigma-algebra
generated by the canonical inclusion of A into $\Omega$.\newline

\noindent \textbf{Exercise 8}.\label{exercise08_sol_doc01-03} Let $\Omega $ be a nonempty
set and let $\mathcal{D}$ a collection of subsets of $\Omega $. We are going
to justify the existence of $a(\mathcal{D})$.\\

\noindent a) Let $(\mathcal{C})_{(t\in T)}$ a non empty family of
algebras on $\Omega $. Show that the intersection of these algebras :

\begin{equation*}
\mathcal{C}_{\infty }=\bigcap\limits_{t\in T}\mathcal{C}_{t}.
\end{equation*}

\bigskip \noindent is an algebra.\newline

\noindent b) (Application). Consider the family $(\mathcal{C})_{(t\in
T)}$ of algebras including $\mathcal{D}$. Say why this class is not empty.
Show that the intersection $\mathcal{C}_{\infty }$ is the smallest algebra
including $\mathcal{D}$.\newline

\noindent c) Show that $a(\mathcal{D})=\mathcal{D}$ if and only if $%
\mathcal{D}$ is an algebra.\newline

\bigskip

\bigskip \noindent \textbf{Solution}.\newline

\noindent (a) We skip the proof that $\mathcal{C}_{\infty }$ is an
algebra since we will do it \ later in Exercise 9 below. From that solution,
you'll see how to adapt in this case.\newline

\noindent (b) Here, the collection of classes algebras of subsets  of $\Omega $ including $\mathcal{D}$ is not empty because the power set $\mathcal{P}(\Omega )$ is in the collection. So the intersection is
meaningful and $\mathcal{C}_{\infty }$ is an algebra. It contains $\mathcal{D%
}$ since all the classes forming the intersection contain it. Next, any
algebra including $\mathcal{D}$ is part of the intersection and then
contains $\mathcal{C}_{\infty }.$ Then $\mathcal{C}_{\infty }$ is the
smallest algebra including $\mathcal{D}$.\newline

\noindent (c) If $a(\mathcal{D})=\mathcal{D},$ clearly $\mathcal{D}$%
, as it is equal to a($\mathcal{D}\mathbb{)},$ is an algebra. Inversely, if $%
\mathcal{D}$\ is an algebra, it is surely the smallest algebra including
itself and then \ $\mathcal{D}=a(\mathcal{D})$.\newline

\bigskip \noindent \textbf{Exercise 9}. \label{exercise09_sol_doc01-03} Let $\Omega$ be a nonempty
set and let $\mathcal{D}$ a collection of subsets of $\Omega$. We are going
to justify the existence of $\sigma(\mathcal{D})$.\newline

\noindent a) Let $(\mathcal{A})_{(t\in T)}$ a non empty family of $%
\sigma $-algebras on $\Omega $. Show that the intersection of these $\sigma $%
-algebras :

\begin{equation*}
\mathcal{A}_{\infty }=\bigcap\limits_{t\in T}\mathcal{A}_{t}.
\end{equation*}

\noindent is a $\sigma$-algebra.\newline

\noindent b) (Application). Consider the family $(\mathcal{A})_{(t\in
T)}$ of $\sigma $-algebras including $\mathcal{D}$. Say why this class is
not empty. Show that the intersection $\mathcal{A}_{\infty }$ is the
smallest $\sigma $-algebra including $\mathcal{D}$

\noindent c) Show that $a(\mathcal{D})=\mathcal{D}$ if and only if $%
\mathcal{D}$ is a $\Omega $-algebra.\newline

\bigskip \noindent \textbf{Solution}.\newline

\noindent With respect to Exercise 8, the only point to prove is
Point (a). Points (b) and (c) are proved exactly as in Exercise 8.\newline

\noindent (a) Let us show that $\mathcal{C}_{\infty }$ is a
sigma-algebra. Let us Check the three points.\newline

\noindent (i) $\Omega \in \mathcal{C}_{\infty }.$ The sets $\Omega $
and $\emptyset $ are in each $\mathcal{C}_{t},$ $t\in T,$ since they are
sigma-algebras. So, for example, $\Omega $ $\in \mathcal{C}_{\infty}$.%
\newline

 \noindent (ii) Let $A\in \mathcal{C}_{\infty }.$ Prove that $%
A^{c}\in \mathcal{C}_{\infty }.$ But $A\in \mathcal{C}_{\infty }$ means 
\begin{equation*}
\forall (t\in T),A\in \mathcal{C}_{t},
\end{equation*}

\bigskip \noindent which , by using that : $A\in \mathcal{C}%
_{t}\Longrightarrow A^{c}\in \mathcal{C}_{t}$ since $\mathcal{C}_{t}$ is a
sigma-algebra, implies 
\begin{equation*}
\forall (t\in T),A^{c}\in \mathcal{C}_{t}.
\end{equation*}

\bigskip \noindent This entails, by definition,%
\begin{equation*}
A^{c}\in \mathcal{C}_{\infty}.
\end{equation*}

\bigskip \noindent (iii) Let $A_{n}\in \mathcal{C}_{\infty },$ $n\geq 1.$
Prove that $\cup _{n\geq 1}A_{n}\in \mathcal{C}_{\infty }.$ But $A_{n}\in 
\mathcal{C}_{\infty },$ for all $n\geq 1$ means%
\begin{equation*}
\forall (n\geq 1),A_{n}\in \mathcal{C}_{\infty},
\end{equation*}

\bigskip \noindent which implies, by the definition of $\mathcal{C}_{\infty
} $,%
\begin{equation*}
\forall (n\geq 1),\forall (t\in T),A_{n}\in \mathcal{C}_{t}.
\end{equation*}

\bigskip \noindent Since the ranges of the two quantifiers $\forall (t\in T)$
and $\forall (n\geq 1)$ are independent one from the other, we may swap them
to get\newline

\begin{equation*}
\forall (t\in T)\text{ }\forall (n\geq 1)\text{ },\text{ }A_{n}\in \mathcal{C}_{t},
\end{equation*}

\bigskip \noindent and use the fact that each $\mathcal{C}_{t}$ is a
sigma-algebra to have
\begin{equation*}
\forall (t\in T)\text{ },\text{ }(\cup _{n\geq 1}A_{n})\in \mathcal{C}_{t},
\end{equation*}

\bigskip \noindent and this implies $(\cup _{n\geq 1}A_{n})$ $\in \mathcal{C}%
_{\infty }$. Question (a) is entirely proved.\newline

\noindent (b) Here, the collection of classes sigma-algebras of subsets of $\Omega $ including $\mathcal{D}$ is not empty because the power
set $\mathcal{P}(\Omega )$ in the collection. So the intersection is meaningful and $\mathcal{C}_{\infty }$ is a sigma-algebra. It contains $%
\mathcal{D}$ since all the classes forming the intersection contain it. Next, any sigma-algebra including $\mathcal{D}$ is part of the intersection
and then contains $\mathcal{C}_{\infty }.$ Then $\mathcal{C}_{\infty }$ is
the smallest sigma-algebra including.\newline

\noindent (c) If $\sigma (\mathcal{D})=\mathcal{D},$ clearly $\mathcal{D}$, as it is equal to $\sigma (\mathcal{D}\mathbb{)},$ is a
sigma-algebra. Conversely, if $\mathcal{D}$\ is a sigma-algebra, it is surely the smallest sigma-algebra including itself and then \ $\mathcal{D}=\sigma (%
\mathcal{D})$.\newline

\bigskip \noindent \textbf{Exercise 10}. \label{exercise11_sol_doc01-03} Let $\Omega $ be a nonempty set and let $\mathcal{D}$ a collection of subsets of $\Omega $. Show that

\begin{equation*}
\sigma (\mathcal{D})=\sigma (a(\mathcal{D})).
\end{equation*}

\bigskip \noindent \textbf{Solution}.\newline

\noindent (a) Direct inclusion. $a(\mathcal{D})$ contains $\mathcal{%
D}$. So the smallest sigma-algebra including $\mathcal{D},$ that is $\sigma
(a(\mathcal{D})),$ contains $\mathcal{D}$. Since $\sigma (a(\mathcal{D}))$\
is a sigma-algebra including $\mathcal{D}$, it contains the smallest
sigma-algebra including $\mathcal{D}$, that is : \ $\sigma (\mathcal{D}%
)\subset \sigma (a(\mathcal{D}))$.\\

\noindent (b) Indirect inclusion : $\sigma (\mathcal{D})$ is an
algebra including $\mathcal{D}$\ (do not forget that a sigma-algebra is also
an algebra), so its contains the smallest algebra including $a(\mathcal{D})$%
, that is : $a(\mathcal{D})$ $\subset \sigma (\mathcal{D}).$ Now $\sigma (%
\mathcal{D})\ $is a sigma-algebra including $a(\mathcal{D})$, then it
contains $\sigma (a(\mathcal{D}))$, that is : $\sigma (a(\mathcal{D}%
))\subset \sigma (\mathcal{D}).$\newline

\noindent Putting together the two points leads to the desired
conclusion.\newline

\bigskip \noindent \textbf{Exercise 11} (Extensions).  \label{exercise11_sol_doc01-03}In the spirit of
Exercises 8 and 9, justify the existence of smallest monotone class
generated by some nonempty class $\mathcal{D}$ denoted $m(\mathcal{D})$, and
that of smallest Dynkin system generated by some nonempty class $\sigma (%
\mathcal{D})$ denoted $d(\mathcal{D})$.

\bigskip \noindent \emph{Solution}.\newline

\noindent \emph{(a) Generated monotone class}.\newline

\noindent Let $\mathcal{D}$\ a non-empty class of subsets of $%
\Omega $\ and denote $\left\{ \mathcal{M}_{t},t\in T\right\} $ the
collection of monotone classes including \ $\mathcal{D}$. This class is
never empty since it contains the power set. Define

\begin{equation*}
\mathcal{M}_{\infty }=\bigcap\limits_{t\in T}\mathcal{M}_{t}.
\end{equation*}

\bigskip \noindent At the light of Exercises 8 and 9, it suffices to prove that $%
\mathcal{M}_{\infty }$ is a monotone class to get that, it is itself, the
smallest monotone class including $\mathcal{D}$, that is : $\mathcal{M}%
_{\infty }=m(\mathcal{M}).$\newline

\noindent Let us prove that $\mathcal{M}_{\infty }$\ is monotone.
We already know that it is not empty (including $\mathcal{D}$\ that is not
empty). We have to check these two points\newline

\noindent (i) Consider a non-decreasing sequence $\left(
A_{n}\right) _{n\geq 1}\subset \mathcal{M}_{\infty }.$ Do we have $\cup
_{n\geq 1}A_{n}\in \mathcal{M}_{\infty }?$ But $\left( A_{n}\right) _{n\geq
1}\subset \mathcal{M}_{\infty }$ means 
\begin{equation*}
\forall (n\geq 1),A_{n}\in \mathcal{M}_{\infty }.
\end{equation*}

\bigskip \noindent which implies, by the definition of $\mathcal{M}_{\infty} $,
\begin{equation*}
\forall (n\geq 1),\forall (t\in T),A_{n}\in \mathcal{M}_{t}.
\end{equation*}

\bigskip \noindent Since the ranges of the two quantifiers $\forall (t\in T)$
and $\forall (n\geq 1)$ are independent one from the other, we may swap them
to get

\begin{equation*}
\forall (t\in T)\text{ }\forall (n\geq 1)\text{ },\text{ }A_{n}\in \mathcal{M%
}_{t},
\end{equation*}

\bigskip \noindent and use the fact that each $\mathcal{M}_{t}$ is a monotone%
\begin{equation*}
\forall (t\in T)\text{ },\text{ }(\cup _{n\geq 1}A_{n})\in \mathcal{C}_{t},
\end{equation*}

\bigskip \noindent and this implies $(\cup _{n\geq 1}A_{n})$ $\in \mathcal{M}_{\infty }$. Question (i) is answered.\newline

\bigskip \noindent (ii) Consider a non-increasing sequence $\left(
A_{n}\right) _{n\geq 1}\subset \mathcal{M}_{\infty }.$ Do we have $\cap
_{n\geq 1}A_{n}\in \mathcal{M}_{\infty }?$ But $\left( A_{n}\right) _{n\geq
1}\subset \mathcal{M}_{\infty }$ means 
\begin{equation*}
\forall (n\geq 1),A_{n}\in \mathcal{M}_{\infty},
\end{equation*}

\bigskip \noindent which implies, by the definition of $\mathcal{M}_{\infty}$,
\begin{equation*}
\forall (n\geq 1),\forall (t\in T),A_{n}\in \mathcal{M}_{t}.
\end{equation*}

\bigskip \noindent Since the ranges of the two quantifiers $\forall (t\in T)$
and $\forall (n\geq 1)$ are independent one from the other, we may swap them
to get

\begin{equation*}
\forall (t\in T)\text{ }\forall (n\geq 1)\text{ },\text{ }A_{n}\in \mathcal{M%
}_{t},
\end{equation*}

\bigskip \noindent and use the fact that each $\mathcal{M}_{t}$ is a monotone%
\begin{equation*}
\forall (t\in T)\text{ },\text{ }(\cap _{n\geq 1}A_{n})\in \mathcal{C}_{t},
\end{equation*}

\bigskip \noindent and this implies $(\cap _{n\geq 1}A_{n})$ $\in \mathcal{M}%
_{\infty }$. Question (ii) is answered.\newline

\noindent \textbf{(a) Generated Dynkyn system}.\newline

\noindent Let $\mathcal{D}$\ a non-empty class of subsets of $%
\Omega $\ and denote $\left\{ \mathcal{G}_{t},t\in T\right\} $ the
collection of Dynkin systems including \ $\mathcal{D}$. This class is never
empty since in contains the power set. Define

\begin{equation*}
\mathcal{G}_{\infty }=\bigcap\limits_{t\in T}\mathcal{G}_{t}.
\end{equation*}

\noindent It will suffice to prove that \ $\mathcal{G}_{\infty }$
is Dynkin system, to get that, it is itself, the smallest Dynkin system
including $\mathcal{D}$, that is : $\mathcal{G}_{\infty }=d(\mathcal{D})$.%
\newline

\noindent Let us prove that $\mathcal{D}_{\infty }$\ is Dynkin
system. We already know that it is not empty (including $\mathcal{D}$\ that
is not empty). We have to check three points :

\noindent (i) $\Omega \in \mathcal{G}_{\infty }.$ But $\Omega $ and 
$\emptyset $ are in each $\mathcal{G}_{t},$ $t\in T,$ since they are Dynkin
systems. So, for example, $\Omega $ $\in \mathcal{G}_{\infty}$.\newline

\noindent (ii) Let $A\in \mathcal{G}_{\infty },B\in \mathcal{G}%
_{\infty }$ with $A\subset B.$ Prove that $B\setminus A\in \mathcal{G}%
_{\infty }.$ $\ $\ But $A\in \mathcal{G}_{\infty }$ and $B\in \mathcal{G}%
_{\infty }$ with $A\subset B$ means 

\begin{equation*}
\forall (t\in T),A\in \mathcal{G}_{t}\text{ and }B\in \mathcal{G}_{t}\text{
with }A\subset B,
\end{equation*}

\bigskip \noindent which implies, by the fact that each $\mathcal{G}_{t}$ is
a Dynkin system, 
\begin{equation*}
\forall (t\in T),B\setminus A\in \mathcal{G}_{t}.
\end{equation*}

\bigskip \noindent This entails, by definition of $\mathcal{G}_{\infty }$,%
\begin{equation*}
B\setminus A\in \mathcal{G}_{\infty }.
\end{equation*}

\bigskip \noindent and Question (ii) is closed.\newline

\noindent (iii) Consider a pairwise disjoint $\left( A_{n}\right)
_{n\geq 1}\subset \mathcal{G}_{\infty }.$ Do we have $\cup _{n\geq
1}A_{n}\in \mathcal{G}_{\infty }?$ But $\left( A_{n}\right) _{n\geq
1}\subset \mathcal{G}_{\infty } $ means 
\begin{equation*}
\forall (n\geq 1),A_{n}\in \mathcal{G}_{\infty }.
\end{equation*}

\bigskip \noindent which implies, by the definition of $\mathcal{G}_{\infty
} $,%
\begin{equation*}
\forall (n\geq 1),\forall (t\in T),A_{n}\in \mathcal{G}_{t}.
\end{equation*}

\bigskip \noindent Since the ranges of the two quantifiers $\forall (t\in T)$
and $\forall (n\geq 1)$ are independent one from the other, we may swap them
to get

\begin{equation*}
\forall (t\in T)\text{ }\forall (n\geq 1)\text{ },\text{ }A_{n}\in \mathcal{G%
}_{t},
\end{equation*}

\bigskip \noindent and use the fact that each $\mathcal{G}_{t}$ is a Dynkin
system class%
\begin{equation*}
\forall (t\in T)\text{ },\text{ }(\cup _{n\geq 1}A_{n})\in \mathcal{G}_{t},
\end{equation*}

\bigskip \noindent and this implies $(\cup _{n\geq 1}A_{n})$ $\in \mathcal{G}%
_{\infty }$. Question (iii) is answered.\newline

\bigskip \noindent \textbf{Exercise 12}. \label{exercise12_sol_doc01-03}Let $\mathcal{D}$\ be a Dynkin
System. Show that $\mathcal{D}$ contain the limits of its non-decreasing
sequences.\newline

\bigskip \noindent \textbf{Solution}. For a non-decreasing sequence $%
(A_{n})_{n\geq 0})$ of elements of , you easily have

\begin{equation*}
\bigcup\limits_{n\geq 1}A_{n}=A_{0}+\sum_{j=1}^{\infty }(A_{j}\setminus
A_{j-1}).
\end{equation*}

\bigskip \noindent Since the $(A_{j}\setminus A_{j-1})^{\prime }s$ are in $%
\mathcal{D}$ by Points (ii), we get that $\ \bigcup\limits_{n\geq 1}A_{n}$
is also a sum of pairwise disjoint elements of \ $\mathcal{D}$. And Point
(iii) implies 
\begin{equation*}
\bigcup\limits_{n\geq 1}A_{n}\in \mathcal{D}.
\end{equation*}

\bigskip 
\bigskip \noindent \textbf{Exercise 13}. \label{exercise13_sol_doc01-03} Let $X$ be an application from $%
\Omega _{1}$ to $\Omega _{2}$ and suppose that $\Omega _{2}$ is endowed with
an $\sigma $-algebra $\mathcal{A}_{2}$. Consider this class of subsets of $%
\Omega _{1}$. Define the inverse image of $B$ in $\Omega _{2}$ by $X$ as
follows

\begin{equation*}
X^{-1}(B)=\{{\omega \in \Omega _{1},X(\omega )\in B\}}.
\end{equation*}

\bigskip \noindent Show that the inverse image preserves all sets operations in the following sense :\\

\noindent (i) $X^{-1}(\Omega_{2})=\Omega_{1}$, \ \  $X^{-1}(\emptyset)=\emptyset $.\newline

\noindent (ii) $\forall B\subset \Omega _{2},X^{-1}(B^{c})=X^{-1}(B)^{c}$.\newline

\noindent (iii) $\forall B_{1}\subset \Omega _{2},B_{2}\subset \Omega _{2}$,

\begin{equation*}
X^{-1}(B_{1}\cup B_{2})=X^{-1}(B_{1})\cup X^{-1}(B_{2})
\end{equation*}

\bigskip \noindent and

\begin{equation*}
X^{-1}(B_{1}\cap B_{2})=X^{-1}(B_{1})\cap X^{-1}(B_{2}).
\end{equation*}%

\bigskip \noindent (iv) $\forall B_{1}\subset \Omega _{2},B_{2}\subset \Omega _{2},$
\begin{equation*}
X^{-1}(B_{1}\setminus B_{2})=X^{-1}(B_{1})\setminus X^{-1}(B_{2})
\end{equation*}

\bigskip \noindent and
\begin{equation*}
X^{-1}(B_{1}\Delta B_{2})=X^{-1}(B_{1})\Delta X^{-1}(B_{2}).
\end{equation*}

\bigskip \noindent (v) $\forall B_{1}\subset \Omega _{2},B_{2}\subset \Omega
_{2},$ $B_{1}$ and $B_{2}$ disjoint%
\begin{equation*}
X^{-1}(B_{1}+B_{2})=X^{-1}(B_{1})+X^{-1}(B_{2}).
\end{equation*}

\bigskip \noindent (vi) $\forall (B_{n})_{n\geq 0}\subset \mathcal{P}(\Omega
_{2}),$ 
\begin{equation*}
X^{-1}(\cup _{n\geq 0}B_{n})=\cup _{n\geq 0}X^{-1}(B_{n})
\end{equation*}

\bigskip \noindent and

\begin{equation*}
X^{-1}(\cap _{n\geq 0}B_{n})=\cap _{n\geq 0}X^{-1}(B_{n}).\newline
\end{equation*}

\bigskip \noindent (vii) $\forall (B_{n})_{n\geq 0}\subset \mathcal{P}(\Omega _{2})$, pairwise disjoint,

\begin{equation*}
X^{-1}\left( \sum\limits_{n\geq 0}{}B_{n}\right) =\sum\limits_{n\geq
0}{}X^{-1}\left( B_{n}\right) 
\end{equation*}

\bigskip 
\noindent \textbf{Solutions}.\\

\noindent  (i) We have 
\begin{equation*}
X^{-1}(\Omega _{2})=\{ \omega \in \Omega_1, \ X(\omega) \in \Omega_{2} \}.
\end{equation*}

\noindent We clearly have $\{\omega \in \Omega _{1}, \ X(\omega )\in \Omega_2\}=\Omega_{1}$ since $X$ is an application (mapping). Next%
\begin{equation*}
X^{-1}(\emptyset )=\{{\omega \in \Omega _{1},X(\omega)\in \emptyset \}}
\end{equation*}

\noindent and, clearly $\{{\omega \in \Omega _{1},X(\omega )\in \emptyset \}}$ since $X
$ is an application and then $X(\omega )$ exists for any $\omega \in \Omega
_{1}$.\\

\noindent (ii) Let $\forall B\subset \Omega _{2}.$ We have%
\begin{eqnarray*}
X^{-1}(B^{c}) &=&\{{\omega \in \Omega _{1},X(\omega )\in B}^{c}{\}} \\
&=&\{{\omega \in \Omega _{1},X(\omega )\notin B\}} \\
&=&\{{\omega \in \Omega _{1},\omega \notin X^{-1}(B)\}} \\
&=&X^{-1}(B)^{c}.
\end{eqnarray*}

\bigskip \noindent (iiia) Let $B_{1}\subset \Omega _{2},B_{2}\subset \Omega _{2}$. We have
\begin{eqnarray*}
X^{-1}(B_{1}\cup B_{2}) &=&\{{\omega \in \Omega _{1},X(\omega )\in B_{1}\cup
B_{2}\}} \\
&=&\{{\omega \in \Omega _{1},X(\omega )\in B_{1}}\text{ or }{X(\omega )\in
B_{2}\}} \\
&=&\{{\omega \in \Omega _{1},X(\omega )\in B_{1}\}\cup }\{{\omega \in \Omega
_{1},}\text{ }{X(\omega )\in B_{2}\}} \\
&=&X^{-1}(B_{1})\cup X^{-1}(B_{2}).
\end{eqnarray*}

\bigskip \noindent (iiib) Let $B_{1}\subset \Omega _{2},B_{2}\subset \Omega _{2}.$ we may prove
the formula as we did it for Point (iiia). But we are going to combine Points (ii) and (iiia). We have%
\begin{eqnarray*}
X^{-1}(B_{1}\cap B_{2}) &=&\left( \left( X^{-1}(B_{1}\cap B_{2})\right)
^{c}\right) ^{c} \\
&=&\left( X^{-1}(\left( B_{1}\cap B_{2})\right) ^{c}\right) ^{c} \\
&=&\left( X^{-1}(B_{1}^{c}\cap B_{2}^{c})\right) ^{c}\text{ \ \ by Point (ii)%
} \\
&=&\left( X^{-1}(B_{1}^{c})\cup X^{-1}(B_{2}^{c})\right) ^{c}\text{ by Point
(iiia)} \\
&=&\left( X^{-1}(B_{1}^{c}))^{c}\cap (X^{-1}(B_{2}^{c})\right) ^{c} \\
&=&X^{-1}(B_{1})\cap X^{-1}(B_{2})\text{ by Point (ii)} \\
&=&X^{-1}(B_{1})\cap X^{-1}(B_{2}).
\end{eqnarray*}%

\bigskip \noindent (iva) Let $B_{1}\subset \Omega _{2},B_{2}\subset \Omega_{2}$. We have

\begin{eqnarray*}
X^{-1}(B_{1}\setminus B_{2}) &=&X^{-1}(B_{1}\cap B_{2}^{c}) \\
&=&X^{-1}(B_{1})\cap X^{-1}(B_{2}^{c})\text{ by Point (iiib)} \\
&=&X^{-1}(B_{1})\cap X^{-1}(B_{2})^{c}\text{ by Point (ii)} \\
&=&X^{-1}(B_{1}) \setminus X^{-1}(B_{2}).
\end{eqnarray*}

\noindent (ivb) Next
\begin{eqnarray*}
X^{-1}(B_{1}\Delta B_{2}) &=&X^{-1}((B_{1} \setminus B_{2})\cup
(B_{2} \setminus B_{1})) \\
&=&X^{-1}(B_{1} \setminus B_{2})\cup X^{-1}(B_{2} \setminus B_{1})\text{ by
Point (iiia)} \\
&=&X^{-1}(B_{1}) \setminus X^{-1}(B_{2})\text{ }\cup \text{ }%
X^{-1}(B_{2}) \setminus X^{-1}(B_{1})\text{ by Point (iva)} \\
&=&X^{-1}(B_{1})\Delta X^{-1}(B_{2}).
\end{eqnarray*}

\bigskip \noindent (v)  Let $B_{1}\subset \Omega _{2},B_{2}\subset \Omega
_{2},$ $B_{1}$ and $B_{2}$ disjoint. Then%
\begin{eqnarray*}
X^{-1}(B_{1})\cap X^{-1}(B_{2}) &=&X^{-1}(B_{1}\cap B_{2})\text{ by Point
(iiib)} \\
&=&X^{-1}(\emptyset )=\emptyset.
\end{eqnarray*}

\noindent Then
\begin{eqnarray*}
X^{-1}(B_{1}+B_{2}) &=&X^{-1}(B_{1})\cup X^{-1}(B_{2}) \\
&=&X^{-1}(B_{1})+X^{-1}(B_{2}).
\end{eqnarray*}

\bigskip \noindent (via) Let $(B_{n})_{n\geq 0}\subset \mathcal{P}(\Omega
_{2}).$ Then 
\begin{eqnarray*}
X^{-1}(\cup _{n\geq 0}B_{n}) &=&\{{\omega \in \Omega _{1},X(\omega )\in \cup
_{n\geq 0}B_{n}\}} \\
&=&\{{\omega \in \Omega _{1},(\exists n\geq 0),}\text{ }{X(\omega )\in
B_{n}\}} \\
&=&\cup _{n\geq 0}\{{\omega \in \Omega _{1},(\exists n\geq 0),}\text{ }{%
X(\omega )\in B_{n}\}} \\
&=&\cup _{n\geq 0}X^{-1}(B_{n}).
\end{eqnarray*}

\noindent (vib) Let $(B_{n})_{n\geq 0}\subset \mathcal{P}(\Omega _{2}).$ Then%
\begin{eqnarray*}
X^{-1}(\cap _{n\geq 0}B_{n}) &=&\left( \left( X^{-1}(\cap _{n\geq
0}B_{n})\right) ^{c}\right) ^{c} \\
&=&\left( X^{-1}(\cup _{n\geq 0}B_{n}^{c})\right) ^{c}\text{ by Point (ii)}
\\
&=&\left( \cup _{n\geq 0}X^{-1}(B_{n}^{c})\right) ^{c}\text{ by Point (via)}
\\
&=&\cap _{n\geq 0}(X^{-1}(B_{n}^{c})^{c} \\
&=&\cap _{n\geq 0}X^{-1}(B_{n}).\newline
\text{ by Point (ii)}
\end{eqnarray*}

\bigskip \noindent (vii) $\forall (B_{n})_{n\geq 0}\subset \mathcal{P}%
(\Omega _{2}),$ pairwise disjoint. By combining Points (iiia) and (v), we
easily have

\begin{equation*}
X^{-1}\left( \sum\limits_{n\geq 0}{}B_{n}\right) =\sum\limits_{n\geq
0}{}X^{-1}\left( B_{n}\right) .
\end{equation*}

\bigskip \noindent \textbf{Exercise 14}. \label{exercise14_sol_doc01-03}\\

\noindent Let $\mathcal{F}$ be an arbitrary collection of disjoint elements of an algebra $\mathcal{C}$ of subsets of $\Omega$. Let $\mathcal{S}$ be
the collection finite unions (sums) of elements of $\mathcal{F}$.\\

\noindent Show that $\mathcal{S}$ is stable under finite unions.\\

\bigskip \noindent \textbf{SOLUTIONS}.\\

\noindent Let $A=\sum_{i \in I} A_i \in \mathcal{S}$ with $A_i \in \mathcal{F}$, $i \in I$, and $A=\sum_{j \in J} B_j \in \mathcal{S}$ with $B_j \in \mathcal{F}$, $j \in J$. Let us define

$$
I_0=\{i \in I, \ A_i \notin \{B_j, \ j\in J\} \} \ and \ J_0=\{J \in J, \ B_j \notin \{A_i, \ i \in I\} \}.
$$

\bigskip \noindent We may alternatively say : $I_0$ is the subset of indexes $i\in I$ such that $A_i$ is one of the $B_j$, $j\in J$, and as well, $J_0$ is the subset of indexes $j\in J$ such that $B_j$ is one of the $A_i$, $i\in I$. Since all the sets $A_i$ and $B_j$, $(i\in I, \ j\in J)$, are mutually disjoints, we have

$$
A \cup B= \sum_{i \in I} A_i + \sum_{j \in (J\setminus J_0)} B_j \in \mathcal{S}. \ \blacksquare 
$$

\newpage
\noindent \LARGE \textbf{DOC 01-04 : Measurable sets. Exercises on $\lambda$ and $\pi$ systems}. \label{doc01-04}\\

 \Large

\bigskip \noindent \textbf{Abstract} Monotone and Dynkin classes are important tools
for theoretical demonstrations in Measure Theory. It is important to be
familiar with them. Their importance reside in Exercises 3 and 4 below.\newline

\bigskip \noindent You will need these basic and easy facts. A $\sigma $-algebra is
an algebra, is a Dynkin system and is a monotone class. An algebra is a $\pi 
$-system.\newline

\bigskip \noindent \textbf{Exercise 1}. \label{exercise01_doc01-04} Let $\mathcal{C}$ be an
algebra and a monotone class. Show that $\mathcal{C}$ is a sigma-algebra.%
\newline

\bigskip \bigskip \noindent \textbf{Exercise 2}. \label{exercise02_doc01-04} Let $\mathcal{D}$ be a $\pi 
$-class and a Dynkin system. Show that $\mathcal{D}$\ is a sigma-algebra.%
\newline

\bigskip \noindent \textbf{Exercise 3}. \label{exercise03_doc01-04} Let $\mathcal{C}$ be an algebra of
subsets of $\Omega $. Prove that $\sigma (\mathcal{C})=m(\mathcal{C})$.%
\newline

\bigskip \noindent \textbf{Proposed solution (in Lo\`{e}ve)}

\bigskip \noindent Denote $\mathcal{M}=m(\mathcal{C})$\\

\bigskip \noindent Define for any subset $A$ of $\Omega$

$\mathcal{M}_{A}=\{B\in \mathcal{M},AB^{c}\in \mathcal{M},A^{c}B\in \mathcal{M},AB\in \mathcal{M}\}$.

\bigskip \noindent \textbf{A1.}  Show that for any $A\in \mathcal{M}$, $\mathcal{M}_{A}$ is a monotone class containing $A$.%
\newline

\noindent \textbf{A2.(i)} Show that for any $A\in \mathcal{C}$, $\mathcal{C} \subset \mathcal{M}_{A}$. Deduce from this, that 
 $\mathcal{M}_{A}=\mathcal{M}$.\\

\noindent \textbf{A2.(ii)} Let $B \in \mathcal{M}$. Use the symmetry of the roles of $A$ and $B$ in the definition $B \in \mathcal{M}_{A}$ and 
A2.(i) to prove that $A \in \mathcal{M}_{B}$ for any $A \subset \mathcal{C}$, and then $\mathcal{C} \subset \mathcal{M}_{B}$. Conclude that
$\mathcal{M}_{B}=\mathcal{M}$.\newline

\bigskip \noindent \textbf{A3} Conclude that for any $(A,B) \in \mathcal{M}^2$, 
$$
AB \in \mathcal{M}, \text{ } AB^{c} \in \mathcal{M}, \text{ } A^{c}B \in \mathcal{M}
$$

\bigskip \noindent and prove from this, $\mathcal{M}$ is a $\sigma $-algebra.\newline

\bigskip \noindent \textbf{A4.} Conclude.\newline

\bigskip \bigskip \noindent \textbf{Exercise 4}. \label{exercise04_doc01-04} Let $G$ be a $\pi $-system
of subsets of $\Omega$. Show that : $\sigma (\mathcal{G})=d(\mathcal{G}$).
\newline

\noindent \textbf{Proposed solution}.\newline

\noindent \textbf{Notation : } : Denote $\mathcal{D}$ as the Dynkin system
generated by $\mathcal{G}$. Put $\mathcal{D}_{A}=\{B\in \mathcal{D},AB\in 
\mathcal{D}\}$.\newline

\bigskip \noindent \textbf{B1.} Show that $\mathcal{D}_{A}$ is a Dynkin
system for any $A\in \mathcal{D}.$\newline

\noindent \textbf{B2.} Prove that for any $A\in \mathcal{G}$, $\mathcal{D}%
_{A}$ is a Dynkin system including $\mathcal{G}$. Conclude that : $\mathcal{D%
}_{A}=\mathcal{D}$.\newline

\noindent \textbf{B3.} Using the symmetrical roles of $A$ and $B$ in the
definition of $B\in \mathcal{D}_{A}$, show that, for any $B\in \mathcal{D}$, 
$\mathcal{G}\subset \mathcal{D}_{B}$ and then $\mathcal{D}_{B}=\mathcal{D}$.%
\newline

\noindent \textbf{B4.} Show then that $\mathcal{D}$ is stable by
intersection and then is a $\sigma $-algebra.\newline

\noindent \textbf{B5.} Conclude.\newline

\newpage
\noindent \LARGE \textbf{DOC 01-05 : Chapter 01 : Measurable sets. Exercises on $\lambda$ and $\pi$ systems, \textit{with solutions}}. \label{doc01-05}\\

 \Large

\bigskip \bigskip \noindent \textbf{Abstract} Monotone and Dynkin classes are important tools
for theoretical demonstrations in Measure Theory. It is important to be
familiar with them. Their importance reside in exercises 3 and 4.\newline

\noindent You will need these basic and easy facts. A $\sigma $-algebra is
an algebra, is a Dynkin system and is a monotone class. An algebra is a $\pi 
$-system.\newline

\bigskip \bigskip \noindent \textbf{Exercise 1}. \label{exercise01_sol_doc01-05} Let $\mathcal{C}$ be an
algebra and a monotone class. Show that $\mathcal{C}$ is a sigma-algebra.%
\newline

\bigskip \noindent  \textbf{Solution}. Since $\mathcal{C}$ is an algebra, we already have these
two points :\\

\bigskip \noindent (i) $\Omega \in $ $\mathcal{C}$\\

\bigskip \noindent (ii) $\forall A\in \mathcal{C}$, $\forall A^{c}\in \mathcal{C}$.\\

\bigskip \noindent To prove that $\mathcal{C}$ is a sigma-algebra, we only have the check, for
example that it is stable under countable unions, that is :

\bigskip \noindent (iii) $\forall \left( A_{n}\right) _{n\geq 1}\subset \mathcal{C}$, $%
\cup _{n\geq 1}A_{n}\in \mathcal{C}$.\\

\bigskip \noindent Let $\left( A_{n}\right) _{n\geq 1}\subset \mathcal{C}$. We use the property
that says that a union is the unions of the partial finite unions :%
\begin{equation*}
\bigcup\limits_{n\geq 1}A_{n}=\bigcup\limits_{n\geq 1}\left(
\bigcup\limits_{m=1}^{n}A_{m}\right) .
\end{equation*}

\bigskip \noindent Denote%
\begin{equation*}
B_{n}=\bigcup\limits_{m=1}^{n}A_{m}
\end{equation*}

\bigskip \noindent so that%
\begin{equation*}
\bigcup\limits_{n\geq 1}A_{n}=\bigcup\limits_{n\geq 1}B_{n}
\end{equation*}

\bigskip \noindent For each $n\geq 1,$ $B_{n}$ is in $\mathcal{C}$\ since $\mathcal{C}$\ is an
algebra. The sequence $(B_{n})_{n\geq 1}$ is non-decreasing. Since $\mathcal{C%
}$\ is a monotone class, we have 
\begin{equation*}
\bigcup\limits_{n\geq 1}B_{n}\in \mathcal{C}
\end{equation*}%

\bigskip \noindent and next 
\begin{equation*}
\bigcup\limits_{n\geq 1}A_{n}=\bigcup\limits_{n\geq 1}B_{n}\in \mathcal{C}.
\end{equation*}

\bigskip \noindent Conclusion. $\mathcal{C}$ is a sigma-algebra.

\bigskip \bigskip \noindent \textbf{Exercise 2}. \label{exercise02_sol_doc01-05} Let $\mathcal{D}$ be a $\pi 
$-class and a Dynkin system. Show that $\mathcal{D}$\ is a sigma-algebra.\\

\bigskip \noindent \textbf{Solution}. Since $\mathcal{D}$ is a Dynkin system, we already have the first point.

\bigskip \noindent \textbf{(i)} $\Omega \in $ $\mathcal{C}$.\\

\bigskip \noindent Next, let us check that\\

\bigskip \noindent (ii) $\forall A\in \mathcal{D}$, $\forall A^{c}\in \mathcal{D}$.\\

\bigskip \noindent Answer. Let $A\in \mathcal{D}$. Since $\mathcal{D}$ is a Dynkin system, since $\Omega
\in \mathcal{D}$ and \ $A\in \mathcal{D}$ with $A\subset \Omega $, we get $%
\Omega \setminus A=\Omega \cap A^{c}=A^{c}\in \mathcal{D}.$\\

\bigskip \noindent Finally, we have to check the last point :\\

\bigskip \noindent (iii) $\forall \left( A_{n}\right) _{n\geq 1}\subset \mathcal{D}$, $\cup
_{n\geq 1}A_{n}\in \mathcal{D}$.\\

\bigskip \noindent Answer. Now, you already know how to render an arbitrary sets union to a sets
summation. You write%
\begin{equation*}
\bigcup\limits_{n\geq 1}A_{n}=\sum\limits_{n\geq 1}B_{n},
\end{equation*}

\bigskip \noindent where
\begin{equation*}
B_{1}=A_{1}\text{ and for }n\geq 2,\text{ }B_{n}=A_{1}^{c}...A_{n-1}^{c}A_{n}
\end{equation*}

\bigskip \noindent The $B_{n}$ are disjoints and they belong to $\mathcal{D}$ since this latter
is stable under complementation (by Point (ii) above) and under finite
intersection by the assumption that $\mathcal{D}$\ is a $\pi$-system. Then $%
\cup _{n\geq 1}A_{n}=\sum_{n\geq 1}B_{n}$ belongs to $\mathcal{D}$\ since $%
\mathcal{D}$, as a Dynkin system, is stable under countable summations of
its subsets.\\

\noindent Conclusion  :  $\mathcal{D}$\ \ is a sigma-algebra.\\

\bigskip \noindent \textbf{Exercise 3}. \label{exercise03_sol_doc01-05} Let $\mathcal{C}$ be an algebra of
subsets of $\Omega $. Prove that $\sigma (\mathcal{C})=m(\mathcal{C})$.%
\newline

\bigskip \noindent \textbf{Proposed solution (in Lo\`{e}ve)}.\\

\noindent Denote $\mathcal{M}=m(\mathcal{C})$\\

\noindent Put also $\mathcal{M}_{A}=\{B\in \mathcal{M},AB^{c}\in \mathcal{M},A^{c}B\in \mathcal{M},AB\in \mathcal{M}\}$.\\

\noindent \textbf{A1.}  Show that for any $A\in \mathcal{M}$, $\mathcal{M}_{A}$ is a monotone class containing $A$.\\

\noindent \textbf{A2.(i)} Show that for any $A\in \mathcal{C}$, $\mathcal{C} \subset \mathcal{M}_{A}$. Deduce from this, that 
 $\mathcal{M}_{A}=\mathcal{M}$.\\

\noindent \textbf{A2.(ii)} Let $B \in \mathcal{M}$. Use the symmetry of the roles of $A$ ad $B$ in the definition $B \in \mathcal{M}_{A}$ and 
A2.(i) to prove that $A \in \mathcal{M}_{B}$ for any $A \subset \mathcal{C}$, and then $\mathcal{C} \subset \mathcal{M}_{B}$. Conclude that
$\mathcal{M}_{B}=\mathcal{M}$.\newline

\noindent \textbf{A3} Conclude that for any $(A,B) \in \mathcal{M}^2$, 
$$
AB \in \mathcal{M}, \text{ } AB^{c} \in \mathcal{M}, \text{ } A^{c}B \in \mathcal{M}
$$

\bigskip \noindent and prove from this, $\mathcal{M}$ is a $\sigma $-algebra.\newline

\bigskip \noindent \textbf{A4.} Conclude.\newline

\bigskip \noindent \textbf{Solution}.\\

\noindent \textbf{A1.} Let $A\in \mathcal{M}$. Let us show that $\mathcal{M}_{A}$ is monotone.\\

\noindent We begin to check that $\mathcal{M}_{A}$ is not empty since $A\in \mathcal{M}%
_{A}$. We get that $A\in \mathcal{M}_{A}$ by checking : $AA^{c}=\emptyset
\in \mathcal{C\subset M}$, $AA^{c}=\emptyset \in \mathcal{C\subset M}$ and $%
AA=A\mathcal{\in M}$ by assumption.\\

\noindent Next, let us check that $\mathcal{M}_{A}$ contains limits of it non-decreasing sequences $(B_{n})_{n\geq 1}$.\\

\noindent If $(B_{n})_{n\geq1}\subset \mathcal{M}_{A}$ and $(B_{n})_{n\geq 1}$ is non-decreasing, we have
for each $n\geq 1$ 
\begin{equation*}
AB_{n}^{c}\in \mathcal{M},A^{c}B_{n}\in \mathcal{M},AB_{n}\in \mathcal{M}
\end{equation*}

\bigskip \noindent Since $\mathcal{M}$ is a monotone class, we obtain by letting $n\uparrow
+\infty $ that

\begin{equation*}
\mathcal{M}\ni AB_{n}^{c}\searrow AB^{c}\in \text{\ }\mathcal{M},\text{ \ }%
\mathcal{M}\ni A^{c}B_{n}\nearrow A^{c}B\in \text{\ }\mathcal{M},\text{ \ }%
\mathcal{M}\ni AB_{n}\nearrow AB\in \text{\ }\mathcal{M}
\end{equation*}%

\bigskip \noindent  We have that $B=\cup _{n\geq 1}B_{n}\in \mathcal{M}_{A}$\\

\noindent Now let us check that $\mathcal{M}_{A}$ contains limits of it non-increasing
sequences $(B_{n})_{n\geq 1}.$ Now, if $(B_{n})_{n\geq 1}\subset \mathcal{M}%
_{A}$ and $(B_{n})_{n\geq 1}$ is non-increasing, we have for each $n\geq 1$ 

\begin{equation*}
AB_{n}^{c}\in \mathcal{M},A^{c}B_{n}\in \mathcal{M},AB_{n}\in \mathcal{M}
\end{equation*}

\bigskip \noindent Since $\mathcal{M}$ is a monotone class, we obtain by letting $n\uparrow
+\infty $ that

\begin{equation*}
\mathcal{M}\ni AB_{n}^{c}\nearrow AB^{c}\in \text{\ }\mathcal{M},\text{ \ }%
\mathcal{M}\ni A^{c}B_{n}\searrow A^{c}B\in \text{\ }\mathcal{M},\text{ \ }%
\mathcal{M}\ni AB_{n}\searrow AB\in \text{\ }\mathcal{M}.
\end{equation*}

\bigskip \noindent We get that $B=\cap _{n\geq 1}B_{n}\in \mathcal{M}_{A}$.\\

\noindent We conclude that each $\mathcal{M}_{A}$ is a monotone class containing $A$.\\

\bigskip \noindent \textbf{A2. (i)} Let $A\in \mathcal{C}$. Any element $B\in \mathcal{C}$ satisfies the
condition of belonging to $\mathcal{M}_{A}$\ \ : $AB^{c}\in \mathcal{%
C\subset M},A^{c}B\in \mathcal{C\subset M},AB\in \mathcal{C\subset M}$,
because $A$ and $B$ are in the algebra $\mathcal{C}$. So%
\begin{equation*}
\mathcal{C\subset M}_{A}\subset \mathcal{M}.
\end{equation*}

\bigskip \noindent Since $\mathcal{M}_{A}$\ is a monotone class including $\mathcal{C}$, then
it includes the monotone class generated by \ $\mathcal{C}$, that is $m(%
\mathcal{C)=M\subset M}_{A}.$ We arrive at%
\begin{equation*}
\mathcal{M\subset M}_{A}\subset \mathcal{M}.
\end{equation*}

\bigskip \noindent We have 
\begin{equation*}
\forall A\in \mathcal{C}\text{, }\mathcal{M}_{A}=\mathcal{M}
\end{equation*}

\bigskip \noindent \textbf{A2. (ii)}. Let $A\in \mathcal{M}$. Any element $B$ of $\mathcal{C}$\ is such that $%
\mathcal{M}_{B}=$\ $\mathcal{M}$. For element $B$ of $\mathcal{C}$, $A\in 
\mathcal{M}_{B}$, that is : $B^{c}A\in \mathcal{C\subset M},BA^{c}\in 
\mathcal{C\subset M},BA\in \mathcal{C}$. Now remark that the role of A and B
are symmetrical in the definition of $A\in \mathcal{M}_{B}$. You conclude
that, For element $B$ of $\mathcal{C}$, $A\in \mathcal{M}_{B}$, and then $%
B\in \mathcal{M}_{A}$. \ You get that $\mathcal{C\subset M}_{A}.$ From there
you conclude exactly as in A2.(i) that

\begin{equation*}
\forall A\in \mathcal{M}\text{, }\mathcal{M}_{A}=\mathcal{M}.
\end{equation*}

\bigskip \noindent \textbf{A3}. What does imply A2? That $\ $for any $(A,B)\in \mathcal{M}^{2}$, $\mathcal{M}_{A}=\mathcal{M}$ and then $B\in \mathcal{M}_{A}$, which implies that%
\begin{equation}
AB^{c}\in \mathcal{M},\text{ }A^{c}B\in \mathcal{M},\text{ }AB\in \mathcal{M}.
\label{CG}
\end{equation}

\bigskip \noindent We may derive from this that is $\mathcal{M}$\ an algebra by checking the three points.\\

\bigskip \noindent (i) $\Omega \in \mathcal{M}$. This comes form $\Omega \in \mathcal{C\subset} \mathcal{M}$ by definition of $\mathcal{M}=m(\mathcal{C}).$\\

\bigskip \noindent (ii) Let $A\in \mathcal{M}.$ Take $B=\Omega ,$ (\ref{CG}) implies $A^{c}B=A^{c}\in \mathcal{M}$.\\

\bigskip \noindent (iii) Let $A\in \mathcal{M}$ and $B\in \mathcal{M}.$  (\ref{CG}) implies $
AB\in \mathcal{M}$.\\

\bigskip \noindent At this point, we got that $\mathcal{M}$ is an algebra and is monotone. By Exercise 1, $%
\mathcal{M}\ $is a sigma-algebra.\\

\bigskip \noindent \textbf{A4.} Since \ $\mathcal{M}=m(\mathcal{C)}$ is an algebra and is also a
monotone class, it is a sigma-algebra by Exercise 1. Besides, it includes $%
\mathcal{C}$. \ Then it includes $\sigma (\mathcal{M})$. And we arrive at 
\begin{equation*}
\sigma (\mathcal{C})\subset m(C\mathcal{)}\text{.}
\end{equation*}

\bigskip \noindent But $\sigma (\mathcal{C)}$, as a sigma-algebra, is also a monotone class
including $\mathcal{C}$. Then it includes $m(\mathcal{C)}$, that is%
\begin{equation*}
m(C\mathcal{)\subset }\text{ }\sigma (\mathcal{C}).
\end{equation*}%
Putting this together gives the final conclusion%
\begin{equation*}
\sigma (\mathcal{C})=\sigma (\mathcal{C})
\end{equation*}%

\bigskip \noindent for any algebra $\mathcal{C}$.\ 

\bigskip 

\bigskip \bigskip \noindent \textbf{Exercise 4}.  \label{exercise04_sol_doc01-05} Let $G$ be a $\pi $-system
of subsets of $\Omega $. We have : $\sigma (\mathcal{G})=(\mathcal{G}$).
\newline

\noindent Proposed \textbf{Solution}.\newline

\noindent \textbf{Notation : } : Denote $\mathcal{D}$ as the Dynkin system
generated by $\mathcal{G}$. Put $\mathcal{D}_{A}=\{B\in \mathcal{D},AB\in 
\mathcal{D}\}$.\newline

\bigskip \noindent \textbf{B1.} Show that $\mathcal{D}_{A}$ is a Dynkin system for any $A\in \mathcal{D}.$\\

\noindent \textbf{B2.} Prove that for any $A\in \mathcal{G}$, $\mathcal{D}%
_{A}$ is a Dynkin system including $\mathcal{G}$. Conclude that : $\mathcal{D%
}_{A}=\mathcal{D}$.\newline

\noindent \textbf{B3.} Using the symmetrical roles of $A$ and $B$ in the
definition of $B\in \mathcal{D}_{A}$, show that, for any $B\in \mathcal{D}$, 
$\mathcal{G}\subset \mathcal{D}_{B}$ and then $\mathcal{D}_{B}=\mathcal{D}$.%
\newline

\noindent \textbf{B4.} Show then that $\mathcal{D}$ is stable by
intersection and then is a $\sigma $-algebra.\newline

\noindent \textbf{B5.} Reminding that a $\sigma $-algebra is a Dynkin
system, show that $\sigma (\mathcal{G})=\mathcal{D}$.\newline

\noindent \textbf{B6.} Show that if $\mathcal{D}$ is a Dynkin system, then $(%
\mathcal{G}\subset \mathcal{D})\Rightarrow (\sigma (\mathcal{G})\subset 
\mathcal{D})$.\newline

\bigskip \noindent \textbf{Solution}.\\

\noindent \textbf{B1.}  Let us check that is a $\mathcal{D}_{A}$ is a Dynkyn system.\\ 

\noindent \textbf{B1.a} Check that $\Omega \in \mathcal{D}_{A}$.\\

\bigskip \noindent \textit{Answer}. Since $A\in \mathcal{D}_{A}$ , we have $\Omega \cap A=A\in \mathcal{D}$ and
then $\Omega \in \mathcal{D}_{A}$.\\

\noindent \textbf{B1.b} Let $(B,C)\in \mathcal{D}_{A}^{2}$ \ with $C\subset B$.
Check that $B\setminus C\in \mathcal{D}_{A}$.\\

\noindent \textit{Answer}. We have $A\cap
(B\setminus C)=(A\cap B)\setminus (A\cap C).$ But , $(A\cap B)\in \mathcal{%
D}$ since $B\in \mathcal{D}_{A}$, and $(A\cap C)\in \mathcal{D}$ since $C\in 
\mathcal{D}_{A}$, and $(A\cap C)\subset (A\cap B),$ and $\mathcal{D}$ is a
Dynkin system. Then 

\begin{equation*}
A\cap (B\setminus C)=(A\cap B)\setminus (A\cap C)\in \mathcal{D},\text{ }
\end{equation*}

\bigskip \noindent and then 
\begin{equation*}
B\setminus C\in \mathcal{D}_{A}.
\end{equation*}

\bigskip \noindent \textbf{B1.c} Let $(B)_{\geq 1}$ a pairwise disjoint sequence of elements of $\mathcal{%
D}_{A}.$ Check that $\sum_{n\geq 1}B_{n}\in \mathcal{D}_{A}$. 

\bigskip \noindent \textit{Answer}. We have for each $n\geq 1,$ $AB_{n}\in \mathcal{D}$. And
\begin{equation*}
A\bigcap \left(\sum_{n\geq 1}B_{n}\right) =\sum_{n\geq 1}AB_{n}.
\end{equation*}

\bigskip \noindent Since all the $AB_{n}$ are in $\mathcal{D}$ \ and $\mathcal{D}$ is a Dynkin
system, we get

\begin{equation*}
A\bigcap \left( \sum_{n\geq 1}B_{n}\right) =\sum_{n\geq 1}AB_{n}\in \mathcal{D}
\end{equation*}

\bigskip \noindent and thus
\begin{equation*}
\left( \sum_{n\geq 1}B_{n}\right) \in \mathcal{D}_{A}.
\end{equation*}

\bigskip \noindent \textbf{B2. (i)} Let $A\in \mathcal{G}$. For any $B\in \mathcal{G}$, we also have $%
AB\in \mathcal{G\subset D}$, since $\mathcal{G}$ is a $\pi $-system. \ Then
For any $B\in \mathcal{G}$, $B\in $ $\mathcal{D}_{A}.$ Thus : \ $\mathcal{%
G\subset D}_{A}$ whenever $A\in \mathcal{G}.$

\bigskip \noindent \textbf{B2.(ii)} Take B2.(i) and B1 into account and conclude like that : for any $%
A\in \mathcal{G}$, $\mathcal{D}_{A}$ - which is a part of $\mathcal{D}$ by
definition - is a Dyskyn system including $\mathcal{G}$. Hence it necessarily
includes the Dynkin system generated by $\mathcal{G}$, that is : \ $\mathcal{%
D}=d(\mathcal{G)\subset D}_{A}\subset \mathcal{D}.$ Then%
\begin{equation*}
\mathcal{D}_{A}=\mathcal{D},
\end{equation*}

\bigskip \noindent for all $A\in \mathcal{G}$.\\

\bigskip \noindent \textbf{B3.(i)} By the symmetry of the roles of $A$ and $B$ in the definition  of $%
B\in \mathcal{D}_{A}$, we may say this : Let $B\in \mathcal{D}.$ For any $%
A\in \mathcal{G}$, $\mathcal{D}_{A}=\mathcal{D}$ and then $B\in \mathcal{D}=%
\mathcal{D}_{A}$, and then $AB\in \mathcal{D},$ and then $A\in \mathcal{D}%
_{B}.$ So, $\mathcal{G}\subset \mathcal{D}_{B}.$ 

\bigskip \noindent \textbf{B3.(ii)} It follows from this and from B1 that : $\mathcal{D}_{B}$ - which is
a part of $\mathcal{D}$ by definition - is a Dyskyn system including $%
\mathcal{G}$. Then it necessarily includes the Dynkin system generated by $%
\mathcal{G}$, that is : \ $\mathcal{D}=d(\mathcal{G)\subset D}_{B}\subset 
\mathcal{D}.$ Then%
\begin{equation*}
\mathcal{D}_{B}=\mathcal{D}.
\end{equation*}

\bigskip \noindent \textbf{B4}. Let us show that $\mathcal{D}$ is a $\pi $-system?\\

\bigskip \noindent \textit{Answer}. Let $(A,B)\in \mathcal{D}%
^{2}.$ But Let $B\in \mathcal{D}=\mathcal{D}_{A}$ and then, by definition, $%
AB\in \mathcal{D}$. So $\mathcal{D}$ is a $\pi $-system. By Exercise 2, $%
\mathcal{D}$ a sigma-algebra since it is both a $\pi $-system and a Dynkin
system.\\

\bigskip \noindent \textbf{B5}. Conclusion. On one hand, $\mathcal{D}=d(\mathcal{G})$ is a sigma-algebra
including $\mathcal{G}$, then $\sigma (\mathcal{G})\subset \mathcal{D}=d(%
\mathcal{G})$. \ On the other hand, $\sigma (\mathcal{G})$ is a Dynkyn
system including $\mathcal{G}$, then $d(\mathcal{G})=\mathcal{D}\subset
\sigma (\mathcal{G})$ . Then%
\begin{equation*}
d(\mathcal{G})=\sigma (\mathcal{G}).
\end{equation*}

\newpage
\noindent \LARGE \textbf{DOC 01-06 : More Exercises}. \label{doc01-06}\\

 \Large

\bigskip \noindent \textbf{Exercise 1}. \label{exercise01_doc01-06} Le $\Omega$ be a nonempty set and set

\begin{equation*}
\mathcal{A}=\{A\subset \Omega ,\ \ A\ \ is\ \ finite\ \ \text{or\newline
}A^{c}\ \ is\ finite\}.
\end{equation*}

\bigskip \noindent (a) Show that $\mathcal{A}$ is an algebra.\newline

\noindent (b) Show that $\mathcal{A}$ is a $\sigma$-algebra if $%
\Omega$ is finite.\newline

\noindent (c) If $\Omega $ is infinite, consider an infinite subset 
$C=\{{\omega _{0},\omega _{1},\omega _{2},..\}}$ and set $A_{n}=\{{\omega
_{2n}\}}$, $n\geq 0$. Show that the $A_{n}$ are in $\mathcal{A}$ and their
union is not.\newline

\noindent (d) Conclude that $\mathcal{A}$ is a $\sigma $-algebra 
\textbf{if and only if} $\Omega $ is finite.\newline

\bigskip \noindent \textbf{Exercise 2}. \label{exercise02_doc01-06} \textbf{1}.\newline
Let $\Omega =\{\omega _{1},...,\omega _{2n}\}$ , $n\geq 2$, be a space of
even cardinality. Let $\mathcal{A}$ the collection of all subsets of $%
\Omega $ having an even cardinality.\newline

\bigskip \noindent (a) Show that $\mathcal{A}$ is a Dynkin system.\newline

 \noindent (b) Use $A=\{\omega _{1},\omega _{2}\}$ and $%
B=\{\omega_{1},\omega _{3}\}.$ to show that $\mathcal{A}$ is not stable by
finite intersection.?\\

\noindent (c) Make a conclusion.\newline

\bigskip \noindent \textbf{Exercise 3.} \label{exercise03_doc01-06} \textbf{Product spaces}. Consider
two measurable spaces $( \Omega_{i}, \mathcal{A}_{i} )$, $i=1,2$ and 
\begin{equation*}
\mathcal{S}=\{A \times B, A \in \mathcal{A}_{1}, B \in \mathcal{A}_{2} \}
\end{equation*}

\noindent S Let $A\in \sigma (\mathcal{S})$. Let $\omega _{1}$ a fixed element
of $\Omega _{1}$. Define the section $A$ at $\omega _{1}$ by : 
\begin{equation*}
A_{\omega _{1}}=\{\omega _{2}\in \Omega _{2},(\omega _{1},\omega _{2})\in
A\}\subset \Omega _{2}.
\end{equation*}

\noindent S Define likely the section $A$ at $\omega_{2}$ for a fixed $%
\omega_{2} \in \Omega_{2}$ as a subset of $\Omega_{1}$.\newline

\noindent We have to show that the sections of a measurable subset $A$ of $%
\Omega =\Omega _{1}\times \Omega _{2}$ with respect to the product $\sigma $%
-algebra $\mathcal{A}=\mathcal{A}_{1}\otimes \mathcal{A}_{2}$ are measurable.%
\newline

\bigskip \noindent \textbf{(a)} Operations on the sections. Prove these formulas for subsets of 
$\Omega =\Omega _{1}\times \Omega _{2}.$ (We fix $\omega _{1}\in \Omega _{1}$
here. The same results are obtained for a fixed $\omega _{2}\in \Omega _{2}).
$.\\

\noindent (1) $\Omega _{\omega _{1}}=\Omega _{2};$ $\emptyset _{\omega _{1}}=\emptyset$\\

\noindent (2) $(A^{c})_{\omega _{1}}=(A_{\omega _{1}})^{c}$\\

\noindent (3) $(A\backslash B)_{\omega _{1}}=A_{\omega _{1}}\backslash B_{\omega _{1}}$\\

\noindent (4) $\left( \cup _{n\geq 0}A_{n}\right) _{\omega _{1}}=\cup _{n\geq 0}(A_{n})_{\omega _{1}}$\\

\noindent (5)  $\left( \cap _{n\geq 0}A_{n}\right) _{\omega _{1}}=\cap _{n\geq
0}(A_{n})_{\omega _{1}}$\\

\noindent (6)  $\left( \sum_{n\geq 0}A_{n}\right) _{\omega _{1}}=\sum_{n\geq
0}(A_{n})_{\omega _{1}}$\\

\noindent (7) etc.\\

\noindent \textbf{(b)} Find consider sections of elements $\mathcal{S}$. Determine $A_{\omega
_{1}}$ when $A=A_{1}\times A_{2}$ is element of $\mathcal{S}$ with $A_{1}\in 
\mathcal{A}_{1}$ and $A_{2}\in \mathcal{A}_{2}$.\newline

\noindent \textit{\textit{Hint}} : Explain why
\begin{equation*}
A_{\omega _{1}}=\left\{ 
\begin{tabular}{lll}
$A_{2}$ & if & $\omega _{1}\in A_{1}$ \\ 
$\emptyset $ & if & $\omega _{1}\notin A_{1}$%
\end{tabular}%
\right. 
\end{equation*}

\bigskip \noindent Conclude that $A_{\omega _{1}}$ is measurable in $\Omega _{2}$.\newline

\noindent \textbf{(c)} Consider $\mathcal{D}=\{A\in \mathcal{A}_{1}\otimes 
\mathcal{A}_{2},A_{\omega _{1}}\text{is measurable}\}$. Show that $\mathcal{D%
}$ is a $\sigma $-algebra containing $\mathcal{S}$. \newline

\noindent \textbf{(d)} Conclude.\newline

\newpage
\noindent \LARGE \textbf{DOC 01-07 : More  Exercises - with solutions}. \label{doc01-07}\\

\Large

\bigskip \noindent \textbf{Exercise 1}. \label{exercise01_sol_doc01-07} Le $\Omega$ be a nonempty set and set

\begin{equation*}
\mathcal{A}=\{A\subset \Omega ,\ \ A\ \ is\ \ finite\ \ \text{or\newline
}A^{c}\ \ is\ finite\}.
\end{equation*}

\bigskip \noindent (a) Show that $\mathcal{A}$ is an algebra.\newline

\noindent (b) Show that $\mathcal{A}$ is a $\sigma$-algebra if $\Omega$ is finite.\newline

\noindent (c) If $\Omega $ is infinite, consider an infinite subset 
$C=\{{\omega _{0},\omega _{1},\omega _{2},..\}}$ and set $A_{n}=\{{\omega
_{2n}\}}$, $n\geq 0$. Show that the $A_{n}$ are in $\mathcal{A}$ and their
union is not.\newline

\noindent (d) Conclude that $\mathcal{A}$ is a $\sigma $-algebra 
\textbf{if and only if} $\Omega $ is finite.\newline

\bigskip \noindent \textbf{Solution}.\newline

\bigskip \noindent (a) Let us check that $\mathcal{A}$ is an algebra.\newline

\bigskip \noindent (i) $\Omega $ belongs to $\mathcal{A}$ since its
complement is empty and then finite.\newline

\bigskip \noindent (ii) The roles of $A$ and $A^{c}$ are symmetrical in the
definition of $\mathcal{A}.$ So $A$ belongs to $\mathcal{A}$ whenever $A^{c}$
belongs to $\mathcal{A}$ and vice-verse.\newline

\bigskip \noindent (iii) Let us show that $\mathcal{A}$ is stable under
finite unions. Let $A$ and $B$ be two elements of $\mathcal{A}.$ we have two
case. Either $A$ and $B$ are finite sets and then $A\cup B$ is finite and
hence $A\cup B$ belongs to $\mathcal{A}$. Or one of them is not finite.
Suppose that $A$ is not finite. Since $A$ is in $\mathcal{A}$ and is not
finite, this implies that $A^{c}$ is finite (otherwise A would not be in $%
\mathcal{A}$). Then%
\begin{equation*}
(A\cup B)^{c}=A^{c}\cap B^{c}\subset A^{c}.
\end{equation*}

\bigskip \noindent Thus $(A\cup B)^{c}$ is finite as part of a finite set.
Hence $A\cup B\in \mathcal{A}$. Then in all possibles cases, $A\cup B$
belongs to $\mathcal{A}$ whenever $A$ and $B$ are in $\mathcal{A}.$\newline

\bigskip \noindent (b) Since $\mathcal{A}$ is already an algebra, if
suffices to show that is is stable under countable unions or intersections.
But if $\Omega $ is a finite set, the collection of subsets is also finite
and so, a countable union of elements of $\mathcal{A}$ is merely a finite
union of element of $\mathcal{A}$ and then remains in $\mathcal{A}.$ Thus,
it is clear that $\mathcal{A}$ is an sigma-algebra if $\Omega $ is finite.\ 

\bigskip \noindent (c) Clearly the sets $A_{n}$ are finite and then, they
belong to $\mathcal{A} $. But their union

\begin{equation*}
A=\cup _{n\geq 0}A_{n}=\{\omega _{0},\omega _{2},...\}
\end{equation*}

\bigskip \noindent is neither finite nor co-finite since%
\begin{equation*}
\{\omega _{1},\omega _{3},...\}\subset A^{c}.
\end{equation*}

\bigskip \noindent (d) The conclusion is : $\mathcal{A}$ is a sigma-algebra
if and only if $\Omega $ is finite.\\

\bigskip 
\bigskip \noindent \textbf{Exercise 2}. \label{exercise02_sol_doc01-07} \newline
Let $\Omega =\{\omega _{1},...,\omega _{2n}\}$ , where $n$ is finite fixed
integer such that $n\geq 2$, be a space of even cardinality. Let $\mathcal{A}
$ the collection of all subsets of $\Omega $ having an even cardinality.%
\newline

\bigskip \noindent (a) Show that $\mathcal{A}$ is a Dynkin system.\newline

\bigskip \noindent (b) Use $A=\{\omega _{1},\omega _{2}\}$and $B=\{\omega
_{1},\omega _{3}\}.$ to show that $\mathcal{A}$ is not stable by finite
intersection?\\

\bigskip \noindent (c) Make a conclusion.\newline

\bigskip \noindent \textbf{Solution}.\\

\bigskip \noindent It is important here to remark that $n$ is finite. So any union of subsets
of $\Omega $ is a finite \textbf{union}. Now let us show that $\mathcal{A}$\
is a Dynkin system.\\

\bigskip \noindent (i) $\Omega \in \mathcal{A}$ ? We have to check that $\Omega $ has an even
cardinality. This is true since the cardinality of $\Omega $ is\ $2n$, and
then even.\\

\bigskip \noindent (ii) Let $A$ and $B$ in $\mathcal{A}$\ \ with \ $A\subset B.$ We have to
check that $B\backslash A$ is in $\mathcal{A}$. But $A\in \mathcal{A}$ and $%
B\in \mathcal{A}$ mean that their cardinals $\#(A)$ and $\#(B)$ are both
even. Also, since  $A\subset B$ holds, we have%
\begin{equation*}
\#(B\backslash A)=\#(B)-\#(A)
\end{equation*}

\bigskip \noindent Since $\#(B\backslash A)$ is difference of even integers, it is also even.
Thus \ $(B\backslash A)\in \mathcal{A}.$

\bigskip \noindent (iii) Let $(A_{n})_{n\geq 0}$ a pairwise disjoint elements of \ $%
\mathcal{A}$. We have to prove that $A=\sum_{n}A_{n}\in \mathcal{A}.$ But as
said earlier, since \ $\Omega $ is finite, any union reduces to a finite
union. So $\sum_{n}A_{n}$ reduce to a finite sum $\sum_{1\leq k\leq p}B_{k},$
where $B_{k}\in \mathcal{A}$ and $p$ is finite. Then

\begin{equation*}
\#(A)=\sum_{1\leq k\leq p}\#(B_{k})
\end{equation*}

\bigskip \noindent and $\#(A)$ is even as a finite sum of even integers. Thus $A=\sum_{n}A_{n}\in \mathcal{A}$.\\

\noindent We may conclude that $\mathcal{A}\ $is a Dynkin system.\\

\noindent (b) : Now, since $\Omega $ has at least four elements, we see that by the
notation, \ $\{\omega _{1},\omega _{2},\omega _{3}\}\subset \Omega .$
Clearly $A=\{\omega _{1},\omega _{2}\}$and $B=\{\omega _{1},\omega _{3}\}$
have even cardinality and then, are in . But their intersection $A\cap $ $%
B=\{\omega _{1}\}$ has an odd cardinality (1) and then does not lie in $%
\mathcal{A}$.\\

\noindent (c) : Conclusion : A Dyskin system is not necessarily stable by finite
intersection, and then, it is not necessarily an algebra.\\

\noindent Remark : You will see in the homework 01 that a Dynkyn system which is also a $\pi $-system, is an algebra and then a sigma-algebra.\\

\bigskip 
\bigskip \noindent \textbf{Exercise 3.}  \label{exercise03_sol_doc01-07}  \textbf{Product spaces}. Consider
two measurable spaces $( \Omega_{i}, \mathcal{A}_{i} )$, $i=1,2$ and 
\begin{equation*}
\mathcal{S}=\{A \times B, A \in \mathcal{A}_{1}, B \in \mathcal{A}_{2} \}.
\end{equation*}

\bigskip \noindent Let $A\in \sigma (\mathcal{S})$. Let $\omega _{1}$ a fixed element
of $\Omega _{1}$. Define the section $A$ at $\omega _{1}$ by : 
\begin{equation*}
A_{\omega _{1}}=\{\omega _{2}\in \Omega _{2},(\omega _{1},\omega _{2})\in
A\}\subset \Omega _{2}.
\end{equation*}

\bigskip \noindent Define likely the section $A$ at $\omega_{2}$ for a fixed $%
\omega_{2} \in \Omega_{2}$ as a subset of $\Omega_{1}$.\newline

\noindent We have to show that the sections of a measurable subset $A$ of $%
\Omega =\Omega _{1}\times \Omega _{2}$ with respect to the product $\sigma $%
-algebra $\mathcal{A}=\mathcal{A}_{1}\otimes \mathcal{A}_{2}$ are measurable.%
\newline

\bigskip \noindent \textbf{(a)} Operations on the sections. Prove these formulas for subsets of 
$\Omega =\Omega _{1}\times \Omega _{2}.$ (We fix $\omega _{1}\in \Omega _{1}$
here. The same results are obtained for a fixed $\omega _{2}\in \Omega _{2}).
$.\\

\bigskip \noindent (1) $\Omega _{\omega _{1}}=\Omega _{2};$ $\emptyset _{\omega _{1}}=\emptyset$\\

\noindent (2) $(A^{c})_{\omega _{1}}=(A_{\omega _{1}})^{c}$\\

\noindent (3) $(A\backslash B)_{\omega _{1}}=A_{\omega _{1}}\backslash B_{\omega _{1}}$\\

\noindent (4) $\left( \cup _{n\geq 0}A_{n}\right) _{\omega _{1}}=\cup _{n\geq 0}(A_{n})_{\omega _{1}}$\\

\noindent (5)  $\left( \cap _{n\geq 0}A_{n}\right) _{\omega _{1}}=\cap _{n\geq 0}(A_{n})_{\omega _{1}}$\\

\noindent (6)  $\left( \sum_{n\geq 0}A_{n}\right) _{\omega _{1}}=\sum_{n\geq 0}(A_{n})_{\omega _{1}}$\\

\noindent (7) etc.\\

\bigskip \noindent \textbf{(b)} Find consider sections of elements $\mathcal{S}$. Determine $A_{\omega
_{1}}$ when $A=A_{1}\times A_{2}$ is element of $\mathcal{S}$ with $A_{1}\in 
\mathcal{A}_{1}$ and $A_{2}\in \mathcal{A}_{2}$.\newline

\noindent \textit{Hint} : Explain why
\begin{equation*}
A_{\omega _{1}}=\left\{ 
\begin{tabular}{lll}
$A_{2}$ & if & $\omega _{1}\in A_{1}$ \\ 
$\emptyset $ & if & $\omega _{1}\notin A_{1}$.
\end{tabular}
\right. 
\end{equation*}

\bigskip \noindent Conclude that $A_{\omega _{1}}$ is measurable in $\Omega _{2}$.\newline

\noindent \textbf{(c)} Consider $\mathcal{D}=\{A\in \mathcal{A}_{1}\otimes 
\mathcal{A}_{2},A_{\omega _{1}}\text{is measurable}\}$. Show that $\mathcal{D%
}$ is a $\sigma $-algebra containing $\mathcal{S}$. \newline

\noindent \textbf{(d)} Conclude.\newline

\bigskip \noindent \textbf{Solution}.\\

\noindent \textbf{(a)} We only proof some of these properties. We only have to apply the definition.\\

\noindent (1) $\Omega _{\omega _{1}}=\{\omega _{2}\in \Omega _{2},(\omega _{1},\omega
_{2})\in \Omega \}=\{\omega _{2}\in \Omega _{2},\omega _{1}\in \Omega _{1}$ $%
and$ $\omega _{2}\in \Omega _{2}\}=\{\omega _{2}\in \Omega _{2}\}=\Omega
_{2}.$\\

\noindent $\emptyset _{\omega _{1}}=\{\omega _{2}\in \Omega _{2},(\omega _{1},\omega
_{2})\in \emptyset \}.$ We can not find \ $\omega _{2}\in \Omega
_{2},(\omega _{1},\omega _{2})\in \emptyset .$ So $\emptyset _{\omega
_{1}}=\emptyset$.\\

\noindent (2) $(A_{\omega _{1}})^{c}=\{\omega _{2}\in \Omega _{2},(\omega _{1},\omega
_{2})\notin A\}=\{\omega _{2}\in \Omega _{2},(\omega _{1},\omega _{2})\in
A^{c}\}=(A^{c})_{\omega _{1}}$.\\

\noindent (3) We have

\begin{eqnarray*}
A_{\omega _{1}}\backslash B_{\omega _{1}} &=&\{\omega _{2}\in \Omega
_{2},(\omega _{1},\omega _{2})\in A\text{ and }(\omega _{1},\omega
_{2})\notin B\} \\
&=&\{\omega _{2}\in \Omega _{2},(\omega _{1},\omega _{2})\in A\text{ and }%
(\omega _{1},\omega _{2})\in B^{c}\}= \\
&=&\{\omega _{2}\in \Omega _{2},(\omega _{1},\omega _{2})\in A\}\cap
\{\omega _{2}\in \Omega _{2},(\omega _{1},\omega _{2})\in B^{c}\} \\
&=&A_{\omega _{1}}\cap (B_{\omega _{1}})^{c}\text{ } \\
&=&A_{\omega _{1}}\cap B_{\omega _{1}}^{c}\text{ (by Point 2)} \\
&=&A_{\omega _{1}}\backprime B_{\omega _{1}}^{c}.
\end{eqnarray*}

\bigskip \noindent (4) We have
\begin{eqnarray*}
\left( \cup _{n\geq 0}A_{n}\right) _{\omega _{1}} &=&\{\omega _{2}\in \Omega
_{2},(\omega _{1},\omega _{2})\in \cup _{n\geq 0}A_{n}\} \\
&=&\{\omega _{2}\in \Omega _{2},\exists n\geq 0,(\omega _{1},\omega _{2})\in
A_{n}\} \\
&=&\cup _{n\geq 0}\{\omega _{2}\in \Omega _{2},(\omega _{1},\omega _{2})\in
A_{n}\}\text{ (by definition of union)} \\
&=&\cup _{n\geq 0}(A_{n})_{\omega _{1}}.
\end{eqnarray*}

\bigskip \noindent (5) Combine (2) and (4) to get it.\\

\noindent (6) : To get (6) from (4), you only need to see that if $A$ and $B$ are two
disjoint subsets of $\Omega $, then  for any $\omega _{1}\in \Omega _{1},$
one has $A_{\omega _{1}}\cap B_{\omega _{1}}=\emptyset .$ Indeed if $A$ and $%
B$ are disjoint, find some $\omega _{2}\in A_{\omega _{1}}\cap B_{\omega
_{1}}$ leads to $(\omega _{1},\omega _{2})\in A\cap B,$ which is impossible.\\

\noindent \textbf{(b)} Let  $A=A_{1}\times A_{2}$ be an  element of $\mathcal{S}$ with $A_{1}\in 
\mathcal{A}_{1}$ and $A_{2}\in \mathcal{A}_{2}$. We have 
\begin{equation*}
A_{\omega _{1}}=\{\omega _{2}\in \Omega _{2},\omega _{1}\in A_{1}\text{ }and%
\text{ }\omega _{2}\in A_{2}\}.
\end{equation*}

\bigskip \noindent Discuss on whether  $\omega _{1}\in A_{1}$ or not. If $\omega _{1}\notin
A_{1},$ you can not find a $\omega _{2}\in \Omega _{2}$ such that $\omega
_{1}\in A_{1}$ $and$ $\omega _{2}\in A_{2}$. If $\omega _{1}\notin A_{1},$
then all the elements $\omega _{2}$ of $A_{2}$ satisfy the condition : $%
\omega _{2}\in \Omega _{2},\omega _{1}\in A_{1}$ $and$ $\omega _{2}\in A_{2}$
and only them satisfy it. So $A_{\omega _{1}}=A_{2}.$\\

\noindent In each case, the section ($\emptyset$ or $A_2$) is measurable.\\

\noindent \textbf{(c)}. You may easily prove that $\mathcal{D}$\ is a sigma-algebra with the
properties in (a)\\

\noindent \textbf{(d)} Conclusion : $\mathcal{D\subset A}_{1}\otimes \mathcal{A}_{2}$ by its
definition. Now $\mathcal{D}$\ is a sigma-algebra including $\mathcal{S}$\
by Question (b). Then \ includes the sigma-algebra generated by $\mathcal{S}$%
, that is $\mathcal{A}_{1}\otimes \mathcal{A}_{2}\subset \mathcal{D}$. We
arrive at $\mathcal{A}_{1}\otimes \mathcal{A}_{2}=\mathcal{D}$. conclude
that for any measurable subset of $\Omega ,$ its sections are measurable.

%\chapter{Introduction to Measurable Applications}
\chapter{Measurable applications} \label{02_applimess}

\bigskip
\noindent \textbf{Content of the Chapter}\\

\begin{table}[htbp]
	\centering
		\begin{tabular}{llll}
		\hline
		Type& Name & Title  & page\\
		\hline
		S & Doc 02-01 & Measurable Applications - A summary& \pageref{doc02-01}\\
		S & Doc 02-02 & What you cannot ignore on Measurability  & \pageref{doc02-02}\\
		S& Doc 02-03 & Exercises on Measurable Applications   & \pageref{doc02-03}\\
		D& Doc 03-04 & Exercises on Measurable Applications with solutions& \pageref{doc02-04}\\
		\hline
		\end{tabular}
\end{table}

\newpage

\noindent \LARGE \textbf{DOC 02-01 - Measurable Applications - A summary}. \label{doc02-01}\\
\Large

\noindent \textbf{Abstract}. Measurable applications are applications preserving
measurability of sets by means of inverse images. Before you continue, you
are asked to go back to \textit{Exercise 13 in Doc 01-02} (page \pageref{exercise13_doc01-02}) and to its solutions
in \textit{Doc 01.03 in Chapter \ref{01_setsmes}} (page \pageref{exercise13_sol_doc01-03}), to revise the properties of inverse
images, in particular to get acquainted with the fact that the inverse image
preserves all sets operations.\newline

\bigskip \noindent \textbf{(02.01)} \textbf{Definition of measurable applications}. \\

\noindent An application $X$ from the measurable space $(\Omega _{1},\mathcal{A}_{1})$
to the measurable space $(\Omega _{2},\mathcal{A}_{2})$, is measurable,
with respect to the two $\sigma $-algebras $\mathcal{A}_{1}$ and $\mathcal{A}%
_{2}$, and we say that $X$ is ($\mathcal{A}_{1}$-$\mathcal{A}_{2}$)-measurable, if and only if 

\begin{equation}
\forall (B\in \mathcal{A}_{2}),\text{ }X^{-1}(B)\in \mathcal{A}_{1},  \tag{M}
\end{equation}

\noindent meaning : the reciprocal image of a measurable subset of \ $\Omega _{2}$\ is
a measurable subset of $\Omega _{1}$.\\

\noindent \textbf{NB}. In measure theory, It is usual to use this notation : $%
X^{-1}(B)=(X\in B)$. Never forget the definition of this set : $(X\in
B)=\{\omega \in \Omega _{1},X(\omega )\in B\} \subset \Omega_{1}$\newline

\bigskip \noindent \textbf{Condition (M) is equivalent to} : 
\begin{equation}
\mathcal{A}_{X}=\{(X\in B),B\in \mathcal{A}_{2}\}\subset \mathcal{A}_{1}. 
\tag{MS}
\end{equation}

\bigskip \noindent \textbf{(02.02)} \textbf{Weakening the definition}.\\

\noindent Fortunately, in many cases, $\mathcal{A}_{2}$ is generated by a small class,
say $\mathcal{C}$, that is $\mathcal{A}_{2}=\sigma (\mathcal{C})$. In these
cases, the application $X$ is measurable if and only if  
\begin{equation}
\forall (B\in \mathcal{C}),\text{ }X^{-1}(B)\in \mathcal{A}_{1}.  \tag{WMC}
\end{equation}

\bigskip \noindent \textbf{Catch this} : the condition (WMC) is called a weak measurability
condition. It is very useful in many situations. Immediate examples are
provided in Points (02.04) and (02.05) below.\\

\bigskip 
\noindent \textbf{(02.03)}. \textbf{Measurability of composition of applications}.\\

\noindent Let $X$ be a measurable application from $(\Omega _{1},\mathcal{A}_{1})$ to 
$(\Omega _{2},\mathcal{A}_{2})$, and let $Y$ be a measurable application
form $(\Omega _{2},\mathcal{A}_{2})$ to $(\Omega _{3},\mathcal{A}_{3})$.\\

\noindent Then the composition of applications (function of function) $Z=Y(X)=YX$ defined from $(\Omega
_{1},\mathcal{A}_{1})$ to $(\Omega _{3},\mathcal{A}_{3})$ by 
\begin{equation*}
\Omega _{1}\ni \omega \hookrightarrow Z(\omega )=Y(X(\omega ))
\end{equation*}

\bigskip \noindent is measurable.\newline

\noindent (Reminder) : In terms of inverse images, we have 
$$
Y(X)^{-1}(C)=X^{-1}(Y^{-1}(C)).
$$

\bigskip \noindent \textbf{(02.04)} \textbf{Simple application to continuous functions}.\newline

\noindent Let $X$ be an application from the topological $(E_{1},\mathcal{T}_{1})$ to the topological space $(E_{2},\mathcal{T}_{2})$, where $\mathcal{T}_{1}$ and $\mathcal{T}_{2}$ are the classes of open sets of the two topological spaces. By definition, $X$ is
continuous if and only if 
\begin{equation}
\forall (B\in \mathcal{T}_2),\text{ } X^{-1}(B)\in \mathcal{T}_{1},  \tag{C}
\end{equation}

\bigskip \noindent .\\

\noindent Now let $\mathcal{B}(E_1)=\sigma(\mathcal{T}_{1})$ and $\mathcal{B}(E_2)=\sigma(\mathcal{T}_{2})$ be the corresponding Borel $\sigma$-algebras on $E_1$ and $E_2$ respectively. Then
(C) implies

\begin{equation}
\forall (B\in \mathcal{T}_2),\text{ } X^{-1}(B)\in \mathcal{B}(E_1),
\end{equation}

\bigskip \noindent since $\mathcal{T}_{1} \subset \mathcal{B}(E_1)$. By
(WMC), $X$ is measurable with respect to the Borel $\sigma$-algebra.\newline

\noindent \textbf{Conclusion} : A continuous application is measurable with respect
to the Borel fields. As a consequence, the class of measurable functions is
broader than the class of continuous functions.\newline

\bigskip 

\bigskip \noindent \textbf{(02.05)} \textbf{Product Space. Measurability with respect the product $\sigma $-algebra}.\\

\noindent Let be given a product space \noindent\ endowed with his product $\sigma $-algebra \ 
\begin{equation*}
(\prod_{1\leq i\leq k}\Omega _{i},\bigotimes_{1\leq i\leq k}\mathcal{A}_{i}).
\end{equation*}

\bigskip \noindent \textbf{Recall} : \ $\bigotimes_{1\leq i\leq k}\mathcal{A}_{i}$ is generated by the
semi-algebra of the class of cylinders
\begin{equation*}
\mathcal{S}=\{A_{1}\times A_{2}\times ...\times A_{k},\text{ }\forall (1\leq
j\leq k),\text{ }A_{j}\in \mathcal{A}_{j}\}.
\end{equation*}

\bigskip \noindent \textbf{\textbf{(02.06)} (I) Projections}.\\

\bigskip \noindent \textbf{(02.06a) Definitions}. For $1\leq i\leq k,$ the $jth$ projection is defined by%
\begin{equation*}
\begin{tabular}{lllll}
$\pi _{i}$ & $:$ & $(\prod_{1\leq j\leq k}\Omega _{j},\bigotimes_{1\leq
j\leq k}\mathcal{A}_{j})$ & $\longmapsto $ & $\Omega _{i}$ \\ 
&  & $(x_{1},...,x_{i-1},\mathbf{x}_{i},x_{i+1},...,x_{k})$ & $%
\hookrightarrow $ & $x_{i}$.
\end{tabular}
\end{equation*}

\bigskip \noindent  \textbf{Exemple : first projection}.
\begin{equation*}
\begin{tabular}{lllll}
$\pi _{1}$ & $:$ & $(\prod_{1\leq j\leq k}\Omega _{j},\bigotimes_{1\leq
j\leq k}\mathcal{A}_{j})$ & $\longmapsto $ & $\Omega _{1}$ \\ 
&  & $(x_{1},...,x_{k}$) & $\hookrightarrow $ & $x_{1}$.
\end{tabular}
\end{equation*}

\bigskip
\noindent \textbf{(02.06b) Special case : $k=2$}.
\begin{equation*}
\pi _{1}(x,y)=x\text{ and }\pi _{2}(x,y)=y.
\end{equation*}

\bigskip \noindent \textbf{(02.06c) Measurability of the projections}. For $B\subset \Omega _{i},$ $B\in 
\mathcal{A}_{i}$, we have%
\begin{equation*}
\pi _{i}^{-1}(B)=\Omega _{1}\times ...\times \Omega _{i-1}\times B\times
\Omega _{i+1}\times ...\times \Omega _{k}\in \mathcal{S}.
\end{equation*}

\bigskip \noindent In the case where $k=2$, we have
\begin{equation*}
\pi _{1}^{-1}(A)=A_{1}\times \Omega _{2}\in \mathcal{S}\text{ \ \ and }\pi
_{2}^{-1}(B)=\Omega _{1}\times B\in \mathcal{S}.
\end{equation*}

\bigskip \noindent \textbf{Conclusion} : The projections $\pi _{i},$ $1\leq i\leq k,$ are
measurable with respect to the product $\sigma $-algebra.\\

\noindent \textbf{(02.07) (II) Multi-component applications}. Let $X$ be an application defined on some
measurable space $(E,\mathcal{B)}$ with values in the product space 

$$
(\prod_{1\leq j\leq k}\Omega _{j},\bigotimes_{1\leq j\leq k}\mathcal{A}_{j}).
$$

\bigskip 
\begin{equation*}
\begin{tabular}{lllll}
$X$ & $:$ & $(E,\mathcal{B)}$ & $\longmapsto $ & $(\prod_{1\leq j\leq
k}\Omega _{j},\bigotimes_{1\leq j\leq k}\mathcal{A}_{j})$ \\ 
&  & $\omega $ & $\hookrightarrow $ & $X(\omega )$.
\end{tabular}
\end{equation*}

\bigskip 
\noindent The image $X(\omega )$ has $k$ components in $(\prod_{1\leq i\leq
k}\Omega _{i},\bigotimes_{1\leq i\leq k}\mathcal{A}_{i})$ : \ 
\begin{equation*}
\begin{tabular}{lllll}
$X$ & $:$ & $(E,\mathcal{B)}$ & $\longmapsto $ & $(\prod_{1\leq j\leq
k}\Omega _{j},\bigotimes_{1\leq j\leq k}\mathcal{A}_{j})$ \\ 
&  & $\omega $ & $\hookrightarrow $ & $X(\omega )=(X_{1}(\omega
),X_{2}(\omega ),...,X_{k}(\omega ))$.
\end{tabular}
\end{equation*}

\bigskip 
\noindent Each component  $X_{i}$ is an application from $(E,\mathcal{B)}$
to $(\Omega _{i},\mathcal{A}_{i})$ :

\begin{equation*}
\begin{tabular}{llllll}
$X:$ & $(E,\mathcal{B)}$ & $\longmapsto $ & $(\prod_{1\leq j\leq k}\Omega
_{j},\bigotimes_{1\leq j\leq k}\mathcal{A}_{j})$ & $\longmapsto $ & $(\Omega
_{i},\mathcal{A}_{i})$  \\ 
& $\omega $ & $\hookrightarrow $ & $X(\omega )=(X_{1}(\omega ),X_{2}(\omega
),...,X_{k}(\omega ))$ & $\hookrightarrow $ & $\pi _{i}(X(\omega
))=X_{i}(\omega)$.
\end{tabular}
\end{equation*}

\bigskip 
\noindent \textbf{(02.07a)} Each component is a composition of functions : 
\begin{equation*}
X_{i}=\pi _{i}\circ X,1\leq i\leq k. 
\end{equation*}

\bigskip 

\noindent \textbf{(02.07b)} For each $A_{1}\times A_{2}\times ...\times A_{k}\in \mathcal{S}$,with $A_{j}\in \mathcal{A}_{j}$, $1\leq j\leq k$, we have

\begin{eqnarray}
(X &\in &A_{1}\times A_{2}\times ...\times
A_{k})=((X_{1},X_{2},...,X_{k})\in A_{1}\times A_{2}\times ...\times A_{k}) \notag \\
&=&(X_{1}\in A_{1},X_{2}\in A_{2},...,X_{k}\in A_{k})  \notag \\
&=&\bigcap\limits_{j=1}^{k}X_{i}^{-1}(A_{i}).  \notag
\end{eqnarray}

\noindent that is

\begin{equation}
(X \in A_{1}\times A_{2}\times ...\times A_{k})=\bigcap\limits_{j=1}^{k}X_{i}^{-1}(A_{i}). \tag{invX} 
\end{equation}

\noindent \noindent Let us make the following remarks.\\

\noindent (R1) By \textit{Point (02.07b)} above, each component of $X$ is measurable whenever $X$ is.\\

\noindent (R2) By (InvX) and (WMC), $X$ is measurable whenever all its components are
measurable.\\

\noindent \textbf{Criterion} : $X$ is measurable if and only if each $X_{i}$ is
measurable.\newline

\newpage
\noindent \LARGE \textbf{DOC 02-02 : What you cannot ignore on Measurability}. \label{doc02-02}\\
\Large

\bigskip \noindent \textbf{Abstract}. Measurable applications are applications preserving
measurability of sets by means of inverse images. Before you continue, you
are asked to go back to \textit{Exercise 13 in Doc 01-02 in Chapter \ref{01_setsmes}} (page \pageref{doc01-02}) to revise
the properties of inverse images, in particular to get acquainted with the
fact that the inverse image preserves all sets operations.\newline

\bigskip \noindent \textbf{I - General introduction}.\newline

\noindent \textbf{(02.08) : Definition}. An application $X$ from the measurable
space $(\Omega _{1},\mathcal{A}_{1})$ onto the measurable space $(\Omega
_{2},\mathcal{A}_{2})$, is measurable, with respect to $(\mathcal{A}_{1}),%
\mathcal{A}_{2}))$, if

\begin{equation}
\forall (B\in \mathcal{A}_{2}),\text{ } X^{-1}(B)\in \mathcal{A}_{1}, 
\tag{M}
\end{equation}

\bigskip \noindent It is usual in measure theory to use the following
notation : 
\begin{equation*}
X^{-1}(B)=(X\in B)
\end{equation*}

\noindent \textbf{(02.09) : Criterion}. Fortunately, in many cases, $\mathcal{A}%
_{2}$ is generated by a small class, say $\mathcal{C}$, that is $\mathcal{A}%
_{2}=\sigma (\mathcal{C})$. In these cases, the application $X$ is
measurable if 
\begin{equation}
\forall (B\in \mathcal{C}),\text{ }X^{-1}(B)\in \mathcal{A}_{1},  \tag{M}
\end{equation}

\bigskip \noindent \textbf{(02.10) : Composition of measurable applications}. Let $%
X$ be a measurable application from $(\Omega _{1},\mathcal{A}_{1})$ onto $%
(\Omega _{2},\mathcal{A}_{2})$, and let $Y$ be a measurable application form 
$(\Omega _{2},\mathcal{A}_{2})$ onto $(\Omega _{3},\mathcal{A}_{3})$. Then
the composite function (function of function) $Z=Y(X)=YX$ defined from $%
(\Omega _{1},\mathcal{A}_{1})$ to $(\Omega _{3},\mathcal{A}_{3})$ by 
\begin{equation*}
\Omega _{1}\ni \omega \Rightarrow Z(\omega )=Y(X(\omega ))
\end{equation*}

\noindent is measurable.\newline

\bigskip

\bigskip \noindent \textbf{II - Operations on measurable functions using the
class of simple functions}.\newline

\bigskip \noindent \textbf{(02.11) - Indicator function of a measurable set}. For
any measurable space A in the measurable space $(\Omega ,\mathcal{A}),$ the
indication function%
\begin{equation*}
\begin{tabular}{lllll}
$1_{A}$ & : & $(\Omega ,\mathcal{A})$ & $\longrightarrow $ & $(\mathbb{R},%
\mathcal{B}(\mathbb{R}))$ \\ 
&  & $\omega $ & $\hookrightarrow $ & $1_{A}(\omega )=\left\{ 
\begin{tabular}{lll}
$1$ & $if$ & $\omega \in A$ \\ 
$0$ & $if$ & $\omega \in $
\end{tabular}
\right. $
\end{tabular}
\end{equation*}

\bigskip \noindent is measurable.\newline

\bigskip \noindent \textbf{(02.12) : Elementary functions}. A function 
\begin{equation*}
\begin{tabular}{lllll}
$f$ & : & $(\Omega ,\mathcal{A})$ & $\longrightarrow $ & $(\mathbb{R},
\mathcal{B}(\mathbb{R}))$
\end{tabular}
\end{equation*}

\bigskip \noindent is an elementary function if and only if it is finite and 
\newline

\bigskip \noindent (i) there exists a finite and measurable subdivision of $%
\Omega $%
\begin{equation*}
\Omega =A_{1}+...+A_{p},\text{ }p\geq 1
\end{equation*}

\noindent and\newline

\bigskip \noindent (ii) there exist $p$ real numbers $\alpha _{1},...,\alpha
_{p}$ such that

\bigskip 
\begin{equation}
\forall (1\leq i\leq p),\text{ }\forall \omega \in A_{i},\text{ }f(\omega
)=\alpha _{i}.  \tag{EF1}
\end{equation}

\bigskip \noindent Formula (EF1) is equivalent to%
\begin{equation*}
f=\alpha _{i}\text{ on }A_{i},\text{ }1\leq i\leq p
\end{equation*}

\bigskip \noindent and to
\begin{equation*}
f=\sum_{i=1}^{p}\alpha _{i}1_{A_{i}}.
\end{equation*}

\bigskip \noindent \textbf{Properties}. We have the following properties :\\

\bigskip \noindent \textbf{(02.13) : Operations in $\mathcal{E}$, the class of all elementary functions}.\newline

\noindent If $f$ and $g$ are elementary functions, and $c$ is a
real number, then $f+g, $ $fg,cf$, $\max (f,g)$ and $\min (f,g)$ are
elementary functions. If $g$ does not take the null value, then $f/g$ is an
elementary function.\\

\noindent \textbf{(02.14) : Approximation of a non-negative measurable
function}.\newline

\noindent For any non-negative and non-decreasing function 
\begin{equation*}
\begin{tabular}{lllll}
$f$ & : & $(\Omega ,\mathcal{A})$ & $\longrightarrow $ & $(\overline{\mathbb{%
R}},\mathcal{B}(\overline{\mathbb{R}}))$%
\end{tabular}%
\end{equation*}

\bigskip \noindent there exists a non-decreasing sequence of non-negative
elementary functions $(f_{n})_{n\geq 0}$ such that%
\begin{equation*}
f_{n}\nearrow f\text{ as }n\nearrow +\infty .
\end{equation*}

\bigskip \noindent \textbf{(02.15) : Approximation of an arbitrary measurable
function}.\newline

\bigskip \noindent Let 
\begin{equation*}
\begin{tabular}{lllll}
$f$ & : & $(\Omega ,\mathcal{A})$ & $\longrightarrow $ & $(\overline{\mathbb{%
R}},\mathcal{B}(\overline{\mathbb{R}}))$%
\end{tabular}%
\end{equation*}

\bigskip \noindent be measurable, we have the decomposition%
\begin{equation*}
f=f^{+}-f^{-},
\end{equation*}

\bigskip \noindent where 
\begin{equation*}
f^{+}=\max (f,0)\text{ \ and }f^{-}=\max (-f,0)
\end{equation*}

\bigskip \noindent \textbf{Fact } : $f^{+}$ and $f^{-}$ are positive and measurable. So there exist a non-decreasing sequences of non-negative
elementary functions $(f_{n}^{+})_{n\geq 0}$ and $(f_{n}^{-})_{n\geq 0}$ such that
\begin{equation*}
f_{n}^{+}\nearrow f^{+}\text{ and }f_{n}^{-}\nearrow f^{-}\text{ as }%
n\nearrow +\infty.
\end{equation*}

\bigskip \noindent Then, there exists a sequence of elementary functions $(f_{n})_{n\geq 0}$ \ \ $[f_{n}=f_{n}^{+}-f_{n}^{-}]$ \ such that
\begin{equation*}
f_{n}\longrightarrow f^{-}\text{ as }n\longrightarrow +\infty.
\end{equation*}

\bigskip \noindent \textbf{(02.16) : Operations on real-valued measurable functions}.\newline

\bigskip \noindent If $f$ and $g$ are measurable functions, and $c$ is a real number, then $f+g,$ $fg,cf$, $\max (f,g)$ and $\min (f,g)$ are measurable whenever theses \ functions are well-defined. If $g$ does not take the null value, then $f/g$ is measurable.\newline

\bigskip \noindent \textbf{III : Criteria of measurability for real-valued functions}.\newline

\noindent \textbf{(02.17)} A real-valued function 
\begin{equation*}
\begin{tabular}{lllll}
$f$ & : & $(\Omega ,\mathcal{A})$ & $\longrightarrow $ & $(\overline{\mathbb{R}},\mathcal{B}(\overline{\mathbb{R}}))$
\end{tabular}
\end{equation*}

\bigskip \noindent is measurable if and only if

\begin{equation*}
\forall a\in \mathbb{R},\text{ }(f>a)\in \mathcal{A}
\end{equation*}

\bigskip \noindent if and only if

\begin{equation*}
\forall a\in \mathbb{R},\text{ }(f\geq a)\in \mathcal{A}
\end{equation*}

\bigskip \noindent if and only if

\begin{equation*}
\forall a\in \mathbb{R},\text{ }(f<a)\in \mathcal{A}
\end{equation*}

\bigskip \noindent if and only if

\begin{equation*}
\forall a\in \mathbb{R},\text{ }(f\leq a)\in \mathcal{A}
\end{equation*}

\bigskip \noindent if and only if

\begin{equation*}
\forall (a,b)\in \mathbb{R}^{2},\text{ }(a<f<a)\in \mathcal{A}.
\end{equation*}

\bigskip \noindent And so forth.\newline

\bigskip

\bigskip \noindent \textbf{IV - Examples of real-valued measurable functions
from $\mathbb{R}$ to $\mathbb{R}$}.\newline

\bigskip \noindent \textbf{(02.18)} Right of left continuous functions,
including continuous functions.\newline

\bigskip \noindent \textbf{(02.19)} Function with at most a countable of
discontinuity points having left limits at each point or having right limit
at each point.\newline

\bigskip \noindent \textbf{(02.20)} Monotonic functions.\newline

\bigskip \noindent \textbf{(02.21)} Convex functions (convex functions are
continuous).\newline

\bigskip \noindent \textbf{V : Examples measurable functions defined from a
topological space (E,T) to $\mathbb{R}$}.\newline

\bigskip \noindent \textbf{(02.22)} Continuous functions.\newline

\noindent If $(E,T)$ is a metric space $(E,d)$.\newline

\noindent \textbf{(02.23)} lower semi-continuous functions.\newline

\noindent \textbf{(02.24)} upper semi-continuous functions.\newline

\bigskip \noindent \textbf{VI : Multi-components function}.\newline

\bigskip \noindent Consider an application with values in a product space
endowed with his product $\sigma$-algebra. \ 
\begin{equation*}
X:(E,\mathcal{B}\mapsto (\prod_{1\leq i\leq k}\Omega _{i},\bigotimes_{1\leq
i\leq k}\mathcal{A}_{i}).
\end{equation*}

\bigskip \noindent The image $X(\omega )$ has $k$ components. 
\begin{equation*}
X(\omega )=(X_{1}(\omega ),X_{2}(\omega ),...,X_{k}(\omega )),
\end{equation*}

\bigskip \noindent and each component $X_{i}$ is an application from $(E,%
\mathcal{B}$ to $(\Omega _{i},\mathcal{A}_{i})$.\newline

\bigskip \noindent \textbf{(0.24)Criterion} : $X$ is measurable if and only if
each $X_{i}$ is measurable.\newline

\newpage
\noindent \LARGE \textbf{DOC 02-3 : Measurable Applications - Exercises}. \label{doc02-03}\\
\Large
 
\bigskip \noindent \textbf{Exercise 1}. \label{exercise01_doc02-03} Let $X$ be an application from the
measurable $(\Omega _{1},\mathcal{A}_{1})$ to the measurable space $(\Omega
_{2},\mathcal{A}_{2})$.
\begin{equation*}
\mathcal{A}_{X}=\{(X\in B),B\in \mathcal{A}_{2}\}\subset \mathcal{A}_{1}.
\end{equation*}

\bigskip \noindent (a) Show that $\mathcal{A}_{X}$ is a $\sigma $-algebra on $\Omega_{1}.$\\

\noindent (b) Show that $X:(\Omega _{1},\mathcal{A}_{X})\longmapsto (\Omega _{2},%
\mathcal{A}_{2})$ is measurable (with respect $\mathcal{A}_{X}$ and $\mathcal{A}_{2}).$\newline

\noindent (c) Suppose that $X:(\Omega _{1},\mathcal{A}_{1})\longmapsto (\Omega _{2},\mathcal{A}_{2})$ is measurable. Prove that

\begin{equation}
\mathcal{A}_{X}\subset \mathcal{A}_{2}.  \tag{M}
\end{equation}

\noindent (d) Give a characterization of $\mathcal{A}_{X}.$\newline

\noindent \textbf{Definition}. $\mathcal{A}_{X}$ is the $\sigma $-algebra
generated by $X$ and $\mathcal{A}_{2}$ on $\Omega _{1}.$

\bigskip

\noindent \textbf{Exercise 2}. \label{exercise02_doc02-03} Let $A$ be a subset of the measurable space $%
(\Omega ,\mathcal{A}).$ Consider the Consider the canonical injection 
\begin{equation*}
i_{A}:A\rightarrow (\Omega ,\mathcal{A}),
\end{equation*}

\bigskip \noindent defined by : for all $x\in A,\ i_{A}(x)=x$.\newline

\noindent (a) Determine the $\sigma $-algebra generated by $i_{A}$ and $%
\mathcal{A}$ on $A$, denoted as $\mathcal{A}_{A}.$\newline

\noindent (b) What can you say about the collection of subsets of $A$ : $\{{%
B\cap A,B\in \mathcal{A\}}}$?\newline

\noindent (c) Show that $\mathcal{A}_{A}\subset \mathcal{A}$ if and only if $%
A$ is measurable (that is A$\in \mathcal{A}).$\newline

\bigskip \noindent \textbf{Exercise 3}. \label{exercise03_doc02-03}Let $X$ be a measurable application
from $(\Omega _{1},\mathcal{A}_{1})$ to $(\Omega _{2},\mathcal{A}_{2})$,
and let $Y$ be a measurable application form $(\Omega _{2},\mathcal{A}_{2})$
to $(\Omega _{3},\mathcal{A}_{3})$. Denote by $Z=Y(X)$ the composition of
the applications $X$ and $Y$.\newline

\noindent (a) Show that for any subset C of $\Omega _{3},$ 
\begin{equation*}
Z^{-1}(C)=X^{-1}(Y^{-1}(C)).
\end{equation*}

\bigskip \noindent (b) Show that $Z$ is measurable.

\bigskip \noindent \textbf{Exercise 4}. \label{exercise04_doc02-03} (Weak measurability condition) Let $%
X $ be an application from the measurable $(\Omega _{1},\mathcal{A}_{1})$ to
the measurable space $(\Omega _{2},\mathcal{A}_{2})$ and suppose that $%
\mathcal{A}_{2}$ is generated by a nonempty class $\mathcal{C}$ of subsets $\ \Omega _{2}$%
, that is $\mathcal{A}_{2}=\sigma (\mathcal{C}).$ Suppose that 
\begin{equation}
\forall (B\in \mathcal{C}),\text{ }X^{-1}(B)\in \mathcal{A}_{1},  \tag{WMC}
\end{equation}

\noindent (a) Show that the following class of subsets of $\Omega _{2}$%
\begin{equation*}
\mathcal{A}_{2,c}=\{B\in \mathcal{A}_{2},X^{-}(B)\in \mathcal{A}_{1}\}.
\end{equation*}

\bigskip \noindent satisfies the three properties : \newline

\noindent (i) $\mathcal{A}_{2,c}$ is a $\sigma $-algebra.\newline
\noindent (ii) $\mathcal{A}_{2,c}$ is included in $\mathcal{A}_{2}:$ $%
\mathcal{A}_{2,c}\subset \mathcal{A}_{2}.$\newline
\noindent (iii) $\mathcal{A}_{2,c}$ includes $\mathcal{C}$ : $\mathcal{%
C\subset A}_{2,c}.$\newline

\bigskip \noindent (b) Deduce from this that $\mathcal{A}_{2,c}=\mathcal{A}%
_{2}$ and then%
\begin{equation}
\forall (B\in \mathcal{A}_{1}),\text{ }X^{-1}(B)\in \mathcal{A}_{1}, 
\tag{WMC}
\end{equation}

\noindent that is, $X$ is measurable.\newline

\bigskip \noindent \textbf{Exercise 5}. \label{exercise05_doc02-03} (Applications)\newline

\noindent (a) \textbf{Continuity}. Let $X$ be an application from the
topological $(E_{1},\mathcal{T}_{1})$ to the topological space $(E_{2},%
\mathcal{T}_{2})$, where $\mathcal{T}_{1}$ and $\mathcal{T}_{2}$ are the
classes of open sets of the two topological spaces. By definition, $X$ is
continuous if and only if 
\begin{equation}
\forall (B\in \mathcal{T}_{2}),\text{ }X^{-1}(B)\in \mathcal{T}_{1},  \tag{C}
\end{equation}

\noindent Now let $\mathcal{B}(E_{1})=\sigma (\mathcal{T}_{1})$ and $%
\mathcal{B}(E_{2})=\sigma (\mathcal{T}_{2})$ be the corresponding Borel $%
\sigma $-algebras on $E_{1}$ and $E_{2}$ respectively.\newline

\noindent Show that $X$ is measurable whenever it is continuous.\newline

\bigskip \noindent (b) \textbf{(Multicomponent application)}. Consider two measurable spaces $(\Omega _{i},\mathcal{A}_{i})$, $i=1,2$. Let $\mathcal{A}=\mathcal{A}_{1}\otimes \mathcal{A}_{2}$ be the product sigma-algebra on the product space $\Omega =\Omega _{1}\times \Omega _{2}$. Recall that this $\sigma$-algebra is generated by : 
\begin{equation*}
\mathcal{S}=\{A\times B,A\in \mathcal{A}_{1},B\in \mathcal{A}_{2}\}.
\end{equation*}

\noindent Give a measurability criterion of applications from a measurable
space $(E,\mathcal{B})$ to ($\Omega _{1}\times \Omega _{2},\mathcal{A}%
_{1}\otimes \mathcal{A}_{2})$ using $\mathcal{S}.$\newline

\noindent (c) Let $X$ be a real-valued application defined on some
measurable space 
\begin{equation*}
X:(\Omega ,\mathcal{A})\longmapsto (\mathbb{R},\mathcal{B}(\mathbb{R})).
\end{equation*}

\bigskip \noindent Give four measurability criteria of $X$, each of them based on
one these generations : (1) $\mathcal{B}(\mathbb{R})=\sigma (\{]a,+\infty
\lbrack ,a\in \mathbb{R}\}),$ (2) $\mathcal{B}(\mathbb{R})=\sigma
(\{[a,+\infty \lbrack ,a\in \mathbb{R}\}),$ (3) $\mathcal{B}(\mathbb{R}%
)=\sigma (\{]-\infty ,a[,a\in \mathbb{R}\}),$ (4) $\mathcal{B}(\mathbb{R}%
)=\sigma (\{]-\infty ,a[,a\in \mathbb{R}\}).$

\bigskip \noindent \textbf{Exercise 6}. \label{exercise06_doc02-03} \textbf{Two product spaces}.
Consider two measurable spaces $(\Omega _{i},\mathcal{A}_{i})$, $i=1,2$. Let 
$\mathcal{A}=\mathcal{A}_{1}\otimes \mathcal{A}_{2}$ be the product
sigma-algebra on the product space $\Omega =\Omega _{1}\times \Omega _{2}$.
Recall that this sigma-algebra is generated by : 
\begin{equation*}
\mathcal{S}=\{A\times B,A\in \mathcal{A}_{1},B\in \mathcal{A}_{2}\}.
\end{equation*}

\bigskip \noindent Consider the two projections : $\pi _{1}:(\Omega ,\mathcal{A}%
)\rightarrow (\Omega _{1},\mathcal{A}_{1})$ and $\pi _{2}:(\Omega ,\mathcal{A%
})\rightarrow (\Omega _{2},\mathcal{A}_{2})$ respectively defined by : 
\begin{equation*}
\pi _{1}(\omega _{1},\omega _{2})=\omega _{1}\text{ and }\pi _{2}(\omega_{1},\omega _{2})=\omega _{2}.
\end{equation*}

\bigskip \noindent (a) Prove , for $A\subset \Omega _{1},$ $A\in \mathcal{%
A}_{1}$, and for for $B\subset \Omega _{2},$ $B\in \mathcal{A}_{2}$, that 
\begin{equation}
\pi _{1}^{-1}(A)=A_{1}\times \Omega _{2}\in \mathcal{S}\text{ \ \ and }\pi_{2}^{-1}(B)=\Omega _{1}\times B\in \mathcal{S}.  \label{rectPi2}
\end{equation}

\noindent Deduce from this that $\pi _{1}$ and $\pi _{2}$ are measurable with respect to $\mathcal{A}=\mathcal{A}_{1}\otimes \mathcal{A}_{2}$,  we have \newline

\noindent (b) Deduce from (\ref{rectPi2}) that for $A\times B\in \mathcal{S}, $ with $A\in \mathcal{A}_{1},B\in \mathcal{A}_{2},$ that 
\begin{equation}
A\times B=\pi _{1}^{-1}(A)\cap \pi _{2}^{-1}(B).  \label{rectPi2b}
\end{equation}

\bigskip \noindent (In case of problems, go back \textit{Doc 00-01, Point (0.21)}.\\

\bigskip \noindent (c) Now, let $\mathcal{B}$ be any $\sigma $-algebra on $%
\Omega =\Omega _{1}\times \Omega _{2}.$ And suppose the projections $\pi
_{i}:(\Omega _{1}\times \Omega _{2},\mathcal{B})\mapsto (\Omega _{i},%
\mathcal{A}_{i})$ are measurable with respect to $\mathcal{B}$ and $\mathcal{%
A}_{i}$, for $i=1,2.$ Use (\ref{rectPi2b}) to show that $\mathcal{S}\subset 
\mathcal{B}$ and next%
\begin{equation*}
\mathcal{A}=\mathcal{A}_{1}\otimes \mathcal{A}_{2}\subset \mathcal{B}.
\end{equation*}

\bigskip \noindent (d) Give a second definition of the product $\sigma $-algebra
using the measurability of the projections.\newline

\bigskip \noindent \textbf{Exercise 7}. \label{exercise07_doc02-03} Consider two measurable spaces $%
(\Omega _{i},\mathcal{A}_{i})$, $i=1,2$. Let $\mathcal{A}=\mathcal{A}%
_{1}\otimes \mathcal{A}_{2}$ be the product sigma-algebra on the product
space $\Omega =\Omega _{1}\times \Omega _{2}$. Recall that this
sigma-algebra is generated by : 
\begin{equation*}
\mathcal{S}=\{A\times B,A\in \mathcal{A}_{1},B\in \mathcal{A}_{2}\}.
\end{equation*}

\bigskip \noindent Let $X$ be an application defined on some measurable space $(E,%
\mathcal{B)}$ with values in the product space $\Omega =\Omega _{1}\times
\Omega _{2}$ endowed with the product space $\mathcal{A}=\mathcal{A}%
_{1}\otimes \mathcal{A}_{2}:$

\bigskip 
\begin{equation*}
\begin{tabular}{lllll}
$X$ & $:$ & $(E,\mathcal{B)}$ & $\longmapsto $ & $(\Omega _{1}\times \Omega
_{2},\mathcal{A}_{1}\otimes \mathcal{A}_{2})$ \\ 
&  & $\omega $ & $\hookrightarrow $ & $X(\omega )=(X_{1}(\omega
),X_{2}(\omega ))$.
\end{tabular}
\end{equation*}

\bigskip \noindent (a) Express each $X_{i}$ with respect to $X$ and $\pi _{i},$ $%
i=1,2. $ Deduce from that that each component $X_{i}$\ is measurable
whenever $X$ is measurable.\newline

\noindent (b) Show that for $A\times B\in \mathcal{S},$ with $A\in 
\mathcal{A}_{1},B\in \mathcal{A}_{2}$, we have that 
\begin{equation}
(X\in A\times B)=X_{1}^{-1}(A)\cap X_{2}^{-1}(B).
\end{equation}

\noindent Deduce from this that $X$ is measurable whenever its components $%
X_{1}$ and $X_{2}$ are measurable.\newline

\noindent (c) Conclude on the measurability of $X$ with respect to that of
its components.

\bigskip \noindent \textbf{NB}. Exercises 6 and 7 may be extended to an
arbitrary product of $k\geq 2$ spaces. The statements are already given in
\textit{Point 02.05 of Doc 02-01, page \pageref{doc02-01}}.\\

\bigskip \noindent \textbf{Exercise 8}. \label{exercise08_doc02-03}. Let $\Omega_0$ be a measurable set in the measurable space $(\Omega, \mathcal{A})$. Let 
$(\Omega_0, \mathcal{A}_{\Omega_0})$ be the induced measurable space on $\Omega_0$. \\

\noindent Let $(E,\mathcal{B})$ a measurable space, $f : (\Omega_0, \mathcal{A}_{\Omega_0}) \mapsto (E,\mathcal{B})$ be a measurable application and $A$ be an element of $E$.\\

\noindent Show that the extension :

\begin{equation*}
\tilde{f}(\omega)=\left\{ 
\begin{tabular}{lll}
$f(\omega)$ & if & $\omega \in \Omega_0$ \\ 
$A$ & if  & $\omega \notin \Omega_0$ \end{tabular}
\right. .
\end{equation*}

\bigskip \noindent is measurable as an application of $(\Omega, \mathcal{A})$ onto $(E,\mathcal{B})$.\\

\noindent \textbf{NB}. This simple exercise is very helpful in many situations, particularly when dealing with convergence of measurable applications.\\

\newpage
\noindent \LARGE \textbf{DOC 02-04 : Measurable Applications : exercises with solutions}. \label{doc02-04}\\
\Large
 
\bigskip \noindent \textbf{Exercise 1}. \label{exercise01_sol_doc02-04} Let $X$ be an application from the
measurable $(\Omega _{1},\mathcal{A}_{1})$ to the measurable space $(\Omega
_{2},\mathcal{A}_{2}).$%
\begin{equation*}
\mathcal{A}_{X}=\{(X\in B),B\in \mathcal{A}_{2}\}\subset \mathcal{A}_{1}.
\end{equation*}

\bigskip \noindent (a) Show that $\mathcal{A}_{X}$ is a $\sigma $-algebra on $\Omega_{1}.$\\

\noindent (b) Show that $X:(\Omega _{1},\mathcal{A}_{X})\longmapsto (\Omega _{2},%
\mathcal{A}_{2})$ is measurable (with respect $\mathcal{A}_{X}$ and $%
\mathcal{A}_{2}).$\newline

\noindent (c) Suppose that $X:(\Omega _{1},\mathcal{A}_{1})\longmapsto
(\Omega _{2},\mathcal{A}_{2})$ is measurable. Prove that

\begin{equation}
\mathcal{A}_{X}\subset \mathcal{A}_{2}.  \tag{M}
\end{equation}

\bigskip  \noindent (d) Give a characterization of $\mathcal{A}_{X}.$\newline

\noindent \textbf{Definition}. $\mathcal{A}_{X}$ is the $\sigma $-algebra
generated by $X$ and $\mathcal{A}_{2}$ on $\Omega _{1}.$\\

\bigskip \noindent \textbf{SOLUTIONS}.\\

\noindent Question (a). See Exercise 6 in Doc 01-02 and its solution in Doc 01-03, pages \ref{exercise06_doc01-02}, \ref{exercise06_sol_doc01-03}.\\

\bigskip \noindent Question (b). By definition, for any $B\in \mathcal{A}_{2},(X\in
B)=X^{-1}(B)\in \mathcal{A}_{X}.$ So $X$ est $\mathcal{A}_{X}-\mathcal{A}_{2}
$ measurable.\\

\bigskip \noindent Question (c) Let $X:(\Omega _{1},\mathcal{A}_{1})\longmapsto (\Omega _{2},%
\mathcal{A}_{2})$ measurable. Then for any $B\in \mathcal{A}_{2},$ $(X\in
B)=X^{-1}(B)\in \mathcal{A}_{1}.$ It comes that all the elements $\mathcal{A}%
_{X}$ are in $\mathcal{A}_{1}.$\\

\bigskip \noindent Question (d) $\mathcal{A}_{X}$ is the smallest $\sigma $-algebra on $\Omega
_{1}$ rendering $X$ measurable.\\

\bigskip

\noindent \textbf{Exercise 2}. \label{exercise02_sol_doc02-04} Let $A$ be a subset of the measurable space $%
(\Omega ,\mathcal{A}).$ Consider the Consider the canonical injection 
\begin{equation*}
i_{A}:A\rightarrow (\Omega ,\mathcal{A}),
\end{equation*}

\bigskip \noindent defined by : for all $x\in A,\ i_{A}(x)=x$.\newline

\noindent (a) Determine the $\sigma $-algebra generated by $i_{A}$ and $%
\mathcal{A}$ on $A$, denoted as $\mathcal{A}_{A}.$\newline

\noindent (b) What can you say about the collection of subsets of $A$ : $\{{%
B\cap A,B\in \mathcal{A\}}}$?\newline

\noindent (c) Show that $\mathcal{A}_{A}\subset \mathcal{A}$ if and only if $%
A$ is measurable (that is $A\in \mathcal{A}).$\newline

\bigskip \noindent \textbf{SOLUTIONS}.\\

\noindent Question (a). For any $B\subset \Omega ,$

\begin{eqnarray*}
i_{A}^{-1}(B) &=&\{\omega \in A,i_{A}(\omega )\in B\} \\
&=&\{\omega \in A,\omega \in B\} \\
&=&A\cap B.
\end{eqnarray*}

\bigskip \noindent Then, by definition 
\begin{eqnarray*}
\mathcal{A}_{A} &=&\{i_{A}^{-1}(B),B\in \mathcal{A\}} \\
&=&\{A\cap B,B\in \mathcal{A\}}\text{.}
\end{eqnarray*}

\bigskip \noindent Question (b) The collection $\{A\cap B,B\in \mathcal{A\}}$ is then the $%
\sigma $-algebra generated by the canonical injection $i_{A}.$\\

\noindent Question (c) It is clear that if $A$ is measurable (that is $A\in \mathcal{A}%
)$, then all the $A\cap B,B\in \mathcal{A},$ are measurable (in \ $\mathcal{A%
}$). Now, if $\mathcal{A}_{A}\subset \mathcal{A},$ the for the particular
case $B=\Omega $, we get $A=A\cap \Omega \in \mathcal{A}.$

\bigskip

\bigskip \noindent \textbf{Exercise 3}. \label{exercise03_sol_doc02-04} Let $X$ be a measurable application
from $(\Omega _{1},\mathcal{A}_{1})$ to $(\Omega _{2},\mathcal{A}_{2})$,
and let $Y$ be a measurable application form $(\Omega _{2},\mathcal{A}_{2})$
to $(\Omega _{3},\mathcal{A}_{3})$. Denote by $Z=Y(X)$ the composition of
the applications $X$ and $Y$.\newline

\noindent (a) Show that for any subset C of $\Omega _{3},$ 
\begin{equation*}
Z^{-1}(C)=X^{-1}(Y^{-1}(C)).
\end{equation*}

\bigskip \noindent (b) Show that $Z$ is measurable.\\

\noindent \textbf{SOLUTIONS}.\\

\noindent Question (a) For any $C\subset \Omega _{3},$%
\begin{eqnarray*}
Z^{-1}(C) &=&\{\omega \in \Omega _{1},Z(\omega )\in C\} \\
&=&\{\omega \in \Omega _{1},Y(X(\omega ))\in C\}.
\end{eqnarray*}

\bigskip \noindent But, by the definition of the inverse image set,%
\begin{equation*}
Y(X(\omega ))\in C\Longleftrightarrow X(\omega )\in Y^{-1}(C)=B.
\end{equation*}%

\noindent So
\begin{eqnarray*}
Z^{-1}(C) &=&\{\omega \in \Omega _{1},Z(\omega )\in C\} \\
&=&\{\omega \in \Omega _{1},Y(X(\omega ))\in C\} \\
&=&\{\omega \in \Omega _{1},X(\omega )\in B\} \\
&=&X^{-1}(B)=X^{-1}(Y^{-1}(C)).
\end{eqnarray*}

\bigskip \noindent Question (b). Let $C\in \mathcal{A}_{3}.$\ Since $Y$ is measurable, we have $%
B=Y^{-1}(C)$ is measurable (on \ $\Omega _{2}$). And, since $X$ is
measurable, $X^{-1}(B)$ is measurable (on $\Omega _{1}$). But \ $%
X^{-1}(B)=Z^{-1}(C).$ So for any $C\in \mathcal{A}_{3}$, $Z^{-1}(C)\in 
\mathcal{A}_{1}$. Hence $Z$ is measurable as a mapping of $\Omega _{1}$ into 
$\Omega _{3}.$

\bigskip \noindent \textbf{Exercise 4}. \label{exercise04_sol_doc02-04} (Weak measurability condition) Let $X
$ be an application from the measurable $(\Omega _{1},\mathcal{A}_{1})$ to
the measurable space $(\Omega _{2},\mathcal{A}_{2})$ and suppose that $%
\mathcal{A}_{2}$ is generated by a nonempty class $\mathcal{C}$ of subsets $\ \Omega _{2}$%
, that is $\mathcal{A}_{2}=\sigma (\mathcal{C}).$ Suppose that 
\begin{equation}
\forall (B\in \mathcal{C}),\text{ }X^{-1}(B)\in \mathcal{A}_{1},  \tag{WMC}
\end{equation}

\noindent (a) Show that the following class of subsets of $\Omega _{2}$%
\begin{equation*}
\mathcal{A}_{2,c}=\{B\in \mathcal{A}_{2},X^{-}(B)\in \mathcal{A}_{1}\}.
\end{equation*}

\bigskip \noindent satisfies the three properties : \newline

\noindent (i) $\mathcal{A}_{2,c}$ is a $\sigma $-algebra.\newline
\noindent (ii) $\mathcal{A}_{2,c}$ is included in $\mathcal{A}_{2}:$ $%
\mathcal{A}_{2,c}\subset \mathcal{A}_{2}.$\newline
\noindent (iii) $\mathcal{A}_{2,c}$ includes $\mathcal{C}$ : $\mathcal{%
C\subset A}_{2,c}.$\newline

\bigskip \noindent (b) Deduce from this that $\mathcal{A}_{2,c}=\mathcal{A}%
_{2}$ and then%
\begin{equation}
\forall (B\in \mathcal{A}_{1}),\text{ }X^{-1}(B)\in \mathcal{A}_{1}, 
\tag{WMC}
\end{equation}

\noindent that is, $X$ is measurable.\newline

\bigskip 
\noindent \textbf{SOLUTIONS}.\\

\noindent Question (a).\\

\noindent Point (i). Use the properties of the inverse images of sets (Doc 01-02,
Exercise 13, solutions in Doc 01-03). We are going to prove that $\mathcal{A}%
_{2,c}$ is a $\sigma $-algebra on $\Omega _{2}.$ Let us check the following
points.\\

\noindent (1) $\Omega _{2}\in $ $\mathcal{A}_{2,c}$? Yes since $X^{-1}(\Omega
_{2})=\Omega _{1}\in \mathcal{A}_{1}.$\\

\noindent (2) Let $B\in \mathcal{A}_{2,c}.$ Do we have $B^{c}\in \mathcal{A}_{2,c}?$
Yes since $X^{-1}(B)\in \mathcal{A}_{1}$ and $X^{-1}(B^{c})=\left(
X^{-1}(B)\right) ^{c}$ . Then $X^{-1}(B^{c})\in \mathcal{A}_{1}$ and thus, $%
B\in \mathcal{A}_{2,c}.$\\

\noindent (3) Let $(B_{n})_{n\geq 1}\subset \mathcal{A}_{2,c}.$ Do we have $\cup
_{n\geq 1}B_{n}\in \mathcal{A}_{2,c}?$ Yes since%
\begin{equation*}
X^{-1}(\cup _{n\geq 1}B_{n})=\cup _{n\geq 1}X^{-1}(B_{n})\in \mathcal{A}_{1}
\end{equation*}

\bigskip \noindent since, by assumption, each $X^{-1}(B_{n})\in \mathcal{A}_{1}.$ Thus $\cup _{n\geq 1}B_{n}\in 
\mathcal{A}_{2,c}.$

\bigskip \noindent Point (ii). We have $\mathcal{A}_{2,c}\subset \mathcal{A}_{2}$ by
construction.\\

\bigskip \noindent Point (iii). We have $\mathcal{C\subset A}_{2,c}$ by the assumption (WMC).\\

\bigskip \noindent Question (b). Since  $\mathcal{A}_{2,c}$ is a $\sigma $-algebra including $%
\mathcal{C}$ then it includes $\sigma (\mathcal{C})=\mathcal{A}_{2}.$
Now, since $\mathcal{A}_{2,c}$ is in $\mathcal{A}_{2}$ by construction, we
get $\mathcal{A}_{2,c}=\mathcal{A}_{2}.$ We get 
\begin{equation}
\forall (B\in \mathcal{A}_{2}),\text{ }X^{-1}(B)\in \mathcal{A}_{1}.  \tag{M}
\end{equation}

\noindent This means that $X$ is measurable. So Condition (WMC) implies (M).
Reversely, (M) obviously implies (WMC). We conclude that $X$ is mesaurable
if and only if (WMC) holds.

\bigskip

\bigskip \noindent \textbf{Exercise 5}. \label{exercise05_sol_doc02-04} (Applications)\newline

\noindent (a) \textbf{Continuity}. Let $X$ be an application from the
topological $(E_{1},\mathcal{T}_{1})$ to the topological space $(E_{2},%
\mathcal{T}_{2})$, where $\mathcal{T}_{1}$ and $\mathcal{T}_{2}$ are the
classes of open sets of the two topological spaces. By definition, $X$ is
continuous if and only if 
\begin{equation}
\forall (B\in \mathcal{T}_{2}),\text{ }X^{-1}(B)\in \mathcal{T}_{1},  \tag{C}
\end{equation}

\noindent Now let $\mathcal{B}(E_{1})=\sigma (\mathcal{T}_{1})$ and $%
\mathcal{B}(E_{2})=\sigma (\mathcal{T}_{2})$ be the corresponding Borel $%
\sigma $-algebras on $E_{1}$ and $E_{2}$ respectively.\newline

\noindent Show that $X$ is measurable whenever it is continuous.\newline

\bigskip \noindent (b) \textbf{(Multi-component application)}. Consider two
measurable spaces $(\Omega _{i},\mathcal{A}_{i})$, $i=1,2$. Let $\mathcal{A}=%
\mathcal{A}_{1}\otimes \mathcal{A}_{2}$ be the product $\sigma$-algebra on the
product space $\Omega =\Omega _{1}\times \Omega _{2}$. Recall that this
$\sigma$ is generated by : 
\begin{equation*}
\mathcal{S}=\{A\times B,A\in \mathcal{A}_{1},B\in \mathcal{A}_{2}\}.
\end{equation*}

\bigskip \noindent State a measurability criterion of applications from a measurable space $(E,\mathcal{B})$ to ($\Omega
_{1}\times \Omega _{2},\mathcal{A}_{1}\otimes \mathcal{A}_{2})$ using $\mathcal{S}.$\newline
using $\mathcal{S}.$\newline

\noindent (c) Let $X$ be a real-valued application defined on some
measurable space 
\begin{equation*}
X:(\Omega ,\mathcal{A})\longmapsto (\mathbb{R},\mathcal{B}(\mathbb{R})).
\end{equation*}

\noindent Give four measurability criteria of $X$, each of them based on
one these generations : (1) $\mathcal{B}(\mathbb{R})=\sigma (\{]a,+\infty
\lbrack ,a\in \mathbb{R}\}),$ (2) $\mathcal{B}(\mathbb{R})=\sigma
(\{[a,+\infty \lbrack ,a\in \mathbb{R}\}),$ (3) $\mathcal{B}(\mathbb{R}%
)=\sigma (\{]-\infty ,a[,a\in \mathbb{R}\}),$ (4) $\mathcal{B}(\mathbb{R}%
)=\sigma (\{]-\infty ,a[,a\in \mathbb{R}\}).$\\

\bigskip \textbf{SOLUTIONS}.\\

\noindent Question (a) Let $X$ be continuous. Use the criterion (WMC) based on $%
\mathcal{T}_{2}$ that generates $\mathcal{B}(E_{2}).$ We have by continuity,
for any $G\in \mathcal{T}_{2}$ (open in $E_{2}),$%
\begin{equation*}
X^{-1}(G)\in \mathcal{T}_{1}\subset \mathcal{B}(E_{1}).
\end{equation*}

\bigskip \noindent So $X$ is measurable.\\

\noindent Question (b) By (WMC), an application $X$ from a measurable space $(E,\mathcal{B})$ to
($\Omega _{1}\times \Omega _{2},\mathcal{A}_{1}\otimes \mathcal{A}_{2})$ is
measurable if and only if for any $A\times B,$ with $A\in \mathcal{A}_{1}$
and $B\in \mathcal{A}_{2},$ \newline
\begin{equation*}
(X\in A\times B)\in \mathcal{B}.
\end{equation*}

\bigskip \noindent Since $X=(X_{1},X_{2})$, we have : $X$ is measurable if and only if 
\begin{equation*}
((X_{1},X_{2})\in A\times B)=X_{1}^{-1}(A)\cap X_{2}^{-1}(B)\in \mathcal{B}.
\end{equation*}

\bigskip \noindent Question (c). As usual, we denote%
\begin{equation*}
X^{-1}(]a,+\infty \lbrack )=\{\omega \in \Omega ,X(\omega )>a)=(X>a)
\end{equation*}%
\begin{equation*}
X^{-1}([a,+\infty \lbrack )=\{\omega \in \Omega ,X(\omega )\geq a)=(X>a)
\end{equation*}%
\begin{equation*}
X^{-1}(]-\infty ,a[)=\{\omega \in \Omega ,X(\omega )<a)=(X<a)
\end{equation*}%

\noindent and
\begin{equation*}
X^{-1}(]-\infty ,a])=\{\omega \in \Omega ,X(\omega )\leq a)=(X\leq a).
\end{equation*}

\bigskip \noindent Based on these four generations for the usual $\sigma $-algebra of $%
\mathbb{R}$, we have the following criteria : $X:(\Omega ,\mathcal{A}%
)\longmapsto (\mathbb{R},\mathcal{B}(\mathbb{R}))$ is measurable if and only 
\begin{equation*}
\forall (a\in \mathbb{R}),(X>a)\in \mathcal{A}
\end{equation*}%
$\Longleftrightarrow $%
\begin{equation*}
\forall (a\in \mathbb{R}),(X\geq a)\in \mathcal{A}
\end{equation*}%
$\Longleftrightarrow $%
\begin{equation*}
\forall (a\in \mathbb{R}),(X>a)\in \mathcal{A}
\end{equation*}%
$\Longleftrightarrow $%
\begin{equation*}
\forall (a\in \mathbb{R}),(X>a)\in \mathcal{A}
\end{equation*}

\bigskip \noindent \textbf{Exercise 6}. \label{exercise06_sol_doc02-04} \textbf{Two product spaces}.
Consider two measurable spaces $(\Omega _{i},\mathcal{A}_{i})$, $i=1,2$. Let 
$\mathcal{A}=\mathcal{A}_{1}\otimes \mathcal{A}_{2}$ be the product
sigma-algebra on the product space $\Omega =\Omega _{1}\times \Omega _{2}$.
Recall that this sigma-algebra is generated by : 
\begin{equation*}
\mathcal{S}=\{A\times B,A\in \mathcal{A}_{1},B\in \mathcal{A}_{2}\}.
\end{equation*}

\bigskip \noindent Consider the two projections : $\pi _{1}:(\Omega ,\mathcal{A}%
)\rightarrow (\Omega _{1},\mathcal{A}_{1})$ and $\pi _{2}:(\Omega ,\mathcal{A%
})\rightarrow (\Omega _{2},\mathcal{A}_{2})$ respectively defined by : 
\begin{equation*}
\pi _{1}(\omega _{1},\omega _{2})=\omega _{1}\text{ and }\pi _{2}(\omega
_{1},\omega _{2})=\omega _{2}.
\end{equation*}

\bigskip \noindent (a) Prove, for $A\subset \Omega _{1},$ $A\in \mathcal{%
A}_{1}$, and for for $B\subset \Omega _{2},$ $B\in \mathcal{A}_{2}$, that 
\begin{equation}
\pi _{1}^{-1}(A)=A_{1}\times \Omega _{2}\in \mathcal{S}\text{ \ \ and }\pi
_{2}^{-1}(B)=\Omega _{1}\times B\in \mathcal{S}.  \label{rectPi2}
\end{equation}

\noindent Deduce from this, that $\pi _{1}$ and $\pi _{2}$ are measurable
with respect to $\mathcal{A}=\mathcal{A}_{1}\otimes \mathcal{A}_{2}$, we have :\newline

\noindent (b) Deduce from (\ref{rectPi2}) that for $A\times B\in \mathcal{S}%
, $ with $A\in \mathcal{A}_{1},B\in \mathcal{A}_{2},$ that 
\begin{equation}
A\times B=\pi _{1}^{-1}(A)\cap \pi _{2}^{-1}(B).  \label{rectPi2b}
\end{equation}

\bigskip \noindent (In case of problems, go back to Tutorial 00, Point 0.21
of Doc 00-01).\newline

\bigskip \noindent (c) Now, let $\mathcal{B}$ be any $\sigma $-algebra on $%
\Omega =\Omega _{1}\times \Omega _{2}.$ And suppose the projections $\pi
_{i}:(\Omega _{1}\times \Omega _{2},\mathcal{B})\mapsto (\Omega _{i},%
\mathcal{A}_{i})$ are measurable with respect to $\mathcal{B}$ and $\mathcal{%
A}_{i}$, for $i=1,2.$ Use (\ref{rectPi2b}) to show that $\mathcal{S}\subset 
\mathcal{B}$ and next%
\begin{equation*}
\mathcal{A}=\mathcal{A}_{1}\otimes \mathcal{A}_{2}\subset \mathcal{B}.
\end{equation*}

\bigskip \noindent (d) Give a second definition of the product $\sigma $-algebra
using the measurability of the projections.\newline

\bigskip \noindent \textbf{SOLUTIONS}.\\

\noindent Question (a) By definition, for $A\subset \Omega _{1},$%
\begin{eqnarray*}
\pi _{1}^{-1}(A) &=&\{(\omega _{1},\omega _{2})\in \Omega _{1}\times \Omega
_{2},\pi _{1}(\omega _{1},\omega _{2})\in A\} \\
&=&\{(\omega _{1},\omega _{2})\in \Omega _{1}\times \Omega _{2},\omega
_{1}\in A\} \\
&=&\{\omega _{1}\in \Omega _{1},\omega _{2}\in \Omega _{2},\omega _{1}\in A\}
\\
&=&\{\omega _{1}\in \Omega _{1}\cap A,\omega _{2}\in \Omega _{2}\} \\
&=&\{\omega _{1}\in A,\omega _{2}\in \Omega _{2}\} \\
&=&A\times \Omega _{2}.
\end{eqnarray*}

\bigskip \noindent Likely, you can prove, for $B\subset \Omega _{2},$%
\begin{equation*}
\pi _{2}^{-1}(B)=\Omega _{1}\times B.
\end{equation*}

\noindent Then for any $A\in \mathcal{A}_{1}$ and for any $B\in \mathcal{A}_{2},$%
\begin{equation*}
\pi _{1}^{-1}(A)=A\times \Omega _{2}\in \mathcal{S}\subset \mathcal{A}=%
\mathcal{A}_{1}\otimes \mathcal{A}_{2}
\end{equation*}

\bigskip \noindent and 

\bigskip 
\begin{equation*}
\pi _{2}^{-1}(B)=\Omega _{1}\times B\in \mathcal{S}\subset \mathcal{A}=%
\mathcal{A}_{1}\otimes \mathcal{A}_{2}.
\end{equation*}

\bigskip \noindent Then the projections are measurable with respect to $\mathcal{A}=\mathcal{A}%
_{1}\otimes \mathcal{A}_{2}.$\\

\bigskip \noindent Question (b) Let $A\times B\in \mathcal{S},$ with $A\in \mathcal{A}_{1}$ and 
$B\in \mathcal{A}_{2}.$ By Question (a) and the property on the intersection
of rectangles,%
\begin{eqnarray*}
\pi _{1}^{-1}(B)\cap \pi _{2}^{-1}(B) &=&(A\times \Omega _{2})\cap (\Omega
_{1}\times B) \\
&=&(A\cap \Omega _{1})\times (\Omega _{2}\cap B) \\
&=&A\times B.
\end{eqnarray*}

\bigskip \noindent Then if, for a $\sigma $-algebra $\mathcal{B}$, $\pi _{1}:(\Omega _{1}\times
\Omega _{2},\mathcal{B})\mapsto (\Omega _{1},\mathcal{A}_{1})$ and $\pi
_{2}:(\Omega _{1}\times \Omega _{2},\mathcal{B})\mapsto (\Omega _{2},%
\mathcal{A}_{2})$ are measurable, then for $A\times B\in \mathcal{S},$ with $%
A\in \mathcal{A}_{1}$ and $B\in \mathcal{A}_{2},$%
\begin{equation*}
\pi _{1}^{-1}(A)\in \mathcal{B}\text{ and }\pi _{2}^{-1}(B)\in \mathcal{B},
\end{equation*}%

\bigskip \noindent and thus
\begin{equation*}
A\times B=\pi _{1}^{-1}(B)\cap \pi _{2}^{-1}(B)\in \mathcal{B}.
\end{equation*}

\bigskip \noindent This implies that $\mathcal{B}$\ includes $\mathcal{S}$ and then

\begin{equation*}
\mathcal{A}=\mathcal{A}_{1}\otimes \mathcal{A}_{2}=\sigma (\mathcal{S})\subset \mathcal{B}\text{.}
\end{equation*}

\bigskip \noindent \textbf{Conclusion} : any $\sigma $-algebra $\mathcal{B}$ rendering the
projections $\pi _{1}:(\Omega _{1}\times \Omega _{2},\mathcal{B})\mapsto
(\Omega _{1},\mathcal{A}_{1})$ and $\pi _{2}:(\Omega _{1}\times \Omega _{2},\mathcal{B})\mapsto (\Omega _{2},\mathcal{A}_{2})$ measurable, includes the product $\sigma $-algebra.\\

\bigskip \noindent Question (c) The product $\sigma $-algebra $\mathcal{A}=\mathcal{A}_{1}\otimes \mathcal{A}_{2}$ on the product space $\Omega _{1}\times \Omega_{2}$ is the smallest $\sigma $-algebra for which the projections are
measurable.\\

\bigskip \noindent \textbf{Exercise 7}. \label{exercise07_sol_doc07-04} Consider two measurable spaces $%
(\Omega _{i},\mathcal{A}_{i})$, $i=1,2$. Let $\mathcal{A}=\mathcal{A}%
_{1}\otimes \mathcal{A}_{2}$ be the product $\sigma$-algebra on the product
space $\Omega =\Omega _{1}\times \Omega _{2}$. Recall that this
$\sigma$-algebra is generated by : 
\begin{equation*}
\mathcal{S}=\{A\times B,A\in \mathcal{A}_{1},B\in \mathcal{A}_{2}\}
\end{equation*}

\bigskip \noindent Let $X$ be an application defined on some measurable space $(E,%
\mathcal{B)}$ with values in the product space $\Omega =\Omega _{1}\times
\Omega _{2}$ endowed with the product $\sigma$-algebra $\mathcal{A}=\mathcal{A}%
_{1}\otimes \mathcal{A}_{2}:$

\bigskip 
\begin{equation*}
\begin{tabular}{lllll}
$X$ & $:$ & $(E,\mathcal{B)}$ & $\longmapsto $ & $(\Omega _{1}\times \Omega
_{2},\mathcal{A}_{1}\otimes \mathcal{A}_{2})$ \\ 
&  & $\omega $ & $\hookrightarrow $ & $X(\omega )=(X_{1}(\omega
),X_{2}(\omega ))$
\end{tabular}
\end{equation*}

\bigskip \noindent (a) Express each $X_{i}$ with respect to $X$ and $\pi _{i},$ $%
i=1,2. $ Deduce from that that each component $X_{i}$\ is measurable
whenever $X$ is measurable.\newline

\noindent (b) Show that,t for $A\times B\in \mathcal{S},$ with $A\in 
\mathcal{A}_{1},B\in \mathcal{A}_{2},$ we have
\begin{equation}
(X\in A\times B)=X_{1}^{-1}(A)\cap X_{2}^{-1}(B).
\end{equation}

\noindent Deduce from this, that $X$ is measurable whenever its components $%
X_{1}$ and $X_{2}$ are measurable.\newline

\noindent (c) Conclude on the measurability of $X$ with respect to that of
its components.

\bigskip \noindent \textbf{NB}. Exercises 6 and 7 may be extended to an
arbitrary product of $k\geq 2$ spaces. The statements are already given
in \textit{Point 02.05 of Doc 02-01, page \pageref{doc02-01}}.\\

\bigskip \noindent \textbf{SOLUTIONS}.\\
 
\noindent Question (a). It is clear that for $X=(X_{1},X_{2}),$ $X_{1}=\pi
_{1}(X_{1},X_{2})=\pi _{1}(X)$ and $X_{2}=\pi _{2}(X_{1},X_{2})=\pi _{2}(X).$
Then, if $X$ is measurable,  $X_{1}$ and $X_{2}$ are compositions of
measurable applications, and then, they are measurable.\\

\bigskip \noindent  Question (b). Let $X_{1}$ and $X_{2}$ be measurable. By Question (b) of
Exercise 5 above, we have \ for any $A\times B\in \mathcal{S},$ with $A\in 
\mathcal{A}_{1}$ and $B\in \mathcal{A}_{2},$ that

\begin{equation*}
((X_{1},X_{2})\in A\times B)=X_{1}^{-1}(A)\cap X_{2}^{-1}(B)\in \mathcal{B.}
\end{equation*}

\bigskip \noindent And then $X=(X_{1},X_{2})$ is measurable.\\

\bigskip \noindent \textbf{Exercise 8}. \label{exercise08_sol_doc02-04}. Let $\Omega_0$ be a measurable set in the measurable space $(\Omega, \mathcal{A})$. Let 
$(\Omega_0, \mathcal{A}_{\Omega_0})$ be the induced measurable space on $\Omega_0$. \\

\noindent Let $(E,\mathcal{B})$ a measurable space, $f : (\Omega_0, \mathcal{A}_{\Omega_0}) \mapsto (E,\mathcal{B})$ be a measurable application and $A$ be an element of $E$.\\

\noindent Show that the extension :

\begin{equation*}
\tilde{f}(\omega)=\left\{ 
\begin{tabular}{lll}
$f(\omega)$ & if & $\omega \in \Omega_0$ \\ 
$A$ & if  & $\omega \notin \Omega_0$.
\end{tabular}
\right. 
\end{equation*}

\bigskip \noindent is measurable as an application of $(\Omega, \mathcal{A})$ onto $(E,\mathcal{B})$.\\

\noindent \textbf{NB}. This simple exercise is very helpful in many situations, particularly when dealing with convergence of measurable applications.\\

\bigskip \noindent \textbf{Solution}.\\

\noindent By definition, we have

$$
\mathcal{A}_{\Omega_0}=\{B \cap \Omega_0, \ B\in \mathcal{A} \},
$$ 

\bigskip \noindent and since $\Omega_0$ is measurable, we also have
$$
\mathcal{A}_{\Omega_0}=\{B, \ B\in \mathcal{A} \ and \ B \subset A\}.
$$ 

\bigskip \noindent Now, let $C \in \mathcal{B}$, we have

$$
(\tilde{f} \in C)=(\tilde{f} \in C) \cap \Omega_0 + (\tilde{f} \in C) \cap (\Omega_0)^c.
$$

\bigskip \noindent By definition of $\tilde{f}$ and using its values both on $\Omega_0$ and on $\Omega_0^c$, this yields

$$
(\tilde{f} \in C)=(f \in C) \cap \Omega_0 + (f \in C) \cap (\Omega_0)^c.
$$

\bigskip \noindent But $(f \in C) \cap \Omega_0 \in \mathcal{A}$ since $(f \in C) \in \mathcal{A}_{\Omega_0} \subset \mathcal{A}$ and $(f\in C)$ is either the empty set or the full set $\Omega$.\\

\noindent Hence $(\tilde{f} \in C) \in \mathcal{C}$ for any $C \in \mathcal{B}$. In conclusion $\tilde{f}$ is measurable.\\

%\chapter{Measurability of Sets and Application is Usual Spaces}
\chapter{Measurability is Usual Spaces} \label{03_setsmes_applimes_cas_speciaux}

\bigskip
\noindent \textbf{Content of the Chapter}\\

\begin{table}[htbp]
	\centering
		\begin{tabular}{llll}
		\hline
		Type& Name & Title  & page\\
		\hline
		S & Doc 03-01 & Borel sets  and $\mathbb{R}^d$-valued Measurable functions& \pageref{doc03-01}\\
		S & Doc 03-02 & Elementary functions  & \pageref{doc03-02} \\
		S& Doc 03-03 & Measurable applications - A summary  & \pageref{doc03-03}\\
		D& Doc 03-04 & Exercises on Measurability 02 01& \\
		&   & Borel sets and Borel functions.   &  \pageref{doc03-04}\\
		D& Doc 03-05 & Exercises on Measurability 02& \\
		&   & Elementary functions, Marginal functions and applications   &  \pageref{doc03-05}\\
		D& Doc 03-06 & Exercises on Measurability 03& \\
		&   & Classical functions in Real Analysis  &  \pageref{doc03-06}\\
		SD& Doc 03-07 & Exercises on Measurability 01 with solutions & \pageref{doc03-07}\\
		SD& Doc 03-08 & Exercises on Measurability 02 with solutions & \pageref{doc03-08}\\
		SD& Doc 03-09 & Exercises on Measurability 03 with solutions & \pageref{doc03-09}\\
		\hline
		\end{tabular}
\end{table}

\newpage
\noindent \LARGE \textbf{DOC 03-01 : Borel sets $\mathbb{R}^d$ and $\mathbb{R}^d$-valued Measurable functions}. \label{doc03-01}\\
\Large

\bigskip \noindent \noindent Part I - The space $\mathbb{R}$.\newline

\noindent \textbf{(03.01)} The usual $\sigma $-algebra on $\overline{\mathbb{R}}$, denoted $\mathbb{B}(\mathbb{R})$, is generated by the whole class of intervals of $\mathbb{R}$ denoted $\mathcal{I}$. And we know that $\mathbb{I}$ is a semi-algebra. But we have more detailed things.\newline

\bigskip \noindent \textbf{(03.02)} Actually $\mathcal{B}(\mathbb{R})$ is generated by any class of one particular type of intervals.\newline

\bigskip \noindent Define :\newline

\begin{equation*}
\mathcal{I}(1)=\{]-\infty ,a[,a\in \mathbb{R}\}\text{ , }\mathcal{I}%
(2)=\{]-\infty ,a],\ \ a\in \overline{\mathbb{R}}\};
\end{equation*}

\bigskip 
\begin{equation*}
\mathcal{I}(3)=\{]-\infty ,a],a\in \mathbb{R}\}\text{ , }\mathcal{I}%
(4)=\{]-\infty ,a],\ \ a\in \overline{\mathbb{R}}\}
\end{equation*};

\bigskip 
\begin{equation*}
\mathcal{I}(5)=\{]a,+\infty \lbrack ,a\in \mathbb{R}\},\mathcal{I}%
(6)=\{[a,+\infty \lbrack, \ \  a\in \overline{\mathbb{R}}\}
\end{equation*}

\bigskip 
\begin{equation*}
\mathcal{I}(7)=\{[a,+\infty \lbrack ,a\in \mathbb{R}\}\text{ , }\mathcal{I}%
(8)=\{]a,+\infty \lbrack, \ \ a\in \overline{\mathbb{R}}\};
\end{equation*}

\bigskip 
\begin{equation*}
\mathcal{I}(9)=\{]a,b],a,b\in \mathbb{R}\}\ ,\ \mathcal{I}%
(10)=\{]a,b], \ \ (a,b)\in \overline{\mathbb{R}}^2\};
\end{equation*}

\bigskip 
\begin{equation*}
I(11)=\{]a,b[,a,b\in \mathbb{R}\}\ ,\ I(12)=\{]a,b[, (a,b)\in \overline{\mathbb{R}}^2\};
\end{equation*}

\bigskip 
\begin{equation*}
\mathcal{I}(13)=\{[a,b[,a,b\in \mathbb{R}\}\text{ , }\mathcal{I}%
(14)=\{[a,b[, (a,b)\in \overline{\mathbb{R}}^2\};
\end{equation*}

\bigskip 
\begin{equation*}
\mathcal{I}(15)=\{[a,b],a,b\in \mathbb{R}\}\text{ , }\mathcal{I}%
(16)=\{[a,b], (a,b) \in \overline{\mathbb{R}}^2\}.
\end{equation*}

\bigskip \noindent \textbf{Important Fact}. Any of these classes generates $\mathcal{B}(\mathbb{R})$. \newline

\bigskip \noindent \textbf{(03.03)} \textit{Borel application}. A real-valued and measurable application is called a Borel function. By using the criterion of measurability for the different generation classes in Point \textit{03.02}, we have a number of propositions, each of them characterizing measurability of functions

\begin{equation*}
X:(\Omega ,\mathcal{A})\rightarrow (\mathbb{R},\mathbb{B}(\mathbb{R})).
\end{equation*}

\bigskip \noindent Before we give some of these criteria, we introduce the notation :

\begin{equation*}
\text{For }B\subset \overline{\mathbb{R}},X^{-1}(B)=\{\omega \in
\Omega,X(\omega )\in B\}\equiv (X\in B).
\end{equation*}

\bigskip \noindent For example, for $a$, $b$ $\in \mathbb{R}$,

\begin{equation*}
(X<a)=X^{-1}([-\infty ,a[), \ \ (X\leq a)=X^{-1}([-\infty ,a]),
\end{equation*}

\begin{equation*}
(a<X\leq b)=X^{-1}(]a,b]), \ \ (a<X < b)=X^{-1}(]a,b[), \ \ etc.
\end{equation*}

\bigskip \noindent So, the real-valued application $X$ is measurable if and only if

\begin{equation*}
(C1):\text{For all }a\in \mathbb{R},(X<a)\in A,
\end{equation*}

\bigskip \noindent if and only if

\begin{equation*}
(C1):\text{For all }a\in \mathbb{R},(X\leq a)\in A,
\end{equation*}

\bigskip \noindent if and only if

\begin{equation*}
(C1):\text{For all }a\in \mathbb{R},(X>a)\in A,
\end{equation*}

\bigskip \noindent if and only if

\begin{equation*}
(C1):\text{For all }a\in \mathbb{R},(X\geq a)\in A,
\end{equation*}

\bigskip \noindent etc. (See Solutions of \textit{Exercise 2 in DOC 03-04} below, \pageref{exercise02_doc03-04}, in Doc 03-07, page \pageref{exercise02_sol_doc03-07}, for a more detailed
list of criteria).\newline

\bigskip \noindent \textbf{Part II - Working in the extended real line }$%
\overline{\mathbf{\mathbb{R}}}$.\newline

\noindent \textbf{Part II (a) - The space} $\overline{\mathbf{\mathbb{R}}}$.%
\newline

\noindent \textbf{(03.04)} We usually work on the extended real line 
\begin{equation*}
\overline{\mathbb{R}}=[-\infty ,+\infty ]=\mathbb{R\cup \{-\infty },\mathbb{%
+\infty \}}.
\end{equation*}

\bigskip \noindent The extended real line is obtained by adjoining $-\infty$ and +$\infty $ to $\mathbb{R}$. To define the sigma-algebra on $\overline{\mathbb{R}}$, we consider the elements of $\mathcal{B}_{\infty }(\mathbb{R})$ and their unions with $\{+\infty \}$, $\{-\infty \}$ and with $\{-\infty,+\infty \},$ that is

\bigskip

\begin{equation*}
\mathcal{B}_{\infty }(\overline{\mathbb{R}}\mathbb{)}=\{A,A\cup \{+\infty
\},A\cup \{+\infty \},A\cup \{-\infty ,+\infty \}\}.
\end{equation*}

\bigskip \noindent We have these remarkable points.\newline

\noindent \textbf{(03.05)} $\mathcal{B}_{\infty }(\overline{\mathbb{R}})$ this is smallest sigma-algebra for which the Borel sets of \ $\mathbb{R}$, $\{-\infty \},\{+\infty \},$ and $\{-\infty ,+\infty\}$ are measurable. You will not have any problem to show that each of the following classes of subsets also generates $\mathcal{B}_{\infty }(\overline{\mathbb{R}})$.\newline

\noindent \textbf{(03.06)} - Generation of $\mathcal{B}_{\infty }(\overline{\mathbb{R}}\mathbb{)}$. The same classes $\mathcal{I}(1),...,\mathcal{I}(14)$ given above generate $\mathcal{B}_{\infty }(\overline{\mathbb{R}}\mathbb{)}$. \newline

\bigskip \noindent Part II (b) \textbf{General Borel Functions}.\newline

\bigskip \noindent By using \textit{Point (03.06)} above, we conclude that :\newline

\noindent \textbf{(03.07)} For any real-valued function, possibly taking infinite values,

\begin{equation*}
X:(\Omega ,\mathcal{A})\rightarrow (\overline{\mathbb{R}},\mathcal{B}%
_{\infty }(\overline{\mathbb{R}}))
\end{equation*}

\bigskip \noindent is measurable if and only if (C1) holds, that is,

\begin{equation*}
\text{For all }a\in \mathbb{R},(X<a)\in \mathcal{A},
\end{equation*}

\bigskip \noindent if and only if (C2) holds, that is,

\begin{equation*}
\text{For all }a\in \mathbb{R},(X\leq a)\in \mathcal{A},
\end{equation*}

\bigskip \noindent if and only if (C3) holds, that is,

\begin{equation*}
\text{For all }a\in \mathbb{R},(X>a)\in \mathcal{A},
\end{equation*}

\bigskip \noindent if and only if (C4) holds, that is,

\begin{equation*}
\text{For all }a\in \mathbb{R},(X\geq a)\in \mathcal{A},
\end{equation*}

\bigskip \noindent etc.\newline

\bigskip \noindent \textbf{Part III}. Application to measurability of envelops and limits of sequences of real-valued applications.\newline

\bigskip \noindent Let $\{X_{n}\}_{n\geq 0}$ be a sequence of measurable real-valued applications defined on the same measurable space : 
\begin{equation*}
X_{n}:(\Omega ,\mathcal{A})\rightarrow (\overline{\mathbb{R}},\mathcal{B}%
_{\infty }(\overline{\mathbb{R}}))
\end{equation*}

\bigskip \noindent \textbf{(03.08)} Mesurability of the supremum : We use (C2) to prove that $\sup_{n\geq 0}X_{n}$ is measurable because for all $a\in \mathbb{R}$ 
\begin{equation*}
(\sup_{n\geq 0}X_{n}\leq a)=\bigcap\limits_{n\geq 0}(X_{n}\leq a).
\end{equation*}

\bigskip \noindent \textbf{(03.09)} Measurability of the infinum : We use (C1) to prove that $\inf_{n\geq 0}X_{n}$ because for all $a\in \mathbb{R}$

\begin{equation*}
(\inf_{n\geq 0}X_{n}<a)=\bigcup\limits_{n\geq 0}(X_{n}<a).
\end{equation*}

\bigskip \noindent \textbf{(03.09)} Measurability of the limit inferior : We
have 
\begin{equation*}
\underline{\lim }X_{n}=\inf_{n\geq 0}\sup_{m\geq n}X_{n}:=\inf_{n\geq 0}Y_{n},
\end{equation*}

\bigskip \noindent where $Y_{n}=\sup_{m\geq n}X_{n}$. Then $Y_{n}$ is measurable as the supremum of a sequence measurable functions and next $\underline{\lim }X_{n}$ is measurable as the infinum of a sequence of measurable functions.\newline

\bigskip \noindent Measurability of the limit superior : We have 
\begin{equation*}
\overline{\lim }X_{n}=\sup_{n\geq 0}\inf_{m\geq n}X_{n}:=\sup_{n\geq 0}Z_{n},
\end{equation*}

\bigskip \noindent where $Z_{n}=\inf_{m\geq n}X_{n}$. Then $Z_{n}$ is measurable as the infimum of a sequence measurable functions and $\overline{\lim }X_{n}$ next is measurable as the infimum of a sequence of measurable functions.\newline

\noindent \textbf{Remark}. The point above includes the finite case, that is $\max(X,Y)$ and $\min(X,Y)$ are finite whenever $X$ and $Y$ are.\\ 

\bigskip \noindent \textbf{(03.10)} \textbf{Measurability of the limit}. If the sequence $X_{n}$ admits a limit $X$ for every $\omega $, this limit is equal to the superior and the inferior limit, and so, it is measurable.\newline

\noindent \textbf{Rule : Limits of measurable applications are measurable}.%
\newline

\bigskip \noindent \textbf{Part IV}. Examples of measurable sets of $\overline{\mathbb{R}}$.\newline

\noindent It is not easy to find non-measurable sets on by simple means. Examples of non measurable sets come after very complicated mathematics.\newline

\noindent Any singleton $\{a\}$, $a\in \overline{\mathbb{R}}$, is measurable since $\{a\}=[a,a]$ is an interval.\newline

\noindent Finite subsets of $\overline{\mathbb{R}}$ are measurable as finite unions of singletons.\newline

\noindent Countable subsets $\overline{\mathbb{R}}$ are measurable as countable unions of singletons.\newline

\noindent The set of rationals $\mathbb{Q}$ is measurable since it is countable. The set of irrationals is also measurable as the complement of $Q$.\newline

\noindent etc.\newline

\bigskip \noindent\textbf{Part V} The spaces $\mathbb{R}^{k}$ and $\overline{%
\mathbb{R}}^{k}$, $k\geq 1$.\newline

\noindent \textbf{(a) Case of \ $\mathbb{R}^{k}$}.\newline

\noindent \textbf{(03.11)} The Borel $\sigma $-algebra on $\mathbb{R}^{k}$ is generated by the class $\mathcal{I}_{k}$ of intervals in $\mathbb{R}^{k}$. An interval in $\mathbb{R}^{k}$ is any set defined by 
\begin{equation*}
(a,b)=\prod\limits_{j=1}^{k}(a_{i},b_{i}),
\end{equation*}

\bigskip \noindent where $a=(a_{1},...,a_{k})\leq b=(b_{1},...,b_{k})$ [meaning that $a_{i}\leq b_{i},1\leq i\leq k]$ . As for the special case $k=1,$ we may consider all types of interval. But we restrict ourselves to two specific classes that are frequently used :

\begin{equation*}
\mathcal{I}_{k}(1)=\{]a,b]=\prod\limits_{j=1}^{k}]a_{i},b_{i}],\ \ (a,b)\in \mathbb{R}^{k},a\leq b\}
\end{equation*}

\bigskip \noindent and

\begin{equation*}
\mathcal{I}_{k}(2)=\{]-\infty ,a]=\prod\limits_{j=1}^{k}]-\infty,a_{i}], \ \ a\in \overline{\mathbb{R}}^{k}\}.
\end{equation*}

\bigskip \noindent We have the following facts.\newline

\noindent \textbf{(03.12)} The classes $\mathcal{I}_{k},$ $\mathcal{I}_{k}(1)$ and $\mathcal{I}_{k}(2)$ generate the same $\sigma $-algebra. The $\sigma$-algebra generated by each of them is the usual $\sigma $-algebra on $\mathbb{R}^{k}$ denoted $\mathcal{B}$($\mathbb{R}^{k})$.\newline

\noindent \textbf{(03.13)} $\mathcal{I}_{k}(2,c)=\{\overline{\mathbb{R}}^{k}\}\cup \mathcal{I}_{k}(2,c)$ is a semi-algebra.\newline

\noindent \textbf{(03.14)} The usual $\sigma $-algebra on $\mathbb{R}^{k}$, denoted $\mathcal{B}$($\mathbb{R}^{k}),$ is actually the product $\sigma$-sigma algebra of $(\mathbb{R},\mathcal{B}(\mathbb{R}))$ by itself $k$ times : $(\mathbb{R}^{k},\mathcal{B}(\mathbb{R}^{k}))=(\mathbb{R}^{k},\mathcal{B}(\mathbb{R})^{\otimes k})$.\newline

\noindent \textbf{(03.15)} The usual $\mathcal{B}(\mathbb{R}^{k})$ is also the Borel $\sigma $-algebra on $\mathbb{R}^{k}$ endowed with one of the three equivalent metric defined for : $a=(a_{1},...,a_{k})$ and $b=(b_{1},...,b_{k})$

\begin{equation*}
d_{e}(a,b)=\sqrt{\sum_{j=1}^{k}(a_{i}-b_{i})^{2}}\text{,}
\end{equation*}

\begin{equation*}
d_{m}(a,b)=\sum_{j=1}^{k}\left\vert a_{i}-b_{i}\right\vert \text{,}
\end{equation*}

\bigskip \noindent and

\begin{equation*}
d_{\infty }(a,b)=\max_{1\leq i\leq k}\left\vert a_{i}-b_{i}\right\vert.
\end{equation*}

\bigskip \noindent These three metrics are respectively called the Euclidean, Manhattan
and max (or uniform) metric.\newline

\bigskip \noindent \textbf{(b) Case of \ $\overline{\mathbb{R}}^{k}$}.\newline

\noindent We endow $\overline{\mathbb{R}}^{k}$ with the product $\sigma $-algebra $\mathcal{B}_{\infty }(\overline{\mathbb{R}})$ $k$ times and get the measurable space ($\overline{\mathbb{R}}^{k},\mathcal{B}_{\infty }(\overline{\mathbb{R}})^{\otimes k}).$ We denote $\mathcal{B}_{\infty }(\overline{\mathbb{R}}^{k})=\mathcal{B}_{\infty }(\overline{\mathbb{R}})^{\otimes k}$.\newline

\noindent We have the following fact.\newline

\noindent \textbf{(03.16)} On ($\overline{\mathbb{R}}^{k},\mathcal{B}(\overline{\mathbb{R}})^{\otimes k}),$ \ $\mathcal{B}(\mathbb{R}^{k})\subset \mathcal{B}(\overline{\mathbb{R}})^{\otimes k}$.\newline

\noindent \textbf{(c) Measurability of functions} $X=(X_{1},...,X_{k})$ from a measurable space $(\Omega ,\mathcal{A})$ to $\overline{\mathbb{R}}^{k}$.\newline

\noindent \textbf{(03.17)} Such a function considered as a function with values in $\mathbb{R}^{k}$ endowed that $\mathcal{B}(\mathbb{R}^{k})$\ or in \ $\overline{\mathbb{R}}^{k}$endowed with $\mathcal{B}(\overline{\mathbb{R}})^{\otimes k}$ is measurable if and only if each component $X_{i}$, $1\leq i\leq k,$ is measurable as an application from $(\Omega ,\mathcal{A})$ to $(\mathbb{R},\mathcal{B}(\mathbb{R})$ or to $(\overline{\mathbb{R}},\mathcal{B}_{\infty }(\overline{\mathbb{R}})),$ is measurable.\\

\bigskip \noindent \textbf{(d) Examples of Measurable sets in $\overline{\mathbb{R}}^{k}$}.\\

As for the special class $k=1,$ non-measurable sets are not easy to get in $\overline{\mathbb{R}}^{k}$. Almost all the classical subsets we know, are measurable. This is not a reason to get for granted that every subset of $\overline{\mathbb{R}}^{k}$ is measurable in $\overline{\mathbb{R}}^{k}$. The following list of measurable types of subsets of $\overline{\mathbb{R}}^{k}$ is given only for illustration.\\

\noindent \textbf{(d1)} Intervals in $\overline{\mathbb{R}}^{k}$ are measurable by construction.\\

\noindent \textbf{(d2)} Singletons in $\overline{\mathbb{R}}^{k}$ are measurable as since they are intervals
\begin{equation*}
(a_{1},...,a_{k})=[a_{1},a_{1}]\times ...\times \lbrack a_{k},a_{k}].
\end{equation*}

\bigskip \noindent \textbf{(d3)} Countable sets are measurable.\\

\noindent \textbf{(d4)} $\mathbb{Q}^{k}$ is measurable as a product of measurable sets in $\mathbb{Q}$. $\mathbb{J}^{k}$ is measurable for the same reason ($\mathbb{Q}$ is the set of rational numbers, $\mathbb{J}$\ is the set if irrational numbers in $\mathbb{R}$).\\

\noindent \textbf{(d5)} Hyperplanes are measurable. Hyperplanes are linear subspaces of $\mathbb{R}^{k}$\ of dimension $k-1$. Any hyperplane $H$ may be defined by an equation like : there exist real numbers $a_{1},a_{2},..,a_{k}$ such that $a_{1}\times a_{2}\times ...\times a_{k}\neq 0$ and a real number $b$ such that  

\begin{equation*}
H=\{(x_{1},....,x_{k})\in R^{k},a_{1}x_{1}+....+a_{k}x_{k}+b=0\}.
\end{equation*}

\noindent Since the application $h(x_{1},....,x_{k})=a_{1}x_{1}+....+a_{k}x_{k}+b$ is continuous from $\mathbb{R}^{k}$ to $\mathbb{R}$ is continuous and the measurable. And we have $H=h^{-1}(\{0\}).$ Hence $H$ is measurable as the inverse image of the measurable set  $\{0\}$ in $\mathbb{R}$.\\

\noindent \textbf{(d6)} Affine subspaces in $\mathbb{R}^{k}$ are measurable.\\

\noindent To see this, proceed as follows. An affine subspace of $\mathbb{R}^{k}$\ is of the form $a+H=\{a+x,x\in H\}$, where $b\in \mathbb{R}^{k}$ and $H$ is a linear subspace of $\mathbb{R}^{k}.$
By $b+H$ is the inverse image of the translation $t_{-b}(x)=x-b,$ which is continuous and then measurable. It is then clear that $b+H$ is measurable whenever $H$ is. But the only linear subspace of $\mathbb{R}^{k}$ of dimension $k$ is $\mathbb{R}^{k},$ the only linear subspace of $\mathbb{R}^{k}$ of dimension $0$ is $\{0\}$ and finally, any linear subspace of $%
\mathbb{R}^{k}$ of dimension $s$ with $0<s<k$ is intersection of $k-s+1$ hyperplanes. Conclusion : any affine subspace of $\mathbb{R}^{k}$ is measurable.\\

\noindent \textbf{(d7)} Rectangles, cubes, hyper-cubes, squares, etc. are measurable as intervals.\\

\noindent \textbf{(d8)} Parallelograms, parallelepiped, cylinders, etc. are measurable.\\ 

\noindent \textbf{(d8.1)} We may use the equations of all these geometric objects and use their continuity of the functions involved (in the equations) to express these objects as inverse images of measurable sets.\\

\noindent \textbf{(d8.2)} We may also the topological concepts of closure, boundary and interior to establish measurability. You will learn more about this with the course of topology.\\

\noindent \textbf{(d6)} etc..\\

\bigskip \noindent \textbf{(e) Non-measurable sets}.\\

\noindent In \textit{Chapter \ref{04_measures}} (See \textit{Exercise 4 in Doc 04-04}), page \pageref{exercise_04_doc-04-04}), the Vitaly's construction of a non-measurable
set will be given.

\newpage
\noindent \LARGE \textbf{DOC 03-02 : Measurable sets and measurable functions in special spaces - Elementary functions}. \label{doc03-02}\\
\Large

\bigskip \noindent \textbf{03.18} : Definition of an elementary function.\\

\noindent Consider a measurable and finite partition of $\Omega $, that is, $k$ mutually disjoint and measurable subsets $A_{1}$,..., $A_{k}$ such that
\begin{equation*}
\Omega =\sum_{1\leq i\leq k}A_{i}.
\end{equation*}

\bigskip \noindent Let $X$ be an application taking constant and finite values $\alpha _{i}\in \mathbb{R}$ on each $A_{i}$, $1\leq i \leq k$, that is 
\begin{equation*}
X(\omega )=\alpha _{i}\text{ for }\omega \in A_{i} \text{ } 1\leq i \leq k,
\end{equation*}

\bigskip \noindent also denoted by  
\begin{equation}
X=\sum_{1\leq i\leq k}\alpha _{i}1_{A_{i}}.  \label{etag01}
\end{equation}

\noindent Such an application (\ref{etag01}) is called an \textit{elementary} or \textit{simple} function.\\

\bigskip \noindent \textbf{(03.19)} Be careful : A simple function may have several expressions. For example, we get another expression if we split one of the $A_i$.\\

\bigskip \noindent \textbf{(03.20)} There exists a unique expression for an elementary function of the form (\ref{etag01}) with distinct values of $\alpha_i$, $1 \leq i \leq k$.\\ 

\bigskip \noindent \textbf{(00.21)} The class $\mathcal{E}$ of simple functions is a Riesz space, and is an algebra of functions.\\

\noindent Denote by $(\mathcal{E}$ the class of all simple functions. When $(\mathcal{E} ,+,\cdot)$ is a vectorial space. The order $\leq$ is compatible with the vectorial structure (defined by Point (c) below) and $\mathcal{E}$ contains the maximum and the minimum of its couple. We say : $(\mathcal{E} ,+ , \cdot,\leq )$ is a lattice vectorial space or a Riesz space. Also $(\mathcal{E} ,+, \times, \cdot)$ is an algebra of functions (in the algebraic sense). The details of all this are :

\bigskip \noindent (a) Let $X$, $Y$ and $Z$ simple functions, $a$ and $b$ real scalars. We have

\bigskip \noindent (b) $aX+bY \in \mathcal{E}$, $XY \in \mathcal{E}$, $max(A,X) \in \mathcal{E}$, $min(X,Y) \in \mathcal{E}$.\\

\bigskip \noindent (c) If $Y$ is everywhere different of zero, then $X/Y \in \mathcal{E}$.\\

\bigskip \noindent (d) $(X \leq Y) \Rightarrow (X+Z \leq Y+Z)$. If $a\geq 0$, then $X \leq Y \Rightarrow aX \leq bY$.\\

\bigskip \noindent Tools for demonstrations : \\

\bigskip \noindent \textbf{(03.22a)} Let $X=\sum_{1\leq i\leq k}\alpha _{i}1_{A_{i}}$ and $X=\sum_{1\leq j\leq m}\beta _{j}1_{B_{j}}$ two writings of $X$. Then we also have :\\
$$
X=\sum_{(i,j)\in I}\gamma_{ij} \ 1_{A_{i}B_{j}},
$$

\bigskip \noindent where $\gamma_{ij}=\alpha_{i}=\beta_{j}$ and where $I=\{(i,j), 1\leq i \leq k, 1\leq j \leq m, A_{i}B_{j} \neq \emptyset \}$.\\

\bigskip \noindent Recall that we have

$$
\Omega = \sum_{1\leq i\leq k} \sum_{1\leq j\leq m}1_{A_{j}}1_{B_{j}}.
$$

\bigskip \noindent  The partition $\{ A_{j}B_{j}, \text{  } 1\leq i\leq k, \text{  } 1\leq j\leq m \}$ is called superposition of the partitions $\{ A_{i}, \text{  } 1\leq i\leq k \}$ and 
$\{ B_{j}, \text{  } 1\leq j\leq m \}$.\\

\noindent \textbf{(03.22b)} As well, if we have two simple functions $X=\sum_{1\leq i\leq k}\alpha _{i}1_{A_{i}}$ and $Y=\sum_{1\leq j\leq m}\beta _{j}1_{B_{j}}$, we may write them with a common partition :

$$
X=\sum_{(i,j)\in I}\alpha_{i}1_{A_{i}B_{j}}
$$

\bigskip \noindent and

$$
X=\sum_{(i,j)\in I}\beta_{j}1_{A_{i}B_{j}}.
$$

\bigskip \noindent \textbf{(03.23)} The most important key of Measure theory, basis of the famous three steps methods.\\

\bigskip \noindent \textbf{Any positive measurable application is limit of an non-decreasing sequences of positive simple functions}.\\

\bigskip \noindent To prove this consider the following construction. Let $X\geq 0$ be measurable. Pour $n$ fix\'{e}, set the following partition of $\mathbb{R}_{+}$ :

$$
\mathbb{R}_{+}=\sum_{k=1}^{2^{2n}}[\frac{k-1}{2^{n}},\frac{k}{2^{n}}[\text{\ }+\text{ }[2^{n},+\infty ]=\sum_{k=1}^{2^{2n}+1}A_{k}.
$$

\bigskip \noindent Now set the sequence of elementary functions, $n\geq 1$,  
$$
X_{n}=\sum_{k=1}^{2^{2n}}\frac{k-1}{2^{n}}\text{ }1_{(\frac{k-1}{2^{2n}}\leq
X<\frac{k}{2^{n}})}+2^{n}\text{ }1_{(X\geq 2^{2n})}. 
$$

\bigskip \noindent The sequence is $X_{n}$ is non-decreasing and $X_{n}$ tends to $X$ as $n\rightarrow +\infty$.\\

\newpage
\noindent \LARGE \textbf{DOC 03-03 : Measurable applications - A summary}. \label{doc03-03}\\
\Large

\bigskip \noindent \textbf{I - General introduction}.\newline

\noindent \textbf{1 : Definition}. An application $X$ from the measurable space $(\Omega _{1},\mathcal{A}_{1})$ onto the measurable space $(\Omega_{2},\mathcal{A}_{2})$, is measurable, with respect to $(\mathcal{A}_{1}), \mathcal{A}_{2}))$, if

\begin{equation}
\forall (B\in \mathcal{A}_{2}),\text{ } X^{-1}(B)\in \mathcal{A}_{1}, 
\tag{M}
\end{equation}

\noindent It is usual in measure theory to use the following notation : 
\begin{equation*}
X^{-1}(B)=(X\in B).
\end{equation*}

\bigskip \noindent \textbf{2 : Criterion}. Fortunately, in many cases, $\mathcal{A}_{2}$ is generated by a small class, say $\mathcal{C}$, that is $\mathcal{A}_{2}=\sigma (\mathcal{C})$. In these cases, the application $X$ is measurable if 
\begin{equation}
\forall (B\in \mathcal{C}),\text{ }X^{-1}(B)\in \mathcal{A}_{1}.  \tag{M}
\end{equation}

\noindent \textbf{3 : Composition of measurable applications}. Let $X$ be a measurable application from $(\Omega _{1},\mathcal{A}_{1})$ onto $(\Omega _{2},\mathcal{A}_{2})$, and let $Y$ be a measurable application form $(\Omega _{2},\mathcal{A}_{2})$ onto $(\Omega _{3},\mathcal{A}_{3})$. Then the composite function (function of function) $Z=Y(X)=YX$ defined from $%
(\Omega _{1},\mathcal{A}_{1})$ to $(\Omega _{3},\mathcal{A}_{3})$ by 

\begin{equation*}
\Omega _{1}\ni \omega \Rightarrow Z(\omega )=Y(X(\omega ))
\end{equation*}

\bigskip \noindent is measurable.\newline

\noindent \textbf{4 : Criteria of measurability for an application taking a countable number of values}.\newline

\noindent Let $X :(\Omega _{1},\mathcal{A}_{1}) \mapsto (\Omega _{2}, \mathcal{A}_{2})$ be an applications taking a countable number of distinct values $x_{i}$ $(i\in I)$ in $\Omega _{2}$, with $I$ countable. The we have

\begin{equation*}
\Omega_1 =\sum\limits_{i\in I}(X=x_{i}),
\end{equation*}

\bigskip \noindent and for any subset $B \in \Omega _{2}$, we have

\begin{equation}
X^{-1}(B)=\sum\limits_{i\in I,x_{i}\in B}(X=x_{i}).  \label{mesCountable}
\end{equation}

\noindent It follows that $X$ is measurable if and only if 
\begin{equation*}
\forall (i\in I)\text{, }(X=\alpha _{i})\text{ is measurable in } \Omega_{1}.
\end{equation*}

\bigskip \noindent \textbf{II - Operations on real-valued and measurable functions using the class of simple functions}.\newline

\noindent \textbf{1 - Indicator function of a measurable set}. For any measurable space A in the measurable space $(\Omega ,\mathcal{A}),$ the indication function
\begin{equation*}
\begin{tabular}{lllll}
$1_{A}$ & : & $(\Omega ,\mathcal{A})$ & $\longrightarrow $ & $(\mathbb{R},\mathcal{B}(\mathbb{R}))$ \\ 
&  & $\omega $ & $\hookrightarrow $ & $1_{A}(\omega )=\left\{ 
\begin{tabular}{lll}
$1$ & $if$ & $\omega \in A$ \\ 
$0$ & $if$ & $\omega \in $
\end{tabular}
\right. $.
\end{tabular}
\end{equation*}

\bigskip \noindent is measurable.\newline

\bigskip \noindent \textbf{2 : Elementary functions}. A function 
\begin{equation*}
\begin{tabular}{lllll}
$f$ & : & $(\Omega ,\mathcal{A})$ & $\longrightarrow $ & $(\mathbb{R},%
\mathcal{B}(\mathbb{R}))$,
\end{tabular}
\end{equation*}

\bigskip \noindent is an elementary function if and only if it is finite and 
\newline

\bigskip \noindent (i) there exists a finite and measurable subdivision of $\Omega$, that is there exist $p\geq 1$ mutually disjoint and measurable sets $A_{1}$, ..., $A_{p}$, such that 
\begin{equation*}
\Omega =A_{1}+...+A_{p},
\end{equation*}

\bigskip \noindent and\\

\bigskip \noindent (ii) there exist $p$ real and \textbf{finite} numbers $\alpha _{1},...,\alpha _{p}$ such that

\bigskip 
\begin{equation}
\forall (1\leq i\leq p),\text{ }\forall \omega \in A_{i},\text{ }f(\omega
)=\alpha _{i}.  \tag{EF1}
\end{equation}

\bigskip \noindent Formula (EF1) is equivalent to
\begin{equation*}
f=\alpha _{i}\text{ on }A_{i},\text{ }1\leq i\leq p
\end{equation*}

\bigskip \noindent and to
\begin{equation*}
f=\sum_{i=1}^{p}\alpha _{i}1_{A_{i}.}
\end{equation*}

\bigskip \noindent \textbf{Properties}. We have the following properties

\bigskip \noindent \textbf{2.1 : Operations in $\mathcal{E}$, the class of all elementary functions}.\newline

\bigskip \noindent If $f$ and $g$ are elementary functions, and $c$ is a real number, then $f+g, $ $fg,cf$, $\max (f,g)$ and $\min (f,g)$ are elementary functions. If $g$ does not take the null value, then $f/g$ is an elementary function.

\bigskip \noindent \textbf{2.2 : Approximation of a non-negative measurable
function}.\newline

\bigskip \noindent For any non-negative and measurable function 
\begin{equation*}
\begin{tabular}{lllll}
$f$ & : & $(\Omega ,\mathcal{A})$ & $\longrightarrow $ & $(\overline{\mathbb{R}},\mathcal{B}(\overline{\mathbb{R}}))$.
\end{tabular}
\end{equation*}

\bigskip \noindent there exists a non-decreasing sequence of non-negative elementary functions $(f_{n})_{n\geq 0}$ such that
\begin{equation*}
f_{n}\nearrow f\text{ as }n\nearrow +\infty .
\end{equation*}

\bigskip \noindent \textbf{2.3 : Approximation of an arbitrary measurable function}.\newline

\bigskip \noindent Let 
\begin{equation*}
\begin{tabular}{lllll}
$f$ & : & $(\Omega ,\mathcal{A})$ & $\longrightarrow $ & $(\overline{\mathbb{%
R}},\mathcal{B}(\overline{\mathbb{R}}))$.
\end{tabular}
\end{equation*}

\bigskip \noindent be measurable, we have the decomposition%
\begin{equation*}
f=f^{+}-f^{-},
\end{equation*}

\bigskip \noindent where 
\begin{equation*}
f^{+}=\max (f,0)\text{ \ and }f^{-}=\max (-f,0).
\end{equation*}

\bigskip \noindent \textbf{Fact } : $f^{+}$ and $f^{-}$ are positive and measurable. So there exist a non-decreasing sequences of non-negative elementary functions $(f_{n}^{+})_{n\geq 0}$ and $(f_{n}^{-})_{n\geq 0}$ such that
\begin{equation*}
f_{n}^{+}\nearrow f^{+}\text{ and }f_{n}^{-}\nearrow f \text{ as } n\nearrow +\infty .
\end{equation*}

\bigskip \noindent Then, there exists a sequence of elementary functions $(f_{n})_{n\geq 0}$ \ \ $[f_{n}=f_{n}^{+}-f_{n}^{-}]$ \ such that%
\begin{equation*}
f_{n}\longrightarrow f^{-}\text{ as }n\longrightarrow +\infty .
\end{equation*}

\bigskip \noindent \textbf{2.4 : Operations on real-valued measurable functions}.\newline

\bigskip \noindent If $f$ and $g$ are measurable functions, and $c$ is a real number, then $f+g,$ $fg,cf$, $\max (f,g)$ and $\min (f,g)$ are measurable whenever theses \ functions are well-defined. If $g$ does not take the null value, then $f/g$ is measurable.\newline

\bigskip

\bigskip \noindent \textbf{III : Criteria of measurability for real-valued
functions}.\newline

\bigskip \noindent A real-valued function 
\begin{equation*}
\begin{tabular}{lllll}
$f$ & : & $(\Omega ,\mathcal{A})$ & $\longrightarrow $ & $(\overline{\mathbb{R}},\mathcal{B}(\overline{\mathbb{R}}))$.
\end{tabular}
\end{equation*}

\bigskip \noindent is measurable if and only if

\begin{equation*}
\forall a\in \mathbb{R},\text{ }(f>a)\in \mathcal{A}
\end{equation*}

\bigskip \noindent if and only if

\begin{equation*}
\forall a\in \mathbb{R},\text{ }(f\geq a)\in \mathcal{A}
\end{equation*}

\bigskip \noindent if and only if

\begin{equation*}
\forall a\in \mathbb{R},\text{ }(f<a)\in \mathcal{A}
\end{equation*}

\bigskip \noindent if and only if

\begin{equation*}
\forall a\in \mathbb{R},\text{ }(f\leq a)\in \mathcal{A}
\end{equation*}

\bigskip \noindent if and only if

\begin{equation*}
\forall (a,b)\in \mathbb{R}^{2},\text{ }(a<f<a)\in \mathcal{A}.
\end{equation*}

\bigskip \noindent And so forth.\newline

\bigskip

\bigskip \noindent \textbf{IV - Examples of real-valued measurable functions from $\mathbb{R}$ to $\mathbb{R}$}.\newline

\noindent \textbf{(a)} Right of left continuous functions, including continuous functions.\newline

\noindent \textbf{(b)} Function with at most a countable of discontinuity points having left limits at each point or having right limit
at each point.\newline

\noindent \textbf{c)} Monotonic functions.\newline

\noindent \textbf{d)} Convex functions (convex functions are
continuous).\newline

\bigskip \noindent \textbf{6 : Examples measurable functions defined from a topological space (E,T) to $\mathbb{R}$}.\newline

\noindent \textbf{(a)} Continuous functions.\newline

\noindent  If $(E,T)$ is a metric space $(E,d)$.\newline

\noindent \textbf{(b)} lower semi-continuous functions.\newline

\noindent \textbf{(c)} upper semi-continuous functions.\newline

\noindent \textbf{(c)} (d) inferior continuous function and superior continuous functions.\\

\noindent \textbf{When a product space is involved}.\newline

\noindent \textbf{(a) : Multi-component function}.\newline

\noindent Consider an application with values in a product space endowed with his product $\sigma $-algebra. \ 
\begin{equation*}
f:(E,\mathcal{B})\mapsto \left( \prod_{1\leq i\leq k}\Omega_{i},\bigotimes_{1\leq i\leq k}\mathcal{A}_{i}\right) .
\end{equation*}

\bigskip \noindent The image $f(\omega )$ has $k$ components. 
\begin{equation*}
f(\omega )=(f_{1}(\omega ),f_{2}(\omega ),...,f_{k}(\omega )),
\end{equation*}

\bigskip \noindent and each component $f_{i}$ is an application from $(E,\mathcal{B})$ to $(\Omega _{i},\mathcal{A}_{i})$.\newline

\noindent \textbf{Criterion} : $f$ is measurable if and only if each $f_{i}$ is measurable.\newline

\noindent \textbf{(b) : Real-valued function of multiple arguments}.\newline

\noindent \textbf{(b1) Definition} \noindent Now let $f$ be a function defined on the product space with values in $\overline{\mathbb{R}}$ : \ 
\begin{equation*}
\begin{tabular}{lllll}
$f$ & $:$ & $\left( \prod_{1\leq i\leq k}\Omega _{i},\bigotimes_{1\leq i\leq k}\mathcal{A}_{i}\right) $ & $\longmapsto $ & $(\overline{\mathbb{R}},\mathcal{B}_{\infty }(\overline{\mathbb{R}}))$ \\ 
&  & $(\omega _{1,}\omega _{2,}...,\omega _{k})$ & $\hookrightarrow $ & $f(\omega _{1,}\omega _{2,}...,\omega _{k})$.
\end{tabular}
\end{equation*}

\bigskip \noindent The function $f$ has $k$ arguments or $k$ variables. For $i\in \{1,2,...\}$ fixed, we may consider $i-$th partial function $f$ by considering it as a function of $\omega _{i}$ only and by fixing the other variables. Let $(\omega _{1,}\omega _{2,}...,\omega _{i-1},\omega_{i+1},...,\omega _{k})$ fixed and define $f^{(i)}$ the as $i-$th partial function $f$ 
\begin{equation*}
\begin{tabular}{lllll}
$f^{(i)}$ & $:$ & $\left( \Omega _{i},\mathcal{A}_{i}\right) $ & $\longmapsto $ & $(\overline{\mathbb{R}},\mathcal{B}_{\infty }(\overline{%
\mathbb{R}}))$ \\ 
&  & $\mathbf{t}$ & $\hookrightarrow $ & $f^{(i)}(t)=f(\omega _{1,}\omega
_{2,}...,\omega _{i-1},\mathbf{t},\omega _{i+1},...,\omega _{k})$,
\end{tabular}
\end{equation*}

\bigskip \noindent which is defined of the space $\left( \Omega _{i},\mathcal{A}_{i}\right)$. All the arguments, except the $i-$th, are fixed. We have this result.\newline

\noindent \textbf{Rule} : For any real-valued and measurable function of several variables, its partial functions are measurable.\newline

\noindent \textbf{Remarks}.\newline

\noindent (1) Be careful. This rule is valid for real-valued function.\\

\noindent (2) The reverse implication not true, exactly as for the continuity of such functions.\newline

\noindent (3) The proof is based on Property (2.3) of this document and properties of sections of measurable sets (see Doc 01-06, Exercise 3).\\

\bigskip \noindent For $k=2$, have these usual notation : for $x$ fixed,
\begin{equation*}
y\hookrightarrow f_{x}(y)=f(x,y)
\end{equation*}

\bigskip \noindent and for $y$ fixed, 
\begin{equation*}
x\hookrightarrow f_{y}(y)=f(x,y).
\end{equation*}

\bigskip \noindent \textbf{(b2) Special case of elementary functions}.\newline

\noindent If $A$ is a subset of $\left( \prod_{1\leq i\leq k}\Omega_{i},\bigotimes_{1\leq i\leq k}\mathcal{A}_{i}\right) ,$ its indicatorfunction is function of the variables $\omega _{1,}\omega _{2,}...,\omega_{k}$ :

\begin{equation*}
(\omega _{1,}\omega _{2,}...,\omega _{k})\hookrightarrow 1_{A}(\omega_{1,}\omega _{2,}...,\omega _{k}).
\end{equation*}

\bigskip \noindent So the $i-$th partial function is defined by
\begin{equation*}
\omega _{i}\hookrightarrow 1_{A}(\omega _{1,}\omega _{2,}...,\omega_{k})=1_{A_{(\omega _{1,}\omega _{2,}...,\omega _{i-1},\omega_{i+1},...,\omega _{k})}}(\omega _{i}),
\end{equation*}

\bigskip \noindent where $A_{(\omega _{1,}\omega _{2,}...,\omega _{i-1},\omega_{i+1},...,\omega _{k})}$ is the section of $A$ at $A_{(\omega _{1,}\omega_{2,}...,\omega _{i-1},\omega _{i+1},...,\omega _{k})}$.\newline

\bigskip \noindent If $f$ is an elementary function 
\begin{equation*}
f=\sum_{j=1}^{p}\alpha _{j}1_{A_{j}},
\end{equation*}

\bigskip \noindent it comes that $i-$th partial function $f^{(i)}$ is%
\begin{equation*}
\omega _{i}\hookrightarrow f^{(i)}(\omega _{i})=\sum_{i=1}^{p}\alpha_{j}1_{(A_{j})_{(\omega _{1,}\omega _{2,}...,\omega _{i-1},\omega_{i+1},...,\omega _{k})}}(\omega _{i}),
\end{equation*}

\bigskip \noindent where the $(A_{j})_{j(\omega _{1,}\omega _{2,}...,\omega_{i-1},\omega _{i+1},...,\omega _{k})}$ are the section of the $A_{j}$ at $(\omega _{1,}\omega _{2,}...,\omega _{i-1},\omega _{i+1},...,\omega _{k}).$

\bigskip \noindent For $k=2$, you will see clearer. For an indicator function for example, for $\omega _{1}$ fixed (for example), we have the partial function 
\begin{equation*}
\omega _{2}\hookrightarrow 1_{A}(\omega _{1,}\omega )=1_{A_{\omega
_{1}}}(\omega _{2})n,
\end{equation*}

\bigskip \noindent where $_{A_{\omega _{1}}}$ is the section of $A$ at $\omega _{1}$. For the elementary function

\begin{equation*}
f=\sum_{j=1}^{p}\alpha _{j}1_{A_{j}},
\end{equation*}

\bigskip \noindent the partial function for $\omega _{1}$ is 
\begin{equation*}
\omega _{2}\hookrightarrow f^{(i)}(\omega _{2})=\sum_{i=1}^{p}\alpha
_{j}1_{(A_{1})_{(\omega _{1})}}(\omega _{2})
\end{equation*}

\noindent \textbf{Consequence}. From Doc 01-06, Exercise 3, we immediately see the partial elementary functions are measurable. The proof of this rule is given un Exercise 6 in Doc 03-07 below.

\newpage
\noindent \LARGE \textbf{DOC 03-04 : Exercises on Measurability 01 - Borel sets and Borel functions}. \label{doc03-04}\\

\Large 

\bigskip  \noindent  \textbf{Summary : } You will discover the direct properties of measurable real-valued functions through criteria, characterizations, limits and some operations.\\

\noindent \textbf{Exercise 1}. \label{exercise01_doc03-04}\\

\noindent \textbf{(a)} Show that each of these classes generates the same $\sigma $-algebra as the class of all intervals $\mathcal{I}$. Do only one or two cases.

\begin{equation*}
\mathcal{I}(1)=\{]-\infty ,a[,\ a\in \mathbb{R}\}\text{ , }\mathcal{I}(2)=\{]-\infty ,a],\ \ a\in \overline{\mathbb{R}}\}.
\end{equation*}

\bigskip 
\begin{equation*}
\mathcal{I}(3)=\{]-\infty ,a],a\in \mathbb{R}\}\text{ , }\mathcal{I}(4)=\{]-\infty ,a],\ \ a\in \overline{\mathbb{R}}\}
\end{equation*}

\bigskip 
\begin{equation*}
\mathcal{I}(5)=\{]a,+\infty \lbrack, \ \ a\in \mathbb{R}\},\mathcal{I}(6)=\{[a,+\infty \lbrack, \ \ a\in \overline{\mathbb{R}}\}
\end{equation*}

\bigskip
\begin{equation*}
\mathcal{I}(7)=\{[a,+\infty \lbrack,\ \ a\in \mathbb{R}\}\text{ , }\mathcal{I}(8)=\{]a,+\infty \lbrack, \ \ a\in \overline{\mathbb{R}}\}
\end{equation*}

\bigskip 
\begin{equation*}
\mathcal{I}(9)=\{]a,b],a,b\in \mathbb{R}\}\ ,\ \mathcal{I}(10)=\{]a,b],\ \ (a,b)\in \overline{\mathbb{R}}^2\}
\end{equation*}

\bigskip 
\begin{equation*}
I(11)=\{]a,b[,a,b\in \mathbb{R}\}\ ,\ I(12)=\{]a,b[,\ \ (a,b)\in \overline{\mathbb{R}}^2\}
\end{equation*}

\bigskip 
\begin{equation*}
\mathcal{I}(13)=\{[a,b[,\ (a,b) \in \mathbb{R}^2\}\text{ , }\mathcal{I}(14)=\{[a,b[,\ \ (a,b)\in \overline{\mathbb{R}}^2\}
\end{equation*}

\bigskip 
\begin{equation*}
\mathcal{I}(15)=\{[a,b],(a,b) \in \mathbb{R}^2\}\text{ , }\mathcal{I}(16)=\{[a,b],\ \ (a,b)\in \overline{\mathbb{R}}^2\}
\end{equation*}

\bigskip \noindent \textbf{(b)} show that the class $\mathcal{I}(4)=\{]-\infty ,a],\ \ a\in \overline{\mathbb{R}}\}$ is a semi-algebra in $\mathbb{R}$.\\

\bigskip \noindent \textbf{Exercise 2}. \label{exercise02_doc03-04} Consider the extended real line 
\begin{equation*}
\overline{\mathbb{R}}=[-\infty ,+\infty ]=\mathbb{R} \cup \{-\infty,+\infty \},
\end{equation*}

\noindent obtained by adding the points $-\infty $ and +$\infty $ to $\mathbb{R}$. Define 

\begin{equation*}
\mathcal{B}_{\infty }(\overline{\mathbb{R}}\mathbb{)}=\{A,A\cup \{-\infty\},A\cup \{+\infty \},A\cup \{-\infty ,+\infty \},\text{ }A\in \mathcal{B}(\mathbb{R)}\}.
\end{equation*}

\bigskip \noindent \textbf{(a)}. Show that $\mathcal{B}_{\infty }(\overline{\mathbb{R}}\mathbb{)}$\ is the smallest sigma-algebra for which the Borel sets of \ $\mathbb{R}$, $%
\{-\infty \},\{+\infty \}$ are measurable.\\

\noindent \textbf{(b)} Show that each of the classes given in Exercise 1 also generates $\mathcal{B}_{\infty }(\overline{\mathbb{R}}\mathbb{)}.$

\bigskip \noindent \textbf{(c)} Show that any function

\begin{equation*}
X:(\Omega ,\mathcal{A})\rightarrow (\mathbb{R},\mathcal{B}(\mathbb{R}))\text{or }X:(\Omega ,\mathcal{A})\rightarrow (\mathbb{R},\mathcal{B}_{\infty }(\overline{\mathbb{R}})).\text{ }
\end{equation*}

\bigskip \noindent is measurable if and only if

\begin{equation*}
(C1):\text{For all }a\in \mathbb{R},(X<a)\in \mathcal{A},
\end{equation*}

\bigskip \noindent if and only if

\begin{equation*}
(C2):\text{For all }a\in \mathbb{R},(X\leq a)\in \mathcal{A}\text{,}
\end{equation*}

\bigskip \noindent if and only if

\begin{equation*}
(C3):\text{For all }a\in \mathbb{R},(X>a)\in \mathcal{A},
\end{equation*}

\bigskip \noindent if and only if

\begin{equation*}
(C4):\text{For all }a\in \mathbb{R},(X\geq a)\in \mathcal{A},
\end{equation*}

\bigskip \noindent etc.\\

\bigskip \noindent \textbf{Exercise 3}. \label{exercise03_doc03-04}\\

\noindent Let $\{X_{n}\}_{n\geq 0}$ be a sequence of measurable real-valued applications defined on the same measurable space : 
\begin{equation*}
X_{n}:(\Omega ,\mathcal{A})\rightarrow (\overline{\mathbb{R}},\mathcal{B}%
_{\infty }(\overline{\mathbb{R}})).
\end{equation*}

\bigskip \noindent (a) Show that the applications $\sup_{n\geq 0}X_{n}$, $\inf_{n\geq 0}X_{n},$ $\liminf_{n\rightarrow +\infty} X_{n}$ and $\limsup\liminf_{n\rightarrow +\infty} X_{n}$ are measurable.\\

\noindent (b) Show that, if the limit $\lim_{n\rightarrow +\infty} X_{n}(\omega )$ exists for all $\omega \in \Omega $, then the application $\lim_{n\rightarrow +\infty} X_{n}$ is measurable.\\

\noindent \textbf{Remark}. The Point (a) above includes the finite case, that is $\max(X,Y)$ and $\min(X,Y)$ are finite whenever $X$ and $Y$ are.\\ 

\bigskip \noindent \textbf{Exercise 4}. \label{exercise04_doc03-04} Treat the two following questions.\\

\noindent \textbf{(a)}. Show that constant function in are measurable.\\

\noindent \textbf{(b)}. Let $X$ be a measurable real-valued application : 
\begin{equation*}
X:(\Omega ,\mathcal{A})\rightarrow (\overline{\mathbb{R}},\mathcal{B}
_{\infty }(\overline{\mathbb{R}})),
\end{equation*}

\bigskip \noindent and let $c\neq 0$ be a constant in $\mathbb{R}$. Show that $cX$ and $\left\vert X\right\vert $ are measurable.\\

\bigskip \noindent \textbf{Exercise 5}. \label{exercise05_doc03-04} Define on $\mathbb{R}^{k}$ the class $\mathcal{I}_{k}$ of all intervals 
\begin{equation*}
(a,b)=\prod\limits_{j=1}^{k}(a_{i},b_{i}),
\end{equation*}

\bigskip \noindent where $a=(a_{1},...,a_{k})\leq b=(b_{1},...,b_{k})$ [meaning that $a_{i}\leq b_{i},1\leq i\leq k]$. Define also 

\begin{equation*}
\mathcal{I}_{k}(1)=\{]a,b]=\prod\limits_{j=1}^{k}]a_{i},b_{i}],(a,b)\in 
\mathbb{R}^{k},a\leq b\},
\end{equation*}

\begin{equation*}
\mathcal{I}_{k}(2)=\{]-\infty ,a]=\prod\limits_{j=1}^{k}]-\infty
_{i},a_{i}],a\in \overline{\mathbb{R}}^{k}\},
\end{equation*}

\bigskip \noindent and 
\begin{equation*}
\mathcal{I}_{k}(3)=\{]a,b[=\prod\limits_{j=1}^{k}]a_{i},b_{i}[,(a,b)\in 
\mathbb{R}^{k},a\leq b\}.
\end{equation*}

\bigskip \noindent \textbf{(a)}. Show that each of the classes $\mathcal{I}_{k}$ $\mathcal{I}_{k}(1)$, $\mathcal{I}_{k}(2)$ and $\mathcal{I}_{k}(2)$ generates the same $\sigma$-algebra. The $\sigma $-algebra generated by each of them is the usual $\sigma $-algebra on $\mathbb{R}^{k}$ denoted $\mathcal{B}(\mathbb{R}^{k})$.\\

\noindent \textbf{Hints and restrictions}. Restrict yourself to $k=2$ and show only $\sigma(\mathcal{I}_{2})=\sigma(\mathcal{I}_{2}(1))$. You only have to prove that each element of $\mathcal{I}_{2}$ can be obtained from countable sets operations on the elements of $\mathcal{I}_{2}(1))$. \textbf{Show this only for two cases} : $I=[a_{1},a_{2}[\times \lbrack a_{1},b_{2}[$ (a bounded interval) and $J=[a_{1},+\infty \lbrack \times ]-\infty ,b_{2}[$ (an unbounded interval). Express $I$ and $J$ as countable unions and intersections of elements of $\sigma(\mathcal{I}_{2}(1))$.\\

\noindent \textbf{(b)}. Show that $\mathcal{I}_{k}(2,c)=\{\overline{\mathbb{R}}^{k}\}\cup \mathcal{I}_{k}(2,c)$ is a semi-algebra.\\

\noindent \textbf{Hints and restrictions}. Sow this only for $k=2$.\\

\noindent \textbf{(b)} Show that $\mathcal{B}(\mathbb{R}^{k}))$ is the product $\sigma $-algebra corresponding to the product
space $(\mathbb{R},\mathcal{B}(\mathbb{R}))$ by itself $k$ times, that is ($\mathbb{R}^{k},\mathcal{B}(\mathbb{R}^{k}))=(\mathbb{R}^{k},\mathcal{B}(\mathbb{R})^{\otimes k})$.\\

\noindent \textbf{Hint} To show that $\mathcal{B}(\mathbb{R})^{\otimes k}) \subset  \mathcal{B}(\mathbb{R}^{k})$, use the characterization of 
$\mathcal{B}(\mathbb{R})^{\otimes k})$ by the measurability of the projections. \\ 

\noindent \textbf{(d)} Show that $\mathcal{B}(\mathbb{R}^{k})$ is the Borel $\sigma $-algebra on $\mathbb{R}^{k}$ endowed with one of the three equivalent metrics defined
for : $a=(a_{1},...,a_{k})$ and $b=(b_{1},...,b_{k})$
\begin{equation*}
d_{e}(a,b)=\sqrt{\sum_{j=1}^{k}(a_{i}-b_{i})^{2}}\text{,}
\end{equation*}%
\begin{equation*}
d_{m}(a,b)=\sum_{j=1}^{k}\left\vert a_{i}-b_{i}\right\vert \text{,}
\end{equation*}

\bigskip \noindent and

\begin{equation*}
d_{\infty }(a,b)=\max_{1\leq i\leq k}\left\vert a_{i}-b_{i}\right\vert.
\end{equation*}

\bigskip \noindent \textbf{NB}. Before you begin the solution, you are suggested to revise the product spaces DOC 00-14. See also Exercise 2 in DOC 00-02 and its solution in DIC 00-03.\\\

\newpage
\noindent \LARGE \textbf{DOC 03-05 :  Exercises on Measurability 02 - Elementary functions, Marginal functions and applications}. \label{doc03-05}\\

\Large

\bigskip \noindent \textbf{Summary : } You will study the important class of elementary functions, know their properties and show how to approximated measurable real-valued functions
by them.\\
 
\noindent \textbf{Nota Bene.}\newline

\noindent \textbf{(1)} In all this document, the measurable applications we are using are defined on a measurable space $(\Omega ,\mathcal{A})$.\newline

\noindent \textbf{(2)} The expressions of the elementary functions are
always given with the full measurable partition of $\Omega$.\newline

\bigskip \noindent \textbf{Exercise 1}. \label{exercise01_doc03-05} \\

\noindent Consider a measurable and finite partition of $\Omega $, that is, $k$ mutually disjoint and measurable
subsets of $\Omega$, $A_{1}$,..., $A_{k}$ such that 
\begin{equation*}
\Omega =\sum_{1\leq i\leq k}A_{i}
\end{equation*}

\bigskip \noindent and let $\alpha_{1}$, $\alpha_{2}$, ..., and $\alpha_{k}$ be $k\geq 1$ finite real numbers. Let $X$ be the elementary function 
\begin{equation*}
X(\omega )=\alpha _{i}\text{ for }\omega \in A_{i},\text{ }1\leq i\leq k,
\end{equation*}

\bigskip \noindent also denoted by 
\begin{equation}
X=\sum_{1\leq i\leq k}\alpha _{i}1_{A_{i}}.  \label{etag01}
\end{equation}

\noindent \textbf{(a)} Suppose that the values $\alpha _{i},$ $1\leq i\leq k,
$ are distinct. Show that for any \noindent $i\in \{1,2,...,k\},$ 
\begin{equation*}
A_{i}=(X=\alpha _{i}).
\end{equation*}

\bigskip \noindent Deduce from this that if
\begin{equation}
X=\sum_{1\leq j\leq m}\beta _{j}1_{B_{j}}  \label{etag02}
\end{equation}

\bigskip \noindent is an expression of $X$ with distinct values $\beta _{j},$ $1\leq
j\leq \ell$, then, necessarily, $k=\ell $ and there exists a permutation $%
\sigma $ of $\{1,...,k\}$ such that for any $j\in \{1,...,k\}$

\begin{equation*}
\beta _{j}=\alpha _{\sigma (j)}\text{ and }B_{j}=A_{\sigma (j),}
\end{equation*}

\bigskip \noindent meaning that (\ref{etag02}) is simply a re-ordering of (\ref%
{etag01}).\newline

\noindent The expression (\ref{etag01}) of $X$ with distinct values of $%
\alpha _{i},$ is called the canonical expression of $X$.\newline

\noindent \textbf{(b)} Consider an arbitrary expression of an elementary
function : 
\begin{equation}
X=\sum_{1\leq i\leq k}\alpha _{i}1_{A_{i}}.
\end{equation}

\bigskip \noindent Explain how to proceed to get the canonical form.\newline

\noindent \textbf{Hint} : Denote the distinct values of the sequence $\alpha
_{1},..,\alpha _{k} $ as $\beta _{1},...,\beta _{\ell }$ with $\ell \leq k.$
Set $I(j)=\{i,1\leq i\leq k,\alpha _{i}=\beta _{j}\}$ and%
\begin{equation*}
B_{j}=\sum\limits_{i\in I(j)}A_{i}.
\end{equation*}

\bigskip \noindent Conclude!\newline

\noindent \textbf{(c)} Use Questions (a) and (b) to prove that any
elementary function is measurable with the following suggested steps :%
\newline

\noindent \textbf{(c1)} Use the canonical form (\ref{etag01}) of $X.$ Show
that for any subset $B$ of $\mathbb{R}$,

\begin{equation*}
(X\in B)=\sum\limits_{i=1}^{k}(X\in B)\cap A_{i}.
\end{equation*}

\bigskip \noindent \textbf{(c2)} Use Question (a) to establish 
\begin{eqnarray*}
(X &\in &B)=\sum\limits_{i=1}^{k}(X\in B)\cap (X=\alpha _{i}) \\
&=&\sum\limits_{\alpha _{i}\in B}(X=\alpha _{i}).
\end{eqnarray*}

\bigskip \noindent \textbf{(c3)} Conclude that any elementary function $X$ is
measurable.\newline

\noindent \textbf{(d)} Use the same method to prove that for a real-valued
application $X$\ taking at most countable distinct values $\left(
x_{i}\right) _{i\in I}$ where $I$ is countable, we have for any subset of $%
\mathbb{R}$,%
\begin{equation*}
(X\in B)=\sum\limits_{\alpha _{i}\in B}(X=x_{i}).
\end{equation*}

\bigskip \noindent \textbf{Conclude that} : A real-valued application $X$ taking at
most a countable number of distinct values $\left( x_{i}\right) _{i\in I}$
where $I$ is countable is measurable if and only if 
\begin{equation*}
\forall (i\in I),(X=x_{i})\text{ is measurable}.
\end{equation*}

\bigskip \noindent \textbf{(e)} Consider two expressions for an elementary function 
\begin{equation*}
X=\sum_{1\leq i\leq k}\alpha _{i}1_{A_{i}}
\end{equation*}

\noindent and 
\begin{equation*}
X=\sum_{1\leq j\leq \ell }\beta _{j}1_{B_{j}}.
\end{equation*}

\bigskip \noindent Find an expression of $X$ based on the superposition of the two
partitions, that is 
\begin{equation*}
\Omega =\sum_{1\leq i\leq k}\sum_{1\leq j\leq \ell }A_{i}B_{j}.
\end{equation*}

\bigskip \noindent \textbf{(f)} Consider two elementary functions 
\begin{equation*}
X=\sum_{1\leq i\leq k}\alpha _{i}1_{A_{i}}
\end{equation*}

\bigskip \noindent and 
\begin{equation*}
Y=\sum_{1\leq j\leq \ell }\beta _{j}1_{B_{j}}.
\end{equation*}

\bigskip \noindent Provide two expressions for $X$ and $Y$ based on the same
subdivision (\textbf{that} you take as the superposition of the partitions $%
(A_{i})_{1\leq i\leq k}$ and $(B_{i})_{\leq j\leq \ell })$.\newline

\bigskip \noindent \textbf{Exercise 2}. \label{exercise02_doc03-05} (Algebra on the space of elementary
functions $\mathcal{E}$).\newline

\noindent Let $X$, $Y$ and $Z$ be simple functions, $a$ and $b$ be real
scalars.\newline

\noindent By using Question (c) of \textit{Exercise 1}, show the following
properties :\newline

\noindent \textbf{(a1)} $aX+bY\in \mathcal{E}$.\newline

\noindent \textbf{(a2)} $XY\in \mathcal{E}$.\newline

\noindent \textbf{(a3)} If $Y$ is everywhere different of zero, then $X/Y\in 
\mathcal{E}$.\newline

\noindent \textbf{(a5)} $(X\leq Y)\Rightarrow (X+Z\leq Y+Z)$. If $a\geq 0$,
then $X\leq Y\Rightarrow aX\leq bY$.\newline

\noindent \textbf{(a6)} $max(A,X)\in \mathcal{E}$, $min(X,Y)\in \mathcal{E}$.\\

\bigskip \noindent \textbf{Terminology} : (You may skip the remarks below). 
\newline

\noindent \textbf{(A)} Needless to prove, the product of finite functions
and the sum of finite functions are commutative and associative. The product
of finite functions is distributive with the sum of finite functions. Taking
into account these facts, we state the following facts.\newline

\noindent By (a1), $(\mathcal{E},+,\cdot )$ is a linear space (where $\cdot $
stands for the external multiplication by real scalars).\newline

\noindent By (a1) and (a2), and the remarks (A), $(\mathcal{E},+,\times )$
is a ring with unit $1_{\Omega}$ (where $\times $ stands for the product
between two functions), and is an algebra of functions.\newline

\noindent By (a6), $(\mathcal{E},\leq )$ is a lattice space, meaning closed
under finite maximum and finite minimum.\newline

\noindent By (a1), (a5) and (a6), $(\mathcal{E},+,\cdot ,\leq )$ is a Riesz
space, that is a lattice (a6) vector space (a1) such that the order is
compatible with the linear structure.\newline

\bigskip \noindent \textbf{Exercise 3}. \label{exercise03_doc03-05}  Let $X\geq 0$ be
measurable. For $n\geq 1$ fixed, consider the following partition of $\mathbb{R}_{+}$ :

\begin{equation*}
\mathbb{R}_{+}=\sum_{k=1}^{2^{2n}}[\frac{k-1}{2^{n}},\frac{k}{2^{n}}[\text{\ }+\text{ }[2^{n},+\infty ]=\sum_{k=1}^{2^{2n}+1}A_{k}.
\end{equation*}

\bigskip \noindent and set 
\begin{equation*}
X_{n}=\sum_{k=1}^{2^{n}}\frac{k-1}{2^{n}}\text{ }1_{(\frac{k-1}{2^{2n}}\leq
X<\frac{k}{2^{n}})}+2^{n}\text{ }1_{(X\geq 2^{2n})}.
\end{equation*}

\bigskip \noindent \textbf{(a)} Show that for each $n\geq $, $X_{n}$ is an elementary
function.\newline

\noindent \textbf{(b)} by using this implication : for any $\omega \in
\Omega ,$ for any $n\geq 1,$ 
\begin{equation*}
X(\omega )\geq \frac{k-1}{2^{n}},k\leq 2^{2n}
\end{equation*}

\bigskip \noindent implies 
\begin{equation}
X_{n}(\omega )=\left\{ 
\begin{tabular}{lll}
$\frac{j-1}{2^{2n}}$ & if for some $j,1\leq j\leq $ $2^{2n}$ : & $\frac{j-1}{%
2^{2n}}\leq \frac{j-1}{2^{n}}\leq X(\omega )<\frac{j}{2^{n}}$ \\ 
&  &  \\ 
$2^{n}$ & otherwise & 
\end{tabular}
\right. ,
\end{equation}

\noindent and by discussing the two cases ($X(\omega )\geq 2^{n}$ or $\frac{%
k-1}{2^{n}}\leq X(\omega )<\frac{k}{2^{n}}$ for some $k\leq 2^{2n}),$ show
that $X_{n}$ is non-decreasing.\newline

\noindent \textbf{(c)} Show that $X_{n}\rightarrow X$ as $n\rightarrow +\infty $ by proceeding as follows. If $X(\omega )=+\infty ,$ show that $X_{n}(\omega )=2^{n}\rightarrow +\infty .$ If $\frac{k-1}{2^{n}}\leq X(\omega )<\frac{k}{2^{n}}$ for some $k\leq 2^{2n},$ show that for $n$ large enough 
\begin{equation*}
0\leq X(\omega )-X_{n}(\omega )\rightarrow \frac{1}{2^{n}}\rightarrow 0.
\end{equation*}

\bigskip \noindent \textit{\textbf{Conclusion} : Any non-negative measurable
applications is limit of a non-decreasing sequence of non-negative
elementary functions.}\\

\bigskip \noindent \textbf{Exercise 4}. \label{exercise04_doc03-05} Let $X$, $Y$ be measurable \textbf{%
finite} applications, $a$ and $b$ be real scalars.\newline

\noindent \textbf{(a)} Show that any real-valued measurable function $Z,$
possibly taking infinite values, is limit of elementary functions by
considering the decomposition%
\begin{equation*}
Z=\max (Z,0)-\max (-Z,0),
\end{equation*}

\bigskip \noindent where $Z^{+}=\max (Z,0)$ and $Z^{-}=\max (-Z,0)$ are non-negative
functions and are respectively called : positive and negative parts of $Z$.
\ By applying the conclusion of Question (c) of Exercise 3, show that any
real-valued measurable function is limit of elementary functions. \newline

\noindent \textbf{(b)} Deduce from this that $aX+bY$, $XY$, $\max (X,Y)$ and 
$\min (X,Y)$ are measurable. If $Y(\omega )\neq 0$ for all $\omega \in \Omega
$, then $X/Y$ is measurable.\newline

\bigskip \noindent \textbf{Exercise 5}. \label{exercise05_doc03-05} Consider a product space endowed with its product $\sigma $-algebra $\prod_{1\leq i\leq k}\Omega_{i},\bigotimes_{1\leq i\leq k}\mathcal{A}_{i},$ and let  $X$ $\ $be a function defined on from this product space to  in $\overline{\mathbb{R}}$ :
\ 
\begin{equation*}
\begin{tabular}{lllll}
$X$ & $:$ & $\left( \prod_{1\leq i\leq k}\Omega _{i},\bigotimes_{1\leq i\leq
k}\mathcal{A}_{i}\right) $ & $\longmapsto $ & $(\overline{\mathbb{R}},%
\mathcal{B}_{\infty }(\overline{\mathbb{R}}))$ \\ 
&  & $(\omega _{1,}\omega _{2,}...,\omega _{k})$ & $\hookrightarrow $ & $%
X(\omega _{1,}\omega _{2,}...,\omega _{k})$.
\end{tabular}
\end{equation*}

\bigskip \noindent For any fixed $i\in \{1,2,...\}$ , define the  $i-$th partial function $X$ by considering it as a function of $\omega _{i}$ only and by fixing the other variables. Let $(\omega _{1,}\omega _{2,}...,\omega_{i-1},\omega _{i+1},...,\omega _{k})$ fixed and define $X^{(i)}$ the as $i-$th partial function $X$ 
\begin{equation*}
\begin{tabular}{lllll}
$X^{(i)}$ & $:$ & $\left( \Omega _{i},\mathcal{A}_{i}\right) $ & $%
\longmapsto $ & $(\overline{\mathbb{R}},\mathcal{B}_{\infty }(\overline{%
\mathbb{R}}))$ \\ 
&  & $\mathbf{\omega }_{i}$ & $\hookrightarrow $ & $X^{(i)}(t)=f(\omega
_{1,}\omega _{2,}...,\omega _{i-1},\mathbf{\omega }_{i},\omega
_{i+1},...,\omega _{k})$.
\end{tabular}
\end{equation*}

\bigskip \noindent The objective of this exercise is to show that that partial functions of real-valued and measurable functions of several variables are also measurable. The proof for $k=2$ may re-conducted word by word to a general $k\geq 2.$ In the sequel $k=2$. Let us fix $\omega _{1}\in \Omega _{1}$ and consider the partial function
\begin{equation*}
\begin{tabular}{lllll}
$X$ & $:$ & $\left( \Omega _{2},\mathcal{A}_{2}\right) $ & $\longmapsto $ & $%
(\overline{\mathbb{R}},\mathcal{B}_{\infty }(\overline{\mathbb{R}}))$ \\ 
&  & $\mathbf{\omega }_{2}$ & $\hookrightarrow $ & $X^{(i)}=X_{\omega
_{1}}(\omega _{2})=X(\omega _{1},\omega _{2})$.
\end{tabular}
\end{equation*}

\bigskip \noindent \textbf{(a)} Show that for any non-decreasing sequence of applications $X_{n}:(\Omega
_{1}\times \Omega _{2},\mathcal{A}_{i}\otimes \mathcal{A}_{2})\longmapsto (%
\overline{\mathbb{R}},\mathcal{B}_{\infty }(\overline{\mathbb{R}})),$ the
sequence of partial applications $\left( X_{n}\right) _{\omega _{1}}$ is
also non-decreasing.\\

\noindent \textbf{(b)} Show that for any sequence of applications $X_{n}:(\Omega _{1}\times \Omega _{2},\mathcal{A}_{i}\otimes \mathcal{A}_{2})\longmapsto (\overline{\mathbb{R}},\mathcal{B}_{\infty }(\overline{\mathbb{R}})),n\geq 1$, converging to the application $X_{n}:(\Omega _{1}\times \Omega _{2},\mathcal{A}_{i}\otimes \mathcal{A}_{2})\longmapsto (\overline{\mathbb{R}},\mathcal{B}%
_{\infty }(\overline{\mathbb{R}})),$ the sequence of partial applications $\left( X_{n}\right) _{\omega _{1}},n\geq 1,$ is also non-decreasing also converges to the partial function $X_{\omega _{1}}$ of $X$.\\

\noindent \textbf{(b)} Show for $Y=1_{A}:(\Omega _{1}\times \Omega _{2},\mathcal{A}_{i}\otimes \mathcal{A}_{2})\longmapsto (\overline{\mathbb{R}},\mathcal{B}_{\infty }(\overline{\mathbb{R}}))$ is the indicator function of a set $A\subset \Omega_{1}\times \Omega _{2}$, then the partial function
\begin{equation*}
Y_{\omega _{1}}(\omega _{2})=Y(\omega _{1},\omega _{2})=1_{A}(\omega
_{1},\omega _{2})=1_{A_{\omega _{1}}}(\omega _{2}), 
\end{equation*}

\bigskip \noindent is the indicator function of the section $A_{\omega _{1}}$ of $A$ at $\omega_{1}$\\

\noindent \textbf{(c)} Deduce from Question (b) than partial functions of indicator functions
of measurable sets are measurable and partial functions of elementary
functions defined on $(\Omega _{1}\times \Omega _{2},\mathcal{A}_{i}\otimes 
\mathcal{A}_{2})\longmapsto (\overline{\mathbb{R}},\mathcal{B}_{\infty }(%
\overline{\mathbb{R}}))$ are measurable.\\

\noindent \textbf{(d)} Let \ \ \ $X:(\Omega _{1}\times \Omega _{2},\mathcal{A}_{i}\otimes \mathcal{A}_{2})\longmapsto (\overline{\mathbb{R}},\mathcal{B}_{\infty }(\overline{\mathbb{R}}))$ be a measurable application. Use Question (a) of Exercise 4 and consider a sequence $X_{n}:(\Omega _{1}\times \Omega _{2}, \mathcal{A}_{i}\otimes \mathcal{A}_{2})\longmapsto (\overline{\mathbb{R}}, \mathcal{B}_{\infty }(\overline{\mathbb{R}})),n\geq 1,$ of elementary functions converging to $X.$ Combine the previous questions to prove that $%
X$ has measurable partial functions.\\

\newpage
\noindent \LARGE \textbf{DOC 03-06 : Exercises on Measurability 03 - Classical functions in Real Analysis}. \label{doc03-06}\\

\Large 

\bigskip \noindent \textbf{Summary : } You study the measurability of usual types of real-valued functions used in Calculus courses : functions which are continuous, right and left continuous, upper and lower semi-continuous, increasing or decreasing, etc.\\

\bigskip \noindent \textbf{Exercise 1}. \label{exercise01_doc03-06} Le $f$ \ be a non-decreasing function
from $\mathbb{R}$ to $\mathbb{R}$. Show that the set of discontinuity points
is at most countable.\newline

\bigskip \noindent \emph{Proceed to this way}. \noindent Let $D(n)$ the number discontinuity points of $f$ in $]-n,n[$. Recall that a
real number $x$ is discontinuity point of a non-decreasing function $f$ if and only if 
\begin{equation*}
f(x+)-f(x-)>0,
\end{equation*}

\bigskip \noindent where $f(x-)$ is the left-hand limit of $f$ at $x$ and $f(x+)$ is the
right-hand limit. In that case $f(x+)-f(x-)$ is the discontinuity jump. Next
let 
\begin{equation*}
D_{k}(n)=\{x\in D(n),\text{ }f(x+)-f(x-)>1/k\}.
\end{equation*}

\bigskip \noindent Pick $x_{1},...,x_{m}$ from $D_{k}(n)$ and justify the
inequality 

\begin{equation*}
\sum_{1\leq i\leq m}f(x_{i}+)-f(x_{i}-) \leq f(n)-f(-n).
\end{equation*}

\bigskip \noindent Deduce from this that 
\begin{equation*}
m\leq k\times (f(n)-f(-n)),
\end{equation*}

\bigskip \noindent and then $D_{k}(n)$ is finite.\\

\noindent Conclude.\newline

\bigskip \noindent \textbf{Exercise 2}. \label{exercise02_doc03-06} Show that any right-continuous or
left-continuous function $X$ from $\mathbb{R}$ to $\mathbb{R}$ is
measurable.\newline

\noindent Hint : Let $X$ be right-continuous and set for all $n\geq 1$ 
\begin{equation*}
X_{n}^{+}(t)=\sum_{k=-\infty }^{k=+\infty }X\left( \frac{k+1}{2^{n}}\right)
1_{\left] \frac{k}{2^{-n}},\frac{k+1}{2^{-n}}\right] }(t).
\end{equation*}

\bigskip \noindent and for $N\geq 1$

\begin{equation*}
X_{n,N}^{+}(t)=\sum_{k=-N}^{k=N}X\left( \frac{k+1}{2^{n}}\right) 1_{\left] 
\frac{k}{2^{-n}},\frac{k+1}{2^{-n}}\right] }(t).
\end{equation*}

\bigskip \noindent $X_{n,N}^{+}$ is measurable because it is elementary
function (refer to the document on elementary functions).

\bigskip \noindent Show that $X_{n,N}^{+}$ converges to $X_{n}^{+}$ as $N
\rightarrow +\infty$ ($n$ being fixed).\newline

\bigskip \noindent Conclude that each $X_{n}$ is measurable and that $X_{n}$
converges to $X$ as $n \rightarrow +\infty$. Make a final conclude.\newline

\noindent Hint : If $X$ be left-continuous and set for all $n\geq 1$ 
\begin{equation*}
X_{n}^{-}(t)=X_{n}^{+}(t)=\sum_{k=-\infty }^{k=+\infty }X\left( \frac{k}{%
2^{n}}\right) 1_{\left[ \frac{k}{2^{-n}},\frac{k+1}{2^{-n}}\right[ }(t).
\end{equation*}

\bigskip \noindent Next proceed similarly to the case of right-continuous
functions.\\

\bigskip \noindent \textbf{Exercise 3}. \label{exercise03_doc03-06} Let $X$ be an application from from $\mathbb{R}$ to $\mathbb{R}$ such that at each $t\in \mathbb{R}$ the left limit $f(x-)$ exists (resp. at each $t\in \mathbb{R}$ the right-hand limit $f(x+)$ exists.\\
 
\noindent \textbf{(a)} Use know facts on the left continuity of $t\mapsto X(t-)$ and the left-continuity of $(x\mapsto X(t+))$) to show that $X$ is measurable with help of Exercise 2.\\

\noindent \textit{(b)} Suppose that at each $t\in \mathbb{R}$ the left limit $f(x-)$ exists (resp. at each $t\in \mathbb{R}$ the right-hand limit $f(x+)$ exists). Show that the function  $t\mapsto X(t-)$ is left-continuous (resp. $(x\mapsto X(t+))$ is right-continuous). [You may skip this question and consider it as a toplogy exercise. But you are recommended to read the proof].\\

\bigskip \noindent \textbf{Exercise 4}.\label{exercise04_doc03-06} Let $X$ be an application $\mathbb{R}
$ to $\mathbb{R}$ such that the set of discontinuity points $D$ of $X$  is
at most countable and $X$ admits right-hand limits at all points (or $X$
admits  left-hand limits at all points ). Justify the identity 
\begin{equation*}
\forall t\in \mathbb{R}, X(t)=X(t+)\times 1_{D^{c}}(t)+\sum_{t\in D}X(t)\times 1_{\{t\}}(t).
\end{equation*}

\bigskip \noindent and conclude that $X$ is measurable.\\

\bigskip \noindent \textbf{Exercise 5}. \label{exercise05_doc03-06}  Use Exercises 1 and 3 to show that any monotone function $f:\mathbb{R}\rightarrow \mathbb{R}$ is 
measurable.\\

\bigskip \noindent \textbf{Exercise 6}. \label{exercise06_doc03-06} Let $X$ be a convex application $%
\mathbb{R}$ to $\mathbb{R}$.\\

\noindent \textbf{(a)} Use a know fact on the right-continuity of a convex function (as reminded in Question (b) below) and conclude with the help of Exercise 2.\\

\noindent \textbf{(b)} Show that any convex function is continuous. [You may skip this question and consider it as a topology exercise. But you are recommended to read the proof].\\

\bigskip \noindent \textbf{Exercise 7}. \label{exercise07_doc03-06} Let $X$ be a function from a metric
space $(E,d)$ to $\mathbb{R}$. Suppose that if $X$ is upper semi-continuous or lower semi-continuous.\\ 

\noindent \textbf{Reminder}.\\

\noindent \textit{(1)} A function $X:(E,d)\rightarrow \mathbb{R}$ is lower-continuous if and only for each $t\in E$,

\begin{equation*}
(\forall 0<\epsilon <0)(\exists \eta >0),d(s,t)<\eta \Rightarrow X(s)\leq
X(t)+\epsilon. \text{ (LSC) }
\end{equation*}

\bigskip \noindent \textbf{(2)} A function $X:(E,d)\rightarrow \mathbb{R}$ is upper semi-continuous if
and only for each $t\in E$,

\begin{equation*}
(\forall 0<\epsilon <0)(\exists \eta >0),d(s,t)<\eta \Rightarrow X(s)\geq
X(t)-\epsilon. \text{ (LSC) }
\end{equation*}

\bigskip
\noindent (a) Use know facts on semi-continuous functions in Topology, as reminded in Question (b), to show that $X$ is measurable.\\

\noindent (b) Suppose that $X$ is lower semi-continuous. Show that for  $c\in \mathbb{R},(X\leq c)$ is closed. Deduce that if $X$ is upper semi-continuous, then
for any $c\in \mathbb{R},(X\geq c)$ is closed. [You may skip this question and consider it as a topology exercise. But you are recommended to read the proof].\\

\bigskip \noindent Here, we we introduce semi-continuous functions in special forms. It is proved in http://dx.doi.org/10.16929/jmfsp/2017.001, the equivalence of other more classical definitions.
We will find these expressions powerful to deal with Lebesgue-Stieljes integrals.\\

\bigskip \noindent \textbf{Exercise 8}. \label{exercise08_doc03-06} Let $f$ be a function $\mathbb{R}\rightarrow \mathbb{R}.$ Let $x\in \mathbb{R}.$ Define the limit superior and the limit inferior of $f$ at $x$ by

\begin{equation*}
f^{\ast}(x)=\lim \sup_{y\rightarrow x}f(y)=\lim_{\varepsilon \downarrow
0}\sup \{f(y),y\in ]x-\varepsilon ,x+\varepsilon \lbrack \}
\end{equation*}

\noindent and S
 
\begin{equation*}
f_{\ast }(x)=\lim \inf_{y\rightarrow x}f(y)=\lim_{\varepsilon \downarrow
0}\inf \{f(y),y\in ]x-\varepsilon ,x+\varepsilon \}.
\end{equation*}

\bigskip \noindent Question (1) Show that $f^{\ast }(x)$ and $f_{\ast }(x)$ always exist. Hint remark
that $\sup \{f(y),y\in ]x-\varepsilon ,x+\varepsilon \}$ is non-increasing
and $\inf \{f(y),y\in ]x-\varepsilon ,x+\varepsilon \}$ is non-decreasing as 
$\varepsilon $ decreases to zero.\\

\bigskip 
\noindent Question (2) Show that \ $(-f)^{\ast }=-(f_{\ast }).$.\\

\noindent Question (3)Show that $f$ is continuous at $x$ if and only if $f^{\ast }(x)$ and $%
f_{\ast }(x)$, and then : $f$ is continuous if and only if $f=$ $f^{\ast}=f_{\ast }.$.\\

\noindent Question (4) Define : $f$ is upper semi-continuous if and only if at $x$ if and only if $%
f=f^{\ast }$ and $f$ is lower semi-continuous iff $f=f_{\ast }$.\\

\noindent Show that $f^{\ast}$ is upper semi-continuous. Exploit (2) to extend this to $f_{\ast}$ and show
that $f_{\ast }$ is lower semi-continuous.\\

\bigskip Question (5) Now let $f$ be upper semi-continuous. Set for each $n\geq 1$
\begin{equation*}
f_{n}(x)=\sum_{k=-\infty }^{k=+\infty }\sup \{f(z),z\in [k2^{-n},(k+1)2^{-n}[ 1_{]k2^{-n},(k+1)2^{-n}[}(x)+\sum_{s\in D_n}f(s)1_{\{s\}}.
\end{equation*}

\bigskip \noindent Show that $f_{n}$ converges to $f$. Conclude that $f$ is measurable.\\

\bigskip  \noindent  Question (6) Use the opposite argument to show that a lower semi-continuous function is measurable.\\

\bigskip \noindent Question (7) The same proof you use in Question (5) is also valid in the more general case where $I=[a,b]$ or $I=\mathcal{R}$ and for each $n$,
$I_{i,n}=(x_{i,n},x_{i+1,n})$, $i\in J$, are consecutive intervals, with non zero length, which partition $I$ such that

$$
\sup_{i} |x_{i+1,n}-x_{i+1,n}| \rightarrow 0 \ as \ n\rightarrow +\infty. 
$$

\bigskip \noindent By setting $D=\cup_n \cup_i x_{i,n}$, the proof of Question (5) leads to : for any $x\notin D$,

$$
f_{n}(x)=\sum_{i} \sup \{f(z),z\in ]x_{i,n},x_{i+1,n}[\} 1_{]x_{i,n},x_{i+1,n}[}(x) \rightarrow f^{\ast}(x)
$$

\bigskip \noindent 

$$
f_{n}(x)=\sum_{i} \inf \{f(z),z\in ]x_{i,n},x_{i+1,n}[\} 1_{]x_{i,n},x_{i+1,n}[}(x) \rightarrow f_{\ast}(x).
$$

\newpage
\noindent \LARGE \textbf{DOC 03-07 : Discover Exercises on Measurability 01 with Solutions with solutions}. \label{doc03-07}\\

\Large

\bigskip \noindent \textbf{Exercise 1}. \label{exercise01_sol_doc03-07}\\

\noindent \textbf{(a)} Show that each of these classes generates the same 
$\sigma $-algebra as the class of all intervals $\mathcal{I}$. Do only one
or two cases.

\begin{equation*}
\mathcal{I}(1)=\{]-\infty ,a[,\ a\in \mathbb{R}\}\text{ , }\mathcal{I}(2)=\{]-\infty ,a],\ \ a\in \overline{\mathbb{R}}\}.
\end{equation*}

\bigskip 
\begin{equation*}
\mathcal{I}(3)=\{]-\infty ,a],a\in \mathbb{R}\}\text{ , }\mathcal{I}(4)=\{]-\infty ,a],\ \ a\in \overline{\mathbb{R}}\}
\end{equation*}

\bigskip 
\begin{equation*}
\mathcal{I}(5)=\{]a,+\infty \lbrack, \ \ a\in \mathbb{R}\},\mathcal{I}(6)=\{[a,+\infty \lbrack, \ \ a\in \overline{\mathbb{R}}\}
\end{equation*}

\bigskip
\begin{equation*}
\mathcal{I}(7)=\{[a,+\infty \lbrack,\ \ a\in \mathbb{R}\}\text{ , }\mathcal{I}(8)=\{]a,+\infty \lbrack, \ \ a\in \overline{\mathbb{R}}\}
\end{equation*}

\bigskip 
\begin{equation*}
\mathcal{I}(9)=\{]a,b],a,b\in \mathbb{R}\}\ ,\ \mathcal{I}(10)=\{]a,b],\ \ (a,b)\in \overline{\mathbb{R}}^2\}
\end{equation*}

\bigskip 
\begin{equation*}
I(11)=\{]a,b[,a,b\in \mathbb{R}\}\ ,\ I(12)=\{]a,b[,\ \ (a,b)\in \overline{\mathbb{R}}^2\}
\end{equation*}

\bigskip 
\begin{equation*}
\mathcal{I}(13)=\{[a,b[,\ (a,b) \in \mathbb{R}^2\}\text{ , }\mathcal{I}(14)=\{[a,b[,\ \ (a,b)\in \overline{\mathbb{R}}^2\}
\end{equation*}

\bigskip 
\begin{equation*}
\mathcal{I}(15)=\{[a,b],(a,b) \in \mathbb{R}^2\}\text{ , }\mathcal{I}(16)=\{[a,b],\ \ (a,b)\in \overline{\mathbb{R}}^2\}
\end{equation*}

\bigskip \noindent \textbf{(b)} show that the class $\mathcal{I}(4)=\{]-\infty ,a],a\in 
\overline{\mathbb{R}}\}$ is a semi-algebra in $\mathbb{R}$\\

\bigskip \noindent \textbf{Solutions}.\\

\noindent \textbf{Question (a)}. We must give 14 answers. It is too long. We only give two
answers. The others will be very similar.\\

\noindent \textbf{(a1)}. Let us show that $\sigma (\mathcal{I})=\sigma (\mathcal{I}(1)).$ The
inclusion $\sigma (\mathcal{I}(1))\subset $ \ $\sigma (\mathcal{I})$ is
obvious since $\mathcal{I}(1)\subset \mathcal{I}$. To prove that  \ $%
\sigma (\mathcal{I})\subset \sigma (\mathcal{I}(1))$ it will enough to show
that any interval of $\mathbb{R}$\ may be obtained from a countable number of
sets operations on the elements of $\mathcal{I}(1)$. Let us check that in
the following table 

\begin{equation*}
\begin{tabular}{|l|l|l|l|}
\hline
Line & Intervals of $\mathbb{R}$ & Is obtained from $\mathcal{I}(1)$ by  & 
Remarks : Uses\\
\hline 
1 & $]-\infty ,+\infty \lbrack $ & $\bigcup\limits_{n\geq 1}]-\infty ,n]\in
\sigma (\mathcal{I}(1))$ & Definition of $\mathcal{I}(1)$ \\
\hline 
2 & $]-\infty ,a[$ & $\bigcup\limits_{n\geq 1}]-\infty ,a-\frac{1}{n}]\in
\sigma (\mathcal{I}(1))$ & Definition of $\mathcal{I}(1)$ \\
\hline 
3 & $]-\infty ,a]$ & $\mathcal{I}(1)$ & Definition of $\mathcal{I}%
(1)$ \\ 
\hline
4 & $]a,+\infty \lbrack $ & $(]-\infty ,a])^{c}\in \sigma (\mathcal{I}(1))$
& Definition of $\mathcal{I}(1)$ \\
\hline 
5 & $[a,+\infty \lbrack $ & $(]-\infty ,a[)^{c}\in \sigma (\mathcal{I}(1))$
& Line 4 \\
\hline 
7 & $]a,b]$ & $]-\infty ,b]\cap (]-\infty ,a]^{c})\in \sigma (\mathcal{I}(1))
$ & Definition of $\mathcal{I}(1)$ \\
\hline 
8 & $\{a\}$ & $\bigcap\limits_{n\geq 1}]a,a+\frac{1}{n}]\in \sigma (%
\mathcal{I}(1))$ & line 7 \\
\hline 
9 & $]a,b[$ & $]a,b]\backslash \{a\}\in \sigma (\mathcal{I}(1))$ & Line
8 \\
\hline 
10 & $[a,b[$ & $(\{a\}+]a,b])\backslash \{b\}\in \sigma (\mathcal{I}(1))$ & 
Line 8 \\
\hline 
11 & $[a,b]$ & $\{a\}+]a,b]\in \sigma (\mathcal{I}(1))$ & Line 8.\\ \hline

\end{tabular}
\end{equation*}

\bigskip \noindent We checked that any interval of $\mathbb{R}$ belongs to $\sigma (\mathcal{I}%
(1))$. Then $\sigma (\mathcal{I})\subset \sigma (\mathcal{I}(1))$ and next $%
\sigma (\mathcal{I})=\sigma (\mathcal{I}(1))$.\\

\bigskip \noindent \textbf{(a2)} Let us show that $\sigma (\mathcal{I})=\sigma (\mathcal{I}(9))$. We
re-conduct the same method by using the table%

\begin{equation*}
\begin{tabular}{|l|l|l|l|}
\hline
Line & Intervals of $\mathbb{R}$ & Is obtained from $\mathcal{I}(9)$ by  & 
Remarks : Uses \\ 
\hline
1 & $]-\infty ,+\infty \lbrack $ & $\bigcup\limits_{n\geq 1}]-n,n]\in
\sigma (\mathcal{I}(9))$ & Definition of $\mathcal{I}(9)$ \\
\hline 
2 & $\{a\}$ & $\bigcap\limits_{n\geq 1}]a,a+\frac{1}{n}]\in \sigma (%
\mathcal{I}(9))$ & Definition of $\mathcal{I}(9)$ \\ 
\hline
3 & $]-\infty ,a]$ & $\bigcup\limits_{n\geq 1}]-n,a]\in \sigma (\mathcal{I}%
(9))$ & Definition of $\mathcal{I}(9)$ \\ 
\hline
4 & $]-\infty ,a[$ & $]-\infty ,a]\backslash \{a\}\in \sigma (\mathcal{I}(9))
$ & Line 2 \\ 
\hline
5 & $]a,+\infty \lbrack $ & $(]-\infty ,a])^{c}\in \sigma (\mathcal{I}(9))$
& Line 3 \\ 
\hline
6 & $[a,+\infty \lbrack $ & $(]-\infty ,a[)^{c}\in \sigma (\mathcal{I}(9))$
& Line 4 \\ 
\hline
7 & $]a,b]$ & $\in \sigma (\mathcal{I}(9))$ & Definition of $\mathcal{I}(9)$  \\ 
\hline
8 & $]a,b[$ & $]a,b]\backslash \{a\}\in \sigma (\mathcal{I}(9))$ & Line
2 \\ 
\hline
9 & $[a,b[$ & $(\{a\}+]a,b])\backslash \{b\}\in \sigma (\mathcal{I}(9))$ & 
Line 2 \\ 
\hline
10 & $[a,b]$ & $\{a\}+]a,b]\in \sigma (\mathcal{I}(9))$ & Line 2.\\
\hline
\end{tabular}
\end{equation*}

\bigskip \noindent \textbf{Question (b)}. Let us check the three following points.\\

\noindent \textbf{(b1)}. Since we are restricted on $\mathbb{R}$, we have  $\mathbb{R}=]-\infty
,+\infty ]\in \mathcal{I}(4)$.\\

\noindent \textbf{(b2)} Consider two elements of $\mathcal{I}(4),$ $]-a_{1},a_{2}]$ and $%
]b_{1},b_{2}]$. We have
\begin{equation*}
]-a_{1},a_{2}]\cap ]b_{1},b_{2}]=]\max (a_{1},a_{2}),\min (b_{1},b_{2})]\in 
\mathcal{I}(4).
\end{equation*}

\bigskip \noindent \textbf{(b3)}. For any  element $]a,b]$ of $\mathcal{I}(4)$, we have $]a,b]=$ $%
]-\infty ,a]+$ $]b,+\infty \lbrack .$ But we may write on $\mathbb{R}$
: $]-a,b]^{c}=$ $]-\infty ,a]+$ $]b,+\infty ].$ This proves that the
complement of an element of $\mathcal{I}(4)$ (in the case where $a=-\infty$) or is a sum of two elements of $\mathcal{I}(4)$.\\

\bigskip \bigskip \noindent \textbf{Exercise 2}. \label{exercise02_sol_doc03-07} Consider the extended real line 
\begin{equation*}
\overline{\mathbb{R}}=[-\infty ,+\infty ]=\mathbb{R} \cup \{-\infty,+\infty \},
\end{equation*}

\bigskip \noindent obtained by adding the points $-\infty $ and +$\infty $ to $\mathbb{%
R}$. Define 

\begin{equation*}
\mathcal{B}_{\infty }(\overline{\mathbb{R}}\mathbb{)}=\{A,A\cup \{+\infty
\},A\cup \{+\infty \},A\cup \{-\infty ,+\infty \},\text{ }A\in \mathcal{B}(%
\mathbb{R)}\}.
\end{equation*}

\bigskip \noindent \textbf{(a)}. Show that $\mathcal{B}_{\infty }(\overline{\mathbb{R}}\mathbb{)}$\ is
the smallest sigma-algebra for which the Borel sets of \ $\mathbb{R}$, $%
\{-\infty \},\{+\infty \}$ are measurable.\\

\noindent \textbf{(b)} Show that each of the classes given in Exercise 1 also generates $%
\mathcal{B}_{\infty }(\overline{\mathbb{R}}\mathbb{)}.$

\bigskip \noindent \textbf{(c)} Show that any function

\begin{equation*}
X:(\Omega ,\mathcal{A})\rightarrow (\mathbb{R},\mathcal{B}(\mathbb{R}))\text{
or }X:(\Omega ,\mathcal{A})\rightarrow (\mathbb{R},\mathcal{B}_{\infty }(%
\overline{\mathbb{R}})),
\end{equation*}

\bigskip \noindent is measurable if and only if

\begin{equation*}
(C1):\text{For all }a\in \mathbb{R},(X<a)\in \mathcal{A},
\end{equation*}

\bigskip \noindent if and only if

\begin{equation*}
(C2):\text{For all }a\in \mathbb{R},(X\leq a)\in \mathcal{A}\text{,}
\end{equation*}

\bigskip \noindent if and only if

\begin{equation*}
(C3):\text{For all }a\in \mathbb{R},(X>a)\in \mathcal{A},
\end{equation*}

\bigskip \noindent if and only if

\begin{equation*}
(C4):\text{For all }a\in \mathbb{R},(X\geq a)\in \mathcal{A},
\end{equation*}

\noindent etc.\\

\bigskip \noindent \textbf{Solutions}.\\

\noindent \textbf{Question (a)}. It is easy to check that  $\mathcal{B}_{\infty }(\overline{%
\mathbb{R}}\mathbb{)}$ is a $\sigma$-algebra. By taking $A=\emptyset$, we also see that $%
\{+\infty \}$ and $\{-\infty \}$ are measurable and next $\{-\infty ,+\infty
\}$ is measurable. Finally, any $\sigma $-algebra including \ $\mathcal{B}(%
\mathbb{R)}$\ and for which $\{-\infty \}$ and $\{+\infty \}$ are mea0surable
includes $\mathcal{B}_{\infty }(\overline{\mathbb{R}}\mathbb{)}$. Then $%
\mathcal{B}_{\infty}(\overline{\mathbb{R}}\mathbb{)}$\ is the smallest
$\sigma$-algebra for which the Borel sets of \ $\mathbb{R}$, and $\{-\infty
\},\{+\infty \}$ are measurable.\\

\bigskip \noindent \textbf{Question (b)}. Consider the class $\mathcal{I}$ of all intervals of $\mathbb{R}
$. Remark that the elements of $\mathcal{I}$\ are subsets of $\overline{\mathbb{R}}$
and we may consider  $\mathcal{B}_{\infty }^{\prime }(\overline{\mathbb{R}}%
\mathbb{)}$ to be the $\sigma $-algebra generated by $\mathcal{I}$ $\ $on $%
\overline{\mathbb{R}}.$ By definition of $\mathcal{B}_{\infty }(\overline{\mathbb{R}})$, we have 
$\mathcal{B}_{\infty }^{\prime}(\overline{\mathbb{R}})\subset \mathcal{B}_{\infty }(\overline{\mathbb{R}}).$
Conversely, we see that for any $a\in \mathbb{R},]n,+\infty ]=]-\infty
,+a]^{c}\in \mathcal{B}_{\infty }^{\prime }(\overline{\mathbb{R}}\mathbb{)}$
and $[-\infty ,a]=]a,+\infty \lbrack ^{c}\in \mathcal{B}_{\infty }^{\prime }(%
\overline{\mathbb{R}}\mathbb{)}.$ And it is clear that the elements $A$ of  $%
\mathcal{B}(\mathbb{R)}$ are in $\mathcal{B}_{\infty}^{\prime }(\overline{%
\mathbb{R}}\mathbb{)}$ and and for these $A$, we have
\begin{equation*}
A\cup \{+\infty \}=A\bigcup \left( \bigcap\limits_{n\geq 1}]n,+\infty
]\right) \in \mathcal{B}_{\infty }^{\prime }(\overline{\mathbb{R}}),
\end{equation*}

\begin{equation*}
A\cup \{-\infty \}=A\bigcup \left( \bigcap\limits_{n\geq 1}[-\infty
,-n]\right) \in \mathcal{B}_{\infty }^{\prime }(\overline{\mathbb{R}}),
\end{equation*}

\bigskip \noindent and

\begin{equation*}
A\cup \{-\infty ,+\infty \}=\left( A\bigcup \{-\infty \}\right) \bigcup
\left( A\bigcup \{-\infty \}\right) .
\end{equation*}

\bigskip \noindent We may conclude that $\mathcal{B}_{\infty }(\overline{\mathbb{R}})\mathbb{\subset 
}\mathcal{B}_{\infty }^{\prime }(\overline{\mathbb{R}})$. We finally get
that $\mathcal{B}_{\infty }(\overline{\mathbb{R}})\mathbb{=}\mathcal{B}%
_{\infty }^{\prime }(\overline{\mathbb{R}})$.\\

\noindent We conclude that the class  $\mathcal{I}$ of all intervals of $\mathbb{R}$
still generates $\mathcal{B}_{\infty }(\overline{\mathbb{R}}).$ We may
easily adapt this proof to establish that each class given in \textit{Exercise 1} still
generates $\mathcal{B}_{\infty }(\overline{\mathbb{R}})$.\\

\bigskip \noindent \textbf{Question (c)}. Let $X:(\Omega ,\mathcal{A})\rightarrow (\mathbb{R},\mathcal{B}(%
\mathbb{R}))$ or $X:(\Omega ,\mathcal{A})\rightarrow (\overline{\mathbb{R}},%
\mathcal{B}_{\infty }(\overline{\mathbb{R}})).$ In both case, each class $%
\mathcal{I}(\circ )$ given in Exercise 1 generates both $\mathcal{B}(\mathbb{R}%
)$ and $\mathcal{B}_{\infty }(\overline{\mathbb{R}})$. We apply the weak
measurability function (WMC) to each of these classes to get the following measurability
criteria :\\

\noindent \textbf{Criterion (C1) based on }$\mathcal{I}(1)$ : $X$ is measurable iff for
any $a\in \mathbb{R},(X<a)$ is measurable.\\

\noindent \textbf{Criterion (C2) based on }$\mathcal{I}(3)$ : $X$ is measurable iff for
any $a\in \mathbb{R},(X\leq a)$ is measurable.\\

\noindent \textbf{Criterion (C3) based on }$\mathcal{I}(5)$ : $X$ is measurable iff for
any $a\in \mathbb{R},(X>a)$ is measurable.\\

\noindent \textbf{Criterion (C4) based on }$\mathcal{I}(7)$ : $X$ is measurable iff for
any $a\in \mathbb{R},(X\geq a)$ is measurable.\\ 

\noindent \textbf{Criterion (C5) based on }$\mathcal{I}(9)$ : $X$ is measurable iff for
any $(a,b)\in \mathbb{R}^{2}$, $(a<X\leq a)$ is measurable.\\ 

\noindent \textbf{Criterion (C6) based on }$\mathcal{I}(11)$ : $X$ is measurable iff
for any $(a,b)\in \mathbb{R}^{2}$, $(a<X<a)$ is measurable.\\

\noindent \textbf{Criterion (C7) based on }$\mathcal{I}(13)$ : $X$ is measurable iff
for any $(a,b)\in \mathbb{R}^{2}$, $(a\leq X<a)$ is measurable.\\

\noindent \textbf{Criterion (C8) based on }$\mathcal{I}(15).$ $X$ is measurable iff
for any $(a,b)\in \mathbb{R}^{2}$, $(a\leq X\leq a)$ is measurable.\\

\bigskip \noindent \textbf{Exercise 3}. \label{exercise03_sol_doc03-07}\\ \noindent Let $\{X_{n}\}_{n\geq 0}$ be a sequence of
measurable real-valued applications defined on the same measurable space : 
\begin{equation*}
X_{n}:(\Omega ,\mathcal{A})\rightarrow (\overline{\mathbb{R}},\mathcal{B}_{\infty }(\overline{\mathbb{R}})).
\end{equation*}

\bigskip \noindent (a) Show that the applications $\sup_{n\geq 0}X_{n}$, $\inf_{n\geq 0}X_{n},$ $\liminf_{n\rightarrow +\infty} X_{n}$ and $\limsup\liminf_{n\rightarrow +\infty} X_{n}$ are measurable.\\

\noindent (b) Show that, if the limit $\lim_{n\rightarrow +\infty} X_{n}(\omega )$ exists for all $\omega \in \Omega $, then the application $\lim_{n\rightarrow +\infty} X_{n}$ is measurable.\\

\bigskip \noindent \textbf{Solutions}. Let $\{X_{n}\}_{n\geq 0}$ be a sequence of measurable real-valued applications defined on the same measurable space : 
\begin{equation*}
X_{n}:(\Omega ,\mathcal{A})\rightarrow (\overline{\mathbb{R}},\mathcal{B}%
_{\infty }(\overline{\mathbb{R}})).
\end{equation*}

\bigskip \noindent \textbf{Question (a)}. Let us use the criteria found in \textit{Excercise 2}.\\

\noindent \textbf{(a1)}. $\sup_{n\geq 0}X_{n}$ is measurable since for an $a\in 
\mathbb{R}$,
\begin{equation*}
\left( \sup_{n\geq 0}X_{n}\leq a\right) =\bigcap\limits_{n\geq 1}\left(
X_{n}\leq a\right) \in \mathcal{A}.
\end{equation*}

\bigskip \noindent \textbf{(a.2)}. $\inf_{n\geq 1}X_{n}$ is measurable since for an $a\in 
\mathbb{R}$,
\begin{equation*}
\left( \inf_{n\geq 1}X_{n}\geq a\right) =\bigcap\limits_{n\geq 1}\left(
X_{n}\geq a\right) \in \mathcal{A}.
\end{equation*}

\bigskip \noindent \textbf{(a3)} $\lim \inf_{n\rightarrow +\infty }X_{n}$ is measurable by \textbf{(a1)} and
\textbf{(a2)} since 
\begin{equation*}
\liminf_{n\rightarrow +\infty }X_{n}=\sup_{n\geq 1}\left( \inf_{m\geq
n}X_{m}\right).
\end{equation*}

\bigskip \noindent (a4) $\limsup_{n\rightarrow +\infty }X_{n}$ is measurable by \textbf{(a1)} and
\textbf{(a2)} since 
\begin{equation*}
\lim \sup_{n\rightarrow +\infty }X_{n}=\inf_{n\geq 1}\left( \sup_{m\geq
n}X_{m}\right).
\end{equation*}

\bigskip \noindent \textbf{Question (b)}. The limit of the sequence $(X_n)_{n\geq 1}$, if it exists, is equal to the limit superior and to the limit inferior, and then, is measurable by Points \textbf{(a3)} and \textbf{(a4)} above.\\

\bigskip \noindent \textbf{Exercise 4}. \label{exercise04_sol_doc03-07} Treat the two following questions.\\

\noindent \textbf{(a)}. Show that constant function in are measurable.\\

\noindent \textbf{(b)}. Let $X$ be a measurable real-valued application : 
\begin{equation*}
X:(\Omega ,\mathcal{A})\rightarrow (\overline{\mathbb{R}},\mathcal{B}%
_{\infty }(\overline{\mathbb{R}})).
\end{equation*}

\bigskip \noindent and let $c\neq 0$ be a constant in $\mathbb{R}$. Show that $cX$ and $\left\vert X\right\vert $ are measurable.\\

\bigskip \noindent \textbf{Solutions}.\\

\noindent \textbf{Question (a)}. Let $X=c\in \overline{\mathbb{R}}.$ Use Criterion (C1) above.
We have for any $a\in \mathbb{R},$ 
\begin{equation*}
(X\leq a)=\left\{ 
\begin{tabular}{lll}
$\Omega$ & if & $c\leq a$ \\ 
$\emptyset $ & if & $c>a$.
\end{tabular}
\right. 
\end{equation*}

\bigskip \noindent Then $X=c$ is measurable.\\

\noindent \textbf{NB} : It may be useful to know that when $X=c$ is a constant function, we
have for any subset $B$ of 
\begin{equation*}
(X\in B)=\left\{ 
\begin{tabular}{lll}
$\Omega$ & if & $c\in B$ \\ 
$\emptyset $ & if & $c\notin B$.
\end{tabular}
\right.
\end{equation*}

\bigskip \noindent \textbf{Question (b)}.\\

\noindent \textbf{(b1)} For any $a\in \mathbb{R},$%
\begin{equation*}
(cX\geq a)=(X\geq a/c)
\end{equation*}

\bigskip \noindent if $c>0$ and

\begin{equation*}
(cX\geq a)=(X\leq a/c)
\end{equation*}

\bigskip \noindent if $c<0$. So, if $c\neq 0$, $cX$ is measurable if and only $X$ is.\\

\noindent \textbf{(b2)}. Let $X$ be measurable. We have that for any $a\in R$, 
\begin{equation*}
(\left\vert X\right\vert >a)=\left\{ 
\begin{tabular}{lll}
$(X<-a)\cup (X>a)$  & if & $a\geq 0$ \\ 
$\emptyset $ & if & $a<0$,
\end{tabular}
\right. 
\end{equation*}

\bigskip \noindent is measurable. Then $X$ is measurable.\\

\noindent \textbf{NB} : If $c=0$ and $X$ is finite (has its values in $\mathbb{R})$, $cX=0$ and it is
measurable because of Point (a). But if $X$ take infinite values, we may face undetermined forms $0\times \infty $. Later, we will see how to overcome such situations with the use of null sets.\\

\bigskip \noindent \textbf{Exercise 5}. \label{exercise05_sol_doc03-07} Define on $\mathbb{R}^{k}$ the class $\mathcal{I}_{k}$ of all
intervals 
\begin{equation*}
(a,b)=\prod\limits_{j=1}^{k}(a_{i},b_{i}),
\end{equation*}

\bigskip \noindent where $a=(a_{1},...,a_{k})\leq b=(b_{1},...,b_{k})$ [meaning that $a_{i}\leq
b_{i},1\leq i\leq k]$ . Define also 

\begin{equation*}
\mathcal{I}_{k}(1)=\{]a,b]=\prod\limits_{j=1}^{k}]a_{i},b_{i}],(a,b)\in 
\mathbb{R}^{k},a\leq b\},
\end{equation*}%
\begin{equation*}
\mathcal{I}_{k}(2)=\{]-\infty ,a]=\prod\limits_{j=1}^{k}]-\infty
_{i},a_{i}],a\in \overline{\mathbb{R}}^{k}\},
\end{equation*}

\noindent and 
\begin{equation*}
\mathcal{I}_{k}(3)=\{]a,b[=\prod\limits_{j=1}^{k}]a_{i},b_{i}[,(a,b)\in 
\mathbb{R}^{k},a\leq b\}
\end{equation*}

\bigskip \noindent \textbf{(a)}. Show that each of the classes $\mathcal{I}_{k}$ $\mathcal{I}_{k}(1)$, $\mathcal{I}_{k}(2)$ and $\mathcal{I}_{k}(2)$ generates the same $\sigma$-algebra. The $\sigma $-algebra generated by each of them is the usual $\sigma $-algebra on $\mathbb{R}^{k}$ denoted $\mathcal{B}(\mathbb{R}^{k})$.\\

\noindent \textbf{Hints and restrictions}. Restrict yourself to $k=2$ and show only $\sigma(\mathcal{I}_{2})=\sigma(\mathcal{I}_{2}(1))$. You only have to prove that each element of $\mathcal{I}_{2}$ can be obtained from countable sets operations on the elements of $\mathcal{I}_{2}(1))$. \textbf{Show this only for two cases} : $I=[a_{1},a_{2}[\times \lbrack a_{1},b_{2}[$ (a bounded interval) and $J=[a_{1},+\infty \lbrack \times ]-\infty ,b_{2}[$ (an unbounded interval). Express $I$ and $J$ as countable unions and intersections of elements of $\sigma(\mathcal{I}_{2}(1))$.

\noindent \textbf{(b)}. Show that $\mathcal{I}_{k}(2,c)=\{\overline{\mathbb{R}}^{k}\}\cup \mathcal{I}_{k}(2,c)$ is a semi-algebra.\\

\noindent \noindent \textbf{Hints and restrictions}. Sow this only for $k=2$.\\

\noindent \textbf{(b)} Show that $\mathcal{B}(\mathbb{R}^{k}))$ is the product $\sigma $-algebra corresponding to the product
space $(\mathbb{R},\mathcal{B}(\mathbb{R}))$ by itself $k$ times, that is ($\mathbb{R}^{k},\mathcal{B}(\mathbb{R}^{k}))=(\mathbb{R}^{k},\mathcal{B}(\mathbb{R})^{\otimes k})$.\\

\noindent \textbf{Hint} To show that $\mathcal{B}(\mathbb{R})^{\otimes k}) \subset  \mathcal{B}(\mathbb{R}^{k})$, use the characterization of 
$\mathcal{B}(\mathbb{R})^{\otimes k})$ by the measurability of the projections.  

\noindent \textbf{(d)} Show that $\mathcal{B}(\mathbb{R}^{k})$ is the Borel $\sigma $-algebra
on $\mathbb{R}^{k}$ endowed with one of the three equivalent metrics defined
for : $a=(a_{1},...,a_{k})$ and $b=(b_{1},...,b_{k})$
\begin{equation*}
d_{e}(a,b)=\sqrt{\sum_{j=1}^{k}(a_{i}-b_{i})^{2}}\text{,}
\end{equation*}%
\begin{equation*}
d_{m}(a,b)=\sum_{j=1}^{k}\left\vert a_{i}-b_{i}\right\vert \text{,}
\end{equation*}

\bigskip \noindent and
\begin{equation*}
d_{\infty }(a,b)=\max_{1\leq i\leq k}\left\vert a_{i}-b_{i}\right\vert.
\end{equation*}

\bigskip \noindent \textbf{NB}. Before you begin the solution, you are suggested to revise the product spaces DOC 00-14. See also Exercise 2 in DOC 00-02 and its solution in DIC 00-03.\\\

\noindent \textbf{Solutions}.\\

\noindent \textbf{Question (a)}. Usually, properties on $\mathbb{R}^{k}$ are straightforward
extensions of those in $\mathbb{R}^{2}$. So, we restrict ourselves to the
case $k=2.$ we have to show that each interval of $\mathbb{R}^{2}$ is in the 
$\sigma $-algebra generated by $\mathcal{I}_{k}(1)$ or $\mathcal{I}_{k}(2).$
Instead of giving all the possibilities as in the real case ($k=1$), we may
prove this for one bounded interval and one unbounded intervals. The others
cases are straightforward extensions. By the definition of $\mathcal{I}%
_{2}(1)$ par example, we have
\begin{equation*}
\{(a_{1},a_{2})\}=\bigcap\limits_{n\geq 1}]a_{1}-\frac{1}{n},a_{1}]\times
]a_{2}-\frac{1}{n},a_{2}]\in \sigma (\mathcal{I}_{2}(1)).
\end{equation*}

\bigskip \noindent Let $[a_{1},a_{2}[\times \lbrack a_{1},b_{2}[$ be a bounded interval of $%
\mathbb{R}^{k}$, we have
\begin{eqnarray*}
\lbrack a_{1},a_{2}[\times \lbrack b_{1},b_{2}[ &=&\bigcap\limits_{n\geq
1}]a_{1}-\frac{1}{n},a_{2}[\times \lbrack b-\frac{1}{n},b_{2}[ \\
&=&\bigcap\limits_{n\geq 1}\bigcup\limits_{m\geq 1}]a_{1}-\frac{1}{n}%
,a_{2}-\frac{1}{m}]\times \lbrack b-\frac{1}{n},b_{2}-\frac{1}{m}]\in 
\mathcal{I}_{2}(1).
\end{eqnarray*}

\bigskip \noindent Let $[a_{1},+\infty \lbrack \times ]-\infty ,b_{2}[$ be an unbounded
interval of $\mathbb{R}^{k}$. We have

\begin{eqnarray*}
\lbrack a_{1},+\infty \lbrack \times ]-\infty ,b_{2}[
&=&\bigcap\limits_{n\geq 1}[a_{1},n]\times ]-n,b_{2}[ \\
&=&\bigcap\limits_{n\geq 1}\bigcap\limits_{m\geq 1}]a_{1}-\frac{1}{m}%
,n]\times ]-n,b_{2}[ \\
&=&\bigcap\limits_{n\geq 1}\bigcap\limits_{m\geq 1}\bigcup\limits_{p\geq
1}]a_{1}-\frac{1}{m},n]\times ]-n,b_{2}-\frac{1}{p}]\in \mathcal{I}_{2}(1).
\end{eqnarray*}

\bigskip \noindent \textbf{Question (b)}. To show that $\mathcal{I}_{k}(2,c)=\{\overline{\mathbb{R}}%
^{k}\}\cup \mathcal{I}_{k}(2,c)$ is a semi-algebra, check the three points.\\

\noindent (b1) $\overline{\mathbb{R}}^{k}\in \mathcal{I}_{k}(2,c)$ by definition.\\

\noindent (b2) Let $A$ and $B$ in $\mathcal{I}_{k}(2,c)$. If one of them is $\overline{\mathbb{R}}^{k}$, $A\cap B$ is equal to the other and belongs to $\mathcal{I}_{k}(2,c)$. if not, we have an intersection of the form (for k=2), 
$$
(]a_{1}^{\prime },a_{2}^{\prime
}]\times ]b_{1}^{\prime },b_{2}^{\prime }])\cap (]a_{1}^{\prime \prime
},a_{2}^{\prime \prime }]\times ]b_{1}^{\prime \prime },b_{2}^{\prime \prime
}]),
$$ 

\bigskip \noindent which is
\begin{equation*}
(]a_{1}^{\prime },a_{2}^{\prime }]\cap ]a_{1}^{\prime \prime },a_{2}^{\prime
\prime }])\times (]b_{1}^{\prime },b_{2}^{\prime }]\cup ]b_{1}^{\prime
\prime },b_{2}^{\prime \prime }])=]a_{1}^{\prime }\vee a_{1}^{\prime \prime
},a_{2}^{\prime }\wedge a_{2}^{\prime \prime }]\times ]b_{1}^{\prime }\vee
b_{1}^{\prime \prime },b_{2}^{\prime }\wedge b_{2}^{\prime \prime }]\in 
\mathcal{I}_{k}(2).
\end{equation*}

\bigskip \noindent (b3) Let $A$ be in $\mathcal{I}_{k}(2,c)$. If $A=\mathbb{R}^k$, then $A^{c}=\emptyset
=]a,a]\times ]a,a]$ whatever be $a\in \mathbb{R}$. Then $A^{c}\in \mathcal{I}_{k}(2).$
If not, $A$ is of the form $]a_{1},a_{2}]\times ]b_{1},b_{2}]$ and 
\begin{eqnarray*}
A^{c} &=&(]-\infty ,a_{1}]+]a_{2},+\infty ])\times
]b_{1},b_{2}]+]a_{1},a_{2}]\times (]-\infty ,b_{1}]\\
&&+]b_{2},+\infty])+(]-\infty ,a_{1}]+]a_{2},+\infty ])\times (]-\infty
,b_{1}]+]b_{2},+\infty ]) \\
&=&]-\infty ,a_{1}]\times ]b_{1},b_{2}]+]a_{2},+\infty ]\times ]b_{1},b_{2}]
\\
&&+]a_{1},a_{2}]\times ]-\infty ,b_{1}]+]a_{1},a_{2}]\times ]b_{2},+\infty ]
\\
&&+]-\infty ,a_{1}])\times ]-\infty ,b_{1}]+]-\infty ,a_{1}]\times
]b_{2},+\infty ] \\
&&+]a_{2},+\infty ]\times ]-\infty ,b_{1}]+]a_{2},+\infty ]\times
]b_{2},+\infty ]).
\end{eqnarray*}

\bigskip \noindent \textbf{Question (b)}.\\

\noindent The class $\mathcal{I}_{k}$ which generates $\mathcal{B%
}(\mathbb{R}^{k})$ is a measurable rectangle in $(\mathbb{R}^{k}),\mathcal{B}%
(\mathbb{R})^{\otimes k})$. Then $\mathcal{B}(\mathbb{R}^{k})$ is included
in the $\sigma $-algebra generated by the measurable rectangles, that is $\mathcal{B}(\mathbb{R}^{k})\subset \mathcal{B}(\mathbb{R})^{\otimes k}.$ T
show the reverse inclusion, it is enough to show that the
projections 
\begin{equation*}
\pi _{i}:(\mathbb{R}^{k},\mathcal{B}(\mathbb{R}^{k}))\mapsto (\mathbb{R},%
\mathcal{B}(\mathbb{R})),1\leq i\leq k.
\end{equation*}

\bigskip \noindent are measurable (See Exercise 6 in DOC 02-02).\\

\noindent Let us check that for the first projection by using criteria (C1). We have
for any $a\in \mathbb{R}$,
\begin{eqnarray*}
(\pi _{1} &>&a)=]a,+\infty ]\times R^{k-1} \\
&=&\bigcup ]]a,n]\times ]-n,n]^{k-1}\in \sigma (\mathcal{I}_{k}(2))=%
\mathcal{B}(\mathbb{R}^{k}).
\end{eqnarray*}

\bigskip \noindent Since the projections are measurable and since $\mathcal{B}(\mathbb{R}%
)^{\otimes k}$ is the smallest $\sigma $-algebra making the projections
measurable, we conclude that $\mathcal{B}(\mathbb{R})^{\otimes k}\subset 
\mathcal{B}(\mathbb{R}^{k})$.\\

\noindent \textbf{Question (d)}. We suppose that we have already proved that $\mathcal{B}(%
\mathbb{R}^{k})$\ is generated by the class $\mathcal{I}_{k}(3).$ Let us use
the max-metric. An open ball centered at $x=(x_{1},...,x_{k})$ with radisu $%
r>0$ is $b=(b_{1},...,b_{k})$%
\begin{equation*}
B(x,r)=\{y=(y_{1},...,y_{k}),\max_{1\leq i\leq k}\left\vert
x_{i}-y_{i}\right\vert <r\}.
\end{equation*}

\bigskip \noindent We have
\begin{equation*}
B(x,r)=\prod\limits_{i=1}^{k}]x_{i}-r,x_{i}+r[.
\end{equation*}

\bigskip \noindent The class $\mathcal{I}_{k}(3)$ is formed by an open sets of the
topological space $\mathbb{R}^{k}$ endowed with the metric $d_{m}.$Then $%
\sigma (\mathcal{I}_{k}(3))\subset \mathcal{B}(\mathbb{R}^{k}).$ Conversely,
any open set is a union of balls centered at its points. Let $G$ be an open
set. Then we have

\begin{equation*}
G=\bigcup\limits_{x\in G}B(x,r_{x}).
\end{equation*}

\bigskip \noindent The space $\mathbb{R}^{k}$\ has the property that an arbitrary union is
equal to one its countable sub-unions. This implies that we may found a
countable subset $\{x^{_{(j)}},$ $j\in J,\neq (J)\leq Card(\mathbb{N})\}$\
of G such that\    
\begin{equation*}
G=\bigcup\limits_{j\in J}B(x^{_{(j)}},r_{x^{_{(j)}}})=\bigcup\limits_{j\in
J}\left(\prod
\limits_{i=1}^{k}]x_{i}^{_{(j)}}-r_{x^{_{(j)}}},x_{i}^{_{(j)}}+r_{x^{_{(j)}}}[\right) .
\end{equation*}

\bigskip \noindent Since the union is countable, we get $G\in \sigma (\mathcal{I}_{k}(3)).$ We
conclude that $\mathcal{B}(\mathbb{R}^{k})\subset \sigma (\mathcal{I}%
_{k}(3)).$ Finally, we have

\newpage
\noindent \LARGE \textbf{DOC 03-08 : Discover Exercises on Measurability 03 with solutions}. \label{doc03-08}\\

\Large
\bigskip

\noindent \textbf{Nota bene.}\newline

\noindent \textbf{(1)} In all this document, the measurable applications we
are using are defined on a measurable space $(\Omega ,\mathcal{A})$.\newline

\noindent \textbf{(2)} The expressions of the elementary functions are
always given with the full measurable partition of $\Omega$.\newline
 
\bigskip \noindent \textbf{Exercise 1}. \label{exercise01_sol_doc03-08} \\ \noindent Consider a measurable and
finite partition of $\Omega $, that is, $k$ mutually disjoint and measurable
subsets of $\Omega$, $A_{1}$,..., $A_{k}$ such that 
\begin{equation*}
\Omega =\sum_{1\leq i\leq k}A_{i}
\end{equation*}

\bigskip \noindent and let $\alpha_{1}$, $\alpha_{2}$, ..., and $\alpha_{k}$ be $%
k\geq 1$ finite real numbers. Let $X$ be the elementary function 
\begin{equation*}
X(\omega )=\alpha _{i}\text{ for }\omega \in A_{i},\text{ }1\leq i\leq k,
\end{equation*}

\bigskip \noindent also denoted by 
\begin{equation}
X=\sum_{1\leq i\leq k}\alpha _{i}1_{A_{i}}.  \label{etag01}
\end{equation}

\noindent \textbf{(a)} Suppose that the values $\alpha _{i},$ $1\leq i\leq k,
$ are distinct. Show that for any \noindent $i\in \{1,2,...,k\},$ 
\begin{equation*}
A_{i}=(X=\alpha _{i}).
\end{equation*}

\bigskip \noindent Deduce from this that if 
\begin{equation}
X=\sum_{1\leq j\leq m}\beta _{j}1_{B_{j}}  \label{etag02}
\end{equation}

\noindent is an expression of $X$ with distinct values $\beta _{j},$ $1\leq
j\leq \ell$, then, necessarily, $k=\ell $ and there exists a permutation $%
\sigma $ of $\{1,...,k\}$ such that for any $j\in \{1,...,k\}$

\begin{equation*}
\beta _{j}=\alpha _{\sigma (j)}\text{ and }B_{j}=A_{\sigma (j),}
\end{equation*}

\bigskip \noindent meaning that (\ref{etag02}) is simply a re-ordering of (\ref%
{etag01}).\newline

\noindent The expression (\ref{etag01}) of $X$ with distinct values of $%
\alpha _{i},$ is called the canonical expression of $X$.\newline

\noindent \textbf{(b)} Consider an arbitrary expression of an elementary
function : 
\begin{equation}
X=\sum_{1\leq i\leq k}\alpha _{i}1_{A_{i}}.
\end{equation}

\noindent Explain how to proceed to get the canonical form.\newline

\noindent \textbf{Hint} : Denote the distinct values of the sequence $\alpha
_{1},..,\alpha _{k} $ as $\beta _{1},...,\beta _{\ell }$ with $\ell \leq k.$
Set $I(j)=\{i,1\leq i\leq k,\alpha _{i}=\beta _{j}\}$ and%
\begin{equation*}
B_{j}=\sum\limits_{i\in I(j)}A_{i}.
\end{equation*}

\bigskip \noindent Conclude!\newline

\noindent \textbf{(c)} Use Questions (a) and (b) to prove that any
elementary function is measurable with the following suggested steps :%
\newline

\noindent \textbf{(c1)} Use the canonical form (\ref{etag01}) of $X.$ Show
that for any subset $B$ of $\mathbb{R}$,

\begin{equation*}
(X\in B)=\sum\limits_{i=1}^{k}(X\in B)\cap A_{i}.
\end{equation*}

\bigskip \noindent \textbf{(c2)} Use Question (a) to establish 
\begin{eqnarray*}
(X &\in &B)=\sum\limits_{i=1}^{k}(X\in B)\cap (X=\alpha _{i}) \\
&=&\sum\limits_{\alpha _{i}\in B}(X=\alpha _{i}).
\end{eqnarray*}

\bigskip \noindent \textbf{(c3)} Conclude that any elementary function $X$ is
measurable.\newline

\noindent \textbf{(d)} Use the same method to prove that for a real-valued
application $X$\ taking at most countable distinct values $\left(
x_{i}\right) _{i\in I}$ where $I$ is countable, we have for any subset of $%
\mathbb{R}$,
\begin{equation*}
(X\in B)=\sum\limits_{\alpha _{i}\in B}(X=x_{i}).
\end{equation*}

\bigskip \noindent \textbf{Conclude that} : A real-valued application $X$ taking at
most a countable number of distinct values $\left( x_{i}\right) _{i\in I}$
where $I$ is countable is measurable if and only if 
\begin{equation*}
\forall (i\in I),(X=x_{i})\text{ is measurable}.
\end{equation*}

\bigskip \noindent \textbf{(e)} Consider two expressions for an elementary function 
\begin{equation*}
X=\sum_{1\leq i\leq k}\alpha _{i}1_{A_{i}}
\end{equation*}

\bigskip \noindent and 
\begin{equation*}
X=\sum_{1\leq j\leq \ell }\beta _{j}1_{B_{j}}.
\end{equation*}

\noindent Find an expression of $X$ based on the superposition of the two
partitions, that is 
\begin{equation*}
\Omega =\sum_{1\leq i\leq k}\sum_{1\leq j\leq \ell }A_{i}B_{j}.
\end{equation*}

\bigskip \noindent \textbf{(f)} Consider two elementary functions 
\begin{equation*}
X=\sum_{1\leq i\leq k}\alpha _{i}1_{A_{i}}
\end{equation*}

\bigskip \noindent and 
\begin{equation*}
Y=\sum_{1\leq j\leq \ell }\beta _{j}1_{B_{j}}.
\end{equation*}

\bigskip \noindent Provide two expressions for $X$ and $Y$ based on the same
subdivision (\textbf{that} you take as the superposition of the partitions $%
(A_{i})_{1\leq i\leq k}$ and $(B_{i})_{\leq j\leq \ell })$.\newline

\bigskip \noindent \textbf{Solutions}.\newline

\noindent \textbf{Question (a)}. Let $X$ be the elementary function defined
by 
\begin{equation*}
X(\omega )=\alpha _{i}\text{ for }\omega \in A_{i},\text{ }1\leq i\leq k.
\end{equation*}

\bigskip \noindent Since, for $i$ fixed, $X(\omega )=\alpha _{i}$ for $\omega \in
A_{i},$ then $A_{i}\subset (X=\alpha _{i}).$ Now if the values $\alpha
_{1},...,\alpha _{k} $ are distinct, then for $i$ fixed, we cannot find $%
\omega \notin A_{i}$ such that $X(\omega )=\alpha _{i}$ since : 
\begin{equation*}
\omega \notin A_{i}\Longrightarrow \omega \in A_{i}^{c}=\sum\limits_{j\neq
i}A_{j},
\end{equation*}

\bigskip \noindent and then $X(\omega )$ is necessarily equal to one the $\alpha
_{j},1\leq j\neq i\leq k$ that are all different from $\alpha_i.$ In
conclusion for any fixed $i\in \{1,2,...,k\},$ $X(\omega )=\alpha _{i}$ if
and only if $\omega \in A_{i}.$ The first part of the question is answered.%
\newline

\bigskip \noindent As for the second part, suppose that 
\begin{equation*}
X=\sum_{1\leq j\leq \ell}\beta _{j}1_{B_{j}},
\end{equation*}

\noindent is an expression of $X$ with distinct values $\beta _{j},$ $1\leq
j\leq \ell.$ This means that $X$ takes the distinct values $\beta _{j},$ $%
1\leq j\leq \ell.$ Then necessarily $\ell=k$ and the sets $\{\alpha
_{i},1\leq i\leq k\}$ and equal $\{\beta _{j},1\leq j\leq k\}$.\newline

\noindent Therefore, the ordered set $(\alpha _{i},1\leq i\leq k)$ is a
permutation of the ordered set $(\beta _{i},1\leq i\leq k)$. This finishes
the answer of the question.\newline

\bigskip \noindent \textbf{Question (b)}. Let $X$ be the elementary function
defined by 
\begin{equation*}
X(\omega )=\alpha _{i}\text{ for }\omega \in A_{i},\text{ }1\leq i\leq k,
\end{equation*}

\noindent where the $\alpha _{1},...,\alpha _{k}$ are not necessarily
distinct and $A_{1},A_{2},...,A_{k}$ form a measurable partition of $\Omega .
$\ Then we regoup them into distinct values $\beta _{1},...,\beta _{\ell
},\ell \leq k$ so that each $\alpha _{i}$ is one of the $\beta _{j}.$ For
each $j\in \{1,...,\ell \}$, we may pick all the $\alpha _{i}$ equal to $%
\beta _{j}.$ Denote by $I(j)$ the set of the subscripts $i$ such that $%
\alpha _{i}=\beta _{j}$

\begin{equation*}
I(j)=\{i,1\leq i\leq k,\alpha _{i}=\beta _{j}\}.
\end{equation*}

\bigskip \noindent It is clear that the sets $I(1),I(2),...,I(\ell )$ form a
partition of $\{1,2,...,k\}.$ Next for each $j\in \{1,...,\ell \},$ we
regroup the sets $A_{i}$ for which $i\in I(j)$ in

\begin{equation*}
B_{j}=\sum\limits_{i\in I(j)}A_{i}.
\end{equation*}

\bigskip \noindent It is also clear that the $B_{1},B_{2},...,B_{\ell }$ form a
measurable partition of $\Omega ,$ and by construction 
\begin{equation*}
X=\beta _{j}\text{ on }B_{j},\text{ }j\in \{1,2,...,k\}.
\end{equation*}

\noindent Another writing of this is : 
\begin{equation*}
X=\sum_{1\leq j\leq \ell }\beta _{j}1_{B_{j}}.
\end{equation*}

\bigskip \noindent Since the values $\beta _{1},...,\beta _{\ell }$ are distinct, this latter expression is a canonical one and for any $j\in \{1,...,\ell \}$
\begin{equation*}
(X=\beta _{j})=B_{j}.
\end{equation*}

\bigskip \noindent \textbf{Question (c)}. Let us prove that an elementary function is measurable. We saw in Exercise 2 that any elementary function has a canonical expression in the form

\begin{equation*}
X(\omega )=\alpha _{i}\text{ for }\omega \in A_{i},\text{ }1\leq i\leq k,
\end{equation*}

\bigskip \noindent with $A_{i}=(X=\alpha _{i}),$ $1\leq i\leq k$ and then, since the $%
A_{i}$ form a partition of $\Omega ,$%
\begin{equation}
\Omega =\sum\limits_{i=1}^{k}(X=\alpha _{i}).  \label{etag0201}
\end{equation}

\noindent We want to prove that $X^{-1}(B)$ is measurable for any borel $B$
set in $\mathbb{R}$. But for any subset B of $\mathbb{R}$, even non
measurable,
\begin{eqnarray*}
X^{-1}(B) &=&(X\in B) \\
&=&(X\in B)\bigcap \Omega \\
&=&(X\in B)\bigcap \left( (X=\alpha _{i})\right) \\
&=&\sum\limits_{i=1}^{k}\left\{ (X\in B)\bigcap (X=\alpha _{i})\right\} 
\text{ (by \ref{etag0201}).}
\end{eqnarray*}

\bigskip \noindent Now, for each $i\in \{1,...,k\},$ we have two cases :\newline

\noindent \textbf{Case 1}. $\alpha _{i}\notin B.$ Then $(X\in B)\cap
(X=\alpha _{i})=\emptyset $ since we cannot find $\omega \in \Omega $ such
that $X(\omega )\in B$ and $X(\omega )=\alpha _{i}$ and $\alpha _{i}\notin B$%
.\newline

\noindent \textbf{Case 2}. $\alpha _{i}\in B.$ Then $(X\in B)\cap (X=\alpha
_{i})=(X=\alpha _{i})$ since : $\omega \in (X=\alpha _{i})\Longrightarrow
X(\omega )=\alpha _{i}\Longrightarrow X(\omega )\in B$ $($since$\ \alpha
_{i}\in B)\Longrightarrow \omega \in (X\in B).$ Thus $(X=\alpha _{i})\subset
(X\in B) $ and next $(X\in B)\cap (X=\alpha _{i})=(X=\alpha _{i})$.\newline

\noindent We summarize Case 1 and Case 2 into :

\begin{equation*}
(X\in B)\bigcap (X=\alpha _{i})=\left\{ 
\begin{tabular}{lll}
$\emptyset $ & if & $\alpha _{i}\notin B$ \\ 
$(X=\alpha _{i})$ & if & $\alpha _{i}\in B$.
\end{tabular}
\right. 
\end{equation*}

\noindent Going back five lines above, we get 
\begin{equation*}
X^{-1}(B)=\sum\limits_{\alpha _{i}\in B}(X=\alpha _{i})=\sum\limits_{\alpha
_{i}\in B}A_{i}.
\end{equation*}

\bigskip \noindent So $X^{-1}(B)$ is a finite sum of measurable sets.\newline

\noindent \textbf{Question (d)}. Suppose that $X$ takes a countable number
of distinct values $x_{i}$ $(i\in I).$ It is also clear that

\begin{equation*}
\Omega =\sum\limits_{i\in I}(X=x_{i})
\end{equation*}

\bigskip \noindent and for any Borel set $B,$ we have

\begin{equation}
X^{-1}(B)=\sum\limits_{i\in I,x_{i}\in B}(X=x_{i}).  \label{etagGen}
\end{equation}

\noindent We conclude as follows : If 
\begin{equation*}
\forall (i\in I)\text{, }(X=\alpha _{i})\text{ is measurable}
\end{equation*}

\bigskip \noindent then, by (\ref{etagGen}), $X$ is measurable. Conversely if $X$ is
measurable, then the $(X=x_{i})=X^{-1}(\{x_{i}\})$ are measurable.\newline

\bigskip \noindent \textbf{Question (e)}. Suppose that we have two
expressions for an elementary function

\begin{equation*}
X=\sum_{1\leq i\leq k}\alpha _{i}1_{A_{i}}
\end{equation*}

\bigskip \noindent and 
\begin{equation*}
X=\sum_{1\leq j\leq \ell }\beta _{j}1_{B_{j}}.
\end{equation*}

\bigskip \noindent By definition, $\Omega =\sum_{1\leq i\leq k}A_{i}$ and $\Omega
=\sum_{1\leq j\leq \ell }B_{j}$. Then

\begin{eqnarray*}
\Omega &=&\Omega \cap \Omega \\
&=&\left( \sum_{1\leq i\leq k}A_{i}\right) \bigcap \left( \sum_{1\leq j\leq
\ell }B_{j}\right) \\
&=&\sum_{1\leq i\leq k}\sum_{1\leq j\leq \ell }A_{i}B_{j}.
\end{eqnarray*}

\bigskip \noindent But some sets $A_{i}B_{j}$\ may be empty. So we restrict ourselves to the partition

\begin{equation}
\Omega =\sum_{(i,j)\in I}A_{i}B_{j},  \label{partSuper}
\end{equation}

\bigskip \noindent where 
\begin{equation*}
I=\{(i,j)\in \{1,...,k\}\times \{1,...,\ell \},A_{i}B_{j}\neq \emptyset \}.
\end{equation*}

\bigskip \noindent On this subdivision, called superposition of the subdivisions $%
\Omega =\sum_{1\leq i\leq k}A_{i}$ and $\Omega =\sum_{1\leq j\leq \ell }B_{j}
$, we have

\begin{equation*}
X=\alpha _{i}=\beta _{j}=\gamma _{ij}\text{ on }A_{i}B_{j},\text{ }(i,j)\in I
\end{equation*}

\bigskip \noindent and, in an other notation, 
\begin{equation*}
X=\sum_{(i,j)\in I}\gamma _{ij}1_{A_{i}B_{j}}.
\end{equation*}

\bigskip \noindent \textbf{Question (f)}. Suppose we have two elementary
functions 
\begin{equation*}
X=\sum_{1\leq i\leq k}\alpha _{i}1_{A_{i}}.
\end{equation*}

\bigskip \noindent and 
\begin{equation*}
Y=\sum_{1\leq j\leq \ell }\beta _{j}1_{B_{j}}
\end{equation*}

\bigskip \noindent Let us use the superposition of the subdivisions $\Omega
=\sum_{1\leq i\leq k}A_{i}$ and $\Omega =\sum_{1\leq j\leq \ell }B_{j}$
given in (\ref{partSuper}). We surely have

\begin{equation*}
X=\alpha _{i}\text{ on }A_{i}B_{j},\text{ }(i,j)\in I
\end{equation*}

\bigskip \noindent and 
\begin{equation*}
Y=\beta _{j}\text{ on }A_{i}B_{j},\text{ }(i,j)\in I
\end{equation*}

\bigskip \noindent corresponding to 
\begin{equation}
X=\sum_{(i,j)\in I}\alpha _{i}1_{A_{i}B_{j}}  \label{superExpCommonA}
\end{equation}

\noindent and 
\begin{equation}
Y=\sum_{(i,j)\in I}\beta _{j}1_{A_{i}B_{j}}.  \label{superExpCommonB}
\end{equation}

\bigskip \noindent \textbf{Exercise 2}. \label{exercise02_sol_doc03-08} (Algebra on the space of elementary
functions $\mathcal{E}$).\newline

\noindent Let $X$, $Y$ and $Z$ be simple functions, $a$ and $b$ be real
scalars.\newline

\noindent By using Question (c) of \textit{Exercise 1}, show the following
properies :\newline

\noindent \textbf{(a1)} $aX+bY\in \mathcal{E}$.\newline

\noindent \textbf{(a2)} $XY\in \mathcal{E}$.\newline

\noindent \textbf{(a3)} If $Y$ is everywhere different of zero, then $X/Y\in 
\mathcal{E}$.\newline

\noindent \textbf{(a5)} $(X\leq Y)\Rightarrow (X+Z\leq Y+Z)$. If $a\geq 0$,
then $X\leq Y\Rightarrow aX\leq bY$.\newline

\noindent \textbf{(a6)} $max(A,X)\in \mathcal{E}$, $min(X,Y)\in \mathcal{E}$.%
\newline

\bigskip \noindent \textbf{Terminology} : (You may skip the remarks below). 
\newline

\noindent \textbf{(A)} Needless to prove, the product of finite functions
and the sum of finite functions are commutative and associative. The product
of finite functions is distributive with the sum of finite functions. Taking
into account these facts, we state the following facts.\newline

\noindent By (a1), $(\mathcal{E},+,\cdot )$ is a linear space (where $\cdot$ stands for the external multiplication by real scalars).\newline

\noindent By (a1) and (a2), and the remarks (A), $(\mathcal{E},+,\times)$ is a ring with unit $1_{\Omega}$ (where $\times $ stands for the product between two functions), and is an algebra of functions.\newline

\noindent By (a6), $(\mathcal{E},\leq )$ is a lattice space, meaning closed under finite maximum and finite minimum.\newline

\noindent By (a1), (a5) and (a6), $(\mathcal{E},+,\cdot ,\leq )$ is a Riesz
space, that is a lattice (a6) vector space (a1) such that the order is
compatible with the linear structure.\newline

\bigskip \noindent \textbf{Solutions}.\newline

\noindent Let $X$, $Y$ be elementary functions, $a$ and $b$ be real scalars
and suppose that 
\begin{equation*}
X=\sum_{1\leq i\leq k}\alpha _{i}1_{A_{i}}
\end{equation*}

\bigskip \noindent and 
\begin{equation*}
Y=\sum_{1\leq j\leq \ell }\beta _{j}1_{B_{j}}.
\end{equation*}

\bigskip \noindent Let us consider the superposition of \ of the subdivisions $\Omega
=\sum_{1\leq i\leq k}A_{i}$ and $\Omega =\sum_{1\leq j\leq \ell }B_{j}$
given in (\ref{partSuper}) and the expressions of $X$ and $Y$ based on this
superposition as in (\ref{superExpCommonA}) and (\ref{superExpCommonB}). By
working of each element $A_{i}B_{j}$ of this partition, we have

\begin{equation*}
aX+bY=a\alpha _{i}+b\beta _{j}\text{ on }A_{i}B_{j},\text{ }(i,j)\in I,
\end{equation*}

\noindent that is 
\begin{equation*}
aX+bY=\sum_{(i,j)\in I}(a\alpha _{i}+b\beta _{j})1_{A_{i}B_{j},}
\end{equation*}

\bigskip \noindent and 
\begin{equation*}
XY=\alpha _{i}\beta _{j}\text{ on }A_{i}B_{j},\text{ }(i,j)\in I,
\end{equation*}

\bigskip \noindent that is 
\begin{equation*}
XY=\sum_{(i,j)\in I}\alpha _{i}\beta _{j}\text{ }1_{A_{i}B_{j},}
\end{equation*}

\bigskip \noindent and

\begin{equation*}
X/Y=\alpha _{i}/\beta _{j}\text{ on }A_{i}B_{j},\text{ }(i,j)\in I,\text{(if
all the }\beta _{j}\text{\ are different from zero),}
\end{equation*}

\bigskip \noindent that is

\begin{equation*}
X/Y=\sum_{(i,j)\in I}(\alpha _{i}/\beta _{j})\text{ }1_{A_{i}B_{j},}
\end{equation*}

\bigskip \noindent and 
\begin{equation*}
\max (X,Y)=\max (\alpha _{i},\beta _{j})\text{ on }A_{i}B_{j},\text{ }%
(i,j)\in I,
\end{equation*}

\bigskip \noindent that is 
\begin{equation*}
\max (X,Y)=\sum_{(i,j)\in I}\max (\alpha _{i},b_{j})\text{ }1_{A_{i}B_{j},}
\end{equation*}

\bigskip \noindent and 
\begin{equation*}
\min (X,Y)=\min (\alpha _{i},\beta _{j})\text{ on }A_{i}B_{j},\text{ }%
(i,j)\in I,
\end{equation*}

\bigskip \noindent that is

\begin{equation*}
\min (X,Y)=\sum_{(i,j)\in I}\min (\alpha _{i},b_{j})\text{ }1_{A_{i}B_{j}}.
\end{equation*}

\bigskip \noindent \textbf{Exercise 3}. \label{exercise03_sol_doc03-08} Let $X\geq 0$ be
measurable. For $n\geq 1$ fixed, consider the following partition of $\mathbb{%
R}_{+}$ :

\begin{equation*}
\mathbb{R}_{+}=\sum_{k=1}^{2^{2n}}[\frac{k-1}{2^{n}},\frac{k}{2^{n}}[\text{\ 
}+\text{ }[2^{n},+\infty ]=\sum_{k=1}^{2^{2n}+1}A_{k}.
\end{equation*}

\bigskip \noindent and set 
\begin{equation*}
X_{n}=\sum_{k=1}^{2^{n}}\frac{k-1}{2^{n}}\text{ }1_{(\frac{k-1}{2^{2n}}\leq
X<\frac{k}{2^{n}})}+2^{n}\text{ }1_{(X\geq 2^{2n})}.
\end{equation*}

\bigskip \noindent \textbf{(a)} Show that for each $n\geq $, $X_{n}$ is an elementary
function.\newline

\noindent \textbf{(b)} by using this implication : for any $\omega \in
\Omega ,$ for any $n\geq 1,$ 
\begin{equation*}
X(\omega )\geq \frac{k-1}{2^{n}},k\leq 2^{2n}
\end{equation*}

\bigskip \noindent implies 
\begin{equation*}
X_{n}(\omega )=\left\{ 
\begin{tabular}{lll}
$\frac{j-1}{2^{2n}}$ & if $\exists j : 1\leq j\leq $ $2^{2n}$ & and $\frac{j-1}{%
2^{2n}}\leq \frac{j-1}{2^{n}}\leq X(\omega )<\frac{j}{2^{n}}$ \\ 
&  &  \\ 
$2^{n}$ & otherwise & 
\end{tabular}
\right.,
\end{equation*}

\bigskip \noindent and by discussing the two cases ($X(\omega )\geq 2^{n}$ or $\frac{%
k-1}{2^{n}}\leq X(\omega )<\frac{k}{2^{n}}$ for some $k\leq 2^{2n}),$ show
that $X_{n}$ is non-decreasing.\newline

\noindent \textbf{(c)} Show that $X_{n}\rightarrow X$ as $n\rightarrow
+\infty $ by proceeding as follows. If $X(\omega )=+\infty ,$ show that $%
X_{n}(\omega )=2^{n}\rightarrow +\infty .$ If $\frac{k-1}{2^{n}}\leq
X(\omega )<\frac{k}{2^{n}}$ for some $k\leq 2^{2n},$ show that for $n$ large
enough 
\begin{equation*}
0\leq X(\omega )-X_{n}(\omega )\rightarrow \frac{1}{2^{n}}\rightarrow 0.
\end{equation*}

\bigskip \noindent \textit{\textbf{Conclusion} : Any non-negative measurable
applications is limit of a an non-decreasing sequence of non-negative
elementary functions.}

\bigskip \noindent \textbf{Solutions}.\newline

\noindent \textbf{Question (a)}.\newline

\noindent For any fixed $n\geq 0,$ the following dyadic partition of $%
\mathbb{R}_{+}$ 
\begin{equation*}
\mathbb{R}_{+}=\sum_{k=1}^{2^{2n}}[\frac{k-1}{2^{n}},\frac{k}{2^{n}}[\text{\ 
}+\text{ }[2^{n},+\infty ].
\end{equation*}

\bigskip \noindent is obvious. To make it more clear, we list its elements below for $%
k=1,2,....$, 

\begin{equation*}
\lbrack \frac{0}{2^{n}},\frac{1}{2^{n}}[,\text{ \ }[\frac{1}{2^{n}},\frac{2}{2^{n}}[,\text{ \ }[\frac{2}{2^{n}},\frac{3}{2^{n}}[, 
\end{equation*}

\begin{equation*}
..., \ [\frac{2^{2^{2n}}-1}{2^{n}},\frac{2^{2n}}{2^{n}}=2^{n}[, \ [2^{n},+\infty \lbrack.
\end{equation*}

\bigskip \noindent Define 
\begin{equation*}
A_{k}=X^{-1}\left( \left[ \frac{k-1}{2^{n}},\frac{k}{2^{n}}\right[ \right)
=\left( \frac{k-1}{2^{n}}\leq X<\frac{k}{2^{n}}\right) ,k=1,...,2^{2n}
\end{equation*}

\bigskip \noindent and 
\begin{equation*}
A_{2^{2n}+1}=X^{-1}([2^{n},+\infty \lbrack )=\left( X\geq 2^{2n}\right).
\end{equation*}

\bigskip \noindent We have a partition of $\Omega $ 
\begin{equation*}
\sum_{k=1}^{2^{2n}+1}A_{k}=\Omega .
\end{equation*}

\bigskip \noindent Then the application 
\begin{equation*}
X_{n}=\sum_{k=1}^{2^{n}}\frac{k-1}{2^{n}}\text{ }1_{(\frac{k-1}{2^{2n}}\leq
X<\frac{k}{2^{n}})}+2^{n}\text{ }1_{(X\geq 2^{2n})}.
\end{equation*}

\bigskip \noindent is an elementary function whenever $X$ is measurable. An other
expression of $X_{n}$ is

\begin{equation*}
X_{n}(\omega )=\left\{ 
\begin{tabular}{lll}
$\frac{k-1}{2^{n}}$ & if & $\frac{k-1}{2^{n}}\leq X(\omega )<\frac{k}{2^{n}},
$ $k=k=1,...,2^{2n}$ \\ 
$2^{n}$ & if & $X(\omega )\geq 2^{2n}$%
\end{tabular}%
\right. .
\end{equation*}

\bigskip \noindent (b) Let us show that the sequence ($X_{n})_{n\geq 1}$ is
non-increasing. Remark that for any $\omega \in \Omega ,$ for any $n\geq 1,$%
\begin{equation*}
X(\omega )\geq \frac{k-1}{2^{n}},k\leq 2^{2n}
\end{equation*}

\bigskip \noindent implies 
\begin{equation*}
X_{n}(\omega )=\left\{ 
\begin{tabular}{lll}
$\frac{j-1}{2^{n}}$ & if $ \exists j : 1\leq j\leq $ $2^{2n}$ & and $\frac{j-1}{%
2^{2n}}\leq \frac{j-1}{2^{n}}\leq X(\omega )<\frac{j}{2^{n}}$ \\ 
&  &  \\ 
$2^{n}$ & otherwise & 
\end{tabular}%
\right.  \ (RULE 1)
\end{equation*}

\bigskip \noindent Fix $n\geq 1$ and $\omega \in \Omega .$ We want to show that $%
X_{n}(\omega )\leq X_{n+1}(\omega ).$ Let us discuss two cases.\newline

\noindent \textbf{Case 1}. $X(\omega )\geq 2^{n}$. We have $X_{n}(\omega)=2^{n}$ and 
\begin{equation*}
X(\omega )\geq 2^{n}=\frac{2^{2n}}{2^{n}}=\frac{2^{2n+1}}{2^{n+1}}=\frac{%
(2^{2n+1}+1)-1}{2^{n+1}}=\frac{k-1}{2^{n+1}},
\end{equation*}

\bigskip \noindent where 
\begin{equation*}
k=2^{2n+1}+1\leq 2^{2(n+1)}\text{ \ (Check this easily).}
\end{equation*}

\bigskip \noindent Apply (RULE 1) for $n+1$ at the place of $n$ to get%
\begin{equation*}
X_{n+1}(\omega )=\left\{ 
\begin{tabular}{lll}
$\frac{j-1}{2^{n+1}}$ & if  $\exists j : 1\leq j\leq $ $2^{2(n+1)}$ & and $\frac{2^{2n+1}+1}{2^{n+1}}\leq \frac{j-1}{2^{n+1}}\leq X(\omega )<\frac{j}{2^{n+1}}$ \\ 
&  &  \\ 
$2^{n+1}$ & otherwise & 
\end{tabular}
\right..
\end{equation*}

\bigskip \noindent We see that 
\begin{equation*}
X_{n+1}(\omega )=\left\{ 
\begin{tabular}{l}
$\frac{j-1}{2^{n+1}}\geq \frac{2^{2n+1}+1}{2^{2n+1}}>\frac{2^{2n+1}}{2^{2n+1}%
}=2^{n}$ \\ 
or \\ 
$2^{n+1}$.
\end{tabular}
\right. 
\end{equation*}

\bigskip \noindent Then $X_{n+1}(\omega )\geq 2^{n}=X_{n}(\omega).$\newline

\noindent \textbf{Case 2}. We have for some $k,1\leq k\leq 2^{2n},$%
\begin{equation*}
X(\omega )=\frac{k-1}{2^{n}}\text{ and }\frac{k-1}{2^{n}}\leq X(\omega )<\frac{2}{2^{n}}.
\end{equation*}

\bigskip \noindent So we have 
\begin{equation*}
X(\omega )\geq \frac{2(k-1)}{2^{n+1}}=\frac{(2k-1)-1}{2^{n+1}}=\frac{K-1}{%
2^{n+1}},
\end{equation*}

\bigskip \noindent with 
\begin{equation*}
K=2k-1\leq 22^{2n}-1=2^{2n+1}-1\leq 2^{2(n+1).}
\end{equation*}

\bigskip \noindent We use again (RULE 1) to get 
\begin{equation*}
X_{n+1}(\omega )=\left\{ 
\begin{tabular}{lll}
$\frac{j-1}{2^{n+1}}$ & if $\exists  j: 1\leq j\leq $ $2^{2(n+1)}$ & and $\frac{K-1}{2^{n+1}}\leq \frac{j-1}{2^{n+1}}\leq X(\omega )<\frac{j}{%
2^{n+1}}$ \\ 
&  &  \\ 
$2^{n+1}$ & otherwise & 
\end{tabular}
\right. .
\end{equation*}

\noindent This implies 
\begin{equation*}
X_{n+1}(\omega )=\left\{ 
\begin{tabular}{l}
$\frac{j-1}{2^{n+1}}\geq \frac{K-1}{2^{n+1}}=\frac{2k-1-1}{2^{n+1}}=\frac{k-1%
}{2^{n}}=X_{n}(\omega ).$ \\ 
\\ 
$2^{n+1}=2\frac{2^{2n}}{2^{n}}\geq 2\frac{k}{2^{n}}=2X_{n}(\omega )$%
\end{tabular}%
\right. .
\end{equation*}

\bigskip \noindent Then for any $n \geq 1$, $X_{n+1}(\omega )\geq X_{n}(\omega)$.%
\newline

\noindent By putting together the two cases, we get that $n\geq 1$ and $%
\omega \in \Omega ,X_{n+1}(\omega )\geq X_{n}(\omega ).$ The sequence $X_{n}$
is non-decreasing.\newline

\noindent \textbf{Question (c)}. Let us consider two cases for a fixed $%
\omega \in \Omega$.\newline

\noindent \textbf{Case 1}. $X(\omega )=+\infty .$ Then for any $n\geq 1,$ $%
X(\omega )\geq 2^{n} $ and by the definition of $X_{n}$, 
\begin{equation*}
X_{n}(\omega )=2^{n}\rightarrow +\infty =X_{n}(\omega ).
\end{equation*}

\bigskip \noindent \textbf{Case 2}. $X(\omega )$ finite. Recall that $X$ is
non-negative. Since $2^{n}\nearrow +\infty $ and $n\nearrow +\infty ,$ then
for $n$ large enough (For example, take $n\geq (\log X(\omega))/\log 2)+1),$
we have 
\begin{equation*}
X(\omega )\in \lbrack 0,2^{n}[=\sum_{k=1}^{2^{2n}}[\frac{k-1}{2^{n}},\frac{k%
}{2^{n}}[.
\end{equation*}

\bigskip \noindent Then for $n$ large enough, there exists an integer $k(n,\omega )$
such that
\begin{equation*}
\left\{ 
\begin{tabular}{l}
$\frac{k-1}{2^{n}}\leq X(\omega )\leq \frac{k}{2^{n}}$ \\ 
$X_{n}(\omega )=\frac{k-1}{2^{n}}$%
\end{tabular}%
\right. .
\end{equation*}

\bigskip \noindent This implies, for large values of $n$, 
\begin{equation*}
0\leq X(\omega )-X_{n}(\omega )\leq \frac{1}{2^{n}}.
\end{equation*}

\bigskip \noindent This implies that 
\begin{equation*}
X_{n}(\omega )\rightarrow X(\omega )\text{ as }n\rightarrow \infty .
\end{equation*}

\bigskip \noindent \textbf{Exercise 4}. \label{exercise04_sol_doc03-08} Let $X$, $Y$ be measurable \textbf{%
finite} applications, $a$ and $b$ be real scalars.\newline

\noindent \textbf{(a)} Show that any real-valued measurable function $Z,$
possibly taking infinite values, is limit of elementary functions by
considering the decomposition%
\begin{equation*}
Z=\max (Z,0)-\max (-Z,0)
\end{equation*}

\bigskip \noindent where $Z^{+}=\max (Z,0)$ and $Z^{-}=\max (-Z,0)$ are non-negative
functions and are respectively called : positive and negative parts of $Z$.
\ By applying the conclusion of Question (c) of Exercise 3, show that any
real-valued measurable function is limit of elementary functions. \newline

\noindent \textbf{(b)} Deduce from this that $aX+bY$, $XY$, $\max (X,Y)$ and 
$\min (X,Y)$ are measurable. If $Y(\omega )\neq 0$ for all $\omega \in \Omega
$, then $X/Y$ is measurable.\newline

\bigskip \noindent \textbf{Solutions}.\newline

\noindent \textbf{Question (a)}. If $Z:(\Omega ,\mathcal{A})\rightarrow 
\mathbb{R}$ is measurable, then $Z^{+}=\max (Z,0)$ and $Z^{-}=\max (-Z,0)$
are measurable as superior envelops of two measurable functions. By \textit{%
Exercise 4}, $Z^{+}$ is limit of a sequence of elementary functions $%
(g_{n}^{+})_{n\geq 1}$ and $Z^{-}$ is limit of a sequence of elementary
functions $(g_{n}^{-})_{n\geq 1}.$ Then $Z=Z^{+}-Z^{-}$ is limit of the
sequence $g_{n}=g_{n}^{+}-g_{n}^{-},$ which is a sequence of elementary
functions.\newline

\bigskip 

\noindent \textbf{Question (b)}. Let $X$, $Y$ be measurable \textbf{finite}
applications, $a$ and $b$ be real scalars. Then all the operations $aX+bY,$ $%
XY$, $\max (X,Y)$ and $\min (X,Y)$ are well-defined. By Question (a), we may
find two sequences of elementary functions $(X_{n})_{n\geq 1}$ and $%
(Y_{n})_{n\geq 1}$ such that

\begin{equation*}
X_{n}\rightarrow X\text{ and }Y_{n}\rightarrow Y,\text{ as }n\rightarrow
+\infty .
\end{equation*}

\bigskip \noindent Then, for each $n\geq 1,$ the functions $aX_{n}+bY_{n},$ $%
X_{n}Y_{n}$, $\max (X_{n},Y_{n})$ and $\min (X_{n},Y_{n})$ are elementary
functions and they respectively converge to $aX+bY,$ $XY$, $\max (X,Y)$ and $%
\min (X,Y).$ By Exercise 3 in Doc 03-04, the functions $aX+bY,$ $XY$, $\max
(X,Y)$ and $\min (X,Y)$ are measurable.\newline

\noindent If $Y(\omega )\neq 0$ for all $\omega \in \Omega ,$ we may change
each $Y_{n} $ to 
\begin{equation*}
\widetilde{Y}_{n}=\left\{ 
\begin{tabular}{lll}
$Y_{n}$ & if & $Y_{n}\neq 0$ \\ 
1 & if & $Y_{n}=0$%
\end{tabular}%
\right. .
\end{equation*}

\bigskip \noindent Each $\widetilde{Y}_{n}$ is still an elementary function and does
not vanish on $\Omega $. Since for any $\omega \in \Omega ,Y(\omega )\neq 0$
and $Y_{n}(\omega )\rightarrow Y(\omega ),$ it follows that $Y_{n}(\omega
)\neq 0$ when $n$ is large enough and then $\widetilde{Y}_{n}(\omega
)=Y_{n}(\omega )$ for large enough values of $n$ and then

\begin{equation*}
\widetilde{Y}_{n}(\omega )\rightarrow Y(\omega )
\end{equation*}

\bigskip \noindent and then 
\begin{equation*}
\mathcal{E} \ni X_{n}/\widetilde{Y}_{n}\rightarrow X/Y,
\end{equation*}

\bigskip \noindent which is then measurable.\newline

\bigskip \noindent \textbf{Exercise 5}. \label{exercise05_sol_doc03-08} Consider a product space endowed
with its product $\sigma $-algebra $\prod_{1\leq i\leq k}\Omega
_{i},\bigotimes_{1\leq i\leq k}\mathcal{A}_{i},$ and let  $X$ $\ $be a
function defined on from this product space to  in $\overline{\mathbb{R}}$ :
\ 
\begin{equation*}
\begin{tabular}{lllll}
$X$ & $:$ & $\left( \prod_{1\leq i\leq k}\Omega _{i},\bigotimes_{1\leq i\leq
k}\mathcal{A}_{i}\right) $ & $\longmapsto $ & $(\overline{\mathbb{R}},%
\mathcal{B}_{\infty }(\overline{\mathbb{R}}))$ \\ 
&  & $(\omega _{1,}\omega _{2,}...,\omega _{k})$ & $\hookrightarrow $ & $X(\omega _{1,}\omega _{2,}...,\omega _{k})$.
\end{tabular}
\end{equation*}

\bigskip \noindent For any fixed $i\in \{1,2,...\}$ , define the  $i-$th partial function $X$ by considering it as a function of $\omega _{i}$ only and by
fixing the other variables. Let $(\omega _{1}, ... ,\omega_{i-1},\omega _{i+1},...,\omega _{k})$ fixed and define $X^{(i)}$ the as $i-$%
th partial function $X$ 
\begin{equation*}
\begin{tabular}{lllll}
$X^{(i)}$ & $:$ & $\left( \Omega _{i},\mathcal{A}_{i}\right) $ & $%
\longmapsto $ & $(\overline{\mathbb{R}},\mathcal{B}_{\infty }(\overline{%
\mathbb{R}}))$ \\ 
&  & $\mathbf{\omega }_{i}$ & $\hookrightarrow $ & $X^{(i)}(t)=f(\omega
_{1,}\omega _{2,}...,\omega _{i-1},\mathbf{\omega }_{i},\omega
_{i+1},...,\omega _{k})$.
\end{tabular}.
\end{equation*}

\bigskip \noindent The objective of this exercise is to show that that partial functions of
real-valued and measurable functions of several variables are also
measurable. The proof for $k=2$ may re conducted word by word to a general $%
k\geq 2.$ In the sequel $k=2.$ Let us fix $\omega _{1}\in \Omega _{1}$ and
consider the partial function%
\begin{equation*}
\begin{tabular}{lllll}
$X$ & $:$ & $\left( \Omega _{2},\mathcal{A}_{2}\right) $ & $\longmapsto $ & $%
(\overline{\mathbb{R}},\mathcal{B}_{\infty }(\overline{\mathbb{R}}))$ \\ 
&  & $\mathbf{\omega }_{2}$ & $\hookrightarrow $ & $X^{(i)}=X_{\omega
_{1}}(\omega _{2})=X(\omega _{1},\omega _{2})$.
\end{tabular}
\end{equation*}

\bigskip \noindent \textbf{(a)} Show that for any non-decreasing sequence of applications $X_{n}:(\Omega
_{1}\times \Omega _{2},\mathcal{A}_{i}\otimes \mathcal{A}_{2})\longmapsto (%
\overline{\mathbb{R}},\mathcal{B}_{\infty }(\overline{\mathbb{R}})),$ the
sequence of partial applications $\left( X_{n}\right) _{\omega _{1}}$ is
also non decreasing.\\

\noindent \textbf{(b)} Show that for any sequence of applications $X_{n}:(\Omega _{1}\times
\Omega _{2},\mathcal{A}_{i}\otimes \mathcal{A}_{2})\longmapsto (\overline{%
\mathbb{R}},\mathcal{B}_{\infty }(\overline{\mathbb{R}})),n\geq 1,$
converging to the application $X_{n}:(\Omega _{1}\times \Omega _{2},\mathcal{%
A}_{i}\otimes \mathcal{A}_{2})\longmapsto (\overline{\mathbb{R}},\mathcal{B}%
_{\infty }(\overline{\mathbb{R}})),$ the sequence of partial applications $%
\left( X_{n}\right) _{\omega _{1}},n\geq 1,$ is also non-decreasing also
converges to the partial function $X_{\omega _{1}}$ of $X$.\\

\noindent \textbf{(b)} Show for $Y=1_{A}:(\Omega _{1}\times \Omega _{2},\mathcal{A}_{i}\otimes 
\mathcal{A}_{2})\longmapsto (\overline{\mathbb{R}},\mathcal{B}_{\infty }(%
\overline{\mathbb{R}}))$ is the indicator function of a set $A\subset \Omega
_{1}\times \Omega _{2},$ then the partial function%
\begin{equation*}
Y_{\omega _{1}}(\omega _{2})=Y(\omega _{1},\omega _{2})=1_{A}(\omega
_{1},\omega _{2})=1_{A_{\omega _{1}}}(\omega _{2}),
\end{equation*}

\bigskip \noindent is the indicator function of the section $A_{\omega _{1}}$ of $A$ at $\omega_{1}$\\

\noindent \textbf{(c)} Deduce from Question (b) than partial functions of indicator functions
of measurable sets are measurable and partial functions of elementary
functions defined on $(\Omega _{1}\times \Omega _{2},\mathcal{A}_{i}\otimes 
\mathcal{A}_{2})\longmapsto (\overline{\mathbb{R}},\mathcal{B}_{\infty }(%
\overline{\mathbb{R}}))$ are measurable.\\

\noindent \textbf{(d)} Let \ \ \ $X:(\Omega _{1}\times \Omega _{2},\mathcal{A}_{i}\otimes 
\mathcal{A}_{2})\longmapsto (\overline{\mathbb{R}},\mathcal{B}_{\infty }(%
\overline{\mathbb{R}}))$ be a measurable application. Use Question (a) of
Exercise 4 and consider a sequence $X_{n}:(\Omega _{1}\times \Omega _{2},%
\mathcal{A}_{i}\otimes \mathcal{A}_{2})\longmapsto (\overline{\mathbb{R}},%
\mathcal{B}_{\infty }(\overline{\mathbb{R}})),n\geq 1,$ of elementary
functions converging to $X.$ Combine the previous questions to prove that $%
X$ has measurable partial functions.\\

\bigskip \noindent \textbf{Solutions}.\\

\noindent \textbf{Question (a)}. If the sequence $X_{n}:(\Omega _{1}\times \Omega _{2},%
\mathcal{A}_{i}\otimes \mathcal{A}_{2})\longmapsto (\overline{\mathbb{R}},%
\mathcal{B}_{\infty }(\overline{\mathbb{R}}))$ is non-decreasing, that%
\begin{equation*}
\forall (\omega _{1},\omega _{2})\in \Omega _{1}\times \Omega _{2},\text{ }%
\forall (n\geq 1),\text{ }X_{n+1}(\omega _{1},\omega _{2})\leq X_{n}(\omega
_{1},\omega _{2}).
\end{equation*}

\bigskip \noindent Then for $\omega _{1}\in \Omega _{1}$, we also have
\begin{eqnarray*}
&&\forall (\omega _{2}\in \Omega _{2}),\text{ }\forall (n\geq 1),\text{ }%
(X_{n+1})_{\omega _{1}}(\omega _{2})=X_{n+1}(\omega _{1},\omega _{2})\\
&\leq& X_{n}(\omega _{1},\omega _{2})=(X_{n})_{\omega _{1}}(\omega _{2}).
\end{eqnarray*}

\bigskip \noindent And the sequence of partial functions $(X_{n})_{\omega _{1}}$ is
non-decreasing.\\

\noindent \textbf{Question (b)}. If the sequence $X_{n}:(\Omega _{1}\times \Omega _{2},%
\mathcal{A}_{i}\otimes \mathcal{A}_{2})\longmapsto (\overline{\mathbb{R}},%
\mathcal{B}_{\infty }(\overline{\mathbb{R}}))$ converges to $X:(\Omega
_{1}\times \Omega _{2},\mathcal{A}_{i}\otimes \mathcal{A}_{2})\longmapsto (%
\overline{\mathbb{R}},\mathcal{B}_{\infty }(\overline{\mathbb{R}})),$ then 
\begin{equation*}
\forall (\omega _{1},\omega _{2})\in \Omega _{1}\times \Omega _{2},\text{  }%
X_{n+1}(\omega _{1},\omega _{2})\longrightarrow X(\omega _{1},\omega _{2})%
\text{ }as\text{ }n\rightarrow +\infty.
\end{equation*}

\bigskip \noindent Then for $\omega _{1}\in \Omega _{1}$, we also have, as $n\rightarrow +\infty$, 
\begin{eqnarray*}
&&\forall (\omega _{2}\in \Omega _{2}),\text{ }(X_{n+1})_{\omega _{1}}(\omega
_{2})=X_{n+1}(\omega _{1},\omega _{2})\\
&\longrightarrow& X_{n}(\omega_{1},\omega _{2})=(X_{n})_{\omega _{1}}(\omega _{2}).
\end{eqnarray*}

\bigskip \noindent And the sequence of partial functions $(X_{n})_{\omega _{1}}$ converges to $%
(X_{n})_{\omega _{1}}$.\\

\noindent \textbf{Question (c)}. Let $A\subset \Omega _{1}\times \Omega _{2}$ and $\omega
_{1}\in \Omega _{1}.$ We want to show that%
\begin{equation}
\left( 1_{A}\right) _{\omega _{1}}=1_{A_{\omega _{1}}}.  \label{partialIndic}
\end{equation}

\noindent Let us compare the graphs. For any $\omega _{2}\in \Omega _{2}$%
\begin{equation*}
\left( 1_{A}\right) _{\omega _{1}}=1_{A}(\omega _{1},\omega _{2})=\left\{ 
\begin{tabular}{lll}
1 & if & $(\omega _{1},\omega _{2})\in A$ \\ 
0 & if & $(\omega _{1},\omega _{2})\notin A$,
\end{tabular}%
\right. 
\end{equation*}

\bigskip \noindent and
\begin{equation*}
1_{A_{\omega _{1}}}(\omega _{2})=\left\{ 
\begin{tabular}{lll}
1 & si & $\omega _{2}\in A_{\omega _{1}}\Longleftrightarrow (\omega
_{1},\omega _{2})\in A$ \\ 
0 & si & $\omega _{2}\notin A_{\omega _{1}}\Longleftrightarrow (\omega
_{1},\omega _{2})\notin A$.
\end{tabular}
\right. 
\end{equation*}

\bigskip \noindent By comparing these two graphs, we get the our target.\\

\noindent Now let $Y$ be an elementary function, that is 
\begin{equation*}
Y=\sum\limits_{j=1}^{k}\alpha _{j}1_{A_{j}}.
\end{equation*}

\bigskip \noindent with the usual assumptions of the $\alpha_i$ and the $A_j$. By applying the result of indicator functions above, we get for $\omega _{1}\in \Omega _{1}$,

\begin{eqnarray}
Y_{\omega _{1}}(\omega _{2}) &=&Y(\omega _{1},\omega _{2})
\label{partialElem} \\
&=&\sum\limits_{j=1}^{k}\alpha _{j}1_{A_{j}}(\omega _{1},\omega _{2}) 
\notag \\
&=&\sum\limits_{j=1}^{k}\alpha _{j}1_{\left( A_{j}\right) _{\omega
_{1}}}(\omega _{2}).  \notag
\end{eqnarray}

\noindent \textbf{Question (c)}. By (\ref{partialIndic})  and (\ref{partialElem})
by the knowledge that sections of measurable functions are measurable (Doc
01-06, Exercise 3), it becomes clear that partial functions of indicators
functions anf of elementary functions are measurable.\\

\noindent \textbf{Question (d)}.  Let $X:(\Omega _{1}\times \Omega _{2},\mathcal{A}_{i}\otimes \mathcal{A}_{2})\longmapsto (\overline{\mathbb{R}},\mathcal{B}_{\infty }(\overline{\mathbb{R}}))$ be a real-valued measurable application.\\

\noindent By \textit{Question (a) of Exercise 4}, there exists a sequence $X_{n}:(\Omega _{1}\times \Omega _{2},%
\mathcal{A}_{i}\otimes \mathcal{A}_{2})\longmapsto (\overline{\mathbb{R}},
\mathcal{B}_{\infty }(\overline{\mathbb{R}})),n\geq 1,$ of elementary
functions converging to $X.$ Let $\omega _{1}\in \Omega _{1}$ be fixed. By
\textit{Question (c)}, the partial functions $(X_{n})_{\omega _{1}}$ are measurable
and by \textit{Question (a)},
\begin{equation*}
(X_{n})_{\omega _{1}}\rightarrow X_{\omega _{1}}\text{ }n\rightarrow +\infty.
\end{equation*}

\bigskip \noindent Then the partial function $X_{\omega _{1}}$ is measurable as a limit of measurable functions.

\newpage
\noindent \LARGE \textbf{DOC 03-09 :Discover Exercises on Measurability 03 with solutions}. \label{doc03-09}\\
\Large
\bigskip

\bigskip \noindent \textbf{Exercise 1}. \label{exercise01_sol_doc03-09} Le $f$ \ be a non-decreasing function from $\mathbb{R}$ to $\mathbb{R}$. Show that the set of discontinuity points
is at most countable.\newline

\bigskip \noindent \emph{Proceed to this way}. \noindent Let $D(n)$ the number discontinuity points of $f$ in $]-n,n[$. Recall that a real number $x$ is discontinuity point of a non-decreasing function $f$ if and only if
\begin{equation*}
f(x+)-f(x-)>0
\end{equation*}

\bigskip \noindent where $f(x-)$ is the left-hand limit of $f$ at $x$ and $f(x+)$ is the right-hand limit. In that case $f(x+)-f(x-)$ is the discontinuity jump. Next let 
\begin{equation*}
D_{k}(n)=\{x\in D(n),\text{ }f(x+)-f(x-)>1/k\}.
\end{equation*}

\bigskip \noindent Pick $x_{1},...,x_{m}$ from $D_{k}(n)$ and justify the inequality 

\begin{equation*}
\sum_{1\leq i\leq m}f(x_{i}+)-f(x_{i}-)\leq f(n)-f(-n).
\end{equation*}

\bigskip \noindent Deduce from this that 
\begin{equation*}
m\leq k\times (f(n)-f(-n)), 
\end{equation*}

\bigskip \noindent and then $D_{k}(n)$ is finite.\\

\noindent Conclude.\newline

\bigskip \noindent \textbf{Solution}. \label{exercise02_sol_doc03-09}  Since $f$ is non-decreasing, $x
$ is a discontinuity point of $f$ if and only if the discontinuity jump $%
f(x+)-f(x-)$ is positive. Denote by $D$ the set of all discontinuity points
of $f$, and for any $k\geq 1$, denote by $D_{k}$ the set of discontinuity
points such that $f(x+)-f(x-) > 1/k$ and by $D_{k,n}$ the set of
discontinuity points in the interval $]-n,n[$ such that $f(x+)-f(x-) > 1/k$.
We are going to show that $D_{k,n}$ is finite.\newline

\noindent Let us suppose we can find $m$ points $x_1$, ..., $x_m$ in $D_{k,n}
$. Since $f$ is non-decreasing, we may see that the sum of the discontinuity
jumps is less than $f(n)-f(-n)$. [ You may make a simple drawing for $m=3$ and
project the jumps to the y-axis to see this easily]. So

\begin{equation*}
\sum_{1\leq j \leq m} f(x+)-f(x-) \leq f(n)-f(-n). 
\end{equation*}

\bigskip \noindent Since each of these jumps exceeds $1/k$, we have

\begin{equation*}
\sum_{1\leq j \leq m} (1/k) \leq \sum_{1\leq j \leq m} f(x+)-f(x-) \leq
F(n)-F(-n). 
\end{equation*}

\noindent

\begin{equation*}
m/k \leq f(n)-f(-n),
\end{equation*}

\bigskip \noindent That is

\begin{equation*}
m \leq k(f(n)-f(-n)) 
\end{equation*}

\bigskip \noindent We conclude by saying that we cannot have more that $[k(f(n)-f(-n))]$ points in $D_{k,n}$, so $D_{k,n}$ is finite. Since

\begin{equation*}
D=\cup_{n\geq 1} \cup_{k\geq 1} D(k,n),
\end{equation*}

\bigskip \noindent we see that $D$ is countable. The solution is complete.\\

\bigskip \noindent \textbf{Exercise 2}. \label{exercise02_sol_doc03-09}  Show that any right-continuous or left-continuous function $X$ from $\mathbb{R}$ to $\mathbb{R}$ is measurable.\newline

\noindent \textit{Hint} : Let $X$ be right-continuous and set for all $n\geq 1$ 

\begin{equation*}
X_{n}^{+}(t)=\sum_{k=-\infty }^{k=+\infty }X\left( \frac{k+1}{2^{n}}\right) 1_{\left] \frac{k}{2^{-n}},\frac{k+1}{2^{-n}}\right] }(t).
\end{equation*}

\bigskip \noindent and for $N\geq 1$

\begin{equation*}
X_{n,N}^{+}(t)=\sum_{k=-N}^{k=N}X\left( \frac{k+1}{2^{n}}\right) 1_{\left] \frac{k}{2^{-n}},\frac{k+1}{2^{-n}}\right] }(t).
\end{equation*}

\bigskip \noindent $X_{n,N}^{+}$ is measurable because it is elementary function (refer to the document on elementary functions).\\

\bigskip \noindent Show that $X_{n,N}^{+}$ converges to $X_{n}^{+}$ as $N\rightarrow +\infty$ ($n$ being fixed).\newline

\bigskip \noindent Conclude that each $X_{n}$ is measurable and that $X_{n}$ converges to $X$ as $n \rightarrow +\infty$. Make a final conclude.\newline

\noindent \textit{Hint} : If $X$ be left-continuous and set for all $n\geq 1$ 
\begin{equation*}
X_{n}^{-}(t)=X_{n}^{+}(t)=\sum_{k=-\infty }^{k=+\infty }X\left( \frac{k}{2^{n}}\right) 1_{\left[ \frac{k}{2^{-n}},\frac{k+1}{2^{-n}}\right[}(t).
\end{equation*}

\bigskip \noindent Next proceed similarly to the case of right-continuous functions.\\

\noindent \textbf{Solution}. \\

\noindent \textbf{(1) Solution for a right-continuous function}.\\

\noindent The function 
\begin{equation*}
X_{n,N}^{+}(t)=\sum_{k=-N}^{k=N}X\left( \frac{k+1}{2^{n}}\right) 1_{\left] 
\frac{k}{2^{-n}},\frac{k+1}{2^{-n}}\right] }(t),
\end{equation*}

\bigskip \noindent is an elementary function. By Exercise 1 in DOC 03-08, Question (c3), you have proved that $X_{n,N}^{+}$ is measurable. And it is obvious evident that for each 
$t\in \mathbb{R}$,

\begin{equation*}
X_{n,N}^{+}(t)\rightarrow X_{n}^{+}(t)\text{ as }N\rightarrow +\infty,
\end{equation*}

\bigskip
\noindent and then $X_{n}^{+}$ is measurable by Exercise 3, Question (b) in DOC 03-07. Now fix $t\in \mathbb{R}.$ For each $n\geq 1,$ we have this partition of $t\in \mathbb{R},$

\begin{equation*}
\mathbb{R}=\sum_{k=-\infty }^{+\infty }\left] \frac{k}{2^{n}},\frac{k+1}{%
2^{n}}\right] .
\end{equation*}

\bigskip \noindent Then there exists a unique value $k=k(n,t)$ such that,%
\begin{equation*}
\left[ \frac{k(n,t)}{2^{-n}},\frac{k(n+t)+1}{2^{-n}}\right[ =\left[ x_{n}(t)-%
\frac{1}{2^{-n}},x_{n}(t)\right[, 
\end{equation*}

\noindent where
\begin{equation*}
x_{n}(t)=\frac{k+1}{2^{n}}.
\end{equation*}

\bigskip \noindent Hence, by definition,
\begin{equation*}
\left\{ 
\begin{tabular}{l}
$X_{n}^{+}(t)=X(x_{n}(t))$ \\ 
$0\leq x_{n}(t)-t<2^{-n}$%
\end{tabular}%
\right. 
\end{equation*}

\bigskip \noindent Since $x_{n}(t)\rightarrow t$ as $n\rightarrow +\infty $\ from the right and
since $X$ is right-continuous, then 
\begin{equation*}
X_{n}^{+}(t)=X(x_{n}(t))\rightarrow X(t)\text{ as }n\rightarrow +\infty.
\end{equation*}

\bigskip
\noindent Hence $X$ is limit of the sequence of measurable Borel functions $X_{n}^{+}$. Then $X$ is
measurable.\\

\noindent \textbf{(2) Solution for the left-continuous function}. We use the function  
\begin{equation*}
X_{n}^{-}(t)=\sum_{k=-\infty }^{k=+\infty }X\left( \frac{k}{2^{n}}\right) 1_{%
\left[ \frac{k}{2^{-n}},\frac{k+1}{2^{-n}}\right[ }(t),
\end{equation*}

\bigskip \noindent which is limit (as $N\rightarrow +\infty )$ of the elementary functions 

\begin{equation*}
X_{n,N}^{-}(t)=\sum_{k=-N}^{k=N}X\left( \frac{k}{2^{n}}\right) 1_{\left[ 
\frac{k}{2^{-n}},\frac{k+1}{2^{-n}}\right[ }(t).
\end{equation*}

\bigskip \noindent To show that $X_{n}^{-}$ converges to $X$, we remark that for
\begin{equation*}
x_{n}(t)=\frac{k}{2^{n}}\text{ for }.t\in \left[ \frac{k}{2^{-n}},\frac{k+1}{%
2^{-n}}\right[,
\end{equation*}

\bigskip \noindent we have

\begin{equation*}
\left\{ 
\begin{tabular}{l}
$X_{n}^{-}(t)=X(x_{n}(t))$ \\ 
$0\leq t-x_{n}(t)<2^{-n}$.
\end{tabular}
\right. 
\end{equation*}

\bigskip
\noindent From there, the left-continuity leads to $X_{n}^{-}\rightarrow X$ and the
conclusion is made.

\bigskip \noindent \textbf{Exercise 3}. \label{exercise03_sol_doc03-09}  Let $X$ be an application from from $\mathbb{R}$ to $\mathbb{R}$ such that at each $t\in \mathbb{R}$ the left limit $f(x-)$ exists (resp. at each $t\in \mathbb{R}$ the right-hand limit $f(x+)$ exists.\\
 
\noindent \textbf{(a)} Use know facts on the left continuity of $t\mapsto X(t-)$ and the left-continuity of $(x\mapsto X(t+))$) to show that $X$ is measurable with help of Exercise 2.\\

\noindent \textit{(b)} Suppose that at each $t\in \mathbb{R}$ the left limit $f(x-)$ exists (resp. at each $t\in \mathbb{R}$ the right-hand limit $f(x+)$ exists). Show that the function  $t\mapsto X(t-)$ is left-continuous (resp. $(x\mapsto X(t+))$ is right-continuous). [You may skip this question and consider it as a topology exercise. But you are recommended to read the proof].\\

\bigskip \noindent \textbf{Solution}.\\

\noindent \textbf{Question (a)}. If for each $t\in \mathbb{R},$ the right-hand limit $X(t+)$ exists, then by results of real analysis, the function $t\mapsto X(t+)$ is right-continuous. Hence, by Exercise 1, it is measurable.\\

\noindent As well, if for each $t\in \mathbb{R},$ the left-hand limit $X(t-)$ exists,
then by results of real analysis, the function $t\mapsto X(t-)$ is left-continuous. Hence, by Exercise 1, it is measurable.\newline

\noindent \textbf{Question (b)}. Let us prove the above mentioned properties of real analysis. We give the solution only for the right-continuity case.\\

\noindent Assume that the left-limit $X(t+)$ exists at each $t\in \mathbb{R}$. For any fixed $t$, we have

\begin{equation*}
(\forall \epsilon >0), (\exists h_0 >0) : 0 < h<h_0 \Rightarrow,
X(t+h)-\epsilon \leq X(t+) \leq X(t+h)+\epsilon. \text{ (LM) }
\end{equation*}

\bigskip \noindent Now, set $Y(t)=X(t+)$. We have to prove that $Y$ is
right-continuous. Fix $t$ and $\epsilon >0$. From (LM), we have

\begin{equation*}
(\forall 0 < h<h_{0}/2, \forall 0 < \ell<h_{0}/2), X(t+h+\ell)-\epsilon \leq
X(t+) \leq X(t+h+\ell)+\epsilon. \text{ (RC1) }
\end{equation*}

\bigskip \noindent Then let $\ell \downarrow 0$, we get

\begin{equation*}
(\forall 0 <h<h_{0}/2, X((t+h)+)-\epsilon \leq X(t+) \leq
X((t+h)+)+\epsilon. \text{ (RC1) }
\end{equation*}

\bigskip \noindent Through (RC1), we have proved this : For any $\epsilon
>0 $, there exists $h_1=h_{0}/2$ such that

\begin{equation*}
(\forall 0 <h<h_{1}, Y(t+h)-\epsilon \leq Y(t) \leq Y(t+h)+\epsilon. \text{ (RC2) }
\end{equation*}

\bigskip 
\noindent The right-continuity of $Y(t)=X(t+)$ is thus proved.\newline

\noindent To treat the case where $X$ has left-limits, do similarly by taking $%
Z(t)=X(t-)$. Be wise!\newline

\bigskip \noindent \textbf{Exercise 4}. \label{exercise04_sol_doc03-09} Let $X$ be an application $\mathbb{R}$ dans $\mathbb{R}$ such that the set of discontinuity points $D$ of $X$  is
at most countable and $X$ admits right-hand limits at all points (or $X$ admits  left-hand limits at all points ). Justify the identity 
\begin{equation*}
\forall t\in \mathbb{R}, X(t)=X(t+)\times 1_{D^{c}}(t)+\sum_{t\in D}X(t)\times 1_{\{t\}}(t).
\end{equation*}

\bigskip \noindent and conclude that $X$ is measurable.\\

\bigskip \noindent \textbf{Solution}. Suppose that $X$ admits right-hand limits at all points
and that the set of discontinuity points $D$ of $X$ is at most countable.
Denote $D=\{s_{1},s_{2},...\}$. It is clear that
\begin{equation*}
\forall t,X(t)=X(t)1_{D^{c}}(t)+X(t)1_{D}(t).
\end{equation*}

\bigskip \noindent But, we have $X(t+)=X(t)$ on $D^{c}$ and then, we surely get

\begin{equation*}
X(t)1_{D}(t)=\sum_{j\geq 1}X(s_{j})1_{\{s_{j}\}}(t),
\end{equation*}

\bigskip \noindent which is either an elementary function (if $D$ is finite) or a limit of the
sequences of elementary functions

\begin{equation*}
\sum_{j=1}^{N}X(s_{j})1_{\{s_{j}\}}(t),
\end{equation*}

\bigskip \noindent as $N\rightarrow +\infty $ (if $D$ is infinite), and is measurable. Finally

\begin{equation*}
\forall t,X(t)=X(t+)\times 1_{D^{c}}(t)+\sum_{j\geq
1}X(s_{j})1_{\{s_{j}\}}(t),
\end{equation*}

\bigskip \noindent is measurable a product of measurable functions plus a measurable one (See Exercise 2 in Doc 03-08).\\

\bigskip \noindent \textbf{Exercise 5}. \label{exercise05_sol_doc03-09}  Use Exercises 1 and 4 to show that any monotone function $f:\mathbb{R}\rightarrow \mathbb{R}$ is 
measurable.\\

\bigskip \noindent \textbf{Solution}. Let $f:\mathbb{R}\rightarrow \mathbb{R}$ be monotone. Suppose that it is
non-decreasing. From Exercise 1, the set $D$ of discontinuity points of $f$ is at most
countable. Surely $f$ has a left and right limit at each point $t\in \mathbb{%
R}.$ Indeed, fix $t\in \mathbb{R}.$ Since $f$ is non-decreasing, the set%
\begin{equation*}
A_{t}=\{F(x),x<t\}
\end{equation*}

\bigskip \noindent is bounded above by $f(t)$. So by the Archimedean property of $\mathbb{R}$, this set has a finite supremum
\begin{equation*}
\alpha =\sup A_{t}.
\end{equation*}

\noindent S Let us show that $\alpha $ is the left limit $f(t-)$ of $f$ at $t$. First, we have : for any $h>0$
\begin{equation*}
f(t-h)\leq f(t).
\end{equation*}

\noindent S As $h\downarrow 0,$ $f(t-h)$ increases in the large sense. Then $f(t-h)$  has a limit $f(t-)$ which is still less than $f(t)$.
 Since $f$ is non-decreasing,the increasing limit $as$ $h\downarrow 0$ is exactly the limit as $h\rightarrow 0$, which is equal to $f(t-)$.\\

\noindent Conversely, by the definition of the supremum, for any $\varepsilon
_{n}=1/n$, $n\geq 1$, there exists one element $A_{t},$ denoted $f(y_{n})$ [
with $x_{n}<t$ ] that exceeds $\alpha -\varepsilon _{n}$ that is%
\begin{equation*}
f(y_{n})\geq \alpha -\varepsilon _{n}.
\end{equation*}

\noindent S Now either $t-\varepsilon _{n}\leq y_{n}<t$ and we put $x_{n}=y_{n}$ or$%
y_{n}<$ $t-\varepsilon _{n}$ and we take $x_{n}=t-\varepsilon _{n}$ and the
non-decreasingness of $f,$ we have%
\begin{equation*}
\left\{ 
\begin{tabular}{l}
$f(x_{n})\geq \alpha -\varepsilon _{n}$ \\ 
$t-\varepsilon _{n}\leq x_{n}<t$.
\end{tabular}
\right. 
\end{equation*}

\bigskip \noindent Thus $x_{n}\rightarrow t$ by left, we get by letting $n\rightarrow +\infty .$%
\begin{equation*}
f(t-)\geq \alpha .
\end{equation*}

\bigskip \noindent We conclude that 
\begin{equation*}
f(t-)=\sup A_{t}.
\end{equation*}

\bigskip
\noindent A non-decreasing function has left-hand limits everywhere and has at most
number of discontinuity points. It is measurable.\\

\noindent If $f$ is non-increasing, we apply the result we found to $-f$ which is non-decreasing. The function $-f$ is measurable, so is $f$.\\

\bigskip \noindent \textbf{Exercise 6}. \label{exercise06_sol_doc03-09}  Let $X$ be a convex application $%
\mathbb{R}$ to $\mathbb{R}$.\\

\noindent \textbf{(a)} Use a know fact on the right-continuity of a convex function (as reminded in Question (b) below) and conclude with the help of Exercise 2.\\

\noindent \textbf{(b)} Show that any convex function is continuous. [You may skip this question and consider it as a topology exercise. But you are recommended to read the proof].\\

\bigskip \noindent \textbf{Solution}.\\

\noindent \textbf{Question (a)}. By known result of real analysis reminded by Question (b), a convex function is right-continuous. It is then measurable by Exercise 2.\\

\noindent \textbf{Question (b)}. Let $X:\mathbb{R}\rightarrow \mathbb{R}$ be convex. We are going to prove that $X$ is right-continuous.\newline

\noindent By definition, $X$ is convex if and only if

\begin{equation*}
(\forall (s,t) \in R^2, s\leq t), (\forall 1\leq \lambda \leq 1), X(\lambda
s + (1-\lambda)t) \leq \lambda X(s) + (1-\lambda)X(t).
\end{equation*}

\bigskip \noindent We begin to prove this formula :

\begin{equation*}
\forall s < t < u, \frac{X(t)-X(s)}{t-s} \leq \frac{X(u)-X(s)}{u-s} \leq \frac{X(u)-X(t)}{u-t}.  \ \ (FCC1)
\end{equation*}

\bigskip \noindent we may check that for $s < t < u$,

\begin{equation*}
t=\frac{u-t}{u-s} s + \frac{t-s}{u-s} u = \lambda s + (1-\lambda)t,
\end{equation*}

\bigskip \noindent where $\lambda=(u-t)/(u-s) \in ]0,1[$ and $1-\lambda=(t-s)/(u-s)$.
By convexity, we have

\begin{equation*}
X(t)\leq \frac{u-t}{u-s} X(s) + \frac{t-s}{u-s} X(u). \text{ (FC1) }
\end{equation*}

\bigskip
\noindent Multiply all members of (FC1) by $(u-s)$ to get

\begin{equation*}
(u-s) X(t)\leq (u-t) X(s) + (t-s) X(u). \text{ (FC1) }
\end{equation*}

\bigskip
\noindent First, split $(u-t)$ into $(u-t)=(u-s)-(t-s)$ in the first term in the left-hand member to have

\begin{eqnarray*}
&&(u-s) X(t)\leq (u-s) X(s) - (t-s) X(s) +(t-s) X(u)\\
&\Rightarrow& (u-s) (X(t)-X(s)) \leq (t-s) (X(s)-X(u)).
\end{eqnarray*}

\bigskip \noindent This leads to

\begin{equation*}
\frac{X(t)-X(s)}{t-s} \leq \frac{X(u)-X(s)}{u-s}. \text{ (FC2a) }
\end{equation*}

\bigskip \noindent Next, split $(t-s)$ into $(t-s)=(u-s)-(u-t)$ in the second term in the left-hand member to have

\begin{eqnarray*}
&&(u-s) X(t)\leq (u-t) X(s) + (u-s) X(u) - (u-t) X(u)\\
&\Rightarrow& (u-s) (X(t)-X(u)) \leq (u-t) (X(s)-X(u)).
\end{eqnarray*}

\bigskip
\noindent Multiply the last inequality by $-1$ to get

\begin{equation*}
(u-s) (X(u)-X(t)) \geq (u-t) (X(u)-X(s)),
\end{equation*}

\bigskip \noindent and then

\begin{equation*}
\frac{X(u)-X(s)}{u-s} \leq \frac{X(u)-X(t)}{u-t}.  \ \  (FC2b)
\end{equation*}

\bigskip
\noindent Together (FC2a) and (FC2b) prove (FCC).\newline

\noindent Before we go further, we think it is better to write (FCC) as :

\begin{equation*}
\forall t_1 < t_2 < t_3, \frac{X(t_2)-X(t_1)}{t_2 - t_1} \leq \frac{X(t_3)-X(t_1)}{t_3-t_1} \leq \frac{X(t_3)-X(t_2)}{t_3-t_2}. \ \ (FCC)
\end{equation*}

\bigskip
\noindent We also may apply (FCC) to :\newline

\begin{equation*}
\forall r < s < t, \frac{X(s)-X(r)}{r-s} \leq \frac{X(s)-X(r)}{s-r} \leq \frac{X(t)-X(s)}{s-t}.  \ \ (FCC)
\end{equation*}

\bigskip
\noindent We are on the point to conclude. Fix $r < s$. From (FFC1), we may
see that the function 
\begin{equation*}
G(v)=\frac{X(v)-X(t)}{v - t}, v >t,
\end{equation*}

\bigskip
\noindent is increasing since in (FCC1), $G(t) \leq G(u)$ for $t<u$. By (FCC2), $G(v)$ is bounded below by $G(r)$ (that is fixed with $r$ and $t$).\\

\noindent We conclude $G(v)$ decreases to a real number as $v \downarrow t$
and this limit is the roght-derivative of $X$ at $t$ :

\begin{equation*}
\lim_{v \downarrow t} \frac{X(v)-X(t)}{v - t}=X^{\prime}_r(t).\text{ (RL) }
\end{equation*}

\bigskip \noindent Since $X$ has a right-derivative at each point, it is
right-continuous at each point since, by (RL),

\begin{equation*}
X(t+h)-X(t)= h (G(t+h)=h(X^{\prime}_r(t)) + o(1)) \ \ as \ \ h\rightarrow 0.
\end{equation*}

\bigskip
\noindent \textbf{Remark}. You also may prove that the left-derivative
exists in a similar way and conclude that $X$ is left-continuous and
finally, $X$ is continuous.\newline

\bigskip \noindent \textbf{Exercise 7}. \label{exercise07_sol_doc03-09}  Let $X$ be a function from a metric
space $(E,d)$ to $\mathbb{R}$. Suppose that if $X$ is upper semi-continuous or lower semi-continuous.\\ 

\noindent \textbf{Reminder}.\\

\noindent \textit{(1)} A function $X:(E,d)\rightarrow \mathbb{R}$ is lower-continuous if and
only for each $t\in E$,

\begin{equation*}
(\forall 0<\epsilon <0)(\exists \eta >0),d(s,t)<\eta \Rightarrow X(s)\leq X(t)+\epsilon. \text{ (LSC) }
\end{equation*}

\noindent \textbf{(2)} A function $X:(E,d)\rightarrow \mathbb{R}$ is upper semi-continuous if and only for each $t\in E$,

\begin{equation*}
(\forall 0<\epsilon <0)(\exists \eta >0),d(s,t)<\eta \Rightarrow X(s)\geq X(t)-\epsilon. \text{ (LSC) }
\end{equation*}

\bigskip
\noindent (a) Use know facts on semi-continuous functions in Topology, as reminded in Question (b), to show that $X$ is measurable.\\

\noindent (b) Suppose that $X$ is lower semi-continuous. Show that for  $c\in \mathbb{R%
},(X\leq c)$ is closed. Deduce that if $X$ is upper semi-continuous, then
for any $c\in \mathbb{R},(X\geq c)$ is closed. [You may skip this question and consider it as a toplogy exercise. But you are recommended to read the proof].\\

\bigskip \noindent \textbf{Solution}.\\

\noindent \textbf{Question (a)}. By the reminder given in the statement of Question (b), a lower
semi-continuous function satisfies : for $c\in \mathbb{R},(X\leq c).$ Then $X
$ is measurable by the criteria of Exercise 2, Question (c).\\

\noindent Next, an upper semi-continuous function is the opposite of an lower semi-continuous. Then it is measurable.\\

\noindent \textbf{Question (b)}. Let $X$ be lower semi-continuous. Let us how that  for any
real number $c$, the set $F=(X\leq c)$ is closed.\newline

\noindent $X:(E,d)\rightarrow \mathbb{R}$ is lower-continuous at $t$ means : for any

\begin{equation*}
(\forall 0<\epsilon <0)(\exists \eta >0),d(s,t)<\eta \Rightarrow X(s)\leq X(t)+\epsilon. \text{ (LSC) }
\end{equation*}

\bigskip
\noindent Now let $c$ be fixed. Let show that $G=F^c=(X>c)$ is open, and then measurable. Let $t \in F^c=(X>c)$, that is $X(t)>c$. Let us
exhibit an open ball centered at $t$ that is in $G$. Set $\epsilon=X(t)-c>0$
and use (LSC). We finish by checking that $B_d(t,\eta)\subset G$, and this comes from,

\begin{eqnarray*}
s \in B_d(t,\eta) &\Rightarrow& X(s)\leq X(t)+\epsilon > X(t) + X(t)-c=c\\
&\Rightarrow& X(s)>c \Rightarrow s\in G.
\end{eqnarray*}

\bigskip
\noindent So, for real number $c$, $(X\leq c)$ is closed, is then
measurable. Finally $X$ is measurable.\newline

\noindent Now, if $f$ is upper semi-continuous. So $-f$ is upper semi-continuous, and then for any $c \in \mathbb{R}$, the set $(-f \leq -c)$ is closed, that is,  $(f \geq c)$ is closed.\\

\newpage
\bigskip \noindent her we semi-continuous functions in special forms. It is proved in http://dx.doi.org/10.16929/jmfsp/2017.001, the equivalence of other more classical definitions.
We will find these expressions powerful to deal with Lebesgue-Stieljes integrals.\\

\bigskip \noindent \textbf{Exercise 8}. \label{exercise08_sol_doc03-09}  Let $f$ be a function $\mathbb{R}\rightarrow \mathbb{R}.$ Let $x\in \mathbb{R%
}$. Define the limit superior and the limit inferior of $f$ at $x$ by%
\begin{equation*}
f^{\ast }(x)=\lim \sup_{y\rightarrow x}f(y)=\lim_{\varepsilon \downarrow
0}\sup \{f(y),y\in ]x-\varepsilon ,x+\varepsilon \lbrack \}
\end{equation*}

\bigskip \noindent and
 
\begin{equation*}
f_{\ast }(x)=\lim \inf_{y\rightarrow x}f(y)=\lim_{\varepsilon \downarrow
0}\inf \{f(y),y\in ]x-\varepsilon ,x+\varepsilon \}.
\end{equation*}

\bigskip \noindent Question (1). Show that $f^{\ast }(x)$ and $f_{\ast }(x)$ always exist.\\

\noindent \textit{Hint} : remark that $\sup \{f(y),y\in ]x-\varepsilon ,x+\varepsilon \}$ is non-increasing and $\inf \{f(y),y\in ]x-\varepsilon ,x+\varepsilon \}$ is non-decreasing as 
$\varepsilon $ decreases to zero.\\

\bigskip \noindent Question (2). Show that \ $(-f)^{\ast }=-(f_{\ast })$.\\

\noindent Question (3). Show that $f$ is continuous at $x$ if and only if $f^{\ast }(x)$ and $f_{\ast }(x)$, and then : $f$ is continuous if and only if $f=$ $f^{\ast}=f_{\ast}$.\\

\noindent Question (4) Define : $f$ is upper semi-continuous if and only if at $x$ if and only if $f=f^{\ast }$ and $f$ is lower semi-continuous iff $f=f_{\ast }$.\\

\noindent Show that $f^{\ast}$ is upper semi-continuous. Exploit (2) to extend this to $f_{\ast}$ and show that $f_{\ast }$ is lower semi-continuous.\\

\bigskip \noindent Question (5) Now let $f$ be upper semi-continuous. Set for each $n\geq 1$%
\begin{equation*}
f_{n}(x)=\sum_{k=-\infty }^{k=+\infty }\sup \{f(z),z\in [k2^{-n},(k+1)2^{-n}[\} 1_{]k2^{-n},(k+1)2^{-n}[}(x)+\sum_{s\in D_n}f(s)1_{\{s\}}.
\end{equation*}

\bigskip \noindent Show that $f_{n}$ converges to $f$. Conclude that $f$ is measurable.\\

\bigskip  \noindent Question (6). Use the opposite argument to show that a lower semi-continuous function is measurable.\\

\newpage
\noindent \textbf{Solution of Exercise 8}.\\

\noindent Question (1). 
\begin{equation*}
f^{\ast }(x)=\lim \sup_{y\rightarrow x}f(y)=\lim_{\varepsilon \downarrow
0}\sup \{f(y),y\in ]x-\varepsilon ,x+\varepsilon \lbrack \}
\end{equation*}

\noindent But
\begin{equation*}
f_{\varepsilon }^{\ast }(x)=\sup \{f(y),y\in ]x-\varepsilon ,x+\varepsilon
\lbrack \}
\end{equation*}

\bigskip \noindent is decreasing as $\varepsilon \downarrow 0$ and is bounded below by $f(x).$
Then its infimum\ 
\begin{equation*}
\min_{\varepsilon >0}f_{\varepsilon }^{\ast }(x),
\end{equation*}

\bigskip \noindent is its limit as $\varepsilon \downarrow 0.$ So $f^{\ast }(x)$ exists. So
does $f_{\ast }(x)$ since $(-f)^{\ast }=-(f_{\ast }).$\\

\bigskip \noindent Question (2).\\

\noindent The solution is obvious and then omitted.\\

\bigskip  \noindent Question (3). Suppose $f$ is continuous.\\

\noindent  Take an arbitrary sequence $\varepsilon
_{n}\downarrow 0$ as $n\uparrow +\infty .$ We are going to prove that $f_{\varepsilon _{n}}^{\ast }(x)=\sup \{f(y),y\in ]x-\varepsilon _{n},x+\varepsilon _{n}[\}.$ Since $f$ is continuous at $x$, there exists $\eta >0$ that $y\in ]x-\varepsilon _{n},x+\varepsilon _{n}[\Longrightarrow f(x)-1<f(y)<f(x)+1.$\\ 

\noindent Then for $n$ large enough so that $\varepsilon
_{n}<\eta ,$ $f_{\varepsilon _{n}}^{\ast }(x)$ is finite. We may use the
characterization of the supremum in R : For any $\theta >0,$ there exists $%
y_{n}\in ]x-\varepsilon _{n},x+\varepsilon _{n}[$ such that%
\begin{equation*}
f_{\varepsilon _{n}}^{\ast }(x)\leq f(y_{n})<f_{\varepsilon _{n}}^{\ast }(x).
\end{equation*}

\bigskip
\noindent Let $n\uparrow +\infty $ and use that $f(y_{n})\rightarrow f(x)$ since $%
y_{n}\rightarrow x$ and $f$ is continuous at $x,$ to have for all $\theta
>0; $
\begin{equation*}
f^{\ast }(x)\leq f(x)\leq f^{\ast }(x)+\theta .
\end{equation*}

\bigskip \noindent Now let $\theta \downarrow 0$ to have $f^{\ast }(x)=f(x).$ To get $f^{\ast
}(x)=f(x),$ use that $-f$ \ is continuous at $x$ if and only if $f$ is
continuous at $x$. Since $-f$ $\ $is continuous at $x$, use our first
conclusion to get that $(-f)^{\ast }=-f(x),$ that is, by Point 2, $-f_{\ast
}(x)=-f(x)$ $\ $and then $f_{\ast }(x)=f(x).$\\

\noindent Now suppose that $f^{\ast }(x)=$ $f_{\ast }(x).$ Since for all $\varepsilon
>0$, 
\begin{equation*}
f_{\ast _{,}\varepsilon }(x)=\inf \{f(y),y\in ]x-\varepsilon ,x+\varepsilon
\lbrack \}\leq f(x)\leq f_{\varepsilon }^{\ast }(x),
\end{equation*}

\bigskip
\noindent we get by letting $\varepsilon \downarrow 0$
\begin{equation*}
f^{\ast }(x)=f_{\ast }(x)=f(x).
\end{equation*}

\bigskip \noindent Now let $(y_{n})_{n\geq 1}$ be an arbitrary sequence such that $y_{n}\rightarrow
x$ as $n\rightarrow +\infty .$So for any $\varepsilon >0,$ there exists $%
n_{0}\geq 1$ such that

\begin{equation*}
n\geq n_{0}\Longrightarrow ]x-\varepsilon ,x+\varepsilon \lbrack,
\end{equation*}

\bigskip \noindent and then
\begin{equation*}
n\geq n_{0}\Longrightarrow \inf \{f(y),y\in ]x-\varepsilon ,x+\varepsilon
\lbrack \}\leq f(y_{n})\leq \sup \{f(y),y\in ]x-\varepsilon ,x+\varepsilon
\lbrack \}.
\end{equation*}

\bigskip \noindent Thus, we let $n\longrightarrow +\infty $ to get
\begin{eqnarray*}
&&\inf \{f(y),y\in ]x-\varepsilon ,x+\varepsilon \lbrack \} \leq \lim
\inf_{n\rightarrow +\infty }f(y_{n})\\
&\leq& \lim \sup_{n\rightarrow +\infty} f(y_{n})\leq \sup \{f(y),y\in ]x-\varepsilon ,x+\varepsilon \lbrack \}.
\end{eqnarray*}

\bigskip \noindent Finally let $\varepsilon \downarrow 0$ to have

\begin{equation*}
f_{\ast }(x)\leq \lim \inf_{n\rightarrow +\infty }f(y_{n})\leq \lim
\sup_{n\rightarrow +\infty }f(y_{n})\leq f^{\ast }(x).
\end{equation*}

\bigskip \noindent  Since $f^{\ast }(x)=$ $f_{\ast }(x),$ we conclude that the limit of 
$f(y_{n})$ exists when $n\rightarrow +\infty $ and $\inf \{f(y),y\in ]x-\varepsilon ,x+\varepsilon \lbrack \}$,
 
\begin{equation*}
\lim_{n\rightarrow +\infty }f(y_{n})=f^{\ast }(x)=f_{\ast }(x)=f(x).
\end{equation*}

\bigskip \noindent Since for any $(y_{n})_{n\geq 1}$ sequence such that $y_{n}\rightarrow x$ as 
$n\rightarrow +\infty $, we have $\lim_{n\rightarrow +\infty }f(y_{n})=f(x),$ we
conclude that $f$ is continuous at $x.$.\\

\noindent Question (4). Let $g=f^{\ast }$. Let us show that $g^{\ast }=g.$ Let $x$ be
fixed. By definition for $\varepsilon >0,$

\begin{equation*}
g_{\varepsilon }^{\ast }(x)=\sup \{g(z),z\in ]x-\varepsilon ,x+\varepsilon
\lbrack \}=\sup \{f^{\ast }(z),z\in ]x-\varepsilon ,x+\varepsilon \lbrack \}.
\end{equation*}

\bigskip \noindent Suppose that $g^{\ast }(x)$ is finite. So is $g_{\varepsilon }^{\ast }(x)$
for $\varepsilon$ small enough. Use the characterization of the supremum,
for all $\eta >0,$there exists $z_{0}\in ]x-\varepsilon ,x+\varepsilon
\lbrack $ such that

\begin{equation*}
g_{\varepsilon }^{\ast }(x)-\eta <f^{\ast }(z_{0})<g_{\varepsilon }^{\ast}(x).
\end{equation*}

\bigskip \noindent Since $z_{0}\in ]x-\varepsilon ,x+\varepsilon \lbrack ,$ there exists $%
\varepsilon _{0}>0$ such that $z_{0}\in ]z_{0}-\varepsilon
_{0},z_{0}+\varepsilon _{0}[\subset z]x-\varepsilon ,x+\varepsilon \lbrack$.\\

\noindent Remind that $f_{h}^{\ast }(z_{0})\downarrow f^{\ast }(z_{0})$ as $%
h\downarrow 0.$ Then for $0<h\leq \varepsilon _{0},$
\begin{equation*}
f_{h}^{\ast }(z_{0})=\sup \{f(z),z\in ]z_{0}-h,z_{0}+h[\}\leq \sup
\{f(z),z\in ]x-\varepsilon ,x+\varepsilon \lbrack \}=f_{\varepsilon }^{\ast
}(x),
\end{equation*}

\bigskip \noindent so that
\begin{equation*}
g_{\varepsilon }^{\ast }(x)-\eta <f^{\ast }(x).
\end{equation*}

\bigskip \noindent Let $\eta \downarrow 0$ and next $\varepsilon \downarrow 0$ to get%
\begin{equation}
g^{\ast }(x)\leq f^{\ast }(x).  \tag{INEQ}
\end{equation}

\noindent Suppose that these inequality is strict, that is%
\begin{equation*}
g^{\ast }(x)<f^{\ast }(x).
\end{equation*}

\bigskip \noindent We can find $\eta >0$ such that

\begin{equation*}
g^{\ast }(x)<f^{\ast }(x)-\eta.
\end{equation*}

\bigskip \noindent Since $g_{\varepsilon }^{\ast }(x)\downarrow g^{\ast }(x)$, there exists $%
\varepsilon_{0}$ such that

\begin{equation*}
g^{\ast }(x)\leq g_{\varepsilon _{0}}^{\ast }(x)<f^{\ast }(x)-\eta,
\end{equation*}

\bigskip \noindent that is
\begin{equation*}
g^{\ast }(x)\leq \sup \{g(z),z\in ]x-\varepsilon _{0},x+\varepsilon
_{0}[\}<f^{\ast }(x)-\eta, 
\end{equation*}

\bigskip \noindent that is also

\begin{equation*}
g^{\ast }(x)\leq \sup \{f^{\ast }(z),z\in ]x-\varepsilon _{0},x+\varepsilon
_{0}[\}<f^{\ast }(x)-\eta. 
\end{equation*}

\bigskip \noindent But
\begin{equation*}
x\in ]x-\varepsilon _{0},]x+\varepsilon _{0}[,
\end{equation*}

\bigskip \noindent and then
\begin{equation*}
f^{\ast }(x)\leq \sup \{f^{\ast }(z),z\in ]x-\varepsilon _{0},]x+\varepsilon
_{0}[\}<f^{\ast }(x)-\eta .
\end{equation*}

\bigskip \noindent We arrive at the conclusion
\begin{equation*}
f^{\ast }(x)<f^{\ast }(x)-\eta .
\end{equation*}

\bigskip \noindent This is absurd. Hence the inequality (INEQ) is an equality, that is%
\begin{equation*}
g^{\ast }(x)=f^{\ast}(x).
\end{equation*}

\noindent Now suppose that $g^{\ast }(x)$ is infinite. This means that for any $A>0$, we can find in any interval $]x-\varepsilon ,]x-\varepsilon \lbrack$,  a
point $z(\varepsilon )$ such that 
\begin{equation*}
g(z(\varepsilon ))>A,
\end{equation*}

\bigskip \noindent Now fix $A$ and $\varepsilon >0.$ Consider $z(\varepsilon )$ such that the last inequality holds. Since $z(\varepsilon )\in ]x-\varepsilon,]x-\varepsilon \lbrack ,$ there exists $r_{0}>0$ such that

\begin{equation*}
z(\varepsilon )\in ]z(\varepsilon )-r_{0},z(\varepsilon )+r_{0}[\subset
]x-\varepsilon ,x+\varepsilon \lbrack .
\end{equation*}

\bigskip \noindent So, since 
\begin{equation*}
f_{r}^{\ast }(z(\varepsilon ))\downarrow g(z(\varepsilon ))=f_{r}^{\ast
}(z(\varepsilon ))>A,
\end{equation*}

\bigskip \noindent we have for $r$ small enough
\begin{equation*}
f_{r}^{\ast }(z(\varepsilon ))>A/2,
\end{equation*}

\bigskip \noindent that is
\begin{equation*}
\sup \{f(z),z\in ]z(\varepsilon )-r,z(\varepsilon )+r[\}>A/2
\end{equation*}

\bigskip \noindent and again there will be a $u\in z\in ]z(\varepsilon )-r,z(\varepsilon )+r[$ such that
\begin{equation*}
f(u)>A/4.
\end{equation*}

\bigskip \noindent If $r$ is taken as small such that $r\leq r_{0}$, we have found a point $%
u\in ]x-\varepsilon ,]x-\varepsilon \lbrack $ such that%
\begin{equation*}
f(u)>A/4.
\end{equation*}

\bigskip \noindent We have proved that for any $A>0,$ for any $\varepsilon >0,$ we can find a
point in $u\in ]x-\varepsilon ,]x-\varepsilon \lbrack $ such that%
\begin{equation*}
f(u)>A/4
\end{equation*}

\bigskip \noindent and hence

\begin{equation*}
f_{\varepsilon }^{\ast }(x)=\sup \{f^{\ast }(z),z\in ]x-\varepsilon
,]x+\varepsilon \lbrack \}\geq A/4.
\end{equation*}

\bigskip \noindent Let $\varepsilon \downarrow 0$ and next $A\uparrow +\infty $ to get that
\begin{equation*}
g(x)=f^{\ast }(x)=+\infty .
\end{equation*}

\bigskip \noindent Thus 
\begin{equation*}
g^{\ast }(x)=g(x)=+\infty .
\end{equation*}

\bigskip \noindent We have finished and we proved that
\begin{equation*}
g^{\ast }=g.
\end{equation*}

\bigskip \noindent Conclusion : if $g^{\ast}$ is finite, then it is upper semi-continuous. Show that
if $g_{\ast }$ is finite, then it is lower semi-continuous by using $(-f)^{\ast
}=-f_{\ast}.$\\

\bigskip \noindent Question (5). We have two cases.\\

\noindent If 
$$
x\in D=\bigcup_{n\geq 1} D_n,
$$

\bigskip \noindent we have, as $n\rightarrow +\infty$,
$$
f_n(x)=f(x) \rightarrow f(x).
$$

\bigskip \noindent if 

$$
x\notin D=\bigcup_{n\geq 1} D_n,
$$

\bigskip \noindent then for $n\geq 1$, there exists $k=k(n,x)$ such that $x\in]k2^{-n},(k+1)2^{-n}[$. Then there exists $\varepsilon(n)>0$ such that for any 
$0<\varepsilon\leq \varepsilon(n)$, we have

$$
]x-\varepsilon,x+\varepsilon[ \subset x\in]k(n,x)2^{-n},(k(n,x)+1)2^{-n}[,  \text{ (A) }
$$

\bigskip
\noindent and next, 

$$
\sup \{f(z),z\in ]x-\varepsilon,x+\varepsilon[\} \leq \sup \{f(z),z\in [k2^{-n},(k+1)2^{-n}[\}. \text{ (B) }
$$

\bigskip
\noindent As well, for any $\eta>0$, for $2^{-n}<\eta$,  for $\varepsilon(n)>0$ such that (A) holds, we obviously have for any $n\geq 1$, 

$$
x\in]k(n,x)2^{-n},(k(n,x)+1)2^{-n}[ \in ]x-\eta,x+\eta[ ;
$$

\noindent S

$$
\sup \{f(z), z\in ]k(n,x)2^{-n},(k(n,x)+1)2^{-n}[\} \subset \sup \{f(z), z \in ]x-\eta,x+\eta[ \}. \text{ (C) }
$$

\bigskip
\noindent By combining (A) and (B), we have for $\eta>0$, $2^{-n}<\eta$, for $\varepsilon(n)>0$ such that (A) holds, for any $0<\varepsilon\leq \varepsilon(n)$, we have 

$$
f^{*,\varepsilon} \leq f_n(x) \leq f^{*,\eta}.
$$

\bigskip
\noindent Now let $\varepsilon\downarrow 0$ and only after, let $n\rightarrow +\infty$ to have

$$
f^{*} \leq \liminf_{n\rightarrow +\infty} f_n(x) \leq \limsup_{n\rightarrow +\infty} f_n(x) \leq f^{*,\eta}.
$$

\bigskip \noindent Finally, let $\varepsilon\downarrow 0$ to get

$$
f_n(x)\rightarrow f^{*}=f(x).
$$

\bigskip
\noindent Thus, we have $f_n(x)=f(x) \rightarrow f(x)$ as $n\rightarrow +\infty$. We conclude that $f$ is measurable.\\

\bigskip \noindent Question (6). The extension of the result of (5) to lower semi-continuous functions is immediate by using opposite functions.

\part{Measures and Integration}

%\chapter{Measures}
\chapter{Measures} \label{04_measures}

\bigskip
\noindent \textbf{Content of the Chapter}\\

\begin{table}[htbp]
	\centering
		\begin{tabular}{llll}
		\hline
		Type& Name & Title  & page\\
		\hline
		S&Doc 04-01  &Measures - A summary   &  \pageref{doc04-01}\\
		D&Doc 04-02  & Introduction to Measures - Exercises  & \pageref{doc04-02}\\
		S&Doc 04-03  & Carath\'{e}odory's Theorem - Lebesgue-Stieljes Measures   & \pageref{doc04-03}\\
		D&Doc 04-04  & Exercises on Lebesgue-Stieljes Measure  & \pageref{doc04-04}\\
		D&Doc 04-05  & General exercises on Measures  & \pageref{doc04-05}\\
		SD&Doc 04-06 & Introduction to Measures - Exercises with solutions & \pageref{doc04-06}\\
		SD&Doc 04-07 & Exercises on Lebesgue-Stieljes Measure with solutions & \pageref{doc04-07}\\
		SD&Doc 04-08 & General exercises on Measures with solutions  & \pageref{doc04-08}\\
		A&Doc 04-09  & The existence of the Lebesgue-Stieljes measure  & \pageref{doc04-09}\\
		A&Doc 04-10  & Proof the Carath\'{e}odory's Theorem  & \pageref{doc04-10}\\
		A&Doc 04-11  & Application of exterior measures to Lusin's theorem & \pageref{doc04-11}\\
		\hline
		\end{tabular}
\end{table}

\bigskip \noindent \textbf{Recommendation}. The results of Measures are summarized in \textit{Doc 04-01} and \textit{Doc 04-03}. The reader is asked to try to know these results by heart.\\

\newpage
\noindent \LARGE \textbf{DOC 04-01 : Measures - A summary}. \label{doc04-01}\\
\bigskip
\Large

\bigskip \noindent \textbf{Global picture : } Measures are modern and formal generalizations of \textit{intuitive measures} of lengths (in $\mathbb{R}$),
of areas (in $\mathbb{R}^{2})$, of volumes in $\mathbb{R}^{3}$, etc.\newline

\noindent \textbf{Definition}. Let $m$ be an application from a class $\mathcal{F}$ of subsets of $\Omega$ to $\overline{\mathbb{R}}$. The application $m$ is said to be
proper if and only if it exists an element $A$ of $\mathcal{F}$ such that $m(A)$ is finite, that is $|m(A)|<+\infty$.\\

\noindent \textbf{Warning}. In all this book, measures are meant as non-negative measures. Differences of non-negative measures are called signed measures.\\

\bigskip \noindent \textbf{I - DEFINITIONS}.

\bigskip \noindent \textbf{04.01a} (Definition) Let $\mathcal{A}$ be a $\sigma $-algebra of subsets
of $\Omega $. An application $m$ from $\mathcal{A}$ to $\overline{\mathbb{R}}_{{+}}$ is
a called a measure if and only if :\newline

\noindent (i) $m(\emptyset )$=0.\\

\noindent (ii) For all countable family $(A_{n})_{(n\geq I)}$, $I \subset \mathcal{N}$, of pairwise
disjoint elements of $\mathcal{A}$, the $\sigma$-additivity formula holds 
\begin{equation*}
m(\sum_{n\in I}A_{n})=\sum_{n\in I}m(A_{n}).
\end{equation*}

\bigskip \noindent \textbf{04.01b} (A second definition). Alternatively, the application $m$ from $\mathcal{A}$ to $\mathbb{R}_{{+}}$ is
a called a measure if and only if it is proper and $\sigma$-additivity.\\

\bigskip \noindent \textbf{04.02} Finite additivity. An application $m$ from an algebra $\mathcal{C}$ to $\overline{\mathbb{R}}_{+}$ is a said to be finitely additive if and only if
for any finite sequence $(A_{n})_{(0\leq n\leq k)}$ of $k$ pairwise disjoint elements of $\mathcal{C}$, with $k\geq 1$ and fixed, we have 
\begin{equation*}
m(\sum_{0\leq n\leq k}A_{n})=\sum_{0\leq n\leq k}m(A_{n}).
\end{equation*}

\bigskip \noindent \textbf{Remark}. Finite additivity is a special case of $\sigma$-additivity, where in (ii), the set $I$ is finite. Then, $\sigma$-additivity 
includes finite additivity.\\

\bigskip \noindent \textbf{04.03a} $\sigma$-sub-additivity. An application $m$ from $%
\mathcal{A}$ to $\overline{\mathbb{R}}_{+}$ is a said to be $\sigma $-sub-additive if and
only if for all sequence $(A_{n})_{(n\geq 0)}$ of elements of $\mathcal{A}$,
we have 
\begin{equation*}
m(\bigcup\limits_{n\geq 0}A_{n})\leq \sum_{n\geq 0}m(A_{n}).
\end{equation*}

\bigskip \noindent \textbf{04.03b}  As well, an application $m$ from an algebra $\mathcal{A}$ to $\overline{\mathbb{R}}_{+}$ is a said to be finitely sub-additive, or simply sub-additive, if and
only if for any finite sequence $(A_{n})_{(0 \leq n \leq k)}$ $k$  elements of $\mathcal{C}$, with $k\geq 1$ and fixed,
we have 
\begin{equation*}
m(\bigcup\limits_{0 \leq n \leq k}A_{n})\leq \sum_{0\leq n\leq k}m(A_{n}).
\end{equation*}

\bigskip \noindent \textbf{04.04} If $m$ is a measure on $\mathcal{A}$, we say that $(\Omega ,\mathcal{A},m)$ is a measure space.\newline

\newpage
\bigskip \noindent \textbf{II - MEASURES ON AN ALGEBRAS OR ON A SEMI-ALGEBRAS}.\\

\noindent In the definition of a measure on the measurable space  ($\Omega $,$\mathcal{A)}$ in \textbf{Points 04.01} or in \textit{Point 04.02} above, we know that $\sum_{0\leq
n\geq k}A_{n}$ is still in $\mathcal{A}$ whenever the $A_{n}$'s are.\\

\noindent But, when using the notion of measure on \textbf{an algebra} $\mathcal{C}$, or on a \textbf{semi-algebra} $\mathcal{S}$, we are not sure that $\sum_{0\leq n\geq k}A_{n}$ is $\mathcal{C}$ in
whenever the $A_{n}$'s are. In such  restricted cases, we use the following
definition.\\

\noindent \textbf{04.05a} Let $\mathcal{C}$ be an algebra or a semi-algebra of subsets of $\Omega 
$. An application $m$ from $\mathcal{C}$ to $\overline{\mathbb{R}}_{+}$ is a called a
measure if and only if \ $m(\emptyset )=0$ and for any sequence $(A_{n})_{(n\geq 0)}$ of pairwise disjoint elements of $\mathcal{A}$ \textbf{such that} $\sum_{n\geq 0}A_{n}\in \mathcal{C},$ we have  
\begin{equation*}
m(\sum_{n\geq 0}A_{n})=\sum_{n\geq 0}m(A_{n}).
\end{equation*}

\bigskip \noindent \textbf{04.05b} $\sigma$-sub-additivity on. An application $m$ from $\mathcal{C}$ to $\overline{\mathbb{R}}_{+}$ is a said to be $\sigma $-sub-additive if and
only if for all sequence $(A_{n})_{(n\geq 0)}$ of elements of $\mathcal{C}$ \textbf{such that} $\bigcup_{n\geq 0}A_{n}\in \mathcal{C}$,
we have 
\begin{equation*}
m(\bigcup\limits_{n\geq 0}A_{n})\leq \sum_{n\geq 0}m(A_{n}).
\end{equation*}

\bigskip \noindent \textbf{04.05c}  As well, an application $m$ from an algebra $\mathcal{C}$ to $\overline{\mathbb{R}}_{+}$ is a said to be finitely sub-additive, or simply sub-additive, if and
only if for any finite sequence $(A_{n})_{(0 \leq n \leq k)}$ $k$  elements of $\mathcal{C}$, with $k\geq 1$ and fixed, \textbf{such that} $\bigcup_{0\leq n\leq k}A_{n}\in \mathcal{C}$
we have 
\begin{equation*}
m(\bigcup\limits_{0 \leq n \leq k}A_{n})\leq \sum_{0\leq n\leq k}m(A_{n}).
\end{equation*}

\bigskip \noindent \textbf{04.05d} (Alternative definition).\\

\noindent (1) An application $m$ from $\mathcal{C}$ to $\overline{\mathbb{R}}_{+}$ is a said to be $\sigma $-sub-additive if and
only if for all sequence $(A_{n})_{(n\geq 0)}$ of pairwise disjoint elements of $\mathcal{C}$ \textbf{such that} $\sum_{n\geq 0}A_{n}\in \mathcal{C}$,
we have 
\begin{equation*}
m(\sum_{n\geq 0}A_{n})\leq \sum_{n\geq 0}m(A_{n}).
\end{equation*}

\bigskip \noindent (2)  As well, an application $m$ from an algebra $\mathcal{C}$ to $\overline{\mathbb{R}}_{+}$ is a said to be finitely sub-additive, or simply sub-additive, if and
only if for any finite sequence $(A_{n})_{(0 \leq n \leq k)}$ $k$  pairwise disjoint elements of $\mathcal{C}$, with $k\geq 1$ and fixed, \textbf{such that} $\sum_{0\leq n\leq k}A_{n}\in \mathcal{C}$
we have 
\begin{equation*}
m(\sum_{0 \leq n \leq k}A_{n})\leq \sum_{0\leq n\leq k}m(A_{n}).
\end{equation*}

\newpage
\bigskip \noindent \textbf{III - MAIN PROPERTIES OF MEASURES}.\newline

\bigskip \noindent \textbf{04.06} $m$ is non-decreasing : If $A\subset B$, then $%
m(A)\leq m(B)$.\newline

\noindent Besides, If $A\subset B$, and $m(B)$ is finite, that $m(B \setminus A)=m(B)-m(A \cup B)$.\\ 

\bigskip \noindent \textbf{04.07} A measure $m$ finite if and only if only if $m(\Omega)$ is finite if and only if $m$ is bounded.\newline

\bigskip \noindent \textbf{04.08} If a measure $m$ is finite, then for any measurable sets $A$ and $B$, then $m(B \setminus A)=m(B)-m(A\cap B)$.\newline

\bigskip \noindent \textbf{04.09} A measure $m$ is $\sigma$-sub-additive.\newline

\bigskip \noindent \textbf{04.10} A measure $m$ is continuous in the following sense :\\

\bigskip \noindent \textbf{(a)} Let $(A_{n})_{(n\geq 0)}$ be a non decreasing sequence of elements of $\mathcal{A}$ and let us denote $A=\bigcup\limits_{n\geq 0}A_{n}$. Then 
$m(A_{n})$ broadly increases to $m(A)$ as $n\rightarrow +\infty$.\newline

\bigskip \noindent \textbf{(b)} Let $(A_{n})_{(n\geq 0)}$ be a non increasing sequence of
elements of $\mathcal{A}$ and let us denote $A=\bigcap\limits_{n\geq 0}A_{n}$ such :\\

$$
\exists n_{0}\geq 0, \ m(A_{n_{0}})<+\infty. 
$$

\noindent Then $m(A_{n})$ broadly decreases to $m(A)$ as $n\rightarrow +\infty$.\\

\bigskip \noindent \textbf{04.11} Non-negative finite and non-negative countable linear
combinations of measures are measures.\newline

\bigskip \noindent \textbf{04.12} A probability measure (or simply a probability) is
a measure $\mathbb{P}$ such that $\mathbb{P}(\Omega )=1$.  If $\mathbb{P}$
is a probability measure on $\mathcal{A}$, $(\Omega ,\mathcal{A},\mathbb{P})$
is called a probability space.\newline

\newpage
\bigskip \noindent \textbf{IV - Simple examples of measure}.\\

\noindent \textbf{(04.13) The null measure}. On any measurable space$(\Omega, \mathcal{A})$, we have the null measure $z$ which assigns the null value to all measurable set :
$$
\forall A \in \mathcal{A}, \ \ z(A)=0.
$$

\bigskip \noindent \textbf{(04.14a) Dirac measure}.\\

\noindent Consider a measurable space $(\Omega, \mathcal{A})$, one point $\omega_0$ of $\Omega$ and a non-negative real number $\alpha$, the application $\delta_{\omega_0,\alpha}$ defined on $\mathcal{A}$ by

$$
\delta_{\omega_0,\alpha}(A)=\alpha 1_{A}, \ \ A \in \mathcal{A},
$$

\bigskip \noindent is a measure, called the delta-measure concentrated at $\omega_0$ with mass $\alpha$.

\noindent if $\alpha=1$, we simply write $\delta_{\omega_0}$ and call it \textit{the delta-measure concentrated at $\omega_0$}.\\

\noindent \textbf{(04.14b) Counting measure}.\\

\noindent The counting measure is a very important one in this theory from a conceptual point of view. It is usually defined on a countable set 
$\Omega=\{\omega_1, \omega_2, ...\}$ endowed with the $\sigma$-algebra $\mathcal{P}(\Omega)$, like $\mathbb{N}$ the set of nonnegative integers, $\mathbb{Z}$ the set of all integers, $\mathbb{Q}$ the set of rational numbers, etc.\\

\noindent By \textit{Point 04.11}, the sum of all the $\delta_{\omega_j}$'s 
$$
\nu = \sum_{j\geq 1} \delta_{\omega_j}
$$

\bigskip \noindent is a measure such that for any $A \subset \Omega$,

$$
\nu(A) = \sum_{j\geq 1} \delta_{\omega_j}(A).
$$

\bigskip \noindent We may quickly see that for any $A \subset \Omega$,

$$
\nu(A) = Card (A).
$$

\bigskip \noindent The name of counting measure comes from the fact that $\nu$ counts the number of elements of subsets of a countable set $\Omega$.\\

\noindent \textbf{NOTA-BENE}. Except this important measure, other nontrivial measures are either constructed or derived from constructed measures.\\

\newpage
\bigskip \noindent \textbf{V - IMAGE MEASURE AND INDUCED MEASURE}.\\

\noindent Let $(\Omega, \mathcal{A},m)$ be a measure space and let $X : (\Omega, \mathcal{A},m) \mapsto (E, \mathcal{B})$ be a measurable application from $(\Omega, \mathcal{A},m)$ to the measure space $(E, \mathcal{B})$.

\noindent We may \textbf{transfer} the measure $m$ defined on $(\Omega, \mathcal{A})$ to a new one denoted $m_X$ defined on $(E, \mathcal{B})$ by means of the inverse image function $X^{-1}$ of $X$.\\

\noindent \textbf{4.15 Image Measure}. The application $m_X$ defined on $\mathcal{B}$ by
$$
m_X(B) = m(X^{-1}(B)), \ \ B \in \mathcal{B},
$$

\bigskip \noindent is a measure on $(E, \mathcal{B})$ called image measure (of $m$ by $X$).\\

\noindent \textbf{Important remark}. Despite of its simplicity, image measures are very important in probability theory. The image of a probabibility measure $\mathbb{P}$ by a measurable application $X$ is called the \textbf{probability law}, denoted $\mathbb{P}_X$. Finding probability laws of measurable applications is the top objective in that theory.\\

\bigskip \noindent \textbf{4.16 Induced Measure}. (Please, revise \textit{Exercise 7, Doc 01-02 of Chapter 1} before proceeding).\\

\noindent Let $A$ be a measurable set on the measure space $(\Omega, \mathcal{A},m)$. The application $m_A$ defined on $(A,\mathcal{A}_A)$ by
$$
m_A(B) = m(A\cap B), \ \ B \subset A,
$$

\bigskip \noindent is a measure on $(A,\mathcal{A}_A)$ called the induced measure on the induced measurable space.\\

\newpage \noindent \textbf{VI - NULL SETS AND ALMOST-EVERYWHERE PROPERTIES}.\\

\noindent \textbf{(A) Null sets}.\\

\bigskip \noindent \textbf{04.17} A measurable set $N$ is said to be a \textit{m-null
set} if and only if $m(N)=0$.\newline

\bigskip \noindent \textbf{Immediate properties}.\\

\noindent \textbf{(1)} A measurable set that is included in a $m$-null set is a $m$-null set.\\

\noindent \textbf{(2)} A finite or a countable union $m$-null sets is a $m$-null set.\newline

\bigskip \noindent \textbf{(B) Almost everywhere true properties}.\\

\bigskip \noindent \textbf{04.18} Let $\mathcal{P}(\omega )$ a proposition depending on $\omega $. We say that $\mathcal{P}$ is true \textit{almost-everywhere (a.e.)} if and only if the set 
$$
\{\omega ,\mathcal{P}(\omega )\text{ is false}\}
$$ 

\bigskip \noindent is measurable and 

\begin{equation*}
m(\{\omega ,\mathcal{P}(\omega )\text{ is false}\})=0.
\end{equation*}

\bigskip
\noindent This is denoted by : $\mathcal{P}$ is true almost-everywhere, abbreviated into : 

$$
\mathcal{P} \  a.e.
$$

\bigskip \noindent \textbf{Measure Space Reduction Property} \label{msrp} : If each assertion of a countable family of assertions $(\mathcal{P}_n, \ n\geq 0)$ defined on the measure space $(\Omega, \mathcal{A}, m)$ holds almost everywhere, then all the assertions simultaneously hold on a measurable subspace $\Omega_0$, which is complement of a null-set, that is $m(\Omega_0)=0$.\\

\noindent In situations where we have a a countable family of \textit{a.e.} true assertions $(\mathcal{P}_n, \ n\geq 0)$, we may and usually assume the they are all true everywhere, but putting ourselves on the induced measure space $(\Omega_0, \mathcal{A}_{\Omega_0}, m_{\Omega_0})$ (See Point \textit{04.16} above). This change does not affect the measure of any measurable set $A$, since

$$
m(A)=m_{\Omega_0}(A).
$$

\newpage
\noindent \textbf{(C) Equivalence classes of measurable applications}. \label{04_equivcal}\\

\noindent \textbf{04.19a} Denote by $\mathcal{L}_{0}(\Omega, \mathcal{A})$ the class of all real-valued and measurable applications defined on $(\Omega, \mathcal{A})$. We write $\mathcal{L}_{0}$ in short.\\

\noindent Consider two elements $f$ and $g$ of $\mathcal{L}_{0}$, the sets
$$
(f=g)=\{\omega \in \Omega, f(\omega)=g(\omega)\} \ \ and  \ \ (f\neq g)=\{\omega \in \Omega, f(\omega) \neq g(\omega)\}
$$

\bigskip \noindent are measurable.\\ 

\noindent \textbf{04.19b} We may define the binary relation $\mathcal{R}$ on $\mathcal{L}_{0}^{2}$ defined by

$$
f \mathcal{R} g \Leftrightarrow m(f\neq g)=0,
$$

\bigskip \noindent that is, $f$ and $g$ are equal outside of a $m$-null set, denoted also by : $f=g$ almost-everywhere (a.e.)\\

\noindent This is an equivalence relation. The equivalence classes of this relation are :\\

$$
\dot{f}=\{g \in \mathcal{L}_{0}, \ g=f \ a.e. \}.
$$ 

\bigskip \noindent \textbf{04.20} The set of all equivalence classes is denoted

$$
L_{0}(\Omega, \mathcal{A},m).
$$

\newpage
\noindent \textbf{(D) \textbf{Complete measure}.}\\

\noindent We know that a \textbf{measurable} subset of m-null set
is a m-null set. We cannot say anything for that a subset of a m-null if we
do not know if this subset is measurable or not. We have the following
definition.\\

\noindent \textbf{04.21}  A measure space is complete if and only if any subset of a null set
is measurable and then is a null set.\\

\noindent \textbf{04.22} Fortunately, it is possible to extend non-complete measures to complete measures in the following way. Assume that 
$(\Omega ,\mathcal{A},m)$ is not complete. Set 
\begin{equation*}
\mathcal{N}=\{N\subset \Omega ,\exists \text{ }B\text{ }m\text{-null set}%
,N\subset B\}).
\end{equation*}

\bigskip \noindent Define

\begin{equation*}
\widehat{\mathcal{A}}=\{A\cup N,A\in \mathcal{A},N\in \mathcal{N}\},
\end{equation*}

\bigskip \noindent and next define the application 
$$
\widehat{m}:\widehat{\mathcal{A}}%
\rightarrow \mathbb{R}_{+}
$$ 

\bigskip \noindent such that 
$$
\widehat{m}(A\cup N)=m(A),
$$ 

\bigskip \noindent for $A\in \mathcal{A}$, $N\in \mathcal{N}$.\\

\noindent  We have : \\

\noindent (1) $\widehat{\mathcal{A}}$ is a $\sigma$-algebra including $\mathcal{A}$,\\

\noindent (2) $\widehat{m}$ is a measure on $(\Omega,\widehat{\mathcal{A}})$,\\

\noindent (3) the measure space \noindent $(\Omega ,\widehat{\mathcal{A}},\widehat{m})$ is a complete measure space and $m$ is a restriction of $%
\tilde{m}$ on $\mathcal{A}$.\\

\noindent We may and do make the following conclusion : \textbf{Any measure may be considered as a complete measure and, from now on, in any measure space; subsets of null-sets are null sets.}\newline

\newpage
\bigskip \noindent \textbf{VI - SIGMA-FINITE MEASURES AND CONTINUITY OF MEASURES}\\

\bigskip \noindent \textbf{4.23a} Let $m$ be an application of a class $\mathcal{F}$ of subsets of $\Omega$ to $\overline{\mathbb{R}}$ and let $\mathcal{F}_0$ be a sub-class of $\mathcal{F}$. The application $m$ is said to be $\sigma$-finite on $\Omega$ with respect to $\mathcal{F}_0$ if and only if there exists a finite or a countable partition of $\Omega$ into elements of $\mathcal{F}_0$, that
is a sequence $(\Omega _{j})_{j\geq 0}$ of pairwise disjoint elements of $\mathcal{F}_0$ satisfying 
\begin{equation*}
\Omega =\sum_{j\geq 0}\Omega _{j} \ \ \ (CSF1),
\end{equation*}

\bigskip \noindent such that for all $j\geq 0$, $m(\Omega_j)<+\infty$.\newline

\bigskip \noindent \textbf{4.23b} Let $m$ be a proper non-negative and additive application an algebra $\mathcal{C}$ of subsets of $\Omega$. The $\sigma$-finiteness also holds if and only if there exists a finite or a countable cover of $\Omega$ by elements of $\mathcal{C}$, that is a sequence $(\Omega _{j})_{j\geq 0}$ of elements of $\mathcal{C}$ satisfying  
\begin{equation*}
\Omega =\bigcup_{j\geq 0}\Omega _{j}  \ \ \ (CSF2),
\end{equation*}

\bigskip \noindent such that for all $j\geq 0$, $m(\Omega_j)<+\infty$.\newline

\bigskip \noindent \textbf{04.24} Let $m_{1}$ and $m_{2}$ defined on the same measurable space $(\Omega ,%
\mathcal{A})$. The measure $m_{1}$ is continuous is with respect to $m_{2}$
if and only if we have

\begin{equation*}
\forall A\in \mathcal{A},(m_{2}(A)=0)\Rightarrow (m_{1}(A)=0).
\end{equation*}

\bigskip \noindent \textbf{04.25} We have the same definition of the continuity of $\sigma $-additive applications between them.\newline

\bigskip \noindent \textbf{VII - Outer measures}

\bigskip \noindent An application $\mu$\ defined on $\mathcal{P}(\Omega)$ is an exterior measure if and only if it is positive, non-decreasing, $\sigma -$sub-additive and assigns the value $0$ to
$\emptyset$, that is we have the following properties\newline

\bigskip \noindent (i) For any sequence ($A_{n})_{n\geq 1}$ of subsets of $\Omega$
\begin{equation*}
\mu(\bigcup_{n\geq 1}A_{n})\leq \sum_{i=1}^{\infty }m^{0}(A_{i}). 
\end{equation*}

\bigskip \noindent (ii) For any subset $A$ of $\Omega $ 
\begin{equation*}
\mu(A)\geq 0.
\end{equation*}

\bigskip \noindent For every $A\subset B,$%
\begin{equation*}
\mu(A)\leq \mu(B). 
\end{equation*}

\bigskip \noindent Also, we have

\begin{equation*}
\mu(\emptyset )=0. 
\end{equation*}

\bigskip \noindent The The class of $\mu$-measurable sets 

$$
\mathcal{A}=\{A \in \mathcal{P}(\Omega) : \ \forall D\subset \Omega, \ \mu(D)=\mu(AD)+\mu(A^{c}D) \},
$$ 

\bigskip \noindent is a $\sigma$-algebra on which $\mu$ is a measure.\\

\noindent \textbf{NB}. This result is part of the Appendix in Doc 04-10 (see \pageref{doc04-10}).\\

\newpage
\noindent \LARGE \textbf{Doc 04-02 : Introduction to Measures - Exercises}. \label{doc04-02}\\

\Large

\noindent \textbf{PART A : General properties on measures and additive sets applications}.\\

\bigskip \noindent \textbf{Exercise 1}. \label{exercise01_doc04-02} Let $\mathcal{A}$ be $\sigma $-algebra of subsets
of $\Omega $. Let  $m$ be an application from $\mathcal{A}$ to $\overline{\mathbb{R}}_{+}$. Suppose that we have.

\noindent (ii) For all sequence $(A_{n})_{(n\geq 0)}$ of pairwise
disjoint elements of $\mathcal{A}$, the $\sigma$-additivity formula holds 
\begin{equation*}
m(\sum_{n\geq 0}A_{n})=\sum_{n\geq 0}m(A_{n})
\end{equation*}

\bigskip \noindent holds. Show that the two assertions are equivalent : \\

\noindent (i1) $m(\emptyset)$=0.\\

\noindent (i2) $m$ is a proper application constantly equal to $+\infty$, that is : 
$$
\exists A \in \mathcal{A}, \ m(A)<+\infty.
$$

\bigskip \noindent Show why the results you obtain still hold if $m$ is only finitely additive on an algebra.\\

\noindent \textit{Conclusion}. The results of the exercise justify the equivalence of the definitions of a measure in \textit{Points (04-01a) and (04-01b)}.\\

\bigskip \noindent \textbf{Exercise 2} \label{exercise02_doc04-02}. Let $(\Omega,\mathcal{A},m)$ be measure space. Show the following properties.\\

\noindent (1)  $m$ is non-decreasing : If $A\subset B$, then $m(A)\leq m(B)$. Besides, show, if $m(B)$ is finite, that $m(B \setminus A)=m(B)-m(A \cap B)$.\\ 

\noindent \textit{Hint} : Use and justify $B = (A\cap B) + (B\setminus A)$ and use the additivity and the non-negativity of a measure.\\

\noindent (2) A measure $m$ finite if and only if only if $m(\Omega )$ is finite if and only if $m$ is bounded.\\

\noindent \textit{Hint} : Use the non-decreasingness proved in Question (1).\\

\noindent (3) A measure $m$ is $\sigma $-sub-additive.\\

\noindent In a similar manner, show that if $m(\emptyset)=0$, and $m$ is non-negative and additive, then it is sub-additive

\noindent \textit{Hint} : Use the formula of \textit{Exercise 4 in DOC 00-02} (page \pageref{exercise02_doc00-02}) in \textit{Chapter \ref{00_sets}}.\\

\noindent (4) Show that $m$ is continuous in the following sense :\\

\bigskip \noindent \textbf{(a)} Let $(A_{n})_{(n\geq 0)}$ be a non-decreasing sequence
of elements of $\mathcal{A}$ and set $A=\bigcup\limits_{n\geq 0}A_{n}$. Then 
$m(A_{n})$ increases to $m(A)$ as $n\rightarrow +\infty$.\newline

\bigskip \noindent \textbf{(b)} Let $(A_{n})_{(n\geq 0)}$ be a non-increasing sequence of
elements of $\mathcal{A}$ and set $A=\bigcap\limits_{n\geq 0}A_{n}$ such that $m(A_{n_{0}})$ for some $n_{o}\geq 0$. Then $m(A_{n})$ decreases to $m(A)$ as $n\rightarrow +\infty$.\\

\noindent \textit{Hint} : (1) To show (a), remark by the help of a simple diagram, that if $(A_{n})_{(n\geq 0)}$ is a non-decreasing sequence, we have

$$
\bigcup\limits_{n\geq 0}A_{n} = A_0 + \sum_{n\geq 1} \left(A_n\setminus A_{n-1}\right).
$$ 

\bigskip \noindent Next, use the $\sigma$-additivity of $m$ and then the partial sums of the obtained series.\\

\noindent To prove (b), use Point (a) by taking the complements of the $A_n$'s.\\

\bigskip \noindent \textbf{Exercise 3}. \label{exercise03_doc04-02} Let $\mathcal{C}$ be an algebra of subsets of $\Omega$. Let $m$ be a proper and additive application from $\mathcal{C}$ to $\overline{\mathbb{R}}_{+}$.\\

\noindent Show that the two following assertions are equivalent.\\

\noindent (a) For any sequence $(A_{n})_{(n\geq 0)}$ of elements of $\mathcal{C}$ \textbf{such that} $\bigcup_{n\geq 0}A_{n}\in \mathcal{C}$,
we have 
\begin{equation*}
m(\bigcup\limits_{n\geq 0}A_{n})\leq \sum_{n\geq 0}m(A_{n}).
\end{equation*}

\bigskip \noindent (b) For all sequence $(A_{n})_{(n\geq 0)}$ of pairwise disjoint elements of $\mathcal{C}$ \textbf{such that} $\sum_{n\geq 0}A_{n}\in \mathcal{C}$,
we have 
\begin{equation*}
m(\sum_{n\geq 0}A_{n})\leq \sum_{n\geq 0}m(A_{n}).
\end{equation*}

\bigskip \noindent \textbf{Exercise 4}. \label{exercise03_doc04-02} Let $\mathcal{C}$ be an algebra of subsets of $\Omega$.  Let $m$ be a proper non-negative and additive application from $\mathcal{C}$ 
to $\overline{\mathbb{R}}_{+}$. Show that the two assertions are equivalent.\\

\noindent (a) There exists a finite or a countable partition of $\Omega$ into elements of $\mathcal{F}_0$, that is a sequence $(\Omega _{j})_{j\geq 0}$ of pairwise disjoint elements of $\mathcal{F}_0$ satisfying 
\begin{equation*}
\Omega =\sum_{j\geq 0}\Omega _{j} \ \ \ (CSF1),
\end{equation*}

\noindent such that for all $j\geq 0$, $m(\Omega_j)<+\infty$.\\

\noindent (b) There exists a finite or a countable cover of $\Omega$ by elements of $\mathcal{C}$, that is a sequence $(\Omega _{j})_{j\geq 0}$ of elements of $\mathcal{C}$ satisfying  
\begin{equation*}
\Omega =\bigcup_{j\geq 0}\Omega _{j}  \ \ \ (CSF2),
\end{equation*}

\noindent such that for all $j\geq 0$, $m(\Omega_j)<+\infty$.\\

\bigskip \noindent The implication $(a) \Rightarrow (b)$ is needless to prove since $(b)$ is a special case of $(a)$. Let $(b)$ holds. Let $(\Omega _{j})_{j\geq 0}$ be a sequence of elements of 
$\mathcal{C}$ such that $(CSF2)$ holds. Then by \textit{Exercise 4 in Chapter 0, Doc 00-02}, we have 
$$
\Omega = \sum_{n\geq 0} \ \Omega_j=\sum_{j\geq 0} \ \Omega^{\prime}_{j}
$$

\bigskip \noindent with

$$
\Omega^{\prime}_{0}=\Omega^{\prime}_{0}, \ \Omega^{\prime}_{j}=\Omega_{0}^{c}...\Omega_{j-1}^{c}\Omega_{j} \subset \Omega_j, \ j\geq 1.
$$ 

\bigskip \noindent and we also have : 

$$
\forall \ j>0, \ m(\Omega^{\prime}_{j}) \leq m(\Omega_{j}) <+\infty.
$$

\bigskip \noindent The solution is complete.$\blacksquare$\\

\newpage \noindent \textbf{PART II : Measures}.

\bigskip \noindent \textbf{Exercise 5}. \label{exercise05_doc04-02}\\ 

\noindent (1) Let $(m_i)_{i\in I}$ a family of measures on $(\Omega, \mathcal{A})$, where $I\subset \mathbb{N}$. Let $(m\alpha_i)_{i\in I}$ a family of non-negative real numbers such that one of them is positive. Consider the application $m$ defined for  $A \in \mathcal{A}$ by

$$
m(A)=\sum_{i \in I} \alpha_{i}m_{i}(A).
$$

\bigskip \noindent Show that $m=\sum_{i \in I} \alpha_{i}m_{i}$ is a measure on $(\Omega, \mathcal{A})$. Give a conclusion.\\

\noindent Hint : use the Fubini's rule for series of non-negative numbers.\\

\bigskip 
\noindent (2) Consider a measurable space $(\Omega, \mathcal{A})$, one point $\omega_0$ of $\Omega$ and a non-negative real number $\alpha$, show that the application $\delta_{\omega_0,\alpha}$ defined on $\mathcal{A}$ by

$$
\delta_{\omega_0,\alpha}(A)=\alpha 1_{A}, \ \ A \in \mathcal{A},
$$

\bigskip \noindent is a measure, called the delta-measure concentrated at $\omega_0$ with mass $\alpha$.\\

\noindent if $\alpha=1$, we denote

$$
\delta_{\omega_0,\alpha} \equiv \delta_{\omega_0}.
$$

\bigskip \noindent (3) Let $\Omega=\{\omega_1, \omega_2, ...\}$ be a countable space endowed with the $\sigma$-algebra $\mathcal{P}(\Omega)$, which is the power set of $\Omega$. By using Points (1) and (2) below, to show that the application defined by
$$
\nu = \sum_{j\geq 1} \delta_{\omega_j}
$$

\bigskip \noindent is a measure. Show that for any $A \subset \Omega$, we have

$$
\nu(A) = Card(A)
$$

\bigskip \noindent Suggest a name for this measure.\\

\bigskip \noindent \textbf{Exercise 6}. \label{exercise06_doc04-02} Let $(\Omega, \mathcal{A},m)$ be a measure space.\\

\noindent (a) Let $X : (\Omega, \mathcal{A},m) \mapsto (E, \mathcal{B})$ be a measurable application from $(\Omega, \mathcal{A},m)$ to the measure space $(E, \mathcal{B})$. Show that the application defined for $B\ \mathcal{B}$ by
$$
m_X(B) = m(X^{-1}(B)), \ \ B \in \mathcal{B},
$$

\bigskip \noindent is a measure on $(E, \mathcal{B})$ called image measure (of $m$ by $X$).\\

\noindent (b)  Let $A$ be a measurable set on the measure space $(\Omega, \mathcal{A},m)$. Show that the application $m_A$ defined on $(A,\mathcal{A}_A)$ by
$$
m_A(B) = m(A\cap B), \ \ B \subset A,
$$

\bigskip \noindent is a measure on $(A,\mathcal{A}_A)$.\\

\newpage
\noindent \textbf{Part III : null sets, equivalence classes, completeness of measures}.\\

\noindent \textbf{Exercise 7}. \label{exercise07_doc04-02} Let $(\Omega, \mathcal{A},m)$ a measure space.\\

\noindent (1) Show that  : (a) A measurable subset of a $m$-null set is a $m$-null set, (b) a countable union of $m$-null sets is a $m$-null set. Hint : use the $\sigma$-sub-additivity of $m$.\\

\noindent (2) Let $f : (\Omega, \mathcal{A},m) \mapsto \overline{\mathbb{R}}$ a real-valued measurable application.\\

\noindent Show that
$$
(f \text{ finite })=(|f|<+\infty)=\bigcup_{k\geq 1} (|f|\leq k) \text{   (A)}
$$

\bigskip \noindent and then, by taking complements,

$$
(f \text{ infinite })=(|f|=+\infty)=\bigcap_{k\geq 1} (|f| > k) \text{   (B)}.
$$

\bigskip \noindent Show also that

$$
(f=0)=\bigcap_{k\geq 1} (|f| <1/k) \text{   (C)}.
$$

\bigskip \noindent Remark : you may use strict or large inequalities in each formula.\\

\noindent Then the definitions \\

\noindent [ $f$ finite a.e. if and only if $m(f \text{ finite })=0$ ]\\

\noindent and\\

\noindent [ $f$ positive a.e. if and only if $m(f=0)=0$ ]\\

\noindent make sense.\\

\bigskip 
\noindent (3) Let $f_n : (\Omega, \mathcal{A},m) \mapsto \overline(\mathbb{R})$ a sequence real-valued measurable applications.\\

\noindent Show that if each $f_n$ is a.e. finite, they all the $f_n$'s are simultaneously finite outside a null-set.\\

\noindent Show that if each $f_n$ is a.e. positive, they all the $f_n$'s are simultaneously positive outside a null-set.\\ 

\bigskip 
\noindent (4) Let $f,g,h : (\Omega, \mathcal{A},m) \mapsto \mathbb{R}$ three finite real-valued measurable applications.\\

\noindent (a) Based on the triangle inequality $|f-h|\leq |f-g|+|g-h|$, remark that

$$
(|f-g|=0) \cap (|g-h|=0) \subset (|f-h|=0).
$$

\bigskip \noindent (b) and by taking the complements, that

$$
(|f-h|>0) \subset (|f-g|>0) \cup (|g-h|>0).
$$

\bigskip \noindent (c) Deduce from this that : $f=g$ a.e. and $g=h$ a.e., implies that $f=h$ a.e.\\

\noindent (d) Let $\mathcal{L}_0(\Omega, \mathcal{A}, m)$ the space of all real-valued measurable and a.e. finite applications $f_n : (\Omega, \mathcal{A},m) \mapsto \overline(\mathbb{R})$. Define the binary relation $\mathcal{R}$ on $\mathcal{L}_0(\Omega, \mathcal{A}, m)^2$ by :

$$
\forall (f,g) \in \mathcal{L}_0(\Omega, \mathcal{A}, m)^2, \ \ f \mathcal{R}g \ f=g \text{ a.e.}.
$$

\bigskip \noindent Show that $\mathcal{R}$ is an equivalence relation and the equivalence class of each $f \in \mathcal{L}_0(\Omega, \mathcal{A}, m)$ is

$$
\overset{\circ}{f}=\{g \in \mathcal{L}_0(\Omega, \mathcal{A}, m), \  f=g \text{ a.e.}\}.
$$

\bigskip \noindent The quotient set is denoted

$$
L_0(\Omega, \mathcal{A}, m)= \{\overset{\circ}{f}, f \in \mathcal{L}_0(\Omega, \mathcal{A}, m) \}.
$$

\bigskip \noindent \textbf{Exercise 8}. \label{exercise08_doc04-02} A measure space $(\Omega ,\mathcal{A},m)$ is said to be complete if and only if all subsets on $m$-null sets are measurable, and the then null $m$-sets.\\

\noindent Such a notion of complete measure is very practical. Indeed, not being able to measure a set one knows it is part of $m$-null set is pathological. This opposes to naive belief that any part of a negligeable (null set) object should be negligeable too. The following exercise aims at fixing this problem.\\ 

\bigskip \noindent Let $(\Omega ,\mathcal{A},m)$ a measure. Set 
\begin{equation*}
\mathcal{N}=\{N\subset \Omega ,\exists \text{ }B\text{ }m\text{-null set} ,N\subset B\}).
\end{equation*}

\noindent Remark that $\mathcal{N}$ includes the $m$-null sets, $\emptyset$ in particular.\\
 
\bigskip \noindent Define

\begin{equation*}
\widehat{\mathcal{A}}=\{A\cup N,A\in \mathcal{A},N\in \mathcal{N}\})
\end{equation*}

\bigskip \noindent and next define the application 
$$
\widehat{m}:\widehat{\mathcal{A}}\rightarrow \mathbb{R}_{+},
$$ 

\bigskip \noindent such that 
$$
\widehat{m}(A\cup N)=m(A).
$$ 

\bigskip \noindent (a) Show that $\widehat{\mathcal{A}}$ is a $\sigma$-algebra including $\mathcal{A}$.\\

\noindent Hint. Showing that $B \in \widehat{\mathcal{A}} \Leftrightarrow B^c \in \widehat{\mathcal{A}}$ is not easy. Proceed as follows. 

\noindent Consider $B \in \widehat{\mathcal{A}}$ with $B=A\cup N$, $A\in \mathcal{A}$, $N\in \mathcal{N}$. Since $N\in \mathcal{N}$, consider $D \mathcal{A}$, $m(D)=0$ and $N \subset D$.\\

\noindent Justify $N^c=D^c+D\cap N^c$.\\

\noindent Plug this into $B^c=A^c \cap N^c$ to get $B^c=(A^c\cap D^c)+N_0$. Identify $N_0$ and conclude that $B^c \in  \widehat{\mathcal{A}}$.\\

\noindent (b) Show that the application $\widehat{m}$ is well-defined by showing that if 
$$
B \in \widehat{\mathcal{A}}
$$

\bigskip \noindent with two expressions $B=A_1\cup N_1=A_2\cup N_2$, $A_i\in \mathcal{A}$, $N_i\in \mathcal{N}$, $i=1,2$, then
$$
m(A_1)=m(A_2).
$$

\bigskip \noindent Hint : Write $B=(A_1\cap A_2) \cup N_3$, identify $N_3$ and show that $N_3 \in \mathcal{N}$. Next show, by using $A_1=(A_1\cap A_2) \cup (A_1 \cap N_3)$ that $A_1=(A_1\cap A_2) + N_4$ and $N_4 \in \mathcal{N}$. Deduce that $m(A_1)=m(A_1\cap A_2)$. Do the same for $A_2$. Then show that $m(A_1)=m(A_2)$. Conclude.\\

\noindent (c) Show now that $\widehat{m}$ is a measure on $(\Omega, \widehat{\mathcal{A}})$, which is an extension of the measure $m$ on 
$(\Omega, \widehat{\mathcal{A}})$.

\noindent (d) Show  that the measure space  $(\Omega, \widehat{\mathcal{A}}, \widehat{m})$ is complete.\\

\noindent Hint. Consider a $\widehat{m}$-null set $B$ of the form $B=A \cap N$ such that $A\in \mathcal{A}$ and $N\in \mathcal{N}$, and $N \subset D$, $D$ being a $m$-null set. Show that $A$ is a 
$m$-null set and deduce from this, that $C\in \mathcal{N}$ and conclude.\\

\noindent \textbf{Conclusion}. We \textbf{may and do} consider any measure space $(\Omega, \mathcal{A},m)$ as complete, eventually by considering it as a subspace of its completition 
$(\Omega, \widehat{\mathcal{A}}, \widehat{m})$.\\

\newpage
\bigskip \noindent \textbf{Part IV : Determining classes for $\sigma$-finite measures}.\\

\bigskip \noindent \textbf{Exercise 9}. \label{exercise09_doc04-02}  Let $m_{1}$ and $m_{2}$ be two $\sigma$-finite measures on an algebra $\mathcal{C}$ of subsets of $\Omega$. Find a countable partition of $\Omega$ into elements of $\mathcal{C}$ in the form

$$
\Omega=\sum_{j\geq 0} \Omega_j, (\Omega_j \in \mathcal{C}, j\geq 0),
$$

\bigskip \noindent such that $m_{1}$ and $m_{2}$ are both finite on each $\Omega_j$, $j\geq 0$.\\

\bigskip \noindent \textbf{Exercise 10}. \label{exercise10_doc04-02} Let $(\Omega, \mathcal{A})$ a measurable space such that $\mathcal{A}$ is generated by an algebra $\mathcal{C}$ of subsets of $\Omega$. Show that two finite measure which are equal on $\mathcal{C}$, are equal (on $\mathcal{A}$).\\

\noindent Hint : Let $m_1$ and $m_2$ be two finite measures on $(\Omega, \mathcal{A})$ equal on $\mathcal{C}$. Set 
$$
\mathcal{M}=\{A \in \mathcal{A}, m_1(A)=m_2(A) \}.
$$

\bigskip \noindent Show that is a monotone class including $\mathcal{C}$. Apply \textit{Exercise 3 in DOC 01-04} to conclude.\\

\bigskip \noindent \textbf{Exercise 11}. \label{exercise11_doc04-02} (The $\lambda$-$\pi$ Lemma) Let $(\Omega, \mathcal{A})$ a measurable space such that $\mathcal{A}$ is generated by an $\pi$-system $\mathcal{P}$ containing $\Omega$. Show that two finite measures which equal on $\mathcal{P}$, are equal (on $\mathcal{A}$).

\noindent Hint : Let $m_1$ and $m_2$ be two finite measures on $(\Omega, \mathcal{A})$ equal on $\mathcal{P}$. Set 
$$
\mathcal{D}=\{A \in \mathcal{A}, m_1(A)=m_2(A) \}.
$$

\bigskip \noindent Show that is a $\lambda$-system (Dynkin system), which includes $\mathcal{P}$. Apply \textit{Exercise 4 in DOC 01-04} to conclude.\\

\bigskip \noindent \textbf{Exercise 12}. \label{exercise12_doc04-02} (Extension of Exercise 10). Let $(\Omega, \mathcal{A})$ a measurable space such that $\mathcal{A}$ is generated by an algebra $\mathcal{C}$ of subsets of $\Omega$. Show that two measures which are equal on $\mathcal{C}$ and $\sigma$-finite on $\mathcal{C}$, are equal (on $\mathcal{A}$).\\

\noindent Hint : Let $m_1$ and $m_2$ be two $\sigma$-finite on $\mathcal{C}$. Use a partition like to the one obtained in Exercise 8. Decompose each measure $m_i$, $i=1,2$ into

$$
m_i(A)=\sum_{j\geq 0} m_i(A\cap \Omega_{j}), \ A \in \mathcal{A}.
$$

\bigskip \noindent Next, apply Exercise 10 to conclude.\\

\newpage

\noindent \textbf{PART V : Extensions}.\\

\bigskip \noindent\textbf{Exercise 13}. \label{exercise13_doc04-02} Let $m$ be a proper and non-negative application defined from an algebra $\mathcal{C}$\ of subsets of $\Omega$ to $\overline{\mathcal{R}}_{+}$. Show that if $m$ is additive and $\sigma$-sub-additive with respect to $\mathcal{C}$, then $m$ is $\sigma$-additive on $\mathcal{C}$.\\

\bigskip \noindent\textbf{Exercise 14}. \label{exercise14_doc04-02} Let $m : \mathcal{C} \rightarrow \overline{\mathbb{R}}_{+}$, be a proper non-negative and finite application. Suppose that $m$ is finitely-additive on $\mathcal{C}$ and continuous at $\emptyset$, i.e., if $(A_{n})_{n\geq 1}$\ is a sequence of elements of $\mathcal{C}$ non-increasing to $\emptyset$, then $m(A_{n})\downarrow 0$, as $n\uparrow +\infty$. Show that $m$ is $\sigma$-additive on $\mathcal{C}$.\\

\bigskip \noindent\textbf{Exercise 15}. \label{exercise15_doc04-02} Let $\mathcal{S}$ be semi-algebra of subsets of $\Omega$. Let $m : \mathcal{S} \rightarrow \overline{\mathbb{R}}_{+}$ be a proper, non-negative and additive application. We already know that $\mathcal{C}=a(\mathcal{S})$ is the collection of all finite sums of elements of $\mathcal{S}$, that is :

$$
a(\mathcal{S})=\mathcal{C}=\{A_1+A_2+...+A_k, \ k\geq 1, \ A_1 \in \mathcal{S}, \ A_2 \in \mathcal{S}, ..., A_k \in \mathcal{S}\}.
$$

\bigskip \noindent Define the application $\widehat{m} : \mathcal{C} \rightarrow \overline{\mathbb{R}}_{+}$ as follows. For any $A \in \mathcal{C}$, with $A=A_1+A_2+...+A_k$, $k\geq 1$, $A_1 \in \mathcal{S}$, 
$ A_2 \in \mathcal{S}$, ..., $A_k \in \mathcal{S}$, set
 
$$
\widehat{m}(A)=\widehat{m}(A_1+A_2+...+A_k)=m(A_1)+m(A_2)+...+m(A_k), \ \ (DEF)
$$

\bigskip \noindent Question (1). Show that the application $\widehat{m}$ is a well-defined.\\

\noindent Question (2). Show that $\widehat{m}$ is an extension of $m$ as a non-negative, proper and additive application with the additive properties.\\ 

\noindent Question (3). Show that $m$ is $\sigma$-additive on $\mathcal{S}$ if and only if $\widehat{m}$ is $\sigma$-additive on $\mathcal{C}$.\\
 
\noindent Question (4). Show that $m$ is $\sigma$-sub-additive on $\mathcal{S}$ if and only if $\widehat{m}$ is $\sigma$-sub-additive on $\mathcal{C}$.\\

\noindent Question (5). Show that $m$ is $\sigma$-finite with respect to $\mathcal{S}$ if and only if  $\widehat{m}$ is $\sigma$-finite with respect to $\mathcal{C}$.\\

\newpage
\noindent \LARGE \textbf{DOC 04-03 : Carath\'{e}odory's Theorem and application to the Lebesgue-Stieljes Measures}. \label{doc04-03}\\
\Large

\bigskip \noindent This document states main construction tool of measure : the Caratheodory's Theorem. Here, what is important is to understand it and to learn how to use it, in the nontrivial case of Lebesgue measure's construction. It is recommended that the learners discuss it about in a session of about thirty minutes.\\

\bigskip \noindent \textbf{Part I : Statement of Carath\'{e}odory's Theorem}.\\

\bigskip \noindent Carath\'{e}odory's Theorem allows extension of the measure on the algebra to the $\sigma$-algebra it generates. It is stated as follows.\\

\bigskip \noindent\textbf{04.26 Main form of Caratheodory's Theorem}.\\

\bigskip \noindent Let $m$ be a proper and non-negative application defined from an algebra $\mathcal{C}$\ of subsets of $\Omega$ to $\overline{\mathbb{R}}_{+}$. If $m$ is 
$\sigma$-additive with respect to $\mathcal{C}$, then it is extensible to a measure $m^{0}$ on the $\sigma$-algebra $\mathcal{A}=\sigma(\mathcal{C})$ generated by $\mathcal{C}$. If $m$ is $\sigma$-finite on $\mathcal{C}$, the extension is unique and $\sigma$-finite.\\

\bigskip \noindent Some times, it is not easy to prove that the application $m$ is $\sigma$-additive with respect to $\mathcal{C}$. We use the following weaker form.\\

\bigskip \noindent\textbf{04.27 Weaker versions}.\\

\noindent Using \textit{Exercise 13 in Doc 04-02} allows to replace the $\sigma$-additivity of $m$ by the additivity and the $\sigma$-sub-additivity to give the second version below. As well, for finite measures, the results in \textit{Exercise 14 in Doc 04-02} lead to third version.\\

\noindent\textbf{04.27a Version 2}.\\

\bigskip \noindent Let $m$ be a proper and non-negative application defined from an algebra $\mathcal{C}$\ of subsets of $\Omega$ to $\overline{\mathcal{R}}_{+}$. If $m$ is additive and 
$\sigma$-sub-additive with respect to $\mathcal{C}$, then it is extensible to a measure $m^{0}$ on the $\sigma$-algebra $\mathcal{A}=\sigma(\mathcal{C})$ generated by $\mathcal{C}$. If $m$ is $\sigma$-finite on $\mathcal{C}$, the extension is unique and $\sigma$-finite.\\

\bigskip \noindent When $m$ is finite, we have to the third one.\\

\bigskip \noindent \textbf{04.27b Version 3}. Let $m : \mathcal{C} \rightarrow \overline{\mathbb{R}}_{+}$, be a proper non-negative and finite application. Suppose that $m$ is finitely-additive on $\mathcal{C}$ and continuous at $\emptyset$, i.e., if (A$_{n})_{n\geq 1}$\ is a sequence of subsets of $\mathcal{C}$ non-increasing to $\emptyset$, then $m(A_{n})\downarrow 0$, as $n\uparrow +\infty$,

\bigskip \noindent Then, $m$ extensible to a measure on $\sigma(\mathcal{C})$.\\

\bigskip \noindent\textbf{04.28 Construction of a measure from a semi-algebra}.\\

\noindent In some interesting cases, we may begin to construct a proper non-negative and additive application on a semi-algebra $\mathcal{S}$ and to extend it on the algebra $\mathcal{C}=a(\mathcal{S})$ generated $\mathcal{S}$ and further, to extend it to $\sigma$-algebra $\mathcal{A}$ generated by both $\mathcal{C}$ and $\mathcal{S}$.\\

\noindent The method is explained by the results in \textit{Exercise 14 in Doc 04-02}. A proper, non-negative and additive application $m$ is defined on $\mathcal{S}$. It is uniquely extended to  a proper, non-negative and additive application $\widehat{m}$ on $\mathcal{C}$. If $\widehat{m}$ is $\sigma$-sub-additive, it extensible to a measure on by the second version. But $\widehat{m}$ is  $\sigma$-sub-additive if and only if $m$ is. This gives the advanced version of Caratheodory's theorem.\\

\bigskip \noindent\textbf{Advanced form of Caratheodory's Theorem}. Let $\mathcal{S}$ be semi-algebra of subsets of $\Omega$. Let $m : \mathcal{S} \rightarrow \overline{\mathbb{R}}_{+}$ be a proper non-negative and additive. It it is $\sigma$-sub-additive, then it is extensible to a measure on $\mathcal{A}=a(\mathcal{S})$.\\

\noindent If $m$ is $\sigma$-finite with respect to $\mathcal{S}$, then the extension is unique and $\sigma$-finite.\\
   
 \newpage
\bigskip \noindent \textbf{Part II : Main application to the Lebesgue-Stieljes Measure on $\mathbb{R}$}.\\

\noindent Let us begin to define distribution functions.\\

\noindent \textbf{04.29 Definition} A non-constant function $F:\mathbb{R}$ $\longmapsto$ $\mathbb{R}$ is a  if and only if \\

\noindent (1) it assigns non-negative lengths to intervals, that is for any real numbers $a$ abd $b$ such that $a\leq b$, we have
$$
\Delta_{a,b} F=F(b)-F(a) \geq 0. \text{         } (\Delta F)
$$

\bigskip \noindent (2) it is right-continuous at any point $t\in \mathbb{R}$.\\

\noindent \textbf{Warning and Remark} The condition ($\Delta F$) incidentally means that $F$ is non-decreasing on $\mathbb{R}$. \textbf{But}, keep using the saying that \textit{assigns non-negative lengths to intervals} or \textit{intervals of $\mathbb{R}$ have non-negative $F$-lengths}. This approach will help you to move smoothly to higher dimension.\\

\noindent The application of the Caratheodory theorem gives the following\\

\noindent \textbf{Theorem-Definition}.\\ 

\noindent \textbf{04.30 (a)} For any distribution function $F$, there exists a unique measure on 
$(\mathbb{R}, \mathcal{B}(\mathbb{R}))$ such that such that any real numbers $a$ and $b$ satisfying $a\leq b$, we have

\begin{equation}
\lambda_{F}(]a,b])=\Delta F(]a,b])=F(b)-F(a).  \label{LS1}
\end{equation}

\bigskip \noindent This measure is called the Lebesgue-Stieljes measure associated with the distribution function $F$.\\

\noindent By applying this to the special case where $F(x)=x$, we get the Lebesgue measures on $(\mathbb{R}, \mathcal{B}(\mathbb{R}))$ as follows.\\

\noindent \textbf{04.30 (b)} The Lebesgue measure on $\lambda$, is the unique measure on $(\mathbb{R}, \mathcal{B}(\overline{\mathbb{R}}))$ such that any real numbers $a$ and $b$ satisfying $a\leq b$, we have 

\begin{equation}
\lambda(]a,b])=b-a.  \label{LS2}
\end{equation}

\noindent The Lebesgue measure on $\mathbb{R}$ is the measure that assigns to an interval its length.\\

\newpage
\bigskip \noindent \textbf{Part III : The product measure space $\mathbb{R}$}.\\

\noindent Consider a finite product space endowed with its product $\sigma$-algebra
$$
(\Omega, \mathcal{A}) \equiv \left(\prod_{1\leq i\leq k}\Omega_{i},\bigotimes_{1\leq i\leq k}\mathcal{A}_{i}\right).
$$

\bigskip \noindent Suppose we have $\sigma$-finite measure $m_i$ on each factor $(\Omega_i, \mathcal{A}_i)$.\\

\noindent There exists a unique $m$ measure on the product space $(\Omega, \mathcal{A})$ characterized by its values on the semi-algebra 
$\mathcal{S}$ of the measurable cylinders 

$$
A=A_1 \times A_2 \times \cdots \times A_k = \prod_{i=1}^{i=k} A_i,  \ A_i \in \mathcal{A}_i,
$$ 

\bigskip \noindent which are

\begin{equation}
m\left(\prod_{i=1}^{i=k} A_i\right) = \prod_{i=1}^{i=k} m_{i}\left(A_i\right). \label{prodMeasure}
\end{equation}

\bigskip \noindent The measure $m$ defined in Formula (\ref{prodMeasure}) has the following features.\\

\noindent (1) $m$ is called the product measure of the measures $m_i$, $i=1,...,k$.\\

\noindent (2) $m$ is unique and $\sigma$-finite if each measure $m_i$, $i=1,...,k$, is $\sigma$-finite.\\

\noindent (3) The application $m$ defined on $\mathcal{S}$ by Formula (\ref{prodMeasure}) is additive and then extensible to an additive application on 
$\mathcal{C}=\sigma(\mathcal{S})$\\

\noindent (4) The extension the application $m$, as defined on $\mathcal{S}$, to a measure on $\mathcal{A}$, is possible by using the Caratheodory's Theorem when the factors are $\mathbb{R}$, or some special cases.\\
 
\noindent (5) There exists a very beautiful and general proof using integration, even when the $m_i$ are not $\sigma$-finite. We will favor this proof that will be given in special chapter.\\

\noindent (6) The product measure is very important in probability theory in relation with the notion of independence.

\newpage
\noindent \LARGE \textbf{Doc 04-04 Measures - Exercises on Lebesgue-Stieljes Measure}. \label{doc04-04}\\
\Large

\bigskip \noindent \textbf{Exercise 1}. \label{exercise01_doc04-04} Let $(\overline{\mathbb{R}}, \mathcal{B}(\overline{\mathbb{R}}), \lambda)$ the Lebesgue measure space of $\mathbb{R}$.\\

\noindent For any Borel set $A$ in $\overline{\mathbb{R}}$, for $x\in \mathbb{R}$ and $0 \neq \gamma \in \mathbb{R}$, Define

$$
A+x =\{ y+x, \ y\in A\} \text{ and } \gamma A =\{ \gamma y, \ y\in A\}.
$$

\bigskip \noindent Define the translation application $t_x$ of vector $x$ defined by : $t_x(y)=x+y$ and the linear application $\ell_{\gamma}$ of coefficient $\gamma$ defined by $\ell_{\gamma}(y)=\gamma y$.\\

\noindent Question (a) Let $A$ be Borel set. Show that $A+x=t_{x}^{-1}(A)$ and $\gamma A=\ell_{1/\gamma}^{-1}(A)$ and deduce that $A+x$ and $\gamma A$ are Borel sets.\\

\noindent Question (b) Show that the Lebesgue measure is translation-invariant, that id for any $x\in \mathbb{R}$, for any Borel set $A$, we have

$$
\lambda(A+x)=\lambda(A).
$$

\bigskip \noindent \textit{Hint} : Fix $x\in \mathbb{R}$. Consider the application $A \hookrightarrow \lambda_{1}(A)=\lambda(A+x)$. Show that $\lambda_{1}$ is an inverse measure of the Lebesgue measure. Compare $\lambda_{1}$ and $\lambda$ of the sigma-algebra

$$
\mathcal{S}=\{]a,b], -\infty \leq a \leq b \leq +\infty\}.
$$

\bigskip \noindent Conclude.\\

\noindent Question c) Show for any $\gamma \in \mathbb{R}$, for any Borel set $A$, we have

$$
\gamma A)=|\gamma| \lambda(A).
$$

\bigskip \noindent \textit{Hint}. The case $\gamma=0$ is trivial. For $\gamma\neq 0$, proceed as in Question (a).\\

\bigskip \noindent \textbf{Exercise 2}. \label{exercise02_doc04-04} Let $(\mathbb{R}, \mathcal{B}(\mathbb{R}))$ be the Borel space.\\

\noindent Question (a) Consider the Lebesgue measure $\lambda$ on $\mathbb{R}$.\\

\noindent (a1) Show that for any $x \in \mathbb{R}$, $\lambda(\{x\})=0$.\\

\noindent (a2) Show that for countable $A$ subset of $\mathbb{R}$, $\lambda(\{A\})=0$.\\

\noindent Question (b) Consider a Lebesgue-Stieljes measure $\lambda_F$ on $\mathbb{R}$ where $F$ is a distribution function.\\

\noindent (b1) Show that for any $x \in \mathbb{R}$, 

$$
\lambda_F(\{x\})=F(x)-F(x-0),
$$

\bigskip \noindent where $F(x-0)$ is the left-hand limit of $F$ at $x$.\\

\noindent Deduce from this that $\lambda_F(\{x\})$ if and only if $F$ is continuous at $x$.\\

\noindent (b2) Show that for any $\lambda_F(\{x\})\neq 0$ for at most a countable number of real numbers.\\

\noindent Hint. Combine this with \textit{Exercise 1 in DOC 03-06}. \\

\bigskip \noindent \textbf{Exercise 3}. \label{exercise03_doc04-04}  The objective of this exercise is to show that the Lebesgue measure is the unique measure on $(\mathbb{R}, \mathcal{B}(\mathbb{R}))$ which is translation invariant and which assigns the unity value to the interval $]0,1]$.\\

\noindent Let $(\mathbb{R}, \mathcal{B}(\mathbb{R}), m)$ be a measure space.\\

\noindent Question (a) Show that if $m$ is translation-invariant and $m(]0,1])=1$, then $m$ is the Lebesgue measure $\mathbb{R}$.\\

\noindent Hint. Use the translation-invariance of $m$ to prove that $m(]0,1/q])=1/q$, for any integer $q\geq 1$. Next, prove that $m(]0,p/q])=p/q$ for $p\geq 1$, $q\geq 1$, $p/q\leq 1$. Next, Using the continuity of $m$ on $(0,1)$ and the density of the set of rational numbers in $\mathbb{R}$, show that $m(]0,b])=b$ for any real number $b\in (0,1)$. Finally, by translation-invariance of $m$, proceed to your conclusion.\\

\noindent Question (b) Deduce from Question (a) that : if $m$ assigns to open intervals positive values and is translation-invariant, then there exists $\gamma >0$, $m=\gamma \lambda$.\\

\noindent Hint. Put $\gamma=m(]0,1])$ and proceed.\\

\newpage
\bigskip \noindent \textbf{Exercise 4}. \label{exercise04_doc04-04}  \label{exercise_04_doc-04-04} Vitali's construction of a non measurable subset on $\mathbb{R}$.\\

\noindent In this exercise, a non-measurable set on $\overline{\mathbb{R}}$ is constructed. The exercise focuses on $]0,1]$ on which is uses the
circular translation, in place of the straight translation by x : $t\hookrightarrow t+x,$ on the whole real line.\\

\noindent Define the circular translation on $]0,1]$ for any fixed $x$, $0<x<1$ by

\begin{equation*}
\begin{tabular}{llll}
$T_{x}:$ & $]0,1]$ & $\longrightarrow $ & $]0,1]$ \\ 
& $y$ & $\hookrightarrow $ & $\left\{ 
\begin{tabular}{lll}
$z=y+x$ & if & $y+x\leq 1$ \\ 
z=$x+y-1$ & if & $y+x>1$%
\end{tabular}%
\right.$.
\end{tabular}
\end{equation*}

\bigskip \noindent and denote 

$$
T_{x}(A)=A\oplus x.
$$

\bigskip \noindent Question (a) Show that for any $x$, $0<x<1$, the application from  $]0,1]$ to $]0,1]$, defined by

\begin{equation*}
\begin{tabular}{llll}
$T_{-x}:$ & $]0,1]$ & $\longrightarrow $ & $]0,1]$ \\ 
& $z$ & $\hookrightarrow $ & $\left\{ 
\begin{tabular}{lll}
$y=z-x$ & if & $z-x>0$ \\ 
$y=z-x+1$ & if & $z-x\leq 0$%
\end{tabular}%
\right. $.
\end{tabular}
\end{equation*}

\bigskip \noindent is the inverse of $T_x$. We denote

$$
T_{-x}(A)=A\ominus x.
$$

\bigskip \noindent Question (b) Show that for $0<x<1$, for any Borel set in $]0,1]$, we have
$$
A\oplus x=T_{-x}^{-1}(A)
$$

\bigskip \noindent and conclude $A\oplus x$ is measurable and the graph

\begin{equation*}
A\hookrightarrow \lambda_{x}(A)=\lambda (A\oplus x)
\end{equation*}

\noindent defines a measure on the Borel sets in $]0,1]$.

\noindent Question (c) Show that $\lambda_{x}=\lambda$.\\

\noindent Hint. Consider an interval $]y, y^{\prime}]$ in $]0,1]$. Use the following table to give  $]y, y^{\prime}]\oplus x$ with respect to the positions of $y+x$ and $y^{\prime}+x$ with respect to the unity.

\begin{equation*}
\begin{tabular}{|l|l|l|l|l|}
\hline
case & y+x & $y^{\prime}+x$ & $]y,y^{\prime }]\oplus x$ & $\lambda (]y,y^{\prime }]\oplus x)$ \\ \hline
1 & $y+x\leq 1$ & $y^{\prime }+x\leq 1$ &  &  \\ 
\hline
2 & $y+x\leq 1$ & $y^{\prime }+x>1$ &  &  \\ 
\hline
3 & $y+x>1$ & $y^{\prime }+x\leq 1$ & impossible & impossible \\ 
\hline
4 & $y+x>1$ & $y^{\prime }+x>1$ &  & \\
\hline
\end{tabular}
\end{equation*}

\noindent \textbf{Extension}. By a very similar method, we also may see that $A\hookrightarrow \lambda_{-x}(A)=\lambda (A\ominus x)$ is a measure and that $\lambda_{-x}=\lambda$.\\

\noindent We may unify the notation by writing $A\ominus x=A\oplus (-x)$ for $-1<x<0$. We have proved that

\begin{equation*}
\forall (-1<x<1),\forall A\in B(]0,1]),\lambda (A\oplus x)=\lambda (A).
\end{equation*}

\bigskip \noindent Question (d).\\

\noindent Consider the binary relation $\mathcal{R}$ on $]0,1]$ defined by%
\begin{equation*}
\forall (y,z)\in ]0,1]^{2},y\mathcal{R}z\Longleftrightarrow y-z\in \mathbb{Q} \cap ]0,1],
\end{equation*}

\bigskip \noindent that is $y$ and $z$ are in relation if and only if their $y-z$ difference is
a rational number, that is also : there exists $r\in \mathbb{Q},$ $y=z+r$
with $-1<r<1$.

\noindent (d1) Show that $\mathcal{R}$ is an equivalence relation.\\

\noindent Form a set $H$ by choosing one element in each equivalence class (apply the choice axiom).\\

\noindent (d2) Show that for any $0\leq r<1,H\oplus r\subset ]0,1]$.\\

\noindent (d3) Show that we have

\begin{equation*}
]0,1]\subset \sum\limits_{r\in Q\cap ]0,1[}H\oplus r.
\end{equation*}

\bigskip \noindent (d4) Show that for $0\leq r_{1}\neq r_{2}<1,H\oplus r_{1}$ and $H\oplus r_{2}$ are disjoint, that is
\begin{equation*}
(H\oplus r_{1})\cap (H\oplus r_{2})=\emptyset .
\end{equation*}

\bigskip \noindent Conclude that 

\begin{equation*}
]0,1]=\sum\limits_{r\in Q\cap ]0,1[}H\oplus r.
\end{equation*}

\bigskip \noindent (d5) Suppose that $H$ is measurable. Apply the $\sigma$-additivity of $\lambda$ and show that that this leads to an absurdity.\\

\noindent Conclusion. Did you get a nonempty non-measurable set on $\mathbb{R}$?\\

\newpage
\noindent \LARGE \textbf{DOC 04-05 : Measures - General exercises}. \label{doc04-05}\\
\Large
 
\bigskip \noindent \textbf{Exercise 1}. \label{exercise01_doc04-05} Let $(E, \mathcal{B}(E),m)$ be a finite measure ($(E, \mathcal{B}(E))$ is a Borel space. Show that any other measure that is equal to $m$ on the class of of open sets, or on the class of closed sets is equal to $m$ on $\mathcal{B}(E)$. Conclusion : A finite measure on a Borel space us regular, $i.e$, it is determined by its values on the open sets, by its values on the closed sets.\\

\bigskip \noindent \textbf{Exercise 2}. \label{exercise02_doc04-05}  (Approximation of a finite measure by its values on open and closed sets in a metric space). Let $(S,d,m)$ be a finite measure on the metric space $(S,d)$ endowed with its Borel $\sigma$-algebra.  \noindent Let $A \in \mathcal{B}(S)$. Show the assertion $\mathcal{P}(A)$ : for any $\varepsilon>0$, there exist a closed set
$F$ and an open set $G$ such that : $F \subset A \subset G$ and $m(G \setminus F)<\varepsilon$.\\

\noindent \textit{Hint}. Proceed as follows.\\

\noindent Question (a) Let $A$ be closed set. Fix $\delta>0$ and set  the $\delta$-dilation or $\delta$-enlargement of $G$ as 

\begin{equation*}
G_{\delta}=\{x,\text{ }d(x,A)<\delta \}.
\end{equation*}

\bigskip \noindent (a) Show that $G_{\delta}$ is an open set (or you may skip this as a result of topology).\\

\noindent Let $(\delta_n)_{n\geq 1}$ be sequence of positive number such that $\delta_n \downarrow +\infty$ as $n \uparrow+\infty$.\\

\noindent (b) By using the continuity of $m$ with the sequence $(G_{\delta_n})_{n\geq 1}$, and by noticing that $A \subset G_{\delta_n}$ for any $n$, show that for any $\varepsilon>0$, there exists $\delta>0$ such that

\begin{equation*}
m(G_{\delta} \setminus A)<\varepsilon.
\end{equation*}

\bigskip \noindent (c) Define by $\mathcal{D}$ the class of measurable sets $A$ such that for any $\varepsilon>0$, there exist a closed set
$F$ and an open set $G$ such that : $F \subset A \subset G$ and $m(G \setminus F)<\varepsilon$.\\

\noindent Show that $\mathcal{D}$ is a $\sigma$-algebra including the class $\mathcal{F}$ of closed sets. Conclude that the assertion $\mathcal{P}$  holds.\\

\noindent Hint. To show that $\mathcal{D}$ is stable under countable union, consider a sequence $A_n \in \mathcal{D}$, $n\geq 0$ and next
$$
A =\bigcup_{n\geq 0} A_n. 
$$

\bigskip \noindent Fix $\varepsilon>0$. By definition, for each $n\geq 0$, there exist a closed set$H_n$ and an open set $G_n$ such that : $H_n \subset A_n \subset G_n$ and $m(H_n \setminus F_n)<\varepsilon/2^{n+1}$.\\

\noindent Denote

$$
F_p=\bigcup_{1\leq n \leq p} H_n, \ p\geq 0, \ H=\bigcup_{n\geq 0} H_n \text{ and } G=\bigcup_{n\geq 0} G_n. 
$$

\bigskip \noindent By using the continuity of $m$, explain that there exists $p_0\geq 0$ show that
$$
m(H \setminus F_{p_0}) < \varepsilon.
$$

\bigskip \noindent Set $F=F_{p_0}$ and justify :

$$
G \setminus F= (H \setminus F) + (G \setminus H).
$$

\bigskip \noindent Finally, use \textit{Formula 00.12 in Doc 00-01} or \textit{Exercise 1, Doc 00-03} in Chapter 0 to conclude.\\

\bigskip \noindent \textbf{Exercise 3}. \label{exercise03_doc04-05}  (Borel-Cantelli Lemma) Let $(\Omega,  \mathcal{A}, m)$ be a measure space.\\

\noindent Question (a) Let $(A_n)_{n\geq 0}$ be a sequence of measurable sets on $(\Omega,  \mathcal{A})$ such that

$$
\sum_{n\geq 0} m(A_n) <+\infty.
$$

\bigskip \noindent Show that $\limsup_{n \rightarrow +\infty} A_n$ is a $m$-null set, that is, 

$$
m\left(\limsup_{n \rightarrow +\infty} A_n\right)=0.
$$ 

\bigskip \noindent \textit{Hint} : Recall

$$
\limsup_{n \rightarrow +\infty} A_n=\lim_{n\uparrow +\infty} \downarrow B_n,
$$

\bigskip \noindent with

$$
B_n=\bigcup_{p\geq n} A_p.
$$

\bigskip \noindent Use the $\sigma$-sub-additivity of $m$ in $m(B_n)$ and use a remarkable property of convergent series.\\

\noindent Question (b) Suppose that the measure is finite with $0<m(\Omega)=M<+\infty$. Let $(A_n)_{n\geq 0}$ be a sequence of measurable sets on $(\Omega,  \mathcal{A})$ such that :

$$
\forall \ n\geq 0, \ \lim_{r\rightarrow +\infty} M^r exp(-\sum_{p=n}^{r} m(A_p))=0.  \label{  (SUMINF) }
$$ 

\bigskip \noindent and for any integers $0\leq r <s$,

$$
m\left( \bigcap_{n=r}^{r=s} A_{n}^c\right) \leq  \prod_{n=r}^{r=s} m(A_{n}^c).  \label{  (SUBIND) }
$$

\bigskip \noindent Show that 

$$
m\left(\limsup_{n \rightarrow +\infty} A_n\right)=M.
$$ 

\noindent \textbf{Useful remark}. In Probability Theory, the condition ((SUBIND) is implied by an independence condition and the condition (SUMINF) is equivalent to
$$
\sum_{n\geq 0} m(A_n)=+\infty, \ M=1.
$$

\bigskip \noindent \textit{Hint}. Show that $m(\liminf_{n \rightarrow +\infty} A_{n}^{c})=0$. Define

$$
C_n=\bigcap_{n \geq p} A_p^c, n\geq 0
$$

\bigskip \noindent and for $n\geq 0$ fixed and $r \geq n$,

$$
C_{n,r}=\bigcap_{p=n}^{p=r} A_p^c.
$$

\bigskip \noindent Remark that $C_n \nearrow \liminf_{n \rightarrow +\infty} A_{n}^{c}$ and for $n\geq 0$ fixed, $C_{n,r} \nearrow C_n$.

\noindent Combine (SUBIND) with the inequality : $1-x\leq e^{-x}$, $0\leq \leq <1$ to get

$$
m(C_{n,r}) \leq = M^{-n} M^{r} \exp\left( - M^{-1} \sum_{n=r}^{r=s} m(A_{n}\right).
$$

\bigskip \noindent Continue and conclude.\\

\newpage
\noindent \LARGE \textbf{Doc 04-06 : Introduction to Measures - Exercises with solutions}. \label{doc04-06}\\
\Large

\noindent \textbf{PART A : General properties on measures and additive sets applications}.\\

\bigskip \noindent \textbf{Exercise 1}. \label{exercise01_sol_doc04-06} Let $\mathcal{A}$ be $\sigma $-algebra of subsets
of $\Omega $. Let  $m$ be an application from $\mathcal{A}$ to $\overline{\mathbb{R}}_{+}$. Suppose that we have :\\

\noindent (ii) For all sequence $(A_{n})_{(n\geq 0)}$ of pairwise
disjoint elements of $\mathcal{A}$, the $\sigma$-additivity formula holds 
\begin{equation*}
m(\sum_{n\geq 0}A_{n})=\sum_{n\geq 0}m(A_{n})
\end{equation*}

\bigskip \noindent holds. Show that the two assertions are equivalent :\\

\noindent (i1) $m(\emptyset)$=0.\\

\noindent (i2) $m$ is a proper application constantly equal to $+\infty$, that is : 
$$
\exists A \in \mathcal{A}, \ m(A)<+\infty.
$$

\bigskip \noindent Show why the results you obtain still hold if $m$ is only finitely additive on an algebra.\\

\noindent \textit{Conclusion}. The results of the exercise justify the equivalence of the definitions of a measure in \textit{Points (04-01a) and (04-01b)}.\\

\noindent \textbf{SOLUTIONS}. We have to prove the equivalence between $\{(i1), (ii)\}$ and $\{(i2), (ii)\}$. But $(i1)$ obviously implies $(i2)$ and hence $\{(i1), (ii)\} \Rightarrow \{(i2), (ii)\}$. It remains to prove that $\{(i2), (ii)\} \Rightarrow \{(i1), (ii)\}$. Now suppose that $\{(i2)$ and $(ii)\}$. We only need to prove that $(i1)$ holds. But by $(ii)$, there exists $A$ in $\mathcal{A}$ such that
$m(A)$ is finite. Since $A$ and $\emptyset$ are two elements of $\mathcal{A}$ and $m$ is $\sigma$-additive and hence additive, we we have

$$
m(A)=m(A+\emptyset)=m(A)+m(\emptyset).
$$

\bigskip \noindent Since $m(A)$ is finite, we get

$$
m(\emptyset)=m(A)-m(A)=0.
$$

\bigskip \noindent So we have the equivalence between $\{(i1), (ii)\}$ and $\{(i2), (ii)\}$. We only use the finite additivity. So the result is also valid in the present from : For an additive application on an algebra, the assertion that $m$ is proper is equivalent to $m(\emptyset)=0$. $\blacksquare$\\

\bigskip \noindent \textbf{Exercise} 2. \label{exercise02_sol_doc04-06} Let $(\Omega,\mathcal{A},m)$ be measure space. Show the following properties.\\

\noindent (1)  $m$ is non-decreasing : If $A\subset B$, then $m(A)\leq m(B)$. Besides, show, if $m(B)$ is finite, that $m(B \setminus A)=m(B)-m(A \cap B)$.\\ 

\noindent Hint : Use and justify $B = (A\cap B) + (B\setminus A)$ and use the additivity and the non-negativity of a measure.\\

\noindent (2) A measure $m$ finite if and only if only if $m(\Omega )$ is finite if and only if $m$ is bounded.\\

\noindent Hint : Use the non-decreasingness proved in Question (1).\\

\noindent (3) A measure $m$ is $\sigma $-sub-additive.\\

\noindent In a similar manner, show that if $m(\emptyset)=0$, and $m$ is non-negative and additive, then it is sub-additive

\noindent Hint : Use the formula if Exercise 4 in Chapter 0, DOC 00-02.\\

\noindent (4) Show that $m$ is continuous in the following sense :\\

\bigskip \noindent \textbf{(a)} Let $(A_{n})_{(n\geq 0)}$ be a non-decreasing sequence
of elements of $\mathcal{A}$ and set $A=\bigcup\limits_{n\geq 0}A_{n}$. Then 
$m(A_{n})$ increases to $m(A)$ as $n\rightarrow +\infty$.\newline

\bigskip \noindent \textbf{(b)} Let $(A_{n})_{(n\geq 0)}$ be a non-increasing sequence of
elements of $\mathcal{A}$ and set $A=\bigcap\limits_{n\geq 0}A_{n}$ such that $m(A_{n_{0}})$ for some $n_{o}\geq 0$. Then $m(A_{n})$ decreases to $m(A)$ as $n\rightarrow +\infty$.\\

\noindent \textit{Hint} : (1) To show (a), remark by the help of a simple diagram, that if $(A_{n})_{(n\geq 0)}$ is a non-decreasing sequence, we have

$$
\bigcup\limits_{n\geq 0}A_{n} = A_0 + \sum_{n\geq 1} \left(A_n\setminus A_{n-1}\right).
$$ 

\bigskip \noindent Next, use the $\sigma$-additivity of $m$ and then the partial sums of the obtained series.\\

\noindent To prove (b), use Point (a) by taking the complements of the $A_n$'s.\\

\bigskip \noindent \textbf{SOLUTIONS}.\\

\noindent \textbf{(1)}. Let $(A,B) \in \mathcal{A}^2$. With the help of a simple diagram, we can see that $A = (A\cap B) + (B\setminus A)$. By additivity, we have
$$
m(B)=m(A\cap B) + m(B\setminus A).
$$

\bigskip \noindent From that formula, we see that, if $A \subset B$, we have $m(B)=m(A) + m(B\setminus A)$ and next, $m(B)\geq m(A)$ since $m(B\setminus A)\geq 0$.\\

\noindent Next, we mat derive from it, if $m(B)$ is finite, that $m(B\setminus)$ is finite by the first conclusion and finally, we have

$$
m(B\setminus A)=m(B) - m(A\cap B).
$$

\bigskip \noindent \textbf{(2)}. Let $m$ be finite. Then $m(\Omega)$ is finite and, by Point (1) of this Exercise, we have for all $A \in \mathcal{A}$, $m(A)\leq  m(\Omega)$. Hence, $m$ is bounded by 
$m(\Omega)$.\\

\noindent In the other side, if $m$ is bounded, say by $M<+\infty$, we get that $m(A)\leq  M$ for all $A \in \mathcal{A}$. Hence $m$ is finite.\\

\noindent Conclusion : The measure $m$ finite if and only if it is bounded, and in that case, $m(\Omega)$ is the maximum value of $m$.\\

\bigskip \noindent \textbf{(3)}. Let $(A_n)_{n\geq 0}$ be a sequence of elements of $\mathcal{A}$. By the results of \textit{Exercise 4 in Chapter 0, Doc 00-02}, we have
$$
\sum_{n\geq 0} \ A_{n}=\sum_{n\geq 0} \ B_{n}
$$

\bigskip \noindent with

$$
B_0=A_0, \ B_n=A_{0}^{c}...A_{n-1}^{c}A_{n} \subset A_n, \ n\geq 1.
$$

\bigskip \noindent By $\sigma$-additivity and next by non-decreasingness of $m$, we have

\begin{eqnarray*}
m(\sum_{n\geq 0} \ A_{n})&=&m(\sum_{n\geq 0} \ B_{n})\\
&=&\sum_{n\geq 0} \ m(B_{n})\\
&\leq & \sum_{n\geq 0} \ m(A_{n}),
\end{eqnarray*}

\bigskip \noindent which concludes the solution of this point. $\blacksquare$\\

\bigskip \noindent \textbf{(4)}. \textbf{Part (a)}.\\

\noindent Let $(A_n)_{n\geq 0}$ be a non-decreasing sequence of elements of $\mathcal{A}$ and let us denote $A=\bigcup_{n\geq 0} A_n$. Let us use again the results of \textit{Exercise 4 in Chapter 0, Doc 00-02}, to have
$$
\sum_{n\geq 0} \ A_{n}=\sum_{n\geq 0} \ B_{n}
$$

\bigskip \noindent with

$$
B_0=A_0, \ B_n=A_{0}^{c}...A_{n-1}^{c}A_{n} \subset A_n, \ n\geq 1.
$$

\bigskip \noindent For $n\geq 1$, by the non-decreasingness of the $A_n$'s, all the $A_{j}^{c}$'s, $j=0,..,n-1$, contains $A_{n-1}^{c}$, so that we have for all $n\neq 1$ 

$$
B_0=A_0, \ B_n=A_{n-1}^{c}A_{n}=A_n \setminus A_{n-1} \subset A_n, \ n\geq 1.
$$

\bigskip \noindent and we get

$$
A=A_0 + \sum_{n\geq 1} (A_n \setminus A_{n-1}),
$$

\bigskip \noindent which we write as

$$
A=A_0 + \sum_{k\geq 1} (A_k \setminus A_{k-1}).
$$

\bigskip \noindent We may have directly derived this formula with the help of a simple diagram. Hence, we get by $\sigma$-additivity and next by partial sums,

$$
m(A)=m(A_0) + \sum_{n\geq 1} m(A_k \setminus A_{k-1})=\lim_{n \rightarrow +\infty} (m(A_0) + \sum_{k=1}^{n} m(A_k \setminus A_{k-1})).
$$

\bigskip \noindent We discuss over two possible cases.\\

\noindent \textbf{Case 1 : all the $A_k$'s, $k\geq 1$, have finite measures}. We remark that for each $k\geq 1$, we have $A_{k-1}\subset A_k$ and $m(A_k)$ is finite. Then by Point (1) of this exercise and
by denoting $a_n=m(A_0) + \sum_{k=1}^{n} m(A_k \setminus A_{k-1})$, we get

\begin{eqnarray*}
a_n&=&m(A_0) + \sum_{k=1}^{n} m(A_k \setminus A_{k-1})\\
&=&m(A_0) + \sum_{k=1}^{n} m(A_k) - m(A_{k-1})\\
&=&m(A_0) + (m(A_1)-m(A_0)) + (m(A_2)-m(A_1))+\cdots +(m(A_n)-m(A_{n-1}))\\
&=&m(A_n).
\end{eqnarray*}

\bigskip \noindent Hence, by combining the two latter formulas, we have

$$
m(A)=\lim_{n \rightarrow +\infty} m(A_n).
$$

\bigskip \noindent \textbf{Case 2 : there exists $n_0\geq 0$ such that $m(A_{n_0})=+\infty$}. It follows that $m(A)=+\infty$ since $A_{n_0} \subset A$ and also, $m(A_n)=+\infty$ for all $n\geq n_0$ for the sale reason. We have that :

$$
\exists n_0\geq 0, \ \ \forall \ n\geq n_0, \ m(A_n)=m(A)=+\infty.
$$

\bigskip \noindent This simply means that $m(A)=\lim_{n \rightarrow +\infty} m(A_n)$.\\

\noindent In both cases, we have $m(A)=\lim_{n \rightarrow +\infty} m(A_n)$ and the sequence $(m(A_n))_{n\geq 0}$ is non-decreasing. This concludes the proof of Part (a) of Point (4).$\square$\\

\noindent \textbf{Part (b)}. Let $(A_n)_{n\geq 0}$ be a non-increasing sequence of elements of $\mathcal{A}$ and let us denote $A=\bigcap_{n\geq 0} A_n$ and let $n_0$ be an integer such that $m(A_{n_0})$ is finite. Consider the sequence $(A_{n_0} \setminus A_{n})_{(n\geq n_0)}$ which is non-decreasing, since for each $n\geq n_0$ we have $A_{n_0} \setminus A_{n}=A_{n_0}A_{n}^{c}$, to 
$A_{n_0} \setminus A$. By applying Part (a) of this Point, we have

$$
m(A_{n_0} \setminus A_{n}) \nearrow m(A_{n_0} \setminus A), \ as \ n \nearrow +\infty.
$$

\bigskip \noindent which, by non-increasingness and by Point (1) of the this exercise,

$$
m(A_{n_0}) - m(A_{n}) \nearrow m(A_{n_0}) - m(A), \ as \ n \nearrow +\infty.
$$

\bigskip \noindent and by dropping the finite number $m(A_{n_0})$ and by taking the opposite numbers, we get

$$
m(A_{n}) \searrow  m(A), \ as \ n \nearrow +\infty.
$$

\bigskip \noindent The solution is complete. $\square$\\

\noindent \bigskip \noindent \textbf{Exercise 3}. \label{exercise03_sol_doc04-06}  Let $\mathcal{C}$ be an algebra of subsets of $\Omega$. Let $m$ be a proper and additive application from $\mathcal{C}$ to $\overline{\mathbb{R}}_{+}$.\\

\noindent Show that the two following assertions are equivalent.\\

\noindent (a) For any sequence $(A_{n})_{(n\geq 0)}$ of elements of $\mathcal{C}$ \textbf{such that} $\bigcup_{n\geq 0}A_{n}\in \mathcal{C}$, we have 
\begin{equation*}
m(\bigcup\limits_{n\geq 0}A_{n})\leq \sum_{n\geq 0}m(A_{n}).
\end{equation*}

\bigskip \noindent (b) For all sequence $(A_{n})_{(n\geq 0)}$ of pairwise disjoint elements of $\mathcal{C}$ \textbf{such that} $\sum_{n\geq 0}A_{n}\in \mathcal{C}$,
we have 
\begin{equation*}
m(\sum_{n\geq 0}A_{n})\leq \sum_{n\geq 0}m(A_{n}).
\end{equation*}

\bigskip \noindent The implication $(a) \Rightarrow (b)$ is needless to prove since $(b)$ is a special case of $(a)$. Let $(b)$ holds. Let $(A_{n})_{(n\geq 0)}$ be a sequence of elements of $\mathcal{C}$ \textbf{such that} $A = \bigcup_{n\geq 0}A_{n}\in \mathcal{C}$.  Then by \textit{Exercise 4 in Chapter 0, Doc 00-02}, we have 
$$
A = \sum_{n\geq 0} \ A_{n}=\sum_{n\geq 0} \ B_{n}
$$

\bigskip \noindent with

$$
B_0=A_0, \ B_n=A_{0}^{c}...A_{n-1}^{c}A_{n} \subset A_n, \ n\geq 1.
$$

\bigskip \noindent Hence the sequence $(B_{n})_{(n\geq 0)}$ remains in $\mathcal{C}$ since $\mathcal{C}$ is an algebra. By applying Part (b) and the non-decreasingness of $m$ due to its addivitity and to its
non-negativity, we have

\begin{eqnarray*}
m(\sum_{n\geq 0} \ A_{n})&=&m(\sum_{n\geq 0} \ B_{n})\\
&\leq &\sum_{n\geq 0} \ m(B_{n})\\
&\leq & \sum_{n\geq 0} \ m(A_{n}),
\end{eqnarray*}

\bigskip \noindent which puts an end to the solution.$\blacksquare$\\

\bigskip \noindent \textbf{Exercise 4}. \label{exercise04_sol_doc04-06}  Let $\mathcal{C}$ be an algebra of subsets of $\Omega$.  Let $m$ be a proper non-negative and additive application from $\mathcal{C}$ to $\overline{\mathbb{R}}_{+}$. Show that the two assertions are equivalent.\\

\noindent (a) There exists a finite or a countable partition of $\Omega$ into elements of $\mathcal{F}_0$, that is a sequence $(\Omega _{j})_{j\geq 0}$ of pairwise disjoint elements of $\mathcal{F}_0$ satisfying 
\begin{equation*}
\Omega =\sum_{j\geq 0}\Omega _{j} \ \ \ (CSF1),
\end{equation*}

\bigskip \noindent such that for all $j\geq 0$, $m(\Omega_j)<+\infty$.\\

\noindent (b) There exists a finite or a countable cover of $\Omega$ by elements of $\mathcal{C}$, that is a sequence $(\Omega _{j})_{j\geq 0}$ of elements of $\mathcal{C}$ satisfying  
\begin{equation*}
\Omega =\bigcup_{j\geq 0}\Omega _{j}  \ \ \ (CSF2),
\end{equation*}

\bigskip \noindent such that for all $j\geq 0$, $m(\Omega_j)<+\infty$.\\

\bigskip \noindent The implication $(a) \Rightarrow (b)$ is needless to prove since $(b)$ is a special case of $(a)$. Let $(b)$ holds. Let $(\Omega _{j})_{j\geq 0}$ be a sequence of elements of 
$\mathcal{C}$ such that $(CSF2)$ holds. Then by \textit{Exercise 4 in Chapter 0, Doc 00-02}, we have 
$$
\Omega = \sum_{n\geq 0} \ \Omega_j=\sum_{j\geq 0} \ \Omega^{\prime}_{j}
$$

\bigskip \noindent with

$$
\Omega^{\prime}_{0}=\Omega^{\prime}_{0}, \ \Omega^{\prime}_{j}=\Omega_{0}^{c}...\Omega_{j-1}^{c}\Omega_{j} \subset \Omega_j, \ j\geq 1.
$$ 

\bigskip \noindent and we also have : 

$$
\forall \ j>0, \ m(\Omega^{\prime}_{j}) \leq m(\Omega_{j}) <+\infty.
$$

\bigskip \noindent The solution is complete. $\square$\\

\newpage \noindent \textbf{PART II : Measures}.

\bigskip \noindent \textbf{Exercise 5}. \label{exercise05_sol_doc04-06} \\ 

\noindent (1) Let $(m_i)_{i\in I}$ a family of measures on $(\Omega, \mathcal{A})$, where $I\subset \mathbb{N}$. Let $(m\alpha_i)_{i\in I}$ a family of non-negative real numbers such that one of them is positive. Consider the application $m$ defined for  $A \in \mathcal{A}$ by

$$
m(A)=\sum_{i \in I} \alpha_{i}m_{i}(A).
$$

\bigskip \noindent Show that $m=\sum_{i \in I} \alpha_{i}m_{i}$ is a measure on $(\Omega, \mathcal{A})$. Give a conclusion.\\

\noindent \textit{Hint} : use the Fubini's rule for sums non-negative series.\\

\bigskip 
\noindent (2) Consider a measurable space $(\Omega, \mathcal{A})$, one point $\omega_0$ of $\Omega$ and a non-negative real number $\alpha$, show that the application $\delta_{\omega_0,\alpha}$ defined on $\mathcal{A}$ by

$$
\delta_{\omega_0,\alpha}(A)=\alpha 1_{A}, \ \ A \in \mathcal{A},
$$

\bigskip \noindent is a measure, called the delta-measure concentrated at $\omega_0$ with mass $\alpha$.\\

\noindent if $\alpha=1$, we denote

$$
\delta_{\omega_0,\alpha} \equiv \delta_{\omega_0}.
$$

\bigskip \noindent (3) Let $\Omega=\{\omega_1, \omega_2, ...\}$ be a countable space endowed with the $\sigma$-algebra $\mathcal{P}(\Omega)$, which is the power set of $\Omega$. By using Points (1) and (2) below, to show that the application defined by
$$
\nu = \sum_{j\geq 1} \delta_{\omega_j}
$$

\bigskip \noindent is a measure. Show that for any $A \subset \Omega$, we have

$$
\nu(A) = Card(A)
$$

\bigskip \noindent Suggest a name for this measure.

\bigskip\noindent \textbf{SOLUTIONS}.\\

\noindent \textbf{Question (1)}.  It is clear that $m$ non-negative as a linear combination of non-negative application with non-negative coefficient. It is also  proper since

$$
m(\emptyset) = \sum_{i\in I} \alpha_i m_i(\emptyset)= \sum_{i\in I} \ 0 = 0.
$$

\bigskip \noindent It remains to show that $m$ is $\sigma$-additive. Le $(A_{n})_{(n\geq 0)}$ be a sequence of pairwise disjoint elements of $\mathcal{A}$. We have

$$
m(\sum_{n \geq 0} A_n) = \sum_{i\in I} \alpha_i  m_i(\sum_{n \geq 0} A_n),
$$

\bigskip \noindent which gives, since each $m_i$ is a measure,

$$
m(\sum_{n \geq 0} A_n) = \sum_{i\in I} \alpha_i  \sum_{n \geq 0} m_i(A_n),
$$

\bigskip \noindent which, by the Fubini's formula for non-negative series, leads to

$$
m(\sum_{n \geq 0} A_n) = \sum_{n \geq 0} \left(\sum_{i\in I} \alpha_i  m_i(A_n)\right).
$$

\bigskip \noindent By using the definition of $m$, we arrive at

$$
m(\sum_{n \geq 0} A_n) = \sum_{n \geq 0} m(A_n).
$$

\bigskip \noindent \textbf{Question (2)}. It is also clear that application $\delta_{\omega_0,\alpha}$ is non-negative and it is proper since

$$
\delta_{\omega_0,\alpha}(\emptyset)= \alpha 1_{\emptyset}(\omega_0)=0.
$$

\bigskip \noindent To prove it is $\sigma$-finite, consider a sequence $(A_{n})_{(n\geq 0)}$ of pairwise disjoint elements of $\mathcal{A}$. We have to prove that

$$
\delta_{\omega_0,\alpha}(\sum_{n \geq 0} A_n) = \sum_{n \geq 0} \delta_{\omega_0,\alpha}(A_n). (CMS)
$$

\bigskip \noindent We face two possibilities.\\

\noindent \textit{Case 1} : $\omega_0$ is in any of the $A_n$'s. Hence, by definition of $\delta_{\omega_0,\alpha}$ for any $n\geq 0$, $\delta_{\omega_0,\alpha}(A_n)=0$. In that case, $\omega_0$ is not in $\sum_{n \geq 0} A_n$ and thus, $\delta_{\omega_0,\alpha}(\sum_{n \geq 0} A_n)=0$. We conclude that Formula (CMS) below holds since both members are zeros.\\

\noindent \textit{Case 2} : $\omega_0$ is in one of the $A_n$'s. Since the $A_n$'s are disjoint, this means that it in one and only one of the $A_n$'s. Let us name it $n_0$. This implies, by definition, that
$\delta_{\omega_0,\alpha}(A_{n_0})=\alpha$ and that for any $n\neq n_0$ and $n\geq 0$, we have $\delta_{\omega_0,\alpha}(A_n)=0$. Hence the left-member of Formula (CMS) us $\alpha$. \\

\noindent In that case, we surely have that $\omega_0$ is in $\sum_{n \geq 0} A_n$ and thus, $\delta_{\omega_0,\alpha}(\sum_{n \geq 0} A_n)=\alpha$. We see that the two members of the equality in Formula 
(CMS) are equal to $\alpha$.\\

\noindent In conclusion, Formula (CMS) holds and thus $\delta_{\omega_0,\alpha}$ is $\sigma$-additive and becomes a measure.\\

\noindent \textbf{Question (3)}. By combining the two first points, we see that $\nu$ is a measure. And we have for any part $A$ of $\Omega$,

$$
\nu(A) = \sum_{j\geq 1} 1_A(\omega_j).
$$

\bigskip \noindent It is clear that $\nu(A)$ is the cardinality of $A$ since $A$ is a subset of $\{\omega_1, \omega_1, ...\}$ and the formula just above is a count of the $\omega_j$ in $A$.\\

\noindent The application $\nu$ is called the counting measure on the countable space $\Omega$ endowed with the power set as a $\sigma$-algebra.$\blacksquare$\\

\bigskip \noindent \textbf{Exercise 6}. \label{exercise06_sol_doc04-06} Let $(\Omega, \mathcal{A},m)$ be a measure space.\\

\noindent (a) Let $X : (\Omega, \mathcal{A},m) \mapsto (E, \mathcal{B})$ be a measurable application from $(\Omega, \mathcal{A},m)$ to the measure space $(E, \mathcal{B})$. Show that the application defined for $B\ \mathcal{B}$ by
$$
m_X(B) = m(X^{-1}(B)), \ \ B \in \mathcal{B},
$$

\bigskip \noindent is a measure on $(E, \mathcal{B})$ called image measure (of $m$ by $X$).\\

\noindent (b)  Let $A$ be a measurable set on the measure space $(\Omega, \mathcal{A},m)$. Show that the application $m_A$ defined on $(A,\mathcal{A}_A)$ by
$$
m_A(B) = m(A\cap B), \ \ B \subset A,
$$

\bigskip \noindent is a measure on $(A,\mathcal{A}_A)$.\\

\bigskip \noindent \textbf{SOLUTIONS}.\\

\noindent Part (a). This exercise relies on the properties of the inverse image $X^{-1}$ as listed in \textit{Exercise 13, Chapiter 2, Doc 01-02}.\\

\noindent It is clear that $m_X$ is non-negative and we also have that $m_X(\emptyset)=m(X^{-1}(\emptyset))=m(\emptyset)=0$. Now, let $(B_{n})_{(n\geq 0)}$ be a sequence of pairwise disjoint elements 
of $\mathcal{B}$. We have

$$
m_X\left(\sum_{n\geq 0} B_n\right)=m\left(X^{-1}(\sum_{n\geq 0} B_n)\right) = m\left(\sum_{n\geq 0} X^{-1}(B_n)\right)
$$

\bigskip \noindent which, by the $\sigma$-additivity of $m$ applied to  $X^{-1}(B_n)$ that belong to $\mathcal{A}$, leads to

$$
m_X\left(\sum_{n\geq 0} B_n\right)= \sum_{n\geq 0} m\left(X^{-1}(B_n)\right) =\sum_{n\geq 0} m_X(B_n).
$$

\bigskip \noindent Part (b). \noindent It is clear that $m_A$ well-defined and non-negative. It is proper since  $m_A(\emptyset)=m(A \cap \emptyset)=0$. Remark that $\mathcal{A}_A$ is only the collection of
measurable sets (in $\mathcal{A}$) which are parts of $A$. So the $\sigma$-additivity of $m_A$ directly derives from that of $m$.

\newpage
\noindent \textbf{Part III : null sets, equivalence classes, completeness of measures}.\\

\noindent \textbf{Exercise 7}. \label{exercise07_sol_doc04-06} Let $(\Omega, \mathcal{A},m)$ a measure space.\\

\noindent (1) Show that  : (a) A measurable subset of a $m$-null set is a $m$-null set, (b) a countable union of $m$-null sets is a $m$-null set. Hint : use the $\sigma$-sub-additivity of $m$.\\

\bigskip 
\noindent (2) Let $f : (\Omega, \mathcal{A},m) \mapsto \overline{\mathbb{R}}$ a real-valued measurable application.\\

\noindent Show that
$$
(f \text{ finite })=(|f|<+\infty)=\bigcup_{k\geq 1} (|f|\leq k) \text{   (A)}
$$

\bigskip \noindent and then, by taking complements, we get

$$
(f \text{ infinite })=(|f|=+\infty)=\bigcap_{k\geq 1} (|f| > k) \text{   (B)}.
$$

\bigskip \noindent Show also that

$$
(f=0)=\bigcap_{k\geq 1} (|f| <1/k) \text{   (C)}.
$$

\bigskip \noindent \textit{Remark} : you may use strict or large inequalities in each formula.\\

\noindent Then the definitions \\

\noindent [ $f$ finite a.e. if and only if $m(f \text{ finite })=0$ ]\\

\noindent and\\

\noindent [ $f$ positive a.e. if and only if $m(f=0)=0$ ]\\

\noindent make sense.\\

\bigskip \noindent (3) Let $f_n : (\Omega, \mathcal{A},m) \mapsto \overline(\mathbb{R})$ a sequence real-valued measurable applications.\\

\noindent Show that if each $f_n$ is a.e. finite, they all the $f_n$'s are simultaneously finite outside a null-set.\\

\noindent Show that if each $f_n$ is a.e. positive, they all the $f_n$'s are simultaneously positive outside a null-set.\\ 

\bigskip \noindent (4) Let $f,g,h : (\Omega, \mathcal{A},m) \mapsto \mathbb{R}$ three finite real-valued measurable applications.\\

\noindent (a) Based on the triangle inequality $|f-h|\leq |f-g|+|g-h|$, remark that

$$
(|f-g|=0) \cap (|g-h|=0) \subset (|f-h|=0).
$$

\bigskip \noindent (b) and by taking the complements, that

$$
(|f-h|>0) \subset (|f-g|>0) \cup (|g-h|>0).
$$

\bigskip \noindent (c) Deduce from this that : $f=g$ a.e. and $g=h$ a.e., implies that $f=h$ a.e.\\

\noindent (d) Let $\mathcal{L}_0(\Omega, \mathcal{A}, m)$ the space of all real-valued measurable and a.e. finite applications $f_n : (\Omega, \mathcal{A},m) \mapsto \overline(\mathbb{R})$. Define the binary relation $\mathcal{R}$ on $\mathcal{L}_0(\Omega, \mathcal{A}, m)^2$ by :

$$
\forall (f,g) \in \mathcal{L}_0(\Omega, \mathcal{A}, m)^2, \ \ f \mathcal{R}g \ f=g \text{ a.e.}.
$$

\bigskip \noindent Show that $\mathcal{R}$ is an equivalence relation and the equivalence class of each $f \in \mathcal{L}_0(\Omega, \mathcal{A}, m)$ is

$$
\overset{\circ}{f}=\{g \in \mathcal{L}_0(\Omega, \mathcal{A}, m), \  f=g \text{ a.e.}\}.
$$

\bigskip \noindent The quotient set is denoted

$$
L_0(\Omega, \mathcal{A}, m)= \{\overset{\circ}{f}, f \in \mathcal{L}_0(\Omega, \mathcal{A}, m) \}.
$$

\bigskip \noindent \textbf{SOLUTIONS}.\\

\noindent Question (1). (a) Let $N$ be a $m$-null set and $A$ a measurable subset of $N$. We have $0\leq m(A) \leq m(N)=0$. Hence $m(A)=0$ and $A$ is a $m$-null set.\\

\noindent (b) Let $(N_n)_{(n\geq 0)}$ a countable family of null sets with respect to the measure $m$. Thus $N=\cup_{n\geq 0} N_n$ is measurable. By $\sigma$-sub-additivity of $m$, we have

$$
0\leq m(N) \leq \sum_{n\geq 0} m(N_n)=0,
$$

\bigskip \noindent which implies that $m(N)=0$ and thus, $\cup_{n\geq 0} N_n$ is a $m$-null set.$\blacksquare$\\

\noindent Question (2). We remind the following properties of real numbers :

$$
x=0 \ \Leftrightarrow \ (\forall \ k\geq 1, \ |x|\leq 1/k) \ \ \ \ (EQ1),
$$

$$
(x \text{ is finite }) \ \Leftrightarrow \ (\exists \ k\geq 1, \ |x|\leq k) \ \ \ \ (EQ2)
$$

\bigskip \noindent and

$$
(x \text{ is infinite }) \   \Leftrightarrow \ (\forall \ k\geq 1, \ |x|>k). \ \ \ \ (EQ3)
$$

\bigskip \noindent With this recall, we see that the formulas :

$$
(f \text{ finite })=(|f|<+\infty)=\bigcup_{k\geq 1} (|f|\leq k) \text{   (A)},
$$

$$
(f \text{ infinite })=(|f|=+\infty)=\bigcap_{k\geq 1} (|f| > k) \text{   (B)}
$$

\bigskip \noindent and

$$
(f=0)=\bigcap_{k\geq 1} (|f| <1/k) \text{   (C)}
$$

\bigskip \noindent are simple sets versions of the equivalences $(EQ1)$, $(EQ2)$ and $(EQ3)$. By the way, they show that if $f$ is measurable, then the sets

$$
(f \text{ finite }), \ \ \ (f \text{ infinite }) \ \ \ and  \ \ \ (f=0)
$$

\bigskip \noindent are measurable. This justifies the definitions introduced.\\

\noindent Question (3). (a) Let us denote $N_n=(f_n \ infinite)$ for each $n\geq 0$. If all the $f_n$'s are finite \textit{a.e.}, then we have $m(N_n)=0$. Put

$$
N=\cup_{n\geq 0} N_n.
$$

\bigskip \noindent Then $N$ is a null-set and each of the $f_n$'s is finite \textit{a.e.} on the complement of $N$.\\

\noindent (b) Let us denote $N_n=(f_n \leq 0)$ for each $n\geq 0$. If all the $f_n$'s are positive \textit{a.e.}, then we have $m(N_n)=0$. Put

$$
N=\cup_{n\geq 0} N_n.
$$

\bigskip \noindent Then $N$ is a null-set and each of the $f_n$'s is positive \textit{a.e.} on the complement of $N$.\\

\noindent Question (4).\\  

\noindent Question (4-a). This part is obvious.\\

\noindent Question (4-b). This is obvious too.\\

\noindent Question (4-c). Since $f$ and $g$ are finite measurable functions, we have $(f=g)=(|f-g|=0)$. Then $f=g$ \textit{a.e} and $g=h$ \textit{a.e} mean $m(|f-g|>0)=0$ and 
$m(|g-h|>0)=0$. So, if $f=g$ \textit{a.e} and $g=h$ \textit{a.e}, we get by (b) that $m(|f-h|>0)=0$, that is $f=h$ \textit{a.e.}

\noindent Question (4-d). The binary relation $\mathcal{R}$ is obviously reflexive and symmetrical. Its transitivity comes from part (c) just above. Thus it is an equivalence relation.

\bigskip \bigskip \noindent \textbf{Exercise 8}. \label{exercise08_sol_doc04-06} A measure space $(\Omega ,\mathcal{A},m)$ is said to be complete if and only if all subsets on $m$-null sets are measurable, and the then null $m$-sets.\\

\noindent Such a notion of complete measure is very pratical. Indeed, not being able to measure a set one knows it is part of $m$-null set is pathological. This opposes to naive belief that any part of a negligeable (null set) object should be negligeable too. The following exercise aims at fixing this problem.\\ 

\noindent Let $(\Omega ,\mathcal{A},m)$ be a measure. Set 
\begin{equation*}
\mathcal{N}=\{N\subset \Omega ,\exists \text{ }B\text{ }m\text{-null set} ,N\subset B\}).
\end{equation*}

\noindent Remark that $\mathcal{N}$ includes the $m$-null sets, $\emptyset$ in particular.\\
 
\bigskip \noindent Define

\begin{equation*}
\widehat{\mathcal{A}}=\{A\cup N,A\in \mathcal{A},N\in \mathcal{N}\})
\end{equation*}

\bigskip \noindent and next define the application 
$$
\widehat{m}:\widehat{\mathcal{A}}\rightarrow \mathbb{R}_{+}
$$ 

\bigskip \noindent such that 
$$
\widehat{m}(A\cup N)=m(A).
$$ 

\bigskip \noindent (a) Show that $\widehat{\mathcal{A}}$ is a $\sigma$-algebra including $\mathcal{A}$.\\

\noindent Hint. Showing that $B \in \widehat{\mathcal{A}} \Leftrightarrow B^c \in \widehat{\mathcal{A}}$ is not easy. Proceed as follows. 

\noindent Consider $B \in \widehat{\mathcal{A}}$ with $B=A\cup N$, $A\in \mathcal{A}$, $N\in \mathcal{N}$. Since $N\in \mathcal{N}$, consider $D \mathcal{A}$, $m(D)=0$ and $N \subset D$.\\

\noindent Justify $N^c=D^c+D\cap N^c$.\\

\noindent Plug this into $B^c=A^c \cap N^c$ to get $B^c=(A^c\cap D^c)+N_0$. Identify $N_0$ and conclude that $B^c \in  \widehat{\mathcal{A}}$.\\

\noindent (b) Show that the application $\widehat{m}$ is well-defined by showing that if 
$$
B \in \widehat{\mathcal{A}}
$$

\bigskip \noindent with two expressions $B=A_1\cup N_1=A_2\cup N_2$, $A_i\in \mathcal{A}$, $N_i\in \mathcal{N}$, $i=1,2$, then
$$
m(A_1)=m(A_2)
$$

\bigskip \noindent \textit{Hint} : Write $B=(A_1\cap A_2) \cup N_3$, identify $N_3$ and show that $N_3 \in \mathcal{N}$. Next show, by using $A_1=(A_1\cap A_2) \cup (A_1 \cap N_3)$ that $A_1=(A_1\cap A_2) + N_4$ and $N_4 \in \mathcal{N}$. Deduce that $m(A_1)=m(A_1\cap A_2)$. Do the same for $A_2$. Then show that $m(A_1)=m(A_2)$. Conclude.\\

\noindent (c) Show now that $\widehat{m}$ is a measure on $(\Omega, \widehat{\mathcal{A}})$, which is an extension of the measure $m$ on 
$(\Omega, \widehat{\mathcal{A}})$.\\

\noindent (d) Show  that the measure space  $(\Omega, \widehat{\mathcal{A}}, \widehat{m})$ is complete.\\

\noindent \textit{Hint}. Consider a $\widehat{m}$-null set $B$ of the form $B=A \cap N$ such that $A\in \mathcal{A}$ and $N\in \mathcal{N}$, and $N \subset D$, $D$ being a $m$-null set. Show that $A$ is a 
$m$-null set and deduce from this, that $C\in \mathcal{N}$ and conclude.\\

\noindent \textbf{Conclusion}. We \textbf{may and do} consider any measure space $(\Omega, \mathcal{A},m)$ as complete, eventually by considering it as a subspace of its completition 
$(\Omega, \widehat{\mathcal{A}}, \widehat{m})$.\\

\bigskip \noindent \textbf{SOLUTIONS}.\\

\noindent Question (a). Let us begin to see that by taking $N=\emptyset$ in the definition of $\widehat{\mathcal{A}}$, we get all the elements of $\mathcal{A}$. Then

$$
\mathcal{A} \subset \widehat{\mathcal{A}}.
$$

\bigskip \noindent Let us check that $\widehat{\mathcal{A}}$ is a $\sigma$-algebra. By taking $A=\Omega$ in the definition of $\widehat{\mathcal{A}}$, we get $\Omega \cup N=\Omega$, whatever be $N \in \mathcal{N}$. Also, $\widehat{\mathcal{A}}$ is stable under countable union for the following reason. If $(A_n \cup N_n)_{(n\geq 0)}$ is a sequence of elements of $\widehat{\mathcal{A}}$ where for each $n\geq 0$,
$A_n \in \mathcal{A}$ and $N_n \in \mathcal{N}$, we have, by associativity of the union,

$$
\bigcup_{n\geq 0} \biggr(A_n\cup N_n\biggr) = \left( \bigcup_{n\geq 0} A_n \right) \bigcup \left(\bigcup_{n\geq 0} N_n\right),
$$

\bigskip \noindent which is of the form $A \cup N$, with $A=\bigcup_{n\neq 0} A_n \in \mathcal{A}$ and $N=\bigcup_{n\neq 0} N_n \in \mathcal{N}$.\\

\noindent To finish, let us prove that $\widehat{\mathcal{A}}$ is stable under complements. Consider $B \in \widehat{\mathcal{A}}$ with $B=A\cup N$, $A\in \mathcal{A}$, $N\in \mathcal{N}$.\\

\noindent Since $N\in \mathcal{N}$, there exists, by definition of $\mathcal{N}$,  $D \mathcal{A}$ such that $m(D)=0$ and $N \subset D$. Since $D^c \subset N^c$, we have 
$N^c = D^c + (N^c \setminus D^c)$, that is 

$$
N^c=D^c+D\cap N^c.
$$

\bigskip \noindent Upon plugging  this into $B^c=A^c \cap N^c$, we get

$$
B^c= A^c \cap N^c= A^c \cap (D^c+D\cap N^c)=(A^c \cap D) + ( D \cap A^c \cap N^c).
$$ 

\bigskip \noindent We have $A_0=A^c \cap D \in \mathcal{A}$ and $N_0=D \cap A^c \cap N^c \subset D$ so that $N_0 \in \mathcal{N}$, and $B^c=A_0 \cup N_0$. In conclusion, $B^c \in \widehat{\mathcal{A}}$.\\

\noindent We conclude that $\widehat{\mathcal{A}}$ is a $\sigma$-algebra containing $\mathcal{A}$.\\

\noindent Question (b). Let us show that is well-defined, that is the value of $\widehat{m}(B)$ does not depend on one particular expression of $B$ as an element of $\widehat{\mathcal{A}}$. Suppose that 
$B \in \widehat{\mathcal{A}}$ with two writings : $B=A_1\cup N_1=A_2\cup N_2$, $A_i\in \mathcal{A}$, $N_i\in \mathcal{N}$, $i=1,2$. For $\widehat{m}(B)$ to be well-defined, we have to prove that

$$
m(A_1)=m(A_2)
$$

\bigskip \noindent By taking the intersection of $B$ by itself, we get

$$
B=x(A_1 \cap N_2)
$$

\bigskip \noindent where

$$
N_3=(A_1 \cap N_2) \cup (N_1 \cap A_2) \cup (N_1 \cap N_2).
$$

\bigskip \noindent Hence, for $i=1$ or $i=2$, 

$$
A_i=A_i \cap B=(A_1 \cap A_2) \cup (N_3 \cap A_i)
$$

\bigskip \noindent which gives

$$
A_i=(A_1 \cap A_2) + ((N_3 \cap A_i) \cap (A_1 \cap A_2)^c).
$$

\noindent Since each $N_i$ - $i=1,2$ - is in $\mathcal{N}$, there exists $D_i$ $m$-null set such that $N_i \subset D_i$, $i=1,2$. It is clear that $N_3$ and next 
$N_4=((N_3 \cap A_i) \cap (A_1 \cap A_2)^c)$ are in $D=D_1 \cup D_2$, which is a $m$-null set. Hence, for each $i=1,2$,

$$
A_i \setminus (A_1 \cap A_2) \subset D,
$$

\bigskip \noindent and, by Point (1) of Exercise 2,

$$
0\leq m(A_i \setminus (A_1 \cap A_2))=m(A_i) - m(A_1 \cap A_2)\leq m(D)=0.
$$

\bigskip \noindent We get

$$
m(A_1)=m(A_2)=m(A_1 \cap A_2),
$$

\noindent which was the target. In conclusion, $\widehat{m}$ is defines an application.$\square$\\

\noindent Question (c). Let us show that $\widehat{m}$ is a measure. First, it is non-negative. It is also proper. To see this, remark that the empty set has only on expression as an element of
 $\widehat{\mathcal{A}}$, that is $\emptyset=\emptyset \cup \emptyset$, with $A=\emptyset$ and $N=\emptyset$. Thus $\widehat{m}(\emptyset)=0$. Let $B_n=A_n \cup N_n$, $n\geq 0$, a sequence of pairwise disjoint elements of $\widehat{\mathcal{A}}$ where for each $n\geq 0$, $A_n \in \mathcal{A}$ and $N_n \in \mathcal{N}$. We have

$$
\bigcup_{n\geq 0} (A_n\cup N_n) = \left( \bigcup_{n\geq 0} A_n \right) \bigcup \left(\bigcup_{n\geq 0} N_n\right),
$$

\bigskip \noindent which is

$$
\sum_{n\geq 0} B_n = \left( \sum_{n\geq 0} A_n \right) \bigcup N,
$$

\bigskip \noindent where

$$
N=\left(\bigcup_{n\neq 0} N_n\right) \in \mathcal{N}.
$$

\bigskip \noindent We get

\begin{eqnarray*}
\widehat{m}(\sum_{n\neq 0} B_n) &=& m(\sum_{n\neq 0} A_n)\\
&=& \sum_{n\neq 0} m(A_n)\\
&=& \widehat{m}(A_n \cup N_n)\\
&=& \widehat{m}(B_n),
\end{eqnarray*}

\bigskip \noindent which proves that $\widehat{m}$ is a measure. Now, let $A\in \mathcal{A}$. We have

$$
m(A)=\widehat{m}(A \cup \emptyset)=\widehat{m}(A).
$$ 

\bigskip \noindent Since $\widehat{m}$ and $m$ coincide on $\mathcal{A}$, we see that $\widehat{m}$ is an extension of $m$ on $\widehat{\mathcal{A}}$.\\

\noindent Question (d). Let $B$ be a $\widehat{m}$-null set and  Let $C$ be a part of $B$. By definition, $B$ is of the form $B=A \cap N$ such that $A\in \mathcal{A}$ and $N$ is part of a $m$-null set $D$. Since we have $\widehat{m}(B)=0$, it follows that $m(A)=0$. Hence $A$ is a $m$-null set and, since $C \subset \cup N$, we have that $C\subset A \cap D$. Since $A \cap D$ is a $m$-null set, we conclude that $C$ is in $\mathcal{N}$, and hence belongs to $\widehat{\mathcal{A}}$. Thus, the measure $\widehat{m}$ is complete.

\newpage
\bigskip \noindent \textbf{Part IV : Determining classes for $\sigma$-finite measures}.\\

\bigskip \noindent \textbf{Exercise 9}. \label{exercise09_sol_doc04-06}  Let $m_{1}$ and $m_{2}$ be two $\sigma$-finite measures on an algebra $\mathcal{C}$ of subsets of $\Omega$. Find a countable partition of $\Omega$ into elements of $\mathcal{C}$ in the form

$$
\Omega=\sum_{j\geq 0} \Omega_j, (\Omega_j \in \mathcal{C}, j\geq 0),
$$

\bigskip \noindent such that $m_{1}$ and $m_{2}$ are both finite on each $\Omega_j$, $j\geq 0$.\\

\bigskip \noindent \textbf{SOLUTIONS}. Each $m_{i}$ ($i=1,2$) is $\sigma$-finite on $\mathcal{C}$. Then there exists, for each $i \in \{1,2\}$, a partition of $\Omega$

$$
\{\Omega^{(i)}_j, j\geq 0\} \subset \mathcal{C} 
$$

\bigskip \noindent such that, for $i \in \{1,2\}$, for all $j\geq 0$, we have  $m_j(\Omega^{(i)}_j)<+\infty$. But we have

$$
\Omega = \sum_{h\geq 0} \Omega^{(1)}_h \ \ and \ \Omega = \sum_{k\geq 0} \Omega^{(2)}_k,
$$ 

\bigskip \noindent and by taking the intersection of $\Omega$ by itself, we get

$$
\Omega = \sum_{h\geq 0} \sum_{k\geq 0} \Omega^{(1)}_h \cap \Omega^{(1)}_k.
$$

\bigskip \noindent Hence, the family

$$
\{\Omega^{(1)}_h \cap \Omega^{(2)}_k, \ h\geq 0, \ k\geq 0\} 
$$

\bigskip \noindent is a subclass of $\mathcal{C}$ since the latter is stable under finite intersection. It is also a countable partition of $\Omega$ on $\mathcal{C}$. By the non-decreasingness of the $m_i$'s, we also have for all $h\geq 0$ and $k \geq 0$,

$$
m_1(\Omega^{(1)}_h \cap \Omega^{(2)}_k) \leq m_1(\Omega^{(1)}_h) <+\infty
$$

\bigskip \noindent and

$$
m_2(\Omega^{(1)}_h \cap \Omega^{(2)}_k) \leq m_2(\Omega^{(k)}_k) <+\infty.
$$

\bigskip \noindent By denoting the countable set $\{\Omega^{(1)}_h \cap \Omega^{(2)}_k, \ h\geq 0, \ k\geq 0\}$ as $\{\Omega_j, \ j\geq 0\}$, we conclude the solution.\\

\bigskip \noindent \textbf{Exercise 10}. \label{exercise10_sol_doc04-06}  Let $(\Omega, \mathcal{A})$ a measurable space such that $\mathcal{A}$ is generated by an algebra $\mathcal{C}$ of subsets of $\Omega$. Show that two finite measure which are equal on $\mathcal{C}$, are equal (on $\mathcal{A}$).\\

\noindent \textit{Hint} : Let $m_1$ and $m_2$ be two finite measures on $(\Omega, \mathcal{A})$ equal on $\mathcal{C}$. Set 
$$
\mathcal{M}=\{A \in \mathcal{A}, m_1(A)=m_2(A) \}.
$$

\bigskip \noindent Show that is a monotone class including $\mathcal{C}$. Apply \textit{Exercise 3 in DOC 01-04} to conclude.\\

\bigskip \noindent \textbf{SOLUTIONS}. By applying first \textit{Point (4) of Exercise 2} of this series, we easily see that $\mathcal{M}$ is monotone. In applying the continuity of the $m_i$'s for non-increasing sequences, we do not have to worry about the finiteness of $m_i(A_n)$ for $n$ large enough, because the $m_i$'s are finite here.\\

\noindent By assumption, $\mathcal{M}$ contains $\mathcal{C}$. So, it contains the monotone class generated by $\mathcal{C}$, which, by \textit{Exercise 3 in DOC 01-04}, is equal to  $\sigma(\mathcal{C})=\mathcal{A}$. By construction, $\mathcal{M}$ is in $\mathcal{A}$. As a conclusion, we have

$$
\{A \in \mathcal{A}, m_1(A)=m_2(A) \}=\mathcal{A}.
$$ 

\bigskip \noindent This means that $m_1=m_2$ on $\mathcal{A}$.\\

\bigskip \noindent \textbf{Exercise 11}. \label{exercise11_sol_doc04-06}  (The $\lambda$-$\pi$ Lemma) Let $(\Omega, \mathcal{A})$ a measurable space such that $\mathcal{A}$ is generated by an $\pi$-system $\mathcal{P}$ containing $\Omega$. Show that two finite measures which equal on $\mathcal{P}$, are equal (on $\mathcal{A}$).\\

\noindent \textit{Hint} : Let $m_1$ and $m_2$ be two finite measures on $(\Omega, \mathcal{A})$ equal on $\mathcal{P}$. Set 
$$
\mathcal{D}=\{A \in \mathcal{A}, m_1(A)=m_2(A) \}.
$$

\bigskip \noindent Show that is a $\lambda$-system (Dynkin system), which includes $\mathcal{P}$. Apply \textit{Exercise 4 in DOC 01-04} to conclude.\\

\bigskip \noindent \textbf{SOLUTIONS}. Let us show that $\mathcal{D}$ is a Dynkin system.\\

\noindent (i) Let us show that $\Omega \in \mathcal{P}$. Since $\Omega$ is $\mathcal{P}$ and since $m_1=m_2$ on $\mathcal{P}$, we have $m_1(\Omega)=m_2(\Omega)$ and hence : $\Omega \in \mathcal{P}$.

\noindent (ii) Let us show that if $A \in \mathcal{D}$, $B \in \mathcal{D}$ and $A \subset B$, we have $B \setminus A \in \mathcal{D}$. Indeed, suppose that $A \in \mathcal{D}$, $B \in \mathcal{D}$ and $A \subset B$. By \textit{Point (1) of Exercise 2} of this series, we have for each $o \in \{1,2\}$,

$$
m_i(B \setminus A)=m_i(B)-m_i(A)
$$

\bigskip \noindent and since $A$ and $B$ are in $\mathcal{D}$, we have

$$
m_1(B \setminus A)=m_1(B)-m_1(A) = m_2(B)-m_2(A)=m_2(B \setminus A).
$$

\bigskip \noindent Thus $m_1(B \setminus A)=m_2(B \setminus A)$ and, by this, $B \setminus A \in \mathcal{D}$.\\

\noindent (iii) Let us show that $\mathcal{D}$ is stable under countable disjoint unions. Let $(A_n)_{(n\geq n)}$ a sequence of pairwise disjoint elements of $\mathcal{D}$ and let

$$
A =\sum_{n\geq 0} A_n.
$$

\bigskip \noindent By $\sigma$-additivity of the $m_i$'s and by using that each $A_n \mathcal{D}$, we have

$$
m_1(A)=\sum_{n\geq 0} m_1(A_n) =  \sum_{n\geq 0} m_2(A_n)=m_2(A).
$$

\bigskip \noindent Thus, $A \in \mathcal{D}$.\\

\noindent In total, $\mathcal{D}$ is a Dynkin system including the $\pi$-system $\mathcal{P}$. So it contains the Dynkin system generated by $\mathcal{P}$, which, by \textit{Exercise 4 in DOC 01-04}, is equal to  $\sigma(\mathcal{P})=\mathcal{A}$. By construction, $\mathcal{D}$ is in $\mathcal{P}$. As a conclusion, we have

$$
\{A \in \mathcal{A}, m_1(A)=m_2(A) \}=\mathcal{A}.
$$ 

\bigskip \noindent This means that $m_1=m2$ on $\mathcal{A}$.\\

\bigskip \noindent \textbf{Exercise 12}. \label{exercise12_sol_doc04-06}  (Extension of Exercise 10). Let $(\Omega, \mathcal{A})$ a measurable space such that $\mathcal{A}$ is generated by an algebra $\mathcal{C}$ of subsets of $\Omega$. Show that two measures which are equal on $\mathcal{C}$ and $\sigma$-finite on $\mathcal{C}$, are equal (on $\mathcal{A}$).

\noindent \textit{Hint} : Let $m_1$ and $m_2$ be two $\sigma$-finite on $\mathcal{C}$. Use a partition like to the one obtained in Exercise 8. Decompose each measure $m_i$, $i=1,2$ into

$$
m_i(A)=\sum_{j\geq 0} m_i(A\cap \Omega_{j}), \ A \in \mathcal{A}.
$$

\bigskip \noindent Next, apply Exercise 10 to conclude.\\

\bigskip \noindent \textbf{SOLUTIONS}. Since $m_1$ and $m_2$ are two measures which are $\sigma$-finite on $\mathcal{C}$, we use \textit{Exercise 8} to get a partition on $\mathcal{C}$ :

$$
\Omega=\sum_{j\geq 0} \Omega_j
$$

\bigskip \noindent such that $m_{1}$ and $m_{2}$ are both finite on each $\Omega_j$, $j\geq 0$. So, for each $i \in \{1, 2\}$, for each $A \in \mathcal{A}$, we have

$$
A=\sum_{j\geq 0} \Omega_j \cup A,
$$

\bigskip \noindent and next

$$
m_1(A)=\sum_{j\geq 0} m_1(\Omega_j \cup A) \ \ and m_2(A)=\sum_{j\geq 0} m_2(\Omega_j \cup A).
$$

\bigskip \noindent For each $j\geq 0$, the application

$$
m_{1,j}(A)=m_1(\Omega_j \cup A), \ \ A \in \mathcal{A}. 
$$

\bigskip \noindent defines the induced measure of $m_1$ on $\Omega_j$. Each $m_{1,j}$ is also a measure on $\mathcal{A}$, which is finite and we have that

$$
m_1= \sum_{j\geq 0} m_{1,j}. \ \ \ (D1)
$$

\bigskip \noindent We do the same for $m_2$ to get

$$
m_2= \sum_{j\geq 0} m_{2,j}, \ \ \ (D2)
$$

\bigskip \noindent where each $m_{2,j}$ is a finite measure on $\mathcal{A}$. Now, for any $A\in \mathcal{C}$, we have that for any $j\geq 0$, $\Omega_j \cup A \in \mathcal{C}$ and by the equality of 
$m_1$ and $m_2$ on $\mathcal{C}$ :

$$
m_{1,j}(A)=m_1(\Omega_j \cup A)=m_2(\Omega_j \cup A)=m_{2,j}(A).
$$

\bigskip \noindent We conclude as follows. For each $j\geq 0$, the measures  $m_{1,j}$ and $m_{2,j}$ are finite and equal on $\mathcal{C}$. So, by \textit{Exercise 9}, they equal. Finally, Formulas (D1) and (D2) together imply that $m_1=m_2$.

\newpage

\noindent \textbf{PART V : Extensions}.\\

\bigskip \noindent\textbf{Exercise 13}. \label{exercise13_sol_doc04-06}  Let $m$ be a proper and non-negative application defined from an algebra $\mathcal{C}$\ of subsets of $\Omega$ to $\overline{\mathcal{R}}_{+}$. Show that if $m$ is additive and $\sigma$-sub-additive with respect to $\mathcal{C}$, then $m$ is $\sigma$-additive on $\mathcal{C}$.\\

\bigskip \noindent \textbf{SOLUTIONS}. Assume that $m$ is a proper and non-negative application defined from an algebra $\mathcal{C}$\ of subsets of $\Omega$ to $\overline{\mathcal{R}}_{+}$ such that
$m$ is additive and $\sigma$-sub-additive with respect to $\mathcal{C}$. Let us show it is $\sigma$-additive.\\

\noindent Let $(A_{n})_{n\geq 1}$\ is a sequence of pairwise disjoint elements of $\mathcal{C}$ and set $A=\sum_{n\geq 0} A_n$. We have, for any $k\geq 1$,

$$
A=\sum_{n=0}^{k} A_n + \sum_{n=k+1}^{k} A_n. 
$$

\bigskip \noindent Let us denote, for $k\geq 1$,

$$
B_k=\sum_{n=k+1}^{k} A_n.
$$

Since S $\mathcal{C}$ is an algebra, the set  $\sum_{n=0}^{k} A_n$ is in $\mathcal{C}$, so is its complement, which is $B_k$. By applying the $\sigma$-sub-additivity of $m$ on $\mathcal{C}$, we have

$$
m(B_k) \leq \sum_{n=k+1}^{k} m(A_n).
$$ 

\bigskip \noindent By additity of$m$ on $\mathcal{C}$, we have

$$
m(A)=\sum_{n=0}^{k} m(A_n) + m(B_k). \ \ (D)
$$

\bigskip \noindent We see that the second term in the left member of the formula below is bounded by the tail of a series. If the series is convergent, we will be able to conclude. So, we may consider two cases : \\

\noindent \textit{Case 1}. The series $\sum_{n=0}^{\infty} m(A_n)$ is convergent. Then its tail $\sum_{n=k+1}^{\infty} m(A_n)$ converges to zero as $k\rightarrow +\infty$. We get

$$
0\leq \limsup_{k\rightarrow +\infty} m(B_k) \leq \limsup_{k\rightarrow +\infty} \sum_{n=k+1}^{k} m(A_n)=0,
$$

\bigskip \noindent so that

$$
\lim_{k\rightarrow +\infty} m(B_k)=0.
$$

\bigskip \noindent This leads to

\begin{eqnarray*}
m(A) &=&\lim_{k\rightarrow +\infty} \sum_{n=0}^{k} m(A_n) +  m(B_k)\\
&=& \lim_{k\rightarrow +\infty} \sum_{n=0}^{k} m(A_n) + \lim_{k\rightarrow +\infty} m(B_k)\\
&=& \lim_{k\rightarrow +\infty} \sum_{n=0}^{k} m(A_n).
\end{eqnarray*} 

\bigskip \noindent \textit{Case 2}. The series $\sum_{n=0}^{\infty} m(A_n)$ is divergent, $\sum_{n=0}^{\infty} m(A_n)=+\infty$. From Equation (D) above, we have for all $k\geq 1$,

$$
m(A) \geq \sum_{n=0}^{k} m(A_n).
$$

\bigskip \noindent By letting $k\rightarrow +\infty$, we get $m(A)=+\infty$ and hence

$$
m(A)=\sum_{n=0}^{\infty} m(A_n)=+\infty.
$$

\bigskip \noindent In both cases, we obtained
$$
m(A) = \sum_{n=0}^{\infty} m(A_n).
$$

\noindent So, $m$ is $\sigma$-additive.\\

\bigskip \noindent\textbf{Exercise 14}. \label{exercise14_sol_doc04-06}  Let $m : \mathcal{C} \rightarrow \overline{\mathbb{R}}_{+}$, be a proper non-negative and finite application. Suppose that $m$ is finitely-additive on $\mathcal{C}$ and continuous at $\emptyset$, i.e., if $(A_{n})_{n\geq 1}$\ is a sequence of elements of $\mathcal{C}$ non-increasing to $\emptyset$, then $m(A_{n})\downarrow 0$, as $n\uparrow +\infty$. Show that $m$ is $\sigma$-additive on $\mathcal{C}$.\\

\bigskip \noindent \textbf{SOLUTIONS}. Assume that $m$ is a proper and non-negative application defined from an algebra $\mathcal{C}$\ of subsets of $\Omega$ to $\overline{\mathcal{R}}_{+}$ such that
$m$ is additive and continuous at $\emptyset$. Let us show it is $\sigma$-additive.\\

\noindent Let $(A_{n})_{(n\geq 0)}$\ is a sequence of pairwise disjoint elements of $\mathcal{C}$ and set $A=\sum_{n\geq 0} A_n$. We have, for any $k\geq 1$,

$$
A=\sum_{n=0}^{k} A_n + \sum_{n=k+1}^{k} A_n. 
$$

\bigskip \noindent Let us denote, for $k\geq 1$,

$$
B_k=\sum_{n=k+1}^{k} A_n.
$$

\bigskip \noindent Since $\mathcal{C}$ is an algebra, the set  $\sum_{n=0}^{k} A_n$ is in $\mathcal{C}$, so is its complement, which is $B_k$. By applying the $\sigma$-sub-additivity of $m$ on $\mathcal{C}$, we have

$$
m(B_k) \leq \sum_{n=k+1}^{k} m(A_n).
$$ 

\bigskip \noindent By additity of $m$ on $\mathcal{C}$, we have

$$
m(A)=\sum_{n=0}^{k} m(A_n) + m(B_k). \ \ (D)
$$

\bigskip \noindent It is clear that $B_k$ decreases to $\emptyset$. To see that, remark that if we increment $k$, we loose the subset $A_{k+1}$. So the sequence $(B_k)$ is non-increasing. Remark also that all the $B_k$ are in $A$. Each element of $A$ is in one and only one $B_k$, say in $A_{k_0}$. So this element is not in $B_{k_0}$, and by this not the intersection of the $B_k$. So this intersection is empty since all the elements of $A$ will drop out.\\

\bigskip \noindent By Formula (D) and by the finiteness of $m$, we have that $m(B_k)=m(A)-\sum_{n=0}^{k} m(A_n)$ is finite for all $k$. The continuity of $m$ at $\emptyset$ implies that $m(B_k) \downarrow 0$ as 
$k \uparrow +\infty$. Thus, we have

\begin{eqnarray*}
m(A) &=&\lim_{k\rightarrow +\infty} \sum_{n=0}^{k} m(A_n) +  m(B_k)\\
&=& \lim_{k\rightarrow +\infty} \sum_{n=0}^{k} m(A_n) + \lim_{k\rightarrow +\infty} m(B_k)\\
&=& \lim_{k\rightarrow +\infty} \sum_{n=0}^{k} m(A_n).
\end{eqnarray*} 

\bigskip \noindent This concludes the solution.\\

\bigskip \noindent\textbf{Exercise 15}. \label{exercise15_sol_doc04-06}  Let $\mathcal{S}$ be semi-algebra of subsets of $\Omega$. Let $m : \mathcal{S} \rightarrow \overline{\mathbb{R}}_{+}$ be a proper, non-negative and additive application. We already know that $\mathcal{C}=a(\mathcal{S})$ is the collection of all finite sums of elements of $\mathcal{S}$, that is :

$$
a(\mathcal{S})=\mathcal{C}=\{A_1+A_2+...+A_k, \ k\geq 1, \ A_1 \in \mathcal{S}, \ A_2 \in \mathcal{S}, ..., A_k \in \mathcal{S}\}.
$$

\bigskip \noindent Define the application $\widehat{m} : \mathcal{C} \rightarrow \overline{\mathbb{R}}_{+}$ as follows. For any $A \in \mathcal{C}$, with $A=A_1+A_2+...+A_k$, $k\geq 1$, $A_1 \in \mathcal{S}$, 
$ A_2 \in \mathcal{S}$, ..., $A_k \in \mathcal{S}$, set
 
$$
\widehat{m}(A)=\widehat{m}(A_1+A_2+...+A_k)=m(A_1)+m(A_2)+...+m(A_k), \ \ (DEF)
$$

\bigskip \noindent Question (1). Show that the application $\widehat{m}$ is a well-defined.\\

\noindent Question (2). Show that $\widehat{m}$ is an extension of $m$ as a non-negative, proper and additive application with the additive properties.\\ 

\noindent Question (3). Show that $m$ is $\sigma$-additive on $\mathcal{S}$ if and only if $\widehat{m}$ is $\sigma$-additive on $\mathcal{C}$.\\
 
\noindent Question (4). Show that $m$ is $\sigma$-sub-additive on $\mathcal{S}$ if and only if $\widehat{m}$ is $\sigma$-sub-additive on $\mathcal{C}$.\\

\noindent Question (4). Show that $m$ is $\sigma$-finite with respect to $\mathcal{S}$ if and only if  $\widehat{m}$ is $\sigma$-finite with respect to $\mathcal{C}$.\\

\noindent Question (1). The application $\widehat{m}$ is well-defined, that is the definition in Formula (DEF) is coherent, if the definition of $\widehat{m}(A)$ in Formula (DEF) does not depend on the
represent of $A \ \mathcal{C}$. To see this, let us consider two expressions of $B$.\\

$$
A=A_1+A_2+...+A_k, \ k\geq 1, A_1 \in \mathcal{S}, \ A_2 \in \mathcal{S}, \cdots, A_k \in \mathcal{S}
$$

\bigskip \noindent and

$$
A=B_1+B_2+...+B_h, \ k\geq 1, B_1 \in \mathcal{S}, \ B_2 \in \mathcal{S}, \cdots, B_k \in \mathcal{S}.
$$

\bigskip \noindent The coherence of the definition is proved if we establish that
$$
\sum_{i=1}^{k} m(A_i) = \sum_{i=j}^{h} m(B_j). \ \ (QC)
$$ 

\bigskip \noindent For any any $j \in \{1,...,h\}$, we have (by taking the intersection of $B_j$ with $A=A_1+A_2+...+A_k$) :

$$
B_j= \sum_{i=1}^{k} B_j \cup A_i.
$$

\bigskip \noindent In the same manner, we also have for any any $i \in \{1,...,k\}$, we have (by taking the intersection of $B_j$ with $A=A_1+A_2+...+A_k$) :

$$
A_i= \sum_{j=1}^{h} A_i \cup B_j.
$$

\bigskip \noindent By additivity of $m$, we have

$$
m(A_i)= \sum_{j=1}^{h} m(A_i \cup B_j) \ \ (SU1) \ \ and \ \ m(B_j)= \sum_{i=1}^{k} m(B_j) \cup A_i \ \ (SU2).
$$

\bigskip \noindent Let us combine all this to get

\begin{eqnarray*}
\sum_{i=1}^{k} m(A_i)&=&\sum_{i=1}^{k} m(A_i)\\
&=& \sum_{i=1}^{k} \sum_{j=1}^{h} m(A_i \cup B_j) \ (SU1)\\
&=& \sum_{j=1}^{h} \sum_{i=1}^{k} m(A_i \cup B_j) \ \ (Fubini)\\
&=& \sum_{j=1}^{h} m(B_j).  \ (SU2) 
\end{eqnarray*}

\bigskip \noindent We are done.\\

\noindent Question (2). First, it is clear that, by construction, $\widehat{m}$ and $m$ coincide on $\mathcal{S}$. Next, obviously $\widehat{m}$ is proper and non-negative. To show that is additive, consider two disjoint elements of $\mathcal{C}$ :
$$
A=A_1+A_2+...+A_k, \ k\geq 1, A_1 \in \mathcal{S}, \ A_2 \in \mathcal{S}, \cdots, A_k \in \mathcal{S}
$$

\bigskip \noindent and

$$
B=B_1+B_2+...+B_h, \ k\geq 1, B_1 \in \mathcal{S}, \ B_2 \in \mathcal{S}, \cdots, B_k \in \mathcal{S}.
$$

\bigskip \noindent Since the $A_i$'s and the $B_j$'s are disjoint, we have

$$
A+B=A_1+A_2+...+A_k + B_1+B_2+...+B_h,
$$

\bigskip \noindent which, by the definition of $\widehat{m}$, implies that

$$
\widehat{m}(A+B)=m(A_1)+m(A_2)+...+m(A_k) + m(B_1)+(B_2)+...+(B_h),
$$

\bigskip \noindent which is

$$
\widehat{m}(A+B)=\widehat{m}(A) + \widehat{m}(B),
$$

\bigskip \noindent which was the target.\\

\noindent Question (3). We only need to show that $\widehat{m}$ is $\sigma$-additive if $m$ is, the other implication is obvious. Suppose that $m$ is $\sigma$-additive. Let $A$ be an element of
$\mathcal{C}$ which is decomposed into elements $A_n$ in $\mathcal{C}$, 

$$
A=\sum_{n\geq 0} A_n, \ A_n \in \mathcal{C}, \ n\geq 0.
$$

\bigskip \noindent We have to prove that 

$$
m(A)=\sum_{n\geq 0} m(A_n).
$$

\bigskip \noindent By definition of $\mathcal{C}$, $A$ is of the form

$$
A=A_1+A_2+...+A_h, \ k\geq 1, A_1 \in \mathcal{S}, \ A_2 \in \mathcal{S}, \cdots, A_k \in \mathcal{S}
$$

\bigskip \noindent and each $A_n$ is of the form

$$
A_n=A_{n,1}+A_{n,2}+...+A_{n,h(n)}, \ h(n)\geq 1, A_{n,1} \in \mathcal{S}, \ A_{n,2} \in \mathcal{S}, \cdots, A_{n,h(n)} \in \mathcal{S}.
$$
 
\bigskip \noindent We have

$$
A=\sum_{n\geq 0} \sum_{j=1}^{h(n)} A_{n,j}
$$

\bigskip \noindent and for each $i \in \{1,...,n\}$, we get by taking the intersection of $A_i$ with $\sum{n\geq 0} \sum_{j=1}^{h(n)} A_{n,j}$ that

$$
A_i=\sum_{n\geq 0} \sum_{j=1}^{h(n)} A_{n,j} \cap A_i.
$$

\bigskip \noindent By the $\sigma$-additivity of $m$ on $\mathcal{S}$, we get

$$
m(A_i)=\sum_{n\geq 0} \sum_{j=1}^{h(n)} m(A_{n,j} \cap A_i). \ (SU1)
$$

\bigskip \noindent Also, for each $n\geq 0$, $\j \in \{1,...,h(n)\}$,

$$
A_{n,j}=A_{n,j} \cap A = \sum_{i=1}^{h} A_{n,j} \cap A_i,
$$

\bigskip \noindent which, by the additivity of $m$ on $\mathcal{S}$, leads to

$$
A_{n,j}=\sum_{j=1}^{h(n)} m(A_{n,j} \cap A_i). \ (SU2)
$$

\bigskip \noindent Now, by the additivity of $\widehat{m}$ on $\mathcal{C}$, we have

\begin{eqnarray*}
\widehat{m}(A)&=&\sum_{i=1}^{h} \widehat{m}(A_i)\\
&=&\sum_{i=1}^{h} m(A_i)\\
&=&\sum_{i=1}^{h} \sum_{n\geq 0} \sum_{j=1}^{h(n)} m(A_{n,j} \cap A_i) \ (SU1)\\
&=& \sum_{n\geq 0} \sum_{j=1}^{h(n)} \sum_{i=1}^{h} m(A_{n,j} \cap A_i) \ (Fubini)\\ 
&=& \sum_{n\geq 0} \sum_{j=1}^{h(n)} m(A_{n,j}) \ (SU2)\\ 
&=& \sum_{n\geq 0}  \widehat{m}(A_n),
\end{eqnarray*}

\bigskip \noindent which ends the solution.\\

\noindent Question (4). We only need to show that $\widehat{m}$ is $\sigma$-sub-additive if $m$ is, the other implication is obvious. Suppose that $m$ is $\sigma$-sub-additive. Let $A$ be an element of
$\mathcal{C}$ which is decomposed into elements $A_n$ in $\mathcal{C}$, 

$$
A=\sum_{n\geq 0} A_n, \ A_n \in \mathcal{C}, \ n\geq 0.
$$

\bigskip \noindent We have to prove that 

$$
m(A) \leq \sum_{n\geq 0} m(A_n).
$$

\bigskip \noindent By definition of $\mathcal{C}$, $A$ is of the form

$$
A=A_1+A_2+...+A_h, \ k\geq 1, A_1 \in \mathcal{S}, \ A_2 \in \mathcal{S}, \cdots, A_k \in \mathcal{S},
$$

\bigskip \noindent and each $A_n$ is of the form

$$
A_n=A_{n,1}+A_{n,2}+...+A_{n,h(n)}, \ h(n)\geq 1, A_{n,1} \in \mathcal{S}, \ A_{n,2} \in \mathcal{S}, \cdots, A_{n,h(n)} \in \mathcal{S}.
$$
 
\bigskip \noindent We have

$$
A=\sum_{n\geq 0} \sum_{j=1}^{h(n)} A_{n,j}
$$

\bigskip \noindent and for each $i \in \{1,...,n\}$, we get by taking the intersection of $A_i$ with $\sum{n\geq 0} \sum_{j=1}^{h(n)} A_{n,j}$ that

$$
A_i=\sum_{n\geq 0} \sum_{j=1}^{h(n)} A_{n,j} \cap A_i.
$$

\bigskip \noindent By the $\sigma$-sub-additivity of $m$ on $\mathcal{S}$, we get

$$
m(A_i) \leq \sum_{n\geq 0} \sum_{j=1}^{h(n)} m(A_{n,j} \cap A_i). \ (SU1)
$$

\noindent Also, for each $n\geq 0$, $\j \in \{1,...,h(n)\}$,

$$
A_{n,j}=A_{n,j} \cap A = \sum_{i=1}^{h} A_{n,j} \cap A_i,
$$

\bigskip \noindent which, by the additivity of $m$ on $\mathcal{S}$, leads to

$$
A_{n,j}=\sum_{j=1}^{h(n)} m(A_{n,j} \cap A_i). \ (SU2)
$$

\bigskip \noindent Now, by the additivity of $\widehat{m}$ on $\mathcal{C}$, we have

\begin{eqnarray*}
\widehat{m}(A)&=&\sum_{i=1}^{h} \widehat{m}(A_i)\\
&=&\sum_{i=1}^{h} m(A_i)\\
&\leq &\sum_{i=1}^{h} \sum_{n\geq 0} \sum_{j=1}^{h(n)} m(A_{n,j} \cap A_i) \ (SU1)\\
&=& \sum_{n\geq 0} \sum_{j=1}^{h(n)} \sum_{i=1}^{h} m(A_{n,j} \cap A_i) \ (Fubini)\\ 
&=& \sum_{n\geq 0} \sum_{j=1}^{h(n)} m(A_{n,j}) \ (SU2)\\ 
&=& \sum_{n\geq 0}  \widehat{m}(A_n),
\end{eqnarray*}

\bigskip \noindent which ends the solution.\\

\noindent Question (5). This question is direct and is based on the fact that a countable sum of elements of elements of $\mathcal{C}$ is a countable elements of $\mathcal{C}$ and vice-versa. The non-decreasingness of $m$ completes the solution.

\newpage
\noindent \LARGE \textbf{Doc 04-07 : Measures - Exercises on Lebesgue Measure \textit{with solutions} }. \label{doc04-07}\\
\Large

\bigskip \noindent \textbf{Exercise 1}. \label{exercise01_sol_doc04-07} Let $(\overline{\mathbb{R}}, \mathcal{B}(\overline{\mathbb{R}}), \lambda)$ the Lebesgue measure space of $\mathbb{R}$.\\

\noindent For any Borel set $A$ in $\overline{\mathbb{R}}$, for $x\in \mathbb{R}$ and $0 \neq \gamma \in \mathbb{R}$, Define

$$
A+x =\{ y+x, \ y\in A\} \text{ and } \gamma A =\{ \gamma y, \ y\in A\}.
$$

\bigskip \noindent Define the translation application $t_x$ of vector $x$ defined by : $t_x(y)=x+y$ and the linear application $\ell_{\gamma}$ of coefficient $\gamma$ defined by $\ell_{\gamma}(y)=\gamma y$.\\

\noindent Question (a) Let $A$ be Borel set. Show that $A+x=t_{x}^{-1}(A)$ and $\gamma A=\ell_{1/\gamma}^{-1}(A)$ and deduce that $A+x$ and $\gamma A$ are Borel sets.\\

\noindent Question (b) Show that the Lebesgue measure is translation-invariant, that id for any $x\in \mathbb{R}$, for any Borel set $A$, we have

$$
\lambda(A+x)=\lambda(A).
$$

\bigskip \noindent \textit{Hint} : Fix $x\in \mathbb{R}$. Consider the application $A \hookrightarrow \lambda_{1}(A)=\lambda(A+x)$. Show that $\lambda_{1}$ is an inverse measure of the Lebesgue measure. Compare $\lambda_{1}$ and $\lambda$ of the sigma-algebra

$$
\mathcal{S}=\{]a,b], -\infty \leq a \leq b \leq +\infty\}.
$$

\bigskip \noindent Conclude.\\

\noindent Question c) Show for any $\gamma \in \mathbb{R}$, for any Borel set $A$, we have

$$
\gamma A)=|\gamma| \lambda(A).
$$

\bigskip \noindent \textbf{Hint}. The case $\gamma=0$ is trivial. For $\gamma\neq 0$, proceed as in Question (a).\\

\bigskip \noindent \textbf{Solutions}.\\

\noindent Question (a). \noindent By definition of the inverse image, we have
\begin{eqnarray*}
t_{-x}^{-1}(A)&=&\{y \in \mathbb{R},\ t_{-x}(y)=y-x \in A\}\\
&=&\{y \in \mathbb{R},\ \exists t \in A, y-x=t\}\\
&=&\{y \in \mathbb{R},\ \exists t \in A, y=x+t\}\\
&=&\{x+t, t \in A\}\\
&=&A+x.
\end{eqnarray*}

\bigskip \noindent We also have

\begin{eqnarray*}
\ell_{1/\gamma}^{-1}(A)&=&\{y \in \mathbb{R},\ \ell_{1/\gamma}(y)=y/\gamma \in A\}\\
&=&\{y \in \mathbb{R},\ \exists t \in A, y/\gamma=t\}\\
&=&\{y \in \mathbb{R},\ \exists t \in A, y=t \gamma\}\\
&=&\{t \gamma, t \in A\}\\
&=& \gamma A.
\end{eqnarray*}

\bigskip \noindent Since for each $x\in \mathbb{R}$ and $\gamma\neq 0$, the application $t_{-x}$ and $\ell_{1/\gamma}$ are continuous from $\mathbb{R}$ to $\mathbb{R}$. Thus, they are measurable and the inverse images of the measurable set $A$, $t_{-x}^{-1}(A)$ and $\ell_{1/\gamma}^{-1}(A)$ are measurable. So are $x+A$ and $\gamma A$.\\

\noindent Question (b). For each $x\in \mathbb{R}$, the application

$$
\lambda_{1}(A)=\lambda(A+x), \ A\in \mathcal{A},
$$

\bigskip \noindent is the image measure of $\lambda$ by $t_{-x}$. Let us show that it coincides with $\lambda$ on  
$$
\mathcal{S}=\{]a,b], -\infty \leq a \leq b \leq +\infty\}.
$$

\bigskip \noindent For any $]a,b]$ such that $-\infty \leq a \leq b \leq +\infty$, $]a,b]+x=]a+x,b+x]$, whether $a$ and $b$ are finite or not. If one of them is infinite ($a=-\infty$ or $b=+\infty$), we have
$$
\lambda_{1}(]a,b])=\lambda(]a,b])=+\infty.
$$ 

\bigskip \noindent If not, we have
$$
\lambda_{1}(]a,b])=\lambda(]a+x,b+x])=(b+x)-(a+x)=b-a=\lambda(]a+x,b+x]).
$$

\bigskip \noindent The two measures are $\sigma$-finite and equal on $\mathcal{S}$. By \textit{Exercise 12 of Doc 04-02 of this Chapter}, these measure are equal on $\sigma(\mathcal{S})=\mathcal{B}(\mathbb{R})$.\\

\noindent Question (b). For each $0\neq \gamma \in \mathbb{R}$, the application

$$
\lambda_{1}(A)=\lambda(\gamma A), \ A\in \mathcal{A},
$$

\bigskip \noindent is the image measure of $\lambda$ by $\ell_{1/\gamma}$. Also the application

$$
\lambda_{2}(A)=|\gamma| \lambda(A)
$$

\bigskip \noindent is also a measure (See \textit{Exercise 5 above in Doc 04002 of this chapter}). Let us show that  $\lambda_1$ and $\lambda_2$ coincide on  
$$
\mathcal{S}=\{]a,b], -\infty \leq a \leq b \leq +\infty\}.
$$

\bigskip \noindent We have to do it for $\gamma<0$ and for $\gamma>0$. Let us do it just for $\gamma<0$. For any $]a,b]$ such that $-\infty \leq a \leq b \leq +\infty$, 
$\gamma [\gamma b, \gamma a[$, whether $a$ and $b$ are finite or not. If one of them is infinite ($a=-\infty$ or $b=+\infty$), we have

$$
\lambda_{1}(]a,b])=\lambda_2(]a,b])=+\infty.
$$ 

\bigskip \noindent If not, we have
$$
\lambda_{1}(]a,b])=\lambda([\gamma b, \gamma a[)=\gamma a -\gamma b=-\gamma(b-a)=\|\gamma\| \lambda(]a,b])=\lambda_2(]a,b]).
$$

\bigskip \noindent The two measures are $\sigma$-finite and equal on $\mathcal{S}$. By \textit{Exercise 12 of Doc 04-02 of this Chapter}, these measure are equal on $\sigma(\mathcal{S})=\mathcal{B}(\mathbb{R})$.\\

\bigskip \noindent \textbf{Exercise 2}. \label{exercise02_sol_doc04-07}  Let $(\mathbb{R}, \mathcal{B}(\mathbb{R}))$ be the Borel space.\\

\noindent Question (a) Consider the Lebesgue measure $\lambda$ on $\mathbb{R}$.\\

\noindent (a1) Show that for any $x \in \mathbb{R}$, $\lambda(\{x\})=0$.\\

\noindent (a2) Show that for countable $A$ subset of $\mathbb{R}$, $\lambda(\{A\})=0$.\\

\noindent Question (b) Consider a Lebesgue-Stieljes measure $\lambda_F$ on $\mathbb{R}$ where $F$ is a distribution function.\\

\noindent (b1) Show that for any $x \in \mathbb{R}$, 

$$
\lambda_F(\{x\})=F(x)-F(x-0),
$$

\bigskip \noindent where $F(x-0)$ is the left-hand limit of $F$ at $x$.\\

\noindent Deduce from this that $\lambda_F(\{x\})$ if and only if $F$ is continuous at $x$.\\

\noindent (b2) Show that for any $\lambda_F(\{x\})\neq 0$ for at most a countable number of real numbers.\\

\noindent Hint. Combine this with \textit{Exercise 1 in DOC 03-06}. \\  

\bigskip \noindent \textbf{Solutions}.\\

\noindent Question (a).\\

\noindent (a1). For all $x\in\mathbb{R}$, we have
$$
\{x\} = \bigcap_{n\geq 1} ]x-1/n, x].
$$

\bigskip \noindent Actually, we have $]x-1/n, x]\downarrow {x}$. By the continuity of the measure, we have $\lambda(]x-1/n, x])=1/n \downarrow \lambda(\{x\})$. So, $\lambda({x})=0$. We also may say that

$$
\forall \delta>0, \{x\} \subset ]x-\delta, x+\delta], 
$$ 

\bigskip \noindent which, by the non-decreasingness of the measure, implies that

$$
\forall \delta>0, 0\leq \lambda(\{x\}) \leq 2\delta. 
$$ 

\bigskip \noindent By letting $\delta \downarrow 0$, we get that $\lambda(\{x\})=0$.\\

\noindent (a2). If $A$ is countable, it is a countable union of the singletons of its elements, which are $\lambda$-null sets. So it is a  $\lambda$-null set by \textit{Exercise 7 in Doc 04-02 of this chapter}.\\

\noindent Question (b).\\

\noindent (b1). For all $x\in\mathbb{R}$, from 

$$
\{x\} = \bigcap_{n\geq 1} ]x-1/n, x].
$$

\bigskip \noindent We get that $\lambda_F(]x-1/n, x])=F(x)-F(x-1/n) \downarrow \lambda_F(\{x\})$. Since $F$ is non-decreasing and $F(x)$ is finite, the limit of $F(x-h)$ as $h$ tends to zero by positive values exists, is bounded below by $F(x)$ and  is the monotone limit of $F(x-h)$ as $h \downarrow 0$. By definition, this limit is the left-limit $F(x-0)$ at $x$. As a sub-limit, we have 
$F(x-1/n) \downarrow F(x-0)$ as $n\uparrow +\infty$. We conclude that

$$
\lambda_F(\{x\})=F(x)-F(x-0).
$$

\bigskip \noindent Since $F$ is right-continuous, the fact $\lambda_F(\{x\})=0$ implies that $F(x-0)=F(x)$, that is $F$ is left-continuous and thus, continuous. Thus, $F$ is continuous at $x$ if and only if
$\lambda_F(\{x\})=0$.\\

\noindent (b2). By \textit{Exercise 1 in Doc 03-06 in Chapter 3}, the monotone function $F$ has at most a countable number of discontinuity points. Hence $\lambda_F(\{x\})\neq 0$ for at most a countable of real numbers $x$.\\

\bigskip \noindent \textbf{Exercise 3}. \label{exercise03_sol_doc04-07} The objective of this exercise is to show that the Lebesgue measure is the unique measure on $(\mathbb{R}, \mathcal{B}(\mathbb{R}))$ which is translation invariant and which assigns the unity value to the interval $]0,1]$.\\

\noindent Let $(\mathbb{R}, \mathcal{B}(\mathbb{R}), m)$ be a measure space.\\

\noindent Question (a) Show that if $m$ is translation-invariant and $m(]0,1])=1$, then $m$ is the Lebesgue measure $\mathbb{R}$.\\

\noindent Hint. Use the translation-invariance of $m$ to prove that $m(]0,1/q])=1/q$, for any integer $q\geq 1$. Next, prove that $m(]0,p/q])=p/q$ for $p\geq 1$, $q\geq 1$, $p/q\leq 1$. Next, Using the continuity of $m$ on $(0,1)$ and the density of the set of rational numbers in $\mathbb{R}$, show that $m(]0,b])=b$ for any real number $b\in (0,1)$. Finally, by translation-invariance of $m$, proceed to your conclusion.\\

\noindent Question (b) Deduce from Question (a) that : if $m$ assigns to open intervals positive values and is translation-invariant, then there exists $\gamma >0$, $m=\gamma \lambda$.\\

\noindent Hint. Put $\gamma=m(]0,1])$ and proceed.\\

\bigskip \noindent \textbf{Solution}.\\

\noindent Question (a). \noindent Let us proceed by steps.\\

\bigskip \noindent Step 1. Let $q \in \mathbb{Q}$ and $q\geq 1$. Let us show that $m(]0,1/q])=1/q$. We have

$$
]0,1]=\sum_{i=1}^{q} \left] \frac{i-1}{q}, \  \frac{i}{q} \right] =: \sum_{i=1}^{q} I(i,q). \ (F1)
$$

\bigskip \noindent Each $I(i,q)=\left] \frac{i-1}{q}, \ \frac{i}{q} \right]$, for $1\leq i\leq q$, is a translation of $]0,1/q]$, that is

$$
\left] \frac{i-1}{q}, \ \frac{i}{q} \right]=\left] 0, \ \frac{1}{q} \right] + \frac{i-1}{q}.
$$

\noindent So, by translation-invariance, we have $m(I(i,q))=m(]0,1/q])$ for all $1\leq i\leq q$, and by additivity of $m$, we get from (F1)

$$
m(]0,1])= \sum_{i=1}^{q} m(I(i,q))=q m(]0,1/q]).
$$

\bigskip \noindent This gives that $m(]0,1/q])=1/q$.\\

\noindent Step 2. Let $(p,q)\in \mathbb{Q}^2$, $p\geq 1$, $q\geq 1$, $p/q \leq 1$. Let us us show that $m(]0,p/q])=p/q$. We have

$$
]0,p/q]=\sum_{i=1}^{p} I(i,q),
$$

\bigskip \noindent which gives

$$
m(]0,p/q])=\sum_{i=1}^{p} m(I(i,q))=p m(]0,1/q])=p/q.
$$

\bigskip \noindent Step 3. Consider $x \in ]0,1]$. Let us show that $m(]0,x])=x$. if $x=1$, we have $m(]0,x])=x$ by assumption. If $0<x<1$, there exists a sequence 
$(r_n)_{(n\geq 1)} \subset \mathbb{Q}\ \cap ]x,1]$ such that $r_n \downarrow x$ as $n\uparrow +\infty$. Thus, $]0,r_n] \downarrow ]0,x]$ and hence

$$
m(]0,x])=\lim_{n \rightarrow +\infty} m(]0,r_n]) = \lim_{n \rightarrow +\infty} r_n=x.
$$ 

\bigskip \noindent Step 4. For any integer $p$, $m(]0,p])=p$ since $]0,p]$ is the sums of the $](i-1),i]$'s for $i=1$ to $p$, which are translation on $]0,1]$.\\

\noindent Final step. By translation-invariance, we deduce from  for any $c \in \mathbb{R}$ and any $x\in ]0,1]$,

$$
m(]c,c+x])=m(]0,x]+c)=m(]0,x])=x.
$$

$$
m(]a,b]=m(]0,b-a])
$$

\bigskip \noindent Now, consider the decomposition (b-a)=[b-a]+h, where $[b-a]$ is the integer part of $b-a$, that is the greatest integer less or equal to $b-a$, and $0\leq h <1$. Thus, we use the decomposition
$]0,[b-a]+h]=]0,[b-a]]+][b-a],[b-a]+h]$ to have

$$
m(]a,b]=m(]0,b-a])=m(]0,[b-a])+m(][b-a],[b-a]+h])=[b-a]+h=(b-a).
$$

\bigskip \noindent This concludes the solution of Question(a).\\

\noindent Question (b). By putting $m_0=m/m(]0,1])$, $m_0$ satisfies the assumptions in Question (a), and $m_0$ is is the Lebesgue measure. So we have $m=m/m(]0,1]) \lambda$.

\bigskip \noindent The next exercise may be skipped at a first reading. Instead, you may read to solution and use the knowledge in.\\

\newpage
\bigskip \noindent \textbf{Exercise 4}. \label{exercise04_sol_doc04-07}  Vitali's construction of a non measurable subset on $\mathbb{R}$.\\

\noindent In this exercise, a non-measurable set on $\overline{\mathbb{R}}$ is constructed. The exercise focuses on $]0,1]$ on which is uses the
circular translation, in place of the straight translation by x : $t\hookrightarrow t+x,$ on the whole real line.\\

\noindent Define the circular translation on $]0,1]$ for any fixed $x$, $0<x<1$ by

\begin{equation*}
\begin{tabular}{llll}
$T_{x}:$ & $]0,1]$ & $\longrightarrow $ & $]0,1]$ \\ 
& $y$ & $\hookrightarrow $ & $\left\{ 
\begin{tabular}{lll}
$z=y+x$ & if & $y+x\leq 1$ \\ 
z=$x+y-1$ & if & $y+x>1$%
\end{tabular}%
\right. $,
\end{tabular}
\end{equation*}

\bigskip \noindent and denote 

$$
T_{x}(A)=A\oplus x.
$$

\bigskip \noindent Question (a) Show that for any $x$, $0<x<1$, the application from  $]0,1]$ to $]0,1]$, defined by

\begin{equation*}
\begin{tabular}{llll}
$T_{-x}:$ & $]0,1]$ & $\longrightarrow $ & $]0,1]$ \\ 
& $z$ & $\hookrightarrow $ & $\left\{ 
\begin{tabular}{lll}
$y=z-x$ & if & $z-x>0$ \\ 
$y=z-x+1$ & if & $z-x\leq 0$%
\end{tabular}%
\right. $%
\end{tabular}%
\end{equation*}%

\bigskip \noindent is the inverse of $T_x$. We denote

$$
T_{-x}(A)=A\ominus x.
$$

\bigskip \noindent Question (b) Show that for $0<x<1$, for any Borel set in $]0,1]$, we have
$$
A\oplus x=T_{-x}^{-1}(A)
$$

\bigskip \noindent and conclude $A\oplus x$ is measurable and the graph

\begin{equation*}
A\hookrightarrow \lambda_{x}(A)=\lambda (A\oplus x)
\end{equation*}

\bigskip \noindent defines a measure on the Borel sets in $]0,1]$.

\noindent Question (c) Show that $\lambda_{x}=\lambda$.\\

\noindent Hint. Consider an interval $]y, y^{\prime}]$ in $]0,1]$. Use the following table to give  $]y, y^{\prime}]\oplus x$ with respect to the positions of $y+x$ and $y^{\prime}+x$ with respect to the unity.

\begin{equation*}
\begin{tabular}{|l|l|l|l|l|}
\hline
case & y+x & $y^{\prime}+x$ & $]y,y^{\prime }]\oplus x$ & $\lambda (]y,y^{\prime }]\oplus x)$ \\ \hline
1 & $y+x\leq 1$ & $y^{\prime }+x\leq 1$ &  &  \\ 
\hline
2 & $y+x\leq 1$ & $y^{\prime }+x>1$ &  &  \\ 
\hline
3 & $y+x>1$ & $y^{\prime }+x\leq 1$ & impossible & impossible \\ 
\hline
4 & $y+x>1$ & $y^{\prime }+x>1$ &  & \\
\hline
\end{tabular}
\end{equation*}

\bigskip \noindent \textbf{Extension}. By a very similar method, we also may see that $A\hookrightarrow \lambda_{-x}(A)=\lambda (A\ominus x)$ is a measure and that $\lambda_{-x}=\lambda$.\\

\noindent We may unify the notation by writing $A\ominus x=A\oplus (-x)$ for $-1<x<0$. We have proved that

\begin{equation*}
\forall (-1<x<1),\forall A\in B(]0,1]),\lambda (A\oplus x)=\lambda (A).
\end{equation*}

\bigskip \noindent Question (d).\\

\noindent Consider the binary relation $\mathbb{R}$ on $]0,1]$ defined by%
\begin{equation*}
\forall (y,z)\in ]0,1]^{2},y\mathcal{R}z\Longleftrightarrow y-z\in \mathbb{Q} \cap ]0,1],
\end{equation*}

\bigskip \noindent that is $y$ and $z$ are in relation if and only if their difference $y-z$ is a rational number, that is also : there exists $r\in \mathbb{Q},$ $y=z+r$
with $-1<r<1$.\\

\noindent (d1) Show that $\mathcal{R}$ is an equivalence relation.\\

\noindent Form a set $H$ by choosing one element in each equivalence class (apply the choice axiom).\\

\noindent (d2) Show that for any $0\leq r<1,H\oplus r\subset ]0,1]$.\\

\noindent (d3) Show that we have

\begin{equation*}
]0,1]\subset \sum\limits_{r\in Q\cap ]0,1[}H\oplus r.
\end{equation*}

\bigskip \noindent (d4) Show that for $0\leq r_{1}\neq r_{2}<1,H\oplus r_{1}$ and $H\oplus r_{2}$ are
disjoint, that is
\begin{equation*}
(H\oplus r_{1})\cap (H\oplus r_{2})=\emptyset .
\end{equation*}

\bigskip \noindent Conclude that 

\begin{equation*}
]0,1]=\sum\limits_{r\in Q\cap ]0,1[}H\oplus r.
\end{equation*}

\noindent (d5) Suppose that $H$ is measurable. Apply the $\sigma$-additivity of $\lambda$ and show that that this leads to an absurdity.\\

\noindent Conclusion. Did you get a nonempty non-measurable set on $\mathbb{R}$?\\

\bigskip \noindent \textbf{Solution}. The solution is given in a self-written document in the next page.\\

\newpage
\noindent \textbf{Vitali's construction of a non measurable subset on $\mathbb{R}$}.\\

\noindent We are going to exhibit a non-measurable subset on $\mathbb{R}$ when endowed with the usual sigma-algebra. We will focus on $]0,1]$ on which we use the circular translation, in place of the straight translation by x : $t\hookrightarrow t+x$ on the whole real line.\\

\bigskip \noindent We begin to define the circular translation and show that the Lebesgue measure is invariant by circular translations. Next we use the translation-invariance of the Lebesgue measure to exhibit a non-measurable subset of $\mathbb{R}$.\\

\bigskip \noindent \textbf{1 - Circular translation}. Let $\mathbb{R}\ni x>0,0<x<1.$ Define the
circular translation on $]0,1]$ by%
\begin{equation*}
\begin{tabular}{llll}
$T_{x}:$ & $]0,1]$ & $\longrightarrow $ & $]0,1]$ \\ 
& $y$ & $\hookrightarrow $ & $\left\{ 
\begin{tabular}{lll}
$z=y+x$ & if & $y+x\leq 1$ \\ 
z=$x+y-1$ & if & $y+x>1$%
\end{tabular}
\right. $.
\end{tabular}
\end{equation*}

\bigskip \noindent We have to understand how this circular translation is defined. If $y+x$ is still in $]0,1]$, we let is as such. But, if is exceeds $1$, that is $y+x>1$, we associate $y$ to $x+y-1$. That is : we want to $y$ to $y+x$. When we reach $1$, we go back to $0$ and continue our process and we stop at $y+x-1$ : the total increase will be $(1-y)$ plus $(x+y-1)$, that is exactly $x$.\\ 

\bigskip \noindent We may simply write
\begin{equation*}
T_{x}(y)=(y+x)1_{(y\leq 1-x)}+(y+x-1)1_{(y>1-x)}.
\end{equation*}

\bigskip \noindent This function is a bijection and its inverse is $T_{-x}$ defined by 
\begin{equation*}
\begin{tabular}{llll}
$T_{-x}:$ & $]0,1]$ & $\longrightarrow $ & $]0,1]$ \\ 
& $z$ & $\hookrightarrow $ & $\left\{ 
\begin{tabular}{lll}
$y=z-x$ & if & $z-x>0$ \\ 
$y=z-x+1$ & if & $z-x\leq 0$%
\end{tabular}%
\right. $,
\end{tabular}
\end{equation*}

\bigskip \noindent that is
\begin{equation*}
T_{-x}(z)=(z-x)1_{(z>x)}+(z-x+1)1_{(z\leq x)}.
\end{equation*}

\bigskip \noindent \textbf{Image measure of the Lebesgue Measure by $T_{-x}$}.\\

\noindent The precedent section implies that for any Borel set $A\subset ]0,1]$, its translation 
\begin{equation*}
T_{x}(A)=T_{-x}^{-1}(A)
\end{equation*}

\bigskip \noindent is a Borel as the inverse image of $A$ by $T_{-x}^{-1}$. We denote
\begin{equation*}
T_{x}(A)=A\oplus x
\end{equation*}

\bigskip \noindent and we may define the image measure $\lambda_{x}$ of $\lambda$ by $T_{-x}^{-1}$ :

\begin{equation*}
A\hookrightarrow \lambda_{x}(A)=\lambda (A\oplus x).
\end{equation*}

\bigskip \noindent \textbf{Fact}. We are going to show that $\lambda_{x}=\lambda$ for all $0\leq x <1$. To prove this, it will be enough to show that $\lambda_{x}$ assigns to each interval $]y,y^{\prime }] \subset ]0,1]$ the value $]y,y^{\prime}-y$, which is a characterization of the Lebesgue measure.

\bigskip \noindent You may check this using the table below :

\begin{equation*}
\begin{tabular}{|l|l|l|l|l|}
\hline
case & y+x & $y^{\prime}+x$ & $]y,y^{\prime }]\oplus x$ & $\lambda (]y,y^{\prime }]\oplus x)$ \\ 
\hline
1 & $y+x\leq 1$ & $y^{\prime }+x\leq 1$ & $]y+x,y^{\prime }+x]$ & $y\prime
-y=\lambda (]y,y^{\prime }])$ \\ 
\hline
2 & $y+x\leq 1$ & $y^{\prime }+x>1$ & $]y+x,1]+]0,y^{\prime }+x-1]$ & $y\prime
-y=\lambda (]y,y^{\prime }])$ \\ 
\hline
3 & $y+x>1$ & $y^{\prime }+x\leq 1$ & impossible & impossible \\ 
\hline
4 & $y+x>1$ & $y^{\prime }+x>1$ & $]y+x-1,y^{\prime }+x-1]$ & $y\prime
-y=\lambda (]y,y^{\prime }])$\\
\hline
\end{tabular}
\end{equation*}

\bigskip \noindent Remarks.\\

\noindent  (1) \textit{Case 3 }: Since $0<y\leq y^{\prime }\leq 1,$ we cannot have $%
y^{\prime }+x\leq 1$ if $y+x>1$.\\

\noindent (2) \textit{Cases 1 and 4} are easy to see. Let us have a close sight at Case 2. On the whole line, $]y,y^{\prime }]+x=]y+,y^{\prime }+x]$ with the straight
translation. Now, since $y^{\prime }+x$ exceeds 1, the circular translation gives first the interval $]y+x,1]$ and we transfer the exceedance $x+y\prime
-1$ to form the interval $]0,y^{\prime }+x-1].$ We may see that $y^{\prime}+x-1\leq $ $y+x,$ otherwise we would have ($y^{\prime}+x-1>y+x\Longleftrightarrow y^{\prime }-y>1$) and this is impossible. Then
the intervals $]0,y^{\prime }+x-1]$ and\ $]y+x,1]$ are disjoint and its Lebesgue measure $y^{\prime }-y$.\\

\noindent We successfully checked that $\lambda_x(]y,y^{\prime}])=y^{\prime}-y$ for all intervals $]y,y^{\prime }] \subset ]0,1]$, with 
$y \leq y^{\prime}$. Then $\lambda_{x}=\lambda$ for all $0\leq x <1$.\\

\noindent So far, we have considered translations by $0\leq x <1$. We may do the same for negative $x$, that is for $-1 < x\leq 0$. Let us do it but with less details. consider a circular translation by a negative number and denote for 

\begin{equation*}
T_{-x}(A)=A\ominus x.
\end{equation*}

\bigskip \noindent We may also check that the measure defined by $\lambda_{-x}(A)=\lambda (A\ominus x)$ is the Lebesgue measure itself by checking 
the following table

\begin{equation*}
\begin{tabular}{|l|l|l|l|l|}
\hline
Cases & $y-x$ & $y^{\prime}-x$ & $]y,y^{\prime }]\ominus x$ & $\lambda (]y,y^{\prime }]\ominus x)$
\\ 
\hline
1 & $y-x>0$ & $y^{\prime }-x>0$ & $]y-x,y^{\prime }-x]$ & $y^{\prime }-y=\lambda
(]y,y^{\prime }])$ \\ 
\hline
2 & $y-x\geq 0$ & $y^{\prime }+x<0$ & $]0,y^{\prime }-x]+]y-x+1,1]$ & $y^{\prime
}-y=\lambda (]y,y^{\prime }])$ \\ 
\hline
3 & $y-x>1$ & $y^{\prime }-x\geq 0$ & impossible & impossible \\ 
\hline
4 & $y-x\leq 0$ & $y^{\prime }-x\leq 0$ & $]y+x+1,y^{\prime }+x+1]$ & $y^{\prime
}-y=\lambda (]y,y^{\prime }])$\\
\hline
\end{tabular}
\end{equation*}

\bigskip \noindent Remark again that the straight translation on the whole line gives $]y,y^{\prime }]-x=]y-x,y^{\prime }-x]$. In this case, the part $]y-x,0]$ is outside $]0,1]$
and the circular translation applies to it and turns it into $]y-x+1,1]$. Remark also that $y^{\prime }-x\leq y-x+1,$ otherwise we would have $y^{\prime }-x>y-x+1\Longleftrightarrow y^{\prime }>y+1\Longrightarrow y^{\prime }>1$, which is impossible.\\

\noindent We use the same argument to see that
\begin{equation*}
\forall A\in B(]0,1]),\lambda (A\circleddash x)=\lambda (A).
\end{equation*}

\bigskip \noindent Finally, we have

\begin{equation*}
\forall 0\leq x<1,\forall A\in B(]0,1]),\lambda (A\oplus x)=\lambda(A\circleddash x)=\lambda (A).
\end{equation*}

\bigskip \noindent To unify the notation, we put for $x\leq 0,$%
\begin{equation*}
A\oplus x=A\circleddash (-x)
\end{equation*}

\bigskip \noindent and we may say
\begin{equation*}
\forall (-1<x<1),\forall A\in B(]0,1]),\lambda (A\oplus x)=\lambda (A).
\end{equation*}

\bigskip \noindent \textbf{2 - Construction of the non-measurable set $H$}.\\

\noindent Consider the binary relation $\mathcal{R}$ on $]0,1]$ defined by 
\begin{equation*}
\forall (y,z)\in ]0,1]^{2},y\mathcal{R}z\Longleftrightarrow y-z\in \mathbb{Q} \cap ]0,1],
\end{equation*}

\bigskip \noindent that is $y$ and $z$ are in relation if and only if their $y-z$ difference is a rational number, that is also : there exists $r\in \mathbb{Q}$, $y=z+r$ with $-1<r<1.$ It is easy to see that 
$\mathcal{R}$ is an equivalence relation. Now form the set $H$ by choosing one element in each equivalence class (apply the choice axiom).\\

\noindent We have to see that 

\begin{equation*}
]0,1]=\sum\limits_{r\in Q\cap ]0,1[}H\oplus r.
\end{equation*}

\bigskip \noindent We will have this after we have checked these three points :\\

\noindent (i) For any $0\leq r<1,H\oplus r\subset ]0,1]$.\\

\noindent  (ii) We have

\begin{equation*}
]0,1]\subset \sum\limits_{r\in Q\cap ]0,1[}H\oplus r.
\end{equation*}

\bigskip \noindent (iii) For $0\leq r_{1}\neq r_{2}<1,H\oplus r_{1}$ and $H\oplus r_{2}$ are disjoint, that is
\begin{equation*}
(H\oplus r_{1})\cap (H\oplus r_{2})=\emptyset .
\end{equation*}

\bigskip\noindent \textbf{Proofs of the three points}.\\

\noindent Point (i) holds by the definition of circular translations.\\

\noindent Point (ii). let $z\in ]0,1],$ then its class is represented by one element its $H$, say $x$. Then $z$ and $x$ are in relation, that is, $z-x=r\in ]-1,1[$ and $r$ rational. Thus  $z=x+r=x\oplus r\in H\oplus r.$ Point (ii) is proved.\\

\noindent Point (iii). Suppose we can find one element in the intersection between $H\oplus r_{1}$ and $H\oplus r_{2}$ with $r_{1}\neq r_{2}.$ So $z$ may be written as $z=x_{1}+r_{1}$ with $x\in H$ and $z=x_{2}+r_{2}$ with $x_{2}\in H.$ We have that $x_{1}$ and $x_{2}$ are not equal since this would lead to $r_{1}=r_{2}$, and this is contrary to the assumptions. Now suppose that $x_{1}$ and $x_{2}$ are different. We have $z=x_{1}+r_{1}=x_{2}+r_{2}$ leads to $x_{1}-x_{2}=r_{1}-r_{2}.$ So  $x_{1}$ and $x_{2}$ are in relation. We have that $x_{1}$ and $x_{2}$ different and belong to the same equivalence class. This is not possible since we have exactly one point of each equivalence in $H$. So $x_{1}$ and $x_{2}$ neither equal nor equal. This is absurd. So the intersection between $H\oplus r_{1}$ and $H\oplus r_{2}$ is empty when $r_{1}\neq r_{2}$.\\

\noindent Put together these three points to conclude.\\

\noindent \textbf{3 - Conclusion : H is not measurable}.\\

\noindent Suppose H is measurable. So the sets $H\oplus r,$ $r\in Q\cap ]0,1[,$ form a countable family of measurable sets, and then%
\begin{equation*}
\lambda (]0,1])=\sum\limits_{r\in Q\cap ]0,1[}\lambda (H\oplus r).
\end{equation*}

\bigskip \noindent Since $\lambda (H\oplus r)=\lambda (H)$ $\in \lbrack 0,1]$ for $r\in \mathbb{%
Q}\cap ]0,1[,$ we get%
\begin{equation*}
1=Card(\mathbb{Q}\cap ]0,1[)\times \lambda (H)=+\infty \times \lambda (H).
\end{equation*}

\bigskip \noindent There exists no real number $\lambda (H)$ $\in \lbrack 0,1]$ such that $%
1=+\infty \times \lambda (H).$\\

\noindent So the assumption that $H$ is measurable is impossible.\\

\noindent \textbf{\LARGE Conclusion : the set $H$ is not measurable}.\\

\newpage
\noindent \LARGE \textbf{Doc 04-08 : Measures - General exercises \textit{with solutions}  }. \label{doc04-08}\\
\Large

\bigskip \noindent \textbf{Exercise 1}. \label{exercise01_sol_doc04-08} Let $(E, \mathcal{B}(E),m)$ be a finite measure ($(E, \mathcal{B}(E))$ is a Borel space. Show that any other measure that is equal to $m$ on the class of of open sets, or on the class of closed sets is equal to $m$ on $\mathcal{B}(E)$. Conclusion : A finite measure on a Borel space us regular, $i.e$, it is determined by its values on the open sets, by its values on the closed sets.\\

\bigskip \noindent \textbf{Solution}. The $\sigma$-algebra $\mathcal{B}(E)$ is generated by the class of open sets which is a $\pi$-class containing the full space $E$. So, by the $\pi-\lambda$ lemma (See \textit{Exercise 11 in Doc 04-02 of this Chapter}), two finite measure equal on the class of open sets are equal. The reasoning is the same for the class of closed sets, which is also a $\pi$-class containing the full space $E$ and which generates $\mathcal{B}(E)$.\\

\bigskip \noindent \textbf{Exercise 2}.  \label{exercise02_sol_doc04-08} (Approximation of a finite measure by its values on open and closed sets in a metric space). Let $(S,d,m)$ be a finite measure on the metric space $(S,d)$ endowed with its Borel $\sigma$-algebra.  \noindent Let $A \in \mathcal{B}(S)$. Show the assertion $\mathcal{P}(A)$ : for any $\varepsilon>0$, there exist a closed set
$F$ and an open set $G$ such that : $F \subset A \subset G$ and $m(G \setminus F)<\varepsilon$.\\

\noindent \textit{Hint}. Proceed as follows.\\

\noindent Question (a) Let $A$ be closed set. Fix $\delta>0$ and set the $\delta$-dilation or $\delta$-enlargement of $G$ as 

\begin{equation*}
G_{\delta}=\{x,\text{ }d(x,A)<\delta \}.
\end{equation*}

\bigskip \noindent (a) Show that $G_{\delta}$ is an open set (or you may skip this as a result of topology).\\

\noindent Let $(\delta_n)_{n\geq 1}$ be sequence of positive number such that $\delta_n \downarrow +\infty$ as $n \uparrow+\infty$.\\

\noindent (b) By using the continuity of $m$ with the sequence $(G_{\delta_n})_{n\geq 1}$, and by noticing that $A \subset G_{\delta_n}$ for any $n$, show that for any $\varepsilon>0$, there exists $\delta>0$ such that

\begin{equation*}
m(G_{\delta} \setminus A)<\varepsilon.
\end{equation*}

\bigskip \noindent (c) Define by $\mathcal{D}$ the class of measurable sets $A$ such that for any $\varepsilon>0$, there exist a closed set
$F$ and an open set $G$ such that : $F \subset A \subset G$ and $m(G \setminus F)<\varepsilon$.\\

\noindent Show that $\mathcal{D}$ is a $\sigma$-algebra including the class $\mathcal{F}$ of closed sets. Conclude that the assertion $\mathcal{P}$  holds.\\

\noindent Hint. To show that $\mathcal{D}$ is stable under countable union, consider a sequence $A_n \in \mathcal{D}$, $n\geq 0$ and next
$$
A =\bigcup_{n\geq 0} A_n. 
$$

\bigskip \noindent Fix $\varepsilon>0$. By definition, for each $n\geq 0$, there exist a closed set$H_n$ and an open set $G_n$ such that : $H_n \subset A_n \subset G_n$ and $m(H_n \setminus F_n)<\varepsilon/2^{n+1}$.\\

\noindent Denote

$$
F_p=\bigcup_{1\leq n \leq p} H_n, \ p\geq 0, \ H=\bigcup_{n\geq 0} H_n \text{ and } G=\bigcup_{n\geq 0} G_n. 
$$

\bigskip \noindent By using the continuity of $m$, explain that there exists $p_0\geq 0$ show that

$$
m(H \setminus F_{p_0}) < \varepsilon.
$$

\bigskip \noindent Set $F=F_{p_0}$ and justify :

$$
G \setminus F= (H \setminus F) + (G \setminus H).
$$

\bigskip \noindent Finally, use \textit{Formula 00.12 in Doc 00-01} or \textit{Exercise 1, Doc 00-03} in Chapter 0 to conclude.\\

\bigskip \noindent \textbf{Solution}. We fix $\varepsilon>0$ for once. \\

\noindent (a/b) We begin to establish $\mathcal{P}$ when $A$ is a closed set. Using results of topology, $G_{\delta_n}$ is an open set for all $n\geq 0$. If $\delta_n \downarrow 0$ ans $n \uparrow +\infty$, $G_{\delta_n} \downarrow \overline{A}=A$, where $\overline{A}$ is the closure of $A$ which is still $A$ since $A$ is closed. Since $m$ is finite, its continuity implies that $m(G_{\delta_n}) \downarrow m(A)$. Thus, there exists $n_0$ such that for all $n\geq n_0$, we have

$$
m(G_{\delta_n}) - m(A)=m(m(G_{\delta_n} \setminus A) < \varepsilon.
$$ 

\bigskip \noindent By taking $G$ as any of the $m(G_{\delta_n})$, $n\geq n_0$, we get a closed set $F=A$ and an open set $G$, such that $F \subset A \subset G$ such that $m(G \setminus F)<\varepsilon$. Conclusion : the assertion $\mathcal{P}(A)$ if $A$ is closed.\\

\noindent (c) Denote

$$
\mathcal{D}=\{A \in \mathcal{S}, \mathcal{P}(A) \ holds\}.
$$

\bigskip \noindent (c1) The class $\mathcal{F}$ of closed sets is contained in $\mathcal{D}$ : $\mathcal{F} \subset \mathcal{D}$.\\

\noindent (c2) Let us show that $\mathcal{D}$ is a $\sigma$-algebra.\\

\noindent (i) $S \in $. Take $G=F=S$ to see this.\\

\noindent (ii) $\mathcal{D}$ is stable by taking complements. Indeed, let $A \in \mathcal{D}$. So, there exist $G$ open set, $F$ closed set such that $F \subset A \subset G$ such that $m(G \setminus F)<\varepsilon$. But we have

$$
F_0=G^c \subset A^c \subset G_0=F^c
$$

\bigskip \noindent and

\begin{eqnarray*}
m(G_0 \ F_0)&=&m(G_0)-m(F_0)= (m(\Omega)-m(F)) \ - \ (m(\Omega)-m(G))\\
&=&m(G)-m(F)=m(G\setminus F)< \varepsilon, 
\end{eqnarray*}

\bigskip \noindent where we were able to drop $m(\Omega)$ since $m(\Omega)$ is finite. Since $G_0$ and $F_0$ are respectively open, we get that $A^c \in \mathcal{D}$.\\

\noindent (iii) Let $(A_n)_{(n\geq 0)}$ be a sequence in $\mathcal{D}$. Let us show that $A=\cup_{n\geq 0} A_n$ is in $\mathcal{D}$. Now for each $n\geq 0$ fixed, we there exists $H_n$ closed, $G_n$ open such that $H_n \subset A_n \subset G_n$ and $m(G_n \setminus H_n)<\varepsilon/2^{n+1}$. Let us denote

$$
F_p=\bigcup_{1\leq n \leq p}, \ p\geq 0, \ H=\bigcup_{n\geq 0} F_n \text{ and } G=\bigcup_{n\geq 0} G_n. 
$$

\bigskip \noindent $G$ is open. But we do not know if $H$ is closed. But for sure, for each $p\geq 0$, $F_p$ s closed and $F_p \uparrow H$ as $p\uparrow +\infty$. Si, by the continuity of the measure and the the finiteness of $m(H)$, we have that there exists a value $p_0$ such that

$$
m(F)-m(F_{p_0})=m(H \setminus F_{p_0}) <\varepsilon/2.  (AP1)
$$

\bigskip \noindent Now, it is clear that $F=F_{p_0}$ is closed, that $G$ is open and that $F \subset A \subset G$. It remains to show that $m(G \setminus F)<\varepsilon$. Since $F \subset H$, we have

$$
F^c= (H \setminus F) + H^c,
$$

\bigskip \noindent so that

$$
G \setminus F=G F^c= G (H \setminus F) + G H^c=G H F^c + G H^c \subset H F^c + G H^c=(H \setminus F) + (G\setminus H).
$$

\bigskip \noindent where we used that $G H F^c \subset H F^c$, and that $H F^c$ and $G H^c$ are still disjoint.  We also use that in all the operations setminus ($\setminus$), the right set is contained in the left set. We get

$$
m(G \setminus F) \leq m(H \setminus F) + m(G\setminus H). \ \ (AP2)
$$

\bigskip \noindent By using \textit{Formula 00.12 in Doc 00-01 in Chapter 0}, we have
  
$$
G \setminus H = \left( \bigcup_{n\geq 0} G_n \right) \setminus \left( \bigcup_{n\geq 0} F_n \right) \subset \bigcup_{n\geq 0} (G_n \setminus F_n),
$$

\bigskip \noindent and next,

\begin{eqnarray*}
m(G \setminus H) &\leq& m\left( \bigcup_{n\geq 0} (G_n \setminus F_n) \right)\\
&\leq& \sum_{n\geq 0} m(G_n \setminus F_n)\\
&< & \sum_{n\geq 0} \varepsilon/2^{n+1}=\varepsilon/2.
\end{eqnarray*}

\bigskip \noindent So we have obtained that

$$
m(G \setminus H) < \varepsilon/2. \ \ (AP3)
$$

\bigskip \noindent By plugging (AP1) and (AP3) in (AP2), we get the result, $i.e$, $m(G \setminus F)<\varepsilon/2$.\\

\noindent (c3) Conclusion : By construction, $\mathcal{D}$ is in $\mathcal{B}(S)$. And $\mathcal{D}$ contains $\mathcal{B}(S)$ since $\mathcal{D}$ is a $\sigma$-algebra containing $\mathcal{F}$ which is a generator of $\mathcal{B}(S)$. So, we have the equality $\mathcal{D}=\mathcal{B}(S)$, which implies that $\mathcal{P}(A)$ holds for all $A \in \mathcal{B}(S)$.\\

\bigskip \noindent \textbf{Exercise 3}.  \label{exercise03_sol_doc04-08} (Borel-Cantelli Lemma) Let $(\Omega,  \mathcal{A}, m)$ be a measure space.\\

\noindent Question (a) Let $(A_n)_{n\geq 0}$ be a sequence of measurable sets on $(\Omega,  \mathcal{A})$ such that

$$
\sum_{n\geq 0} m(A_n) <+\infty.
$$

\bigskip \noindent Show that $\limsup_{n \rightarrow +\infty} A_n$ is a $m$-null set, that is, 

$$
m\left(\limsup_{n \rightarrow +\infty} A_n\right)=0.
$$ 

\bigskip \noindent \textit{Hint}. Recall

$$
\limsup_{n \rightarrow +\infty} A_n=\lim_{n\uparrow +\infty} \downarrow B_n,
$$

\bigskip \noindent with

$$
B_n=\bigcup_{p\geq n} A_p.
$$

\bigskip \noindent Use the $\sigma$-sub-additivity of $m$ in $m(B_n)$ and use a remarkable property of convergent series.\\

\noindent Question (b) Suppose that the measure is finite with $0<m(\Omega)=M<+\infty$. Let $(A_n)_{n\geq 0}$ be a sequence of measurable sets on $(\Omega,  \mathcal{A})$ such that :

$$
\forall \ n\geq 0, \ \lim_{r\rightarrow +\infty} M^r exp(-\sum_{p=n}^{r} m(A_p))=0.  \label{  (SUMINF) }
$$ 

\bigskip \noindent and for any integers $0\leq r <s$,

$$
m\left( \bigcap_{n=r}^{r=s} A_{n}^c\right) \leq  \prod_{n=r}^{r=s} m(A_{n}^c).  \label{  (SUBIND) }
$$

\bigskip \noindent Show that 

$$
m\left(\limsup_{n \rightarrow +\infty} A_n\right)=M.
$$ 

\bigskip \noindent \textbf{Useful remark}. In Probability Theory, the condition ((SUBIND) is implied by an independence condition and the condition (SUMINF) is equivalent to
$$
\sum_{n\geq 0} m(A_n)=+\infty, \ M=1.
$$

\bigskip \noindent \textit{Hint}. Show that $m(\liminf_{n \rightarrow +\infty} A_{n}^{c})=0$. Define

$$
C_n=\bigcap_{n \geq p} A_p^c, n\geq 0
$$

\bigskip \noindent and for $n\geq 0$ fixed and $r \geq n$,

$$
C_{n,r}=\bigcap_{p=n}^{p=r} A_p^c.
$$

\bigskip \noindent Remark that $C_n \nearrow \liminf_{n \rightarrow +\infty} A_{n}^{c}$ and for $n\geq 0$ fixed, $C_{n,r} \nearrow C_n$.

\noindent Combine (SUBIND) with the inequality : $1-x\leq e^{-x}$, $0\leq \leq <1$ to get

$$
m(C_{n,r}) \leq = M^{-n} M^{r} \exp\left( - M^{-1} \sum_{n=r}^{r=s} m(A_{n}\right).
$$

\bigskip \noindent Continue and conclude.\\

\bigskip \noindent \textbf{Solution}.\\

\noindent Question (a). Denote

$$
B_n=\bigcup_{p\geq n} A_p, \ n\geq 0.
$$

\bigskip \noindent We have
$$
\limsup_{n \rightarrow +\infty} A_n=\lim_{n\uparrow +\infty} \downarrow B_n.
$$

\bigskip \noindent Besides, we have $\sigma$-sub-additivity, for $n\geq 0$,
$$
m(B_n) \leq \sum_{p\geq n} m(A_p). \ (BD)
$$

\bigskip \noindent This all the $m(B_n)$'s are finite. The continuity of the measure $m$ leads to

$$
m(\limsup_{n \rightarrow +\infty} A_n)=\lim_{n\rightarrow +\infty} m(B_n). (L)
$$

\bigskip \noindent But $\sum_{p\geq n} m(A_n)$ is the tail of the convergent series $\sum_{n\geq 0} m(A_n)$. Hence $\sum_{p\geq n} m(A_p) \rightarrow 0$ as $n \rightarrow +\infty$. Aplying this to Formula
(BD) gives that $m(B_n)$ also tends to zero. Next, we get through Formula (L) that 
that

$$
m(\limsup_{n \rightarrow +\infty} A_n)=0.
$$

\bigskip \noindent Question (b). It is equivalent to show that $m(\liminf_{n \rightarrow +\infty} A_{n}^{c})=0$. Denote

$$
C_n=\bigcap_{p\geq n} A_p^c, \ n\geq 0.
$$

\bigskip \noindent Since $C_n \uparrow \liminf_{n \rightarrow +\infty} A_{n}^{c}$, we have by the continuity of the measure $m$,

$$
m(\liminf_{n \rightarrow +\infty} A_{n}^{c})=\lim_{n\uparrow +\infty} m(C_n).  \ \ (L)
$$

\noindent But for each $n\geq 0$, we have

$$
C_{n,r}\bigcap_{p=n}^{r} A_{p}^c \downarrow C_n. 
$$

\bigskip \noindent Since $m$ is finite we get by the continuity of the measure $m$ that, for all $n\geq 0$,

$$
m(C_n) =\lim_{r \uparrow +\infty} m(C_{n,r}).
$$

\bigskip \noindent We have by assumption that for $n$ fixed and for $r>n$,

\begin{eqnarray*}
m(C_{n,r}) &\leq& \prod_{n=r}^{r=s} m(A_{n}^c)\\
&=& M^{r-n} \prod_{n=r}^{r=s} m(A_{n}^c)/M\\
&=& M^{r} M^{r} \prod_{n=r}^{r=s} (1-m(A_{n})/M)\\
\end{eqnarray*}

\bigskip \noindent By using the elementary inequality : $1-x \leq exp(-x)$ for $0\leq x \leq 1$ each $m(A_{n})/M$, we get for $n$ fixed and for $r>n$,

$$
m(C_{n,r}) \leq = M^{-n} M^{r} \exp \left( - M^{-1} \sum_{n=r}^{r=s} m(A_{n}\right).
$$

\bigskip \noindent By Condition (SUMINF), $m(C_{n,r})\rightarrow 0$ as $r\rightarrow +\infty$. Hence $m(C_n)=0$ for all $n \geq 0$ and then, by Formula (L), we have

$$
m(\liminf_{n \rightarrow +\infty} A_{n}^{c})=0,
$$

\bigskip \noindent which upon taking the complements, we get

$$
m(\limsup_{n \rightarrow +\infty} A_{n})=M.
$$

\newpage
\noindent \LARGE \textbf{Doc 04-09 : Extensions of Measures - Carath\'{e}odory's Theorem : Application to the existence of the Lebesgue-Stieljes measure }. \label{doc04-09}\\
\bigskip
\Large

\bigskip \noindent  Let $m$ be a measure on the Borel sets of $\mathbb{R}$\ such that for any  any $t\in \mathbb{R}$,\ 
\begin{equation}
F_{m}(t)=m(]-\infty ,t])\text{ is finite.}  \label{1.001a}
\end{equation}
The function $F_{m}:\mathbb{R}\longmapsto $\ $\mathbb{R}$ defined by 
\begin{equation}
t\mapsto F_{m}(t)  \label{1.00ab}
\end{equation}

\noindent has the following properties.

\begin{claim}
$F_{m}$\ assigns non-negative lengths to intervals, that is., 
\begin{equation}
t\leq s\Longrightarrow \Delta F_{m}(]t,s])=F_{m}(s)-F_{m}(t)\geq 0. \label{1.002}
\end{equation}
\end{claim}

\bigskip \noindent Here we have to pay attention to the terminology. We say that $\Delta F_{m}(]t,s]=F_{m}(s)-F_{m}(t)$ is the length of the interval $]t,s]$ by $%
F_{m}.$ It is a pure coincidence that \ref{1.002} means that $F_{m}$ is non-decreasing in this particular case. In higher order spaces, this notion of increasingness will disappear while the notion of positive lengths will be naturally extended to that positive areas and generally to positive volumes.\\

\bigskip \noindent Now the proof of \ref{1.002} is immediately seen by remarking that
\begin{equation*}
(t\leq s)\Rightarrow (]-\infty ,t])\subseteq ]-\infty ,s])
\end{equation*}

\bigskip \noindent and (\ref{1.002}) holds.

\begin{claim}
$F_{m}$\ is right-continuous, that is, 
\begin{equation*}
F_{m}(t^{n})\downarrow F_{m}(t)
\end{equation*}%
as $t^{n}\downarrow t$ 
\begin{equation*}
\Leftrightarrow (\forall (1\leq i\leq
k),t_{i}^{n}\downarrow t_{i}).
\end{equation*}
\end{claim}

\bigskip \noindent  \textbf{Proof}.\\

\bigskip \noindent Let $t^{n}\downarrow t.$ We get\ 
\begin{equation*}
]-\infty ,t^{n}]=\text{ }\downarrow \text{ }]-\infty ,t].
\end{equation*}

\bigskip \noindent Since $m(]-\infty ,t^{n}])$\ is finite for all $n$'s, we get by right-continuity of the measure $m,$ that 
\begin{equation*}
F_{m}(t^{n})=m(]-\infty ,t^{n}])\downarrow m(]-\infty ,t])=F_{m}(t).
\end{equation*}

\bigskip \noindent Suppose now that $m$ is a probability measure, that is $m(\mathbb{R})=1$, we
have two additional points

\begin{claim}
$F_{m}(t)\downarrow 0$ as $t\downarrow -\infty $ \ and $F_{m}(t)\uparrow 1$
as $t\uparrow +\infty .$
\end{claim}

\bigskip \noindent These two points result from the continuity of the measure $m$. First 

\begin{equation*}
]-\infty ,t] \downarrow \emptyset \text{ as } t\downarrow -\infty \text{ and } m(]-\infty ,t])<+\infty \text{ for all } t\in \mathbb{R}
\end{equation*}

\bigskip \noindent implies

\begin{equation*}
m(]-\infty ,t]) \downarrow m(\emptyset)=0 \\\  as \\\ t \downarrow -\infty.
\end{equation*}

\bigskip \noindent Next 
\begin{equation*}
]-\infty ,t]\uparrow \mathbb{R}\text{ as }t\uparrow +\infty 
\end{equation*}

\noindent implies

\begin{equation*}
m(]-\infty ,t])\uparrow m(\mathbb{R})=0\text{ as }t\uparrow -\infty .
\end{equation*}

\bigskip \noindent  In the sequel, we make the following notation : $F(-\infty)=\lim_{t\downarrow -\infty }F(t)$ \ and $F(+\infty )=\lim_{t\uparrow
+\infty }F(t).$ We summarize our study theses two definitions.

\begin{definition}
A non-constant function $F:\mathbb{R}$ $\longmapsto$ $\mathbb{R}$ is a distribution
function if and only if \\

\noindent (1) it assigns non-negative lengths to intervals\\

\noindent and \\

\noindent (2) it is right-continuous at any point $t\in \mathbb{R}$.
\end{definition}

\bigskip \noindent This definition is very broad. $F$ is not necessarily non-negative. It is not required that $F(-\infty )=0$. The second definition is more restrictive

\bigskip

\begin{definition}
A \textbf{non-negative} function $F:\mathbb{R}\longmapsto $\ $\mathbb{R}$ is a probability
distribution function if and only if

\noindent (1) it assigns non-negative lengths to intervals,\\

\noindent (2) it is right-continuous\\

\noindent and\\

\noindent (3) $F(-\infty )=0$ \ and $F(+\infty )=1$.\\
\end{definition}

\bigskip \noindent  We conclude by saying that a measure $m$ which assigns finite values to bounded above intervals generates the distribution functions $F_{m}$. A probability measure generates the probability distribution functions $F_{m}$.\\ 

\bigskip \noindent Inversely, let a distribution function $F$ be given. From the distribution function $F+c$,  it is possible to create a measure on the semi-algebra 
\begin{equation*}
\mathcal{I}_{1}=\{]a,b],a\leq b\},
\end{equation*}

\bigskip \noindent which generates the usual sigma-algebra $\mathcal{B}(\mathbb{R})$ by

\begin{equation}
m_{F}(]a,b])=\Delta F(]a,b])=F(b)-F(a).  \label{def_m1}
\end{equation}

\bigskip \noindent Let us show that $m_{F}$ is additive on $\mathcal{I}_{1}$ . The only possibility to have that an element of is sum of two elements of is to split on element of into two others, that is%
\begin{equation*}
]a,b]=]a,c]+]c,b]\text{, where }c\text{ is finite and }a<c<b,
\end{equation*}

\bigskip \noindent and then

\begin{equation*}
m_{F}(]a,c])+m_{F}(]c,b])=F(c)-F(a)+F(b)-F(c)=F(b)-F(a)=m_(]a,b]).
\end{equation*}

\bigskip \noindent Besides $m_{F}$ is $\sigma$-finite on since we have

$$
\mathbb(R)=\sum_{k\in \mathbb{Z}} ]k,k+1],
$$

\bigskip \noindent and for any $k\in \mathbb{Z}$, we have $m_F(]k,k+1]) <+\infty$.\\

\bigskip \noindent From this, we may apply the advanced version of Caratheodory's Theorem (see Point 4.28 of Doc 04-03) to uniquely extend $m_{F}$ to a measure on $\mathcal{R}$, still denoted by $m_{F}$. We will do this below.\\

\bigskip \noindent The distribution function $F_{m_{F}}$\ generated by \ $m_{F}$ is obtained by letting $a$ go to $-\infty $ in (\ref{def_m1}) that gives%
\begin{equation*}
F_{m_{F}}(b)=m_{F}(]-\infty ,b])=F(b)-F(-\infty),
\end{equation*}

\noindent which is exactly $\ F$ if $F(-\infty )=0$. And, we have the following conclusion.\\

\noindent Let $\mathcal{M}_{0}$ the class of measures $m$ on $\mathbb{R}$, such that
\begin{equation*}
\forall u\in \mathbb{R},m(]-\infty ,u])<+\infty 
\end{equation*}

\noindent and $\mathcal{F}_{0}$ the class of distribution functions $F$ on $\mathbb{R}$ such that 
\begin{equation*}
F(-\infty )=0.
\end{equation*}

\bigskip \noindent There is a one-to-one mapping between $\mathcal{M}_{0}$ and $\mathcal{F}_{0}$
in the following way 
\begin{equation*}
\begin{tabular}{lll}
$F(u)=m(]-\infty ,u])$ & $\longleftrightarrow $ & $m(]a,b])=\Delta F(a,b)$.
\end{tabular}
\end{equation*}

\bigskip \noindent In particular, there is a one-to-one mapping between the $\mathcal{P}$ of probability measures on $\mathbb{R}$ and the class $\mathcal{F}_{P}$ of probability distribution functions $\mathbb{R}$. Probability measures are characterized by their distribution functions.\\

\newpage
\noindent \textbf{Existence of the Lebesgue-Stieljes measures}.\\

\noindent We have to prove the theorem

\begin{theorem} \label{ls}
Let $F$ be a distribution function on $\mathbb{R}$\ then the application $m_{F}$
that is defined on the semi-algebra 
\begin{equation*}
\mathcal{S}_{k}=\{]a,b]=\prod_{i=1}^{k}]a_{i},b_{i}],\text{ }(a,b)\in (\overline{\mathbb{R}})^{2}\}
\end{equation*}

\bigskip \noindent by
 
\begin{equation*}
\left\{ 
\begin{array}{cccc}
m_{F}: & \mathcal{S}_{k} & \mapsto  & \mathbb{R}_{+} \\ 
& ]a,b] & \hookrightarrow  & \Delta _{a,b}F%
\end{array}
\right.
\end{equation*}

\bigskip \noindent is proper, non-negative, additive and $\sigma$-finite on $\mathcal{S}$, and is uniquely extensible to a measure on $\mathcal{B}(\mathbb{R})$.
\end{theorem}

\bigskip \noindent \textbf{Proof of Theorem \ref{ls}}.\\

\noindent We already know $m_{F}$ is proper, non-negative, additive and $\sigma$-finite on $\mathcal{S}$. Thus, by advanced version of Caratheodory's Theorem (see Point 4.28 of Doc 04-03) to uniquely extend $m_{F}$ to a measure on $\mathcal{R}$, it will be enough to show that $m_F$ is $\sigma$-sub-finite on $\mathcal{S}$. Since $m_F$ is already non-negative, additive, we use the characterization of $\sigma$-sub-additivity as established in \textit{Exercise 4 in Doc 04-02}.\\

\noindent We also get that $m_F$ extended to a non non-negative and additive application still denoted by $m_F$ so that $m_F$ is non-decreasing and sub-additive so $\mathcal{S}$ and on $\mathcal{C}$.\\

\noindent So we have to prove that : if the interval $]a,b]$, with $a\in \overline{\mathbb{R}}$, $b\in \overline{\mathbb{R}}$, decomposed like
\begin{equation*}
]a,b]=\sum_{n\geq 1}]a_n,\ b_n].
\end{equation*}

\bigskip \noindent then we have

\begin{equation}
m_{F}(]a,b]) \leq\sum_{n\geq 1}m_{F}(]a_{n},b_{n}]).  \label{sousadditivite1}
\end{equation}

\noindent To prove that we proceed by two steps.\\

\noindent Step 1 : $a$ and $b$ are finite.\\

\noindent if $a=b$, the left-side member is zero and Formula (\ref{sousadditivite1}) is abvious. So, we may suppose that $a<b$. Then we can find $\delta>0$ such that
\begin{equation}
\lbrack a+\delta ,b]\subseteq ]a,b].
\end{equation}

\bigskip \noindent All the numbers $a_{n}$ and $b_{n}$\ have finite components since  we have $a\leq a_{n}\leq b_{n}\leq b$ for all $n\geq 1$.\\

\noindent Since $F$ is right-continuous at $b_{n}$ we have
\begin{equation*}
m_F(]a_{n},b_{n}+\varepsilon ])=F(b^{n}+\varepsilon)-F(a_n)\downarrow F(b^{n})-F(a_n) \ as \ \varepsilon \downarrow 0.
\end{equation*}

\bigskip \noindent So for any $\eta >0$, there exists $\varepsilon_n>0$\ such that  
\begin{equation}
m_F(]a_{n},b_{n}])\leq m_F(]a_{n},b_{n}+\varepsilon^{n}])\leq
m(]a_{n},b_{n}])+\eta /2^{-n}.  \label{lj08}
\end{equation}

\noindent Then, we have
\begin{equation}
\lbrack a+\delta ,b]\subseteq ]a,b]\subseteq \bigcup_{n\geq
1}]a_{n},b_{n}]\subseteq \bigcup_{n\geq 1}]a_{n},b_{n}+\varepsilon_{n}/2[\subseteq \bigcup_{n\geq 1}]a_{n},b_{n}+\varepsilon_{n}]
\label{lj09}
\end{equation}

\bigskip \noindent and next 
\begin{equation}
\lbrack a+\delta ,b]\subseteq \bigcup_{n\geq 1}]a_{n},b_{n}+\varepsilon_{n}/2[.  \label{lj10}
\end{equation}

\bigskip \noindent The compact set $[a+\delta ,b]$\ is included in the union of open sets $]a_{n},b_{n}+\varepsilon_{n}/2[$. By the Heine-Borel property, we may extract from this union a finite union $\bigcup_{1\leq j\leq p}]a_{n_{j}},b_{n_{j}}+\varepsilon_{n_{j}}/2[$ such that 
\begin{equation*}
\lbrack a+\delta ,b]\subseteq \bigcup_{1\leq j\leq
p}]a_{n_{j}},b_{n_{j}}+\varepsilon_{n_{j}}/2[.
\end{equation*}

\bigskip \noindent Now, we get,  as $\delta \downarrow 0,$\ 
\begin{equation*}
]a,b]\subseteq \bigcup_{1\leq j\leq p}]a^{n_{j}},b^{n_{j}}+\varepsilon^{n_{j}}].
\end{equation*}

\bigskip \noindent Since $m_F$ is non-decreasing and sub-additive on $\mathcal{C}$, we have

\begin{eqnarray*}
m_{F}(]a,b])&\leq &\sum_{i=1}^{k}m_{F}(]a_{n_{j}},b_{n_{j}}+\varepsilon_{n_{j}}])\\
&\leq& \sum_{1\leq i\leq p}\left\{ m(]a_{n_{j}},b_{n_{j}}])+\eta 2^{n_{j}}\right\}\\
&\leq& \sum_{1\leq i\leq p}m(]a^{n},b^{n}])+\eta\\
&\leq& \sum_{n\geq 1}m(]a_{n},b_{n}])+\eta .
\end{eqnarray*}

\bigskip \noindent Since this is true for an arbitrary value of $\eta >0,$ we may let $\eta
\downarrow 0$ to have
\begin{equation}
m(]a,b])\leq \sum_{n\geq 1}m(]a_{n},b_{n}]).  \label{lj11}
\end{equation}

\bigskip \noindent \textbf{Step 2}. General where $a\leq b$, $a$ may be $-\infty$ and/or $b$ may be $+\infty$. We proceed as follows. Suppose we have 
\begin{equation*}
]a,b]=\sum_{n\geq 1}]a_n,b_n].
\end{equation*}

\bigskip \noindent Then for any $k>0$, we have

\begin{equation*}
]a,b] \cap ]-k,k]=\sum_{n\geq 1}]a_n,b_n] \cap ]-k,k].
\end{equation*}

\noindent Since $\mathcal{S}$ is stable under finite intersection, the sets $]a,b] \cap ]-k,k]$ and $]a_n, \ b_n] \cap ]-k,k]$'s remain in $\mathcal{S}$. By Step 1 and 
by the non-decreasingness of $m_F$, we have

\begin{eqnarray*}
m_F(]a,b] \cap ]-k,k]) &\leq & \sum_{n\geq 1} m_F(]a_n,b_n] \cap ]-k,k])\\
&\leq & \sum_{n\geq 1} m_F(]a_n,b_n]). \ \ (S1)
\end{eqnarray*}

\bigskip \noindent But it is not difficult to see that $]a,b] \cap ]-k,k])=F(min(b,k))-F(min(a,-k) \uparrow ]a,b])$ and, for $k$ large enough,
$$
m_F(]a,b] \cap ]-k,k])=F(min(b,k))-F(min(a,-k)) \uparrow F(b)-F(a)=m_F(]a,b]).
$$

\noindent To see that, you may consider the three cases : $a=-\infty$ and $b$ finite, $a$ finite, $b=+\infty$, $a=-\infty$ and $b=+\infty$.\\

\noindent By letting $k\uparrow +\infty$ in Formula (S1) above, we get

\begin{eqnarray*}
m_F(]a,b] \sum_{n\geq 1} m_F(]a_n,b_n]).
\end{eqnarray*}

\bigskip \noindent The solution is complete.\\

\newpage
\noindent \LARGE \textbf{Doc 04-10 : Appendix : Proof the Carath\'{e}odory's Theorem}. \label{doc04-10}\\
\bigskip
\Large

\noindent This text gives the full details of the proof of the theorem. The students may skip this for a first reading. It is recommended to read the document with the help of an assistant. This document also introduces to outer measures.\newline

\bigskip \noindent \textbf{Part I : Statement of Carath\'{e}odory's Theorem}.
\newline

\bigskip \noindent Carath\'{e}odory's Theorem allows extension of the measure on the algebra to the $\sigma$-algebra it generates. It is stated as follows.\newline

\bigskip \noindent\textbf{Main form of Caratheodory's Theorem}.\newline

\bigskip \noindent Let $m$ be a non-negative application not constantly equal to $+\infty$ on an algebra $\mathcal{C}$\ of subsets of $\Omega$. If $m$ is $\sigma$-additive, then it is extendable to measure $m^{0}$ on the $\sigma$ -algebra $\mathcal{A}=\sigma(\mathcal{C})$ generated by $\mathcal{C}$. If $m$
is $\sigma$-finite on $\mathcal{C}$, the extension is unique.\newline

\bigskip \noindent We also have other versions.\newline

\bigskip \noindent \textbf{Version 2}. Let $m : \mathcal{C} \rightarrow \mathbb{R}$, be a non-negative application not constantly equal to +$\infty$ such that \newline

\bigskip \noindent (i) $m$ is finitely-additive on $\mathcal{C}$, i.e. 
\begin{equation*}
\forall (A,B)\in \mathcal{C}^{2},\text{ m(A+B)=m(A)+m(B)}
\end{equation*}

\bigskip \noindent and\newline

\bigskip \noindent (ii) $m$ is $\sigma$-sub-additive on $\mathcal{C}$, i.e., if an element $A$ of $\mathcal{C}$ is decomposed as a sum of elements $A_{i}$ \ of $\mathcal{C}$, 
\begin{equation*}
A=\sum_{i=1}^{\infty }A_{i}.
\end{equation*}

\bigskip \noindent Then, we have
\begin{equation*}
m(A)\leq \sum_{i=1}^{\infty }m(A_{i}).
\end{equation*}

\bigskip \noindent Then, $m$ is extendable to a measure on $\sigma(\mathcal{C})$.\newline

\bigskip \noindent \textbf{Version 3}. Let $m : \mathcal{C} \rightarrow \mathbb{R}$, be a non-negative and finite application not constantly equal to +$\infty$ such that \newline

\bigskip \noindent (i) $m$ is finitely-additive on $\mathcal{C}$, i.e. 
\begin{equation*}
\forall (A,B)\in \mathcal{C}^{2},\text{ m(A+B)=m(A)+m(B)}
\end{equation*}

\bigskip \noindent and\newline

\bigskip \noindent (ii) $m$ is continuous at $\emptyset ,$\ i.e., if $(A_{n})_{n\geq 1}$\ is a sequence of subsets of $\mathcal{C}$ non-increasing to $\emptyset$, then $m(A_{n})\downarrow 0$, as $n\uparrow +\infty$.\\

\noindent Then, $m$ extendable to a measure on $\sigma(\mathcal{C})$.\newline

\bigskip \noindent The proof of the theorem will be derived from outer measures. So, we are going to study outer measures before we come back to the proof.

\newpage \bigskip \noindent \textbf{II - Outer measures}.

\bigskip \noindent An application $m^{0}$\ defined on $\mathcal{P}(\Omega)$ is an exterior measure if and only if it is non-negative, non-decreasing, $\sigma-$sub-additive and takes the null value at $\emptyset$, that is \newline

\bigskip \noindent (i) For any sequence (A$_{n})_{n\geq 1}$of subsets of $\Omega$
\begin{equation*}
m^{0}(\bigcup_{n\geq 1}A_{n})\leq \sum_{i=1}^{\infty }m^{0}(A_{i}).
\end{equation*}

\bigskip \noindent (ii) For any subset A of $\Omega $ 
\begin{equation*}
m^{0}(A)\geq 0 .
\end{equation*}

\bigskip \noindent For every $A\subset B$,
\begin{equation*}
m^{0}(A)\leq m^{0}(B). 
\end{equation*}

\begin{equation*}
m^{0}(\emptyset)=0.
\end{equation*}

\bigskip \noindent Now, we are going to construct a $\sigma$-algebra $\mathcal{A}^0$ of subsets of $(\Omega,\mathcal{A})$ such that the restriction of $m^0$ on $\mathcal{A}$ is a measure. To begin with, we define $m^0$-measurability. \newline

\bigskip \noindent A subset $A$ of $\Omega$ is said to be $m^{0}$-measurable if and only if for every subset $D$ of $\Omega$, 
\begin{equation}
m^{0}(D)\geq m^{0}(AD)+m^{0}(A^{c}D).  \label{cara01a}
\end{equation}

\bigskip \noindent This formula is enough for establishing the $m^{0}$-measurability, but we have more. Because of sub-additivity, we have 
\begin{equation*}
m^{0}(D)=m^{0}(AD+A^{c}D)\leq m^{0}(AD)+m^{0}(A^{c}D)
\end{equation*}

\bigskip \noindent so that that for every subset $A$ $m^{0}$-measurable, we have
\begin{equation}
m^{0}(D)=m^{0}(AD)+m^{0}(A^{c}D).  \label{cara01b}
\end{equation}

\bigskip \noindent Now consider the $\mathcal{A}^{0}$ the collection of all $m^{0}$-measurable subsets. This is our candidate. We are going to proof that $\mathcal{A}^{0}$ is a $\sigma$-algebra. But we ill need this intermediate result.\newline

\bigskip \noindent \textbf{(P)} If $A$ and $B$ are two disjoints elements of $\mathcal{A}^{0}$ then $A+B \in \mathcal{A}^{0}$ and $m^{0}(A+B)=m^{0}(A)+m^{0}(B)$.\newline

\noindent By using the fact that $A \in \mathcal{A}^{0}$, we have for any subset $D$ of $\Omega$, we have 

\begin{equation*}
m^{0}(D(A+B))=m^{0}(A(A+B)D)+m^{0}((A+B)A^{c}D). 
\end{equation*}

\bigskip \noindent Since $B \subset A^c$, this gives

\begin{equation*}
m^{0}(D(A+B))=m^{0}(AD)+m^{0}(BD). \text{ } (F0) 
\end{equation*}

\bigskip \noindent We want to prove that $A+B \in \mathcal{A}^{0}$ and this is equivalent to for all $D$ 
\begin{equation*}
m^{0}(D)=m^{0}((A+B)D)+m^{0}((A+B)^{c}D), \text{ } (F1) 
\end{equation*}

\bigskip \noindent that is 
\begin{equation*}
m^{0}(D)=m^{0}((A+B)D)+m^{0}(A{c}B{c}D). \text{ } (F2) 
\end{equation*}

\bigskip \noindent By using (F1), the left hand of (F2) is

\begin{equation*}
m^{0}(AD)+m^{0}(BD)+m^{0}(A^{c}B^{c}D) = m^{0}(AD)+ (
m^{0}(BD)+m^{0}(A^{c}B^{c}D) ) 
\end{equation*}

\bigskip \noindent Since $B \subset A^c$, $BD=A^{c}BD$ and next $m^{0}(BD)+m^{0}(A^{c}B^{c}D)= m^{0}(A^{c}DB)+m^{0}(A^{c}DB^{c})=m^{0}(A^{c}D)$. The left member is 
\begin{equation*}
m^{0}(AD)+m^{0}(A^{c}D)=m^{0}(D). 
\end{equation*}

\bigskip \noindent So $(F2)$ is true. Then $A+B \in \mathcal{A}^{0}$. Now, the definition of $A \in \mathcal{A}^{0}$ with $D=A+B$ to get

\begin{equation*}
m^{0}(A+B)=m^{0}((A+B)A)+m^{0}((A+B)A^c) = m^{0}(A)+m^{0}(B) 
\end{equation*}

\bigskip \noindent Conclusion : $m^{0}$ is additive in $A \in \mathcal{A}^{0}$.\newline

\bigskip \noindent Now, we are going to prove that $\mathcal{A}^{0}$ is a sigma-algebra. It will suffice to prove that is is a Dynkin-system stable
under intersection.\newline

\bigskip \noindent We begin to prove that $\mathcal{A}^{0}$ is a Dynkyn system. We are going to check three points \newline

\bigskip \noindent (i) $\Omega$ belongs to $\mathcal{A}^{0}$ since for all subsets $D$, (F2) reduces to $m^{0}(D) \geq m^{0}(D)$ for $A=\Omega$. We
prove also that $\empty \in \mathcal{A}^{0}$ i the same manner.\newline

\bigskip \noindent (ii) If $A$ and $B$ are elements of $\mathcal{A}^{0}$ such that $A\subset B$, then $B \backslash A \in \mathcal{A}^{0}$. Since $(B
\backslash A)^c=A+B^c$, we have to check that for any subset $D$,

\begin{equation*}
m^{0}(D)=m^{0}(Ba^{c}D)+m^{0}((A+B^c)D). \text{ } (F4) 
\end{equation*}

\bigskip \noindent It is easy to get by the symmetry of the definition of a $m^0$-measurable set $A$ in (F2), that $A$ and $A^c$ measurable at the same time. So here, $B^c$ is $m^0$-measurable. So we may use (F0) in the last term of the left-hand in (F4) to get

\begin{eqnarray*}
m^{0}(D)&=&m^{0}(BA^{c}D)+m^{0}(AD)+m^{0}(B^{c}D)\\
&=&m^{0}(AD)+ \{ m^{0}(A^{c}BD)+m^{0}(B^{c}D) \}. (F5) 
\end{eqnarray*}

\bigskip \noindent Since $A\subset B$, we have $B^c\subset A^c$ and then $B^{c}D=A^{c}B^{c}D$ and next, the expression between curly brackets is such
that

\begin{equation*}
m^{0}(A^{c}BD)+m^{0}(B^{c}D)=m^{0}(A^{c}BD)+m^{0}(A^{c}B^{c}D).
\end{equation*}

\bigskip \noindent By using the definition of $m^0$-measurability of $B$, this gives $m^{0}(A^{c}D)$. Using this (F5) and using again the definition of $m^0$-measurability of $A$ both show that the left-hand is $m^{0}(D)$. Finally (F4) holds, meaning that $B \backslash A \in \mathcal{A}^{0}$.\\

\bigskip \noindent (iii) Let $(A_n)_{n\geq 0} \in \mathcal{A}^{0}$, pairwise disjoint. We have ti check that $A=\sum_{n\geq}A_n \in \mathcal{A}^{0}$.\newline

\bigskip \noindent To see this, put $B_n=\sum_{i=0}^{i=n}A_n$ for $n\geq 0$. Each $B_n$ belongs to $\mathcal{A}^{0}$ by the property (P) and for any $D$, 

\begin{equation*}
m^{0}(D) \geq m^{0}(AB_n)+m^{0}(A_{n}^{c}D).
\end{equation*}

\bigskip \noindent By using Formula (F0) and by remarking that $A^c \subset B_{n}^{c}$ and the fact that $m^0$ is non-decreasing, we  get

\begin{equation*}
m^{0}(D)\geq D)+m^{0}(A^{c}D). 
\end{equation*}

\bigskip \noindent We may let $n$ go to $+\infty$ get

\begin{equation*}
m^{0}(D)\geq \sum_{i\geq 0}m^{0}(A_{i}D)+m^{0}(A^{c}D). (F6) 
\end{equation*}

\bigskip \noindent Now $\sigma$-sub-additivity,

\begin{equation*}
m^{0}(D)\geq m^{0}(\sum_{i\geq 0}(A_{i}D)+m^{0}(A^{c}D) 
\end{equation*}

\bigskip \noindent that is, for any $A \in \mathcal{A}^{0}$, we have

\begin{equation*}
m^{0}(D)\geq m^{0}(AD)+m^{0}(A^{c}D).
\end{equation*}

\bigskip \noindent By taking taking $D=A$ in (F5) leads to

\begin{equation*}
m^{0}(A)\geq m^{0}(A). \ \square
\end{equation*}

\bigskip \noindent This and the $\sigma$-sub-additivity lead to

\begin{equation*}
m^{0}(\sum_{i=0}^{n}m^{0}(A_{i})\geq \sum_{i\geq 0}m^{0}(A_{i}D).
\end{equation*}

\bigskip \noindent We get that $\mathcal{A}^{0}$ is a Dynkin-system and $m^0$ is $\sigma$-additive on $\mathcal{A}^{0}$.\\

\bigskip \noindent (iv) We have to check that $\mathcal{A}^{0}$ is stable under finite intersection. \newline

\bigskip \noindent To prove that, assume that $A$ and $B$ are elements of $\mathcal{A}^{0}$ and consider a subset $D$. By definition of $A \in \mathcal{A}^{0}$, we have 
\begin{equation*}
m^{0}(D)=m^{0}(AD)+m^{0}(A^{c}D).
\end{equation*}

\bigskip \noindent We expand each member of left-hand by using the definition of $B \in \mathcal{A}^{0}$ to have 
\begin{equation*}
m^{0}(D)=\left\{
m^{0}(ADB)+m^{0}(ADB^{c})\right\}+\{m^{0}(A^{c}DB)+m^{0}(A^{c}DB^{c})\}. 
\end{equation*}

\bigskip \noindent Now, use sub-additivy for the three last terms in the left-hand to get

\begin{equation*}
m^{0}(D) \geq m^{0}(ADB)+m^{0}( \{ (AB^{c})\cup (A^{c}B)\cup (A^{c}B^{c})\}D). 
\end{equation*}

\bigskip \noindent We have that $(AB^{c})\cup (A^{c}B)\cup(A^{c}B^{c})=(AB)^c$. How to see this : you say that the elements $(AB)^c$ are the elements of $\Omega$ which do not belongs to both of $A$ and $B$.
You split this three disjoints subsets : the subset of elements of $\Omega$ belonging to $A$ and not belonging to $B$, the subset of elements of $\Omega$ belonging to $B$ and not belonging to $A$, and the subset of elements of $\Omega$ not belonging to $A$ and not belonging to $B$.\newline

\bigskip \noindent Final conclusion : $\mathcal{A}^{0}$ is a Dynkin-system stable under intersection. The it is a sigma-algebra. Finally $(\Omega,\mathcal{A}^{0}, m^0)$ is measure space.\\

\bigskip \noindent We are now going to apply this to prove the theorem of Caratheodory.

\newpage

\begin{center}
\textbf{Proof of theorem of Caratheodory.}
\end{center}

\bigskip \noindent We will now define an exterior measure from $m$ satisfying the hypotheses of the theorem.

\bigskip \noindent \textbf{Definition of an exterior measure}.\newline

\bigskip \noindent Let $m$ be an application defined on the algebra $\mathcal{C}$ satisfying properties (i) and (ii). We define for a subset $A$ of $\Omega$, 
\begin{equation*}
m^{0}(A)=\inf \{\sum_{i=1}^{\infty }m(A_{i}),\text{ where (A}_{i})_{n\geq
1}\subseteq \mathcal{C},\text{ A}\subseteq \bigcup_{n\geq 1}A_{n}\}.
\end{equation*}

\noindent S This means that : $m^{0}(A)$\ is the infimum of the numbers $\sum_{i=1}^{\infty }m(A_{i}),$\ for every cover of $A$ by a countable family $\bigcup_{n\geq 1}A_{n}$ of $\mathcal{C}$. This number always exists since $\Omega \in \mathcal{C}$, and $\{\Omega \}$ is a cover for every subset $A$. We are going to show that $m^{0}$\ is an exterior measure through steps (I), ..., (IV).\newline

\bigskip \noindent (I) m$^{0}$ is non-negative and non(decreasing.\newline

\bigskip \noindent The non-negativity is obvious. So also is the non-decreasingness since if A$\subset B$\ and $\bigcup_{n\geq 1}A_{n} $ is a covering $B$ by elements of $\mathcal{C}$, it will also be a covering in  of $A$ by elements of $\mathcal{C}$ and by, by definition, 
\begin{equation*}
m^{0}(A)\leq \sum_{n\geq 1}m(A_{i})
\end{equation*}

\bigskip \noindent Taking the infinum over the coverings of $B$ by elements of $\mathcal{C}$, we have 
\begin{equation*}
m^{0}(A)\leq m^{0}(B).
\end{equation*}

\bigskip \noindent (II) We show that if $A\in \mathcal{C},$\ then $m^{0}(A)=m(A).$\ In fact if $A\in \mathcal{A}$, $\{A\}$ is a covering of itself and belongs to $\mathcal{C}$. So by definition, 
\begin{equation*}
m^{0}(A)\leq m(A).
\end{equation*}

\bigskip \noindent Conversely, let $A\in \mathcal{C}$, $A_{n})_{n\geq 1}$\ $\subseteq \mathcal{C}$, a covering of $A$. We have 
\begin{equation*}
\text{A}\subseteq \bigcup_{n\geq 1}A_{n}\Longrightarrow A=A\cap
(\bigcup_{n\geq 1}A_{n})=\bigcup_{n\geq 1}A\cap A_{n}.
\end{equation*}

\bigskip \noindent From hypothesis (iv) 
\begin{equation*}
m(A)\leq \sum_{i=1}^{\infty }m(A\cap A_{i})
\end{equation*}

\bigskip \noindent and by non-decreasingness of $m$ on $\mathcal{C}$, 
\begin{equation*}
m(A)\leq \sum_{i=1}^{\infty }m(A\cap A_{i})\leq \sum_{i=1}^{\infty }m(A_{i}).
\end{equation*}

\bigskip \noindent Taking the infimum on the left-hand side on all the coverings of $A$ by elements of $\mathcal{C}$, we get
\begin{equation*}
m(A)\leq m^{0}(A).
\end{equation*}

\bigskip \noindent So we have 
\begin{equation*}
\forall (A\in \mathcal{C}{\large ),m}^{0}(A)=m(A).
\end{equation*}

\bigskip \noindent (III) m($\emptyset )=m(\emptyset )=0.$

\bigskip \noindent (IV) We show that m$^{0}$\ is $\sigma $-sub-additive$.$
Let A=$\bigcup_{n\geq 1}A_{n},$\ we show that 
\begin{equation}
m^{0}(A)\leq \sum_{i=1}^{\infty }m^{0}(A_{i})  \label{cara02}
\end{equation}

\bigskip \noindent We remark that if $m^{0}(A_{i})=+\infty$, then formula (\ref{cara02}) is true. Now, suppose that all the m$^{0}(A_{i})$ are finite. Let $\varepsilon >0$. By the characterization of a finite infimum there exists a covering $A_{i}$\ by elements of $\mathcal{C}$, (A$_{i,n})_{n\geq 1}\subseteq \mathcal{C},$\ such that 
\begin{equation*}
\sum_{n=1}^{\infty }m(A_{i,n})\leq m^{0}(A_{i})+\varepsilon 2^{-i}.
\end{equation*}

\bigskip \noindent But 
\begin{equation*}
\bigcup_{i\geq 1}A_{i}\subseteq \bigcup_{i\geq 1,n\geq 1}A_{i,n}.
\end{equation*}

\bigskip \noindent Hence the elements of the left form a covering by the elements
of $\mathcal{C}$. So we get
\begin{equation*}
m^{0}(\bigcup_{i\geq 1}A_{i})\leq \sum_{i\geq 1}\sum_{ni\geq 1}m(A_{i,n}).
\end{equation*}

\bigskip \noindent Now we get 
\begin{equation*}
m^{0}(\bigcup_{i\geq 1}A_{i})\leq \sum_{i\geq 1}\sum_{ni\geq
1}m(A_{i,n})\leq \sum_{i\geq 1}(m^{0}(A)_{i}+\varepsilon 2^{-i})\leq
\sum_{i\geq 1}m^{0}(A_{i})+\varepsilon
\end{equation*}

\bigskip \noindent for every $\varepsilon >0.$\ So for $\varepsilon
\downarrow 0,$\ we have 
\begin{equation*}
m^{0}(\bigcup_{i\geq 1}A_{i})\leq \sum_{i\geq 1}m^{0}(A_{i}).
\end{equation*}

\bigskip \noindent Hence the $\sigma -$sub-additivity of $m^{0}$ is proved.\\

\bigskip \noindent Finally, we show that $\mathcal{C} \subset \mathcal{A}^{0}$. Let $A\in \mathcal{C}$, we are going to show that for all $D$ 
\begin{equation*}
m^{0}(D)\geq m^{0}(AD)+m^{0}(A^{c}D).
\end{equation*}

\bigskip \noindent For $m^{0}(D)=\infty$, there is nothing to do. Let $m^{0}(D)<\infty$. For all $\varepsilon >0$, we could find a covering $\bigcup_{i\geq 1}A_{i}$\ in $\mathcal{C}$\ of D, such that 
\begin{equation*}
m^{0}(D)+\varepsilon \geq \sum_{i\geq 1}m(A_{i}).
\end{equation*}

\bigskip \noindent But, $m$ is additive on $\mathcal{C}$. So 
\begin{equation*}
m^{0}(D)+\varepsilon \geq \sum_{i\geq 1}(m(AA_{i})+m(A^{c}A_{i})).
\end{equation*}

\bigskip \noindent But the $AA_{i}^{\prime}$'s  are in $\mathcal{C}$ and cover AD and the A$^{c}A_{i}^{\prime }s$\ are in $\mathcal{C}$ and cover $A^{c}D.$\ So by definition of $m^{0}$, we have
\begin{equation*}
\sum_{i\geq 1}m(AA_{i})\geq m^{0}(AD)
\end{equation*}

\bigskip \noindent and 
\begin{equation*}
\sum_{i\geq 1}m(A^{c}A_{i})\geq m^{0}(A^{c}D).
\end{equation*}

\bigskip \noindent By combining the last three formulas, we have for all $\varepsilon >0$, \begin{equation*}
m^{0}(D)+\varepsilon \geq m^{0}(AD)+m^{0}(A^{c}D)
\end{equation*}

\bigskip \noindent Hence, taking $\varepsilon \downarrow 0$, yields $A\in ^{0}$. We conclude by saying that $\mathcal{C}\subset \mathcal{A}^{0}$, which is a $\sigma -$algebra. So, we have

\begin{equation*}
\sigma (\mathcal{C})\subset \mathcal{A}^{0}.
\end{equation*}

\bigskip \noindent We conclude by saying that $m^{0}$\ is a measure on $\mathcal{A}^{0},$\ containing $\sigma(\mathcal{C})$\ and is equal to m on $\mathcal{C}$. So $m^{0}$\ is a measure measure, extension of m on $\sigma(\mathcal{C}).$\ if m is $\sigma -$finite, the extension is $\sigma$-finite as well. And since two $\sigma -$finite measures on an algebra are equal on the sigma-algebra it generates, we see that the extension is unique, once $m$ is $\sigma -$finite.\\

\newpage
\noindent \LARGE \textbf{Doc 04-11 Application of the exterior measure to Lusin's theorem}. \label{doc04-11}\\
\bigskip
\Large

\bigskip

\noindent \textbf{Introduction} By applying the proof of the Caratheodory measure for the Lebesgues measure on $\mathbb{R}$, we get this interesting property which approximates any Borel set with finite measure to a sum on intervals.\\

\noindent \textbf{Lusin's theorem}.\\

\noindent Let be a Borel subset $A$ of $\mathbb{R}$ such that its Lebesgues measure is finite, that is $\lambda (E)<+\infty$. Then, for any $\varepsilon >0,$ there exists a finite union $K$ of bounded intervals such that 
$$
\lambda (E\Delta K)<\varepsilon.
$$

\bigskip \noindent \textbf{Solution}.  We have to go back to the proof of the Theorem of Caratheodory which
justified the existence which the Lebesgue measure.\\

\noindent  Recall that the Lebesgue measure is uniquely defined on the class of intervals
\begin{equation*}
\mathcal{I}=\{]a,b],-\infty \leq a\leq b<+\infty \}
\end{equation*}

\bigskip \noindent by
\begin{equation*}
\lambda (]a,b])=b-a.
\end{equation*}

\bigskip \noindent This class $\mathcal{I}$ is a semi-algebra and the algebra $\mathcal{C}$ its generates is the class of finite sums of disjoints intervals. It is
easily proved that $\lambda $ is additive on $\mathcal{I}$ and is readily extensible to a an additive application on $\mathcal{C}$, that we always denote by $\lambda$. The broad extension of $\lambda $ to a measure on a $\sigma$-algebra $\mathcal{A}^{0}$ including $\mathcal{C}$ may be done by the method of the outer measure, defined as follows, for any subset
 $A$ of $\mathbb{R}$ :

\begin{equation}
\lambda ^{0}(A)=\inf \{\sum_{n=0}^{\infty }\lambda (A_{n}),\text{ }A\subset
\bigcup\limits_{n=0}^{\infty }A_{n},A_{n}\in \mathcal{C}\}.  \label{funct.topo.002}
\end{equation}

\bigskip \noindent A subset of $\mathbb{R}$ is $\lambda ^{0}$-measurable if and only for any subset$D$ of $\mathbb{R}$, we have
\begin{equation*}
\lambda ^{0}(A)=\lambda ^{0}(AD)+\lambda ^{0}(AD^{c}).
\end{equation*}

\bigskip \noindent By denoting $\mathcal{A}^{0}$ \ the set of \ $\lambda ^{0}$-measurable subsets of $\mathbb{R}$, is proved that $(\Omega ,\mathcal{A}^{0},\lambda ^{0})$ is a measurable space and $\mathcal{C} \subset \mathcal{A}^{0}.$ This measurable space, surely, includes the measurable space $(\Omega ,\mathcal{B}(\mathbb{R}),\lambda ^{0})$ since $\mathcal{B}(R)=\sigma (\mathcal{C})\subset \mathcal{A}^{0}.$ The measure $\lambda ^{0}$ is the unique extension of $\lambda $ to a measure on $\mathcal{B}(\mathbb{R})$, still denoted by $\lambda$.\\

\noindent By using \ref{funct.topo.002} and the characterization of the infinum based on the fact that $\lambda (E)$ is finite, we conclude that for any $\varepsilon>0,$ there exists a union  $\bigcup_{n\geq 1} A_{n}$, formed by elements $A_n$ in $\mathcal{C}$ and covering $E$, such that we have 

\begin{equation}
\sum_{n=0}^{\infty }\lambda (A_{n})<\lambda (E)+\varepsilon /2.  \label{funct.topo.003}
\end{equation}

\noindent But for each $n\geq 1$,  $A_{n}$ is a finite sum of intervals - elements of $\mathcal{I}$ -  of the form 

$$
A_{n}=\sum_{1\leq j\leq p(n)}I_{n,j},
$$ 

\bigskip \noindent with 

$$
\lambda (A_{n})=\sum_{1\leq j\leq p(n)}\lambda (I_{n,j}),
$$

\bigskip \noindent since $\lambda $ is additive on $\mathcal{C}$. We may rephrase this by saying that  for any $\varepsilon >0,$ there exists a union of intervals $I_{k}$, covering $E$, such that 
\begin{equation*}
\sum_{k=1}^{\infty }\lambda (I_{k})<\lambda (E)+\varepsilon /2.
\end{equation*}

\bigskip \noindent Here, we see that none of these intervals is unbounded. Otherwise, for one them, we would have  $\lambda (I_{k})=+\infty $ and Formula (\ref{funct.topo.003}) would be
impossible. Since\ $\sum_{k=1}^{\infty }\lambda (I_{k})$ is finite, we may find $k_{0}$ such that
\begin{equation}
0\leq \ \sum_{k=1}^{\infty }\lambda (I_{k})-\ \sum_{k=1}^{k_{0}}\lambda
(I_{k})=\sum_{k\geq k_{0}+1}\lambda (I_{k})<\varepsilon /2.  \label{8.004}
\end{equation}

\noindent Let us put $K=\cup _{1\leq j\leq k_{0}} I_j$. We finally have
\begin{equation*}
K\backslash E=KE^{c}\subset \left( \bigcup\limits_{k=0}^{\infty
}I_{k}\right) E^{c}=\left( \bigcup\limits_{k=0}^{\infty }I_{k}\right)
\backslash E,
\end{equation*}

\bigskip \noindent with
\begin{equation*}
E\subset \left( \bigcup\limits_{k=0}^{\infty }I_{k}\right) =\left(\bigcup\limits_{n=0}^{\infty }A_{n}\right).
\end{equation*}

\bigskip \noindent Then, by \ref{funct.topo.002}, we have
\begin{equation*}
\lambda (K\backslash E)\leq \lambda \left( \bigcup\limits_{k=0}^{\infty
}I_{k}\right) -\lambda (E)\leq \left( \sum_{k=0}^{\infty }\lambda
(A_{k})\right) -\lambda (E)<\varepsilon /2.
\end{equation*}

\bigskip \noindent Next, by denoting $J=\bigcup\limits_{k=0}^{\infty }I_{k},$ we have $E\subset J$ and 
\begin{equation*}
K^{c}=\left( \bigcup\limits_{k=0}^{\infty }I_{k}\backslash
\bigcup\limits_{k=0}^{k_{0}}I_{k}\right) +J^{c},
\end{equation*}

\bigskip \noindent and then
\begin{equation*}
E\backslash K=EK^{c}=E\left( \bigcup\limits_{k=0}^{\infty }I_{k}\backslash
\bigcup\limits_{k=0}^{k_{0}}I_{k}\right) +EJ^{c}=E\left(
\bigcup\limits_{k=0}^{\infty }I_{k}\backslash
\bigcup\limits_{k=0}^{k_{0}}I_{k}\right) \subset
\bigcup\limits_{k_{0}+1}^{\infty }I_{k}
\end{equation*}

\bigskip \noindent and
\begin{equation*}
\lambda (E\backslash K)\leq \lambda \left( \bigcup\limits_{k_{0}+1}^{\infty
}I_{k}\right) \leq \sum_{k\geq k_{0}+1}\lambda (I_{k})<\varepsilon /2.
\end{equation*}

\bigskip \noindent We conclude that
\begin{equation*}
\lambda (E\Delta K)<\varepsilon .
\end{equation*}

%\chapter{Integration}
\chapter{Integration} \label{05_integration}

\begin{table}[htbp]
	\centering
		\begin{tabular}{llll}
		\hline
		Type& Name & Title  & page\\
		\hline
		S & Doc 05-01 & Integration with respect to a measure - A summary   & \pageref{doc05-01}\\
		S & Doc 05-02 & Integrals of elementary functions - Exercises & \pageref{doc05-02}\\
		D & Doc 05-03 & Integration with the counting measure - Exercises  & \pageref{doc05-03}\\
		D & Doc 05-04 & Lebesgue/Riemann-Stieljes integrals on $\mathbb{R}$  - Exercises & \pageref{doc05-04}\\
		SD& Doc 05-05 & Integrals of elementary functions - Exercises with solutions & \pageref{doc05-05}\\
		SD& Doc 05-06 & Integration with the counting measure - Solutions & \pageref{doc05-06}\\
		SD & Doc 05-07 & Lebesgue/Riemann-Stieljes integrals on $\mathbb{R}$ - Solutions  & \pageref{doc05-07}\\
		A& Doc 05-08 &  Technical document on the definition of the integral & \\
		&           &  of a function of constant sign & \pageref{doc05-08}\\
		\hline
		\end{tabular}
\end{table}

\newpage
\noindent \LARGE \textbf{Doc 05-01 : Integration with respect to a measure - A summary}. \label{doc05-01}\\
\Large

\bigskip \noindent Summary : The integration with respect to a measure is the most general theory of integration. It generalizes the Riemman integral on compact intervals on $\mathbb{R}$. This measure in an abstract one, more oriented to theoretical purposes.\\

\noindent \textbf{I - Construction of the integral and Main properties}.\\

\bigskip \noindent The main point is that this measure is constructed. This means the the construction process should be mastered and never be forgotten.\\

\bigskip \noindent The construction is achieved in three points. This leads to an easy handling of integration.\\

\bigskip \noindent \textbf{05-01 Definition of the integral with respect to a measure}. Assume that we have a measure space $(\Omega, \mathcal{A},m)$. We are going to construct the integral of a real-value measurable  function $f : (\Omega, \mathcal{A}) \rightarrow \overline{\mathbb{R}}$ with respect to the measure $m$ (that may take infinite values) denoted by
$$
\int f \text{ }dm=\int_{\Omega} f(\omega) \text{ }dm(\omega) = \int_{\Omega} f(\omega) \text{ } m(d\omega).
$$

\bigskip  \noindent \textbf{05-01a Definition of the integral of a non-negative simple function : $f \in \mathcal{E}_{+}$}.\\  

\noindent The integral of a non-negative simple function

$$
f=\sum_{1\leq i\leq k}\alpha _{i}\text{ }1_{A_{i}}, \ (\alpha_i \in \mathbb{R}_{+}, \ A_i \in \mathcal{A}, \ 1\leq i \leq k), \ k\geq 1,
$$

\bigskip  \noindent is defined by

\begin{equation*}
\int f\text{ }dm=\sum_{1\leq i\leq k} \alpha _{i} \ m(A_{i}). \ \ (IPEF)
\end{equation*}

\bigskip \noindent \textbf{Convention - Warning 1} In the definition (IPEF), the product $\alpha _{i} \ m(A_{i})$ is zero whenever $\alpha_i=0$, event if $m(A_i)=+\infty$.\\   

\noindent The definition (IPEF) is coherent. This means that $\int f \ dm$ does not depend on one particular expression of $f$.\\

\bigskip \noindent \textbf{05-01b Definition of the integral for a non-negative measurable function}.\\

\noindent Let $f$ be any non-negative measurable function. By \textit{Point (03-23) in Doc 03-02 in Chapter \ref{03_setsmes_applimes_cas_speciaux}, page \pageref{doc03-02}}, there exists a non-decreasing sequence $(f_n)_{n\geq 0}\subset \mathcal{E}_{+}$ such that
$$
f_n \uparrow f \text{as n } \uparrow +\infty.
$$

\bigskip \noindent  We define
 
$$
\int f \text{ }dm =\lim_{n\uparrow +\infty} \int f_n \text{ }dm  \text{ as n } \uparrow +\infty. \ (IPF)
$$

\bigskip \noindent  This definition (IPF) is also coherent since it does not depends of the sequence which is used in the definition.

\bigskip \noindent \textbf{05-01c Definition of the integral for a measurable function}.\\

\noindent In the general case, let us use the decomposition of $f$ into its positive part and its negative part :

$$
f=f^{+}-f^{-} \, |f|=f^{+}+f^{-}, \ and \ f^{+}f^{-}=0, \ (PPPN)
$$

\bigskip \noindent where $f^{+}=max(0,f)$ and $f^{-}=max(0,-f)$, which are measurable, form the unique couple of functions such that Formulas (PPPN) holds.\\

\noindent By (05-1b), the numbers 

$$
\int f^+ \ dm \ \ and \int f^- \ dm
$$

\bigskip \noindent exist in $\overline{\mathbb{R}}_{+}$. If one of them is finite, i.e.,

$$
\int f^+ \ dm <+\infty \ \ \textbf{or} \int f^- \ dm <+\infty,
$$

\bigskip \noindent we say that $f$ is quasi-integrable with respect to $m$ and we define

$$
\int f dm = \int f^+ \ dm \ - \ \int f^- \ dm.
$$

\bigskip \noindent \textbf{Warning} The integral of a real-valued and measurable function with respect to a measure $m$ exists \textit{if only if} : either it is of constant sign or the integral of its positive part or its negative part is finite.\\

\bigskip \noindent \textbf{05-02 Integration over a measurable set}.\\

\noindent If $A$ is a measurable subset of $\Omega$ and $1_A f$ is quasi-integrable, we denote

$$
\int_A f \ dm = \int 1_A f \ dm.
$$ 

\bigskip \noindent \textbf{Convention - Warning 2} The function $1_A f$ is defined by $1_A f(\omega)=f(\omega)$ of $\omega \in A$, and $zero$ otherwise.\\

\bigskip \noindent \textbf{05-03 Integrable functions}.\\

\noindent The function $f$ is said to be integrable if and only if the integral $\int f dm$ exist (in  $\mathcal{R}$) and is finite, i.e.,

$$
\int f^+ \ dm <+\infty \ \ \textbf{and} \int f^- \ dm <+\infty,
$$
  
\bigskip \noindent The set of all integrable functions with respect to $m$ is denoted

$$
\mathcal{L}^{1}(\Omega, \mathcal{A},m).
$$

\newpage
\noindent \textbf{II - Main properties}.\\

\noindent Let $f$ and $g$ be measurable functions, such that   $\int f \ dm$, $\int g \ dm$, $f+g$ is defined \textit{a.e.}, 
$\int f+g \ dm$, exists and $\int f \ dm + \int g \ dm$ make senses, $A$ and $B$ be two disjoint measurable sets, and $c$ be a finite non-zero real scalar. We have the following properties and facts.\\

\noindent \textbf{05-04 - Linearity}.\\

$$
\int cf \ dm = c \int f \ dm, \ \ (L1)
$$

\bigskip \noindent which is still valid for $c=0$ if $f$ is integrable,

$$
\int f+g \ dm=\int f \ dm + \int g \ dm, \ \ (L2)
$$

\bigskip \noindent and

$$
\int_{A+B} f \text{ }dm=\int_{A} f \text{ }dm+\int_{B} f \text{ }dm, \text{  } \int cf \text{ }dm=c \int f \text{ }dm, \ \ (L3)
$$

\bigskip \noindent \textbf{05-05 - Order preservation}.\\

$$
f \leq g  \Rightarrow \int f \  dm \leq \int g \ dm, \ \ (O1)
$$

$$
f=g \  a.e. \  \Rightarrow \int f \ dm=\int g \text{ }dm, \ \ (O2)
$$

\bigskip \noindent \textbf{05-06 - Integrability}.

$$
f  \ integrable  \Leftrightarrow |f| \ integrable  \Rightarrow f \ finite \ a.e, \ \ (I1)
$$

$$
|f| \leq g \ integrable \Rightarrow f integrable, \ \ (I2)
$$

\bigskip \noindent and

$$
\biggr((a,b)\in \mathcal{R}^2, \  f \ and \ g \ integrable \biggr) \Rightarrow \biggr( af+bg \ integrable \biggr). \ \ (I3)
$$

\newpage
\noindent \textbf{III - Spaces of integrable functions}.\\

\noindent \textbf{05-07a Space $\mathcal{L}^{p}(\Omega, \mathcal{A},m)$, $p\geq 1$}.\\

\noindent $\mathcal{L}^{p}(\Omega, \mathcal{A},m)$ is the space of all real-valued and measurable functions $f$ defined on $\Omega$ such that $|f|^p$ is integrable with respect to $m$, that is

$$
\mathcal{L}^{p}(\Omega, \mathcal{A},m)=\{f \in \mathcal{L}_0(\Omega, \mathcal{A}) \int |f|^p \ dm  < +\infty \}.
$$

\bigskip \noindent Recall : $\mathcal{L}_0(\Omega, \mathcal{A})$ is the class of all real-valued and measurable functions defined on $\Omega$ as in \textit{Point 4.19a in Doc 04-01 of Chapter \ref{04_measures}}, page \pageref{04_equivcal}.\\

\bigskip \noindent \textbf{05-07b Space of $m$-a.e. bounded functions $\mathcal{L}^{\infty}(\Omega, \mathcal{A},m)$}.\\

\noindent $\mathcal{L}^{\infty}(\Omega, \mathcal{A},m)$ is the class of all real-valued and measurable functions which are bounded $m$-a.e., i.e,

$$
\mathcal{L}^{\infty}(\Omega, \mathcal{A},m)=\{f \in \mathcal{L}_0(\Omega, \mathcal{A}) 0, \ \exists C\in \mathbb{R}_{+}, \ |f|\leq C \ m-a.e \}.
$$

\bigskip \noindent \textbf{05-07c Equivalence classes}.\\

\bigskip \noindent When we consider only the equivalence classes pertaining to the equivalence relation introduced in \textit{Point 4.19b in Doc 04-01 in Chapter \ref{04_measures}}, page \pageref{doc04-01}, we get the quotient spaces

$$
L^{p}(\Omega, \mathcal{A},m)=\mathcal{L}^{p}(\Omega, \mathcal{A},m)/ \mathcal{R}.
$$

\bigskip \noindent that is,

$$
L^{p}(\Omega, \mathcal{A},m)=\{ \dot{f}, f \in  \mathcal{L}_0(\Omega, \mathcal{A}), \ \int |f|^p \ dm \text{ } < +\infty \}
$$

\bigskip \noindent and

$$
L^{\infty}(\Omega, \mathcal{A},m)=\mathcal{L}^{+\infty}(\Omega, \mathcal{A},m)/ \mathcal{R}
$$

\bigskip \noindent that is 

$$
L^{\infty}(\Omega, \mathcal{A},m)=\{ \dot{f}, f \mathcal{L}_0(\Omega, \mathcal{A}),\  |f| \ m-bounded \}.
$$

\bigskip \noindent \textbf{Measure Space Reduction Property} \label{msrpA}. We may complete the first statement this principle (See page \pageref{msrp} in Chapter \ref{04_measures}) by saying : In situations where we have a a countable family of \textit{a.e.} true assertions $(\mathcal{P}_n, \ n\geq 0)$, we may and usually assume the they are all true everywhere, but putting ourselves on the induced measure space $(\Omega_0, \mathcal{A}_{\Omega_0}, m_{\Omega_0})$. This change does not affect neither the measure values of measurable sets nor the values of integrals of real-valued functions whenever they exist (See Solutions of Exercise 10 of Doc 05-02 in Doc 05-05, page \pageref{exercise10_sol_doc05-05}).

\bigskip \noindent \textbf{05-08 Norms on $L^{p}$}.\\

\noindent The space $L^{p}(\Omega, \mathcal{A},m)$, $1\leq p <+\infty$, is equipped with the norm

$$
||f||_{p} = ( \int |f|^p dm )^{1/p}.
$$

\bigskip \noindent The space $L^{+\infty}(\Omega, \mathcal{A},m)$ is equipped with the norm

$$
||f||_{\infty} = \inf \{C \in \mathbb{R}_{+}, \ |f|\leq C \ m-.a.e  \}.
$$

\bigskip \noindent \textbf{05-9 Banach spaces}.\\

\noindent  For all $p \in [1, +\infty]$, $L^{p}(\Omega,\mathcal{A}, m) (+, ., ||f||_{p})$ is a Banach space.\\

\noindent $L^{2}(\Omega,\mathcal{A}, m) (+, ., ||f||_{2})$ is a Hilbert space.\\

\newpage
\noindent \LARGE \textbf{Doc 05-02 : Integrals of elementary functions - Exercises}. \label{doc05-02}\\
\bigskip
\Large

\noindent \textbf{NB} : \textbf{In this work, you will have to use the \textit{demonstration tools} of Part 04.12 (see \textbf{04.12a} and \textbf{04.12b}) of Doc 12.}\\

\noindent Let $(\Omega, \mathcal{A}, m)$ be measure space. Let $\mathcal{E}^{(1)}$ be the class of non-negative simple functions or integrable simple functions.\\

\noindent \textbf{Important Remarks}. The properties which are asked to be proved here are intermediate ones. They are used to prove the final properties which are exposed in 
\textit{Part II in Doc 05-01 Main properties}.

\bigskip  \noindent \textbf{Part I - Properties of the integral on $\mathcal{E}^{(1)}$}.\\

\bigskip  \noindent Let be given a non-negative simple function  
\begin{equation*}
f=\sum_{1\leq i\leq k}\alpha _{i}1_{A_{i}}, \ (\alpha_i \in \mathbb{R}_{+}, \ A_i \in \mathcal{A}, \ 1\leq i \leq k), \ k\geq 1. \ \ (NSF)
\end{equation*}

\bigskip  \noindent define its integral with respect to $m$ by :

\begin{equation}
\int f \ dm=\sum_{1\leq i\leq k}\alpha _{i} \ m(A_{i}). \ \ (IFS)
\end{equation}

\bigskip  \noindent \textbf{Exercise 1}. \label{exercise01_doc05-02}  (Coherence) Prove that this definition is coherent. In order to do this, consider two expressions of $f$, use the expression based on the partition of $\Omega$ formed by the superposition of the two first partitions as in \textit{Point 03.22a of Doc 03-02 of Chapter \ref{03_setsmes_applimes_cas_speciaux}}, page \pageref{doc03-02}.

\bigskip  \noindent \textbf{Exercise 2}. \label{exercise02_doc05-02} Show that Formula (IFS) holds if $f$ is non-positive and if $f$ is integrable. Remark that the elementary function $f$ integrable if and only if, we have 

$$
\forall i \in \{1,...,k\}, \ (\alpha_i\neq 0) \Rightarrow m(A_i) <+\infty.
$$

\bigskip  \noindent \textbf{Exercise 3}. \label{exercise03_doc05-02}  Establish the following properties on $\mathcal{E}{+}$.\\

\noindent Hereafter, $\alpha$ and $\beta$ are two real numbers,  $f$ and $g$ are two elements of $\mathcal{E}^{+}$, $A$ and $B$ are measurable spaces.\\

\bigskip  \noindent \textbf{P1} : The integral is non-negative, that is :

$$
f\geq 0 \Rightarrow \int f \ dm\geq 0.
$$

\bigskip  \noindent \textbf{P2} : We have  

\begin{equation*}
f=0 \ a.e. \Leftrightarrow \int f \ dm=0.
\end{equation*}

\bigskip  \noindent \textbf{P3} If $A$ is a null-set, then 
\begin{equation*}
\int_{A}f\text{ }dm \equiv \int 1_{A} f\text{ }dm=0.
\end{equation*}

\bigskip  \noindent \textbf{P4} : If $a\geq 0$ and $b \geq 0$, then we have 

\begin{equation*}
\int (a g+b h) \ dm=a \int g \ dm+b \int h \ dm.
\end{equation*}

\bigskip  \noindent \textbf{P5} : The integral is non-decreasing, that is :

\begin{equation*}
f\geq g\Rightarrow \int  f  \ dm \geq \int  g \ dm.
\end{equation*}

\bigskip  \noindent \textbf{P6} : Let $A$ and $B$ are disjoint, then 

\begin{equation*}
\int_{A+B} f \ dm=\int_{A} f \ dm+\int_{B} f \ dm.
\end{equation*}

\bigskip  \noindent \textbf{P7} Consider two non-decreasing sequences of simple functions $(f_{n})_{n\geq 0}$ and $(g_{n})_{n\geq 0}$. If 

$$
\lim_{n\uparrow +\infty } f_{n}  \leq g_{k},
$$

\bigskip \noindent then, we have

$$
\lim_{n\uparrow +\infty } \int f_{n} \ dm  \leq \int g_{k} \ dm.
$$

\noindent \textbf{Skip} this part and, please follow the proof in \textit{Doc 05-08} below. Repeat the proof as a homework.\\

\bigskip \noindent \textbf{Exercise 4}. \label{exercise04_doc05-02} (Coherence of the definition of the integrable of the non-negative and measurable function). Let  $(f_{n})_{(n\geq 0)} \subset \mathcal{E}_{+}$ and $(g_{n})_{(n\geq 0)} \subset \mathcal{E}_{+}$ be two non-decreasing sequences of simple functions such that

\begin{equation*}
\lim_{n\uparrow +\infty} h_{n} =\lim_{n\uparrow +\infty} g_{n}=f. \ \ (EDL)
\end{equation*}

\bigskip \noindent Show, by using Property (P7) in Exercise 3 above, that

\begin{equation*}
\lim_{n\uparrow +\infty} \int f_{n} \ dm =\lim_{n\uparrow +\infty} \int  g_{n} \ dm. \ \ (ELS)
\end{equation*}

\bigskip \noindent \textbf{Exercise 5}. \label{exercise05_doc05-02} Extend the Properties (P1) to (P5) to the class of all non-negative and measurable functions.\\

\noindent Here, $a$ and $b$ are two real numbers,  $f$ and $g$ are two non-negative and measurable functions, $A$ and $B$ are disjoint measurable sets. Let us denote

$$
\mathcal{L}_{0}^{+}(\Omega, \mathcal{A})=\{f \in \mathcal{L}_{0}(\Omega, \mathcal{A}), \ f\geq 0\}.
$$

\newpage 
\noindent \textbf{Part III - General functions}.\\

\noindent Hereafter, $\alpha$ and $\beta$ are two real numbers,  $f$ and $g$ are two measurable functions, $A$ and $B$ are disjoint measurable sets.\\

\bigskip \noindent \textbf{Exercise 6}. \label{exercise06_doc05-02} (Positive and negative parts).\\

\noindent Question (a) Show that the unique solution of the functional system

$$
f=h-g, \ h\geq 0, \ g \geq 0, hg=0, \ \ (FS) 
$$
 
\bigskip \noindent is given by $h=f^+=max(f,0)$ and $g=f^-=max(f,0)$. Besides, we have

$$
|f|=f^+ \ \ f^-. \ \ (AV)
$$

\bigskip \noindent and

$$
(f+g)^+ \leq f^+ + g^+ \ \ and  \ \ (f+g)^- \leq f^- + g^-. \ \  (BO)
$$

\bigskip \noindent and

$$
(f \leq g ) \Leftrightarrow (f^+ \leq g^+ \ and \ f^- \geq g^-), \ \ (COM) 
$$

\bigskip \noindent Question (b) Show that we have for any real constant $c$ :

$$
If \ c>0, \ then \ (cf)^+=cf^+ \ and \ (cf)^-=cf^-, \ \ (MPN1)
$$ 

\bigskip \noindent and

$$
If \ c<0, \ then \ (cf)^+=-cf^- \ and \ (cf)^-=-c f^+, \ \ (MPN2)
$$ 

\bigskip \noindent Finally, if $f$ is finite $a.e.$, and if $c=0$, $cf=(cf)^+=(cf)^-=0 \ a.e$\\

\bigskip \noindent Question (c) Suppose that $f$ is finite and $g$ is a function of condtant sign. Show that
 
$$
If \ g\geq 0, \ then \ (gf)^+=gf^+ \ and \ (gf)^-=gf^-, \ \ (MPNF1)
$$ 

\bigskip \noindent and

$$
If \ g \leq 0, \ then \ (gf)^+=-gf^- \ and \ (gf)^-=-g f^+. \ \ (MPNF2)
$$

\bigskip \noindent \textbf{Exercise 7}. \label{exercise07_doc05-02} Show the following points.\\

\noindent Question (a)  

$$
\int |f| dm = \int f^+ \ dm + \int f^- \ dm
$$ 

\bigskip \noindent $f$ is integrable if and only if $|f|$ is.\\

\noindent Question (b) If $\int f \ dm$ exists, then

$$
\left| \int  f \ dm \right| \leq \int |f| \ dm.
$$

\bigskip \noindent Question (c) Show that if $|f|\geq g$ and $g$ is integrable, then $f$ is.\\

\bigskip \noindent \textbf{Exercise 8}. \label{exercise08_doc05-02} (Markov inequality) Prove that for any $\lambda>0$

$$
m(|f| \geq \lambda) \leq \frac{1}{\lambda} \int |f| \ dm.  
$$

\bigskip \noindent \textit{Hint}. Put $A=(|f| \geq \lambda)$. Decompose the integral of $|f|$ over $A$ and $A^c$, drop the part on $A^c$ and use Property (P4) for non-negative measurable functions.\\

\bigskip \noindent \textbf{Exercise 9}. \label{exercise09_doc05-02} Use the Markov inequality to show that an integrable function is finite $a.e.$.\\

\noindent Hint. Use the following result of \textit{Exercise 7 of Doc 04-02 in Chapter \ref{04_measures} (page \pageref{exercise07_doc04-02}} : $(f \ infinite)=(|f|=+\infty)=\bigcap_{k\geq 1} (|f|\geq k)$. Remark that the $(|f|\geq k)$'s form a non-increasing sequence in $k$. Apply the Markov inequality to conclude.\\

\bigskip \noindent \textbf{Exercise 10}. \label{exercise10_doc05-02} Let 

\noindent Question (a) Suppose $\int f \ dm$ exists and let $c$ be a non-zero real number. Show that

$$
\int (cf) \ dm= c \int f \ dm.
$$

\bigskip \noindent Question (b) By taking $c=-1$ in question (a), tell simply why Properties (P2), (P4), (P4), (P5), (P6) are valid for non-positive functions, in other words for function of a constant sign.\\

\bigskip \noindent Question (c) Suppose that $f$ is an arbitrary measurable function and $A$ is a $m$-nul. Show that  $\int_{A} f \ dm$ exists and we have

$$
\int_{A} f \ dm=0.
$$

\bigskip \noindent Question (d) Suppose that $f$ exists and let $A$ and $B$ be two disjoint measurable sets. Show that we have

$$
\int_{A+B} f \ dm=\int_{A} f \ dm + \int_{B} f \ dm.
$$

\bigskip  \noindent (e). Let $f$ and $g$ be two measurable functions such that $f=g$ a.e. Show that their integrals exist or not simultaneously and, of they exist, we have
 
$$
\int f \ dm = \int g \ dm= \int_M f \ dm =\int_M g \ dm,
$$

\bigskip \noindent where $M=(f=g)$.\\

\bigskip  \noindent (f) Suppose that $f$ and $g$ are integrable. Show that 
$$
\int (f+g) \ dm=\int f \ dm + \int g \ dm.
$$

\bigskip \noindent Hint : Consider the measurable decomposition of $\Omega$ :

$$
\Omega =(f\geq 0, g \geq 0)+(f < 0, g \geq 0)+(f\geq 0, g <0)+(f < 0, g < 0)=A+B+C+D.
$$

\bigskip \noindent Each of these sets is, itself, decomposable into two set :
$$
A=A \cap (f+g \geq 0) + A \cap (f+g < 0)=A_1 + A _2.
$$

\bigskip \noindent define likewise $B_i$, $C_i$, $D_i$, $i=1,2$. Show that 
$$
\int_E (f+g) \ dm=\int_E f  dm + \int_E f \ dm.
$$

\bigskip \noindent when $E$ is any of the $A_i$, $B_i$, $C_i$, $D_i$, $i=1,2$ and conclude.\\

\noindent Hint. The key tool is to manage to have a decomposition of one of elements of three terms
$\int_E (f+g) \ dm$, $\int_E f  dm$, $\int_E f \ dm$ into the two others such that these two ones have the same constant signs.\\

\noindent Here, prove this only for one of them.\\

\bigskip  \noindent (g) Suppose $\int f \ dm$, $\int g \ dm$ exist and $f+g$ is defined and $\int f+g \ dm$ exist. Show that 

$$
\int (f+g) \ dm=\int f \ dm + \int f \ dm.
$$

\bigskip \noindent \textbf{Recommendation} This extension of Question (e) is rather technical. We recommend that you read the solution in the Appendix below.\\

\noindent (h) Suppose that $\int f \ dm$, $\int g \ dm$ exist and $f\leq g$. Show that

$$
\int f \ dm \geq \int g \ dm. \ \ (INC)
$$ 

\bigskip \noindent \textbf{Exercise 11}. \label{exercise11_doc05-02} Combine all the results of the different questions of Exercise 9 to establish the final properties of the integrals given in 

\textbf{II - Main properties} of \textit{Doc 05-01 of this chapter}. Use sentences only.\\

\noindent \textbf{From now, you only use these main properties which includes all the partial ones you gradually proved to reach them}.\\

\bigskip \noindent \textbf{Exercise 12}. \label{exercise11_doc05-02}  (Monotone Convergence Theorem of series).\\

\noindent Let $f : \mathbb{Z} \mapsto \overline{\mathbb{R}}_{+}$ be non-negative function and let $f_p : \mathbb{Z} \mapsto \overline{\mathbb{R}}_{+}$,  $p\geq 1$, be a sequence of non-negative functions increasing to $f$ in the following sense :

\begin{equation}
\forall (n\in \mathbb{Z}), \ \ 0\leq f_{p}(n)\uparrow f(n)\text{ as }p\uparrow\infty .  \label{comptage01}
\end{equation}

\bigskip
\noindent Then, we have 
\begin{equation*}
\sum_{n\in \mathbb{Z}}f_{p}(n)\uparrow \sum_{n\in \mathbb{Z}}f(n).
\end{equation*}

\newpage
\noindent \textbf{Appendix}. Solution of Question (g) of Exercise 10 : If $f+g$ is defined a.e., $\int f \ dm$, $\int g \ dm$ and $\int f+g \ dm$ exist and the addition 
$\int f \ dm + \int g \ dm$ makes sense, then

$$
\int f+g \ dm = \int f  dm  + \int g \ dm. \ (AD)
$$

\bigskip \noindent Let us consider three cases.\\

\noindent \textit{Case 1}. $f$ and $g$ both integrable. This case reduces to Question (e) of Exercise 10.\\

\noindent \textit{Case 2}. One of $f$ and $g$ is integrable. It will be enough to give the solution for $f$ integrable.  Let us do it by supposing that $f$ is integrable. We adopt the same prove in the proof
of Question(c) of Exercise 10. Indeed, in all the eight (8), we only need to move $\int f \ dm$ from one member to another, which possible thanks to its finiteness.

\noindent \textit{Case 3}. Both $f$ and $g$ have infinite integrals. If the sum $\int f \ dm + \int g \ dm$ makes sense, that means that the two integrals are both $+\infty$ or both $-\infty$.

$$
\int f^+ \ dm <+\infty, \ \ \int f^- \ dm =+\infty, \ \ \int g^+ \ dm <+\infty \ \ and \ \ \int g^- \ dm =+\infty 
$$

\bigskip \noindent or

$$
\int f^+ \ dm =+\infty, \ \ \int f^- \ dm <+\infty, \ \ \int g^+ \ dm =\infty \ \ and \ \ \int g^- \ dm <+\infty.
$$

\bigskip \noindent Let us continue through two sub-cases.\\

\noindent First, let us remark that the method used in the solution of Question (f) of Exercise 10 is valid until the statements of the eight decompositions at Line marked (APC). From that point, the only problem concerned the moving of a finite quantity from the left member to the second after we applied the integral operator at each member and used the linearity for functions of same constant signs. But this concerned only the cases (3), (4), (5) and (6).\\

\noindent Here, we are going to use the finiteness of both integrals of $f$ and $g$, or both integrals of negative parts.\\

\noindent \textbf{Case (i) where the positives parts of $f$ and $g$ are finite}.\\

\noindent In the three the following three cases, we may directly apply the explained principles.\\

\noindent (3) On $B_1=(f < 0, \ g \geq 0, \ f+g \geq 0)$, we use a sum non-positive functions :

$$
1_E -(f+g)  + 1_E  -g =   1_E -f; \ (In \ this \ case : \ \int 1_E  -g \ dm = - \int_E  g^+ \ dm \ finite).
$$

\bigskip \noindent (4) On $B_2=(f < 0, \ g \geq 0, f+g<0)$, we use a sum non-positive functions :

$$
1_E (f+g)   + 1_E  -g =   1_E -f. \ (In \ this \ case : \ \int 1_E  -g \ dm = - \int_E  g^+ \ dm \ finite).
$$

\bigskip \noindent (6) On $C_2=(f\geq 0, \ g <0, f+g<0)$, we use a sum non-negative functions :

$$
1_E -(f+g) + 1_E  f =  1_E -g. \ (In \ this \ case : \ \int 1_E  f \ dm =  \int_E  f^+ \ dm \ finite).
$$

\bigskip \noindent The case (5) is handled as follows.\\:

\noindent (5) On $C_1=(f \geq 0, \ g <0, \ f+g \geq 0)$, we use a sum non-positive functions : *

$$
1_E -(f+g)  + 1_E g =  1_E  -f. 
$$

\bigskip \noindent Since $\int_E  f^ + \ dm=\int_E  f^+ \ dm$ is finite, we get by the additivity of integrals of functions of same constants signs that $- \int 1_E -(f+g) \ dm$ and 
$\int 1_E g \ dm$ are both finite.\\

\noindent In total, in all these cases, we will be able to move the second term of the first member of the equalities after having applied the integral operator.\\

\noindent \textbf{Case (ii) where the negative parts of $f$ and $g$ are finite}.\\

\noindent Here, we have :

\noindent (3) On $B_1=(f < 0, \ g \geq 0, \ f+g \geq 0)$, we use a sum if non-positive functions :

$$
1_E -(f+g)  + 1_E  f =   1_E -g ; \ (In \ this \ case : \ \int 1_E  f \ dm = - \int_E  f^- \ dm \ finite).
$$

\bigskip \noindent (4) On $B_2=(f < 0, \ g \geq 0, f+g<0)$, we use a sum of non-negative functions :

$$
1_E -(f+g)   + 1_E  -f =  -1_E -g. \ (In \ this \ case : \ \int 1_E  -g \ dm = - \int_E  g^+ \ dm \ finite).
$$

\bigskip \noindent (5) On $C_1=(f \geq 0, \ g <0, \ f+g \geq 0)$, we use a sum of of non-positive functions :

$$
1_E -(f+g)  + 1_E g =  1_E  (-f). \ (In \ this \ case : \ \int 1_E  -g \ dm = - \int_E  g^+ \ dm \ finite).
$$

\bigskip \noindent (6) On $C_2=(f\geq 0, \ g <0, f+g<0)$, we use a sum non-positive functions :

$$
1_E -(f+g) + 1_E  f =  1_E -g.
$$

\bigskip \noindent The conclusion of the treatments of the cases (3), (4) and (6) in case (i) are repeated for (3), (4) and (5) in case (ii). The treatment of (5) in case (i) repeated
for (6) in case (ii). $\blacksquare$

\newpage
\noindent \LARGE \textbf{Doc 05-03  : Integration with the counting measure - Exercises}. \label{doc05-03}\\
\bigskip
\Large

\bigskip

\bigskip \noindent \textbf{Exercise 1}. \label{exercise01_doc05-03} (Integrating with respect to the counting measure).\\

\noindent Define the counting measure $\nu$ on a countable space$\mathbb{I}$ that is endowed with the discrete sigma-algebra, that is $\mathcal{P}(\mathbb{I})$ by
\begin{equation*}
\nu =\sum_{n\in \mathbb{I}}\delta _{n}
\end{equation*}

\bigskip \noindent and for any subset $A$ of $\mathbb{I}$. 
\begin{equation*}
\nu (A)=\sum_{n\in \mathbb{I}}\delta _{n}(A)=\sum_{n\in \mathbb{I}%
}1_{A}(n)=Card(A)
\end{equation*}

\bigskip \noindent  Any function  
\begin{equation*}
f:\mathbb{I}\mapsto \overline{\mathbb{R}}
\end{equation*}

\bigskip \noindent is measurable. Actually, such a function is a real sequence of real numbers $(f(n))_{n\in \mathbb{I}}$.\\

\noindent Your are going to find how the integral of $f$ with respect to $\nu$ is constructed through the four steps. We are going to see that, for any constant-sign function $f$ or integrabke function, we have

$$
\int f \ d\nu =\sum_{n\in \mathbb{I}}f(n). \ \ (COUNTINTEG)
$$

\bigskip \noindent \textit{Question (a) [Step 1]}. Let $f$ be an indicator function : $f=1_{A}$. Show that Formula (COUNTINTEG) holds.\\
 
\noindent \textit{Question (b) [Step 2]} Let $f$ be  simple and non-negative function, of the form 
\begin{equation*}
f=\sum_{1\leq i\leq k}\alpha _{i}1_{A_{i}}
\end{equation*}

\bigskip \noindent  where the $A_{i}$ are subsets of $\mathbb{I}$ and the $\alpha _{i}$ are finite non-negative numbers. Show also that Formula (COUNTINTEG) holds.\\

\bigskip \noindent \textbf{Question (c) [Step 3]}. Let $f$ be non-negative. Consider a sequence of simple functions $(f_{p})_{p\geq 0}$ increasing to $f$, that is 
\begin{equation}
\forall (n\in \mathbb{N}),0\leq \text{ }f_{p}(n)\uparrow f(n)\text{ as }n\uparrow \infty.   \label{comptage01}
\end{equation}

\bigskip \noindent By combining Question (b), \textit{Exercise 12 in Doc 05-02} and the construction on the integral of a non-negative and measurable function, show that  Formula (COUNTINTEG) 
still holds
 
\bigskip \noindent \textbf{Question (d) [Step 4]} Let $f$ be an arbitrary function from $\mathbb{I}$ to $\overline{\mathbb{R}}$. Give the expressions of the positive and the negative parts of $f$.
Show that if $f$ is $\nu$-integrable, then Formula (COUNTINTEG) holds.\\

\noindent \textbf{Question (e)} Let $f$ be an arbitrary function from $\mathcal{Z}$ to $\overline{\mathbb{R}}$. Show that if $f$ is quasi-integrable, then Formula (COUNTINTEG) holds.

\noindent \textbf{Question (f)}. Propose a conclusion of the exercise.\\

\bigskip \noindent \textbf{Exercise 2}. \label{exercise02_doc05-03} Consider the following sequences

$$
f(n)=(-1)^{n} g(n), n\geq 1,
$$

\bigskip \noindent where $g(n)$ is non-increasing and $g(n) \rightarrow 0$ as $n\rightarrow +\infty$. We know from the earlier courses of calculus that 
$\sum_{n \geq 1} f(n)$ exists and is finite.\\

\noindent Now, take $g(n)=1/(n+1)$, $n\geq 0$.\\

\noindent Question : Can you say that $\sum_{n \geq 1} (-1)^{n}/(n+1)$ is the integral of the counting measure on $\mathbb{N}$? Why? \\

\bigskip \noindent \textbf{Exercise 3}. \label{exercise03_doc05-03}  (Abel's rule) Let $(v_n)_{n\geq 0}$ be a totally bounded sequence of real numbers, that is, there exists a real number $A$ such that for any $n\geq 0$, $m\geq n$,

$$
|v_n + v_{n+1} + \cdots + v_{m}| \leq A.
$$

\bigskip \noindent Let $(\varepsilon_n)_{n\geq 0}$ be a sequence of reals number such that 
$$
\sum_{n\geq 1} |\varepsilon_{n}-\varepsilon_{n-1}|<+\infty \ (A1)
$$

\bigskip \noindent and

$$
\lim_{n\rightarrow +\infty} \varepsilon_{n}=0. \ (A2)
$$

\bigskip \noindent Show that the series 
$$
\sum_{n\geq 1} v_{n} \varepsilon_{n}
$$

\bigskip \noindent is convergent.\\

\noindent Hint. Put for $n \leq m$, $R_{n,m}=\sum_{n\leq k \leq m} v_{n} \varepsilon_{n}$ and $V_{n,j}=v_n +  \cdots + v_{n+j}$, $j\geq 0$. Show that :
\begin{eqnarray*}
R_{n,m}&=& \varepsilon_{n} V_{n,0}+ \varepsilon_{n+1} (V_{n,1}-V_{n,0})+\cdots+\varepsilon_{m} (V_{n,m-2}-V_{n,n-m-1}\\
&=&\biggr( V_{n,0} (\varepsilon_{n}-\varepsilon_{n+1}) + \cdots + V_{n,m-m-1} (\varepsilon_{m-n-2}-\varepsilon_{m-n-1}) \biggr)+ \biggr(V_{n,m-n} \varepsilon_{m}\biggr)
\end{eqnarray*}

\bigskip \noindent The second line is obtained by developing the first. From there, use (A1) and (A2) to show that

$$
R_{n,m}=\sum_{n\leq k \leq m} v_{k} \varepsilon_{k} \rightarrow 0 \ as \ (n,m) \rightarrow (+\infty, +\infty). 
$$

\bigskip \noindent Conclude

\newpage
\noindent \LARGE \textbf{Doc 05-04 : Lebesgue/Riemann-Stieljes integrals on $\mathbb{R}$ - Exercises}. \label{doc05-04}\\
\bigskip
\Large 

\bigskip \noindent This document is only an introduction to the comparison between the Riemann-Stieljes Integral and the Lebesgue-Stieljes integral.\\

\noindent Later, we will come back to a more deep document on the same subject after we have mastered the convergence theorems.\\
  
\bigskip \noindent \textbf{Reminder of the definition of the Riemann-Stieljes Integral}.\newline

\noindent The Classical Riemann-Stieljes Integral is defined for bounded functions on compacts intervals. In all this text, $a$ and $b$ are two real numbers  such that $a<b$ and $F: [a,b]\longrightarrow \mathbb{R}$ is non-constant non-decreasing function.\\
 
\noindent \textbf{Definition of the Riemann-Stieljes Integral on a compact set}.\\

\noindent Let $f : \ ]a,b]\longrightarrow \mathbb{R}$ be an arbitrary function. We are going to define the Riemann-Stieljes integral of $f$ on $]a,b]$ associated with $F$ and denoted 
\begin{equation*}
I=\int_{a}^{b}f(x)\text{ }dF(x).
\end{equation*}

\bigskip \noindent We begin to define the Riemann-Stieljes sums. For each $n\geq 1,$ consider a subdivision $\pi_n$ of $]a,b]$, that divides into $]a,b]$ into the $\ell(n)$ sub-intervals

\begin{equation*}
]a,b]=\sum_{i=0}^{\ell (n)-1}]x_{i,n},x_{i+1,n}],
\end{equation*}

\bigskip \noindent with $a=x_{0,n}<x_{1,n}< \cdots < x_{\ell(n),n}=b$. The modulus of the subdivision $\pi _{n}$ is defined by%
\begin{equation*}
m(\pi _{n})=\max_{0\leq i\leq m(n)}(x_{i+1,n}-x_{i+1,n}).
\end{equation*}

\bigskip \noindent Associate to this subdivision $\pi _{n}$ an arbitrary sequence $c_{n}=(c_{i,n})_{1\leq i\leq \ell (n)}$ such that $c_{0,n}\in ]a,x_{1,n}]$ and $c_{i,n}\in \lbrack x_{i,n},x_{i+1,n}],$ $1\leq i\leq \ell(n)-1$.\\

\noindent We may define a sequence of Riemann-Stieljes sum associated to the subdivision $\pi_n$ and the vector $c_n$ in the form 

\begin{equation*}
S_{n}(f,F,a,b,\pi _{n},c_{n})=\sum_{i=0}^{\ell
(n)-1}f(c_{i,n})(F(x_{i+1,n})-F(x_{i,n})).  \label{RSS}
\end{equation*}

\bigskip \noindent Since $f,$ $F$ and $]a,b]$ are fixed here, we may and do drop them in the expression of $S_n$.

\bigskip \noindent \textbf{Definition}. A bounded function $f$ is Riemann-Stieljes integrable with respect to $F$ if there exists a real
number $I$ such that any sequence of Riemann-Stieljes sums $S_{n}(\pi_{n},c_{n})$ converges to $I$ as $n\rightarrow 0$\ \ whenever $m(\pi
_{n})\rightarrow 0$ as $n\rightarrow \infty$, and the number $I$ is called the Riemann-Stieljes of $f$ on $[a,b]$.\\

\noindent If $F$ is the identity function, that is $F(x)=x$, $x\in \mathbb{R}$, $I$ is simply called the Riemann integral of $f$ over $[a,b]$ and the sums defined in Formula (RSS) are called Riemann sums.\\

\bigskip \noindent \textbf{Exercise 1.} \label{exercise01_doc05-04} \\

\noindent Question (a) Let $f : ]a,b]\longrightarrow \mathbb{R}$ be a bounded and measurable function on $[a,b]$. Is-it Lebesque-Stieljes integrable? If yes, give a bound of its integral.\\

\noindent Question (b) Consider the function

$$
g=1_{]a,b]\cap \mathbb{Q}},
$$

\bigskip \noindent where $\mathbb{Q}$ is the set of rational numbers.\\

\noindent Consider Riemann-Stieljes sums associated to $g$ in which the points $c_{i,n}$ are chosen in $\mathbb{Q}$. Justify why you can make this choice. What are the values of these Riemann-Stieljes sums.\\

\noindent Proceed similarly by choosing the points $c_{n,i}$ as irrational numbers. Justify why you can make this choice. What are the values of these Riemann-Stieljes sums.\\

\noindent Conclude that $g$ is not Riemann-Stieljes integrable.\\

\noindent Question (c). Make a first comparison between the two integrals.\\

\bigskip \noindent \textbf{Exercise 2}. \label{exercise02_doc05-04} \\

\noindent We require the following assumption on $F$ : $F$ is of bounded variation on $[a,b]$, that is, there exists a real number$A>0$ such that for any subdivision $\pi=(z_0,...,z_s)$ of 
$[a,b]$, that is $a=z_0<z_1<...<z_s=b$, we have

$$
V(\pi,[a,b])=\sum_{j=0}^{s-1} |F(z_{i+1}-z_{i}| <A.
$$

\bigskip \noindent Question (a) Prove that if $f$ is continuous then the sequences of Riemann-Stieljes (RSS) is Cauchy, and then, the Rieman-integral of $f$ exists on $[a,b]$.\\

\noindent \textit{Hint}. By using the results \textit{Exercise 15 in Doc 04-02 of Chapter \ref{04_measures}}, page \pageref{exercise15_doc04-02}, it is direct to see that the application defined by

$$
\lambda_F(]x,y])=F(b)-F(a)
$$

\bigskip \noindent is a proper and an additive application on the semi-algebra $\mathcal{S}([a,b])=\{{a}, [a,b], ]x,y], \ a\leq x \leq y \leq b\}$.\\

\noindent We stress that the extension theory made in Exercise 15 is still valid expect the properties related to sub-additivity. So, $\lambda_F$ is directly extended to an additive application one the class of finite sums of type $]x,y]$ in $[a,b]$. But, be careful : we cannot have the no-negativity without further assumptions on $F$, what we do not require for the moment. The theory of integral of elementary functions(but based here on intervals of the form $]x,y]$ is also valid. Precisely, if we denote

$$
h= \sum_{i=0}^{\ell (n)-1}f(c_{i,n}) 1_{]x_{i,n}, \ x_{i-1,n}]},
$$

\bigskip \noindent the \textit{integral},

$$
\int h \ d\lambda_F =\sum_{i=0}^{\ell (n)-1}f(c_{i,n}) F(x_{i+1})-F(x_{i,n}).
$$

\bigskip \noindent In the same spirit of the remarks above, this integral defined on the $\mathcal{E}_e$class of functions of the form

$$
\sum_{i=0}^{s-1} \alpha_i 1_{]z_{i+1}, \ z_i]}, \ a=z_0<z_1<...<z_s=b, \alpha_i \in \mathbb{R}, \  0\leq i \leq s-1,
$$

\bigskip \noindent is still well-defined and its linearity holds since we do not face infinite values. In the sequel, many properties in the solutions of Exercises series in Doc 05-05 are used.\\
 
\noindent Now take two Riemann-integrals for integers $p>1$ and $q>1$

\begin{equation*}
S_{p}=\sum_{i=0}^{\ell (p)-1}f(c_{i,p})(F(x_{i+1,p})-F(x_{i,p})).  \label{RSS1}
\end{equation*}

\bigskip \noindent and 

\begin{equation*}
S_{q}=\sum_{i=0}^{\ell(q)-1}f(c_{i,q})(F(x_{i+1,q})-F(x_{i,q})).  \label{RSS2}
\end{equation*}

\bigskip \noindent with, for $r=1,2$, we have $a=x_{0,r}<x_{1,r}<\cdots<x_{\ell(r),r}$, $c_{0,r}\in ]a,x_{1,r}]$ and $c_{i,r}\in \lbrack x_{i,r},x_{i+1,r}],$ $1\leq i\leq \ell (r)-1$,  $\pi_{r}=(x_{0,r}, ..., <x_{\ell(r),r})$.\\

\noindent Denote

$$
h_p= \sum_{i=0}^{\ell (p)-1}f(c_{i,p}) 1_{]x_{i,p}, \ x_{i-1,p}]},
$$

\bigskip \noindent and

$$
h_q= \sum_{i=0}^{\ell (q)-1}f(c_{i,q}) 1_{]x_{i,q}, \ x_{i-1,q}]},
$$

\bigskip \noindent and use the identities

$$
S_{p} = \int f_p \ d\lambda_F \ and \ S_{q} = \int f_q \ d\lambda_F
$$

\bigskip \noindent To lessen the notations, put $A_i=]x_{i,p}, ]x_{i+1,p}$ and $B_j=]x_{j,q}, ]x_{j+1,p}$.\\

\noindent Now, use the techniques of the superpositions of partitions as seen in \textit{Point 03.22a of Doc 03-02 of Chapter \ref{03_setsmes_applimes_cas_speciaux}}, page \pageref{doc03-02} and use the superposition $\{ A_{j}B_{j}, \ (i,j)\in H \}$ of the partitions $\{ A_{i},  0\leq i\leq \ell(p)-1 \}$ and $\{ B_{j}, \ 0\leq j\leq \ell(q)-1 \}$, where $H=\{(i,j), \ 1\leq i\leq \ell(p), \  1\leq j\leq \ell(q), A_iB_j\neq \emptyset \}$. Provide new expressions of $h_p$ and $h_q$ related to the partition $\{ A_{j}B_{j}, \ (i,j)\in H\}$ of $]a,b]$.\\

\noindent Finally, use the additivity of $\lambda_F$, in particular \textit{Formulas (DA), (DA1) and (DA2) in the solution of Exercise 15 in Doc 05-05},  to establish that

$$
|S_p-S_q|\leq A \sup_{(s,t)\in [a,b]^2, |s-t|\leq m(\pi_p) \vee m(\pi_s)} |f(t)-f(s)|
$$

\bigskip \noindent Deduce from this, that the sequence $(S_n)_{n\geq 1}$ is Cauchy if $m(\pi_n)$ tends to zero as $n \rightarrow +\infty$.\\

\newpage
\noindent \LARGE \textbf{Doc 05-05 : Integrals of elementary functions - Exercises with solutions}. \label{doc05-05}\\
\bigskip
\Large

\noindent \textbf{NB} : \textbf{In this work, you will have to use the \textit{demonstration tools} of Part 04.12 (see \textbf{04.12a} and \textbf{04.12b}) of Doc 12.}\\

\noindent Let $(\Omega, \mathcal{A}, m)$ be measure space. Let $\mathcal{E}^{(1)}$ be the class of non-negative simple functions or integrable simple functions.\\

\noindent \textbf{Important Remarks}. The properties which are asked to be proved here are intermediate ones. They are used to prove the final properties which are exposed in \textbf{Part II}
 in \textbf{Doc 05-01} (page \pageref{doc05-01}) if this chapter..

\bigskip  \noindent \textbf{Part I - Properties of the integral on $\mathcal{E}^{(1)}$}.\\

\bigskip  \noindent Let be given a non-negative simple function  
\begin{equation*}
f=\sum_{1\leq i\leq k}\alpha _{i}1_{A_{i}}, \ (\alpha_i \in \mathbb{R}_{+}, \ A_i \in \mathcal{A}, \ 1\leq i \leq k), \ k\geq 1. \ \ (NSF)
\end{equation*}

\bigskip  \noindent define its integral with respect to $m$ by :

\begin{equation}
\int f \ dm=\sum_{1\leq i\leq k}\alpha _{i} \ m(A_{i}). \ \ (IFS)
\end{equation}

\bigskip  \noindent \textbf{Exercise 1}. \label{exercise01_sol_doc05-05} (Coherence) Prove that this definition is coherent. In order to do this, consider two expressions of $f$, use the expression based on the partition of $\Omega$ formed by the superposition of the two first partitions as in \textit{Point 03.22a of Doc 03-02 of Chapter \ref{03_setsmes_applimes_cas_speciaux}}, page \pageref{doc03-02}.

\bigskip \noindent \textbf{Solutions}. Suppose that $f$ has an other expression

$$
f=\sum_{1\leq j\leq m}\beta _{j}1_{B_{j}}, \ \ (\beta_j \in \mathbb{R}_{+}, \ B_j \in \mathcal{A}, \ 1\leq j \leq \ell), \ \ell\geq 1.
$$

\bigskip \noindent S The definition of the integral of the non-negative elementary is well done if we prove that 

\begin{equation}
\sum_{1\leq i\leq k}\alpha _{i} \ m(A_{i}) = \sum_{1\leq i\leq \ell} \beta_{j} \ m(B_{j}). \ \ (CINSF)
\end{equation}

\noindent Let us use the superposition $\{ A_{j}B_{j}, \text{  } 1\leq i\leq k, \text{  } 1\leq j\leq \ell \}$ of the partitions $\{ A_{i},  1\leq i\leq k \}$ and 
$\{ B_{j}, \ 1\leq j\leq \ell \}$, as seen in \textit{Point 03.22a of Doc 03-02 of Chapter \ref{03_setsmes_applimes_cas_speciaux}}, page \pageref{doc03-02}, to have

$$
X=\sum_{(i,j)\in I}\gamma_{ij} \ 1_{A_{i}B_{j}},
$$

\bigskip \noindent S where $I=\{(i,j), \ 1\leq i \leq k, \  1\leq j \leq k, \ A_iB_j \neq \emptyset \}$, $\gamma_{i,j}=\alpha_i=\beta_j$ for $(i,j) \in I$ and arbitrary (say zero) for $(i,j) \notin I$. \noindent It is then clear, by that for each $i \ \{1,...,k\}$,

$$
A_i=\sum_{j=1}^{\ell} A_{i}B_{j}=\sum_{1\leq j \leq \ell, \ (i,j) \in I} A_{i}B_{j}.
$$

\bigskip \noindent S Hence, we have for each $i \ \{1,...,k\}$,

$$
m(A_i)=\sum_{1\leq j \leq \ell, \ (i,j) \in I} m(A_{i}B_{j}). \ \ (DA)
$$

\bigskip \noindent S Now, we have

$$
\sum_{1\leq i\leq k} \sum_{1\leq i\leq \ell} \gamma_{ij} \ m(A_i B_{j})=\sum_{(i,j) \in I}  \gamma_{ij} \ m(A_i B_{j}),
$$

\bigskip \noindent S since $m(A_i B_{j})=0$ if $(i,j)$ is not in $I$. We get, by Fubini's formula for sums of finite number of real numbers,

$$
\sum_{1\leq i\leq k} \sum_{1\leq i\leq \ell}  \ m(A_i B_{j})=\sum_{i=1}^{k} \sum_{1\leq j \leq \ell, \ (i,j) \in I} \gamma_{ij} m(A_{i}B_{j}),
$$

\bigskip \noindent S which, by Formula (DA) above and by additivity of $m$, leads to

\begin{eqnarray*}
\sum_{1\leq i\leq k} \sum_{1\leq i\leq \ell} \gamma_{ij} \ m(A_i B_{j})&=&\sum_{i=1}^{k} \sum_{1\leq j \leq \ell, \ (i,j) \in I} \alpha_i m(A_{i}B_{j})\\
&=&\sum_{i=1}^{k} \alpha_i \sum_{1\leq j \leq \ell, \ (i,j) \in I}  m(A_{i}B_{j})\\
&=&\sum_{i=1}^{k} \alpha_i m(\sum_{1\leq j \leq \ell, \ (i,j) \in I} A_{i}B_{j}),\\
&=&\sum_{i=1}^{k} \alpha_i m(A_i),
\end{eqnarray*}

\bigskip \noindent S so that we get

$$
\sum_{i=1}^{k} \alpha_i m(A_i)=\sum_{1\leq i\leq k} \sum_{1\leq i\leq \ell} \gamma_{i,j} \ m(A_i B_{j}). \ \ (DA1)
$$

\bigskip \noindent S Since the roles of the partitions $\{ A_{i},  1\leq i\leq k \}$ and $\{ B_{j}, \ 1\leq j\leq \ell \}$ on one side, and that of the sets of numbers $\{ \alpha_{i},  1\leq i\leq k \}$ and 
$\{\beta_{j}, \ 1\leq j\leq \ell \}$ on the other side, are symmetrical, we also have

$$
\sum_{j=1}^{\ell} \beta_j m(B_j)=\sum_{1\leq i\leq k} \sum_{1\leq i\leq \ell} \gamma_{i,j} \ m(A_i B_{j}). \ \ (DA2)
$$

\bigskip \noindent S We may and do conclude by putting together Formulas (DA1) and (DA2).\\

\bigskip  \noindent \textbf{Exercise 2}. \label{exercise02_sol_doc05-05} Show that Formula (IFS) holds if $f$ is non-positive and if $f$ is integrable. Remark that the elementary function $f$ integrable if and only if, we have 

$$
\forall i \in \{1,...,k\}, \ (\alpha_i\neq 0) \Rightarrow m(A_i) <+\infty.
$$

\bigskip \noindent \textbf{Solutions}.\\

\noindent According to \textit{Point 05-01c of Doc 05-01 of this chapter} (page \pageref{doc05-01}), we have to check the integral of the positive and the negative parts. Here, in the special case of an elementary function, we may denote

$$
I^{+}=\{i, \ 1\leq i\leq k, \ \alpha_i\geq 0 \} \ \ and I^{+}=\{i, \ 1\leq i\leq k, \ \alpha_i\leq 0 \}.
$$

\bigskip \noindent S It is clear that we have

$$
f^{+}=\sum_{ i \in I^{+}} \alpha _{i}1_{A_{i}} \ \ and \ \ \int f^{+}\ dm=\sum_{ i \in I^{+}} \alpha _{i} m(A_{i})
$$

\bigskip \noindent S and

$$
f^{-}=\sum_{ i \in I^{-}} (-\alpha _{i}) 1_{A_{i}} \ \ and \ \ \int f^{-}\ dm=- \sum_{ i \in I^{-}} \alpha _{i} m(A_{i}).
$$

\bigskip \noindent S We see that, for $\int f^{+}\ dm$ to finite, we should have $\alpha _{i} m(A_{i})$, $i \in I^{+}$, whenever $m(A_{i})=+\infty$ by \textbf{}\textbf{Convention - Warning 1} in the earlier lines of \textit{Doc 05-01 of this chapter} (page \pageref{doc05-01}). This is possible only if $\alpha_i=0$. The same being true for the negative part, we have that the elementary function $f$ is integrable if and only all the $m(A_i)$'s are finite or $m(A_i)=+\infty$ if and only if $\alpha_i=0$ for all $1\leq i \leq k$.

\bigskip  \noindent \textbf{Exercise 3}. \label{exercise03_sol_doc05-05} Establish the following properties on $\mathcal{E}{+}$.\\

\noindent Hereafter, $\alpha$ and $\beta$ are two real numbers,  $f$ and $g$ are two elements of $\mathcal{E}^{+}$, $A$ and $B$ are measurable spaces.\\

\bigskip  \noindent \textbf{P1} : The integral is non-negative, that is :

$$
f\geq 0 \Rightarrow \int f \ dm\geq 0.
$$

\bigskip  \noindent \textbf{P2} : We have  

\begin{equation*}
f=0 \ a.e. \Leftrightarrow \int f \ dm=0.
\end{equation*}

\bigskip  \noindent \textbf{P3} If $A$ is a null-set, then 
\begin{equation*}
\int_{A}f\text{ }dm \equiv \int 1_{A} f\text{ }dm=0.
\end{equation*}

\bigskip  \noindent \textbf{P4} : If $a\geq 0$ and $b \geq 0$, then we have 

\begin{equation*}
\int (a g+b h) \ dm=a \int g \ dm+b \int h \ dm.
\end{equation*}

\bigskip  \noindent \textbf{P5} : The integral is non-decreasing, that is :

\begin{equation*}
f\geq g\Rightarrow \int  f  \ dm \geq \int  g \ dm.
\end{equation*}

\bigskip  \noindent \textbf{P6} : Let $A$ and $B$ are disjoint, then 

\begin{equation*}
\int_{A+B} f \ dm=\int_{A} f \ dm+\int_{B} f \ dm.
\end{equation*}

\bigskip  \noindent \textbf{P7} Consider two non-decreasing sequences of simple functions $(f_{n})_{n\geq 0}$ and $(g_{n})_{n\geq 0}$. If 

$$
\lim_{n\uparrow +\infty } f_{n}  \leq g_{k},
$$

\bigskip \noindent S then, we have

$$
\lim_{n\uparrow +\infty } \int f_{n} \ dm  \leq \int g_{k} \ dm.
$$

\bigskip \noindent S \textbf{Skip} this part and, please follow the proof in \textit{Doc 05-08} below. Repeat the proof as a homework.\\

\bigskip \noindent \textbf{Solutions}.\\

\noindent \textbf{(P1)}. If $f$ is a non-negative elementary function written as in Formula (NSF) above, then we have

$$
\forall i, \ 1\leq i \leq k, \ \alpha_i\geq 0,
$$

\bigskip \noindent S which, by Formula (IFS) above, implies that

$$
\int f \ dm \geq 0.\ \square
$$

\bigskip \noindent \textbf{(P2)} Let $f$ be a element of $\mathcal{E}^{+}$ written as in Formula (NSF). Denote

$$
I_1=\{i, \ 1\leq i \leq k, \ \alpha_i =0 \} \and \ I_2=\{i, \ 1\leq i \leq k, \ \alpha_i > 0 \}.
$$

\bigskip \noindent S It is clear that the set where $f$ is not null is given by

$$
(f \neq 0)=\sum_{i \in I_2} A_i
$$

\bigskip \noindent S and, by Formula (ISF) which valid on $\mathcal{E}_{+}$, we have

$$
\int f \ dm =\sum_{i \in I_2} \alpha_i m(A_i).
$$

\bigskip \noindent S (a) Let us suppose first tha $f=0 \ a.e$. Hence, we have 

$$
m(f \neq 0)=m(\sum_{i \in I_2} A_i)=0,
$$

\bigskip \noindent S it comes that $m(A_i)=0$ for each $i \in I_2$ and next, that

$$
\int f \ dm =\sum_{i \in I_2} \alpha_i m(A_i)=0. \square
$$ 

\bigskip \noindent Let us suppose that $\int f \ dm=0$. This means that

$$
\int f \ dm =\sum_{i \in I_2} \alpha_i m(A_i)=0.
$$ 

\bigskip \noindent S Since $\alpha_i>0$ for all $i \in I_2$, we have $m(A_i)=0$ for all $i \in I_2$. Hence, since 

$$
(f \neq 0)=\sum_{i \in I_2} A_i,
$$

\bigskip \noindent S we get

$$
(f \neq 0)=0,
$$

\bigskip \noindent S that is : $f=0$, $a.e$.\\

\noindent The conclusion is made by putting together (a) and (b).\\

\bigskip \noindent \textbf{(P3)} Let $A$ be a null set and let $f$ be a element of $\mathcal{E}_{+}$ written as in Formula (NSF). Then we have

$$
1_A f= \sum_{i=1}^{k} \alpha_i 1_{AA_i}.
$$

\bigskip \noindent S But for all $i$, $1 \leq i \leq k$, $m(AA_i) \leq m(A)=0$. Thus, we get

$$
\int 1_A f \ dm = \sum_{i=1}^{k} \alpha_i m(AA_i)=0. \ \square
$$

\bigskip \noindent \textbf{(P4)} Let $a \geq 0$, $b\geq 0$ and let $f$ be a element of $\mathcal{E}_{+}$ written as in Formula (NSF). Let $g$ another element of $\mathcal{E}_{+}$ written as

$$
f=\sum_{1\leq j\leq m}\beta _{j}1_{B_{j}}, \ \ (\beta_j \in \mathbb{R}_{+}, \ B_j \in \mathcal{A}, \ 1\leq j \leq \ell), \ \ell\geq 1.
$$

\bigskip \noindent S Let us use the superposition $\{ A_{j}B_{j}, \  1\leq i\leq k, \ 1\leq j\leq \ell \}$ of the partitions $\{ A_{i},  1\leq i\leq k \}$ and 
$\{ B_{j}, \ 1\leq j\leq \ell \}$, as seen in \textit{Point 03.22a of Doc 03-02 of Chapter \ref{03_setsmes_applimes_cas_speciaux}}, page \pageref{doc03-02}, to have

$$
a f = \sum_{(i,j) \in I} a \alpha_i 1_{A_i B_{j}} \ and \  b g= \sum_{(i,j) \in I} b \beta_j 1_{A_i B_{j}}.
$$

\bigskip \noindent S where $I=\{(i,j), \ 1\leq i \leq k, \  1\leq j \leq k, \ A_iB_j \neq \emptyset \}$, which leads to 

$$
a f +b g= \sum_{(i,j) \in I} (a \alpha_i + b \beta_j) 1_{A_i B_{j}},
$$

\bigskip \noindent S which, by definition, implies

$$
\int (a f +b g) \ dm = \sum_{(i,j) \in I} (a \alpha_i + b \beta_j) m(A_i B_{j}). \ \  (SIFS)
$$

\bigskip \noindent S This, in turn, implies

\begin{eqnarray*}
\int (a f +b g) \ dm &=& a \sum_{i=1}^{k} \alpha_i \sum_{1\leq j \leq \ell, \ (i,j) \ I} m(A_i B_{j})\\
&+& \sum_{j=1}^{\ell} \beta_j \sum_{1\leq i \leq k, \ (i,j) \ I} m(A_i B_{j})
\end{eqnarray*}

\bigskip \noindent S By using Formula (DA) in Exercise 1 for the composition of $A_i$ in the first line, and for the decomposition of $B_j$ in the second line, we get 

$$
\int (a f +b g) \ dm = a \sum_{i=1}^{k} \alpha_i m(A_i) + \sum_{j=1}^{\ell} m(B_{j}),
$$

\bigskip \noindent S that is

$$
\int (a f +b g) \ dm = a \int f \ dm + b \int  g \ dm. \ \square
$$

\bigskip \noindent \textbf{(P5)} Let $f \geq g$. We have $f= g + (f-g)$. Since $g$ and $f-g$ are non-negative elementary functions, Property (P4) just above implies

$$
\int f \ dm= \int g \ dm + \int (f-g)\ dm.
$$

\bigskip \noindent S Next, by Property (P1), $\int (f-g)\ dm \geq 0$. Hence,

$$
\int f \ dm \geq \int g \ dm.
$$

\bigskip \noindent \textbf{(P1)} Let $A$ and $B$ two disjoint measurable subsets of $\Omega$ and let $f$ be a element of $\mathcal{E}^{+}$ written as in Formula (NSF). We have

$$
1_{A+B} f = 1_{A} f + 1_{B} f.
$$

\bigskip \noindent S All the functions $1_{A+B} f$, $1_{A} f$ and $1_{B} f$ are in $\mathcal{E}^{+}$. By Property (P4), we have

$$
\int 1_{A+B} f \ dm= \int  1_{A} f \ dm + \int  1_{B} f \ dm,
$$

\bigskip \noindent S that is

$$
\int_{A+B} f \ dm= \int_{A} f \ dm + \int_{B} f \ dm.
$$

\bigskip \noindent \textbf{(P7)} You are strongly suggested to follow the  proof of this property in \textit{Doc 05-08} below and to repeat the proof as a homework.\\

\newpage

\bigskip  \noindent \textbf{Part 2 - Non-negative functions}.\\

\noindent Hereafter, $\alpha$ and $\beta$ are two real numbers,  $f$ and $g$ are two non-negative and measurable functions, $A$ and $B$ are measurable spaces.\\

\noindent \textbf{Exercise 4}. \label{exercise04_sol_doc05-05} (Coherence of the definition of the integrable of the non-negative and measurable function). Let  $(f_{n})_{(n\geq 0)} \subset \mathcal{E}_{+}$ and 
$(g_{n})_{(n\geq 0)} \subset \mathcal{E}_{+}$ be two non-decreasing sequences of simple functions such that

\begin{equation*}
\lim_{n\uparrow +\infty} h_{n} =\lim_{n\uparrow +\infty} g_{n}=f. \ \ (EDL)
\end{equation*}

\bigskip \noindent S Show, by using Property (P7) in Exercise 3 above, that

\begin{equation*}
\lim_{n\uparrow +\infty} \int f_{n} \ dm =\lim_{n\uparrow +\infty} \int  g_{n} \ dm. \ \ (ELS)
\end{equation*}

\bigskip \noindent \textbf{Solutions}.\\

\noindent It is clear that Formula (EDL) and the non-decreasingness of the sequences $(f_{n})_{(n\geq 0)}$ and $(g_{n})_{(n\geq 0)}$, that we have, for any $k\geq 0$,

\begin{equation*}
\lim_{n\uparrow +\infty}  f_n \geq g_{k}.
\end{equation*}

\bigskip \noindent S By letting $k\uparrow +\infty$ at right, we get

\begin{equation*}
\lim_{n\uparrow +\infty}  f_n \geq \lim_{k\uparrow +\infty} g_{k},
\end{equation*}

\bigskip \noindent S which is

\begin{equation*}
\lim_{n\uparrow +\infty}  f_n \geq \lim_{n\uparrow +\infty} g_{n}.
\end{equation*}

\bigskip \noindent S Since the roles of the sequences $(f_{n})_{(n\geq 0)}$ and $(g_{n})_{(n\geq 0)}$ are symmetrical here, we get

\begin{equation*}
\lim_{n\uparrow +\infty}  f_n = \lim_{n\uparrow +\infty} g_{n}.
\end{equation*}

\bigskip \noindent \textbf{Exercise 5}. \label{exercise05_sol_doc05-05} Extend the Properties (P1) to (P5) to the class of all non-negative and measurable functions.\\

\noindent Here, $a$ and $b$ are two real numbers,  $f$ and $g$ are two non-negative and measurable functions, $A$ and $B$ are disjoint measurable sets. Let us denote

$$
\mathcal{L}_{0}^{+}(\Omega, \mathcal{A})=\{f \in \mathcal{L}_{0}(\Omega, \mathcal{A}), \ f\geq 0\}.
$$

\bigskip \noindent \textbf{Solutions}.\\

\noindent Throughout the solutions of this exercise,  $(f_{n})_{(n\geq 0)} \subset \mathcal{E}_{+}$ and $(g_{n})_{(n\geq 0)} \subset \mathcal{E}_{+}$ are non-decreasing sequences such that

\begin{equation*}
f=\lim_{n\uparrow +\infty}  f_n \ \ and  \ \ g=\lim_{n\uparrow +\infty} g_{n}. \ (LF1)
\end{equation*}

\bigskip \noindent S and by definition,

$$
\int f \ dm=\lim_{n\uparrow +\infty}  \int f_n \ dm \ \ and \ \ \int g \ dm=\lim_{n\uparrow +\infty}  \int g_n \ dm. \ \ (LF2)
$$

\bigskip \noindent \textbf{(P1)} Since the numbers $\int f_n \ dm$ are non-negative by property (P1) of Exercise 1, we have, by Formula (LF2) above, that $\int f \ dm\geq 0.$ $\square$\\

\bigskip \noindent \textbf{(P2)}\\

\noindent (a) Suppose that $m(f \neq 0)=0$. Since the functions $f$ and $f_n$, $n\geq 0$, are non-negative, we have :$(f \neq 0)=(f>0)$ and $(f \neq 0)=(f_n>0)$, $n\geq 0$. Since $f_n\leq f$, for each $n\geq 0$, we have 

$$
\forall \ n\geq 0, \ (f_n \neq 0) \subset (f \neq 0). \ \  (F1)
$$

\bigskip \noindent S Hence $m(f \neq 0)=0$ implies that $m(f_n \neq 0)=0$ for all $n\geq 0$, that is $f_n=0 \ a.e.$ for all $n\geq 0$. By Property (P2) of Exercise 3, we have $\int f_n \ dm =0$ for all $n\geq 0$. It comes, by definition that

$$
\int f \ dm=0.
$$

\bigskip \noindent S (b) Suppose that $\int f \ dm=0$. We may see that

$$
(f>0)= \bigcup_{n\geq 0} (f_n>0). \ (F2)
$$  

\bigskip \noindent S Indeed, let  $\omega \in (f>0)$, that is $f(\omega)>0$. Since $f_n(\omega)\uparrow f(\omega)$, as $n \uparrow +\infty$, there exits some $n_0$ such that $f_n(^\omega)>0$. This proves the direct inclusion in Formula (F2). The indirect one comes from Formula (F1).\\

\noindent By Formula (LF2), and by the non-negativity of the integrals used here, $\int f \ dm=0$ implies that $\int f_n \ dm=0$ for all $n\geq 0$, which implies, by Property (P2) of Exercise 3, that
$f_n=0 \ a.e$ for all $n\geq 0$. By Formula (F2), we finally get that $f=0 \ a.e.$\\

\noindent We will a simpler ans nicer method to prove this in the general case.\\

\bigskip \noindent \textbf{(P3)} Let $A$ be a null set. We also have

\begin{equation*}
1_A f=\lim_{n\uparrow +\infty}  1_Af_n \ \ and  \ \ \int 1_A f \ dm=\lim_{n\uparrow +\infty}  \int 1_A f_n \ dm,
\end{equation*}

\bigskip \noindent S by the fact that the $1_A f_n$'s form a non-decreasing sequence of non-negative functions in $\mathcal{E}_{+}$, converging to $1_A$ and by definition of the integral of non-negative measurable function. Since by Property (P3) of Exercise 3, all the 
$\int 1_A f_n \ dm$'s are zeros, we get that

$$
\int_A f \ dm=0. \ \square
$$

\bigskip \noindent \textbf{(P4)} We have

\begin{equation*}
(af+bg)=\lim_{n\uparrow +\infty}  (a f_n+b g_n) \ \ and  \ \ \int (a f+b g) \ dm=\lim_{n\uparrow +\infty}  \int (a f_n+b g_n) \ dm,
\end{equation*}

\bigskip \noindent S by the fact that the $(a f_n + b g_n)$'s form a non-decreasing sequence of non-negative functions in $\mathcal{E}_{+}$, converging to $af+bg$ and by definition of the integral of non-negative measurable function. Since by Property (P4) of Exercise 3, we get

$$
\int (a f+b g) \ dm=a \lim_{n\uparrow +\infty}  \int f_n \ dm + b \lim_{n\uparrow +\infty}  \int g_n \ dm=a \int f \ dm + b \int g \ dm. \ \square
$$

\bigskip \noindent \textbf{(P5)} This is prove prove (P4) exactly as in Exercise 3.

\bigskip \noindent \textbf{(P6)} Here again we gave

\begin{equation*}
1_{A+B} f=\lim_{n\uparrow +\infty}  1_{A+B} f_n \ \ and  \ \ \int 1_{A+B} f \ dm=\lim_{n\uparrow +\infty}  \int 1_{A+B} f_n \ dm,
\end{equation*}

\bigskip \noindent S by the fact that the $1_{A+B} f_n$'s form a non-decreasing sequence of non-negative functions in $\mathcal{E}_{+}$, converging to $1_{A+B} f$ and by definition of the integral of non-negative measurable function. Since by Property (P4) and (P6) of Exercise 3, we get

$$
\int 1_{A+B} f \ dm=\lim_{n\uparrow +\infty}  \int 1_{A} f_n \ dm + \lim_{n\uparrow +\infty}  \int 1_{B} f_n \ dm=\int_A f \ dm + \int_B f \ dm.
$$

\bigskip \noindent S All  the properties (P1)-(P5) are established. $\blacksquare$\\

\newpage 
\noindent \textbf{Part III - General functions}.\\

\noindent Hereafter, $\alpha$ and $\beta$ are two real numbers,  $f$ and $g$ are two measurable functions, $A$ and $B$ are disjoint measurable sets.\\

\bigskip \noindent \noindent \textbf{Exercise 6}. \label{exercise06_sol_doc05-05} (Positive and negative parts).\\

\noindent Question (a) Show that the unique solution of the functional system

$$
f=h-g, \ h\geq 0, \ g \geq 0, hg=0, \ \ (FS) 
$$
 
\bigskip \noindent S is given by $h=f^+=max(f,0)$ and $g=f^-=max(f,0)$. Besides, we have

$$
|f|=f^+ \ \ f^-. \ (AV)
$$

\bigskip \noindent S and

$$
(f+g)^+ \leq f^+ + g^+ \ \ and  \ \ (f+g)^- \leq f^- + g^-. \ \  (BO)
$$

\bigskip \noindent S and

$$
(f \leq g ) \Leftrightarrow (f^+ \leq g^+ \ and \ f^- \geq g^-), \ \ (COM) 
$$

\noindent Question (b) Show that we have for any real constant $c$ :

$$
If \ c>0, \ then \ (cf)^+=cf^+ \ and \ (cf)^-=cf^-, \ \ (MPN1)
$$ 

\bigskip \noindent S and

$$
If \ c<0, \ then \ (cf)^+=-cf^- \ and \ (cf)^-=-c f^+, \ \ (MPN2)
$$ 

\bigskip \noindent S Finally, if $f$ is finite $a.e.$, and if $c=0$, $cf=(cf)^+=(cf)^-=0 $, \textit{a.e}\\

\bigskip \noindent Question (c) Suppose that $f$ is finite and $g$ is a function of constant sign. Show that
 
$$
If \ g\geq 0, \ then \ (gf)^+=gf^+ \ and \ (gf)^-=gf^-, \ \ (MPNF1)
$$ 

\bigskip \noindent S and

$$
If \ g \leq <0, \ then \ (gf)^+=-gf^- \ and \ (gf)^-=-g f^+, \ \ (MPNF2)
$$ 

\bigskip \noindent \textbf{Solutions}.\\

\noindent Question (a).\\

\noindent (a1) Let us begin to show that $f^+$ and $f^-$ are solutions of the system (SF). Since we work on functions, we have to compare their graphs by establishing the relations for
all $\omega \in \Omega$. It is enough to work with three cases.\\

\noindent If $f(\omega)=0$, we have $f^+(\omega)=f^-(\omega)=0$ and (SF) holds for $h=f^+$ and $g=f^-$ and (AV) also holds.\\

\noindent If $f(\omega)>0$, we have $f^+(\omega)=f(\omega)$ and $f^-(\omega)=0$. You may easily check that (SF) holds for $h=f^+$ and $g=f^-$ and (AV) also holds.\\

\noindent If $f(\omega)<0$, we have $f^+(\omega)=0$ and $f^-(\omega)=-f(\Omega)$. You may easily check that (SF) holds for $h=f^+$ and $g=f^-$ and (AV) also holds.\\

\bigskip \noindent To prove Formulas (BO), consider the four cases :\\

\noindent \textit{Case 1} : $f(\omega)\geq 0$ and $g(\omega)\geq 0$ : We have

$$
\biggr((f(\omega)+g(\omega))^+=f(\omega)+g(\omega), \ f^+(\omega)=f(\omega), \ g^+(\omega)=g(\omega) \biggr) 
$$

\bigskip \noindent S and 

$$
\biggr( (f(\omega)+g(\omega))^-=0, \ f^-(\omega)=0, \  g^-(\omega)=0    \biggr).
$$

\bigskip \noindent S By looking at the details in the each of the couples of big parentheses, we easily see that the relations in Formula (BO) hold for these $\omega$'s.\\

\noindent \textit{Case 2} : $f(\omega)\leq 0$ and $g(\omega) \leq 0$. We have

$$
\biggr((f(\omega)+g(\omega))^+=0, \ f^+(\omega)=0, \ g^+(\omega)=0 \biggr)
$$

\bigskip \noindent S and

$$
\biggr( (f(\omega)+g(\omega))^-=-f(\omega)-g(\omega), \ f^-(\omega)=-f(\omega), \  g^-(\omega)=-g(\omega)    \biggr). 
$$

\bigskip \noindent S We conclude similarly.\\

\noindent \textit{Case 3 }: One of the $f(\omega)$ and $g(\omega)$ is non-negative, the other non-positive, say $f(\omega)\geq 0$ and $g(\omega) \leq 0$. Here, we have to consider two sub-cases.\\ 

\noindent \textit{Sub-case 31} : $f(\omega)+g(\omega)\geq 0$. We have

$$
\biggr((f(\omega)+g(\omega))^+=f(\omega)+g(\omega), \ f^+(\omega)=f(\omega), \ g^+(\omega)=0 \biggr)
$$

\bigskip \noindent S and

$$
\biggr( (f(\omega)+g(\omega))^-=0, \ f^-(\omega)=0, \  g^-(\omega)=-g(\omega)    \biggr).
$$

\bigskip \noindent S The conclusion is also easily obtained. But remark that form the details in the first couple of big parentheses, $((f(\omega)+g(\omega))^+=f(\omega)+g(\omega)$ is less that $f(\omega)$ since 
$g(\omega)$ is non-positive and $f^+(\omega)+f^-(\omega)=f(\omega)$, since $f^-(\omega)=0$.\\

\noindent \textit{Sub-case 32} : $f(\omega)+g(\omega)< 0$. We have

$$
\biggr((f(\omega)+g(\omega))^+=0, \ f^+(\omega)=f(\omega), \ g^+(\omega)=0 \biggr) 
$$

\bigskip \noindent S and

$$
\biggr( (f(\omega)+g(\omega))^-=(-g(\omega))-f(\omega), \ f^-(\omega)=0, \  g^-(\omega)=-g(\omega)    \biggr).
$$

\bigskip \noindent S We  conclude as in the case 31 by saying that the second couple parentheses leads to : $(f(\omega)+g(\omega))^-=(-g(\omega))-f(\omega)$ is less that $-g(\omega)$ $-f(\omega)$ is non-positive and $f^-(\omega)+f^-(\omega)=-g(\omega)$ since $f^-(\omega)=0$.\\

\bigskip \noindent The prove of Formula (COM) is on the fact that the function $max(.,.)$ is coordinate-wise non-decreasing. Thus $f \leq g$ implies that $max(f,0) \leq max(g,0)$, that is $f^+ \leq g^+$. we also have $f \leq g$ implies $-f \geq -g$ which leads to$max(f,0) \geq max(g,0)$, that is $f^- \geq g^-$. The reverse implication is straightforward.\\

\bigskip \noindent (a2) Next suppose that $h$ and $g$ satisfy (SF).\\

\noindent Suppose that $f(\omega)=0$. We have $h(\omega)-g(\omega)=0$ and $h(\omega)g(\omega)=0$. Hence, if $h(\omega)$ is not zero, so is $g(\omega)$ and  
$f(\omega)=h(\omega)$ would not be not null. Similarly, if $g(\omega)$ is not zero, so is $h(\omega)$ and we would have 
$f(\omega)=-g(\omega)$ not null. Thus, we necessarily have

$$
h(\omega)=0  \ \ and \ \ g(\omega)=0.
$$

\bigskip \noindent S Suppose that $f(\omega)>0$. Since we have $h(\omega)g(\omega)=0$, one of $h(\omega)$ and $g(\omega)$ is zero. And in cannot be $h(\omega)$, since we would have $f(\omega)<0$. Thus we have
$$
h(\omega)=f(\omega) \ \ and \ \ g(\omega)=0.
$$

\bigskip \noindent S Suppose that $f(\omega)<0$. Since we have $h(\omega)g(\omega)=0$, one of $h(\omega)$ and $g(\omega)$ is zero. And in cannot be $h(\omega)$, since we would have $f(\omega)>0$. Thus we have
$$
h(\omega)=0 \ \ and \ \ g(\omega)=-f(\omega).
$$

\bigskip \noindent S By comparing with the results of part (a1) of this question, we may see that $h$ and $f^+$ on one side and $g$ and $f-$ have the same values in the three possible regions $f(\omega)=0$, $f(\omega)<0$ and $f(\omega)>0$. S we have the equality of theirs graphs, that is

$$
h=f^+ \ \ and \ \ g=f^-.
$$

\bigskip \noindent Question (b).\\

\noindent For $c>0$, we have $cf=cf^+ \ - \ cf^-$. Now $h=cf^+$ and $g=cf^-$ are non-negative and $hg=0$. We have the conclusion of formula (MPN1) by applying Question (a).\\

\noindent For $c<0$, we have $cf=cf^+ \ - \ cf^-=(-cf^-) - (-cf^+)$. Now $h=-cf^-$ and $g=-cf^+$ are non-negative and $hg=0$. We have the conclusion of formula (MPN2) by applying Question (a).\\

\noindent For $c=0$, we may have a determination problem for $\omega$'s such that $|f(\omega)|=+\infty$. We require at least that $f$ is finite $a.e$ and the sequel is obvious.

\bigskip Question (C). Proceed exactly as in Question (a).\\

\bigskip \noindent \textbf{Exercise 7}. \label{exercise07_sol_doc05-05} Show the following points.\\

\noindent Question (a)  

$$
\int |f| dm = \int f^+ \ dm + \int f^- \ dm.
$$ 

\bigskip \noindent S Question (b) $f$ is integrable if and only if $|f|$ is.\\

\noindent Question (c) If $\int f \ dm$ exists, then

$$
\left| \int  f \ dm \right| \leq \int |f| \ dm.
$$

\bigskip \noindent S Question (d) Show that if $|f|\geq g$ and $g$ is integrable, then $f$ is.\\

\bigskip \noindent \textbf{Solutions}.\\

\noindent Question (a). This is a direct implication of Property (P4) of Exercise 3 concerning some of non-negative measurable functions and Formula (AV) of Exercise 6.\\

\noindent Question (b). By definition, $f$ is integrable if and only both $\int f^+ \ dm$ and $\int f^- \ dm$ are finite. Since these two numbers are non-negative, they are both finite if and only if their sum is finite. So, $f$ is integrable if and only if

$$
\int f^+ \ dm +  \int f^- \ dm= \int |f| \ dm <+\infty.
$$

\bigskip \noindent S Question (c). Suppose that we have $|f|\geq g$. We may apply Property (P5) for non-negative measurable functions to have

$$
\int |f| \ dm \leq \int g \ dm.
$$

\bigskip \noindent S So, if $g$ is integrable, so is $f$, which is equivalent to the integrability of $f$ according to Question (b).

\bigskip \noindent \textbf{Exercise 8}. \label{exercise08_sol_doc05-05} (Markov inequality) Prove that for any $\lambda>0$

$$
m(|f| \geq \lambda) \leq \frac{1}{\lambda} \int |f| \ dm.  
$$

\noindent Hint. Put $A=(|f| \geq \lambda)$. Decompose the integral of $|f|$ over $A$ and $A^c$, drop the part on $A^c$ and use Property (P4) for non-negative measurable functions.\\

\bigskip \noindent \textbf{Solutions}.\\

\noindent Put $A=(|f| \geq \lambda)$. We have

$$
\int |f| \ dm =\int_{A+A^c} |f| \ dm = \int_{A} |f| \ dm + \int_{A^c} |f| \ dm \geq \int_{A} |f| \ dm, 
$$

\bigskip \noindent S since the dropped part $\int_{A^c} |f| \ dm$ is non-negative. By we also have $1_A|f| \geq \lambda 1_A$. We may use Property (P5) for non-negative measurable functions to get

$$
\int |f| \ dm  \geq \int_{A} |f| \ dm \geq \int_{A} \lambda \ dm= \int \lambda 1_A\ dm=\lambda m(A). 
$$

\bigskip \noindent S By replacing $A$ by its value, we get

$$
\lambda m(|f| \geq \lambda) \geq \int |f| \ dm.
$$

\bigskip \noindent S which was the target. $\square$\\

\bigskip \noindent \textbf{Exercise 9}. \label{exercise09_sol_doc05-05} Use the Markov inequality to show that an integrable function is finite $a.e.$.\\

\noindent Hint. Use the following result of \textit{Exercise 7 of Doc 04-02 in Chapter \ref{04_measures}, page \pageref{doc04-02}} : $(f \ infinite)=(|f|=+\infty)=\bigcap_{k\geq 1} (|f|\geq k)$. Remark that the $(|f|\geq k)$'s form a non-increasing sequence in $k$. Apply the Markov inequality to conclude.\\

\bigskip \noindent \textbf{Solutions}.\\

\noindent From \textit{Exercise 7 of Doc 04-02 in Chapter \ref{04_measures}, page \pageref{doc04-02}}, we have

$$
(f \ infinite)=(|f|=+\infty)=\bigcap_{k\geq 1} (|f|\geq k).
$$

\bigskip \noindent S The $(|f|\geq k)$'s form a non-increasing sequence in $k$, since for any $k\geq 1$, $(|f|\geq k+1) \subset (|f|\geq k)$. Now by Markov's inequality, we have for all $k\geq 1$,

$$
m(|f|\geq k) \geq \frac{1}{k} \inf |f| \ dm. \ \ (MA)
$$

\bigskip \noindent S So, if $f$ is integrable, the numbers $m(|f|\geq k)$ are finite. Thus, by the continuity of the measure $m$, and since $(|f|\geq k)\downarrow (f \ infinite)$ as $k \uparrow +\infty$, we have

$$
m(|f|\geq k)\downarrow m(f \ infinite) \ as k\uparrow +\infty.
$$

\bigskip \noindent S Formula (MA) implies that $m(|f|\geq k)\downarrow 0 \ as k\uparrow +\infty$. We conclude that

$$
m(f \ infinite)=0,
$$

\bigskip \noindent S which was the target. $\square$\\

\bigskip \noindent \textbf{Exercise 10}. \label{exercise10_sol_doc05-05} Let 

\noindent Question (a) Suppose $\int f \ dm$ exists and let $c$ be a non-zero real number. Show that

$$
\int (cf) \ dm= c \int f \ dm.
$$

\bigskip \noindent Question (b) By taking $c=-1$ in question (a), tell simply why Properties (P2), (P4), (P4), (P5), (P6) are valid for non-positive functions, in other words for function of a constant sign.\\

\bigskip \noindent Question (c) Suppose that $f$ is an arbitrary measurable function and $A$ is a $m$-nul. Show that  $\int_{A} f \ dm$ exists and we have

$$
\int_{A} f \ dm=0.
$$

\bigskip \noindent Question (d) Suppose that $f$ exists and let $A$ and $B$ be two disjoint measurable sets. Show that we have

$$
\int_{A+B} f \ dm=\int_{A} f \ dm + \int_{B} f \ dm.
$$

\bigskip  \noindent (e). Let $f$ and $g$ be two measurable functions such that $f=g$ a.e. Show that their integrals exist or not simultaneously and, of they exist, we have
 
$$
\int f \ dm = \int g \ dm= \int_M f \ dm =\int_M g \ dm
$$

\bigskip \noindent S where $M=(f=g)$.\\

\noindent (f) Suppose that $f$ and $g$ are integrable. Show that 
$$
\int (f+g) \ dm=\int f \ dm + \int g \ dm.
$$

\bigskip \noindent Hint : Consider the measurable decomposition of $\Omega$ :

$$
\Omega =(f\geq 0, g \geq 0)+(f < 0, g \geq 0)+(f\geq 0, g <0)+(f < 0, g < 0)=A+B+C+D
$$

\bigskip \noindent Each of these sets is, itself, decomposable into two set :
$$
A=A \cap (f+g \geq 0) + A \cap (f+g < 0)=A_1 + A _2.
$$

\bigskip \noindent define likewise $B_i$, $C_i$, $D_i$, $i=1,2$. Show that 
$$
\int_E (f+g) \ dm=\int_E f  dm + \int_E f \ dm.
$$

\bigskip \noindent when $E$ is any of the $A_i$, $B_i$, $C_i$, $D_i$, $i=1,2$ and conclude.\\

\noindent Hint. The key tool is to manage to have a decomposition of one of elements of three terms
$\int_E (f+g) \ dm$, $\int_E f  dm$, $\int_E f \ dm$ into the two others such that these two ones have the same constant signs.\\

\noindent Here, prove this only for one of them.\\

\bigskip  \noindent (g) Suppose $\int f \ dm$, $\int g \ dm$ exist and $f+g$ is defined and $\int f+g \ dm$ exist. Show that 

$$
\int (f+g) \ dm=\int f \ dm + \int f \ dm.
$$

\bigskip \noindent S \textbf{Recommendation} This extension of Question (e) is rather technical. We recommend that you read the solution in the Appendix below.\\

\bigskip  \noindent (h) Suppose that $\int f \ dm$, $\int g \ dm$ exist and $f\leq g$. Show that

$$
\int f \ dm \geq \int g \ dm. \ \ (INC)
$$

\bigskip \noindent \textbf{Solutions}.\\

\noindent Question (a).\\

\noindent Suppose $c>0$. If $f$ is quasi-integrable, one of $\int f^+ \ dm$ and $\int f^- \ dm$ is finite. So one of $c \int f^+ \ dm$ and $ c\int f^- \ dm$ is finite.\\

\noindent By Exercise 6, Question (b), and by Property (P4) for non-negative measurable functions, we have the following.\\

$$
\int (cf)^+ \ dm =c \int f^+\ dm \ and \ \int (cf)^- \ dm =c \int f^- \ dm.
$$

\bigskip \noindent S By combining these two points, we say that $fc$ is also quasi-integrable. Also, if $\int f^+\ dm - \int f^- \ dm$ makes sense, we get

$$
\int (cf)^+ \ dm - \int (cf)^- \ dm =c \left( \int f^+\ dm - \int f^- \ dm \right),
$$

\bigskip \noindent S which $\int cf \ dm = c\int f \ dm$.\\

\noindent Suppose $c<0$. We have the same remarks on the quasi-integrability of $cf$ and we have by the same arguments,

$$
\int (cf)^+ \ dm =-c \int f^-\ dm \ and \ \int (cf)^- \ dm =- c \int f^+ \ dm
$$

\bigskip \noindent S and next

$$
\int cf \ dm =\int (cf)^+ \ dm - \int (cf)^- \ dm =-c  \int f^-\ dm - (-c \int f^+ \ dm),
$$

$$
=-c \left( \int f^-\ dm - \int f^+ \ dm \right)= c \left( \int f^+\ dm - \int f^- \ dm \right).
$$

\bigskip \noindent S And we get the same conclusion.\\

\noindent Question (b) Let $f$ be a measurable non-positive function, we get by Question (a)

$$
\int f \ dm = - \int (-f) \ dm.
$$

\bigskip \noindent S From this formula, it is obvious that all Properties (P2) to (P6), which already hold for measurable non-negative functions are automatically
transferred to measurable non-positive functions through the non-negative through $(-f)$.\\

\noindent Question (c). By Exercise 6, Question (c), we have

$$
(1_Af)^+ \ dm = 1_A f^+ \ \  and \ (1_Af)^- \ dm = 1_A f^-.
$$

\bigskip \noindent S We get by definition and by Property (P3) for non-negative measurable functions,

$$
\int (1_A f) \ dm=\int (1_A f)^+ \ dm - \int (1_A f)^- \ dm=  \int_A f^+ \ dm - \int_A f^- \ dm=0. 
$$

\bigskip \noindent S Question (d). Since $f$ is quasi-integrable, one of $\int f^+ \ dm$ and $\int f^- \ dm$ is finite. Suppose here that is $\int f^+ \ dm$ is finite. By Property (P5) for non-negative measurable functions, we also have that $\int_{A+B} f^+ \ dm$, $\int_A f^+ \ dm$ and $\int_B f^+ \ dm$ are finite, so that $1_{A+B} f$, $1_A f$ and $1_B f$ are quasi-integrable.\\

\noindent By Exercise 6, Question (c), we have

$$
(1_{A+B}f)^+  = 1_{A+B} f^{+}=1_{A} f^{+} + 1_{B} f^+ \ \  and \ (1_{A+B}f)^-  = 1_{A+B} f^{-}=1_{A} f^{-} + 1_{B} f^-.
$$

\bigskip \noindent S By Property (P4) and (P6) for non-negative measurable functions, we have

$$
\int (1_{A+B}f)^+ \ dm =\int_A f^+ \ dm + \int_B f^+ \ dm \ \ and \ \ \int (1_{A+B}f)^- \ dm =\int_A f^- \ dm + \int_B f^- \ dm.
$$

\bigskip \noindent S Since the numbers $\int_A f^+ \ dm$ and $\int_B f^+ \ dm$, the following operation 

$$
\int (1_{A+B}f) \ dm - \int (1_{A+B}f)^- \ dm =\left( \int_A f^+ \ dm - \int_A f^- \ dm\right) + \left(\int_B f^+ \ dm - \int_B f^- \ dm\right),
$$

\bigskip \noindent S which is

$$
\int_{A+B} f \ dm = \int_{A} f \ dm  + \int_{B} f \ dm. \ \square.
$$

\bigskip \noindent Question (e). Put $N=(f \neq g)$. If $f=g$ $a.e$, we have $m(f \neq g)=0$. Since $f=g$ on $N^c$, we may denote

$$
A=\int_{N^c} f^+ \ dm=\int_{N^c} g^+ \ dm, \ C=\int_{N^c} f^- \ dm=\int_{N^c} g^- \ dm.
$$

\bigskip \noindent S By Question (c), we have

$$
\int f^+ \ dm=\int_N f^+ \ dm + \int_{N^c} f^+ \ dm=\int_{N^c} f^+ \ dm=A.
$$

\bigskip \noindent S By doing the same for $f^-$, $g^+$ and $g^-$, we get

$$
\int f^+ \ dm=A, \ \int ^+ \ dm=A, \int f^- \ dm=B \ and \int f^- \ dm=b.
$$

\bigskip \noindent S This clearly implies that $f$ and $g$ are quasi-integral at the same time, or are not at the same time. So, if they are quasi-integrable, we have

$$
\int f \ dm = A - B = \int g dm
$$

\bigskip \noindent S and clearly

$$
A - B=\int_{N^c} f^+ \ dm - \int_{N^c} f^- \ dm =\int_{N^c} f \ dm.
$$

\bigskip \noindent S We also have

$$
A + B=\int_{N^c} g^+ \ dm - \int_{N^c} g^- \ dm =\int_{N^c} g \ dm.
$$

\bigskip \noindent S We conclude that

$$
\int f \ dm = \int g \ dm = \int f \ dm = \int_{N^c}  f \ dm = \int_{N^c} g \ dm
$$

\bigskip \noindent S and

$$
\int f^- \ dm=\int_N f^+ \ dm + \int_{N^c} f^+ \ dm=\int_{N^c} f^+ \ dm.
$$

\bigskip \noindent Question (f).\\

\noindent Suppose that $f$ and $g$ are integrable. So $f$ and $g$ are a.e finite and then $f+g$ is a.e. finite -(say, $f$ and $g$ are equal on $M$ with a $m$-null complement). By Question (c), it will enough to prove that

$$
\int_{M} f+g \ dm = \int_{M} f \ dm + \int_{M} g \ dm.
$$ 

\bigskip \noindent S In the sequel, we suppose that $N^c=\Omega$ since this will not effect the solution. Define

$$
\Omega =(f\geq 0, g \geq 0)+(f < 0, g \geq 0)+(f\geq 0, g <0)+(f < 0, g < 0)=A+B+C+D.
$$

\bigskip \noindent S Next, let us decompose $A$ into
$$
A=A \cap (f+g \geq 0) + A \cap (f+g < 0)=A_1 + A _2.
$$

\bigskip \noindent S and define likewise $B_i$, $C_i$, $D_i$, $i=1,2$.\\

\noindent Based on Question (d) and since any $h \in \{f, g, f+g\}$ is quasi-integrable, we have

\begin{eqnarray*}
&& \int h \ dm=\int_{A_1} h \ dm+\int_{A_2} h \ dm+\int_{2_1} h \ dm+\int_{B_2} h \ dm\\
&+&\int_{C_1} h \ dm+\int_{C_2} h \ dm+\int_{D_1} h \ dm+\int_{A_2} h \ dm.
\end{eqnarray*}

\bigskip \noindent S So, to prove that $\int f+g \ dm=\int f \ dm + \int f \ dm$, it will be enough to show that 
$$
\int_E (f+g) \ dm=\int_E f  dm + \int_E f \ dm. \ \ (SUM)
$$

\bigskip \noindent S when $E$ is any of the $A_i$, $B_i$, $C_i$, $D_i$, $i=1,2$ and conclude.\\

\noindent We are going to do this for some $E$'s and give the method to complete all the cases.\\

\noindent Let us pick $E=C_2=((f\geq 0, g <0, f+g<0)$. We have the 

$$
1_E (f+g) = 1_E (f+g) + 1_E g.
$$

\bigskip \noindent S Since  all the function are non-positive, we have Formula (SUM).\\

\noindent Let us pick $E=B_2=(f < 0, g \geq 0, f+g\geq 0)$. We use 

$$
1_E (f+g) +  1_E (-f) = 1_E g.
$$

\bigskip \noindent S Since all the function are non-positive, we have

$$
\int_E (f+g) \ dm - \int_E f \ dm = \int_E g \ dm.
$$

\bigskip \noindent S Since $\int_E f \ dm$ is finite, we may add it both members of the equality to get Formula (SUM).\\

\noindent \textbf{(APC) Here are all the cases}.\\

\noindent (1) On $A_1=(f\geq 0, \ g \geq 0, \ f+g \geq 0)$, we use a sum non-negative functions :

$$
1_E (f+g) \ dm  =  1_E  f + 1_E g.
$$

\bigskip \noindent S (2) $A_2=(f\geq 0, \ g \geq 0, \ f+g<0)$ is empty. We drop it.\\

\noindent (3) On $B_1=(f < 0, \ g \geq 0, \ f+g \geq 0)$, we use a sum non-negative functions :

$$
1_E -(f+g)  + 1_E  (-f)=   1_E g.
$$

\bigskip \noindent S (4) On $B_2=(f < 0, \ g \geq 0, f+g<0)$, we use a sum non-positive functions :

$$
1_E -(f+g)   + 1_E  (-f)=   + 1_E g.
$$

\bigskip \noindent S (5) On $C_1=(f\geq 0, \ g <0, \ f+g \geq 0)$, we use a sum non-positive functions :

$$
1_E -(f+g)  + 1_E (-f) =  1_E  (-g).
$$

\bigskip \noindent S (6) On $C_2=(f\geq 0, \ g <0, f+g<0)$, we use a sum non-positive functions :

$$
1_E (f+g) + 1_E  (-f)=  1_E g.
$$

\bigskip \noindent S (7) $D_1=(f < 0, \ g < 0, \ f+g \geq 0)$ is empty. We drop it.\\

\noindent (8) On $D_2=(f < 0, \ g < 0, f+g<0)$, we use a sum negative functions :

$$
1_E (f+g) \ dm  =  1_E  f + 1_E g.
$$

\bigskip \noindent S In any of these case, the quantity to move from one member to an other, if it exists, is finite. Each case leads to Formula (SUM).\\

\bigskip \noindent Question (g). As suggested, read the solution in the Appendix below.\\

\bigskip \noindent Question (h). Let us discuss on the finiteness of $\int f \ dm$.\\

\noindent Case $\int f \ dm$ finite. Hence $f$ finite \textit{a.e} and we have $g=f + (g-f)$ and we apply Question (g) to get

$$
\int g \ dm = \int f \ dm + \int g-f \ dm \geq \int f \ dm,
$$

\bigskip \noindent S since $\int g-f \ dm \geq 0$ based on the fact that $g-f \geq 0$ \textit{a.e.}\\

\noindent Case $\int f \ dm$ infinite. We have :\\

\noindent \textit{Either} $\int f^+ \ dm=+\infty$and $\int f^- \ dm$ is finite and by Formula (COMP) in ,
$$
f^+ \leq g^+ \Rightarrow +\infty=\int f^+ \ dm \leq \int g^+ \ dm,
$$

\noindent that is $\int g^+ \ dm=+\infty$ and next, $\int g^- \ dm$ is finite since $g$ is quasi-integrable. Formula (INC) holds since the the members are both $\infty$.

\noindent \textit{Or}  $\int f^- \ dm$ is finite and $\int f^- \ dm=+\infty$ finite. Hence $\int f \ dm=-\infty$ is always less than $\int g  \ dm$ which exists in $\overline{\mathcal{R}}$.\\

\noindent By putting together all these cases, we get Formula (INC).\\

\bigskip \noindent \textbf{Exercise 11}. \label{exercise11_sol_doc05-05}  Combines all the results of the different questions of Exercise 9 to establish the final properties of the integrals given in the section \textbf{II}, 
(\textit{Main properties}) in \textit{Doc 05-01 of this chapter} (page \pageref{doc05-01}). Use sentences only.\\

\noindent \textbf{From now, you only use these main properties which includes all the partial ones you gradually proved to reach them}.\\

\bigskip \noindent\textbf{Solutions}.\\

\noindent (1) Formulas in Point (05.04) of Doc 05-01 (page \pageref{doc05-01}). They correspond to Questions (a), (d), (f) and (g) in \textit{Exercise 10} above.\\

\noindent (2) Formulas in Point (05.05) of Doc 05-01 (page \pageref{doc05-01}). Formula (O1) comes from Question (h) in \textit{Exercise 9}. Formula (O2) comes from Question (e) in \textit{Exercise 10}.\\

\noindent (2) Formulas in Point (05.06) of Doc 05-01 (page \pageref{doc05-01}). Formula (I1) and (I2) comes Exercise 7. To get Formula (I3), combine Question (a), (f) and (g) and observe that if $f$ and $g$ are integrable, then by Exercise 10, $af+bg$ is defined \textit{a.e} and $a \int \ dm f +b \int g \int \ dm f$ makes sense. If $a$ and $b$ are non-zero, then (I3) holds in virtue of Question(g). Since $f$ and $g$ are finite \textit{a.e}, the scalars $a$ and $b$ may take the value zero and Formula (I2) remains true.

\bigskip \noindent \textbf{Exercise 12}. \label{exercise12_sol_doc05-05}  (Monotone Convergence Theorem of series).\\

\noindent Let $f : \mathbb{Z} \mapsto \overline{\mathbb{R}}_{+}$ be non-negative function and let $f_p : \mathbb{Z} \mapsto \overline{\mathbb{R}}_{+}$,  $p\geq 1$, be a sequence of non-negative functions increasing to $f$ in the following sense :

\begin{equation}
\forall (n\in \mathbb{Z}), \ \ 0\leq f_{p}(n)\uparrow f(n)\text{ as }p\uparrow\infty.   \label{comptage01}
\end{equation}

\bigskip
\noindent Then, we have 
\begin{equation*}
\sum_{n\in \mathbb{Z}}f_{p}(n)\uparrow \sum_{n\in \mathbb{Z}}f(n).
\end{equation*}

\bigskip \noindent \textbf{Solution}.\\

\noindent Let us study two cases.\\

\noindent \textbf{Case 1}. There exists $n_{0} \in \mathbb{Z}$ such that $f(n_{0})=+\infty$. Then  
\begin{equation*}
\sum_{n\in \mathbb{Z}}f(n)=+\infty.
\end{equation*}

\bigskip \noindent S Since $f_{p}(n_{0})\uparrow f(n_{0})=+\infty$, then for any $M>0$, there exits $P_0\geq 1$, such that  
\begin{equation*}
p>P_0 \Rightarrow f_{p}(n_{0})>M.
\end{equation*}

\bigskip \noindent S Thus, we have 
\begin{equation*}
\forall (M>0),\exists (P_0 \geq 1),\text{ }p>P_0 \Rightarrow \sum_{n\in \mathbb{Z}%
}f_{p}(n)>M.
\end{equation*}

\bigskip \noindent S By letting $p\rightarrow \infty$ in the latter formula, we have : for all $M>0$,

 \begin{equation*}
\sum_{n\in \mathbb{Z}}f_{p}(n) \geq M.
\end{equation*}

\bigskip S By letting $\rightarrow +\infty$, we get

\begin{equation*}
\sum_{n\in \mathbb{Z}}f_{p}(n)\uparrow +\infty =\sum_{n\in \mathbb{Z}}f(n).
\end{equation*}

\bigskip \noindent S The proof is complete in the case 1.\\

\noindent \textbf{Case 2}. Suppose that que $f(n)<\infty$ for any $n\in \mathbb{Z}$. By \ref{comptage01}, it is clear that 
\begin{equation*}
\sum_{n\in \mathbb{Z}}f_{p}(n)\leq \sum_{n\in \mathbb{Z}}f(n).
\end{equation*}

\bigskip \noindent S The left-hand member is non-decreasing in $p$. So its limit is a monotone limit and it always exists in $\overline{\mathbb{R}}$ and we have  
\begin{equation*}
\lim_{p\rightarrow +\infty}\sum_{n\in \mathbb{Z}}f_{p}(n)\leq \sum_{n\in \mathbb{Z}}f(n).
\end{equation*}

\bigskip \noindent S Now, fix an integer $N>1$. For any $\varepsilon >0$, there exists $P_{N}$ such that  
\begin{equation*}
p>P_{N}\Rightarrow (\forall (-N\leq n\leq N),\text{ }f(n)-\varepsilon
/(2N+1)\leq f_{p}(n)\leq f(n)+\varepsilon /(2N+1)
\end{equation*}

\bigskip \noindent S and thus,  
\begin{equation*}
p>P_{N}\Rightarrow \sum_{-N\leq n\leq N}f_{p}(n)\geq \sum_{-N\leq n\leq N}f(n)-\varepsilon.
\end{equation*}

\bigskip \noindent S Thus, we have
\begin{equation*}
p>P_{N}\Rightarrow \sum_{n\in }f_{p}(n)\geq \sum_{-N\leq n\leq N}f_{p}(n)\geq
\sum_{-N\leq n\leq N}f(n)-\varepsilon, 
\end{equation*}

\bigskip \noindent S meaning that for $p>P_{N}$,

\begin{equation*}
\sum_{-N\leq n\leq N}f(n)-\varepsilon \leq \sum_{-N\leq n\leq N}f_{p}(n).
\end{equation*}

\noindent  S and then, for $p>P_{N}$,

\begin{equation*}
\sum_{-N\leq n\leq N}f(n)-\varepsilon \leq \sum_{n\in \mathbb{Z}}f_{p}(n).
\end{equation*}

\bigskip \noindent S By letting  $p\uparrow \infty$ first and next, $N\uparrow \infty$, we get for any $\varepsilon >0,$%
\begin{equation*}
\sum_{n\in \mathbb{Z}}f(n)-\varepsilon \leq \lim_{p\rightarrow +\infty}\sum_{n\in \mathbb{Z}}f_{p}(n).
\end{equation*}

\noindent Finally, by letting $\varepsilon \downarrow 0$, we get 
\begin{equation*}
\sum_{n\in \mathbb{Z}}f(n)\leq \lim_{n}\sum_{n\in \mathbb{Z}}f_{p}(n).
\end{equation*}

\bigskip \noindent S We conclude that  
\begin{equation*}
\sum_{n\in \mathbb{Z}}f(n)=lim_{p\uparrow+\infty} \sum_{n\in \mathbb{Z}}f_{p}(n).
\end{equation*}

\bigskip \noindent S The case 2 is now closed and the proof of the theorem is complete.\\

\newpage
\bigskip \noindent \textbf{Appendix}. Solution of Question (g) of Exercise 10 : If $f+g$ is defined a.e., $\int f \ dm$, $\int g \ dm$ and $\int f+g \ dm$ exist and the addition 
$\int f \ dm + \int g \ dm$ makes sense, then

$$
\int f+g \ dm = \int f  dm  + \int g \ dm. \ (AD)
$$

\bigskip \noindent S Let us consider three cases.\\

\noindent Case 1. $f$ and $g$ both integrable. This case reduces to Question (e) of Exercise 10.\\

\noindent Case 2. One of $f$ and $g$ is integrable. It will be enough to give the solution for $f$ integrable.  Let us do it by supposing that $f$ is integrable. We adopt the same prove in the proof
of Question(c) of Exercise 10. Indeed, in all the eight (8), we only need to move $\int f \ dm$ from one member to another, which possible thanks to its finiteness.

\noindent Case 3. Both $f$ and $g$ have infinite integrals. If the sum $\int f \ dm + \int g \ dm$ makes sense, that means that the two integrals are both $+\infty$ or both $-\infty$.

$$
\int f^+ \ dm <+\infty, \ \ \int f^- \ dm =+\infty, \ \ \int g^+ \ dm <+\infty \ \ and \ \ \int g^- \ dm =+\infty 
$$

\bigskip \noindent S or

$$
\int f^+ \ dm =+\infty, \ \ \int f^- \ dm <+\infty, \ \ \int g^+ \ dm =\infty \ \ and \ \ \int g^- \ dm <+\infty.
$$

\bigskip \noindent S Let us continue through two sub-cases.\\

\noindent First, let us remark that the method used in the solution of Question (f) of Exercise 10 is valid until the statements of the eight decompositions at Line marked (APC). From that point, the only problem concerned the moving of a finite quantity from the left member to the second after we applied the integral operator at each member and used the linearity for functions of same constant signs. But this concerned only the cases (3), (4), (5) and (6).\\

\noindent Here, we are going to use the finiteness of both integrals of $f$ and $g$, or both integrals of negative parts.\\

\noindent \textbf{Case (i) where the positives parts of $f$ and $g$ are finite}.\\

\noindent In the three the following three cases, we may directly apply the explained principles.\\

\noindent (3) On $B_1=(f < 0, \ g \geq 0, \ f+g \geq 0)$, we use a sum non-positive functions :

$$
1_E -(f+g)  + 1_E  -g =   1_E -f; \ (In \ this \ case : \ \int 1_E  -g \ dm = - \int_E  g^+ \ dm \ finite).
$$

\noindent (4) S On $B_2=(f < 0, \ g \geq 0, f+g<0)$, we use a sum non-positive functions :

$$
1_E (f+g)   + 1_E  -g =   1_E -f. \ (In \ this \ case : \ \int 1_E  -g \ dm = - \int_E  g^+ \ dm \ finite).
$$

\bigskip \noindent S (6) On $C_2=(f\geq 0, \ g <0, f+g<0)$, we use a sum non-negative functions :

$$
1_E -(f+g) + 1_E  f =  1_E -g. \ (In \ this \ case : \ \int 1_E  f \ dm =  \int_E  f^+ \ dm \ finite).
$$

\bigskip \noindent S The case (5) is handled as follows.\\:

\noindent (5) On $C_1=(f \geq 0, \ g <0, \ f+g \geq 0)$, we use a sum non-positive functions : 

$$
1_E -(f+g)  + 1_E g =  1_E  -f. 
$$

\bigskip \noindent S Since $\int_E  f^ + \ dm=\int_E  f^+ \ dm$ is finite, we get by the additivity of integrals of functions of same constants signs that $- \int 1_E -(f+g) \ dm$ and 
$\int 1_E g \ dm$ are both finite.\\

\noindent In total, in all these cases, we will be able to move the second term of the first member of the equalities after having applied the integral operator.\\

\noindent \textbf{Case (ii) where the negative parts of $f$ and $g$ are finite}.\\

\noindent Here, we have :\\

\noindent (3) On $B_1=(f < 0, \ g \geq 0, \ f+g \geq 0)$, we use a sum if non-positive functions :

$$
1_E -(f+g)  + 1_E  f =   1_E -g ; \ (In \ this \ case : \ \int 1_E  f \ dm = - \int_E  f^- \ dm \ finite).
$$

\bigskip \noindent S (4) On $B_2=(f < 0, \ g \geq 0, f+g<0)$, we use a sum of non-negative functions :

$$
1_E -(f+g)   + 1_E  -f =  -1_E -g. \ (In \ this \ case : \ \int 1_E  -g \ dm = - \int_E  g^+ \ dm \ finite).
$$

\bigskip \noindent S (5) On $C_1=(f \geq 0, \ g <0, \ f+g \geq 0)$, we use a sum of of non-positive functions :

$$
1_E -(f+g)  + 1_E g =  1_E  (-f). \ (In \ this \ case : \ \int 1_E  -g \ dm = - \int_E  g^+ \ dm \ finite).
$$

\bigskip \noindent S (6) On $C_2=(f\geq 0, \ g <0, f+g<0)$, we use a sum non-positive functions :

$$
1_E -(f+g) + 1_E  f =  1_E -g.
$$

\bigskip \noindent S The conclusion of the treatments of the cases (3), (4) and (6) in case (i) are repeated for (3), (4) and (5) in case (ii). The treatment of (5) in case (i) repeated
for (6) in case (ii). $\blacksquare$

\newpage
\noindent \LARGE \textbf{Doc 05-06  : Integration with the counting measure - Exercises with solutions}. \label{doc05-06}\\
\bigskip
\Large

\bigskip

\bigskip \noindent \textbf{Exercise 1}. \label{exercise01_sol_doc05-06} (Integrating with respect to the counting measure).\\

\noindent Define the counting measure $\nu$ on a countable space$\mathbb{I}$ that is endowed with the discrete sigma-algebra, that is $\mathcal{P}(\mathbb{I})$ by
\begin{equation*}
\nu =\sum_{n\in \mathbb{I}}\delta _{n}
\end{equation*}

\bigskip \noindent and for any subset $A$ of $\mathbb{I}$. 
\begin{equation*}
\nu (A)=\sum_{n\in \mathbb{I}}\delta _{n}(A)=\sum_{n\in \mathbb{I}%
}1_{A}(n)=Card(A)
\end{equation*}

\bigskip \noindent  Any function  
\begin{equation*}
f:\mathbb{I}\mapsto \overline{\mathbb{R}}
\end{equation*}

\bigskip \noindent is measurable. Actually, such a function is a real sequence of real numbers $(f(n))_{n\in \mathbb{I}}$.\\

\noindent Your are going to find how the integral of $f$ with respect to $\nu$ is constructed through the four steps. We are going to see that, for any constant-sign function $f$ or integrabke function, we have

$$
\int f \ d\nu =\sum_{n\in \mathbb{I}}f(n). \ \ (COUNTINTEG)
$$

\bigskip \noindent \textit{Question (a) [Step 1]}. Let $f$ be an indicator function : $f=1_{A}$. Show that Formula (COUNTINTEG) holds.\\
 
\noindent \textit{Question (b) [Step 2]} Let $f$ be  simple and non-negative function, of the form 
\begin{equation*}
f=\sum_{1\leq i\leq k}\alpha _{i}1_{A_{i}}
\end{equation*}

\bigskip \noindent  where the $A_{i}$ are subsets of $\mathbb{I}$ and the $\alpha _{i}$ are finite non-negative numbers. Show also that Formula (COUNTINTEG) holds.\\

\bigskip \noindent \textbf{Question (c) [Step 3]}. Let $f$ be non-negative. Consider a sequence of simple functions $(f_{p})_{p\geq 0}$ increasing to $f$, that is 
\begin{equation}
\forall (n\in \mathbb{N}),0\leq \text{ }f_{p}(n)\uparrow f(n)\text{ as }n\uparrow \infty.   \label{comptage01}
\end{equation}

\bigskip \noindent By combining Question (b), \textit{Exercise 12 in Doc 05-02} and the construction on the integral of a non-negative and measurable function, show that  Formula (COUNTINTEG) 
still holds
 
\bigskip \noindent \textbf{Question (d) [Step 4]} Let $f$ be an arbitrary function from $\mathbb{I}$ to $\overline{\mathbb{R}}$. Give the expressions of the positive and the negative parts of $f$.
Show that if $f$ is $\nu$-integrable, then Formula (COUNTINTEG) holds.\\

\noindent \textbf{Question (e)} Let $f$ be an arbitrary function from $\mathcal{Z}$ to $\overline{\mathbb{R}}$. Show that if $f$ is quasi-integrable, then Formula (COUNTINTEG) holds.

\noindent \textbf{Question (f)}. Propose a conclusion of the exercise.\\

\bigskip \noindent \textbf{Solutions}.\\

\noindent \textbf{Question (a)}.  We have

\begin{equation*}
\int f\text{ }d\nu =\nu (A)=\sum_{n\in \mathbb{I}}1_{A}(n)=\nu
(A)=\sum_{n\in \mathbb{I}}f(n).
\end{equation*}

\bigskip \noindent \textbf{Question (b)}. By applying Question (a), we have

\begin{equation*}
\int fd\nu =\sum_{1\leq i\leq k}\alpha _{i}\nu (A_{i})=\sum_{1\leq i\leq
k}\alpha _{i}\sum_{n\in \mathbb{I}}1_{A_{i}}(n)=\sum_{n\in \mathbb{I}%
}\sum_{1\leq i\leq k}\alpha _{i}1_{A_{i}}(n)
\end{equation*}
\begin{equation*}
=\sum_{n\in \mathbb{I}}f(n).
\end{equation*}

\bigskip \noindent \textbf{Question (c)}. Let $f$ be non-negative.  Consider a sequence of simple functions $(f_{p})_{p\geq 0}$ increasing to $f$, that is 
\begin{equation}
\forall (n\in \mathbb{N}),0\leq \text{ }f_{p}(n)\uparrow f(n)\text{ as }n\uparrow
\infty   \label{comptage01}
\end{equation}

\noindent By Question (b), we have

$$
\inf f_p \ dm = \sum_{n\in \mathbb{I}}f_p(n)
$$

\bigskip \noindent and, by the construction of the integral, we have

$$
\sum_{n\in \mathbb{I}}f_p(n) = \inf f_p \ dm  \nearrow \int f \ dm \ as \ n \ \nearrow +\infty.
$$

\bigskip \noindent By using \textit{Exercise 12 in Doc 05-02}, we get from the formula above that (COUNTINTEG) holds.\\

\noindent \textbf{Question (d)}. We have

\begin{equation*}
f^{+}(n)=\left\{ 
\begin{array}{c}
f(n)\text{ if }f(n)\geq 0 \\ 
0\text{ otherwise}%
\end{array}%
\right. \text{ \ \ \ \ \ \ and \ \ \ \ \ \ \ }f^{-}(n)=\left\{ 
\begin{array}{c}
0\text{ otherwise } \\ 
-f(n)\text{ if }f(n)\leq 0%
\end{array}%
\right..
\end{equation*}

\bigskip \noindent So, by construction, $f$ is integrable if and only if $\int f^+ \ d\nu$ and $\int f^- \ d\nu$ are finite, that is, by the properties of the integral, 
$\int |f| \ d\nu=\int f^- \ d\nu +\int f^- \ d\nu$ is finite. By Question (c), we have

$$
\int f^+ \ d\nu = \sum_{n\in \mathbb{I}}f^+(n)= \sum_{f(n) \geq 0}f(n)= \sum_{f(n) \geq 0} |f(n)| \ \ (S1)
$$

\bigskip \noindent and

$$
\int f^- \ d\nu = \sum_{n\in \mathbb{I}}f^+(n)= \sum_{f(n) \leq 0} (-f(n))= \sum_{f(n) \leq 0} |f(n)|. \ \ (S2)
$$

\bigskip \noindent Hence, $f$ is $\nu$-integrable if and only if

$$
\sum_{f(n) \geq 0} |f(n)| + \sum_{f(n) \geq 0} |f(n)| \ \ (S) 
$$

\bigskip \noindent is finite. You will not have any difficulty to see that the sum in (S) is that of the absolute values of $f(n)$ over $\mathbb{I}$. Finally, $f$ is integrable if and only if

$$
\int |f| d\nu =\sum_{n \in \mathbb{I}} |f(n)| <+\infty 
$$

\bigskip \noindent and thus

\begin{eqnarray*}
\int f d\nu &=& \int f^+ d\nu - \int f^- d\nu\\
&=&\biggr( \sum_{f(n) \geq 0}f(n) \biggr) - \biggr( \sum_{f(n) \leq 0} (-f(n))\biggr)\\
&=&\biggr( \sum_{f(n) \geq 0}f(n) \biggr) + \biggr( \sum_{f(n) \leq 0} f(n)\biggr)\\
&=&\sum_{n \in \mathbb{I}} f(n),\\
\end{eqnarray*}

\bigskip \noindent that is, 

$$
\int f d\nu=\sum_{n \in \mathbb{I}} f(n)
$$

\bigskip \noindent and hence, (COUNTINTEG) holds.\\

\noindent \textbf{Question (e)}. $f$ quasi-integrable means that $\int f^+ \ dm$ or $\int f^- \ dm$ are finite. If both are finite, (COUNTINTEG) holds by Question (d). If one of them, only, if finite, we get by combining Formula (S1) and (S2) in the solution of Question (d) that (COUNTINTEG) holds for both members being positive or negative infinity.\\

\noindent \textbf{Question (f)}. The conclusion is the following.\\

\noindent A series $\sum_{n\in \mathcal{Z}} f(n)$ is the integral $\int f \ d\nu$ of $f$ with respect to the counting measure $\nu$ on 
$\mathcal{Z}$ if $f$ is the sums of positive terms  or the sum of negative terms is finite.\\

\bigskip \noindent \textbf{Exercise 2}. \label{exercise02_sol_doc05-06} Consider the following sequences

$$
f(n)=(-1)^{n} g(n), n\geq 1,
$$

\bigskip \noindent where $g(n)$ is non-increasing and $g(n) \rightarrow 0$ as $n\rightarrow +\infty$. We know from the earlier courses of calculus that 
$\sum_{n \geq 1} f(n)$ exists and is finite.\\

\noindent Now, take $g(n)=1/(n+1)$, $n\geq 0$.\\

\noindent Question : Can you say that $\sum_{n \geq 1} (-1)^{n}/(n+1)$ is the integral of the counting measure on $\mathbb{N}$? Why?

\bigskip \noindent \textbf{Solution}.\\

\noindent Let use the details in the solution of Question (d) in \textit{Exercise 1}. Here $\mathbb{I}=\mathbb{N}$. Formulas (S1) and (S2) in the solution of Question (d) become

$$
\int f^+ \ d\nu =  \sum_{f(n) \geq 0} f(n)=\sum_{n \geq 0} 1/(2n+1)=+\infty
$$

\bigskip \noindent and

$$
\int f^- \ d\nu =  \sum_{f(n) \leq 0} -f(n)=\sum_{n \geq 0} 1/(2n+3)=+\infty.
$$

\bigskip \noindent Hence $g$ does not have an integral with respect to $\nu$. So we cannot write (COUNTINTEG) for $g$.\\

\bigskip \noindent \textbf{Exercise 3}. \label{exercise03_sol_doc05-06}  (Abel's rule) Let $(v_n)_{n\geq 0}$ be a totally bounded sequence of real numbers, that is, there exists a real number $A$ such that for any $n\geq 0$, $m\geq n$,

$$
|v_n + v_{n+1} + \cdots + v_{m}| \leq A.
$$

\bigskip \noindent Let $(\varepsilon_n)_{n\geq 0}$ be a sequence of reals number such that 
$$
\sum_{n\geq 1} |\varepsilon_{n}-\varepsilon_{n-1}|<+\infty \ (A1)
$$

\bigskip \noindent and

$$
\lim_{n\rightarrow +\infty} \varepsilon_{n}=0. \ \ (A2)
$$

\bigskip \noindent Show that the series 

$$
\sum_{n\geq 1} v_{n} \varepsilon_{n}
$$

\bigskip \noindent is convergent.\\

\noindent Hint. Put for $n \leq m$, $R_{n,m}=\sum_{n\leq k \leq m} v_{n} \varepsilon_{n}$ and $V_{n,j}=v_n +  \cdots + v_{n+j}$, $j\geq 0$. Show that : 
\begin{eqnarray*}
R_{n,m}&=& \varepsilon_{n} V_{n,0}+ \varepsilon_{n+1} (V_{n,1}-V_{n,0})+\cdots+\varepsilon_{m} (V_{n,m-2}-V_{n,n-m-1}\\
&=&\biggr( V_{n,0} (\varepsilon_{n}-\varepsilon_{n+1}) + \cdots + V_{n,m-m-1} (\varepsilon_{m-n-2}-\varepsilon_{m-n-1}) \biggr)+ \biggr(V_{n,m-n} \varepsilon_{m}\biggr)
\end{eqnarray*}

\bigskip \noindent The second line is obtained by developing the first. From there, use (A1) and (A2) to show that

$$
R_{n,m}=\sum_{n\leq k \leq m} v_{k} \varepsilon_{k} \rightarrow 0 \ as \ (n,m) \rightarrow (+\infty, +\infty). 
$$

\bigskip \noindent Conclude

\bigskip \noindent \textbf{Solution}. \\

\noindent To show that the series $\sum_{n \geq 0} v_{n} \varepsilon_{n}$ is convergent, it is enough enough that show that the sequence of partial sums 
$(\sum_{0 \leq k \leq n} v_{k} \varepsilon_{k})_{n\geq 0}$ is Cauchy, that is

$$
R_{n,m}=\sum_{n\leq k \leq m} v_{k} \varepsilon_{k} \rightarrow 0 \ as \ (n,m) \rightarrow (+\infty, +\infty). \ (A3)
$$

\bigskip \noindent Let us fix $0 \leq n \leq m$ and defines $V_{n,j}=v_n +  \cdots + v_{n+j}$, $j\geq 0$. We see that

$$
R_{n,m}=\varepsilon_{n} V_{n,0}+ \varepsilon_{n+1} (V_{n,1}-V_{n,0})+\cdots+\varepsilon_{m} (V_{n,m-2}-V_{n,n-m-1}.
$$

\bigskip \noindent By developing the right member of the above equation, we get

$$
R_{n,m} =\biggr( V_{n,0} (\varepsilon_{n}-\varepsilon_{n+1}) + \cdots + V_{n,m-m-1} (\varepsilon_{m-n-2}-\varepsilon_{m-n-1}) \biggr)+ \biggr(V_{n,m-n} \varepsilon_{m}\biggr)
$$

\bigskip \noindent We get, by using the total boundedness of the sequence $(v_n)_{n\geq 0}$, that

$$
|R_{n,m}| \leq A\biggr( \biggr(\sum_{n \leq k \leq m} |\varepsilon_{k+1}-\varepsilon_{k}|\biggr) + \varepsilon_{m} \biggr). \ (A3)
$$

\bigskip \noindent By (A1), the partial sums of the sequence $(\varepsilon_{k+1}-\varepsilon_{k})_{k \geq 0}$ is Cauchy, and by combining with (A1), we get -A3) from (A4). $\square$

\newpage
\noindent \LARGE \textbf{Doc 05-07 : Lebesgue/Riemann-Stieljes integrals on $\mathbb{R}$ - Exercises with solutions}. \label{doc05-07}\\
\bigskip
\Large 

\bigskip \noindent This document is only an introduction to the comparison between the Riemann-Stieljes Integral and the Lebesgue-Stieljes integral.\\

\noindent Later, we will come back to a more deep document on the same subject after we have mastered the convergence theorems.\\

\bigskip \noindent \textbf{Reminder of the definition of the Riemann-Stieljes Integral}.\newline

\noindent The Classical Riemann-Stieljes Integral is defined for bounded functions on compacts intervals. In all this text, 
$a$ and $b$ are two real numbers  such that $a<b$ and $F: [a,b]\longrightarrow \mathbb{R}$ is non-constant non-decreasing function.\\

\bigskip
 
\noindent \textbf{Definition of the Riemann-Stieljes Integral on a compact set}.\\

\noindent Let $f : \ ]a,b]\longrightarrow \mathbb{R}$ be an arbitrary function. We are going to define the Riemann-Stieljes integral of $f$ on $]a,b]$ associated with $F$ and denoted 
\begin{equation*}
I=\int_{a}^{b}f(x)\text{ }dF(x).
\end{equation*}

\bigskip \noindent We begin to define the Riemann-Stieljes sums. For each $%
n\geq 1,$ consider a subdivision $\pi_n$ of $]a,b]$, that divides into $]a,b]$ into
the $\ell(n)$ sub-intervals

\begin{equation*}
]a,b]=\sum_{i=0}^{\ell (n)-1}]x_{i,n},x_{i+1,n}],
\end{equation*}

\bigskip \noindent with $a=x_{0,n}<x_{1,n}< \cdots < x_{\ell(n),n}=b$. The modulus of the subdivision $\pi _{n}$ is defined by%
\begin{equation*}
m(\pi _{n})=\max_{0\leq i\leq m(n)}(x_{i+1,n}-x_{i+1,n}).
\end{equation*}

\bigskip \noindent Associate to this subdivision $\pi _{n}$ an arbitrary sequence $c_{n}=(c_{i,n})_{1\leq i\leq \ell (n)}$ such that $c_{0,n}\in ]a,x_{1,n}]$ and $c_{i,n}\in \lbrack x_{i,n},x_{i+1,n}],$ $1\leq i\leq \ell
(n)-1$.\\

\noindent We may define a sequence of Riemann-Stieljes sum associated to the subdivision $\pi_n$ and the vector $c_n$ in the form 

\begin{equation*}
S_{n}(f,F,a,b,\pi _{n},c_{n})=\sum_{i=0}^{\ell
(n)-1}f(c_{i,n})(F(x_{i+1,n})-F(x_{i,n})).  \label{RSS}
\end{equation*}

\bigskip \noindent Since $f,$ $F$ and $]a,b]$ are fixed here, we may and do drop them in the expression of $S_n$.

\bigskip \noindent \textbf{Definition}. A bounded function $f$ is Riemann-Stieljes integrable with respect to $F$ if there exists a real
number $I$ such that any sequence of Riemann-Stieljes sums $S_{n}(\pi_{n},c_{n})$ converges to $I$ as $n\rightarrow 0$\ \ whenever $m(\pi
_{n})\rightarrow 0$ as $n\rightarrow \infty$, and the number $I$ is called the Riemann-Stieljes of $f$ on $[a,b]$.\\

\noindent If $F$ is the identity function, that is $F(x)=x$, $x\in \mathbb{R}$, $I$ is simply called the Riemann integral of $f$ over $[a,b]$ and the sums defined in Formula (RSS) are called Riemann sums.\\

\bigskip \noindent \textbf{Exercise 1.} \label{exercise01_sol_doc05-07}\\

\noindent Question (a) Let $f : ]a,b]\longrightarrow \mathbb{R}$ be a bounded and measurable function on $[a,b]$. Is-it Lebesque-Stieljes integrable? If yes, give a bound of its integral.\\

\noindent Question (b) Consider the function

$$
g=1_{]a,b]\cap \mathbb{Q}},
$$

\bigskip \noindent where $\mathbb{Q}$ is the set of rational numbers.\\

\noindent Consider Riemann-Stieljes sums associated to $g$ in which the points $c_{i,n}$ are chosen in $\mathbb{Q}$. Justify why you can make this choice. What are the values of these Riemann-Stieljes sums.\\

\noindent Proceed similarly by choosing the points $c_{n,i}$ as irrational numbers. Justify why you can make this choice. What are the values of these Riemann-Stieljes sums.\\

\noindent Conclude that $g$ is not Riemann-Stieljes integrable.\\

\noindent Question (c). Make a first comparison between the two integrals.\\

\bigskip \noindent \textbf{Solution}.\\

\noindent Question (a). Here we consider the induced Lebesgue measure $\lambda_{[a,b]}$  of the Lebesgue measure on $[a, \ b]$ which endowed with the induced with the induced Borel $\sigma$-algebra
$$
\mathcal{B}([a,b])=\{B \subset [a,b], B\in \mathcal{B}(\mathbb{R})\}.
$$

\bigskip \noindent One easily shows that $\mathcal{B}([a,b])$ is still generated by the semi-algebra $\mathcal{S}([a,b])=\{{a}, [a,b], ]x,y], \ a\leq x \leq y \leq b\}$. The measure $\lambda_{[a,b]}$ is finite. Thus, any constant function $M=M1_{[a,b]}$ is integrable, since

$$
\int M \ d\lambda_{[a,b]}=M \lambda_{[a,b]}([a,b])=M(b-a).
$$

\bigskip \noindent Next, by Formula (I2) of \textit{Point 05.06 of Doc 05-01 of this chapter}, any measurable and bounded function on $[a,b]$ is integrable and its integral is bounded by $(b-a)$ times 
its bound.\\

\noindent Question (b). Consider an arbitrary subdivision of $\pi_n=(a=x_{0,n}<x_{1,n}< \cdots x_{\ell(n),n}=b)$. Since the set of rational numbers $\mathbb{Q}$ is dense in $\mathbb{R}$, we may choose in any interval $]x_{i,n}, \ x_{i+1,n}[$  a rational number $e_{i,n}$, $0\leq i \leq \ell(n)-1$. With such a choice for the vector $(e_n)=(e_{0,n},...,e_{\ell(n)-1,n})$, we have

$$
\forall n\geq 1, \ S_{n}(g,F,a,b,\pi _{n},e_{n})=0.
$$

\bigskip \noindent Since the set of irrational numbers $\mathbb{J}$ is dense in $\mathbb{R}$, we may choose in any interval $]x_{i,n}, \ x_{i+1,n}[$  an irrational number $d_{i,n}$, $0\leq i \leq \ell(n)-1$. With such a choice for the vector $(d_n)=(d_{0,n},...,d_{\ell(n)-1,n})$, we have

$$
\forall n\geq 1, \ S_{n}(g,F,a,b,\pi _{n},d_{n})=b-a.
$$

\bigskip \noindent Hence, the sequence of Riemann sums $S_{n}(f,F,a,b,\pi _{n},c_{n})$ does not converges as $m(\pi_n)$ converges to as as $n\rightarrow +\infty$ for all choices of the sequence $(c_n)_{n\geq 1}$.
Thus, there is no Riemann-integral for $g$.\\

\noindent Question (c). The conclusion is the following : There exists at least on measurable and bounded function on $[a,b]$ which is Lebesgue and not Riemann.\\

\bigskip \noindent \textbf{Exercise 2}. \label{exercise02_sol_doc05-07}\\

\noindent We require the following assumption on $F$ : $F$ is of bounded variation on $[a,b]$, that is, there exists a real number$A>0$ such that for any subdivision $\pi=(z_0,...,z_s)$ of 
$[a,b]$, that is $a=z_0<z_1<...<z_s=b$, we have

$$
V(\pi,[a,b])=\sum_{j=0}^{s-1} |F(z_{i+1}-z_{i}| <A.
$$

\bigskip \noindent Question (a) Prove that if $f$ is continuous then the sequences of Riemann-Stieljes (RSS) is Cauchy, and then, the Rieman-integral of $f$ exists on $[a,b]$.

\noindent \textit{Hint}. By using the results \textit{Exercise 15 in Doc 04-02 of Chapter \ref{04_measures}, page \pageref{exercise15_doc04-02}}, it is direct to see that the application defined by

$$
\lambda_F(]x,y])=F(b)-F(a)
$$

\bigskip \noindent is a proper and an additive application on the semi-algebra $\mathcal{S}([a,b])=\{{a}, [a,b], ]x,y], \ a\leq x \leq y \leq b\}$.\\

\noindent We stress that the extension theory made in that Exercise 15 is still valid expect the properties related to sub-additivity. So, $\lambda_F$ is directly extended to an additive application one the class of finite sums of type $]x,y]$ in $[a,b]$. But, be careful : we cannot have the no-negativity without further assumptions on $F$, what we do not require for the moment. The theory of integral of elementary functions(but based here on intervals of the form $]x,y]$ is also valid. Precisely, if we denote

$$
h= \sum_{i=0}^{\ell (n)-1}f(c_{i,n}) 1_{]x_{i,n}, \ x_{i-1,n}]},
$$

\bigskip \noindent the \textit{integral},

$$
\int h \ d\lambda_F =\sum_{i=0}^{\ell (n)-1}f(c_{i,n}) F(x_{i+1})-F(x_{i,n}).
$$

\bigskip \noindent In the same spirit of the remarks above, this integral defined on the $\mathcal{E}_e$class of functions of the form

$$
\sum_{i=0}^{s-1} \alpha_i 1_{]z_{i+1}, \ z_i]}, \ a=z_0<z_1<...<z_s=b, \alpha_i \in \mathbb{R}, \  0\leq i \leq s-1,
$$

\bigskip \noindent is still well-defined and its linearity holds since we do not face infinite values. In the sequel, many properties in the solutions of Exercises series in Doc 05-05 are used.\\
 
\noindent Now take two Riemann-integrals for integers $p>1$ and $q>1$

\begin{equation*}
S_{p}=\sum_{i=0}^{\ell (p)-1}f(c_{i,p})(F(x_{i+1,p})-F(x_{i,p})).  \label{RSS1}
\end{equation*}

\bigskip \noindent and 

\begin{equation*}
S_{q}=\sum_{i=0}^{\ell(q)-1}f(c_{i,q})(F(x_{i+1,q})-F(x_{i,q})).  \label{RSS2}
\end{equation*}

\bigskip \noindent with, for $r=1,2$, we have $a=x_{0,r}<x_{1,r}<\cdots<x_{\ell(r),r}$, $c_{0,r}\in ]a,x_{1,r}]$ and $c_{i,r}\in \lbrack x_{i,r},x_{i+1,r}],$ $1\leq i\leq \ell (r)-1$,  $\pi_{r}=(x_{0,r}, ..., <x_{\ell(r),r})$.\\

\noindent Denote

$$
h_p= \sum_{i=0}^{\ell (p)-1}f(c_{i,p}) 1_{]x_{i,p}, \ x_{i-1,p}]},
$$

\bigskip \noindent and

$$
h_q= \sum_{i=0}^{\ell (q)-1}f(c_{i,q}) 1_{]x_{i,q}, \ x_{i-1,q}]},
$$

\bigskip \noindent and use the identities

$$
S_{p} = \int f_p \ d\lambda_F \ and \ S_{q} = \int f_q \ d\lambda_F
$$

\bigskip \noindent To lessen the notations, put $A_i=]x_{i,p}, ,x_{i+1,p}]$ and $B_j=]x_{j,q}, \ x_{j+1,p}]$.\\

\noindent Now, use the techniques of the superpositions of partitions as seen in \textit{Point 03.22a of Doc 03-02 of Chapter \ref{03_setsmes_applimes_cas_speciaux}}, page \pageref{doc03-02} and use the superposition $\{ A_{j}B_{j}, \ (i,j)\in H \}$ of the partitions $\{ A_{i},  0\leq i\leq \ell(p)-1 \}$ and $\{ B_{j}, \ 0\leq j\leq \ell(q)-1 \}$, where $H=\{(i,j), \ 1\leq i\leq \ell(p), \  1\leq j\leq \ell(q), A_iB_j\neq \emptyset \}$. Provide new expressions of $h_p$ and $h_q$ related to the partition $\{ A_{j}B_{j}, \ (i,j)\in H\}$ of $]a,b]$.\\

\bigskip \noindent Finally, use the additivity of $\lambda_F$, in particular \textit{Formulas (DA), (DA1) and (DA2) in the solutions of Exercise 1 in Doc 05-05}, page  \pageref{exercise01_sol_doc05-05}, to establish that

$$
|S_p-S_q|\leq A \sup_{(s,t)\in [a,b]^2, |s-t|\leq m(\pi_p) \vee m(\pi_s)} |f(t)-f(s)|.
$$

\bigskip \noindent Deduce from this that the sequence $(S_n)_{n\geq 1}$ is Cauchy if $m(\pi_n)$ tends to zero as $n \rightarrow +\infty$.\\

\bigskip \noindent \textbf{Solution}.\\

\noindent We have to that the sequences of Riemann-Stieljes sums 
\begin{equation*}
S_{n}=\sum_{i=0}^{\ell(n)-1}f(c_{i,n})(F(x_{i+1,n})-F(x_{i,n})), \ n\geq 1,  \label{RSS}
\end{equation*}

\bigskip \noindent is Cauchy fay any choice of the vectors $c_n$, $n \geq 1$, as the modulus $m(\pi_n)$ goes to zero as $n$ increases indefinitely. We are going to follow the indications in the statements.\\

\noindent With the notations introduced, we have

$$
h_p= \sum_{i=0}^{\ell (q)-1}f(c_{i,q}) 1_{A_i}, \ and \ h_q= \sum_{i=0}^{\ell (q)-1}f(c_{i,q}) 1_{A_i}.
$$

\bigskip \noindent Using the superposition of the two partitions leads to

$$
h_p= \sum_{(i,j) \in H} f(c_{i,p}) 1_{A_iB_j}, \ and \ h_q= \sum_{(i,j) \in H} f(c_{i,q}) 1_{A_iB_j}.
$$

\bigskip \noindent This implies that for $s\in A_iB_j$, $h_p(s)=f(c_{i,p})$ and $h_q(s)=f(c_{i,q})$ with $((c_{i,p},c_{i,p})^2 \in A_iB_j=]x_{i,p},x_{i+1,p}[ \cap ]x_{i,q},x_{i+1,q}[$, for $(i,j) \in H$. Hence, for each $(i,j)\in H$,

\begin{eqnarray*}
|f(c_{i,p})-f(c_{i,q})| &\leq& \sup_{(s,t)\in [a,b]^2, |s-t|\leq |x_{i+1,p}-x_{i,p}|, |s-t|\leq |x_{i+1,q}-x_{i,q}|} |f(t)-f(s)|\\
&\leq& \sup_{(s,t)\in [a,b]^2, |s-t|\leq m(\pi_p) \vee m(\pi_q)} |f(t)-f(s)|=\delta(f,m(\pi_p) \vee m(\pi_q)),
\end{eqnarray*}

\bigskip \noindent where $\delta(f,r))$ stands for the uniform continuity modulus of $f$ with radius $r>0$.\\

\noindent Next, as in Formulas (DA1) and (DA2) in the solutions of \textit{Exercise 1 of Doc 05-05} (See page \pageref{exercise01_sol_doc05-05}),

$$
S_p=\int h_p \ d\lambda_F= \sum_{(i,j) \in H} f(c_{i,p}) \lambda_F(A_iB_j)) \ and \  S_q=\int h_q \ d\lambda_F= \sum_{(i,j) \in H} f(c_{i,q}) \lambda_F(A_iB_j)),
$$

\bigskip \noindent so that

$$
S_p-S_q= \sum_{(i,j) \in H} f(c_{i,p})-f(c_{i,q})- \lambda_F(A_iB_j)),
$$

\bigskip \noindent which, by combining the previous results with Formulas (DA)  in the solutions of \textit{Exercise 1 of Doc 05-05}, implies

\begin{eqnarray*}
|S_p-S_q| &leq&  \delta(f,m(\pi_p) \vee m(\pi_q))  \sum_{(i,j) \in H}  \lambda_F(A_iB_j))\\
&=& \delta(f,m(\pi_p) \vee m(\pi_q))  \sum_{i=0}^{\ell(p)-1} \sum_{j : (i,j) \in H}  \lambda_F(A_iB_j))\\
&=& \delta(f,m(\pi_p) \vee m(\pi_q))\sum_{i=0}^{\ell(p)-1}  |F(x_{i+1,p}-F(x_{i,p})|\\
&\leq &  A \delta(f,m(\pi_p) \vee m(\pi_q)).
\end{eqnarray*}

\bigskip \noindent From there, we say that a continuous function $f$ on the compact set $[a,b]$ is uniformly continuous and thus, the uniform continuity modulus $\delta(f,r)$ tends to zero with $r$. Thus, as $m(\pi_p)$ goes to zero as $p \rightarrow +\infty$ and $m(\pi_q)$ goes to zero as $q \rightarrow +\infty$, we have

$$
S_p-S_q \rightarrow 0,
$$

\bigskip \noindent which proves that $(S_n)$ is Cauchy.

\newpage
\noindent \LARGE \textbf{Doc 05-08 : Technical document on the definition of the integral of a function of constant sign}. \label{doc05-08}\\
\bigskip
\Large

\bigskip \noindent \textbf{Property P6}.\\

\bigskip \noindent  Let $0\leq h_{n}\uparrow ,$\ and (h$_{n})_{n\geq 1}\subset \mathcal{E}$\ such that there exists 0$\leq $g$\in \mathcal{E}$\ and  
\begin{equation*}
\text{ }\lim_{n}h_{n}\geq g,
\end{equation*}%
then
\begin{equation*}
\int h_{n}\text{ }dm\geq \int g\text{ }dm.  (EXT)
\end{equation*}%

\bigskip \noindent \textbf{Proof}.\\

\begin{equation*}
M=\max g\text{ and\ \ }r=\min g.
\end{equation*}%

\bigskip \noindent Remind that $g$ is a simple function. So it takes a finite number of values that are all non-negative. Let $M$ and $r$ be respectively maximum and the minimum if those values.\\

\noindent Suppose that (EXT) is proved for $r>0$. We are going to prove that we may extend it to the case
$r=0$. Indeed, set $A=(g>0)$. If $A$ is empty, then $g=0$ and (EXT) is true. Otherwise, the function

\begin{equation*}
g_{1}=g\text{ }1_{A}
\end{equation*}%

\bigskip \noindent is in $\mathcal{E}$, is positive with a nonzero minimum. We we may apply (EXT) to the
sequence ($h_{n}1_{A})$\ and to $g_{1}=g1_{A}$ since they are simple functions such that we have
\begin{equation*}
\lim_{n}h_{n}1_{A}\geq g_{1}\text{, }\min g_{1}>0.
\end{equation*}

\bigskip \noindent Now, using the fact that $g=0$ on $A^{c}$, we get
\begin{equation*}
\lim_{n}\int h_{n}1_{A}\text{ }dm\geq \int g_{1}\text{ }dm=\int g1_{A}\text{ 
}dm=\int g1_{A}\text{ }dm+\int g1_{A^{c}}dm=\int g\text{ }dm.
\end{equation*}

\bigskip \noindent But   
\begin{equation*}
\int h_{n}\text{ dm}\geq \int h_{n}1_{A}
\end{equation*}%

\bigskip \noindent we arrive at  
\begin{equation*}
\lim_{n}\int h_{n}\text{ dm}\geq \lim_{n}\int h_{n}1_{A}\geq \int g\text{ }dm.
\end{equation*}

\bigskip \noindent Then, it will be enough to show (EXT) for $\min g>0$.\\

\bigskip \noindent We use the same reasoning with respect to the maximum $M$ 
Suppose that (EXT) holds for $M<\infty$ We may extend to the case $M=+\infty$ as follows
Set $A=(g<+\infty)$ If $A$ is empty, then $g=+\infty$ and   
\begin{equation*}
h_{n}\uparrow +\infty 1_{\Omega }
\end{equation*}

\bigskip \noindent and, by definition,  
\begin{equation*}
\lim \int h_{n}\text{ }dm=+\infty =\int g\text{ }dm,
\end{equation*}

\bigskip \noindent so that (EXT) is still true, since $\int g dm=+\infty$. If $A$ is nonempty, set 
$B=A^{c}=(g=+\infty)$, and for all $c>0$,   
\begin{equation*}
g_{c}=g1_{A}+c1_{A^{c}}.
\end{equation*}

\bigskip \noindent It is clear $g_{c}$ is in $\mathcal{E}$. Further 
\begin{equation*}
g\geq g_{c}\geq 0
\end{equation*}

\bigskip \noindent and the maximum of $g_{c}$ is one of the finite value of $g$, so that  
\begin{equation*}
M_{c}=\max \text{ g}_{c}<+\infty.
\end{equation*}%

\bigskip \noindent Then  
\begin{equation*}
\lim_{n}h_{n}\geq g\geq g_{c}\text{ et }\max g_{c}<\infty.
\end{equation*}

\bigskip \noindent We can apply (EXT) to ($h_{n})$ and to $g_{c}$ :  
\begin{equation*}
\lim_{n}\int h_{n}\text{ }dm\geq \int g_{c}\text{ }dm\text{=}\int
g1_{A}+\int c1_{B}=\int g1_{A}+cm(g=+\infty).
\end{equation*}

\bigskip \noindent From this, we let $c\uparrow +\infty$, and get  

\begin{eqnarray*}
\lim_{n}\int h_{n}\text{ }dm&\geq& \int g1_{A}\text{ }dm+\infty \text{ }m(g=+\infty )\\
&=&\int g1_{A}\text{ }dm+\int +\infty \text{ }1_{A^{c}}\text{ }dm\\
&=&\int g1_{A}\text{ }dm+\int g1_{A^{c}}\text{ }dm\\
&=&\int g\text{ }dm.
\end{eqnarray*}

\bigskip \noindent This is the extension of (EXT) to $M=+\infty$.\\

\noindent Here again, is enough to show (EXT) for $\max g<+\infty$.\\

\noindent When we put all this together, we say that is enough to prove (EXT) when

\begin{equation*}
\max g<+\infty \text{ and }\min g>0.
\end{equation*}

\bigskip \noindent At this stage, we consider two cases.\\

\bigskip \noindent \textbf{Case : $m(\Omega )<+\infty$}.\\

\noindent Since  
\begin{equation*}
\min g=r>0,
\end{equation*}

\bigskip \noindent we can find an $\varepsilon$ such that  
\begin{equation*}
0<\varepsilon <r.
\end{equation*}

\bigskip \noindent Then we set  
\begin{equation*}
A_{n}=(h_{n}>g-\varepsilon).
\end{equation*}

\bigskip \noindent Since, for any $\omega \in \Omega $%
\begin{equation*}
h_{n}(\omega )\uparrow h(\omega )\geq g(\omega )>g(\omega )-\varepsilon, 
\end{equation*}

\bigskip \noindent there exists an integer $N$ such that  
\begin{equation*}
n\geq N\Rightarrow h_{n}(\omega )>g(\omega).
\end{equation*}

\bigskip \noindent Then
\begin{equation*}
\omega \in \bigcup_{n\geq 1}(h_{n}>g-\varepsilon )=\bigcup_{n\geq 1}A_{n}.
\end{equation*}

\bigskip \noindent But the the sequence $A_{n}$ is non-decreasing in $n$ and this leads to  
\begin{equation*}
\Omega =\bigcup_{n\geq 1}A_{n}=\lim_{n}A_{n}.
\end{equation*}

\bigskip \noindent Now, by the continuity of the measure,   
\begin{equation*}
m(A_{n})\uparrow m(\Omega )
\end{equation*}

\bigskip \noindent and since $m$ is bounded, we get  
\begin{equation*}
m(A_{n}^{c})\downarrow 0.
\end{equation*}

\bigskip \noindent The following inequalities involving non-negative simple functions 
\begin{equation*}
h_{n}\geq 1_{A_{n}}\text{ }h_{n}\geq (g-\varepsilon )1_{A_{n}}
\end{equation*}

\bigskip \noindent and the  fact that the integral is monotone of $\mathcal{E}$, imply  
\begin{eqnarray*}
\int h_{n}\text{ }dm&\geq& \int 1_{A_{n}}\text{ }h_{n}\text{ }dm\geq \int (g-\varepsilon )1_{A_{n}}\text{ }dm\\
&=&\int g\text{ }1_{A_{n}}\text{ }dm-\int \varepsilon 1_{A_{n}}\text{ }dm\\
&=&\int g\text{ }1_{A_{n}}\text{ }dm-\varepsilon \text{ }m(A_{n})\\
&=&\int g\text{ }(1-1_{A_{n}^{c}})\text{ }dm-\varepsilon \text{ }m(A_{n})\\
&=&\int g\text{ }dm-\int 1_{A_{n}^{c}}\text{ }dm-\varepsilon \text{ }m(A_{n})\\
&=&\int g\text{ }dm-m(A_{n}^{c})-\varepsilon \text{ }m(A_{n}).
\end{eqnarray*}

\bigskip \noindent We get   
\begin{equation*}
\int h_{n}\text{ }dm\geq \int g\text{ }dm-m(A_{n}^{c})-\varepsilon \text{ }m(A_{n}).
\end{equation*}

\bigskip \noindent Now, we let $\varepsilon \downarrow 0$, to get  
\begin{equation*}
\int h_{n}\text{ }dm\geq \int g\text{ }dm-m(A_{n}^{c}).
\end{equation*}

\bigskip \noindent Finally, let $n\uparrow +\infty$ and have
\begin{equation*}
\int h_{n}\text{ }dm\geq \int g\text{ }dm,
\end{equation*}

\bigskip \noindent which is (EXT).\\

\noindent \textbf{Case : $m(\Omega )=+\infty$}.\\

\bigskip \noindent We have 
\begin{equation*}
\int h_{n}\text{ }dm\geq \int 1_{A_{n}}\text{ }h_{n}\text{ }dm\geq \int (g-\varepsilon )1_{A_{n}}\text{ }dm\geq (m-\varepsilon )m(A_{n})\rightarrow \infty. 
\end{equation*}

\bigskip \noindent Then  
\begin{equation*}
\lim_{n}\int h_{n}\text{ }dm=+\infty \geq \int g\text{ }dm,
\end{equation*}

\bigskip \noindent which is (EXT).\\

\part{More on Measures and Integration}
%\newpage
%\chapter{Convergence Theorems}
%\newpage
\chapter{Introduction to Convergence Theorems and Applications} \label{06_convergence}

\noindent \textbf{Content of the Chapter}

\begin{table}[htbp]
	\centering
		\begin{tabular}{llll}
		\hline
		Type& Name & Title  & page\\
		S& Doc 06-01 &  The two Main Convergence Types - A summary & \pageref{doc06-01}\\
		S& Doc 06-02 & The two Main Convergence Theorems   & \pageref{doc06-02} \\
		S& Doc 06-03  & First Applications of The Convergence Theorems & \pageref{doc06-03} \\
		D& Doc 06-04  & Convergence Types - Exercises & \pageref{doc06-04} \\
		D& Doc 06-05  & Convergence Theorems - Exercises & \pageref{doc06-05} \\
		D& Doc 06-06  & Applications of Convergence Theorems - Exercises & \pageref{doc06-06} \\
		D& Doc 06-07  & Lebesgue-Stieljes/Riemann-Stieljes integrals - Exercises & \pageref{doc06-07} \\
		SD& Doc 06-08  & Convergence Types - Exercises with Solutions& \pageref{doc06-08} \\
		SD& Doc 06-09  & Convergence Theorems - Exercises with Solutions & \pageref{doc06-09} \\
		SD& Doc 06-10  & Applications of Convergence Theorems - Solutions & \pageref{doc06-10} \\
		SD& Doc 06-11  & Lebesgue-Stieljes/Riemann-Stieljes integrals - Solutions & \pageref{doc06-11} \\
		%D& Doc 06-12 & More on Convergence in measure - Exercises  & \\
		%D& Doc 06-12 & More on Convergence in measure - Exercises  with Solutions& \\
		%&  &   & \\
		\hline
		\end{tabular}
\end{table}

\newpage
\noindent \LARGE \textbf{Doc 06-01 : The two Main Convergence Types - A summary}. \label{doc06-01}\\
\bigskip
\Large

\bigskip \noindent In this tutorial, we provide the basic types of convergence and next the fundamental convergence theorems and some of their applications.\\

\bigskip \noindent Before you begin to read this note, we remind you to revise the appendix on limits in $\overline{\mathbb{R}}$, called : \textit{what should not be ignored on limits in $\overline{\mathbb{R}}$}, in Doc 11-01 in Chapter \ref{11_appendix}. You may find it in page \pageref{doc11-01}. A perfect mastering of inferior limits and superior limits is needed to be in ease when working on Convergence questions of any type.\\

\bigskip This chapter is devoted to convergence of sequences of real-valued measurable functions. A measure space $(\Omega ,\mathcal{A},m)$ is given for once in all the
sequel. Without any other explicit assumptions,  $(f_{n})_{n\geq 1}$, $(g_{n})_{n\geq 1}$, $(h_{n})_{n\geq 1}$ stand for sequences of real-valued and measurable functions defined from ($\Omega ,\mathcal{A},m)$\ to $\overline{\mathbb{R}}$ and $f$, $g$, $h$, etc. are real-valued and measurable functions defined on $(\Omega ,\mathcal{A},m)$.

\bigskip \noindent The following notation is adopted :\\

\bigskip \noindent The following fact is important before the definition of the convergence that will comes first.\\

$$
(f_{n} \rightarrow \ f) = \{\omega \in \Omega, \ \ f_{n}(\omega) \ converges \  to \ f(\omega)\}
$$

\noindent and

$$
(f_{n} \nrightarrow \ f) = \{\omega \in \Omega, \ \ f_{n}(\omega) \ does \ not \ converge \  to \ f(\omega)\}
$$

\noindent with

$$
(f_{n} \rightarrow \ f)^c=(f_{n} \nrightarrow \ f).
$$

\bigskip \noindent \textbf{(06.00)}. If $(f_{n})_{n\geq 1}$ is a sequence real-valued and measurable functions defined from ($\Omega ,\mathcal{A},m)$\ to $\overline{\mathbb{R}}$ and if $f$ is also a real-valued
and measurable function, then the set $(f_{n} \nrightarrow \ f)$ us measurable.\\

\bigskip \noindent \textbf{(06.01) - Definitions of the two main types of convergence}.\\

\bigskip \noindent \textbf{(06.01a) Almost-everywhere convergence}.\\

\bigskip \noindent A sequence of measurable applications $f_{n})_{n\geq 1}$\ defined from ($\Omega ,\mathcal{A},m)$\ to $\overline{\mathbb{R}}$ converges almost-everywhere to a measurable function $f:(\Omega ,\mathcal{A},m)$\ $\longmapsto \overline{\mathbb{R}}$ and we
denote
\begin{equation*}
f_{n}\longrightarrow f,\text{ }a.e.,
\end{equation*}

\bigskip \noindent if and only if the elements $\omega \in \Omega $ for which $f_{n}(\omega)$ does not converge to $f(\omega )$ form a null set with respect to $m$ that is 

\begin{equation*}
m(f_n \ \nrightarrow \ f)=0.
\end{equation*}

\bigskip \noindent \textbf{remark}. The \textit{a.e.} convergence may may be defined even if the functions $f_n$'s and $f$ take infinite values on non-null sets (with respect to $m$). But if these functions
are \textsl{a.e.} finite, more precise formulas are given.\\

\bigskip \noindent \textbf{(06.01b) Convergence in measure}.\\

\noindent A sequence of measurable and \textbf{a.e. finite} applications $(f_{n})_{n\geq 1}$\ defined from ($\Omega ,\mathcal{A},m)$\ to $\mathbb{R}$
converges in measure with respect to $m$ to an \textit{a.e} finite measurable function $f:(\Omega ,\mathcal{A},m)$\ $\longmapsto \mathbb{R}$, denoted
\begin{equation*}
f_{n}\rightarrow _{m}f
\end{equation*}

\bigskip \noindent if and only for any $\varepsilon >0,$
\begin{equation*}
m(\left\vert f_{n}-f\right\vert >\varepsilon)\longrightarrow 0\text{ as }
n\longrightarrow +\infty .
\end{equation*}

\noindent \textbf{Remark}. The convergence in measure is only possible if the limit $f$ and the $f_n$'s are \textbf{a.e.} since we need to make the differences $f_n-f$. The \textit{a.e.} finiteness justifies this.\\

\noindent As remarked in the \textit{Measure Space Reduction Property} seen in Doc 05-01 in Chapter  \ref{05_integration} (page \pageref{msrpA}), we may and do suppose that $f$ and the $f_{n}$'s are finite on the whole space.\\

\noindent \textbf{NB}. It is important to notice the inequality in $(| f_{n}-f|>\varepsilon)$ may be strict or not.\\

\bigskip \noindent \textbf{(06.02) - Properties of the \textit{a.e.} convergence}.\\

\bigskip \noindent \textbf{(06.02a) - The \textit{a.e.} convergence is well-defined}.\\

\bigskip \noindent The definition of the \textit{a.e. limit} is possible if the set $(f_{n}\rightarrow f)$ is measurable, as pointed out earlier.\\

\noindent  Set 

\begin{equation*}
A=(f=+\infty ) \bigcap \left( \bigcap_{k\leq 1} \bigcup_{n\geq 1} \bigcap_{p\geq n} (f_{p} \geq k) \right)
\end{equation*}

\bigskip\ 
\begin{equation*}
B=(f=-\infty ) \bigcap \left( \bigcap_{k\leq 1} \bigcup_{n\geq 1} \bigcap_{p\geq n} (f_{p} \leq -k) \right)
\end{equation*}

\bigskip \noindent and \\

\begin{equation*}
C=(|f|<+\infty ) \bigcap \left( \bigcap_{k \leq 1} \bigcup_{n\geq 1} \bigcap_{p\geq n} (|f_{p}-f| \leq 1/k) \right).
\end{equation*}

\bigskip \noindent Hence, we have 

\begin{equation*}
(f_{n}\longrightarrow f)=A+B+C
\end{equation*}

\bigskip \noindent  and that $(f_{n} \nrightarrow f)$ is measurable.\\

\bigskip \noindent \textbf{(06.02b) - \textit{a.e.} uniqueness of \textit{a.e.} limits}.\\

\noindent The \textit{a.e.} limit is \textit{a.e.} unique, meaning that :\\

$$
(f_n \rightarrow f \ a.e. \text{ and } f_n \rightarrow g \ a.e.) \ \ \Rightarrow f=g \ a.e. 
$$

\bigskip \noindent \textbf{(06.02c) - Characterization of the \textit{a.e} convergence if $f$ and the $f_n$'s are \textit{a.e.} finite}.\\

\noindent We use the Measure Space Reduction Property (page \pageref{msrpA}) and do assume that $f$ and the $f_n$'s are finite. Besides, in Point \textit{06.02a}, the sets $A$ and $B$ are null sets. It remains to only check that $C$ is a null set to have the \textit{a.e.} convergence. We have :\\

$$
f_n \rightarrow f \ \textit{a.e.} \ as \ n\rightarrow +\infty
$$

\bigskip \noindent if and only if 

$$
m \left( \bigcap_{k\geq 1}\bigcap_{N\geq 1}\bigcup_{n\geq N}(\left\vert f_{n}-f\right\vert <1/k)  \right) = 0
$$

\bigskip \noindent if and only if, for any $k\geq 1$,

$$
m \left( \bigcup_{N\geq 1}\bigcap_{n\geq N}(\left\vert f_{n}-f\right\vert <1/k)  \right) = 0
$$

\bigskip  \noindent if and only if, for any $\varepsilon >0$,

$$
m \left( \bigcup_{N\geq 1}\bigcap_{n\geq N}(\left\vert f_{n}-f\right\vert \geq \varepsilon)  \right) = 0.
$$

\bigskip \noindent \textbf{(06.02d) - Operations of \textit{a.e} limits}. If we deal with \textit{a.e.} finite functions, it is possible to make operations on between sequences of them. And we have the following properties.\\

\noindent Let $(f_{n})_{n\geq 1}$ and $(g_{n})_{n\geq 1}$ be sequences of \textit{a.e.} finite functions. Let $a$ and $b$ be finite real numbers. Suppose that  $f_{n}\rightarrow f$ $a.e.$ and $g_{n}\rightarrow g$ $a.e.$. Let $H(x,y)$ a continuous function of $(x,y) \in D$, where $D$ is an open set of $\mathbb{R}^2$. We have : \\

\noindent \textbf{(1)} $af_{n}+bg_{n}\rightarrow af+bg$ $a.e$.\\

\noindent \textbf{(2)} $f_{n}g_{n}\rightarrow fg$ $a.e$\\

\noindent \textbf{(3)} If $m(g=0)=0$ (that is $g$ is $a.e$ nonzero), then

$$
f_{n}/g_{n}\rightarrow f/g, \ a.e.
$$

\noindent \textbf{(4)} If $(f_{n}, g_n)_{n\geq 1} \subset D$ \textit{a.e.} and $(f,g) \in D$ \textit{a.e.}, then

$$
H(f_{n},g_{n})\rightarrow H(f,g), \ a.e.\\
$$

\bigskip \noindent \textbf{(06.02e) - \textit{a.e.} Cauchy sequences}. If we deal with \textit{a.e.} finite functions, it is possible to consider Cauchy Theory on sequences of them. And we have the following definition and characterizations.\\

\noindent \textbf{Definition}. A sequence $(f_{n})_{n\geq 1}$ of \textit{a.e.} finite functions is an $m$-\textit{a.e.} Cauchy sequence if and only if

$$
m( f_p-f_q \nrightarrow 0 \ as \ (p,q) \rightarrow (+\infty,+\infty) )=0,
$$

\bigskip \noindent that is, the $\omega$ for which the real sequence $(f_n(\omega))_{n\geq 0}$ is not Cauchy sequence on $\mathbb{R}$ form a $m$-null set.\\

\bigskip \noindent \textbf{Other expressions}. A sequence $(f_{n})_{n\geq 1}$ of \textit{a.e.} finite functions is an $m$-\textit{a.e.} Cauchy sequence :\\

\noindent if and only if for any $k\geq 1$,

$$
m\left( \bigcap_{n\geq 1} \bigcup_{p\geq n} \bigcup_{q \geq n} (|f_p-f_q| > 1/k) \right)=0
$$

\noindent if and only if for any $k\geq 1$,

$$
m\left( \bigcap_{n\geq 1} \bigcup_{p\geq 0}  (|f_{p+n}-f_n| > 1/k) \right)=0
$$

\noindent if and only if for any $\varepsilon >0$,

$$
m\left( \bigcap_{n\geq 1} \bigcup_{p\geq n} \bigcup_{q \geq n} (|f_p-f_q| > \varepsilon) \right)=0
$$

\noindent if and only if for any $\varepsilon>0$,

$$
m\left( \bigcap_{n\geq 1} \bigcup_{p\geq 0}  (|f_{p+n}-f_n| > \varepsilon) \right)=0.
$$

\bigskip \noindent \textbf{Property}. Let  $(f_{n})_{n\geq 1}$ be a sequence of \textit{a.e.} finite functions.\\ 

\noindent $(f_{n})_{n\geq 1}$ is an $m$-\textit{a.e.} Cauchy sequence if and only if $(f_{n})_{n\geq 1}$ converges \textit{a.e.} to an \textit{a.e.} finite function.

\bigskip  \bigskip \noindent \textbf{(06.03) - Properties of the convergence in measure}.\\

\noindent \textbf{(06.03a) - \textit{a.e.} uniqueness of limits in measure}.\\

\noindent The limit in measure is \textit{a.e.} unique, meaning that :\\

$$
(f_n \rightarrow_m f  \text{ and } f_n \rightarrow_m g) \ \ \Rightarrow f=g \ a.e. 
$$

\bigskip \bigskip \noindent \textbf{(06.03b) - First comparison between the two types of limits}.\\

\noindent Let $m$ be a \textbf{finite measure} and $f$ be \textit{a.e.} finite. Let $(f_{n})_{n\geq 1}$ be a sequence of  $m$-\textit{a.e.} functions. If $f_n$ converges to $f$ \textit{a.e.}, then
$f_n$ converges to $f$ in measure.\\

\bigskip \noindent \textbf{(06.03c) - Operations of limits in measure}. We still deal with \textit{a.e.} finite functions.\\

\noindent The operations of limits in measure are linked with advanced notions like boundedness in measure and other advanced types of convergence. This explains why we cannot have the same properties as the \textit{a.e.} convergence. But if the limits in measure are finite real numbers, we have equivalent results. Otherwise, you have to treat each as a a special case.\\

\bigskip \noindent \textbf{(06.03c-A) - General case}.\\

\noindent Let $f_{n}\rightarrow _{m} f$ and $g_{n}\rightarrow _{m}g,$ $a\in \mathbb{R}$. We have :\\

\noindent \textbf{(1)} $f_{n}+g_{n}\rightarrow_{m}f+g$.\\

\noindent \textbf{(2)} $af_{n} \rightarrow af$\\

\bigskip \noindent \textbf{(06.03c-B) - Finite and constant limits in measure}.\\

\noindent Let $(f_{n})_{n\geq 1}$ and $(g_{n})_{n\geq 1}$ be sequences of \textit{a.e.} finite functions. Let $a$ and $b$ be finite real numbers. Suppose that  $f_{n}\rightarrow_m A$ and 
$g_{n}\rightarrow_m B$, where $A$ and $B$ are finite real numbers. Let $H(x,y)$ a continuous function of $(x,y) \in D$, where $D$ is an open set of $\mathbb{R}^2$. We have :\\

\noindent \textbf{(1)} $af_{n}+bg_{n}\rightarrow_m aA+bB$.\\

\noindent \textbf{(2)} $f_{n}g_{n}\rightarrow_m AB$.\\

\noindent \textbf{(3)} If $B \neq 0$,  then

$$
f_{n}/g_{n}\rightarrow_m A/B.
$$

\bigskip \noindent \textbf{(4)} If $(f_{n}, g_n)_{n\geq 1} \subset D$ \textit{a.e.} and $(A,B) \in D$, then

$$
H(f_{n},g_{n})\rightarrow_m H(A,B).
$$

\bigskip \noindent \textbf{(06.03d) - Cauchy sequences in measure}. Here again, we deal with \textit{a.e.} finite functions and consider a Cauchy Theory on sequences of them. And we have the following definition and characterizations.\\

\noindent \textbf{Definition}. A sequence $(f_{n})_{n\geq 1}$ of \textit{a.e.} finite functions is Cauchy sequence in measure if and only if, for any $\varepsilon$,

$$
m( |f_p-f_q| > \varepsilon) \rightarrow 0 \ as \ (p,q) \rightarrow (+\infty,+\infty).
$$

\bigskip \noindent \textbf{Properties}. Let  $(f_{n})_{n\geq 1}$ be a sequence of \textit{a.e.} finite functions.\\ 

\noindent \textbf{P1} \noindent $(f_{n})_{n\geq 1}$ is a Cauchy sequence in measure if and only if $(f_{n})_{n\geq 1}$ converges in measure to an \textit{a.e.} finite function.\\

\noindent \textbf{P2} \noindent If $(f_{n})_{n\geq 1}$ is a Cauchy sequence in measure, then $(f_{n})_{n\geq 1}$ possesses a subsequence $(f_{n_k})_{k\geq 1}$ and an \textit{a.e.} measurable function $f$ such that 
$$
f_{n_k} \rightarrow f \ \textit{a.e.} \ as \ k\rightarrow +\infty, 
$$

\noindent and 

$$
f_n \rightarrow_m f \ as \ n\rightarrow +\infty.
$$

\bigskip \noindent \textbf{(06.04) - Comparison between $a.e$ convergence and convergence in measure}.\\

\noindent \textbf{(06.04b)}.  If $f_{n}\rightarrow f$ $a.e.$ and $m$ is finite, then $f_{n}\rightarrow _{m}f$.\\

\noindent \textbf{Terminology}. From this, $a.e.$ convergences are called strong in opposite of convergence in measure that are called weak.\\

\bigskip \noindent The reverse implication is not true. It is only true for a sub-sequence as follows.\\

\bigskip \noindent \textbf{(06.04b)}. Let $f_{n}\rightarrow_{m}f$. Then, there exists a sub-sequence $(f_{n_{k}})_{k\geq 1}$ of $(f_{n})_{n\geq 1}$ converging to $f$.\\

\bigskip \noindent \textbf{(06.05) - Terminology in the Frame of Probability Theory}.\\

\noindent If the measure space is a probability space, that is $m=\mathbb{P}$ a probability measure, an almost-everywhere \textit{a.e.} property $\mathcal{P}$ becomes an almost-sure \textit{a.s.} property, that is

$$
\mathbb{P}( \mathcal{P} \ true )=1.
$$

\noindent In that case : \\

\noindent (a) The \textit{a.e.} convergence becomes the \textit{a.s.} convergence.\\

\noindent (b) The convergence in measure becomes convergence in probability, denoted as $ \rightarrow_{\mathcal{P}}$.\\

\noindent Since a probability measure is finite, the \textit{a.s.} convergence implies the convergence in probability and the convergence in probability implies the existence of a subsequence converging \textit{a.s.} to the limit in probability.\\

\noindent \LARGE \textbf{Doc 06-02 : The two Main Convergence Theorems}. \label{doc06-02}\\
\bigskip
\Large

\noindent \textbf{NB}. The theorems of this section are among the most important tools not only in modern integration theory but in modern Mathematics. They should be very well mastered.\\

\bigskip \noindent \textbf{(06.06) Monotone Convergence Theorem - MCT}.\\

\bigskip \noindent Let $(f_{n})_{n\geq 0}$ be a non-decreasing sequence of \textit{a.e.} non-negative real-valued and measurable functions. Then
\begin{equation*}
\lim_{n \uparrow +\infty} \int f_{n}dm = \int \lim_{n \uparrow +\infty} f_n \ dm.
\end{equation*}

\noindent \textbf{NB}. This result holds also the name of Beppo-Levi.\\

\bigskip \noindent \textbf{(06.07a) Fatou-Lebesgue Lemma}.\\

\bigskip \noindent \textbf{(06.07a-A)}. Let $(f_{n})_{n\geq 0}$ be a sequence real-valued and measurable functions \textit{a.e.} bounded \textbf{below} by an integrable function $h$, that is 
$h \leq f_n$ \textit{a.e.} for all $n\neq 0$. Then, we have

\begin{equation*}
\liminf_{n\rightarrow +\infty}\int f_{n}dm \geq \int \liminf_{n\rightarrow +\infty }f_{n}dm. 
\end{equation*}

\bigskip \noindent By using the opposites of the functions above, we get an analog related to the limit superior.\\

\bigskip \noindent \textbf{(06.07a-B)}. Let $(f_{n})_{n\geq 0}$ be a sequence real-valued and measurable functions \textit{a.e.} bounded \textbf{above} by an integrable function $h$, that is 
$f_n \leq h$ \textit{a.e.} for all $n\neq 0$. Then, we have

$$
\limsup_{n\rightarrow +\infty}\int f_{n}dm \leq \int \limsup_{n\rightarrow +\infty }f_{n}dm.
$$

\bigskip \bigskip \noindent \textbf{(06.07b) (Dominated Convergence Theorem - DCT)}.\\

 \noindent Let $(f_{n})_{n\geq 0}$ be a sequence real-valued and measurable functions \textit{a.e.} such that :\\

\noindent (1) it bounded  by an integrable function $h$, that is $|f_n| \leq h$ \textit{a.e.} for all $n\geq 0$,\\

\noindent (2) it converges in measure or \textit{a.e.} to a function $f$.\\

\noindent  Then $f$ is integrable and the sequence $(\int f_{n} \ dm)_{n\geq 0}$ has a limit with

$$
\lim_{n\rightarrow +\infty}\int f_{n} \ dm = \int f \ dm.
$$

\bigskip \noindent and

$$
\int |f| \ dm = \int h \ dm.
$$

\bigskip \noindent \textbf{NB}. This result us also called the Lebesgue Dominated Convergence Theorem.\\

\bigskip \bigskip \noindent \textbf{(06.07c) (The Young Dominated Convergence Theorem - YDCT)}.\\

 \noindent Let $(f_{n})_{n\geq 0}$ be a sequence real-valued and measurable functions \textit{a.e.} such that :\\

\noindent (1) there exists a sequence of integrable functions $(h_n)_{n\geq 0}$ such that is $|f_n| \leq h_n$ \textit{a.e.} for all $n\geq 0$,\\

\noindent (2) the sequence $(h_n)_{n\geq 0}$ converges to an integrable $h$ such that 
$$
\lim_{n\rightarrow +\infty}\int h_{n} \ dm = \int h \ dm,
$$

\noindent (3) the sequence $(f_{n})_{n\geq 0}$ converges in measure or \textit{a.e.} to a function $f$.\\

\noindent  Then $f$ is integrable and the sequence $(\int f_{n} \ dm)_{n\geq 0}$ has a limit with

$$
\lim_{n\rightarrow +\infty}\int f_{n} \ dm = \int f \ dm.
$$

\bigskip \noindent and

$$
\int |f| \ dm = \int h \ dm.
$$

\noindent \LARGE \textbf{Doc 06-03 : First Applications of The Convergence Theorems}. \label{doc06-03}\\
\bigskip
\Large

\bigskip \noindent \textbf{(06.08) - Applications of the Monotone Convergence Theorem}.\\

\bigskip \noindent \textbf{(06.08a)}. Let $(f_{n})$ be sequence of \textit{a.e.} non-negative measurable functions. We have
\begin{equation*}
\sum_{n}\left\{ \int f_{n}\text{ }dm\right\} =\int \left\{\sum_{n}f_{n}\right\} \text{ }dm.
\end{equation*}

\bigskip \noindent \textbf{(06.08b)}. Let be $(a_{n})_{n\geq 0}=(\left( a_{n,p}\right) _{p\geq 0})_{p\geq
0}=(a_{n,p})_{(n\geq 0,p\geq 0)}$ a rectangular sequence of non-negative terms. We have 
\begin{equation*}
\sum_{n\geq 0}\sum_{p\geq 0}a_{n,p}=\sum_{p\geq 0}\sum_{n\geq 0}a_{n,p}.
\end{equation*}

\bigskip \noindent \textbf{(06.09) - Applications Dominated Convergence Theorem (DCT)}.\\

\bigskip \noindent \textbf{(06.09a)}. Let $(f_{n})$ be sequence of measurable functions such that%
\begin{equation*}
\int \left\{ \sum_{n}\left\vert f_{n}\right\vert \right\} dm<\infty.
\end{equation*}

\bigskip \noindent We have

\begin{equation*}
\sum_{n}\left\{ \int f_{n}\text{ }dm\right\} =\int \left\{\sum_{n}f_{n}\right\} \text{ }dm.
\end{equation*}

\bigskip \noindent \textbf{((06.09b))}. Let be $(a_{p})_{p\geq 0}=(\left( a_{n,p}\right) _{n\geq 0})_{p\geq
0}=(a_{n,p})_{(n\geq 0,p\geq 0)}$ a rectangular sequence such that%
\begin{equation*}
\sum_{p\geq 0}\sum_{n\geq 0}\left\vert a_{n,p}\right\vert <\infty.
\end{equation*}

\bigskip \noindent We have
\begin{equation*}
\sum_{p\geq 0}\sum_{n\geq 0}a_{n,p}=\sum_{n\geq 0}\sum_{p\geq 0}a_{n,p}.
\end{equation*}

\bigskip \noindent \textbf{(06.09c)}. Let $f(t,\omega ), \ t\in T$ a family of measurable functions indexed by $t\in I$, where $I=]a,b[, \ a<b$, to fix ideas. Let $t_{0}\in I$.\\

\bigskip \noindent (a) (Local continuity) Suppose that for some $\varepsilon >0,$ for $\varepsilon $ small enough
that $]t_{0}-\varepsilon ,t_{0}+\varepsilon \lbrack \subset I),$ there
exists $g$ integrable such that%
\begin{equation*}
\forall (t\in ]t_{0}-\varepsilon ,t_{0}+\varepsilon \lbrack ),\left\vert
f(t,\omega )\right\vert \leq g
\end{equation*}

\bigskip \noindent and, for almost every $\omega$, the function
\begin{equation*}
t\hookrightarrow f(t,\omega)
\end{equation*}

\bigskip \noindent is continuous at $t_{0}$. Then the function%
\begin{equation*}
F(t)=\int f(t,\omega )dm(\omega )
\end{equation*}

\bigskip \noindent is continuous at $t_{0}$ and%
\begin{equation*}
\lim_{t\rightarrow t_{0}}\int f(t,\omega )dm(\omega )=\int
\lim_{t\rightarrow t_{0}}f(t,\omega )dm(\omega )=\int f(t_{0},\omega
)dm(\omega ).
\end{equation*}

\bigskip \noindent (b) (Global continuity) Suppose that there exists $g$ integrable such that
\begin{equation*}
\forall (t\in I),\left\vert f(t,\omega )\right\vert \leq g
\end{equation*}

\bigskip \noindent and, for almost every $\omega ,$ the function%
\begin{equation*}
t\hookrightarrow f(t,\omega )
\end{equation*}

\bigskip \noindent is continuous at each $t\in I$. Then the function%
\begin{equation*}
F(t)=\int f(t,\omega )dm(\omega )
\end{equation*}

\bigskip \noindent is continuous at each $t\in I$ and for any $u\in I$%
\begin{equation*}
\lim_{t\rightarrow u}\int f(t,\omega )dm(\omega )=\int \lim_{t\rightarrow
u}f(t,\omega )dm(\omega )=\int f(u,\omega )dm(\omega ).
\end{equation*}

\bigskip \noindent \textbf{(06.09b)}. Let $f(t,\omega ),t\in T$ a family of measurable functions indexed by $%
t\in I$, where $I=]a,b[,a<b$, to fix ideas. Let $t_{0}\in I$.\\

\bigskip \noindent (a) (Local differentiability) Suppose that for almost every $\omega$, the function%
\begin{equation*}
t\hookrightarrow f(t,\omega )
\end{equation*}

\bigskip \noindent is differentiable at $t_{0}$ with derivative
\begin{equation*}
\frac{df(t,\omega )}{dt}|_{t=t_{0}}=f(t_{0},\omega )
\end{equation*}

\bigskip \noindent and that for some $\varepsilon >0,$ for $\varepsilon $ small enough that $%
]t_{0}-\varepsilon ,t_{0}+\varepsilon \lbrack \subset I),$ there exists $g$
integrable such that, 
\begin{equation*}
\forall (\left\vert h\right\vert \leq \varepsilon ),\left\vert \frac{%
f(t+h,\omega )-f(t,\omega )}{h}\right\vert \leq g.
\end{equation*}

\bigskip \noindent Then the function
\begin{equation*}
F(t)=\int f(t,\omega )dm(\omega )
\end{equation*}

\bigskip \noindent is differentiable $t_{0}$ and%
\begin{equation*}
F^{\prime }(t_{0})\frac{d}{dt}\int f(t,\omega )dm(\omega )|_{t=t_{0}}=\int 
\frac{df(t,\omega )}{dt}|_{t=t_{0}}dm(\omega).
\end{equation*}

\bigskip \noindent (b) (Everywhere differentibility) Suppose that $f(t,\omega )$ is $a.e$ differentiable at each t$\in T$ and
there $g$ integrable such that%
\begin{equation*}
\forall (t\in I),\left\vert \frac{d}{dt}f(t,\omega )\right\vert \leq g.
\end{equation*}

\bigskip \noindent Then, the function
\begin{equation*}
t\hookrightarrow f(t,\omega )
\end{equation*}

\bigskip \noindent is differentiable a each $t\in I$ and at each $t\in I,$ 
\begin{equation*}
\frac{d}{dt}F(t)=\int \frac{d}{dt}f(t,\omega)dm(\omega).
\end{equation*}

\bigskip \noindent \textbf{(06.09e)}. Let
\begin{equation*}
U(t)=\sum_{n}u_{n}(t)
\end{equation*}

\bigskip \noindent a convergent series of functions for $\left\vert t\right\vert <R,R>0.$ Suppose
that   each $u_{n}(t)$ is differentiable and there exists $g=(g_{n})$ integrable that is 
\begin{equation*}
\sum_{n}g_{n}<+\infty 
\end{equation*}

\bigskip \noindent such that
\begin{equation*}
\text{For all }\left\vert t\right\vert <R,\left\vert \frac{d}{dt}
u_{n}(t)\right\vert \leq g.
\end{equation*}

\bigskip \noindent Then, the function $U(t)$ is differentiable and
\begin{equation*}
\frac{d}{dt}U(t)=\sum_{n}\frac{d}{dt}u_{n}(t).
\end{equation*}

\noindent \LARGE \textbf{DOC 06-04 : Convergence Types - Exercises } \label{doc06-04}\\
\bigskip
\Large

\bigskip \noindent \textbf{NB}. Any unspecified limit below is meant as $n \rightarrow +\infty$.\\ 

\noindent \textbf{Exercise 01}. \label{exercises01_doc06-04} (Useful Formulas in $\mathbb{R}$).\\

\noindent (a) Quickly check the following equivalences for an real number $x$ or a sequence or elements of $\overline{\mathbb{R}}$.\\

$$
(\forall \varepsilon>0, |x|\leq \varepsilon) \Leftrightarrow (\forall r \in \mathbb{Q} \ and \ r>0, |x|\leq r). \ \ (EQUIV01)
$$

$$
(\forall \varepsilon>0, |x|\leq \varepsilon) \Leftrightarrow  (\forall k \in \mathbb{N} \ and \ k\geq 1, |x|\leq 1/k). \ \ (EQUIV02)
$$

$$
(\forall A>0, |x|\geq A) \Leftrightarrow  (\forall k \in \mathbb{N} \ and \ k\geq 1, |x| \geq k). \ \ (EQUIV03)
$$

\bigskip \noindent Formula (EQUIV04a) :

$$
(\forall A>0, \ \exists n\geq 0, \ \forall p\geq n : \ |x_p|\geq A)
$$
$$
\Leftrightarrow
$$
$$
(\forall k \in \\mathbb{N} \ and \ k\geq 1, \ \exists n\geq 0, \ \forall p\geq n : \ |x_p|\geq k).
$$

\bigskip \noindent Formula (EQUIV04b) :

$$
(\forall A>0, \ \exists n\geq 0, \ \forall p\geq n : \ |x_p|\leq - A)
$$
$$
\Leftrightarrow
$$
$$
(\forall k \in \mathbb{N} \ and \ k\geq 1, \ \exists n\geq 0, \ \forall p\geq n : \ |x_p|\geq -k).
$$

\bigskip \noindent Formula (EQUIV05) :

$$
(\forall \varepsilon>0, \ \exists n\geq 0, \ \forall p\geq n : \ |x_p|\leq \varepsilon)
$$
$$
\Leftrightarrow
$$
$$
(\forall k \in \mathbb{N} \ and \ k\geq 1, \ \exists n\geq 0, \ \forall p\geq n : \ |x_p|\leq 1/k).
$$

\bigskip \noindent etc.\\

\noindent (b) Let $a$ and $b$ be two real numbers and $\varepsilon>0$. Let $(a_j)_{1\leq j \leq k}$ be a family of $k>0$ real numbers. Let also $(\varepsilon_j)_{1\leq j \leq k}$ be a family of 
positive real numbers such that $\varepsilon_1+\varepsilon_2+ ... + \varepsilon_k = \varepsilon$. By using the negation of the triangle inequality, justify the relations 

$$
(|a+b|>\varepsilon) \Rightarrow (|a|>\varepsilon/2) \ or \ (|b|>\varepsilon/2). \ \ (ROS01)
$$

\bigskip \noindent and

$$
\biggr(\left| \sum_{1\leq j \leq k} a_j \right|>\varepsilon\biggr) \Rightarrow (\exists j \in \{1,...,k\}, \ |a_j|>\varepsilon_j). \ \ (ROS02)
$$

\bigskip \noindent (c) Let $f$ and $g$ be two functions defined on $\Omega$ with values in $\overline{\mathbb{R}}$. Let $(f_j)_{1\leq j \leq k}$ be a family $k>0$ of functions defined on $\Omega$ with values in $\overline{\mathbb{R}}$. real numbers. Let also $(\varepsilon_j)_{1\leq j \leq k}$ be a family of positive real numbers such that 
$\varepsilon_1+\varepsilon_2+ ... + \varepsilon_k = \varepsilon$. By using the previous question, justify the relations

$$
(|f+g|>\varepsilon) \subset (|f|>\varepsilon/2) \cup (|g|>\varepsilon/2). \ \ (OS01)
$$

\bigskip \noindent and

$$
\biggr(\left| \sum_{1\leq j \leq k} f_j \right|>\varepsilon\biggr) \subset \bigcup_{{1\leq j \leq k}} (|f_j|>\varepsilon_j). (OS02)
$$

\bigskip \noindent (d) Let $f$ and $g$ be two measurable functions defined on $\Omega$ with values in $\overline{\mathbb{R}}$. Let $(f_j)_{1\leq j \leq k}$ be a family $k>0$ of easurable functions defined on $\Omega$ with values in $\overline{\mathbb{R}}$. real numbers. Let also $(\varepsilon_j)_{1\leq j \leq k}$ be a family of positive real numbers such that 
$\varepsilon_1+\varepsilon_2+ ... + \varepsilon_k = \varepsilon$. By using the previous question, justify the relations

$$
m(|f+g|>\varepsilon) \leq m(|f|>\varepsilon/2) + m(|g|>\varepsilon/2). \ \ (MOS01)
$$

\bigskip \noindent and

$$
m\biggr(\left| \sum_{1\leq j \leq k} f_j \right|>\varepsilon\textbf{}) \leq \sum_{{1\leq j \leq k}} (|f_j|>\varepsilon_j). (MOS02)
$$

\bigskip \noindent \textbf{NB}. This exercise gives technical and simple facts which are known on $\mathbb{R}$, to be used in the sequel.\\

\bigskip \noindent \textbf{Exercise 02}. \label{exercises02_doc06-04}\\

\noindent Let  $f_{n})_{n\geq 1}$\ be a sequence of measurable functions from($\Omega ,\mathcal{A},m)$\ to $\overline{\mathbb{R}}$ and $f : (\Omega ,\mathcal{A},m)$\ $\longmapsto \overline{\mathbb{R}}$ be a measurable functions.\\

\noindent Question (1). Let us remind the definitions of the following limits, as $\rightarrow +\infty$ :\\

\noindent If $f(\omega)=+\infty$,
$$
(f_n(\omega) \rightarrow f(\omega)=+\infty) \Leftrightarrow (\forall \eta>0, \ exists n\geq 0 : \ \forall p\geq n, \ f_n(\omega)\geq \eta; \ \ (DEF1)
$$

\bigskip \noindent If $f(\omega)=-\infty$,
$$
(f_n(\omega) \rightarrow f(\omega)=-\infty) \Leftrightarrow (\forall \eta>0, \ exists n\geq 0 : \ \forall p\geq n, \ f_n(\omega)\leq -\eta;  \ \ (DEF2)
$$

\bigskip \noindent If $f(\omega) \in \mathbb{R}$,

$$
(f_n(\omega) \rightarrow f(\omega) \in \mathbb{R}) \Leftrightarrow (\forall \varepsilon>0, \ exists n\geq 0 : \ \forall p\geq n, \ |f_n(\omega)-f_n(\omega)|\leq \varepsilon.  \ \ (DEF3)
$$

\bigskip \noindent Set 

\begin{equation*}
A1=(f=+\infty ) \bigcap \left( \bigcap_{\eta>0} \bigcup_{n\geq 1} \bigcap_{p\geq n} (f_{p} \geq \eta) \right),
\end{equation*}

\bigskip\ 
\begin{equation*}
B1=(f=-\infty ) \bigcap \left( \bigcap_{\eta >0} \bigcup_{n\geq 1} \bigcap_{p\geq n} (f_{p} \leq -\eta) \right)
\end{equation*}

\bigskip \noindent and \\

\begin{equation*}
C1=(|f|<+\infty ) \bigcap \left( \bigcap_{\varepsilon >0} \bigcup_{n\geq 1} \bigcap_{p\geq n} (|f_{p}-f| \leq \varepsilon) \right).
\end{equation*}

\noindent By using decomposition

$$
(f_n \rightarrow f)=(f_n \rightarrow f)\cap (f=+\infty) + (f_n \rightarrow f)\cap (f=-\infty) + (f_n \rightarrow f)\cap (|f|<+\infty),
$$

\noindent Show that

$$
(f_n \rightarrow f)=A1 + B1 + C1.
$$

\bigskip \noindent Since the values of $\eta$ and $\varepsilon$ are not countable, can you conclude that $(f_n \rightarrow f)$ is measurable?\\

\noindent Now, by using Exercise 1, show that 

$$
(f_n \rightarrow f)=A + B + C,
$$

\noindent with :

\begin{equation*}
A=(f=+\infty ) \bigcap \left( \bigcap_{k\leq 1} \bigcup_{n\geq 1} \bigcap_{p\geq n} (f_{p} \geq k) \right),
\end{equation*}

\bigskip

\begin{equation*}
B=(f=-\infty ) \bigcap \left( \bigcap_{k\leq 1} \bigcup_{n\geq 1} \bigcap_{p\geq n} (f_{p} \leq -k) \right)
\end{equation*}

\bigskip \noindent and \\

\begin{equation*}
C=(|f|<+\infty ) \bigcap \left( \bigcap_{k \leq 1} \bigcup_{n\geq 1} \bigcap_{p\geq n} (|f_{p}-f| \leq 1/k) \right),
\end{equation*}

\noindent and conclude that $(f_n \rightarrow f)$ is measurable.\\

\bigskip \noindent Question (2). Show that $f_n \rightarrow f$ \textit{a.e.} if and only if the three equations hold :

\begin{equation*}
m\left((f=+\infty ) \bigcap \left( \bigcup_{k\leq 1} \bigcap_{n\geq 1} \bigcup_{p\geq n} (f_{p} < k) \right)\right)=0;
\end{equation*}

\begin{equation*}
m\left((f=-\infty ) \bigcap \left( \bigcup_{k\leq 1} \bigcap_{n\geq 1} \bigcup_{p\geq n} (f_{p} > -k) \right)\right) =0;
\end{equation*}

\bigskip \noindent and \\

\begin{equation*}
m\left( (|f|<+\infty ) \bigcap \left( \bigcup_{k \leq 1} \bigcap_{n\geq 1} \bigcup_{p\geq n} (|f_{p}-f| > 1/k) \right) \right)=0.
\end{equation*}

\bigskip \noindent Question (3). If $f$ is \textit{a.e}, show that $f_n \rightarrow f$ \textit{a.e.} if and only if

\begin{equation*}
m\left( (|f|<+\infty ) \bigcap \left( \bigcup_{k \leq 1} \bigcap_{n\geq 1} \bigcup_{p\geq n} (|f_{p}-f| > 1/k) \right) \right)=0.
\end{equation*}

\noindent \bigskip \noindent Question (3). Let $f$ still be finite \textit{a.e.} Denote

$$
C_{\infty} = \bigcup_{k \leq 1} \bigcap_{n\geq 1} \bigcup_{p\geq n} (|f_{p}-f| > 1/k) =0.
$$

\noindent and for $k\leq 1$,

$$
C_k=\bigcap_{n\geq 1} \bigcup_{p\geq n} (|f_{p}-f| > 1/k).
$$

\bigskip Show that $f_n \rightarrow f$ if and only if, for any $k\geq 1$,

$$
m(C_k)=0.
$$
 
\noindent \bigskip \noindent Question (4). Let $f$ still be finite \textit{a.e.}. Combine the previous result to get that  the results of Exercise 1 to get : $f_n \rightarrow f$ if and only if, for any 
$\varepsilon>0$, $m(D_{\varepsilon})=0$, where

$$
D_\varepsilon=\bigcap_{n\geq 1} \bigcup_{p\geq n} (|f_{p}-f| > \varepsilon).
$$
 
\bigskip \noindent \bigskip \noindent Question (5). Let  $f_{n})_{n\geq 1}$ and $f$ \textit{a.e.} finite and let $m$ be a finite measure. Set

$$
D_{n,\varepsilon}=\bigcup_{p\geq n} (|f_{p}-f| > \varepsilon).
$$ 

\bigskip \noindent Remark that the sequence $(D_{n,\varepsilon})_{n\geq 0}$ is non-increasing and that $(|f_n-f|>\varepsilon) \subset D_{n,\varepsilon}$. Use the continuity of the measure $m$ 
(Point 04.10 of Doc 04-01, page \pageref{doc04-10} in Chapter \ref{04_measures}) to show that

$$
(f_n \rightarrow f) \Rightarrow (f_n \rightarrow_m f).
$$

\bigskip \noindent \textbf{Exercise 3}. \label{exercises03_doc06-04}\\

\noindent \textbf{General Hint}. In each question, use the Measure Space Reduction property (see pages \pageref{msrp} and \pageref{msrpA}) and combine to others assumptions to place yourself in a space
$\Omega_0$, whose complement is a null-set, and on which all the \textit{a.e.} properties hold for each $\omega \in \Omega_0$.\\

\noindent Let $(f_{n})_{n\geq 1}$ and $(g_{n})_{n\geq 1}$ be sequences of real-valued measurable functions such that $f_{n}\rightarrow f$ $a.e.$ and $g_{n}\rightarrow g$ $a.e.$. 
Let $a$ and $b$ be finite real numbers. Let $H(x,y)$ a continuous function of $(x,y) \in D$, where $D$ is an open set of $\mathbb{R}^2$. Show that the following assertion hold true : \\

\noindent \textbf{(1)} The \textit{a.e.} limit of $(f_{n})_{n\geq 1}$, if it exists, is unique \textbf{a.e.}.\\

\noindent \textbf{(2)} If $af+bg$ is \textit{a.e.} defined, then $af_{n}+bg_{n}\rightarrow af+bg$ $a.e$.\\

\noindent \textbf{(3)} If $fg$ is \textit{a.e.} defined, then $f_{n}g_{n}\rightarrow fg$ $a.e$\\

\noindent \textbf{(4)} If $f/g$ is \textit{a.e.} defined, (that is $g$ is $a.e$ nonzero and $f/g \neq \frac{\pm \infty}{\pm \infty}$ \textit{a.e.}), then

$$
f_{n}/g_{n}\rightarrow f/g, \ a.e.
$$

\noindent \textbf{(5)} If $(f_{n}, g_n)_{n\geq 1} \subset D$ \textit{a.e.} and $(f,g) \in D$ \textit{a.e.}, then

$$
H(f_{n},g_{n})\rightarrow H(f,g), \ a.e.\\
$$

\bigskip \bigskip \noindent \textbf{SOLUTIONS}.\\

\noindent (1) Suppose that $f_n$ converges \textit{a.e.} to both $f_1$ and $f_2$. Denote

$$
\Omega_0^c=(f_n \nrightarrow f_{(1)}) \cap (f_n \nrightarrow f_{(2)}).
$$ 

\bigskip \noindent It is clear that $\Omega_0^c$ is a null-set. Besides, for all $\omega \in \Omega_0$,

$$
f_n(\omega) \rightarrow f_{(1)}(\omega) \ \ and \ \ f_n(\omega) \rightarrow f_{(2)}(\omega).
$$

\bigskip \noindent We get, by the uniqueness of the limit in $\overline{\mathbb{R}}$, we get that $f_{(1)}=f_{(2)}$ on $\Omega_0$, that is outside a null-set. Hence $f_{(1)}=f_{(2)}$ \textit{a.e.}. $\square$\\

\bigskip  \noindent (2) By combining the hypotheses, we can find a measurable space sub-space $\Omega_0$ of $\Omega$ whose complement is a $m$-null and, on which, we have $f_{n}\rightarrow f$ $a.e.$ and $g_{n}\rightarrow g$ and $af+bg$ is well defined.  Then for any $\omega \in \Omega_0$, $af_n(\omega)+bg_n(\omega)$ is well defined for large values of $n$ and hence

$$
(af_n+bg_n)(\omega) \rightarrow (af+bg)(\omega).
$$

\bigskip \noindent Hence $af_n+bg_n$ converges to $af+bg$ outside the null-set $\Omega_0^c$. $\square$\\

\bigskip \noindent (3) By combining the hypotheses, we can find a measurable space sub-space $\Omega_0$ of $\Omega$ whose complement is a $m$-null and, on which, we have $f_{n}\rightarrow f$ $a.e.$ and $g_{n}\rightarrow g$ and $fg$ is well defined.  Then for any $\omega \in \Omega_0$, $f_n(\omega) g_n(\omega)$ is well defined for large values of $n$ and hence

$$
(f_n g_n)(\omega) \rightarrow (fg)(\omega).
$$

\bigskip \noindent Hence $f_n g_n$ converges to $fg$ outside the null-set $\Omega_0^c$. $\square$\\

\bigskip \noindent (4) Repeat the same argument.\\

\noindent (5)  By combining the hypotheses, we can find a measurable space sub-space $\Omega_0$ of $\Omega$ whose complement is a $m$-null and, on which, we have $f_{n}\rightarrow f$ $a.e.$ and $g_{n}\rightarrow g$ and all the couples $(f_n, g_n)$'s are in $D$ and $(f,g)$ is in D. Then for any $\omega \in \Omega_0$, $(f_n(\omega), g_n(\omega))$ converges $(f(\omega), g(\omega))$ in $D$ By continuity of $H$, we get

$$
(\forall \omega \in \Omega_), \ H(f_n(\omega), g_n(\omega)) \rightarrow H(f(\omega), g(\omega)).
$$

\bigskip \noindent Hence, $H(f_n, g_n)$ converges to $H(f, g)$ outside the null-set $\Omega_0$. $\blacksquare$\\

\bigskip \noindent \textbf{Exercise 4}. \label{exercises04_doc06-04}\\

\noindent Let $(f_{n})_{n\geq 1}$ be a sequence of \textit{a.e.} finite measurable real-valued functions.\\

\noindent Question (1) By using the definition of a Cauchy sequence on $\mathbb{R}$ and by using the techniques in Exercise 1, show that :   $(f_{n})_{n\geq 1}$ is an $m$-\textit{a.e.} Cauchy sequence if and only if

$$
m\left( \bigcup_{k\geq 1} \bigcap_{n\geq 1} \bigcup_{p\geq n} \bigcup_{q\geq n} (|f_p-f_q|>1/k)\right)=0. \ \ (CAU01)
$$ 

\noindent Question (2). Show that (CAU01) is equivalent to each of the two following assertions.\\

$$
m\left( \bigcup_{k\geq 1} \bigcap_{n\geq 1} \bigcup_{p\geq 0} \bigcup_{q\geq 0} (|f_{n+p}-f_{n+p}| > 1/k)\right)=0. \ \ (CAU02)
$$ 

$$
m\left( \bigcup_{k\geq 1} \bigcap_{n\geq 1} \bigcup_{p\geq 0}  (|f_{n+p}-f_{n}| > 1/k)\right)=0. \ \ (CAU03)
$$

\bigskip \noindent Question (3). Show that $(f_{n})_{n\geq 1}$ is an $m$-\textit{a.e.} Cauchy sequence if and only if $f_n$ converges \textbf{a.e.} to an \textbf{a.e.} measurable function $f$.\\

\bigskip \noindent \textbf{Exercise 5}. \label{exercises05_doc06-04}\\

\noindent In all this exercises, the functions $f_n$, $g_n$, $f$ and $g$ are \textbf{a.e.} finite.\\

\noindent Question (1). Show that the limit in measure is unique \textbf{a.e.}.\\

\noindent Hint : Suppose $f_n \rightarrow_m f_{(1)}$ and $f_n \rightarrow_m f_{(2)}$. By Formula (OS1) in Exercise 1 :
$$
(|f_1-f_2|> 1/k) \subset (|f_n-f_{(1)}|> 1/(2k)) \cup (|f_n-f_{(2)}|> 1/(2k)) 
$$

\noindent Use also (See Formula (C) in Exercise 7 in DOC 04-06, page \pageref{exercise07_sol_doc04-06}, Chapter \ref{04_measures}).\\

$$
(f_1-f_2 \neq 0)= \bigcup_{k\geq 1} (|f_1-f_2|> 1/k).
$$

\bigskip \noindent Question (2). Let $c$ be finite real number. Show that the sequence of \textit{a.e.} constant functions, $f_{n}=c$ \textit{a.e.}, converges in measure to $f=c$.\\

\bigskip \noindent Question (3).  Let $f_{n}\rightarrow _{m} f$ and $g_{n}\rightarrow _{m}g,$ $a\in \mathbb{R}$. Show that\\

\noindent \textbf{(a)} $f_{n}+g_{n}\rightarrow_{m}f+g$.\\

\noindent \textbf{(b)} $af_{n} \rightarrow af$\\

\noindent Question (4).  Suppose that  $f_{n}\rightarrow_m A$ and $g_{n}\rightarrow_m B$, where $A$ and $B$ are finite real numbers. Let $H(x,y)$ a continuous function of $(x,y) \in D$, where $D$ is an open set of $\mathbb{R}^2$. Let $a$ and $b$ be finite real numbers. Show the following properties.\\

\noindent \textbf{(a)} $af_{n}+bg_{n}\rightarrow_m aA+bB$.\\

\noindent \textbf{(b)} $f_{n}g_{n}\rightarrow_m AB$.\\

\noindent \textbf{(c)} If $B \neq 0$,  then

$$
f_{n}/g_{n}\rightarrow_m A/B.
$$

\bigskip \noindent \textbf{(d)} If $(f_{n}, g_n)_{n\geq 1} \subset D$ \textit{a.e.} and $(A,B) \in D$, then

$$
H(f_{n},g_{n})\rightarrow_m H(f,g).
$$

\bigskip \noindent \textbf{Exercise 6}. \label{exercises06_doc06-04}\\

\noindent Let $(f_{n})_{n\geq 1}$ be a sequence of \textit{a.e.} finite and measurable real-valued functions.\\

\noindent Question (a) Suppose that the $(f_{n})_{n\geq 1}$ converges in measure to an \textit{a.e.} finite, measurable and real-valued function $f$. Show that the sequence is
a Cauchy sequence in measure.\\

\noindent Question (b) Suppose that $(f_{n})_{n\geq 1}$ is a Cauchy sequence in measure and that there exists a sub-sequence $(f_{n_k})_{k\geq 1}$ of $(f_{n})_{n\geq 1}$ converging in measure to an \textit{a.e.} function $f$. Show that $(f_{n})_{n\geq 1}$ converges in measure to $f$.

\noindent Question (c) Suppose that $(f_{n})_{n\geq 1}$ is a Cauchy sequence in measure and that there exists a sub-sequence $(f_{n_k})_{k\geq 1}$ of $(f_{n})_{n\geq 1}$ converging \textbf{a.e.} to an \textit{a.e.} function $f$. Suppose in addition that the measure is finite. Show that $(f_{n})_{n\geq 1}$ converges in measure to $f$.\\

\noindent \textsl{Hint}. Combine Question (b) of the current exercise and Question (4) if Exercise 2.\\

\bigskip \noindent \textbf{Exercise 7}. \label{exercises07_doc06-04}\\

\noindent Let $(f_{n})_{n\geq 1}$ be a sequence of \textit{a.e.} finite and measurable real-valued functions.\\

\noindent Question (a). Suppose that $(f_{n})_{n\geq 1}$ is Cauchy in measure. Show that $(f_{n})_{n\geq 1}$ possesses a sub-sequence $(f_{n_k})_{k\geq 1}$ of $(f_{n})_{n\geq 1}$ converging both \textbf{a.e.} and in measure to an \textit{a.e.} function $f$.\\

\noindent \textbf{Detailed Hints}.\\

\noindent \textbf{Step 1}. From the hypothesis : for any $\varepsilon>0$, for all $\eta>0$,

$$
\exists n\geq 1, \forall p\geq n, \ \forall q \geq n, \ m(|f_q - f_q| > \varepsilon) < \eta. \ \ (CAUM01)
$$

\bigskip \noindent Remark (without doing any supplementary work) that Formula (CAUM01) is equivalent to each of (CAUM02), (CAUM03) and (CAUM04), exactly as you already proved that (CAU01) is equivalent to each of (CAU02) and (CAU03) in Exercise 4 : 

$$
\exists n \geq 0, \forall p\geq 0, \ \forall q \geq 0, \ m(|f_{n+p} - f_{n+q}| > \varepsilon) < \eta; \ \ (CAUM02)
$$

$$
\exists n\geq 0, \forall p\geq 0,  \ m(|f_{n+p} - f_n| >\varepsilon) < \eta; \ \ (CAUM03)
$$

$$
\exists n\geq 0, \forall N\geq n, \forall p\geq 0,  \ m(|f_{N+p} - f_N| >\varepsilon) < \eta. \ \ (CAUM04)
$$

\bigskip \noindent From Formula (CAUM04) construct an increasing sequence $(n_k)_{k\geq 1}$, that is

$$
n_1 <n_2 < ... < n_k < n_{k+1} < ...
$$

\bigskip \noindent such that for all $k\geq 1$,

$$
\forall N\geq n_k, \ \forall p\geq 0, \ \ m\left(|f_{N+p} - f_N| > \frac{1}{2^k} \right) < \frac{1}{2^k}. \ (FI)
$$

\bigskip \noindent proceed as follows : For each $k=1$, take $\varepsilon=\eta=1/2$ and take $n_1$ as the value of $n$ you found in (CAUM04). For $k=2$, take $\varepsilon=\eta=1/2^2$ and take $n(2)$ as the value of $n$ you found in (CAUM04) and put $n_2=\max(n_1+1,n(2)$. It is clear that $n_1<n_2$ and, since

$$
\forall N\geq n(2), \  \forall p \geq 0, \ \ m\left(|f_{N+p} - f_{N)}| > \frac{1}{2^2} \right) < \frac{1}{2^2}, \ \ (FI2)
$$

\bigskip \noindent you will have

$$
N\geq n_2, \  \forall p \geq 0, \ \ m\left(|f_{N+p} - f_{N}| > \frac{1}{2^2} \right) < \frac{1}{2^2}, \ \ (FI2P)
$$

\bigskip \noindent since values of $N$ greater that $n_2$ in (FI2P) are also greater that $n(2)$ in $(FI2)$. Next, proceed from $k=2$ to $k=3$ as you did from $k=1$ to $k=2$.\\

\noindent \textbf{Step 2}. Set

$$
A_k=\left( \left| f_{n_{k+1}} - f_{n_{k}} \right| > \frac{1}{2^k}\right), \ k\geq 1, \ and \ B_N=\bigcup_{k\geq N} A_k, \ N\geq 0.
$$

\bigskip \noindent Show that

$$
m(B_N) < \frac{1}{2^{N-1}}.
$$

\bigskip \noindent Fix $\varepsilon>0$ and take $N$ satisfying $2^{-(N-1)} < \varepsilon$. Show that : on $B_n^c$ we have for all $p\geq 0$
$$
\left| f_{n_{N+p}} - f_{n_{N}}\right| \geq \frac{1}{2^{N-1}} \leq \varepsilon. \ \ (FI3)
$$
 
\bigskip \noindent where you compose the difference $f_{n_{N+p}} - f_{n_{N}}$ from the increments $f_{n_{N+i}} - f_{n_{N+i-1}}$, $i=1, ..., p$, and next extend the bound to all the
increments for $k\geq N$.\\

\noindent Next using this with the obvious inclusion

$$
\bigcap_{q\geq 0} \bigcup_{p\geq 0} \left( \left| f_{n_{q+p}} - f_{n_{q}}\right| > \varepsilon \right)  \subset  \bigcup_{p\geq 0} \left( \left| f_{n_{N+p}} - f_{n_{N}}\right| > \varepsilon \right),
$$

\bigskip \noindent show that 

$$
m\left(\bigcap_{q\geq 0} \bigcup_{p\geq 0} \left( \left| f_{n_{N+p}} - f_{n_{N}}\right| \geq \varepsilon \right) \right) \leq m(B_N).
$$

\bigskip \noindent Deduce from this that for all $\varepsilon$,

$$
m\left(\bigcap_{q\geq 0} \bigcup_{p\geq 0} \left( \left| f_{n_{q+p}} - f_{n_{q}}\right| > \varepsilon \right) \right)
\leq m\left(\bigcap_{q\geq 0} \bigcup_{p\geq 0} \left( \left| f_{n_{q+p}} - f_{n_{q}}\right| \geq \varepsilon \right) \right)=0.
$$

\bigskip \noindent Apply this to $\varepsilon=1/r$, $r\geq 1$ integer, compare with Formula (CAU03) in Exercise 4 above, to see that sub-sequence $(f_{n_k})_{k\geq 1}$ is an \textit{a.e.} Cauchy sequence. By  Question 3 in Exercise 4, consider the \textit{a.e} finite function to which $(f_{n_k})_{k\geq 1}$ converges \textit{a.e.}\\

\noindent Consider the measurable sub-space  $\Omega_0$ whose complement is a null-set and on which $(f_{n_k})_{k\geq 1}$ converges \textit{a.e.} to $f$.\\

\noindent \textbf{Step 3}. Consider Formula (FI3). Place yourself on $\Omega_0 \cap B_N^c$ and let $p$ go to $\infty$. What is the formula you get on $\Omega_0 \cap B_N^c$?\\

\noindent Deduce the formula

$$
m\left( \left|f - f_{n_{N}}\right| > \varepsilon \right) \leq m\left( \left|f - f_{n_{N}}\right| \geq  \varepsilon \right) \leq m(B_N),
$$

\bigskip \noindent for all $N\geq 1$  satisfying $2^{-(N-1)} < \varepsilon$. Deduce that $f_{n_{N}}$ converges to $f$ in measure as $N \rightarrow +\infty$.\\

\noindent \noindent \textbf{Step 4}. By using Question (b) in Exercise 6, conclude that $f_{n}$ converges to $f$ in measure as $n \rightarrow +\infty$.\\

\noindent Write your conclusion.\\

\noindent Question (b). Combine this with the results of Exercise 6 above to state that a sequence $(f_{n})_{n\geq 1}$ is a Cauchy in measure if and only if it converges to a an \textit{a.e.} finite and measurable real-valued function.\\

\noindent \LARGE \textbf{DOC 06-05 : Convergence Theorems - Exercises.} \label{doc06-05}\\
\bigskip
\Large

\bigskip \noindent \textbf{NB}. In this document and in the sequel of the textbook, the Measure Space Reduction Principle (see pages \pageref{msrp} and \pageref{msrpA}) which consists in placing ourselves in a sub-measurable space $\Omega_0$ whose complement is a null-set with respect to the working measure, and on which all the countable \textit{a.e.} properties simultaneously hold for $\omega \in \Omega_0$. By working on the induced measure on $\Omega_0$, the values of the measures of sets and the values of the integrals of real-valued and measurable function remain unchanged.\\

\bigskip \noindent \textbf{Exercise 1}. \label{exercise1_doc06-05}

\noindent Question (1). Let $(f_{n})_{n\geq 0}$ be a non-decreasing sequence of \textit{a.e.} non-negative real-valued and measurable functions. Show that
\begin{equation*}
\lim_{n \uparrow +\infty} \int f_{n}dm = \int \lim_{n \uparrow +\infty} f_n \ dm.
\end{equation*}

\noindent \textit{Hint}.\\

\noindent (a) First denote $f=\lim_{n \uparrow +\infty} f_n$.\\

\noindent By  \textit{Point (02.14) in Doc 02-01 in Chapter \ref{02_applimess}} (see page \pageref{doc02-01}), consider for each $n\geq 0$, a non-decreasing sequence of non-negative elementary functions $(f_{n,k})_{k\geq 0}$ such that $f_{n,k} \uparrow f_n$ as $k\uparrow +\infty$.\\

\noindent Next set

$$
g_k= \max_{k \leq n} f_{n,k}, \ n\geq 0.
$$

\bigskip \noindent  By \textit{Point (02.13) in Doc 02-01 in Chapter \ref{02_applimess}} (see page \pageref{doc02-01}), see that the $g_k$ are non-negative elementary functions.\\

\noindent Write $g_k$ and $g_{k+1}$ in extension and say why the sequence $(g_k)_{k\geq 0}$ is non-decreasing.\\

\noindent (b) Show that we have 

$$
\forall 0 \geq n, \ f_{n,k} \leq g_k \leq f_k, (I1)
$$

\bigskip \noindent  and deduce from it, by letting first $k\rightarrow +\infty$ and next $n\rightarrow +\infty$, that

$$
\lim_{k\rightarrow +\infty} g_k=f. (I2)
$$

\bigskip \noindent  Why all the limits you consider exist?\\

\noindent (c) Now by taking the integrals in (I1), by letting first $k\rightarrow +\infty$ and next $n\rightarrow +\infty$, and by using the definition of the integral of 
$\int \lim_{k\rightarrow +\infty} g_k \ dm$ and by raking (I2) into account, conclude that

$$
\lim_{k\rightarrow +\infty} \int f_k \ dm= \int \lim_{k\rightarrow +\infty}  f_k \ dm. \ \square.
$$

\bigskip \noindent  Question (2) Extend the result of Question (1) when $(f_{n})_{n\geq 0}$ is a non-decreasing sequence of functions bounded below \textit{a.e.} by an integrable function $g$, that is : $g \leq f_n$, for all $n\geq 0$. (Actually, it is enough that $g \leq f_n$ for $n\geq n_0$, where $n_0$ is some integer).\\

\noindent \textit{Hint}. Consider $h_n=f_n-g$ and justify that you are in the case of Question (1). Next, do the mathematics.\\

\bigskip \noindent  Question (2) Extend the result of Question (2) when $(f_{n})_{n\geq 0}$ is a non-decreasing sequence of functions such that there exists $n_0\geq 0$, such that $\int f_{n_0} \ dm >-\infty$.\\

\noindent \textit{Hint}. Discuss around the two cases : $\int f_{n_0} \ dm=+\infty$ and $\int f_{n_0} \ dm<+\infty$.\\

\bigskip \noindent \textbf{Exercise 2}. (Fatou-Lebesgue's lemma and Dominated Convergence Theorem) \label{exercise02_doc06-05}\\

\noindent Question (1). (Fatou-Lebesgue's lemma) Let $(f_{n})_{n\geq 0}$ be a sequence real-valued and measurable functions \textit{a.e.} bounded \textbf{below} by an integrable function $h$, that is 
$h \leq f_n$ \textit{a.e.} for all $n\neq 0$. Show that

\begin{equation*}
\liminf_{n\rightarrow +\infty}\int f_{n}dm \geq \int \liminf_{n\rightarrow +\infty }f_{n}dm. 
\end{equation*}

\noindent \textit{Hints}.\\

\noindent (a) Denote for all $n\geq 0$

$$
g_n=\inf_{k\geq n} (f_k-h).
$$

\noindent The $g_n$ are they defined? Why? What are their sign? What is the limit of $g_n$? What is its monotonicity, if any?\\

\noindent Justify and apply the Monotone Convergence Theorem.\\

\noindent Compare $\int \inf_{k\geq n} (f_k-h) \ dm$ with each $\int (f_k-h) \ dm$, $k\geq n$ and deduce a comparison with $\inf_{k\geq n} \int (f_k-h) \ dm$.\\

\noindent Conclude. Show your knowledge on superior limits and inferior limits.\\

\noindent Question (2). Let $(f_{n})_{n\geq 0}$ be a sequence real-valued and measurable functions \textit{a.e.} bounded \textbf{above} by an integrable function $h$, that is 
$f_n \leq h$ \textit{a.e.} for all $n\neq 0$. Then, we have

$$
\limsup_{n\rightarrow +\infty}\int f_{n}dm \leq \int \limsup_{n\rightarrow +\infty }f_{n}dm.
$$

\bigskip \noindent  \textit{Hints}. Take the opposites of all the functions in this question and apply Question (1).\\

\noindent Question (3). Let $(f_{n})_{n\geq 0}$ be a sequence real-valued and measurable functions \textit{a.e.} such that :\\

\noindent (a) it is bounded  by an integrable function $g$, that is $|f_n| \leq g$ \textit{a.e.} for all $n\geq 0$,\\

\noindent (b) it converges in measure or \textit{a.e.} to a function $f$.\\

\noindent  Then $f$ is integrable and the sequence $(\int f_{n} \ dm)_{n\geq 0}$ has a limit with

$$
\lim_{n\rightarrow +\infty}\int f_{n} \ dm = \int f \ dm. \ (D1)
$$

\bigskip \noindent and

$$
\int |f| \ dm = \int g \ dm. \ (D2)
$$

\noindent \textit{Hints}. Do two cases.\\

\noindent Case 1 : $f_n \rightarrow f$ \textit{a.e}. Remark the your sequence is bounded below and above. Apply Questions (1) and (2) and replace $\lim_{n\rightarrow +\infty} f_n$ by its value. Combine the two results at the light of the natural order of the limit inferior and the limit superior to conclude.\\

\noindent Case 2 : $f_n \rightarrow_m f$. Use the\textit{ Prohorov criterion given in Exercise 4 in Doc 11-01 in Chapter \ref{11_appendix} (See page \pageref{doc11-01})} in the following way.\\

\noindent Consider an adherent point $\ell$ of the sequence $(\int f_n \ dm)_{n\geq 0}$, limit a sub-sequence $(\int f_{n_k} \ dm)_{k\geq 0}$.\\

\noindent Since we still have $f_{n_k} \rightarrow_m f$ as $k\rightarrow +\infty$, consider a sub-sequence $(f_n{_{k_j}})_{j\geq 0}$ of $(f_{n_k})_{k\geq 0}$ converging both \textit{a.e.} and in measure to a finite function $h$ (See Exercise 7 in Doc 06-04, page \pageref{doc06-04}). Justify that $h=f$ \textit{a.e.} and apply Case 1 to the sub-sub-sequence $(f_{n_{k_j}})_{j\geq 0}$ and conclude that 

$$
\ell= \int f \ dm.
$$

\bigskip \noindent Apply the Prohorov criterion and conclude.\\

\bigskip \noindent \textbf{Exercise 3}. (Young Dominated Convergence Theorem) \label{exercise-05}\\

\noindent Let $(f_{n})_{n\geq 0}$ be a sequence real-valued and measurable functions \textit{a.e.} such that :\\

\noindent (1) there exists a sequence of integrable functions $(h_n)_{n\geq 0}$ such that is $|f_n| \leq h_n$ \textit{a.e.} for all $n\geq 0$,\\

\noindent (2) the sequence $(h_n)_{n\geq 0}$ converges to an integrable $h$ such that 
$$
\lim_{n\rightarrow +\infty}\int h_{n} \ dm = \int h \ dm,
$$

\bigskip \noindent (3) the sequence $(f_{n})_{n\geq 0}$ converges in measure or \textit{a.e.} to a function $f$.\\

\noindent  Then $f$ is integrable and the sequence $(\int f_{n} \ dm)_{n\geq 0}$ has a limit with

$$
\lim_{n\rightarrow +\infty}\int f_{n} \ dm = \int f \ dm.
$$

\bigskip \noindent and

$$
\int |f| \ dm = \int h \ dm.
$$

\bigskip \noindent \textit{Hint}. At the the method used in Exercise 2 above, you only need to prove an analog of the Fatou-Lebesgue lemma of the form : Let $(f_{n})_{n\geq 0}$ be a sequence real-valued and measurable functions \textit{a.e.} bounded \textbf{below} by a sequence of integrable functions $(h_n)_{n\geq 0}$, which converges to an integrable $h$ so that  
$$
\lim_{n\rightarrow +\infty}\int h_{n} \ dm = \int h \ dm.
$$

\bigskip \noindent Then we have

\begin{equation*}
\liminf_{n\rightarrow +\infty}\int f_{n}dm \geq \int \liminf_{n\rightarrow +\infty }f_{n}dm. (D3)
\end{equation*}

\bigskip \noindent Prove (D3) by repeating that of Question (1) in Exercise 2, but be careful since you will not have the linearity of the limit inferior. The method will lead you to the formula

$$
\int \liminf_{n\rightarrow +\infty} (f_n-h_n) \ dm \leq  \liminf_{n\rightarrow +\infty} \int (f_n-h_n) \ dm. \ \ (Y3)
$$  

\bigskip \noindent From there, use the super-additivity in the left-hand, use in the right-hand the the formula : for all $\eta>0$, there exists $N\geq 0$ such that for all $n\geq N$, 
$$
\int h_n \ dm \geq \int h \ dm.
$$

\bigskip \noindent Combine the two formula $n\geq N$, let $n\rightarrow +\infty$ and next $\eta\rightarrow 0$, to conclude.\\

\noindent \LARGE \textbf{DOC 06-06 : Applications of Convergence Theorems - Exercises} \label{doc06-06}\\
\bigskip
\Large

\textbf{NB}. Recall that a continuous limit (in $\overline{\mathbb{R}}$, in the form 

$$
x(t) \rightarrow y  \ \ as \ \ t \rightarrow s
$$ 

\bigskip \noindent may be discretized and shown to be equivalent to : for any sequence $(t_n)_{n\geq 0}$ such that $t_n\rightarrow s$ as $n\rightarrow +\infty$, we have

$$
x(t_n) \rightarrow y \ \ as \ \ n\rightarrow +\infty.
$$

\bigskip \noindent  Using the discretized from is more appropriate with measurability.\\

\noindent \textbf{Exercise 1}. \label{exercise01_doc06-06} (Convergence of sums of applications)\\

\noindent Question (1). Let $(f_{n})$ be sequence of \textit{a.e.} non-negative measurable functions. Use the MCT and show that 
\begin{equation*}
\sum_{n=0}^{+\infty} \biggr( \int f_{n} \  dm \biggr) =\int \biggr( \sum_{n}^{+\infty} f_{n} \ dm \biggr). \ \ (S01)
\end{equation*}

\noindent Question (2). Let be $(a_{n})_{n\geq 0}=(\left( a_{n,p}\right) _{(n\geq 0,p\geq 0)}$ a rectangular sequence of non-negative real numbers. Use the MCT to the counting measure on $\mathbb{N}$ and show that

\begin{equation*}
\sum_{n\geq 0} \biggr( \sum_{p\geq 0}a_{n,p} \biggr)=\sum_{p\geq 0} \biggr(\sum_{n\geq 0}a_{n,p}\biggr). \ \ (S02)
\end{equation*}

\noindent Question (3) Let $(f_{n})$ be sequence of measurable functions such that%
\begin{equation*}
\int \sum_{n\geq 0} \left\vert f_{n}\right\vert  dm<+\infty.
\end{equation*}

\noindent Use the DCT, take into account Question (1) and establish that 

\begin{equation*}
\sum_{n\geq 0} \biggr( \int f_{n} \ dm \biggr)  =\int \biggr( \sum_{n\geq 0} f_{n} \biggr) \ dm.
\end{equation*}

\noindent Question (4). Let $(\left( a_{n,p}\right) _{(n\geq 0,p\geq 0)}$ a rectangular sequence such that
\begin{equation*}
\sum_{p\geq 0} \sum_{n\geq 0}\left\vert a_{n,p}\right\vert <\infty.
\end{equation*}

\bigskip \noindent Apply Question (3) on the counting measure on $\mathbb{N}$ and conclude
\begin{equation*}
\sum_{p\geq 0} \sum_{n\geq 0} a_{n,p}=\sum_{n\geq 0} \sum_{p\geq 0} a_{n,p}.
\end{equation*}

\noindent Question (5). Let $m$ be a countable sum of measures $m_j$, $j\geq 0$ :

$$
m=\sum_{j\geq 0} m_j.
$$

\noindent Show that for any non-negative and measurable function $f$, we have

$$
\int f \ dm = \sum_{j\geq 0} \int f \ dm_j. \ \ (IS)
$$

\noindent Give simple conditions under which Formula (IS) holds for some for  measurable $f$ function.\\

\noindent \textit{Hints}. Describe only the solution at the light of the four steps method.\\

\bigskip \noindent \textbf{Exercise 2}. \label{exercises02_doc06-06} (Continuity of parametrized functions)\\

\noindent Let $f(t,\omega ), \ t\in T$ a family of measurable functions indexed by $t\in I$, where $I=]a,b[, \ a<b$, to fix ideas. Let $t_{0}\in I$.\\

\bigskip \noindent Question(1). (Local continuity) Suppose that for some $\varepsilon >0,$ for $\varepsilon $ small enough
that $]t_{0}-\varepsilon ,t_{0}+\varepsilon \lbrack \subset I)$, there exists $g$ integrable such that

\begin{equation*}
\forall (t\in ]t_{0}-\varepsilon ,t_{0}+\varepsilon \lbrack ),\left\vert
f(t,\omega )\right\vert \leq g,
\end{equation*}

\bigskip \noindent and, for almost every $\omega$, the function
\begin{equation*}
t\hookrightarrow f(t,\omega)
\end{equation*}

\bigskip \noindent is continuous at $t_{0}$. Then the function
\begin{equation*}
F(t)=\int f(t,\omega )dm(\omega)
\end{equation*}

\bigskip \noindent is continuous at $t_{0}$ and%
\begin{equation*}
\lim_{t\rightarrow t_{0}}\int f(t,\omega )dm(\omega )=\int \lim_{t\rightarrow t_{0}}f(t,\omega )dm(\omega )=\int f(t_{0},\omega) \ dm(\omega).
\end{equation*}

\bigskip \noindent Question(2). (Global continuity) Suppose that there exists an integrable function $g$ such that  
\begin{equation*}
\forall (t\in I),\left\vert f(t,\omega )\right\vert \leq g, \ \textit{a.e.}
\end{equation*}

\bigskip \noindent  and that the function
\begin{equation*}
t\hookrightarrow f(t,\omega )
\end{equation*}

\bigskip \noindent is \textit{a.e.} continuous at each $t\in I$. Then the function
\begin{equation*}
F(t)=\int f(t,\omega)dm(\omega )
\end{equation*}

\bigskip \noindent is continuous at each $t\in I$ and for any $u\in I$
\begin{equation*}
\lim_{t\rightarrow u}\int f(t,\omega )dm(\omega )=\int \lim_{t\rightarrow
u}f(t,\omega )dm(\omega )=\int f(u,\omega )dm(\omega).
\end{equation*}

\bigskip \noindent \textbf{Exercise 3}. \label{exercises03_doc06-06} (Differentiability of parametrized functions)\\

\noindent Let $f(t,\omega ),t\in T$ a family of measurable functions indexed by $t\in I$, where $I=]a,b[,a<b$, to fix ideas. Let $t_{0}\in I$.\\

\bigskip \noindent (a) (Local differentiability) Suppose that for almost every $\omega$, the function%
\begin{equation*}
t\hookrightarrow f(t,\omega )
\end{equation*}

\bigskip \noindent is differentiable at $t_{0}$ with derivative
\begin{equation*}
\frac{df(t,\omega )}{dt}|_{t=t_{0}}=f(t_{0},\omega ).
\end{equation*}

\bigskip \noindent and that for some $\varepsilon >0,$ for $\varepsilon $ small enough that $%
]t_{0}-\varepsilon ,t_{0}+\varepsilon \lbrack \subset I),$ there exists $g$ integrable such that, 
\begin{equation*}
\forall (\left\vert h\right\vert \leq \varepsilon ),\left\vert \frac{f(t+h,\omega )-f(t,\omega )}{h}\right\vert \leq g.
\end{equation*}

\bigskip \noindent Then the function
\begin{equation*}
F(t)=\int f(t,\omega )dm(\omega )
\end{equation*}

\bigskip \noindent is differentiable $t_{0}$ and%
\begin{equation*}
F^{\prime }(t_{0})\frac{d}{dt}\int f(t,\omega )dm(\omega )|_{t=t_{0}}=\int 
\frac{df(t,\omega )}{dt}|_{t=t_{0}}dm(\omega ).
\end{equation*}

\bigskip \noindent (b) (Everywhere differentiability) Suppose that $f(t,\omega )$ is $a.e$ differentiable at each t$\in T$ and
there $g$ integrable such that
\begin{equation*}
\forall (t\in I),\left\vert \frac{d}{dt}f(t,\omega)\right\vert \leq g.
\end{equation*}

\bigskip \noindent Then, the function
\begin{equation*}
t\hookrightarrow f(t,\omega)
\end{equation*}

\bigskip \noindent is differentiable a each $t\in I$ and at each $t\in I,$ 
\begin{equation*}
\frac{d}{dt}F(t)=\int \frac{d}{dt}f(t,\omega )dm(\omega).
\end{equation*}

\noindent (c) Let
\begin{equation*}
U(t)=\sum_{n}u_{n}(t)
\end{equation*}

\bigskip \noindent be a convergent series of functions for $\left\vert t\right\vert <R,R>0.$ Suppose
that  each $u_{n}(t)$ is differentiable and there exists $g=(g_{n})$ integrable, that is 
\begin{equation*}
\sum_{n}g_{n}<+\infty,
\end{equation*}

\bigskip \noindent such that
\begin{equation*}
\text{For all }\left\vert t\right\vert <R,\left\vert \frac{d}{dt}
u_{n}(t)\right\vert \leq g.
\end{equation*}

\bigskip \noindent Then, the function $U(t)$ is differentiable and
\begin{equation*}
\frac{d}{dt}U(t)=\sum_{n}\frac{d}{dt}u_{n}(t).
\end{equation*}

\noindent \LARGE \textbf{Doc 06-07 : Lebesgue-Stieljes and Riemann-Stieljes integrals - Exercise}. \label{doc06-07}\\
\bigskip
\Large

\bigskip \noindent The definition of the Riemann-Stieljes Integral is given in Doc 05-04 in Chapter \ref{05_integration}.\\

\noindent We recommend also to read at least the statement of \textit{Exercise 8 in Doc 03-06} \ in Chapter \ref{03_setsmes_applimes_cas_speciaux}, since we are going to use semi-continuous functions.\\

\noindent In all this text, $a$ and $b$ are two real numbers such that $a<b$ and $F: [a,b]\longrightarrow \mathbb{R}$ is
a non-constant and non-decreasing function.\\ 

\bigskip \noindent \textbf{A - Integral on compact sets}.\\

\noindent \textbf{Exercise 1}. \label{exercise01_doc06-07}\\

\noindent Provide a bounded function $f$ on $[a,b]$, which is Lebesgue-Stieljes integrable with respect to any non-decreasing function $F$, and which is not Riemann-Stieljes with respect to any non-decreasing function $F$.\\

\noindent \textit{Hint}. Consider
\begin{equation*}
f=1_{]a,b]\cap \mathbb{Q}}, 
\end{equation*}

\noindent where $\mathbb{Q}$ is the set of rational numbers.\\

\noindent (a) What is the value of $\int_{]a,b]} g \ d\lambda_F$?\\

\noindent (b) Use Riemann-Stieljes sums 
\begin{equation}
S_{n}=\sum_{i=0}^{\ell(n)-1}f(c_{i,n})(F(x_{i+1,n})-F(x_{i,n})).  \label{rsums}
\end{equation}

\noindent and assign to all the $c_{i,n}$'s rational values, and next, irrational values.\\

\noindent (c) Conclude.\\

\bigskip \noindent \textbf{Exercise 2}. \label{exercises02_doc06-07}\\

\noindent Let $F: [a,b]\longrightarrow \mathbb{R}$ be a non-constant and non-decreasing function.\\

\noindent Show that any bounded function $f$ on $[a,b]$ is Riemann-Stieljes integrable on $[a,b]$ if and only it is $\lambda_F$-\textbf{a.e.} continuous with respect to the Lebesgue measure-Stieljes $\lambda_F$.\\

\noindent \textit{\textbf{A - Hint on the part} : If a bounded function $f$ on $[a,b]$ is Riemann-Stieljes integrable on $[a,b]$, then $f$ is $\lambda_F$-\textit{a.e.} continuous with respect to the Lebesgue-Stieljes measure, and its Riemann-Stieljes integral and its Lebesgue-Stieljes integral are equal}.\\

\noindent We are going to use the functions and superior or inferiors limits defined in \textit{Exercise 8 in Doc 03-06} \ in Chapter \ref{03_setsmes_applimes_cas_speciaux} as follows. Put $\varepsilon >0$, 
\begin{equation*}
f^{\ast ,\varepsilon }(x)=\sup \{f(z);z\in \lbrack a,b]\cap ]x-\varepsilon
,x+\varepsilon \lbrack \},
\end{equation*}

\bigskip \noindent and 
\begin{equation*}
f_{\ast }^{\varepsilon }(x)=\inf \{f(z);z\in \lbrack a,b]\cap ]x-\varepsilon.
,x+\varepsilon \lbrack \}.
\end{equation*}

\bigskip \noindent The functions limit and superior limit functions are given by 
\begin{equation*}
f_{\ast }(x)=\lim_{\varepsilon \downarrow 0}f_{\ast }^{\varepsilon }(x)\leq
\lim_{\varepsilon \downarrow 0}f^{\ast ,\varepsilon }(x)=f^{\ast }(x).
\end{equation*}

\bigskip \noindent According to the terminology introduced in \textit{Exercise 8 in Doc 03-06} \ in Chapter \ref{03_setsmes_applimes_cas_speciaux}, we proved that $f_{\ast }$ is lower semi-continuous and that $f^{\ast}$ is upper semi-continuous, and as such, are measurable.\\

\noindent Consider for each $n\geq 1$, a subdivision of $]a,b]$, that partitions  into $]a,b]$ into the $\ell
(n)$ sub-intervals

\begin{equation*}
]a,b]=\sum_{i=0}^{\ell (n)-1}]x_{i,n},x_{i+1,n}].
\end{equation*}

\bigskip \noindent with modulus
\begin{equation*}
m(\pi _{n})=\max_{0\leq i\leq m(n)}(x_{i+1,n}-x_{i+1,n}) \rightarrow 0 as n\rightarrow +\infty.
\end{equation*}

\noindent Define for each $n\geq 0$, for each $i$, $0\leq i \leq \ell(n)-1$,

\begin{equation*}
m_{i,n}=\inf \{f(z),x_{i,n}\leq z<x_{i+1,n}\}\text{ and }M_{i,n}=\sup
\{f(z),x_{i,n}\leq z<x_{i+1,n}\},
\end{equation*}

\begin{equation*}
h_{n}=\sum_{i=0}^{\ell (n)-1}m_{i,n}1_{]x_{i,n},x_{i+1,n}]}\text{ and }%
H_{n}=\sum_{i=0}^{\ell (n)-1}M_{i,n}1_{]x_{i,n},x_{i+1,n}]}
\end{equation*}

\noindent and\\

\begin{equation*}
D=\bigcup_{n\geq 0}\{x_{0,n},x_{1,n},...,x_{\ell (n),n}\}.
\end{equation*}

\noindent (1)  Explain simple sentences why $f^{\ast}$ and $f_{\ast}$, $h_n$ and $H_n$, $n\geq 0$, are bounded by $M$ and why $D$ is countable. Explain why you can choose $D$ as a set of continuity points
of $F$, based on the countability of the set of discontinuity points of $F$ (See Exercise 1 in Doc DOC 03-06, page \pageref{exercise01_doc03-06}).\\

\noindent (2) Fix $x$ in  $x\in \lbrack a,b]\setminus D$, for any $n\geq 1$, there exists $i$, $0\leq i(n)\leq \ell (n)-1$\ such that $x\in
]x_{i(n),n},x_{i(n)+1,n}[$.\\

\noindent Let $\varepsilon >0$\ such that  
\begin{equation*}
]x-\varepsilon ,x+\varepsilon \lbrack \subseteq ]x_{i(n),n},x_{i(n)+1,n}[.
\end{equation*}

\noindent Show that for any $x\in \lbrack a,b]\setminus D$

\begin{eqnarray*}
h^{n}(x) \leq f_{\ast }^{\varepsilon }(x).
\end{eqnarray*}

\noindent and

\begin{eqnarray*}
H_{n}(x) \geq f_{\ast }^{\varepsilon }(x).
\end{eqnarray*}

\noindent (3) Deduce that for any $x\in \lbrack a,b]\setminus D$,

\begin{equation*}
h_{n}(x)\leq f_{\ast }(x)\leq f^{\ast }\leq H_{n}(x).  \ (B1)
\end{equation*}

\noindent (4) Let $\eta>0$. For $n$ fixed. By using the characterization of the suprema and the infima on $\mathbb{R}$, show that the Lebesgue integral of $h_n$ and $H_n$ can be approximated by an 
Riemann-Stieljes sum in the form

\begin{equation*}
\int H_{n}d\lambda =\sum_{i=0}^{\ell
(n)-1}M_{i,n}(x_{i+1,n}-x_{i,n})<\sum_{i=0}^{\ell
(n)-1}f(d_{i,n})(x_{i+1,n}-x_{i,n})+\eta
\end{equation*}

\noindent and 

\begin{equation*}
\int h_{n}d\lambda =\sum_{i=0}^{\ell(n)-1}m_{i,n}(F(x_{i+1,n})-F(x_{i,n}))>\sum_{i=0}^{\ell
(n)-1}f(c_{i,n})(F(x_{i+1,n})-F(x_{i,n}))-\eta.
\end{equation*}

\noindent (5) Mix up all that to prove that

\begin{eqnarray*}
&& \sum_{i=0}^{\ell (n)-1}f(c_{i,n})\left( F(x_{i+1,n})-F(x_{i,n})\right) -\eta\\
&& \leq \int_{\lbrack a,b]\setminus D}f_{\ast }(x)d\lambda \\
&& \leq \int_{[a,b]\setminus D}f_{\ast }(x)d\lambda \\
&& \leq \sum_{i=0}^{\ell (n)-1}f(d_{i,n})(x_{i+1,n}-x_{i,n})+\eta.
\end{eqnarray*}

\noindent In establishing the latter double inequality, do not forget that $D$ is a $\lambda_F$-null set.\\

\noindent (6) Conclude : (i) by letting  $\eta\rightarrow 0$, and : (ii) by considering that $f$ is continuous at $x$ if and only if $f_{\ast}=f^{\ast}(x)$.\\

\bigskip 
\noindent \textit{\textbf{B - Hint on the part} : If a bounded function $f$ on $[a,b]$ is $\lambda_F$-\textit{a.e.} continuous, then it is Riemann-Stieljes integrable on $[a,b]$ exists, and its 
Riemann-Stieljes integral and its Lebesgue-Stieljes  integral are equal}.\\
 
\noindent (1) Denote by $D\subset [a,b]$ a measurable set on which $f$ is continuous such that $\lambda_F(D)=0$. Consider an arbitrary sequence of Riemann-Stieljes  sums  

$$
S_{n}=\sum_{i=0}^{\ell (n)-1}f(c_{i,n})(F(x_{i+1,n})-F(x_{i,n})),
$$

\noindent for which the sequence of modulus tends to zero with $n$. Justify the double inequality

$$
\int h_{n}d\lambda_F \leq S_n \leq \int H_{n}d\lambda_F
$$

\noindent and next, by denoting $H=[a,b]\setminus D$

$$
\int h_{n} 1_{H} d\lambda_F \leq S_n \leq \int H_{n} h_{n} 1_{H} d\lambda_F.
$$

\noindent By using Question (7) of  \textit{Exercise 8 in Doc 03-06} \ in Chapter \ref{03_setsmes_applimes_cas_speciaux}, what are the limits of $h_{n} 1_{H}$ and 
$H_{n} 1_{H}$.\\

\noindent Use the Dominated Convergence Theorem to get

$$
\int f_{\ast} 1_{H} d\lambda_F \leq  \liminf_{n\rightarrow +\infty} S_n \leq \limsup_{n\rightarrow +\infty} S_n \leq\int f*{\ast} 1_{H} d\lambda_F.
$$

\bigskip \noindent Conclude.\\

\bigskip \noindent \textbf{Exercise 3} \label{exercise01_doc06-07}\\

\noindent Let $F: [a,b]\longrightarrow \mathbb{R}$ be a non-constant and non-decreasing function.\\

\noindent On the whole extended real line $\mathbb{R}$, the Lebesgue-Stieljes integral of a measurable application \textit{a.e.} $f$ is defined on  if and only if

$$
\int_{\mathbb{R}} f^{-} d\lambda_F<+\infty \text{ and }  \int_{\mathbb{R}} f^{+} d\lambda_F<+\infty
$$

\noindent and the integral of $f$ 

$$
\int_{\mathbb{R}} f d\lambda_F=\int_{\mathbb{R}} f^{+} d\lambda_F - \int_{\mathbb{R}} f^{-} d\lambda_F<+\infty
$$

\noindent is finite if and only if the integral of $|f|$ 

$$
\int_{\mathbb{R}} |f| d\lambda_F=\int_{\mathbb{R}} f^{+} d\lambda + \int_{\mathbb{R}} f^{-} d\lambda_F<+\infty
$$

\noindent is.\\

\bigskip \noindent The Riemann-Stieljes integral is not defined on the whole real line. Instead, we have improper integrals obtained by limits.\\

\noindent The improper integral of $f$ on $\mathbb{R}$ 

\begin{equation}
\int_{-\infty }^{+\infty }f(x)\text{ }dF(x)  \label{RS-LB4}
\end{equation}

\noindent is defined as the limit of the integrals of 
 
\begin{equation*}
\int_{a}^{b}f(x)\text{ }dF(x),
\end{equation*}

\noindent as $a\rightarrow -\infty $ and $b\rightarrow +\infty$, independently of $a\rightarrow -\infty $ and $b\rightarrow +\infty$.

\bigskip \noindent Question (a) Show that if $f$ is $\lambda_F$-integrable on $\mathbb{R}$, locally bounded and $\lambda_F$-\textit{a.e.} continuous, then the improper integral of $f$ exists, is finite and is equal to the Lebesgue-Stieljes integral, that is 

$$
\int_{-\infty}^{+\infty} f \ dF(x) = \int f \ d\lambda_F. \ \ (B1)
$$

\bigskip \noindent Question (b) Show that if the improper integral of $|f|$ exists and is finite, then $f$ is Lebesgue-Stieljes integrable and (B1) holds.\\

\bigskip \noindent Question (c). Let $f=1_{\mathbb{Q}}$. Show that the improper integral of $f$ is not defined and that the Lebesgue-Stieljes integral of $f$ exists and give its value.\\

\bigskip \noindent Question (d) Consider here the particular case $F(x)=x$ and $f(x)=(\sin x)/x$, $x\in \mathbb{R}$. Show that the Lebesgue integral of $f$ does not exists, the Lebesgue integral of $|f|$ is infinite, the improper integral of $f$ exists and is infinite and the improper integral of $f$ exists and is finite and cannot be a Lebesgue integral.\\

\noindent Follow this way : \\

\noindent (i) Remark that $f$ is locally continuous and the locally bounded. Use Exercise 2 to get for all $a\leq 0$ and $b\geq 0$, 

$$
\int_{a}^{b} |f(x)| dx)=\int_{[a,b]} |f| \ d\lambda. \ \ (B18)
$$

\bigskip \noindent Let $a\downarrow +\infty$ and $b \uparrow +\infty$ and conclude that 

$$
\int_{-\infty}^{+\infty} |f(x)| dx=\int_{\mathbb{R}} |f(x)| \ d\lambda=+\infty, \ \ (B19)
$$

\bigskip \noindent (ii) Define 

$$
A=\{k\pi, k\in \mathbb{Z}\}, \ B=\sum_{k\ in \mathbb{Z}} ]2k\pi, (2k+1)\pi[ \ and \ B=\sum_{k\ in \mathbb{Z}} ](2k+1)\pi, 2(k+1)\pi[.
$$

\bigskip \noindent Establish that

$$
f=1 \ on \ \{0\}, f=0 \ on \ A\setminus \{0\} , f>0 \ on \ B, \ f<0 \ on \ C.
$$

\bigskip \noindent and

$$
f^+=f 1_{A+C}, \  f^-=f 1_{A+C} \ and \ \lambda(A)=0.
$$

\bigskip \noindent Conclude that $f^+$ and $f^-$ are continuous outside the null set $A$ and are  are locally \textbf{a.e.} continuous and they are also locally bounded. Repeat the argiments used in  Formulas (B18) and (B19) to 

$$
\int_{-\infty}^{+\infty} f^-(x) dx=\int_{\mathbb{R}} f^-(x) \ d\lambda \ \ and \ \ \int_{-\infty}^{+\infty} f^+(x) dx=\int_{\mathbb{R}} f^+(x) \ d\lambda. \ \ (B20)
$$

\bigskip \noindent (ii) Justify

$$
\int_{-\infty}^{+\infty} f^-(x) dx= \sum_{k\ in \mathbb{Z}} \int_{(2k+1)\pi}^{2(k+1)\pi} \frac{-\sin x}{x} \ dx
$$

\bigskip \noindent and

$$
\int_{-\infty}^{+\infty} f^+(x) dx= \sum_{k\ in \mathbb{Z}} \int_{2k\pi}^{(2k+1)\pi} \frac{\sin x}{x} \ dx
$$

\bigskip \noindent and compare the series in the right-hand members and prove the equality of the right-hand members to get

$$
\int_{-\infty}^{+\infty} |f|(x) dx=2 \int_{-\infty}^{+\infty} f^-(x) dx=\int_{-\infty}^{+\infty} f^+(x) dx=+\infty.
$$

\bigskip \noindent Now conclude.\\

\noindent \LARGE \textbf{DOC 06-08 : Convergence Types - Exercises With Solutions} \label{doc06-08}\\
\bigskip
\Large

\bigskip \noindent \textbf{NB}. Any unspecified limit below is meant as $n \rightarrow +\infty$.\\ 

\noindent \textbf{Exercise 01}. \label{exercises01_sol_doc06-08} (Useful Formulas in $\mathbb{R}$).\\

\noindent (a) Quickly check the following equivalences for an real number $x$ or a sequence or elements of $\overline{\mathbb{R}}$.\\

$$
(\forall \varepsilon>0, |x|\leq \varepsilon) \Leftrightarrow (\forall r \in \mathbb{Q} \ and \ r>0, |x|\leq r), \ \ (EQUIV01)
$$

$$
(\forall \varepsilon>0, |x|\leq \varepsilon) \Leftrightarrow  (\forall k \in \mathbb{N} \ and \ k\geq 1, |x|\leq 1/k), \ \ (EQUIV02)
$$

\bigskip \noindent and

$$
(\forall A>0, |x|\geq A) \Leftrightarrow  (\forall k \in \mathbb{N} \ and \ k\geq 1, |x| \geq k). \ \ (EQUIV03)
$$

\bigskip \noindent Formula \textit{(EQUIV04a)} : 
$$
(\forall A>0, \ \exists n\geq 0, \ \forall p\geq n : \ |x_p|\geq A)
$$
$$
\Leftrightarrow
$$
$$
(\forall k \in \mathbb{N} \ and \ k\geq 1, \ \exists n\geq 0, \ \forall p\geq n : \ |x_p|\geq k).
$$

\bigskip \noindent Formula \textit{(EQUIV04b)}
$$
(\forall A>0, \ \exists n\geq 0, \ \forall p\geq n : \ |x_p|\leq - A)
$$
$$\Leftrightarrow$$
$$
 (\forall k \in \mathbb{N} \ and \ k\geq 1, \ \exists n\geq 0, \ \forall p\geq n : \ |x_p|\geq -k)).
$$

\bigskip \noindent Formula \textit{(EQUIV05)}
$$
(\forall \varepsilon>0, \ \exists n\geq 0, \ \forall p\geq n : \ |x_p|\leq \varepsilon)
$$
$$
\Leftrightarrow
$$
$$
(\forall k \in \mathbb{N} \ and \ k\geq 1, \ \exists n\geq 0, \ \forall p\geq n : \ |x_p|\leq 1/k).
$$

\bigskip \noindent etc.\\

\noindent (b) Let $a$ and $b$ be two real numbers and $\varepsilon>0$. Let $(a_j)_{1\leq j \leq k}$ be a family of $k>0$ real numbers. Let also $(\varepsilon_j)_{1\leq j \leq k}$ be a family of 
positive real numbers such that $\varepsilon_1+\varepsilon_2+ ... + \varepsilon_k = \varepsilon$. By using the negation of the triangle inequality, justify the relations 

$$
(|a+b|>\varepsilon) \Rightarrow (|a|>\varepsilon/2) \ or \ (|b|>\varepsilon/2). \ \ (ROS01)
$$

\bigskip \noindent and

$$
\biggr(\left| \sum_{1\leq j \leq k} a_j \right|>\varepsilon\biggr) \Rightarrow (\exists j \in \{1,...,k\}, \ |a_j|>\varepsilon_j). (ROS02)
$$

\bigskip \noindent (c) Let $f$ and $g$ be two functions defined on $\Omega$ with values in $\overline{\mathbb{R}}$. Let $(f_j)_{1\leq j \leq k}$ be a family $k>0$ of functions defined on $\Omega$ with values in $\overline{\mathbb{R}}$. real numbers. Let also $(\varepsilon_j)_{1\leq j \leq k}$ be a family of positive real numbers such that 
$\varepsilon_1+\varepsilon_2+ ... + \varepsilon_k = \varepsilon$. By using the previous question, justify the relations

$$
(|f+g|>\varepsilon) \subset (|f|>\varepsilon/2) \cup (|g|>\varepsilon/2). \ \ (OS01)
$$

\bigskip \noindent and

$$
\biggr(\left| \sum_{1\leq j \leq k} f_j \right|>\varepsilon\biggr) \subset \bigcup_{{1\leq j \leq k}} (|f_j|>\varepsilon_j). \ \ (OS02)
$$

\bigskip \noindent (d) Let $f$ and $g$ be two measurable functions defined on $\Omega$ with values in $\overline{\mathbb{R}}$. Let $(f_j)_{1\leq j \leq k}$ be a family $k>0$ of easurable functions defined on $\Omega$ with values in $\overline{\mathbb{R}}$. real numbers. Let also $(\varepsilon_j)_{1\leq j \leq k}$ be a family of positive real numbers such that 
$\varepsilon_1+\varepsilon_2+ ... + \varepsilon_k = \varepsilon$. By using the previous question, justify the relations

$$
m(|f+g|>\varepsilon) \leq m(|f|>\varepsilon/2) + m(|g|>\varepsilon/2). \ \ (MOS01)
$$

\bigskip \noindent and

$$
m\biggr(\left| \sum_{1\leq j \leq k} f_j \right|>\varepsilon\biggr) \leq \sum_{{1\leq j \leq k}} (|f_j|>\varepsilon_j). \ \ (MOS02)
$$
\bigskip \noindent \textbf{NB}. This exercise gives technical and simple facts very know of $\mathbb{R}$ to be used in the sequel.\\

\bigskip \noindent \textbf{SOLUTIONS}.\\

\noindent \textbf{Question (a)}.\\

\noindent \textbf{EQUIV01}. The direct implication $\rightarrow$ is obvious. To prove the indirect implication $\leftarrow$, suppose that we have

$$
(\forall r \in \mathbb{Q} \ and \ r>0, |x|\leq r).
$$  

\bigskip \noindent Now, let $\varepsilon>0$. Since $\mathbb{Q}$ is dense in $\mathbb{R}$, the interval $]0,\varepsilon[$ intersects with $\mathbb{Q}$, that is there exists $r \in \mathbb{Q} \cap ]0,\varepsilon[$. Then by assumption,

$$
|x|\leq r < \varepsilon.
$$

\bigskip \noindent Hence, we have

$$
(\forall \varepsilon>0, |x|\leq \varepsilon). \ \square
$$

\bigskip \noindent \textbf{EQUIV02}. The direct implication $\rightarrow$ is obvious. To prove the indirect implication $\leftarrow$, suppose that we have

$$
(\forall k \in \mathbb{N} \ and \ k>1, |x|\leq 1/k).
$$  

\bigskip Now, let $\varepsilon>0$. Since $\mathbb{R}$ is Archimedes, by comparing $\varepsilon>0$ and $1$, there exists an integer $k \varepsilon>1$, that is $1/k < \varepsilon>0$. Then by assumption,

$$
|x|\leq 1/k < \varepsilon.
$$

\bigskip \noindent Hence, we have

$$
(\forall \varepsilon>0, |x|\leq \varepsilon). \ \square
$$

\bigskip \noindent \textbf{EQUIV03}. The direct implication $\rightarrow$ is obvious. To prove the indirect implication $\leftarrow$, suppose that we have

$$
(\forall k \in \mathbb{N} \ and \ k>1, |x|\geq k).
$$  

\bigskip Now, let $A>0$. Since $\mathbb{R}$ is Archimedes, by comparing $a=1$ and $A>0$, there exists an integer $k>0$ such that $k \times a > A$, that is $k > A$. Then by assumption,

$$
|x|\geq k > A.
$$

\bigskip \noindent Hence, we have

$$
(\forall A>0, |x|\leq A). \ \square
$$

\bigskip \noindent \textbf{EQUIV04}. The direct implication $\rightarrow$ is obvious. To prove the indirect implication $\leftarrow$, suppose that we have

$$
(\forall k \in \mathbb{N} \ and \ k\geq 1, \ \exists n_0 \geq 0, \ \forall p\geq n_0 : \ |x_p|\geq k).
$$

$$
(\forall A>0, \ \exists n\geq 0, \ \forall p\geq n : \ |x_p|\geq A)
$$
$$
\Leftrightarrow
$$
$$
(\forall k \in \mathbb{N} \ and \ k\geq 1, \ \exists n\geq 0, \ \forall p\geq n : \ |x_p|\geq k). \ \ (EQUIV04)
$$

\bigskip \noindent Now, let $\varepsilon>0$. By the argument used the solution of \textit{(EQUIV03)}, consider an integer $k$ such that $k > A$. Then by assumption, we have

$$
\forall p\geq n_0 : \ |x_p|\geq k > A.
$$

\bigskip \noindent Hence, we have

$$
(\forall A>0, \ \exists n_0\geq 0, \ \forall p\geq n_0 : \ |x_p|\geq A). \ \square
$$

\bigskip \noindent \textbf{EQUIV05}. The direct implication $\rightarrow$ is obvious. To prove the indirect implication $\leftarrow$, suppose that we have

$$
\forall k \in \mathbb{N} \ and \ k\geq 1, \ \exists n_0\geq 0, \ \forall p\geq n_0 : \ |x_p|\leq 1/k.
$$

\bigskip \noindent Now, let $\varepsilon>0$. By the argument used the solution of \textit{(EQUIV02)}, consider an integer $k$ such that $1/k \varepsilon$. Then by assumption, we have

$$
\forall p\geq n_0 : \ |x_p|\leq 1/k \leq \varepsilon.
$$

\bigskip \noindent We get

$$
(\forall \varepsilon>0, \ \exists n\geq 0, \ \forall p\geq n : \ |x_p|\leq \varepsilon). \ \square
$$

\bigskip \bigskip \noindent \textbf{Question (b)}.\\

\noindent \textbf{Formulas (ROS01) and (RS02)}. By the triangle inequality, we have

$$
(|a|\leq \varepsilon/2) \ and \ (|b|\leq \varepsilon/2) \Rightarrow (|a+b|\leq |a|+ |b| \leq \varepsilon/2 +\varepsilon/2=\varepsilon)
$$

\bigskip \noindent which yields

$$
(|a|\leq \varepsilon/2) \ and \ (|b|\leq \varepsilon/2) \Rightarrow (|a+b|\leq \varepsilon. \ \ (CROS01)
$$

\bigskip \noindent The repetition of the triangle inequality $k$ times also gives

\begin{eqnarray*}
&&(\forall j \in \{1,...,k\}, \ |a_j|\leq \varepsilon_j)\\
&\Rightarrow& \left| \sum_{1\leq j \leq k} a_j \right| \leq  \sum_{1\leq j \leq k} |a_j|\\
&\leq&  \sum_{1\leq j \leq k} \varepsilon_j=\varepsilon,
\end{eqnarray*}

\bigskip \noindent which yields

$$
(\forall j \in \{1,...,k\}, \ |a_j|\leq \varepsilon_j) \Rightarrow \left(\left| \sum_{1\leq j \leq k} a_j \right|  \leq \varepsilon\right). \ \ (CROS02).
$$

\bigskip \noindent The formulas (ROS01) and (RSO02) to be proved are the negations of (CROS01) and (CRSO02) respectively.\\

\bigskip \bigskip \noindent \textbf{Question (c)}. The formulas (OS01) and (SO02) to be proved are the sets versions (ROS1) and (RSO2). For example, we say : \\

\begin{eqnarray*}
\omega \in (|f+g|>\varepsilon) &\Rightarrow& |f+g|(\omega)>\varepsilon \ (By \ definition)\\
&\Leftrightarrow& |f(\omega)+g(\omega)|>\varepsilon \\  
&\Rightarrow& |f(\omega)|>\varepsilon/2 \ or \ |g(\omega)|>\varepsilon/2 \ (By \ (ROS1)) \\  
&\Leftrightarrow& \omega \in (|f|>\varepsilon/2) \ or \ \omega \in (|g|>\varepsilon/2)\\
&\Leftrightarrow& \omega \in (|f|>\varepsilon/2) \cap (|g|>\varepsilon/2),
\end{eqnarray*}

\noindent which is equivalent to

$$
(|f+g|>\varepsilon) \subset (|f|>\varepsilon/2) \cap (|g|>\varepsilon/2). \ \square
$$

\bigskip \noindent $\blacksquare$\\

\bigskip \noindent \textbf{Exercise 02}. \label{exercises02_sol_doc06-08}\\

\noindent Let  $f_{n})_{n\geq 1}$\ be a sequence of measurable functions from($\Omega ,\mathcal{A},m)$\ to $\overline{\mathbb{R}}$ and $f : (\Omega ,\mathcal{A},m)$\ $\longmapsto \overline{\mathbb{R}}$ be a measurable functions.\\

\noindent Question (1). Let us remind the definitions of the following limits, as $\rightarrow +\infty$ :\\

\noindent If $f(\omega)=+\infty$,
$$
(f_n(\omega) \rightarrow f(\omega)=+\infty) \Leftrightarrow (\forall \eta>0, \ exists n\geq 0 : \ \forall p\geq n, \ f_n(\omega)\geq \eta; \ \ (DEF1)
$$

\bigskip \noindent If $f(\omega)=-\infty$,
$$
(f_n(\omega) \rightarrow f(\omega)=-\infty) \Leftrightarrow (\forall \eta>0, \ exists n\geq 0 : \ \forall p\geq n, \ f_n(\omega)\leq -\eta;  \ \ (DEF2)
$$

\bigskip \noindent If $f(\omega) \in \mathbb{R}$,

$$
(f_n(\omega) \rightarrow f(\omega) \in \mathbb{R}) \Leftrightarrow (\forall \varepsilon>0, \ exists n\geq 0 : \ \forall p\geq n, \ |f_n(\omega)-f_n(\omega)|\leq \varepsilon.  \ \ (DEF3)
$$

\bigskip \noindent Set 

\begin{equation*}
A1=(f=+\infty ) \bigcap \left( \bigcap_{\eta>0} \bigcup_{n\geq 1} \bigcap_{p\geq n} (f_{p} \geq \eta) \right)
\end{equation*}

\bigskip
 
\begin{equation*}
B1=(f=-\infty ) \bigcap \left( \bigcap_{\eta >0} \bigcup_{n\geq 1} \bigcap_{p\geq n} (f_{p} \leq -\eta) \right)
\end{equation*}

\bigskip \noindent and \\

\begin{equation*}
C1=(|f|<+\infty ) \bigcap \left( \bigcap_{\varepsilon >0} \bigcup_{n\geq 1} \bigcap_{p\geq n} (|f_{p}-f| \leq \varepsilon) \right).
\end{equation*}

\noindent By using decomposition

$$
(f_n \rightarrow f)=(f_n \rightarrow f)\cap (f=+\infty) + (f_n \rightarrow f)\cap (f=-\infty) + (f_n \rightarrow f)\cap (|f|<+\infty),
$$

\noindent show that

$$
(f_n \rightarrow f)=A1 + B1 + C1.
$$

\bigskip \noindent Since the values of $\eta$ and $\varepsilon$ are not countable, can you conclude that $(f_n \rightarrow f)$ is measurable?\\

\noindent Now, by using Exercise 1, show that 

$$
(f_n \rightarrow f)=A + B + C,
$$

\noindent with :

\begin{equation*}
A=(f=+\infty ) \bigcap \left( \bigcap_{k\leq 1} \bigcup_{n\geq 1} \bigcap_{p\geq n} (f_{p} \geq k) \right)
\end{equation*}

\bigskip

\begin{equation*}
B=(f=-\infty ) \bigcap \left( \bigcap_{k\leq 1} \bigcup_{n\geq 1} \bigcap_{p\geq n} (f_{p} \leq -k) \right)
\end{equation*}

\bigskip \noindent and \\

\begin{equation*}
C=(|f|<+\infty ) \bigcap \left( \bigcap_{k \leq 1} \bigcup_{n\geq 1} \bigcap_{p\geq n} (|f_{p}-f| \leq 1/k) \right).
\end{equation*}

\noindent and conclude that $(f_n \rightarrow f)$ is measurable.\\

\bigskip \noindent Question (2). Show that $f_n \rightarrow f$ \textit{a.e.} if and only if the three equations hold :

\begin{equation*}
m\left((f=+\infty ) \bigcap \left( \bigcup_{k\leq 1} \bigcap_{n\geq 1} \bigcup_{p\geq n} (f_{p} < k) \right)\right)=0;
\end{equation*}

\begin{equation*}
m\left((f=-\infty ) \bigcap \left( \bigcup_{k\leq 1} \bigcap_{n\geq 1} \bigcup_{p\geq n} (f_{p} > -k) \right)\right) =0
\end{equation*}

\bigskip \noindent and \\

\begin{equation*}
m\left( (|f|<+\infty ) \bigcap \left( \bigcup_{k \leq 1} \bigcap_{n\geq 1} \bigcup_{p\geq n} (|f_{p}-f| > 1/k) \right) \right)=0.
\end{equation*}

\bigskip \noindent Question (3). If $f$ is \textit{a.e}, show that $f_n \rightarrow f$ \textit{a.e.} if and only if

\begin{equation*}
m\left( (|f|<+\infty ) \bigcap \left( \bigcup_{k \leq 1} \bigcap_{n\geq 1} \bigcup_{p\geq n} (|f_{p}-f| > 1/k) \right) \right)=0.
\end{equation*}

\noindent \bigskip \noindent Question (3). Let $f$ still be finite \textit{a.e.} Denote

$$
C_{\infty} = \bigcup_{k \leq 1} \bigcap_{n\geq 1} \bigcup_{p\geq n} (|f_{p}-f| > 1/k) =0.
$$

\noindent and for $k\leq 1$,

$$
C_k=\bigcap_{n\geq 1} \bigcup_{p\geq n} (|f_{p}-f| > 1/k).
$$

\bigskip Show that $f_n \rightarrow f$ if and only if, for any $k\geq 1$,

$$
m(C_k)=0.
$$
 
\noindent \bigskip \noindent Question (4). Let $f$ still be finite \textit{a.e.}. Combine the previous result to get that  the results of Exercise 1 to get : $f_n \rightarrow f$ if and only if, for any 
$\varepsilon>0$, $m(D_{\varepsilon})=0$, where

$$
D_\varepsilon=\bigcap_{n\geq 1} \bigcup_{p\geq n} (|f_{p}-f| > \varepsilon).
$$
 
\bigskip \noindent \bigskip \noindent Question (5). Let  $f_{n})_{n\geq 1}$ and $f$ \textit{a.e.} finite and let $m$ be a finite measure. Set

$$
D_{n,\varepsilon}=\bigcup_{p\geq n} (|f_{p}-f| > \varepsilon).
$$ 

\bigskip \noindent Remark that the sequence $(D_{n,\varepsilon})_{n\geq 0}$ is non-increasing and that $(|f_n-f|>\varepsilon) \subset D_{n,\varepsilon}$. Use the continuity of the measure $m$ 
(Point 04.10 of Doc 04-01, page \pageref{doc04-10} in Chapter \ref{04_measures}) to show that

$$
(f_n \rightarrow f) \Rightarrow (f_n \rightarrow_m f).
$$

\bigskip\bigskip \noindent \textbf{SOLUTIONS}.\\

\noindent Question (1). Let us use the decomposition

$$
(f_n \rightarrow f)=(f_n \rightarrow f)\cap (f=+\infty) + (f_n \rightarrow f)\cap (f=-\infty) + (f_n \rightarrow f)\cap (|f|<+\infty),
$$

\noindent Now let $\omega \in (f_n \rightarrow f)\cap (f=+\infty)$. From the reminder, since $(f(\omega)=+\infty)$, $f_n(\omega) \rightarrow f(\omega)$ if and only if

$$
\forall \eta>0, \ exists n\geq 0 : \ \forall p\geq n, \ f_n(\omega)\geq \eta,
$$

\bigskip \noindent that is, by using the sets language,

$$
\omega \in \bigcap_{\eta>0} \bigcup_{n\geq 0} \bigcap_{p\geq n} (f_n\geq \eta).
$$

\bigskip It is then clear that

$$
(f_n \rightarrow f)\cap (f=+\infty) =  \left( \bigcap_{\eta>0} \bigcup_{n\geq 0} \bigcap_{p\geq n} (f_n\geq \eta) \right) \bigcap (f=+\infty)=A.
$$

\noindent We use the same method to show that $(f_n \rightarrow f)\cap (f=-\infty)=B1$ and $(f_n \rightarrow f)\cap (|f|<+\infty)=C1$ and reach the conclusion.\\

\noindent \textbf{Unfortunately}, we are not able to say if the sets $A1$, $B1$ are $C1$ are measurable or not since the intersection in $eta$ or in $\varepsilon$ is not a countable one.\\

\noindent By Exercise 1, we may use the Archimedes property of $\mathbb{R}$ to see that we may replace $\eta$ by $k \in \mathbb{N} \setminus \{0\}$ in (DEF1) and (DEF1) above (as in (EQUIV04a) and 
(EQUIV04b) in Exercise 1), and $\varepsilon$ by $1/k$, $k \in \mathbb{N} \setminus \{0\}$ in (DEF1) (as in (EQUIV05) in Exercise 1). The same reasoning gives

$$
(f_n \rightarrow f)= A + B + C.
$$

\bigskip \noindent But $A$, $B$ and $C$ are countable intersections of countable unions of countable intersections of measurable  sets $(f_n\geq k)$, $(f_n\leq -k)$ or $(|f_n-f| \leq  k)$. Hence
$(f_n \rightarrow f)$ is measurable.
 
\noindent Question (2). We still use the same technique and write

$$
(f_n \nrightarrow f)=(f_n \nrightarrow f)\cap (f=+\infty) + (f_n \nrightarrow f)\cap (f=-\infty) + (f_n \nrightarrow f)\cap (|f|<+\infty),
$$

\noindent Next by taking the negations of the definitions of convergence that used the integers $k$, we get

$$
(f_n \nrightarrow f)\cap (f=+\infty) =  \left( \bigcup_{k\geq 1} \bigcap_{n\geq 0} \bigcup_{p\geq n} (f_n < k) \right) \bigcap (f=+\infty)=A0,
$$

$$
(f_n \nrightarrow f)\cap (f=-\infty) =  \left( \bigcup_{k\geq 1} \bigcap_{n\geq 0} \bigcup_{p\geq n} (f_n > -k) \right) \bigcap (f=+\infty)=B0.
$$

\noindent and

$$
(f_n \nrightarrow f)\cap (f \ finite) =  \left( \bigcup	_{k\geq 1} \bigcap_{n\geq 0} \bigcup_{p\geq n} (|f_n-f| > 1/k) \right) \bigcap (f \ finite)=C0.
$$

\noindent Since 

$$
(f_n \nrightarrow f)=A0+B0+C0,
$$

\bigskip \noindent $f_n \rightarrow f$ if and only  $m(A0+B0+C0)=0$ if and only ($m(A0)=0$, $m(B0)=0$ and $m(C0)=0$). $\square$\\

\noindent Question (3). If $f$ is \textit{a.e.} finite, it is clear that $0\leq m(A0) \leq (f=+\infty)=0$ and $0\leq m(B0) \leq (f=1\infty)=0$. We get that : $f_n \rightarrow f$ iff $m(C0)=0$. But since 
the complement of $(f \ finite )$ is a null-set, we get that $m(C0)=0$ if and only

$$
m\left( \bigcup	_{k\geq 1} \bigcap_{n\geq 0} \bigcup_{p\geq n} (|f_n-f| > 1/k) \right)=0.
$$

\bigskip By denoting for $k\geq 1$

$$
C_k=\bigcap_{n\geq 0} \bigcup_{p\geq n} (|f_n-f| > 1/k),
$$

\noindent we have, if $f$ is \textit{a.e.} finite : $f_n \rightarrow f$ if and only if 

$$
m\left( \bigcup	_{k\geq 1} C_k \right)=0,
$$

\bigskip \noindent which is equivalent to $m(C_k)=0$ for all $k\geq 1$.\\

\bigskip \noindent Question (4). By using the Archimedes property of $\mathbb{R}$ as we did in Exercise 1 (in the proof of (EQUIV02)) and by the non-decreasingness of the measure $m$, we quickly get the following equivalence

$$
(\forall \varepsilon, \ m(D_{\varepsilon})=0) \Leftrightarrow (\forall \varepsilon, \ m(C_{k})=0). \square
$$

\bigskip \noindent Question (5). Let the $f_n$ and $f$ be \textit{a.e.} finite. Let $\varepsilon>0$. Suppose that $f_n \rightarrow f$. By Question (4), we have $m(D_{\varepsilon})=0$. Next,
the sequence $(D_{n,\varepsilon})_{n\geq 0}$ is non-decreasing since as $n$ grows, the unions looses members and thus become smaller. The set limit is the intersection of all the 
$D_{n,\varepsilon}$'s, that is $D_{\varepsilon}$. Since the measure is finite, the continuity of $m$ implies

$$
m(D_{n,\varepsilon}) \downarrow m(D_{\varepsilon})=0. \ \ (AP1)
$$ 

\bigskip \noindent It is also clear that for all $n\geq 0$,

$$
(|f_n-f| > \varepsilon) \subset \bigcup_{p\geq n} (|f_p-f| > \varepsilon) 
$$

\bigskip \noindent which implies

$$
m(|f_n-f| > \varepsilon) \leq m(D_{n,\varepsilon}). \ \ (AP2)
$$

\bigskip \noindent By combining taking into account Formula (AP1) and by letting $n\rightarrow +\infty$ in Formula (AP2), we get

$$
m(|f_n-f| > \varepsilon) \rightarrow 0 \ as \ n\rightarrow +\infty,
$$

\bigskip \noindent and this is valid for any $\varepsilon>0$. Hence $f_n \rightarrow_m f$. $\square$\\

\bigskip \noindent \textbf{Exercise 3}. \label{exercises03_sol_doc06-08}

\noindent \textbf{General Hint}. In each question, use the Measure Space Reduction property (see pages \pageref{msrp} and \pageref{msrpA}) and combine to others assumptions to place yourself in a space
$\Omega_0$, whose complement is a null-set, and on which all the \textit{a.e.} properties hold for each $\omega \in \Omega_0$.\\

\noindent Let $(f_{n})_{n\geq 1}$ and $(g_{n})_{n\geq 1}$ be sequences of real-valued measurable functions such that $f_{n}\rightarrow f$ $a.e.$ and $g_{n}\rightarrow g$ $a.e.$. 
Let $a$ and $b$ be finite real numbers. Let $H(x,y)$ a continuous function of $(x,y) \in D$, where $D$ is an open set of $\mathbb{R}^2$. Show that the following assertion hold true : \\

\noindent \textbf{(1)} The \textit{a.e.} limit of $(f_{n})_{n\geq 1}$, if it exists, is unique \textbf{a.e.}.\\

\noindent \textbf{(2)} If $af+bg$ is \textit{a.e.} defined, then $af_{n}+bg_{n}\rightarrow af+bg$ $a.e$.\\

\noindent \textbf{(3)} If $fg$ is \textit{a.e.} defined, then $f_{n}g_{n}\rightarrow fg$ $a.e$\\

\noindent \textbf{(4)} If $f/g$ is \textit{a.e.} defined, (that is $g$ is $a.e$ nonzero and $f/g \neq \frac{\pm \infty}{\pm \infty}$ \textit{a.e.}), then

$$
f_{n}/g_{n}\rightarrow f/g, \ a.e.
$$

\noindent \textbf{(5)} If $(f_{n}, g_n)_{n\geq 1} \subset D$ \textit{a.e.} and $(f,g) \in D$ \textit{a.e.}, then

$$
H(f_{n},g_{n})\rightarrow H(f,g), \ a.e.\\
$$

\bigskip \bigskip \noindent \textbf{SOLUTIONS}.\\

\noindent (1) Suppose that $f_n$ converges \textit{a.e.} to both $f_1$ and $f_2$. Denote

$$
\Omega_0^c=(f_n \nrightarrow f_{(1)}) \cap (f_n \nrightarrow f_{(2)}).
$$ 

\bigskip \noindent It is clear that $\Omega_0^c$ is a null-set. Besides, for all $\omega \in \Omega_0$,

$$
f_n(\omega) \rightarrow f_{(1)}(\omega) \ \ and \ \ f_n(\omega) \rightarrow f_{(2)}(\omega).
$$

\bigskip \noindent We get, by the uniqueness of the limit in $\overline{\mathbb{R}}$, we get that $f_{(1)}=f_{(2)}$ on $\Omega_0$, that is outside a null-set. Hence $f_{(1)}=f_{(2)}$ \textit{a.e.}. $\square$\\

\bigskip  \noindent (2) By combining the hypotheses, we can find a measurable space sub-space $\Omega_0$ of $\Omega$ whose complement is a $m$-null and, on which, we have $f_{n}\rightarrow f$ $a.e.$ and $g_{n}\rightarrow g$ and $af+bg$ is well defined.  Then for any $\omega \in \Omega_0$, $af_n(\omega)+bg_n(\omega)$ is well defined for large values of $n$ and hence

$$
(af_n+bg_n)(\omega) \rightarrow (af+bg)(\omega).
$$

\bigskip \noindent Hence $af_n+bg_n$ converges to $af+bg$ outside the null-set $\Omega_0^c$. $\square$\\

\bigskip \noindent (3) By combining the hypotheses, we can find a measurable space sub-space $\Omega_0$ of $\Omega$ whose complement is a $m$-null and, on which, we have $f_{n}\rightarrow f$ $a.e.$ and $g_{n}\rightarrow g$ and $fg$ is well defined.  Then for any $\omega \in \Omega_0$, $f_n(\omega) g_n(\omega)$ is well defined for large values of $n$ and hence

$$
(f_n g_n)(\omega) \rightarrow (fg)(\omega).
$$

\bigskip \noindent Hence $f_n g_n$ converges to $fg$ outside the null-set $\Omega_0^c$. $\square$\\

\bigskip \noindent (4) Repeat the same argument.\\

\noindent (5)  By combining the hypotheses, we can find a measurable space sub-space $\Omega_0$ of $\Omega$ whose complement is a $m$-null and, on which, we have $f_{n}\rightarrow f$ $a.e.$ and $g_{n}\rightarrow g$ and all the couples $(f_n, g_n)$'s are in $D$ and $(f,g)$ is in D. Then for any $\omega \in \Omega_0$, $(f_n(\omega), g_n(\omega))$ converges $(f(\omega), g(\omega))$ in $D$ By continuity of $H$, we get

$$
(\forall \omega \in \Omega_), \ H(f_n(\omega), g_n(\omega)) \rightarrow H(f(\omega), g(\omega)).
$$

\bigskip \noindent Hence, $H(f_n, g_n)$ converges to $H(f, g)$ outside the null-set $\Omega_0$. $\blacksquare$\\

\bigskip \noindent \textbf{Exercise 4}. \label{exercises04_sol_doc06-08}\\

\noindent Let $(f_{n})_{n\geq 1}$ be a sequence of \textit{a.e.} finite measurable real-valued functions.\\

\noindent Question (1) By using the definition of a Cauchy sequence on $\mathbb{R}$ and by using the techniques in Exercise 1, show that :   $(f_{n})_{n\geq 1}$ is an $m$-\textit{a.e.} Cauchy sequence if and only if

$$
m\left( \bigcup_{k\geq 1} \bigcap_{n\geq 1} \bigcup_{p\geq n} \bigcup_{q\geq n} (|f_p-f_q|>1/k)\right)=0. \ \ (CAU01)
$$ 

\noindent Question (2). Show that (CAU01) is equivalent to each of the two following assertions.\\

$$
m\left( \bigcup_{k\geq 1} \bigcap_{n\geq 1} \bigcup_{p\geq 0} \bigcup_{q\geq 0} (|f_{n+p}-f_{n+p}| > 1/k)\right)=0. \ \ (CAU02)
$$ 

$$
m\left( \bigcup_{k\geq 1} \bigcap_{n\geq 1} \bigcup_{p\geq 0}  (|f_{n+p}-f_{n}| > 1/k)\right)=0. \ \ (CAU03)
$$

\bigskip \noindent Question (3). Show that $(f_{n})_{n\geq 1}$ is an $m$-\textit{a.e.} Cauchy sequence if and only if $f_n$ converges \textbf{a.e.} to an \textbf{a.e.} measurable function $f$.\\

\bigskip \noindent \textbf{SOLUTIONS}.\\

\bigskip \noindent Question (1). By definition, $(f_{n})_{n\geq 1}$ is an $m$-\textit{a.e.} Cauchy sequence if and only if

$$
m( \{\omega \in \Omega, (f_n(\omega))_{n\geq 1} \ is \ not \ a \ Cauchy \ sequence \ in \ \mathbb{R}\} )=0.
$$

\bigskip \noindent Denote $C^c=\{\omega \in \Omega, (f_n(\omega))_{n\geq 1} \ is \ a \ Cauchy \ sequence \ in \ \mathbb{R}\}$. But  $(f_n(\omega))_{n\geq 1}$  is  a  Cauchy  sequence  in $\mathbb{R}$ if and only if

$$
f_p(\omega)-f_q(\omega) \rightarrow 0 \ as \ (p,q)\rightarrow (+\infty,+\infty),
$$

\bigskip \noindent which is equivalent to
 
$$
\forall \varepsilon>0, \ \exists n(\omega) : \ \forall p\geq n, \forall q \geq n), \ |f_n(\omega)-f_q(\omega)|<\varepsilon.
$$

\bigskip \noindent which, by the techniques of Exercise 1, is equivalent to
 
$$
\forall k\geq 1, \ \exists n(\omega) : \ \forall p\geq n, \forall q \geq n), \ |f_n(\omega)-f_q(\omega)|<1/k,
$$

\bigskip \noindent that is,

$$
\omega \in \bigcap_{k\geq 1} \bigcup_{n\geq 0} \bigcap_{p\geq n} \bigcap_{q\geq n} (|f_n(\omega)-f_q(\omega)| \leq 1/k).
$$

\bigskip \noindent Hence, we get

$$
C=\bigcup_{k\geq 1} \bigcap_{n\geq 0} \bigcup_{p\geq n} \bigcup_{q\geq n} (|f_n(\omega)-f_q(\omega)| > 1/k). \ \ ECAU01.
$$

\bigskip \noindent Hence, $(f_{n})_{n\geq 1}$ is an $m$-\textit{a.e.} Cauchy sequence if and only if

$$
m\left( C \right)=0
$$
 
\bigskip \noindent which is Formula (CAU01).\\

\bigskip \noindent Question (3). To obtain (CAU02), it is enough to put $r=p-n$ and $s=q-n$ in Formula (ECAU01) to get

$$
C=\bigcup_{k\geq 1} \bigcap_{n\geq 0} \bigcup_{r\geq 0} \bigcup_{s\geq 0} (|f_{n+r}(\omega)-f_{n+s}(\omega)| \leq 1/k). \ \ ECAU02.
$$

\bigskip \noindent which leads to Formula (CAU02) as equivalent to Formula (CAU01).\\

\noindent Now, it clear that Formula (CAU02) implies (CAU03). To prove the indirect implication, suppose that (CAU03) holds. Denote

$$
D^c=\bigcap_{k\geq 1} \bigcup_{n\geq 0} \bigcap_{p\geq 0}  (|f_{n+p}(\omega)-f_{n}(\omega)| \leq 1/k). \ \ ECAU03.
$$

\bigskip \noindent Let $\omega \in D$. Let $k\geq 1$ be fixed and let us consider a value $n=n(\omega)$ such that for all $r\geq 0$, we have $|f_{n+p}(\omega)-f_{n}(\omega)| \leq 1/k$. Now, for all $p`geq 0$ and $q\geq 0$,

$$
|f_{n+p}(\omega)-f_{n+q}(\omega)| \leq |f_{n+p}(\omega)-f_{n}(\omega)| + |f_{n+q}(\omega)-f_{n}(\omega)| < 2/k.
$$

\bigskip \noindent Hence

$$
D^c \subset \bigcap_{k\geq 1} \bigcup_{n\geq 0} \bigcap_{p\geq 0} \bigcap_{q\geq 0}  (|f_{n+p}(\omega)-f_{n+q}(\omega)| \leq 2/k),
$$

\bigskip \noindent that is

$$
\bigcup_{k\geq 1} \bigcap_{n\geq 0} \bigcup_{p\geq 0} \bigcup_{q\geq 0}  (|f_{n+p}(\omega)-f_{n+q}(\omega)| \leq 2/k) \subset D^c.
$$

\bigskip \noindent If Formula (CAU03) holds, we have $m(D^c)=0$, and thus for 

$$
E_1=\bigcup_{k\geq 1} \bigcap_{n\geq 0} \bigcup_{p\geq 0} \bigcup_{q\geq 0}  (|f_{n+p}(\omega)-f_{n+q}(\omega)| \leq 2/k),
$$

\bigskip \noindent we have

$$
m(E_1)=0.
$$

\bigskip \noindent The proof is complete if we establish that $E$ may be written also as

$$
E_2=\bigcup_{k\geq 1} \bigcap_{n\geq 0} \bigcup_{p\geq 0} \bigcup_{q\geq 0}  (|f_{n+p}(\omega)-f_{n+q}(\omega)| \leq 1/k).
$$

\bigskip \noindent But you may use the techniques in Exercise 1 to see that both sets $E_1$ and $E_2$ are exactly

$$
\bigcup_{\varepsilon >0} \bigcap_{n\geq 0} \bigcup_{p\geq 0} \bigcup_{q\geq 0}  (|f_{n+p}(\omega)-f_{n+q}(\omega)| \leq \varepsilon).
$$

\bigskip \noindent (Remark : by handling $E_2$, use the Archimedes property between $\varepsilon$ and $2$). $\square$\\

\bigskip \noindent Question (3).\\

\noindent (a) Suppose that $f_n \rightarrow f$ \textit{a.e.} and $f$ is \textit{a.e.} finite. Denote by $\Omega_0$ a measurable set whose complement is a null-set, on which all the $f_n$'s and $f$ are finite and $f_n(\omega) \rightarrow f(\omega)$ for any $\omega \in \Omega_0$. Then for any $\omega \in \Omega_0$, for any $k\geq 1$, there exists a value $n=n(\omega)$ such that for all $p\geq n$,

$$
|f_p(\omega)-f(\omega)| \leq 1/(2k).
$$  

\bigskip \noindent Hence, for all for all $p\geq n$ and $q\geq n$,

$$
|f_p(\omega)-f_q(\omega)| \leq |f_p(\omega)-f(\omega)|+ |f_q(\omega)-f(\omega)| \leq (1/(2k)+1/(2k))=1/k.
$$  

\bigskip \noindent we get, for all $\omega \in \Omega_0$,

$$
\forall k\geq 1, \ \exists n\geq 0, \ \forall p\geq n, \ \forall q\geq : |f_p(\omega)-f_q(\omega)|\leq 1/k.
$$

\bigskip \noindent that is

$$
\omega \in \bigcap_{k\geq 1} \bigcup_{n\geq 0} \bigcap_{p\geq n} \bigcap_{q\geq n}  (|f_{p}(\omega)-f_{q}(\omega)| \leq 1/k),
$$

\bigskip \noindent This leads to
$$
\Omega_0 \subset \bigcap_{k\geq 1} \bigcup_{n\geq 0} \bigcap_{p\geq n} \bigcap_{q\geq n}  (|f_{p}(\omega)-f_{q}(\omega)| \leq 1/k),
$$

\bigskip \noindent that is, upon taking the complements,

$$
\bigcup_{k\geq 1} \bigcap_{n\geq 0} \bigcup_{p\geq n} \bigcup_{q\geq n}  (|f_{p}(\omega)-f_{q}(\omega)| > 1/k) \subset \Omega_0^c.
$$

\bigskip \noindent Since $\Omega_0^c$ is a null-set, we readily get Formula Formula (CAU01), that is $(f_{n})_{n\geq 1}$ is an $m$-\textit{a.e.}  Cauchy sequence.\\

\bigskip \noindent (b) Suppose that $(f_{n})_{n\geq 1}$ is an $m$-\textit{a.e.} Cauchy sequence. Denote

$$
\omega_0=\{\omega \in \Omega : (f_{n}(\omega))_{n\geq 1} \ is \ a \ Cauchy \ sequence \ in \ \mathbb{R} \}.
$$

\bigskip  \noindent Since $\mathbb{R}$ is a complete space, each sequence $(f_{n}(\omega))_{n\geq 1}$ converges to a finite real number $f(\omega)$. If, for each $n\geq 0$, we denote by $\tilde{f}_{n}$ the restriction of $f_n$ on $\Omega_0$,
we just obtained that the sequence $(\tilde{f}_{n})_{n\geq 0}$ of measurable applications defined on the induce measure space $(\Omega_0, \mathcal{A}_{\Omega_0})$ converges to $f$. Hence $f$ is defined on $\Omega_0$ and is $\mathcal{A}_{\Omega_0}$-measurable.\\

\noindent Now, by Exercise 8 in Doc 02-03 i Chapter \ref{02_applimess}, page \pageref{exercise08_doc02-03}, we may extend $f$ to an $\mathcal{A}$-measurable function on $\Omega$ by assigning it an arbitrary
value on $\Omega_0^c$. Let us name that extension by $\tilde{f}$.\\

\noindent Finally the sequence $f_n$ converges to $\tilde{f}$ on $\Omega_0$ whose complement is a null-set. Hence, $f_n$ \textit{a.e.} converges to an \textit{a.e.} finite function. $\blacksquare$\\

\bigskip \noindent \textbf{Exercise 5}. \label{exercises05_sol_doc06-08}\\

\noindent In all this exercises, the functions $f_n$, $g_n$, $f$ and $g$ are \textbf{a.e.} finite.\\

\noindent Question (1). Show that the limit in measure is unique \textbf{a.e.}.\\

\noindent Hint : Suppose $f_n \rightarrow_m f_{(1)}$ and $f_n \rightarrow_m f_{(2)}$. By Formula (OS1) in Exercise 1 :
$$
(|f_1-f_2|> 1/k) \subset (|f_n-f_{(1)}|> 1/(2k)) \cup (|f_n-f_{(2)}|> 1/(2k)) 
$$

\noindent Use also (See Formula (C) in Exercise 7 in DOC 04-06, page \pageref{exercise07_sol_doc04-06}, Chapter \ref{04_measures}).\\

$$
(f_1-f_2 \neq 0)= \bigcup_{k\geq 1} (|f_1-f_2|> 1/k).
$$

\bigskip \noindent Question (2). Let $c$ be finite real number. Show that the sequence of \textit{a.e.} constant functions, $f_{n}=c$ \textit{a.e.}, converges in measure to $f=c$.\\

\bigskip \noindent Question (3).  Let $f_{n}\rightarrow _{m} f$ and $g_{n}\rightarrow _{m}g,$ $a\in \mathbb{R}$. Show that\\

\noindent \textbf{(a)} $f_{n}+g_{n}\rightarrow_{m}f+g$.\\

\noindent \textbf{(b)} $af_{n} \rightarrow af$.\\

\noindent Question (4).  Suppose that  $f_{n}\rightarrow_m A$ and $g_{n}\rightarrow_m B$, where $A$ and $B$ are finite real numbers. Let $H(x,y)$ a continuous function of $(x,y) \in D$, where $D$ is an open set of $\mathbb{R}^2$. Let $a$ and $b$ be finite real numbers. Show the following properties.\\

\noindent \textbf{(a)} $af_{n}+bg_{n}\rightarrow_m aA+bB$.\\

\noindent \textbf{(b)} $f_{n}g_{n}\rightarrow_m AB$.\\

\noindent \textbf{(c)} If $B \neq 0$,  then

$$
f_{n}/g_{n}\rightarrow_m A/B.
$$

\bigskip \noindent \textbf{(d)} If $(f_{n}, g_n)_{n\geq 1} \subset D$ \textit{a.e.} and $(A,B) \in D$, then

$$
H(f_{n},g_{n})\rightarrow_m H(f,g).
$$

\bigskip \bigskip \noindent \textbf{SOLUTIONS}.\\

\noindent Question (1). Suppose that $f_n \rightarrow_m f_{(1)}$ and $f_n \rightarrow_m f_{(2)}$. By Formula (OS1) in Exercise 1, we have for all $k\geq 1$

$$
(|f_1-f_2|> 1/k) \subset (|f_n-f_{(1)}|> 1/(2k)) \cup (|f_n-f_{(2)}|> 1/(2k)), 
$$

\noindent and hence by the sub-additivity of the measure $m$, we have

$$
0\leq m(|f_1-f_2|> 1/k) \leq m(|f_n-f_{(1)}|> 1/(2k)) + m(|f_n-f_{(2)}|> 1/(2k)). 
$$

\bigskip \noindent The right-hand member of the inequality tends to zero as $n\rightarrow +\infty$ by assumptions. By the sandwich theorem, we have for all $k\geq 1$,

$$
m(|f_{(1)}-f_{(2)}|> 1/k)=0.
$$

\bigskip \noindent Since we have  (See Formula (C) in Exercise 7 in DOC 04-06, page \pageref{exercise07_sol_doc04-06}, Chapter \ref{04_measures}),

$$
(f_{(1)} \neq f_{(2)})= \bigcup_{k \geq 1}(|f_{(1)}-f_{(2)}|> 1/k),
$$

\bigskip \noindent we get $m(f_{(1)} \neq f_{(2)})=0$, since $(f_{(1)} \neq f_{(2)})$ is a countable union of null-sets. This means that $f_{(1)}=f_{(2)}$ \textit{a.e.} $\square$\\

\noindent Question (2). Let $f_n=c$ for all $n \geq 0$, where $c$ is a finite real numbers. Set $f=c$ \textit{a.e.} We have for all $n\geq 0$, for all $\varepsilon>0$ :

$$
m(|f_n-f|>\varepsilon)=m( 0>\varepsilon)=m(\emptyset)=0.
$$

\bigskip \noindent Hence $f_n \rightarrow_m f$. $\square$\\

\noindent Question (3).\\

\noindent (a) By Formula (OS1) applied to $|(f_n+g_n)-(f+g)|=|(f_n-f) - (g_n-g)|$, we have for any $\varepsilon>0$ :

$$
(|(f_n+g_n)-(f+g)|> \varepsilon) \subset (|f_n-f|> \varepsilon/2) \cup (|g_n-g|> \varepsilon/2)
$$

\bigskip \noindent which leads to

$$
0\leq m(|(f_n+g_n)-(f+g)|> \varepsilon) \leq m(|f_n-f|> \varepsilon/2) +m(|g_n-g|> \varepsilon/2),
$$

\bigskip \noindent and we may conclude by applying the assumption and by using the sandwich rule.\\

\noindent (b) if $a=0$, $af_n=0$ and $af=0$ both \textit{a.e.} (remind that the $f_n$'s and $f$ are \textit{a.e.} finite) and Question (2) ensures the convergence in measure. next, suppose that $a\neq 0$, we have for any $\varepsilon>0$,

$$
0\leq m(|af_n-af|>\varepsilon)=m(|a||f_n-f|>\varepsilon)=m(|f_n-f|>\varepsilon/|a|).
$$

\bigskip \noindent Since $\varepsilon_1=\varepsilon/|a|$ is positive and finite number, $m(|f_n-f|>\varepsilon/|a|) \rightarrow 0$ by assumption. Applying the sandwich rule concludes the solution. $\square$\\

\noindent Question (4).\\

\noindent Point (a). This is solved by the combination of Points (a) and (b) in the general case in Question (3).\\

\noindent Point (b). We have

$$
f_ng_n-AB=f_ng_n-f_nB +f_nB - AB=f_n (g_n-B) + (B (f_n-f)).
$$
 
\bigskip \noindent By application of Formula (OS1), we get for any $\varepsilon/2$,

$$
0 \leq m(|f_ng_n-AB| > \varepsilon)\leq =m(|f_n| |g_n-B| \varepsilon/2) + m(|B| |f_n-f|>\varepsilon/2). \ \ (F1)
$$

\bigskip  \noindent We show that 
$$
m(|B| |f_n-f|>\varepsilon/2)\rightarrow 0, \ \ (F2)
$$ 

\bigskip \noindent exactly as in Point (b) of Question (3). Let us prove that 
$$
m(|f_n| |g_n-B| \varepsilon/2)\rightarrow 0.
$$

\bigskip \noindent But we have, for any $\eta>0$,

\begin{eqnarray*}
(|f_n| |g_n-B| \varepsilon/2)&=&(|f_n| |g_n-B| \varepsilon/2)\cap (|f_n-A|\leq \eta)\\
&+& (|f_n| |g_n-B| \varepsilon/2) \cap (|f_n-A|>\eta). \ \ (F3)
\end{eqnarray*}

\bigskip \noindent On $(|f_n| |g_n-B| \varepsilon/2)\cap (|f_n-f|\leq \eta)$, we have 
$$
|f_n|< \max(A+\eta, A-\eta)=M\geq 0,
$$

\bigskip which implies, from  that $|f_n| |g_n-B|>\varepsilon/2$, implies that
$A |g_n-B| \varepsilon/2$. We get

\begin{eqnarray*}
m(|f_n| |g_n-B| \varepsilon/2)\cap (|f_n-A|\leq \eta) &\leq& m(A |g_n-B| \varepsilon/2)\cap (|f_n-A|< \eta)\\
&\leq& m(A |g_n-B| \varepsilon/2). \ \ (F4)
\end{eqnarray*}

\bigskip \noindent For $\eta>0$, we show that 

$$
m(A |g_n-B| \varepsilon/2)\rightarrow 0. \ \ (F5)
$$ 

\bigskip \noindent exactly as in Point (b) of Question (3). Remark that

$$
m(|f_n| |g_n-B| \varepsilon/2) \cap (|f_n-A|>\eta) \leq m(|f_n-A|>\eta), \ \ (F6)
$$

\bigskip \noindent and that $m(|f_n-A|>\eta)\rightarrow 0$. Now (F1) becomes

$$
0 \leq m(|f_ng_n-AB| > \varepsilon) \leq m(A |g_n-B| \varepsilon/2)+m(|f_n-A|>\eta)+ m(|B| |f_n-f|>\varepsilon/2),
$$

\bigskip \noindent which, by combining Formulas (F2-F6) and by applying the sandwich rule, finishes the solution. $\square$\\

\bigskip \noindent Point (c). It is enough here to prove that if $f_n \rightarrow_m A\neq 0$, then $1/f_n \rightarrow_m 1/A$ and to combine with the result of Point (b). So suppose that $f_n \rightarrow_m A\neq 0$ with $|A|>0$. Fix $\eta>0$ such that $|A|-\eta>0$. we have for any $\varepsilon>0$,

$$
\left(|\frac{1}{f_n}-\frac{1}{A}| > \varepsilon \right)=\left(\frac{|f_n-A|}{|A||f_n|} > \varepsilon \right). \ \ (F7)
$$

\bigskip\noindent On $(|f_n-A|>\eta)$, we have $|f_n| \leq M=\max(A+\eta, A-\eta)$ and  $M=\max(A+\eta, A-\eta)$ is either $|A|+\eta$ or $|A|-\eta$ and thus positive in any case.\\

\noindent On $(|f_n-A|>\eta) \cap \left(|\frac{1}{|f_n-A|}-\frac{1}{|A||f_n|}| > \varepsilon \right)$, we have

$$
\varepsilon <|\frac{1}{|f_n-A|}-\frac{1}{|A||f_n|}| \leq \frac{|f_n-A|}{M|A|}. 
$$

\bigskip\noindent This leads to

\begin{eqnarray*}
\left(|\frac{1}{f_n}-\frac{1}{A}| > \varepsilon \right)&=&\left(|\frac{1}{f_n}-\frac{1}{A}| > \varepsilon \right) \bigcap (|f_n-A|>\eta)\\
&+& \left(|\frac{1}{f_n}-\frac{1}{A}| > \varepsilon \right) \bigcap (|f_n-A|>\eta)\\
&\subset& (|f_n-A|>\eta) \bigcup \left( \frac{|f_n-A|}{M|A|}>\varepsilon \right).
\end{eqnarray*}

\noindent and next to

$$
\leq m\left(|\frac{1}{f_n}-\frac{1}{A}| > \varepsilon \right) \leq m(|f_n-A|>\eta) + m\left( |f_n-A|> M|A|\varepsilon \right).
$$

\bigskip \noindent From this, we may conclude based on the convergence in measure to $A$.\\

\bigskip \noindent Point (d). By definition, the continuity of $H$ at $(A,B)$ implies that : for any $\varepsilon>0$, there exists $\eta>0$ such that

$$
\forall (x,y)\in D, (|x-A|\leq \eta, \ |y-B|\leq \eta) \Rightarrow |H(x,y)-H(A,B)|\leq \varepsilon. \ \ (CH)
$$

\bigskip \noindent (Notice that we may use strict inequalities of not this definition of continuity). Now fix $\varepsilon>0$ and consider one value of $\eta$ such that Formula (CH) holds.\\

\noindent We have

\begin{eqnarray*}
B_n=(|f_n-A| \leq \eta, \ |g_n-B| \leq \eta)^c&=&(|f_n-A| > \eta, \ |g_n-B| \leq \eta)\\
&+&(|f_n-A| \leq \eta, \ |g_n-B| > \eta)\\
&+&(|f_n-A| > \eta, \ |g_n-B|> \eta). 
\end{eqnarray*}
 
\bigskip \noindent which implies

$$
B_n=(|f_n-A| \leq \eta, \ |g_n-A| \leq \eta)^c \subset (|f_n-A| > \eta) \bigcup (|g_n-B|> \eta). 
$$

\bigskip \noindent First, we have

$$
0 \leq m(B_n) \leq m(|f_n-A| > \eta) + m(|g_n-B|> \eta) \rightarrow 0, \ \ 
$$

\bigskip \noindent by the assumption on the convergence in measure of the sequences $f_n$ and $g_n$. Secondly, by Formula (CH) above, we have on $B_n^c$

$$
(|H(f_n,g_n)-H(A,B)| \leq \varepsilon).
$$

\bigskip \noindent Finally we have

\begin{eqnarray*}
(|H(f_n,g_n)-H(A,B)| > \varepsilon)&=&(|H(f_n,g_n)-H(A,B)| > \varepsilon) \cap B_n^c\\
&+& (|H(f_n,g_n)-H(A,B)| > \varepsilon) \cap B_n,
\end{eqnarray*}

\bigskip \noindent which, at the light of the previous lines, implies that the set $(|H(f_n,g_n)-H(A,B)| > \varepsilon)$ is a subset of

$$
(|H(f_n,g_n)-H(A,B)| > \varepsilon, |H(f_n,g_n)-H(A,B)| \leq \varepsilon ) \bigcup  B_n. 
$$

\bigskip \noindent Since the first term of union is empty, we get

$$
0 \leq m(|H(f_n,g_n)-H(A,B)| > \varepsilon) \leq m(B_n), \ \ (F8)
$$

\bigskip \noindent and we conclude by applying Formula (F8) and the sandwich rule.\\

\bigskip \noindent \textbf{Exercise 6}. \label{exercise06_sol_doc06-08}\\

\noindent Let $(f_{n})_{n\geq 1}$ be a sequence of \textit{a.e.} finite and measurable real-valued functions.\\

\noindent Question (a) Suppose that the $(f_{n})_{n\geq 1}$ converges in measure to an \textit{a.e.} finite, measurable and real-valued function $f$. Show that the sequence is
a Cauchy sequence in measure.\\

\noindent Question (b) Suppose that $(f_{n})_{n\geq 1}$ is a Cauchy sequence in measure and that there exists a sub-sequence $(f_{n_k})_{k\geq 1}$ of $(f_{n})_{n\geq 1}$ converging in measure to an \textit{a.e.} function $f$. Show that $(f_{n})_{n\geq 1}$ converges in measure to $f$.

\noindent Question (c) Suppose that $(f_{n})_{n\geq 1}$ is a Cauchy sequence in measure and that there exists a sub-sequence $(f_{n_k})_{k\geq 1}$ of $(f_{n})_{n\geq 1}$ converging \textbf{a.e.} to an \textit{a.e.} function $f$. Suppose in addition that the measure is finite. Show that $(f_{n})_{n\geq 1}$ converges in measure to $f$.\\

\noindent \textsl{Hint}. Combine Question (b) of the current exercise and Question (4) if Exercise 2.\\

\bigskip \noindent \textbf{SOLUTIONS}.\\

\noindent Question (a). By definition of the convergence in measure of $(f_{n})_{n\geq 1}$ to $f$, we have for all $\varepsilon>0$,

$$
m(|f_n-f|>\varepsilon/2) \rightarrow 0.
$$

\bigskip \noindent By applying Formula (OS1) of Exercise 1, we have

$$
(|f_p-f_q|>\varepsilon) \subset (|f_p-f|>\varepsilon/2) \cup (|f_q-f|>\varepsilon/2),
$$

\bigskip \noindent which, by the sub-additivity of the measure $m$, leads to

$$
0 \leq m(|f_p-f_q|>\varepsilon) \leq m(|f_p-f|>\varepsilon/2) + m(|f_q-f|>\varepsilon/2).
$$

\bigskip  \noindent By using the Sandwich rule to the formula above, we get for any $\varepsilon>0$,

$$
\leq m(|f_p-f_q|>\varepsilon) \rightarrow 0, \ as \ (p,q) \rightarrow (+\infty,+\infty),
$$

\bigskip \noindent which means that $(f_{n})_{n\geq 1}$ us Cauchy sequence in measure. $\square$\\

\bigskip \noindent Question (b). By exploiting Formulation (OS1) in Exercise 1, we get for any $p\geq 0$, for any $k\geq 0$,

$$
0 \leq m(|f_p-f|>\varepsilon) \leq m(|f_p-f_{n_k}|>\varepsilon/2) + m(|f_{n_k}-f|>\varepsilon/2).
$$

\bigskip \noindent Now let $n \rightarrow +\infty$ and $k \rightarrow +\infty$, $m(|f_p-f_{n_k}|>\varepsilon/2)$ goes to zero since $(f_{n})_{n\geq 1}$ is a Cauchy sequence in measure and  $m(|f_{n_k}-f|>\varepsilon/2)$ tends to zero since the sub-sequence $(f_{n_k})_{k\geq 1}$ converges to $f$ . By the Sandwich rule, the solution is over.\\

\noindent Question (c).Since $m$ is supposed finite in this question, Question (4) of Exercise 2 implies that $(f_{n_k})_{k\geq 1}$ converges to $f$ in measure. Hence, we find ourselves in the situation of Question (b) and the conclusion is immediate. $\square$\\

\newpage

\noindent \textbf{Exercise 7}. \label{exercises07_sol_doc06-08}\\

\noindent Let $(f_{n})_{n\geq 1}$ be a sequence of \textit{a.e.} finite and measurable real-valued functions.\\

\noindent Question (a). Suppose that $(f_{n})_{n\geq 1}$ is Cauchy in measure. Show that $(f_{n})_{n\geq 1}$ possesses a sub-sequence $(f_{n_k})_{k\geq 1}$ of $(f_{n})_{n\geq 1}$ converging both \textbf{a.e.} and in measure to an \textit{a.e.} function $f$.\\

\noindent \textbf{Detailed Hints}.\\

\noindent \textbf{Step 1}. From the hypothesis : for any $\varepsilon>0$, for all $\eta>0$

$$
\exists n\geq 1, \forall p\geq n, \ \forall q \geq n, \ m(|f_q - f_q| > \varepsilon) < \eta. \ \ (CAUM01)
$$

\bigskip \noindent Remark (without doing any supplementary work) that Formula (CAUM01) is equivalent to each of (CAUM02), (CAUM03) and (CAUM04), exactly as you already proved that (CAU01) is equivalent to each of (CAU02) and (CAU03) in Exercise 4 : 

$$
\exists n \geq 0, \forall p\geq 0, \ \forall q \geq 0, \ m(|f_{n+p} - f_{n+q}| > \varepsilon) < \eta, \ \ (CAUM02)
$$

$$
\exists n\geq 0, \forall p\geq 0,  \ m(|f_{n+p} - f_n| >\varepsilon) < \eta. \ \ (CAUM03)
$$

$$
\exists n\geq 0, \forall N\geq n, \forall p\geq 0,  \ m(|f_{N+p} - f_N| >\varepsilon) < \eta. \ \ (CAUM04)
$$

\bigskip \noindent From Formula (CAUM04) construct an increasing sequence $(n_k)_{k\geq 1}$, that is

$$
n_1 <n_2 < ... < n_k < n_{k+1} < ...
$$

\bigskip \noindent such that for all $k\geq 1$,

$$
\forall N\geq n_k, \ \forall p\geq 0, \ \ m\left(|f_{N+p} - f_N| > \frac{1}{2^k} \right) < \frac{1}{2^k}. \ (FI)
$$

\bigskip \noindent proceed as follows : For each $k=1$, take $\varepsilon=\eta=1/2$ and take $n_1$ as the value of $n$ you found in (CAUM04). For $k=2$, take $\varepsilon=\eta=1/2^2$ and take $n(2)$ as the value of $n$ you found in (CAUM04) and put $n_2=\max(n_1+1,n(2)$. It is clear that $n_1<n_2$ and, since

$$
\forall N\geq n(2), \  \forall p \geq 0, \ \ m\left(|f_{N+p} - f_{N)}| > \frac{1}{2^2} \right) < \frac{1}{2^2}, \ \ (FI2)
$$

\bigskip \noindent you will have

$$
N\geq n_2, \  \forall p \geq 0, \ \ m\left(|f_{N+p} - f_{N}| > \frac{1}{2^2} \right) < \frac{1}{2^2}, \ \ (FI2P)
$$

\bigskip \noindent since values of $N$ greater that $n_2$ in (FI2P) are also greater that $n(2)$ in $(FI2)$. Next, proceed from $k=2$ to $k=3$ as you did from $k=1$ to $k=2$.\\

\noindent \textbf{Step 2}. Set

$$
A_k=\left( \left| f_{n_{k+1}} - f_{n_{k}} \right| > \frac{1}{2^k}\right), \ k\geq 1, \ and \ B_N=\bigcup_{k\geq N} A_k, \ N\geq 0.
$$

\bigskip \noindent Show that

$$
m(B_N) < \frac{1}{2^{N-1}}.
$$

\bigskip \noindent Fix $\varepsilon>0$ and take $N$ satisfying $2^{-(N-1)} < \varepsilon$. Show that : on $B_n^c$ we have for all $p\geq 0$
$$
\left| f_{n_{N+p}} - f_{n_{N}}\right| \geq \frac{1}{2^{N-1}} \leq \varepsilon. \ \ (FI3)
$$
 
\bigskip \noindent where you compose the difference $f_{n_{N+p}} - f_{n_{N}}$ from the increments $f_{n_{N+i}} - f_{n_{N+i-1}}$, $i=1, ..., p$, and next extend the bound to all the
increments for $k\geq N$.\\

\noindent Next using this with the obvious inclusion

$$
\bigcap_{q\geq 0} \bigcup_{p\geq 0} \left( \left| f_{n_{q+p}} - f_{n_{q}}\right| > \varepsilon \right)  \subset  \bigcup_{p\geq 0} \left( \left| f_{n_{N+p}} - f_{n_{N}}\right| > \varepsilon \right),
$$

\bigskip \noindent show that 

$$
m\left(\bigcap_{q\geq 0} \bigcup_{p\geq 0} \left( \left| f_{n_{N+p}} - f_{n_{N}}\right| \geq \varepsilon \right) \right) \leq m(B_N).
$$

\bigskip \noindent Deduce from this that for all $\varepsilon$,

$$
m\left(\bigcap_{q\geq 0} \bigcup_{p\geq 0} \left( \left| f_{n_{q+p}} - f_{n_{q}}\right| > \varepsilon \right) \right)
\leq m\left(\bigcap_{q\geq 0} \bigcup_{p\geq 0} \left( \left| f_{n_{q+p}} - f_{n_{q}}\right| \geq \varepsilon \right) \right)=0.
$$

\bigskip \noindent Apply this to $\varepsilon=1/r$, $r\geq 1$ integer, compare with Formula (CAU03) in Exercise 4 above, to see that sub-sequence $(f_{n_k})_{k\geq 1}$ is an \textit{a.e.} Cauchy sequence. By  Question 3 in Exercise 4, consider the \textit{a.e} finite function to which $(f_{n_k})_{k\geq 1}$ converges \textit{a.e.}\\

\noindent Consider the measurable sub-space  $\Omega_0$ whose complement is a null-set and on which $(f_{n_k})_{k\geq 1}$ converges \textit{a.e.} to $f$.\\

\noindent \textbf{Step 3}. Consider Formula (FI3). Place yourself on $\Omega_0 \cap B_N^c$ and let $p$ go to $\infty$. What is the formula you get on $\Omega_0 \cap B_N^c$?\\

\noindent Deduce the formula

$$
m\left( \left|f - f_{n_{N}}\right| > \varepsilon \right) \leq m\left( \left|f - f_{n_{N}}\right| \geq  \varepsilon \right) \leq m(B_N),
$$

\bigskip \noindent for all $N\geq 1$  satisfying $2^{-(N-1)} < \varepsilon$. Deduce that $f_{n_{N}}$ converges to $f$ in measure as $N \rightarrow +\infty$.\\

\noindent \noindent \textbf{Step 4}. By using Question (b) in Exercise 6, conclude that $f_{n}$ converges to $f$ in measure as $n \rightarrow +\infty$.\\

\noindent Write your conclusion.\\

\noindent Question (b). Combine this with the results of Exercise 6 above to state that a sequence $(f_{n})_{n\geq 1}$ is a Cauchy in measure if and only if it converges to a an \textit{a.e.} finite and measurable real-valued function.\\

\newpage
\bigskip \noindent \textbf{SOLUTIONS}.\\

\noindent Question (a). \\

\noindent We are going to closely follow the hints.\\

\noindent \textbf{Step 1}. We asked to do nothing about the equivalence between the formulas (CAUM0X). But, at least, let us say this : (CAUM01) is exactly the mathematical expression of the fact that $(f_{n})_{n\geq 1}$ is Cauchy in measure. And it is straightforward that (CAUM01) implies (CAUM04).\\

\noindent The construction of the sequence $n_k$ is clear as described by the hints and there is nothing else to do.

\bigskip \noindent \textbf{Step 2}. From the definition of $A_k$, $k\geq 1$ and $B_N$, $N\geq 0$, we have by the $\sigma$-sub-additivity of $m$,

$$
m(B_N) \leq \sum_{k\geq N} \frac{1}{2^{k}}=\frac{1}{2^{N}} \sum_{k\geq 0} \frac{1}{2^{k}}=\frac{1}{2^{N-1}}. \ \ (FI6)
$$

\bigskip \noindent Fix $\varepsilon>0$ and take $N\geq 1$ satisfying $2^{-(N-1)} < \varepsilon$. For 
$$
\omega \in B_N^c=\bigcap_{k\geq N} A_k^c,
$$

\bigskip \noindent we have for all $k\geq N$,
 
$$
\omega \in A_k^c \Rightarrow  |f_{n_{k+1}} - f_{n_{k}}|(\omega)\leq 2^{-k}.
$$

\bigskip \noindent But by the triangle inequality, we have
 
$$
\left|f_{n_{N+p}} - f_{n_{N}}\right|(\omega) \leq \sum_{1\leq i \leq p} |f_{n_{N+i}} - f_{n_{N+i-1}}|(\omega).
$$
 
\bigskip \noindent But the sum of absolute values of the increments $|f_{n_{N+i}} - f_{n_{N+i-1}}|(\omega)$, $i=1, ..., p$, is sub-sum of the sum of all absolute values of the increments $|f_{n_{k+1}} - f_{n_{k}}|(\omega)$ for $k=N$ (in absolute values). Hence, for all $p\geq 0$,

$$
\left|f_{n_{N+p}} - f_{n_{N}}\right|(\omega) \leq \sum_{k\geq N} \leq |f_{n_{k+1}} - f_{n_{k}}|(\omega) \leq \sum_{k\geq N} \frac{1}{2^{k}}=\frac{1}{2^{N-1}}<\varepsilon.
$$

\bigskip \noindent Hence, $|f_{n_{N+p}} - f_{n_{N}}|<\varepsilon$ on $B_N^c$, that is  $B_N^c \subset (|f_{n_{N+p}} - f_{n_{N}}|<\varepsilon)$. Taking the complements leads to

$$
(|f_{n_{N+p}} - f_{n_{N}}|\geq \varepsilon) \subset B_N. \ \ (FI7)
$$

\bigskip \noindent By using the following obvious inclusion, 
$$
\bigcap_{q\geq 0} \bigcup_{p\geq 0} \left( \left| f_{n_{q+p}} - f_{n_{q}}\right| \geq  \varepsilon \right)  \subset  \bigcup_{p\geq 0} \left( \left| f_{n_{N+p}} - f_{n_{N}}\right| \geq \varepsilon \right),
$$

\bigskip \noindent and, and by taking (FI6) and (FI7) into account, we see that all the sets $\left( \left| f_{n_{N+p}} - f_{n_{N}}\right| \leq  \varepsilon \right)$, $p\geq 0$, are in $B_N$. We get

$$
0 \leq m\left( \bigcap_{q\geq 0} \bigcup_{p\geq 0} \left( \left| f_{n_{q+p}} - f_{n_{q}}\right| \geq \varepsilon \right) \right) \leq m(B_N).
$$

\bigskip \noindent Since this holds for all $N\geq 1$ satisfying $2^{-(N-1)}$, we may let $N \rightarrow +\infty$, which entails that $m(B_N) \rightarrow 0$ and next,

$$
m\left(\bigcap_{q\geq 0} \bigcup_{p\geq 0} \left( \left| f_{n_{q+p}} - f_{n_{q}}\right| \geq \varepsilon \right) \right)=0. \ \ (F8)
$$

\noindent We may and do replace the inequality $\geq$ by the strict one by considering for example, for all $\varepsilon>0$ 

$$
\left(\bigcap_{q\geq 0} \bigcup_{p\geq 0} \left( \left| f_{n_{q+p}} - f_{n_{q}}\right| > \varepsilon \right) \right) \subset \left(\bigcap_{q\geq 0} \bigcup_{p\geq 0} \left( \left| f_{n_{q+p}} - f_{n_{q}}\right| \geq \varepsilon \right) \right).
$$

\bigskip \noindent By using this Formula (F8) with $\varepsilon=1/r$, $r\geq 1$ integer and by comparing with Formula (CAU03) in Exercise 4 above, we see that sub-sequence $(f_{n_k})_{k\geq 1}$ is an \textit{a.e.} Cauchy sequence. By  Question 3 in Exercise 4, we know that it converges \textit{a.e.} to a \textit{a.e.} function, say $f$.\\

\noindent Let us consider the measurable sub-space  $\Omega_0$ whose complement is a null-set and on which $(f_{n_k})_{k\geq 1}$ converges \textit{a.e.} to $f$.\\

\noindent \textbf{Step 3}. Let us place yourself on $\Omega_0 \cap B_N^c$ and let $p$ go to $\infty$. Thus $f_{n_{q+p}}$ converges to $f(\omega)$ and the inequality becomes

$$
\left| f_ - f_{n_{N}}\right| \leq \frac{1}{2^{N-1}} < \varepsilon \ on \ \Omega_0 \cap B_N^c,
$$

\bigskip \noindent which implies $\Omega_0 \cap B_N^c \subset (\left| f_ - f_{n_{N}}\right| < \varepsilon)$, that is upon taking the complements,

$$
(\left| f_ - f_{n_{N}}\right| \geq \varepsilon) \subset B_N \cup \Omega_0,
$$

\bigskip \noindent and hence, $N\geq 1$ satisfying $2^{-(N-1)}$ (that is $N \leq N_0=[1 +\log(/\varepsilon)/\log 2)]$),

$$
m(\left| f_ - f_{n_{N}}\right| > \varepsilon) \leq m(\left| f_ - f_{n_{N}}\right| \geq  \varepsilon) \leq m(B_N \cup \Omega_0)=m(B_N)< \varepsilon,
$$

\bigskip \noindent which establishes that $f_{n_{N}}$ converges to $f$ in measure as $N \rightarrow +\infty$.\\

\noindent \textbf{Conclusion} : A Cauchy sequence in measure possesses a sub-sequence converging both in measure and \textit{a.e.} to an \textit{a.e.} finite function.$\square$\\

\bigskip \noindent Question (b). This question is straightforward and does not need a solution.\\

\noindent 

\noindent \LARGE \textbf{DOC 06-09 : Convergence Theorems - Exercises with Solutions.} \label{doc06-09}\\
\bigskip
\Large

\bigskip \noindent \textbf{NB}. In this document and in the sequel of the textbook, the Measure Space Reduction Principle (see pages \pageref{msrp} and \pageref{msrpA}) which consists in placing ourselves in a sub-measurable space $\Omega_0$ whose complement is a null-set with respect to the working measure, and on which all the countable \textit{a.e.} properties simultaneously hold for $\omega \in \Omega_0$. By working on the induced measure on $\Omega_0$, the values of the measures of sets and the values of the integrals of real-valued and measurable function remain unchanged.\\

\noindent \textbf{Exercise 1}. \label{exercise01_sol_doc06-09}

\noindent Question (1). Let $(f_{n})_{n\geq 0}$ be a non-decreasing sequence of \textit{a.e.} non-negative real-valued and measurable functions. Show that
\begin{equation*}
\lim_{n \uparrow +\infty} \int f_{n}dm = \int \lim_{n \uparrow +\infty} f_n \ dm.
\end{equation*}

\noindent \textit{Hint}.\\

\noindent (a) First denote $f=\lim_{n \uparrow +\infty} f_n$.\\

\noindent By  \textit{Point (02.14) in Doc 02-01 in Chapter \ref{02_applimess}} (see page \pageref{doc02-01}), consider for each $n\geq 0$, a non-decreasing sequence of non-negative elementary functions $(f_{n,k})_{k\geq 0}$ such that $f_{n,k} \uparrow f_n$ as $k\uparrow +\infty$.\\

\noindent Next set

$$
g_k= \max_{k \leq n} f_{n,k}, \ n\geq 0.
$$

\bigskip \noindent  By \textit{Point (02.13) in Doc 02-01 in Chapter \ref{02_applimess}} (see page \pageref{doc02-01}), see that the $g_k$ are non-negative elementary functions.\\

\noindent Write $g_k$ and $g_{k+1}$ in extension and say why the sequence $(g_k)_{k\geq 0}$ is non-decreasing.\\

\noindent (b) Show that

$$
\forall 0 \geq n, \ f_{n,k} \leq g_k \leq f_k, \ \ (I1)
$$

\bigskip \noindent  and deduce from this by letting first $k\rightarrow +\infty$ and next $n\rightarrow +\infty$ that

$$
\lim_{k\rightarrow +\infty} g_k=f. \ \ (I2)
$$

\bigskip \noindent  Why all the limits you consider exist?\\

\noindent (c) Now by taking the integrals in (I1), by letting first $k\rightarrow +\infty$ and next $n\rightarrow +\infty$, and by using the definition of the integral of 
$\int \lim_{k\rightarrow +\infty} g_k \ dm$ and by raking (I2) into account, conclude that

$$
\lim_{k\rightarrow +\infty} \int f_k \ dm= \int \lim_{k\rightarrow +\infty}  f_k \ dm. \ \square.
$$

\bigskip \noindent  Question (2) Extend the result of Question (1) when $(f_{n})_{n\geq 0}$ is a non-decreasing sequence of functions bounded below \textit{a.e.} by an integrable function $g$, that is : $g \leq f_n$, for all $n\geq 0$. (Actually, it is enough that $g \leq f_n$ for $n\geq n_0$, where $n_0$ is some integer).\\

\noindent \textit{Hint}. Consider $h_n=f_n-g$ and justify that you are in the case of Question (1). Next, do the mathematics.\\

\bigskip \noindent  Question (2) Extend the result of Question (2) when $(f_{n})_{n\geq 0}$ is a non-decreasing sequence of functions such that there exists $n_0\geq 0$, such that $\int f_{n_0} \ dm >-\infty$.\\

\noindent \textit{Hint}. Discuss around the two cases : $\int f_{n_0} \ dm=+\infty$ and $\int f_{n_0} \ dm<+\infty$.\\

\bigskip \noindent \textbf{SOLUTIONS}.\\

\noindent (a) The functions $g_k$, as finite maxima of non-negative elementary functions are non-negative elementary functions. We have for all $k\geq 0$,

$$
g_k=\max(f_{0,k}, f_{1,k}, ..., f_{k,k}) \ and \ g_k=\max(f_{0,k+1}, f_{1,k+1}, ..., f_{k,k+1}) \vee f_{k+1,k+1}  
$$

\bigskip \noindent  and for $0\leq n \leq k$, $f_{n,k}\leq f_{1,k+1}$, since the sequence $(f_{n,k})_{k\geq 0}$ in non-decreasing to $f_n$. Thus the sequence is non-decreasing.\\

\noindent Next, the elements $f_{n,k}$ used to form $g_k$ are such that for each $n\leq k$, $f_{n,k}\leq f_n \leq f_k$ by the non-decreasingness of the sequence $(f_h)_{h\geq 0}$. Hence $g_k \leq f_k$. This and the definition of $g_k$ lead to Formula (I2).\\

\noindent All the sequences used here and the sequences of their integrals are non-decreasing. So their limits exist in $\overline{\mathbb{R}}$ and we do need need to justify them.\\

\noindent Now, by the definition of the integrals of non-negative functions, we have

$$
\int f_n \ dm = \lim_{k\rightarrow +\infty} \int f_{n,k} \ dm, \ \ (I4)
$$

\bigskip \noindent 

$$
\int \lim_{k\rightarrow +\infty}  g_k \ dm = \lim_{k\rightarrow +\infty} \int g_k \ dm. \ \ (I5)
$$

\bigskip \noindent  We are going to conclude in the following way :\\

\noindent By letting first $k\rightarrow +\infty$ and next $n\rightarrow +\infty$ in Formula (I2), we obviously get that  
$$
\lim_{k\rightarrow +\infty} g_k=f. \ \ (I6)
$$

\bigskip \noindent  By taking the integrals in Formula (I2) and by letting first $k\rightarrow +\infty$ and next $n\rightarrow +\infty$, we get

$$
\lim_{n\rightarrow +\infty} \int f_n \ dm \leq  \lim_{k\rightarrow +\infty}  \int g_k \ dm \leq \lim_{k\rightarrow +\infty}  \int f_k \ dm.
$$

\bigskip \noindent  But the two extreme terms are equal and the middle term, by Formula (I5) is $\int \lim_{k\rightarrow +\infty}  g_k \ dm$, which by (I6) is $\int f  g_k \ dm$. $\square$\\

\noindent Question (2). If $g\leq f_n$ \textit{a.e.} for each $n \geq 0$, with $g$ integrable, we get $g$ is \textit{a.e.} finite and the $h_n=f_n-g$ are \textit{a.e.} defined and non-negative and form e non-decreasing function. We may and do apply the results of Question (1) and get

$$
\lim_{n\rightarrow +\infty} \int (f_n-g) \ dm =  \int \lim_{n\rightarrow +\infty} (f_n-g) \ dm.
$$ 

\bigskip \noindent  In both sides of the equation, we may and do use the linearity of the integrals and the limits and next from the finite number $\int g \ dm$ from both sides to conclude.\\

\noindent Question (3). We have to situation. Either $\int f_{n_0} \ dm=+\infty$.  Hence, for all $n\geq n_0$, we have $\int f_n \ dm=+\infty$ and as well $\int f \ dm=+\infty$ due to the non-decreasingness of the sequence of the sequence $(f_n)_{n\geq 0}$. We surely have, for $n\geq n_0$, $\int f_n \ dm=+\infty \rightarrow +\infty=\int f \ dm$\\

\noindent Or $\int f_{n_0} \ dm<+\infty$. From there, we apply the result of Question (2) to the sequence $(f_n)_{n \geq n_0}$ which is bounded below by $f_{n_0}$.\\

\bigskip \noindent \textbf{Exercise 2}. (Fatou-Lebesgue's lemma and Dominated Convergence Theorem) \label{exercise02_sol_doc06-09}\\

\noindent Question (1). (Fatou-Lebesgue's lemma) Let $(f_{n})_{n\geq 0}$ be a sequence real-valued and measurable functions \textit{a.e.} bounded \textbf{below} by an integrable function $h$, that is 
$h \leq f_n$ \textit{a.e.} for all $n\neq 0$. Show that

\begin{equation*}
\liminf_{n\rightarrow +\infty}\int f_{n}dm \geq \int \liminf_{n\rightarrow +\infty }f_{n}dm. 
\end{equation*}

\noindent \textit{Hints}.\\

\noindent (a) Denote for all $n\geq 0$

$$
g_n=\inf_{k\geq n} (f_k-h).
$$

\noindent The $g_n$ are they defined? Why? What are their sign? What is the limit of $g_n$? What is its monotonicity, if any?\\

\noindent Justify and apply the Monotone Convergence Theorem.\\

\noindent Compare $\int \inf_{k\geq n} (f_k-h) \ dm$ with each $\int (f_k-h) \ dm$, $k\geq n$ and deduce a comparison with $\inf_{k\geq n} \int (f_k-h) \ dm$.\\

\noindent Conclude. Show your knowledge on superior limits and inferior limits.\\

\noindent Question (2). Let $(f_{n})_{n\geq 0}$ be a sequence real-valued and measurable functions \textit{a.e.} bounded \textbf{above} by an integrable function $h$, that is 
$f_n \leq h$ \textit{a.e.} for all $n\neq 0$. Then, we have

$$
\limsup_{n\rightarrow +\infty}\int f_{n}dm \leq \int \limsup_{n\rightarrow +\infty }f_{n}dm.
$$

\bigskip \noindent  \textit{Hints}. Take the opposites of all the functions in this question and apply Question (1).\\

\noindent Question (3). Let $(f_{n})_{n\geq 0}$ be a sequence real-valued and measurable functions \textit{a.e.} such that :\\

\noindent (a) it is bounded  by an integrable function $g$, that is $|f_n| \leq g$ \textit{a.e.} for all $n\geq 0$,\\

\noindent (b) it converges in measure or \textit{a.e.} to a function $f$.\\

\noindent  Then $f$ is integrable and the sequence $(\int f_{n} \ dm)_{n\geq 0}$ has a limit with

$$
\lim_{n\rightarrow +\infty}\int f_{n} \ dm = \int f \ dm. \ \ (D1)
$$

\bigskip \noindent and

$$
\int |f| \ dm = \int g \ dm. \ \ (D2)
$$

\noindent \textit{Hints}. Do two cases.\\

\noindent Case 1 : $f_n \rightarrow f$ \textit{a.e}. Remark the your sequence is bounded below and above. Apply Questions (1) and (2) and replace $\lim_{n\rightarrow +\infty} f_n$ by its value. Combine the two results at the light of the natural order of the limit inferior and the limit superior to conclude.\\

\noindent Case 2 : $f_n \rightarrow_m f$. Use the\textit{ Prohorov criterion given in Exercise 4 in Doc 11-01 in Chapter \ref{11_appendix} (See page \pageref{doc11-01})} in the following way.\\

\noindent Consider an adherent point $\ell$ of the sequence $(\int f_n \ dm)_{n\geq 0}$, limit a sub-sequence $(\int f_{n_k} \ dm)_{k\geq 0}$.\\

\noindent Since we still have $f_{n_k} \rightarrow_m f$ as $k\rightarrow +\infty$, consider a sub-sequence $(f_n{_{k_j}})_{j\geq 0}$ of $(f_{n_k})_{k\geq 0}$ converging both \textit{a.e.} and in measure to a finite function $h$ (See Exercise 7 in Doc 06-04, page \pageref{doc06-04}). Justify that $h=f$ \textit{a.e.} and apply Case 1 to the sub-sub-sequence $(f_{n_{k_j}})_{j\geq 0}$ and conclude that 

$$
\ell= \int f \ dm.
$$

\bigskip \noindent Apply the Prohorov criterion and conclude.\\

\bigskip \noindent \textbf{SOLUTIONS}.\\

\noindent Question (1). Since $h$ is integrable, it is \textit{a.e} finite and the functions $f_k-h$ are defined and \textit{a.e.} non-negative. Besides the functions $g_n$ are also measurable, non-negative \textit{a.e.}. By the definition of inferior limits, we know that $g_n$ broadly increases to $\liminf_{k\rightarrow +\infty} f_n$. We may and do apply the Monotone Convergence theorem to have

$$
\int \liminf_{n\rightarrow +\infty} (f_n-h) \ dm = \lim_{n\rightarrow +\infty} \int  g_n \ dm. \ \ (I1)
$$

\bigskip \noindent  By definition, we have $g_n \leq (f_k-h)$ for $k\geq n$. We have  $\int g_n \ dm \leq \int (f_k-h) \ dm$ for $k\geq n$. We get

$$
\int g_n \ dm \leq \inf_{k\geq n} \int (f_k-h) \ dm. \ \ (I2)
$$

\noindent The right-hand member broadly increases to $\liminf_{n\rightarrow +\infty} \int (f_n-h) \ dm$. By letting $n\rightarrow +\infty$ in Formula (I2) and by taking Formula (I1) into account, we have

$$
\int \liminf_{n\rightarrow +\infty} (f_n-h) \ dm \leq  \liminf_{n\rightarrow +\infty} \int (f_n-h) \ dm.
$$

\bigskip \noindent  Let us the additivity of the inferior limit (when one the two term is a constant as a sequence) and the linearity of the integral (here, the integrability is important) to have

$$
(\int \liminf_{n\rightarrow +\infty} f_n \ dm) - (\int h \ dm)\leq (\liminf_{n\rightarrow +\infty} \int f_n dm) - (\int h \ dm).
$$  

\bigskip \noindent  Drop the finite number $\int h \ dm$ to conclude.\\

\noindent \textbf{Warning}. The inferior limit operator is over-additive and the inferior limit operator is sub-additive. But when of the two sequence is constant, we have the additivity of both of them.\\

\noindent Question (2). By passing to the opposites $-f_n$ and $-g$, we may apply Question (1) and use easy properties on limits superior to conclude. Once again, you are referred to ... for a round up on limits in $\mathbb{R}$.

\noindent Question (3).\\

\noindent \textit{Case 1}. We have $-g \leq f_n \leq g$, \textit{a.e.}, $n\geq 0$. So, we have both a below integrable bound and an above integrable bound. We may and do use the results of Questions (1) and (2) simultaneously to get

$$
\int \liminf_{n\rightarrow +\infty} f_n dm \leq \liminf_{n\rightarrow +\infty} \int  f_n dm \leq \limsup_{n\rightarrow +\infty} \int  f_n dm \leq \int  \limsup_{n\rightarrow +\infty} f_n dm.
$$

\bigskip \noindent  The first inequality comes from Question (1), the third from Question (2) and the second is the natural order between the limit inferior limit and the limit superior. Since the extreme members are both equal to $\int f \ dm$, by the assumption that $f=\lim_{n\rightarrow +\infty}$, we get

$$
\liminf_{n\rightarrow +\infty} \int  f_n dm = \limsup_{n\rightarrow +\infty} \int  f_n dm=\int f \ dm.
$$
  
\bigskip \noindent  This proves (D1). Next $|f_n| \leq g$, \textit{a.e.}, $n\geq 0$ implies

$$
| \int f_n \ dm | \leq \int |f_n| \ dm \leq \int g \ dm, \ n\geq 0.
$$ 

\bigskip \noindent By letting $n\rightarrow +\infty$ and by applying (D1), we get (D2).\\

\bigskip \noindent \textit{Case 2}. Suppose the $f_n \rightarrow_m f$ and condition (b) holds. We are going to show that the sequence has only adherent point.\\

\noindent Indeed let $\ell \in \overline{\mathbb{R}}$ be an adherent point of $(\int f_n \ dm)_{n\geq 0}$. This, it limit a sub-sequence $(\int f_{n_k} \ dm)_{k\geq 0}$.\\

\noindent We still have $f_{n_k} \rightarrow_m f$ as $k\rightarrow +\infty$. By (See Exercise 7 in Doc 06-04, page \pageref{doc06-04}), there exists a sub-sequence $(f_{n_{k_j}})_{j\geq 0}$ of $(f_ {n_k} )_{k\geq 0}$ converging both \textit{a.e.} and in measure to a finite function $h$. Since we also have $f_{n_{k_j}} \rightarrow_m f$ as $j\rightarrow +\infty$, we get by the \textit{a.e.} uniqueness of the limit in measure that $h=f$ \textit{a.e.}.\\

\noindent The sequence $(f_n{_{k_j}})_{j\geq 1}$ obviously satisfies the conditions of Case 1. We thus get :

$$
\int f_{n_{k_j}} \ dm \rightarrow \int f \ dm \ as \ j\rightarrow +\infty.
$$

\bigskip \noindent But as sub-sequence of $(\int f_{n_k} \ dm)_{k\geq 0}$, $(\int f_n{_{k_j}} \ dm)_{j\geq 0}$ converges to $\ell$, which by the uniqueness of limits in  $\overline{\mathbb{R}}$ leads to \\
$$
\ell= \int f \ dm.
$$

\bigskip \noindent and (D2) holds. Now, as the unique adherent point, $\ell$ is the limit superior and the limit inferior and thus the limit of $(\int f_n \ dm)_{n\geq 0}$. $\square$\\

\bigskip \noindent \textbf{Exercise 3}. (Young Dominated Convergence Theorem) \label{exercise03_sol_doc06-09}\\

\noindent Let $(f_{n})_{n\geq 0}$ be a sequence real-valued and measurable functions \textit{a.e.} such that :\\

\noindent (1) there exists a sequence of integrable functions $(h_n)_{n\geq 0}$ such that is $|f_n| \leq h_n$ \textit{a.e.} for all $n\geq 0$,\\

\noindent (2) the sequence $(h_n)_{n\geq 0}$ converges to an integrable $h$ such that 
$$
\lim_{n\rightarrow +\infty}\int h_{n} \ dm = \int h \ dm,
$$

\bigskip \noindent (3) the sequence $(f_{n})_{n\geq 0}$ converges in measure or \textit{a.e.} to a function $f$.\\

\noindent  Then $f$ is integrable and the sequence $(\int f_{n} \ dm)_{n\geq 0}$ has a limit with

$$
\lim_{n\rightarrow +\infty}\int f_{n} \ dm = \int f \ dm
$$

\bigskip \noindent and

$$
\int |f| \ dm = \int h \ dm.
$$

\bigskip \noindent \textit{Hint}. At the the method used in Exercise 2 above, you only need to prove an analog of the Fatou-Lebesgue lemma of the form : Let $(f_{n})_{n\geq 0}$ be a sequence real-valued and measurable functions \textit{a.e.} bounded \textbf{below} by a sequence of integrable functions $(h_n)_{n\geq 0}$, which converges to an integrable $h$ so that  
$$
\lim_{n\rightarrow +\infty}\int h_{n} \ dm = \int h \ dm.
$$

\bigskip \noindent Then we have

\begin{equation*}
\liminf_{n\rightarrow +\infty}\int f_{n}dm \geq \int \liminf_{n\rightarrow +\infty }f_{n}dm. \ \ (D3)
\end{equation*}

\bigskip \noindent Prove (D3) by repeating that of Question (1) in Exercise 2, but be careful since you will not have the linearity of the limit inferior. The method will lead you to the formula

$$
\int \liminf_{n\rightarrow +\infty} (f_n-h_n) \ dm \leq  \liminf_{n\rightarrow +\infty} \int (f_n-h_n) \ dm. \ \ (Y3)
$$  

\bigskip \noindent From there, use the super-additivity in the left-hand, use in the right-hand the the formula : for all $\eta>0$, there exists $N\geq 0$ such that for all $n\geq N$, 
$$
\int h_n \ dm \geq \int h \ dm
$$

\bigskip \noindent Combine the two formula $n\geq N$, let $n\rightarrow +\infty$ and next $\eta\rightarrow 0$, to conclude.\\

\bigskip \noindent \textbf{SOLUTIONS}. Let us proceed as the hint suggested. We have a  sequence real-valued and measurable functions $(f_{n})_{n\geq 0}$  which is \textit{a.e.} bounded \textbf{below} by a sequence of integrable functions $(h_n)_{n\geq 0}$, which converges to an integrable $h$ so that $\int h_n \ dm$ converges to $\int f \ dm$.\\

\noindent Let us follows the steps of the proof of Question 1 of Exercise 2.\\

\noindent Since each $h_n$, $n\geq 0$ is integrable, it is \textit{a.e} finite and the functions $f_k-h_k$ are defined and \textit{a.e.} non-negative. Besides the functions 
$$
g_n=\inf{h\geq n} f_k-h_k,
$$
 
\bigskip \noindent are also measurable, non-negative \textit{a.e.}. By the definition of inferior limits, we know that $g_n$ broadly increases to $\liminf_{k\rightarrow +\infty} f_n-g_n$. We may and do apply the Monotone Convergence theorem to have

$$
\int \liminf_{n\rightarrow +\infty} (f_n-g_n) \ dm = \lim_{n\rightarrow +\infty} \int  g_n \ dm. \ \ (Y1)
$$

\bigskip \noindent  By definition, we have $g_n \leq (f_k-h_k)$ for $k\geq n$. We have  $\int g_n \ dm \leq \int (f_k-h_k) \ dm$ for $k\geq n$. We get

$$
\int g_n \ dm \leq \inf_{k\geq n} \int (f_k-h_k) \ dm. (Y2)
$$

\bigskip \noindent The right-hand member broadly increases to $\liminf_{n\rightarrow +\infty} \int (f_n-h_k) \ dm$. By letting $n\rightarrow +\infty$ in Formula (I2) and by taking Formula (I1) into account, we have

$$
\int \liminf_{n\rightarrow +\infty} (f_n-h_n) \ dm \leq  \liminf_{n\rightarrow +\infty} \int (f_n-h_n) \ dm. \ \ (Y3)
$$  

\bigskip \noindent Since the limit inferior is over-additive, we have

$$
\liminf_{n\rightarrow +\infty} (f_n-h_n) \geq \liminf_{n\rightarrow +\infty} f_n + \liminf_{n\rightarrow +\infty} (-h_n)=\liminf_{n\rightarrow +\infty} f_n -h, \ \ (Y4)
$$

\bigskip \noindent and thus

\begin{eqnarray*}
\int \liminf_{n\rightarrow +\infty} (f_n-h_n) \ dm &\geq& \int \liminf_{n\rightarrow +\infty} f_n + \liminf_{n\rightarrow +\infty} (-h_n) \ dm\\
&=&\int \liminf_{n\rightarrow +\infty} f_n -h \ dm\\
&=& \int \liminf_{n\rightarrow +\infty} f_n - \int h \ dm, \ \ (Y5)\\
\end{eqnarray*}

\noindent The right-hand member is

$$
\liminf_{n\rightarrow +\infty} (\int f_n \ dm) - (\int h_n \ dm). \ \  (Y6)
$$

\bigskip \noindent  By assumption, we have for any $\eta>0$, there exists $N\geq 0$ such that for all $n\geq N$,

$$
| (\int h_n \ dm) - (\int h_n \ dm)|<\eta.
$$

\bigskip \noindent By combining Formulas (Y3)-(Y6), we get for all $n\geq N$,

\begin{eqnarray*}
\int \liminf_{n\rightarrow +\infty} f_n - \int h \ dm &\leq& \liminf_{n\rightarrow +\infty} \biggr( (\int f_n \ dm) - (\int h \ dm) +\eta \biggr)\\
&=& \biggr(\liminf_{n\rightarrow +\infty}  \int f_n \ dm \biggr) - (\int h \ dm) +\eta. 
\end{eqnarray*}

\noindent By dropping the finite number $\int h \ dm$ from both sides and by letting $\eta\rightarrow 0$, we finish the solution.\\
\noindent \LARGE \textbf{DOC 06-10 : Applications of Convergence Theorems - Exercises with Solutions} \label{doc06-10}\\
\bigskip
\Large

\textbf{NB}. Recall that a continuous limit (in $\overline{\mathbb{R}}$, in the form 

$$
x(t) \rightarrow y  \ \ as \ \ t \rightarrow s,
$$ 

\bigskip \noindent may be discretized and shown to be equivalent to : for any sequence $(t_n)_{n\geq 0}$ such that $t_n\rightarrow s$ as $n\rightarrow +\infty$, we have

$$
x(t_n) \rightarrow y \ \ as \ \ n\rightarrow +\infty.
$$

\bigskip \noindent  Using the discretized from is more appropriate with measurability.\\

\textbf{Exercise 1}. \label{exercise01_sol_doc06-10} (Convergence of sums of applications)\\

\noindent Question (1). Let $(f_{n})$ be sequence of \textit{a.e.} non-negative measurable functions. Use the MCT and show that 
\begin{equation*}
\sum_{n=0}^{+\infty} \biggr( \int f_{n} \  dm \biggr) =\int \biggr( \sum_{n}^{+\infty} f_{n} \ dm \biggr). \ \ (S01)
\end{equation*}

\noindent Question (2). Let be $(a_{n})_{n\geq 0}=(\left( a_{n,p}\right) _{(n\geq 0,p\geq 0)}$ a rectangular sequence of non-negative real numbers. Use the MCT to the counting measure on $\mathbb{N}$ and show that

\begin{equation*}
\sum_{n\geq 0} \biggr( \sum_{p\geq 0}a_{n,p} \biggr)=\sum_{p\geq 0} \biggr(\sum_{n\geq 0}a_{n,p}\biggr). \ \ (S02)
\end{equation*}

\noindent Question (3) Let $(f_{n})$ be sequence of measurable functions such that%
\begin{equation*}
\int \sum_{n\geq 0} \left\vert f_{n}\right\vert  dm<+\infty.
\end{equation*}

\noindent Use the DCT, take into account Question (1) and establish that 

\begin{equation*}
\sum_{n\geq 0} \biggr( \int f_{n} \ dm \biggr)  =\int \biggr( \sum_{n\geq 0} f_{n} \biggr) \ dm.
\end{equation*}

\noindent Question (4). Let $(\left( a_{n,p}\right) _{(n\geq 0,p\geq 0)}$ a rectangular sequence such that
\begin{equation*}
\sum_{p\geq 0} \sum_{n\geq 0}\left\vert a_{n,p}\right\vert <\infty.
\end{equation*}

\bigskip \noindent Apply Question (3) on the counting measure on $\mathbb{N}$ and conclude
\begin{equation*}
\sum_{p\geq 0} \sum_{n\geq 0} a_{n,p}=\sum_{n\geq 0} \sum_{p\geq 0} a_{n,p}
\end{equation*}

\noindent Question (5). Let $m$ be a countable sum of measures $m_j$, $j\geq 0$ :

$$
m=\sum_{j\geq 0} m_j.
$$

\noindent Show that for any non-negative and measurable function $f$, we have

$$
\int f \ dm = \sum_{j\geq 0} \int f \ dm_j. \ \ (IS)
$$

\noindent Give simple conditions under which Formula (IS) holds for some for  measurable $f$ function.\\

\bigskip \noindent \textbf{SOLUTIONS}.\\

\noindent Question (1). Since the $f_n$'s are \textit{a.e.} non-negative, the sequence of functions $S_k=\sum_{n=0}^k f_n$, $k\geq 0$ is non-decreasing and converges to 
$S=\sum_{n=0}^{+\infty} S_n$, $k\geq 0$. The application of the MCT gives 

$$
\int S \ dm = \lim_{k\rightarrow +\infty} \int S_k \ dm,
$$

\noindent that is

$$
\int \sum_{n=0}^{+\infty} S_n =\lim_{k\rightarrow +\infty} \int \sum_{n=0}^k S_k \ dm=\lim_{k\rightarrow +\infty} \sum_{n=0}^k \int f_n \ dm = \sum_{n=0}^{+\infty} \int f_n \ dm. \ \square.
$$

\noindent Question (2). By denoting by $\nu$ the counting measure on $\mathbb{N}$ and by $f_p$, for each $n\geq 0$, the real-valued function defined on $\mathbb{N}$ such that $f_n(p)=a_{n,p}$, we see that 
(S02) is (S01) with $\nu=m$.\\

\noindent Question (3). By Question (1), we have

$$
\int \sum_{n\geq 0} \left\vert f_{n}\right\vert  dm =  \sum_{n\geq 0} \biggr( \int \left\vert f_{n}\right\vert  dm\biggr) <+\infty.
$$

\noindent Hence the functions are bounded by $g=\sum_{n\geq 0} \left\vert f_{n}\right\vert$ which is integrable. since $g$ is \textit{a.e.} finite, the series $\sum_{n\geq 0} f_{n}$ is \textit{a.e.} absolutely converging. We may and so apply the DCT to have

$$
\int S \ dm = \lim_{k\rightarrow +\infty} \int S_k \ dm,
$$

\noindent which concludes the solution.\\

\noindent Question (4). This is the question (3) applied to the counting measure on $\mathbb{N}$.\\

\noindent Question (5). Formula is nothing else than

$$
\int f \ dm =\sum_{j\geq 0} f \ dm_j,
$$

\noindent for all $f=1_A$, where $A \in \mathcal{A}$. This formula is readily extended to non-negative elementary function and next to non-negative and measurable functions $f$ by the definition of the
integral of non-negative and measurable functions and by Exercise 12 on the Monotone Convergence Theorem of non-negative series.\\

\noindent Finally, for a measurable function, if for $\sum_{j\geq 0} f^- \ dm_j$ or $\sum_{j\geq 0} f^+ \ dm_j$ is finite, thus $\int f \ dm$ exists and, by regrouping term by term, we get Formula (IS).\\

\bigskip \noindent \textbf{Exercise 2}. \label{exercises02_sol_doc06-10} (Continuity of parametrized functions)\\

\noindent Let $f(t,\omega ), \ t\in T$ a family of measurable functions indexed by $t\in I$, where $I=]a,b[, \ a<b$, to fix ideas. Let $t_{0}\in I$.\\

\bigskip \noindent Question(1). (Local continuity) Suppose that for some $\varepsilon >0,$ for $\varepsilon $ small enough
that $]t_{0}-\varepsilon ,t_{0}+\varepsilon \lbrack \subset I),$ there
exists $g$ integrable such that%
\begin{equation*}
\forall (t\in ]t_{0}-\varepsilon ,t_{0}+\varepsilon \lbrack ),\left\vert
f(t,\omega )\right\vert \leq g
\end{equation*}

\bigskip \noindent and, for almost every $\omega$, the function
\begin{equation*}
t\hookrightarrow f(t,\omega )
\end{equation*}

\bigskip \noindent is continuous at $t_{0}$. Then the function%
\begin{equation*}
F(t)=\int f(t,\omega )dm(\omega )
\end{equation*}

\bigskip \noindent is continuous at $t_{0}$ and%
\begin{equation*}
\lim_{t\rightarrow t_{0}}\int f(t,\omega )dm(\omega )=\int
\lim_{t\rightarrow t_{0}}f(t,\omega )dm(\omega )=\int f(t_{0},\omega) \ dm(\omega).
\end{equation*}

\bigskip \noindent Question(2). (Global continuity) Suppose that there exists $g$ integrable 
\begin{equation*}
\forall (t\in I),\left\vert f(t,\omega )\right\vert \leq g
\end{equation*}

\bigskip \noindent almost-everywhere and, \textit{almost-everywhere} the function
\begin{equation*}
t\hookrightarrow f(t,\omega )
\end{equation*}

\bigskip \noindent is continuous at each $t\in I$. Then the function%
\begin{equation*}
F(t)=\int f(t,\omega)dm(\omega)
\end{equation*}

\bigskip \noindent is continuous at each $t\in I$ and for any $u\in I$%
\begin{equation*}
\lim_{t\rightarrow u}\int f(t,\omega )dm(\omega )=\int \lim_{t\rightarrow
u}f(t,\omega )dm(\omega )=\int f(u,\omega )dm(\omega).
\end{equation*}

\bigskip \noindent \textbf{SOLUTIONS}.\\

\noindent Question (1). We have to establish that $F(t) \rightarrow F(t_0)$   as $t \rightarrow t_0$. It will be enough to establish that for any sequence $(h_n)_{n\geq 0}$ such that $h_n\rightarrow 0$ as $n\rightarrow +\infty$, we have $F(t_0+h_n) \rightarrow F(t_0)$ as $n\rightarrow +\infty$. But

$$
F(t_0+h_n)-F(t_0)= \int (f(t_0+h_n,\omega)-f(t_0,\omega)) \ dm =: \int Y_n(\omega) \ dm(\omega), \ n\geq 0.
$$

\noindent By the \textit{a.e.} continuity of $f$ at $t_0$, we have $Y_n\rightarrow 0$ \textit{a.e.} and the functions $Y_n$'s are bounded by $g$. Hence, we have by the DCT,
 
$$
F(t_0+h_n)-F(t_0) \rightarrow 0 \ as \ n\rightarrow +\infty.
$$

\noindent \noindent Question (2). Under the conditions, we may apply Question (1) to each $t \in I$.\\
 
\bigskip \noindent \textbf{Exercise 3}. \label{exercises03_sol_doc06-10} (Differentiability of parametrized functions)\\

\noindent Let $f(t,\omega ),t\in T$ a family of measurable functions indexed by $%
t\in I$, where $I=]a,b[,a<b$, to fix ideas. Let $t_{0}\in I$.\\

\bigskip \noindent (a) (Local differentiability) Suppose that for almost every $\omega$, the function%
\begin{equation*}
t\hookrightarrow f(t,\omega )
\end{equation*}

\bigskip \noindent is differentiable at $t_{0}$ with derivative
\begin{equation*}
\frac{df(t,\omega )}{dt}|_{t=t_{0}}=f(t_{0},\omega )
\end{equation*}

\bigskip \noindent and that for some $\varepsilon >0,$ for $\varepsilon $ small enough that $%
]t_{0}-\varepsilon ,t_{0}+\varepsilon \lbrack \subset I),$ there exists $g$ integrable such that, 
\begin{equation*}
\forall (\left\vert h\right\vert \leq \varepsilon ),\left\vert \frac{f(t+h,\omega )-f(t,\omega )}{h}\right\vert \leq g.
\end{equation*}

\bigskip \noindent Then the function
\begin{equation*}
F(t)=\int f(t,\omega )dm(\omega)
\end{equation*}

\bigskip \noindent is differentiable $t_{0}$ and%
\begin{equation*}
F^{\prime }(t_{0})\frac{d}{dt}\int f(t,\omega )dm(\omega )|_{t=t_{0}}=\int 
\frac{df(t,\omega )}{dt}|_{t=t_{0}}dm(\omega).
\end{equation*}

\bigskip \noindent (b) (Everywhere differentiability) Suppose that $f(t,\omega )$ is $a.e$ differentiable at each t$\in T$ and
there $g$ integrable such that
\begin{equation*}
\forall (t\in I),\left\vert \frac{d}{dt}f(t,\omega )\right\vert \leq g.
\end{equation*}

\bigskip \noindent Then, the function
\begin{equation*}
t\hookrightarrow f(t,\omega )
\end{equation*}

\bigskip \noindent is differentiable a each $t\in I$ and at each $t\in I,$ 
\begin{equation*}
\frac{d}{dt}F(t)=\int \frac{d}{dt}f(t,\omega )dm(\omega).
\end{equation*}

\noindent (c) Let
\begin{equation*}
U(t)=\sum_{n}u_{n}(t)
\end{equation*}

\bigskip \noindent be a convergent series of functions for $\left\vert t\right\vert <R,R>0.$ Suppose
that   each $u_{n}(t)$ is differentiable and there exists $g=(g_{n})$ integrable that is 
\begin{equation*}
\sum_{n}g_{n}<+\infty, 
\end{equation*}

\bigskip \noindent such that
\begin{equation*}
\text{For all }\left\vert t\right\vert <R,\left\vert \frac{d}{dt}
u_{n}(t)\right\vert \leq g.
\end{equation*}

\bigskip \noindent Then, the function $U(t)$ is differentiable and
\begin{equation*}
\frac{d}{dt}U(t)=\sum_{n}\frac{d}{dt}u_{n}(t).
\end{equation*}

\bigskip \noindent \textbf{SOLUTIONS}.\\

\noindent Question (1). To show that $F$ is differentiable, we have to establish the existence of the limit of

$$
\frac{F(t_0+h)-F(t_0)}{h} \ as \ h \rightarrow 0.
$$

\bigskip \noindent It will be enough to establish that there exists a unique finite number $D$ such that for any sequence $(h_n)_{n\geq 0}$ such that $h_n\rightarrow 0$ as $n\rightarrow +\infty$, 

$$
\frac{F(t_0+h_n)-F(t_0)}{h_n} \rightarrow D \ as \ n \rightarrow +\infty.
$$

\bigskip \noindent But we have

$$
\frac{F(t_0+h_n)-F(t_0)}{h_n} =\int  \frac{f(t_0+h_n,\omega)-F(t_0,\omega)}{h_n} \ dm =: \int Y_n(\omega) \ dm(\omega), \ n\geq 0.
$$

\bigskip \noindent For large values of $n$ such that $|h_n|<\varepsilon$, we have $|Y_n|\leq g$ \textit{a.e.} and by assumption $Y_n(\omega) \rightarrow \frac{d}{dt}f(t_0,\omega)$, $\omega$-\textit{a.e.}. By the DCT, we have

$$
\frac{F(t_0+h_n)-F(t_0)}{h_n} \rightarrow \int \frac{d}{dt}f(t_0,\omega) \ dm(\omega) \ as \ n \rightarrow +\infty. \ \square
$$

\bigskip \noindent Question (2). If =$f(t,\omega)$ is \textit{a.e.} differentiable at each $t\in I$, we apply the mean value theorem (MVT) at each point $t_0\in I$ :

$$
\left| \frac{f(t_0+h_n,\omega)-F(t_0,\omega)}{h_n} \right| =\left| \frac{d}{dt}f(t_0+\theta_n(\omega) h_n,\omega) \right| \leq g, \ |\theta_n(\omega)|\leq \ n\geq 0.
$$

\bigskip \noindent and get the same conclusion as in Question (1).\\

\noindent Question (3). This is Question (2) applied to the counting measure.\\

\noindent \LARGE \textbf{Doc 06-11 : Lebesgue-Stieljes and Riemann-Stieljes integrals - Exercise}. \label{doc06-11}\\
\bigskip
\Large

\bigskip \noindent The definition of the Riemann-Stieljes Integral is given in Doc 05-04 in Chapter \ref{05_integration}.\\

\noindent We recommend also to read at least the statement of \textit{Exercise 8 in Doc 03-06} \ in Chapter \ref{03_setsmes_applimes_cas_speciaux}, since we are going to use semi-continuous functions.\\

\noindent In all this text, $a$ and $b$ are two real numbers such that $a<b$ and $F: [a,b]\longrightarrow \mathbb{R}$ is
a non-constant and non-decreasing function, \textbf{or} a function of bounded variation on $[a,b]$.\\ 

\bigskip \noindent \textbf{A - Integral on compact sets}.\\

\noindent \textbf{Exercise 1}. \label{exercises01_sol_doc06-10}\\

\noindent Provide a bounded function $f$ on $[a,b]$, which is Lebesgue-Stieljes integrable with respect to any non-decreasing function $F$, and which is not Riemann-Stieljes with respect to any non-decreasing function $F$.\\

\noindent \textit{Hint}. Consider
\begin{equation*}
f=1_{]a,b]\cap \mathbb{Q}}, 
\end{equation*}

\noindent where $\mathbb{Q}$ is the set of rational numbers.\\

\noindent (a) What is the value of $\int_{]a,b]} g \ d\lambda_F$?\\

\noindent (b) Use Riemann-Stieljes sums 
\begin{equation}
S_{n}=\sum_{i=0}^{\ell(n)-1}f(c_{i,n})(F(x_{i+1,n})-F(x_{i,n})).  \label{rsums}
\end{equation}

\noindent and assign to all the $c_{i,n}$'s rational values, and next, irrational values.\\

\noindent Conclude.\\

\bigskip \noindent \textbf{SOLUTIONS}.\\

\noindent (a) $|f|$ is bounded by $g=1_{]a,b]}$ for which we have $\int f \ d\lambda_F=F(b)-F(a)$. Hence $f$ is integrable and 

\begin{eqnarray*}
\int f \ d\lambda_F&=&\lambda_F(]a,b]\cap \mathbb{Q})\\
&=&\sum_{x\in \mathbb{Q}\cap ]a,b]} F(x+0)-F(x-0)\\
&\leq& F(b)-F(a).
\end{eqnarray*}

\bigskip (b) The Riemann-Stieljes sums is zero for irrational values of the $c_{i,n}$'s and $F(b)-F(a)$ for rational values. So the Riemann-Stieljes
integral of $f$ does not exists.\\

\bigskip (c) This concludes the Exercise.\\

\bigskip \noindent \textbf{Exercise 2}. \label{exercise02_sol_doc06-10}\\

\noindent Let $F: [a,b]\longrightarrow \mathbb{R}$ be a non-constant and non-decreasing function.\\

\noindent Show that any bounded function $f$ on $[a,b]$ is Riemann-Stieljes integrable on $[a,b]$ if and only it is $\lambda_F$-\textbf{a.e.} continuous with respect to the Lebesgue measure-Stieljes $\lambda_F$.\\

\noindent \textit{\textbf{A - Hint on the part} : If a bounded function $f$ on $[a,b]$ is Riemann-Stieljes integrable on $[a,b]$, then $f$ is $\lambda_F$-\textit{a.e.} continuous with respect to the Lebesgue-Stieljes measure, and its Riemann-Stieljes integral and its Lebesgue-Stieljes integral are equal}.\\

\noindent We are going to use the functions and superior or inferiors limits defined in \textit{Exercise 8 in Doc 03-06} \ in Chapter \ref{03_setsmes_applimes_cas_speciaux} as follows. Put $\varepsilon >0$, 
\begin{equation*}
f^{\ast ,\varepsilon }(x)=\sup \{f(z);z\in \lbrack a,b]\cap ]x-\varepsilon
,x+\varepsilon \lbrack \}
\end{equation*}

\bigskip \noindent and 

\begin{equation*}
f_{\ast }^{\varepsilon }(x)=\inf \{f(z);z\in \lbrack a,b]\cap ]x-\varepsilon
,x+\varepsilon \lbrack \}.
\end{equation*}

\noindent The functions limit and superior limit functions are given by 
\begin{equation*}
f_{\ast }(x)=\lim_{\varepsilon \downarrow 0}f_{\ast }^{\varepsilon }(x)\leq
\lim_{\varepsilon \downarrow 0}f^{\ast ,\varepsilon }(x)=f^{\ast }(x).
\end{equation*}

\noindent According to the terminology introduced in \textit{Exercise 8 in Doc 03-06} \ in Chapter \ref{03_setsmes_applimes_cas_speciaux}, we proved that $f_{\ast }$ is lower semi-continuous and that $f^{\ast}$ is upper semi-continuous, and as such, are measurable.\\

\noindent Consider for each $n\geq 1$, a subdivision of $]a,b]$, that partitions  into $]a,b]$ into the $\ell
(n)$ sub-intervals

\begin{equation*}
]a,b]=\sum_{i=0}^{\ell (n)-1}]x_{i,n},x_{i+1,n}],
\end{equation*}

\bigskip \noindent with modulus
\begin{equation*}
m(\pi _{n})=\max_{0\leq i\leq m(n)}(x_{i+1,n}-x_{i+1,n}) \rightarrow 0 as n\rightarrow +\infty.
\end{equation*}

\noindent Define for each $n\geq 0$, for each $i$, $0\leq i \leq \ell(n)-1$

\begin{equation*}
m_{i,n}=\inf \{f(z),x_{i,n}\leq z<x_{i+1,n}\}\text{ and }M_{i,n}=\sup
\{f(z),x_{i,n}\leq z<x_{i+1,n}\},
\end{equation*}

\begin{equation*}
h_{n}=\sum_{i=0}^{\ell (n)-1}m_{i,n}1_{]x_{i,n},x_{i+1,n}]}\text{ and }%
H_{n}=\sum_{i=0}^{\ell (n)-1}M_{i,n}1_{]x_{i,n},x_{i+1,n}]}
\end{equation*}

\noindent and\\

\begin{equation*}
D=\bigcup_{n\geq 0}\{x_{0,n},x_{1,n},...,x_{\ell (n),n}\}.
\end{equation*}

\noindent (1)  Explain simple sentences why $f^{\ast}$ and $f_{\ast}$, $h_n$ and $H_n$, $n\geq 0$, are bounded by $M$ and why $D$ is countable. Explain why you can choose $D$ as a set of continuity points
of $F$, based on the countability of the set of discontinuity points of $F$ (See Exercise 1 in Doc DOC 03-06, page \pageref{exercise01_doc03-06}).\\

\noindent (2) Fix $x$ in  $x\in \lbrack a,b]\setminus D$, for any $n\geq 1$, there exists $i$, $0\leq i(n)\leq \ell (n)-1$\ such that $x\in
]x_{i(n),n},x_{i(n)+1,n}[$.\\

\noindent Let $\varepsilon >0$\ such that  
\begin{equation*}
]x-\varepsilon ,x+\varepsilon \lbrack \subseteq ]x_{i(n),n},x_{i(n)+1,n}[.
\end{equation*}

\noindent Show that for any $x\in \lbrack a,b]\setminus D$

\begin{eqnarray*}
h^{n}(x) \leq f_{\ast }^{\varepsilon }(x)
\end{eqnarray*}

\noindent and

\begin{eqnarray*}
H_{n}(x) \geq f_{\ast }^{\varepsilon }(x).
\end{eqnarray*}

\noindent (3) Deduce that for any $x\in \lbrack a,b]\setminus D$,

\begin{equation*}
h_{n}(x)\leq f_{\ast }(x)\leq f^{\ast }\leq H_{n}(x). \ \ (B1)
\end{equation*}

\noindent (4) Let $\eta>0$. For $n$ fixed. By using the characterization of the suprema and the infima on $\mathbb{R}$, show that the Lebesgue integral of $h_n$ and $H_n$ can be approximated by an 
Riemann-Stieljes sum in the form

\begin{equation*}
\int H_{n}d\lambda =\sum_{i=0}^{\ell
(n)-1}M_{i,n}(x_{i+1,n}-x_{i,n})<\sum_{i=0}^{\ell
(n)-1}f(d_{i,n})(x_{i+1,n}-x_{i,n})+\eta
\end{equation*}

\noindent and 

\begin{equation*}
\int h_{n}d\lambda =\sum_{i=0}^{\ell(n)-1}m_{i,n}(F(x_{i+1,n})-F(x_{i,n}))>\sum_{i=0}^{\ell
(n)-1}f(c_{i,n})(F(x_{i+1,n})-F(x_{i,n}))-\eta.
\end{equation*}

\noindent (5) Mix up all that to prove that

\begin{eqnarray*}
&& \sum_{i=0}^{\ell (n)-1}f(c_{i,n})\left( F(x_{i+1,n})-F(x_{i,n})\right) -\eta\\
&& \leq \int_{\lbrack a,b]\setminus D}f_{\ast }(x)d\lambda \\
&& \leq \int_{[a,b]\setminus D}f_{\ast }(x)d\lambda \\
&& \leq \sum_{i=0}^{\ell (n)-1}f(d_{i,n})(x_{i+1,n}-x_{i,n})+\eta.
\end{eqnarray*}

\noindent In establishing the latter double inequality, do not forget that $D$ is a $\lambda_F$-null set.\\

\noindent (6) Conclude : (i) by letting  $\eta\rightarrow 0$, and : (ii) by considering that $f$ is continuous at $x$ if and only if $f_{\ast}=f^{\ast}(x)$.\\

\bigskip 
\noindent \textit{\textbf{B - Hint on the part} : If a bounded function $f$ on $[a,b]$ is $\lambda_F$-\textit{a.e.} continuous, then it is Riemann-Stieljes integrable on $[a,b]$ exists, and its 
Riemann-Stieljes integral and its Lebesgue-Stieljes  integral are equal}.\\
 
\noindent (1) Denote by $D\subset [a,b]$ a measurable set on which $f$ is continuous such that $\lambda_F(D)=0$. Consider an arbitrary sequence of Riemann-Stieljes  sums  

$$
S_{n}=\sum_{i=0}^{\ell (n)-1}f(c_{i,n})(F(x_{i+1,n})-F(x_{i,n})),
$$

\noindent for which the sequence of modulus tends to zero with $n$. Justify the double inequality

$$
\int h_{n}d\lambda_F \leq S_n \leq \int H_{n}d\lambda_F
$$

\noindent and next, by denoting $H=[a,b]\setminus D$

$$
\int h_{n} 1_{H} d\lambda_F \leq S_n \leq \int H_{n} h_{n} 1_{H} d\lambda_F.
$$

\noindent By using Question (7) of  \textit{Exercise 8 in Doc 03-06} \ in Chapter \ref{03_setsmes_applimes_cas_speciaux}, what are the limits of $h_{n} 1_{H}$ and 
$H_{n} 1_{H}$.\\

\noindent Use the Dominated Convergence Theorem to get

$$
\int f_{\ast} 1_{H} d\lambda_F \leq  \liminf_{n\rightarrow +\infty} S_n \leq \limsup_{n\rightarrow +\infty} S_n \leq\int f*{\ast} 1_{H} d\lambda_F.
$$

\noindent Conclude.\\

\bigskip \noindent \textbf{SOLUTIONS}.

\noindent (A). Let us begin to show that if the bounded function $f$ on $[a,b]$ is Riemann-Stieljes integrable on $[a,b]$, then $f$ is almost-surely continuous with respect to the Lebesgue--Stieljes measure $\lambda_F$, and its Riemann-Stieljes integral and its Lebesgue-Stieljes integral are equal. Suppose that $f$, bounded by $M$ is Riemann-Stieljes integrable.\\

\noindent The functions $f_{\ast}$ and $f^{\ast}$ are limits of functions $f_{\ast }^{\varepsilon}$ and $f^{\ast,\varepsilon}$, $\varepsilon>0$, which are extrema of sets bounded by $M$. Thus, they are bounded by $M$. As well, the functions $h_n$ and $H_n$, $n\geq 0$, are bounded by $M$. The set $D$ is also obviously countable as a countable unions of finite sets.\\

\noindent It is clear that for any $x$ in  $x\in \lbrack a,b]\setminus D$, for any $n\geq 1$, there exists $i$, $0\leq i(n)\leq \ell (n)-1$\ such that $x\in
]x_{i(n),n},x_{i(n)+1,n}[$, and thus there exists $\varepsilon_0>0$ such that
\begin{equation*}
]x-\varepsilon, x+\varepsilon[ \subseteq ]x_{i(n),n},x_{i(n)+1,n}[.
\end{equation*}

\noindent This implies that
\begin{eqnarray*}
f^{\ast,\varepsilon_0}&=&\sup \{f(z), z\in ]x-\varepsilon_0, x+\varepsilon_0[\}\leq \sup \{f(z), z\in ]x_{i(n),n},x_{i(n)+1,n}[\}\\
&=&M_{i(n),n}=H_n(x), \ \ (B11)
\end{eqnarray*}

\noindent and

\begin{eqnarray*}
f_{\ast}^{\varepsilon_0}&=&\inf \{f(z), z\in ]x-\varepsilon_0, x+\varepsilon_0[\} \geq \inf \{f(z), z\in ]x_{i(n),n},x_{i(n)+1,n}[\}\\
&=&m_{i(n),n}=h_n(x). \ \ (B12)
\end{eqnarray*}

\bigskip \noindent Remind that the functions $f_{\ast}^{\varepsilon}$'s and $f_{\ast}^{\varepsilon}$'s are respectively non-increasing and non-decreasing in $\varepsilon$. This justifies that (B11) and (B12) still hold when we replace $\varepsilon_0$ by $\varepsilon>0$ such that $\varepsilon\leq \varepsilon_0$. $0<\varepsilon\leq \varepsilon_0$, we have

$$
h_n 1_{\lbrack a,b]\setminus D} \leq f_{\ast}^{\varepsilon} 1_{\lbrack a,b]\setminus D} \leq f^{\ast,\varepsilon} 1_{\lbrack a,b]\setminus D} \leq H_n 1_{\lbrack a,b]\setminus D}. \ \ (B13)
$$

\bigskip \noindent We obtain Formula (B1) by letting $\varepsilon \downarrow 0$.\\

\noindent Since the set of discontinuity points of $F$  is countable (See Exercise 1 in Doc DOC 03-06, page \pageref{exercise01_doc03-06}), we may choose all points of $D$ as continuity points of $F$. This implies that (See Exercise 2 in Doc Doc 04-04, page \pageref{exercise02_doc04-04})

$$
\lambda_F(D)=\sum_{x \in D} F(x)-F(x-)=0.
$$

\noindent  Next, let fix $\eta>0$. Since, each $m_{i,n}$ is the infimum of $\sup \{f(z), z\in ]x_{i(n),n},x_{i(n)+1,n}[\}$ which is finite, there exists 
$c_{i,n} \in ]x_{i(n),n},x_{i(n)+1,n}[$, such that $m_{i,n} \leq f(c_{i,n}) <m_{i,n}-\eta/(F(b)-F(a))$. Hence

\begin{eqnarray*}
\int h_{n}d\lambda_F &=&\sum_{i=0}^{\ell(n)-1}m_{i,n}(F(x_{i+1,n})-F(x_{i,n}))\\
&>&\sum_{i=0}^{\ell(n)-1}(f(c_{i,n})-\eta/(F(b)-F(a)))(F(x_{i+1,n})-F(x_{i,n}))\\
&=&\sum_{i=0}^{\ell(n)-1}f(c_{i,n})(F(x_{i+1,n})-F(x_{i,n}))\\
&-&\frac{\eta}{F(b)-F(a)}\sum_{i=0}^{\ell(n)-1}(F(x_{i+1,n})-F(x_{i,n}))\\
&=&\sum_{i=0}^{\ell(n)-1}f(c_{i,n})(F(x_{i+1,n})-F(x_{i,n}))-\eta.
\end{eqnarray*} 

\noindent As well, since  each $m_{i,n}$ is the supremum of $\sup \{f(z), z\in ]x_{i(n),n},x_{i(n)+1,n}[\}$ which is finite, there exists 
$d_{i,n} \in ]x_{i(n),n},x_{i(n)+1,n}[$, such that $M_{i,n}-\eta/(F(b)-F(a)) < f(d_{i,n}) \leq M$. Hence

\begin{eqnarray*}
\int H_{n}d\lambda_F &=&\sum_{i=0}^{\ell(n)-1}M_{i,n}(F(x_{i+1,n})-F(x_{i,n}))\\
&<&\sum_{i=0}^{\ell(n)-1}(f(d_{i,n})+\eta/(F(b)-F(a)))(F(x_{i+1,n})-F(x_{i,n}))\\
&=&\sum_{i=0}^{\ell(n)-1}f(d_{i,n})(F(x_{i+1,n})-F(x_{i,n}))\\
&+&\frac{\eta}{F(b)-F(a)}\sum_{i=0}^{\ell(n)-1}(F(x_{i+1,n})-F(x_{i,n}))\\
&=&\sum_{i=0}^{\ell(n)-1}f(c_{i,n})(F(x_{i+1,n})-F(x_{i,n}))+\eta.
\end{eqnarray*} 

\noindent Hence, for one hand, we have

\begin{eqnarray*}
&&\sum_{i=0}^{\ell(n)-1}f(c_{i,n})(F(x_{i+1,n})-F(x_{i,n}))-\eta \leq \int h_{n}d\lambda_F\\
&=&\int_{[a,b]\setminus D} h_n \ d\lambda_F \ \ (Since \ D \ is \ a \ \lambda-null \ set)\\
&\leq&\int_{[a,b]\setminus D} f_{\ast}  \ \ d\lambda_F \ (By \ Formula \ (B1))\\
\end{eqnarray*} 

\noindent and, based on the same arguments, we have for another hand, 

\begin{eqnarray*}
\int_{[a,b]\setminus D} H_n \ d\lambda_F &=& \int_{[a,b]\setminus D} H_n \ d\lambda_F\\
&\leq& \int_{[a,b]\setminus D} f^{\ast} \ d\lambda_F\\
&\leq&\sum_{i=0}^{\ell(n)-1}f(d_{i,n})(F(x_{i+1,n})-F(x_{i,n}))-\eta.
\end{eqnarray*} 

\noindent If $f$ is Riemann-Stieljes integrable, the Riemann-Stieljes sums converges to $\int_{a}^{b} f(x) \ dF(x)$ when the sequences of modilii converges to zero. Wet get by letting $n\rightarrow +\infty$ and next $\eta\rightarrow 0$,

$$
\int f(x) dF(x)=\int_{[a,b]\setminus D} f^{\ast} \ d\lambda_F=\int_{[a,b]\setminus D} f^{\ast} \ d\lambda_F. \ \ (B14)
$$

\noindent Since $f_{\ast} \leq f^{\ast}$, the second equality means that

$$
\int 1_H \left| f^{\ast}-f_{\ast} \right| \ d\lambda_F =0,
$$

\noindent where $H=[a,b]\setminus D$. By Property (O2) of Point (05.05) in Doc 05-01 (page \pageref{doc05-01}), we conclude that $1_H f^{\ast}= 1_H f_{\ast}$ $\lambda_F$-\textit{a.e.}, and hence $f^{\ast}=f_{\ast}$ since $H^c$ is a $\lambda_F$-null set on $[a,b]$. But on the set $(f^{\ast}=f_{\ast})$, $f=f^{\ast}=f_{\ast}$ is continuous. Thus $f$ is continuous $\lambda_F$-\textit{a.e.} and becomes (B14)

$$
\int f(x) \ dF(x) = \int f \ d\lambda_F.
$$ 

\bigskip \bigskip 
\noindent (B). Let us show now that if a bounded function $f$ on $[a,b]$ is $\lambda_F$-\textit{a.e.} continuous, then it is Riemann integrable on $[a,b]$ exists, and its Riemann-Stieljes integral and its Lebesgue-Stieljes integral are equal. Suppose that $f$ is bounded by $M$ and is $[a,b]$ is $\lambda_F$-\textit{a.e.} continuous and denote $H=(f \text{continuous}) \subset [a,b]$.\\

\noindent Let us consider an arbitrary sequence of Riemann sums 
$$
S_{n}=\sum_{i=0}^{\ell (n)-1}f(c_{i,n})(F(x_{i+1,n})-F(x_{i,n})), 
$$

\noindent for which the sequence of modulus tends to zero with $n$. By the definition of the $h_n$'s and the $H_n$'s and the expressions of the integrals of elementary functions, we have

$$
\int h_{n}d\lambda_F \leq S_n \leq \int H_{n}d\lambda_F
$$

\noindent and next, since $\lambda_F(H^c)=0$,

$$
\int h_{n} 1_{H} d\lambda_F \leq S_n \leq \int H_{n} h_{n} 1_{H} d\lambda_F.
$$

\noindent By Question (7) of  \textit{Exercise 8 in Doc 03-06} \ in Chapter \ref{03_setsmes_applimes_cas_speciaux}, see \pageref{exercise07_sol_doc03-09},  $h_{n} 1_{H}$ converges to $f_{\ast}1_H$ and
$H_{n} 1_{H}$ to $f_{\ast}1_H$. Since the functions $h_{n} 1_{H}$ and $M_{n} 1_{H}$ are bounded by $g=M1_{[a,b]}$ which is $\lambda_F$-integrable on $[a,b]$, we may and do apply the DCT to get\\

$$
\int f_{\ast} 1_{H} d\lambda_F \leq  \liminf_{n\rightarrow +\infty} S_n \leq \limsup_{n\rightarrow +\infty} S_n \leq\int f*{\ast} 1_{H} d\lambda_F.
$$

\noindent By Hypothesis, $f_{\ast} 1_{H}=f^{\ast} 1_{H}=f 1_H$. Thus we have 

$$
\int f 1_{H} d\lambda_F =  \liminf_{n\rightarrow +\infty} S_n = \limsup_{n\rightarrow +\infty} S_n.
$$

\noindent Any arbitrary sequence of Riemann sums for which the sequence of modulus tends to zero with $n$, converges to the finite  number $\int f 1_{H} d\lambda_F$. So, this number is the Rieamann
integral of $f$ on $[a,b]$ and, since since is co-null set, we have

$$
\int f(x) \ dF(x)=\int f d\lambda_F.
$$

\bigskip \noindent \textbf{Exercise 3}. \label{exercise03_sol_doc06-10}\\

\noindent On the whole extended real line $\mathbb{R}$, the Lebesgue integral of a measurable application \textit{a.e.} $f$ is defined on  if and only if

$$
\int_{\mathbb{R}} f^{-} d\lambda<+\infty \text{ and }  \int_{\mathbb{R}} f^{+} d\lambda<+\infty
$$

\noindent and the integral of $f$ 

$$
\int_{\mathbb{R}} f d\lambda=\int_{\mathbb{R}} f^{+} d\lambda - \int_{\mathbb{R}} f^{-} d\lambda<+\infty
$$

\noindent is finite if and only if the integral of $|f|$ 

$$
\int_{\mathbb{R}} |f| d\lambda=\int_{\mathbb{R}} f^{+} d\lambda + \int_{\mathbb{R}} f^{-} d\lambda<+\infty
$$

\noindent is.\\

\bigskip \noindent (b) The Riemann integral is not defined on the whole real line. Instead, we have improper integrals obtained by limits.

\noindent The improper integral of $f$ on $\mathbb{R}$ 

\begin{equation}
\int_{-\infty }^{+\infty }f(x)\text{ }dF(x)  \label{RS-LB4}
\end{equation}

\noindent is defined as the limit of the integrals of 
 
\begin{equation*}
\int_{a}^{b}f(x)\text{ }dF(x)
\end{equation*}

\noindent as $a\rightarrow -\infty $ and $b\rightarrow +\infty$, independently of $a\rightarrow -\infty $ and $b\rightarrow +\infty$.

\bigskip \noindent Question (a) Show that if $f$ is $\lambda$-integrable on $\mathbb{R}$, locally bounded and $\lambda$-\textit{a.e.} continuous, then the improper integral of $f$ exists, is finite and is equal to the Lebesgue integral, that is 

$$
\int_{-\infty}^{+\infty} f \ dF(x) \neq \int f \ d\lambda. \ \ (B1)
$$

\bigskip \noindent Question (b) Show that if the improper integral of $|f|$ exists and is finite, then $f$ is Lebesgue integrable and (B1) holds.

\bigskip \noindent Question (c). Let $f=1_{\mathbb{Q}}$. Show that the improper integral of $f$ is not defined and that the Lebesgue integral of $f$ exists and has the null value.

\bigskip \noindent Question (d) Consider here the particular case $F(x)=x$ and $f(x)=(\sin x)/x$, $x\in \mathbb{R}$. Show that the Lebesgue integral of $f$ does not exists, the Lebesgue integral of $|f|$ is infinite, the improper integral of $f$ exists and is infinite and the improper integral of $f$ exists and is finite and cannot be a Lebesgue integral.\\

\noindent Follow this way : \\

\noindent (i) Remark that $f$ is locally continuous and the locally bounded. Use Exercise 2 to get for all $a\leq 0$ and $b\geq 0$, 

$$
\int_{a}^{b} |f(x)| dx)=\int_{[a,b]} |f| \ d\lambda. \ \ (B18)
$$

\bigskip \noindent Let $a\downarrow +\infty$ and $b \uparrow +\infty$ and conclude that 

$$
\int_{-\infty}^{+\infty} |f(x)| dx=\int_{\mathbb{R}} |f(x)| \ d\lambda=+\infty, \ \ (B19)
$$

\bigskip \noindent (ii) Define 

$$
A=\{k\pi, k\in \mathbb{Z}\}, \ B=\sum_{k\ in \mathbb{Z}} ]2k\pi, (2k+1)\pi[ \ and \ B=\sum_{k\ in \mathbb{Z}} ](2k+1)\pi, 2(k+1)\pi[.
$$

\bigskip \noindent Establish that

$$
f=1 \ on \ \{0\}, f=0 \ on \ A\setminus \{0\} , f>0 \ on \ B, \ f<0 \ on \ C.
$$

\bigskip \noindent and

$$
f^+=f 1_{A+C}, \  f^-=f 1_{A+C} \ and \ \lambda(A)=0.
$$

\bigskip \noindent Conclude that $f^+$ and $f^-$ are continuous outside the null set $A$ and are  are locally \textbf{a.e.} continuous and they are also locally bounded. Repeat the argiments used in  Formulas (B18) and (B19) to 

$$
\int_{-\infty}^{+\infty} f^-(x) dx=\int_{\mathbb{R}} f^-(x) \ d\lambda \ \ and \ \ \int_{-\infty}^{+\infty} f^+(x) dx=\int_{\mathbb{R}} f^+(x) \ d\lambda. \ \ (B20)
$$

\bigskip \noindent (ii) Justify

$$
\int_{-\infty}^{+\infty} f^-(x) dx= \sum_{k\ in \mathbb{Z}} \int_{(2k+1)\pi}^{2(k+1)\pi} \frac{-\sin x}{x} \ dx
$$

\bigskip \noindent and

$$
\int_{-\infty}^{+\infty} f^+(x) dx= \sum_{k\ in \mathbb{Z}} \int_{2k\pi}^{(2k+1)\pi} \frac{\sin x}{x} \ dx
$$

\bigskip \noindent and compare the series in the right-hand members and prove the equality of the right-hand members to get

$$
\int_{-\infty}^{+\infty} |f|(x) dx=2 \int_{-\infty}^{+\infty} f^-(x) dx=\int_{-\infty}^{+\infty} f^+(x) dx=+\infty.
$$

\bigskip \noindent Now conclude.\\

\bigskip \noindent \textbf{SOLUTIONS}.\\

\noindent Question (a). Suppose that $f$ is $\lambda$-integrable on $\mathbb{R}$, locally bounded and $\lambda_F$-\textit{a.e.} continuous. Thus for each $n\geq 0$, 

$$
\left| \int_{[-n,n]} f \ d\lambda_F \right| \leq \int_{[-n,n]} |f| \ d\lambda_F \leq \int |f| \ d\lambda_F=M<+\infty. \ \ (B16)
$$

\bigskip \noindent Since $f$ is locally bounded and $\lambda_F$-\textit{a.e.} continuous, we have by Exercise 1 that for each couple $(a,b)$, $a\leq b$,

$$
\int_{a}^{b} f(x) dF(x) = \int_{[a,b]} f \ d\lambda_F\leq M<+\infty.
$$

\bigskip \noindent Next, by letting $a\rightarrow -\infty$ and $a\rightarrow +\infty$, the right-hand member converges to the improper integral of $f$ while the left-hand member, by the DCT, converges to the Lebesgue integral of $f$.\\

\bigskip \noindent Question (b) Suppose that the improper integral of $|f|$ exists and is finite, and that $f$ is locally bounded and integrable. Hence, for all $a\leq 0$ and $b\geq 0$, 

$$
\int_{a}^{b} |f(x)| dF(x) \leq \int_{-\infty}^{+\infty} |f(x)| \ dF(x)=:M<+\infty. \ \ (B17a)
$$

\bigskip \noindent Thus, by Exercise 2, we have that for all $a\leq 0$ and $b\geq 0$,

$$
\int_{[a,b]} |f| \ d\lambda_F=\int_{a}^{b} |f(x)| dF(x)<+\infty. \ \ (B17b)
$$

\bigskip \noindent Hence $\int_{[a,b]} f \ d\lambda_F$ is defined. By applying the MCT on (B17b), we get as $a\downarrow +\infty$ and $b \uparrow +\infty $ (we eventually may take the discretized limits),

$$
\int_{\mathbb{R}} |f| \ d\lambda_F\leq M<+\infty.
$$

\bigskip \noindent Now, $f$ is Riemann-Stieljes locally integrable, meaning that for all $-\infty <a<0<b<+\infty$, the integral $\int_{a}^{b} f(x) dF(x)$ exists, we may and do apply based on the following set of facts  

$$
\int_{a}^{b} f(x) dF(x)=\int_{[a,b]} f \ d\lambda_F, \ \ \left|1_{[a,b]} f\right|\leq |f|, \ \ \int |f| \ d\lambda_F<+\infty,
$$ 

\bigskip \noindent to get, as $a\downarrow +\infty$ and $b \uparrow +\infty$, that

$$
\int_{-\infty}^{+\infty} f(x) dx=\int_{\mathbb{R}} f \ d\lambda=M<+\infty.
$$ 

\bigskip \noindent Question (c). Let $f=1_{\mathbb{Q}}$. We know by Exercise 1 that the Riemann-Stieljes integral of $f$ does not exists on $[a,b]$, whatever be $-\infty <a<0<b\leq <+\infty$. Hence, its improper integral does not exists. But $f$ the Lebesgue-Stieljes integral of $f$ exists and has the value

$$
\int f \ d\lambda=\lambda(\mathbb{Q})=\sum_{x\in \mathbb{Q}} F(x+0)-F(x-0).
$$

\bigskip \noindent Question (d). Let $f(x)=(\sin x)/x$, $x\in \mathbb{R}$. we know in Classical Analysis that the proper integral of is absolute value is infinite, that is

$$
\int \left| \frac{\sin x}{x} \right| dx=+\infty.
$$

\bigskip \noindent Surely $f$ is locally continuous and the locally bounded.  Hence, for all $a\leq 0$ and $b\geq 0$, 

$$
\int_{a}^{b} |f(x)| dx=\int_{[a,b]} |f| \ d\lambda. \ \ (B18)
$$

\bigskip \noindent By applying the definition of the improper integral at left and the MCT at right, we get as $a\downarrow +\infty$ and $b \uparrow +\infty$

$$
\int_{-\infty}^{+\infty} |f(x)| dx=\int_{\mathbb{R}} |f(x)| \ d\lambda=+\infty, \ \ (B19)
$$

\bigskip \noindent Now let 

$$
A=\{k\pi, k\in \mathbb{Z}\}, \ B=\sum_{k\ in \mathbb{Z}} ]2k\pi, (2k+1)\pi[ \ and \ C=\sum_{k\ in \mathbb{Z}} ](2k+1)\pi, 2(k+1)\pi[.
$$

\bigskip \noindent We surely have

$$
f=1 \ on \ \{0\}, f=0 \ on \ A\setminus \{0\} , f>0 \ on \ B, \ f<0 \ on \ C.
$$

\noindent and

$$
f^+=f 1_{A+C}, \  f^-=f 1_{A+C} \ and \ \lambda(A)=0.
$$

\bigskip \noindent and hence $f^+$ and $f^-$ are continuous outside the null set $A$ (Actually, they are continuous on each interval involved in the formation of $B$ and $C$). So they are locally \textbf{a.e.} continuous and they are also locally bounded since they are bounded by $f$ which is locally uniformly continuous. So we may repeat analogs of Formulas (B18) and (B19) for them to get

$$
\int_{-\infty}^{+\infty} f^-(x) dx=\int_{\mathbb{R}} f^-(x) \ d\lambda \ \ and \ \ \int_{-\infty}^{+\infty} f^+(x) dx=\int_{\mathbb{R}} f^+(x) \ d\lambda. \ \ (B20)
$$

\bigskip \noindent We justify

$$
\int_{-\infty}^{+\infty} f^-(x) dx= \sum_{k\ in \mathbb{Z}} \int_{(2k+1)\pi}^{2(k+1)\pi} \frac{-\sin x}{x} \ dx
$$

\bigskip \noindent and

$$
\int_{-\infty}^{+\infty} f^+(x) dx= \sum_{k\ in \mathbb{Z}} \int_{2k\pi}^{(2k+1)\pi} \frac{\sin x}{x} \ dx
$$

\bigskip \noindent and see that the two series in the right-hand members are equal term by term and get

$$
\int_{-\infty}^{+\infty} |f|(x) dx=2 \int_{-\infty}^{+\infty} f^-(x) dx=\int_{-\infty}^{+\infty} f^+(x) dx=+\infty,
$$

\bigskip \noindent which leads to

$$
\int_{\mathbb{R}} f^-(x) \ d\lambda=\int_{\mathbb{R}} f^+(x) \ d\lambda=+\infty.
$$

\bigskip \noindent We conclude that the Lebesgue integral of $f$ does not exists, the Lebesgue integral of $|f|$ is infinite, the improper integral of $|f|$ exists and is infinite and the improper integral of $f$ exists and is finite and cannot be a Lebesgue integral.\\

%\chapter{Finite Product Measures}
%\newpage
\chapter{Finite Product Measure Space and Fubini's Theorem} \label{07_finiteProductMeasure}

\noindent \textbf{Content of the Chapter}

\begin{table}[htbp]
	\centering
		\begin{tabular}{llll}
		\hline
		Type& Name & Title  & page\\
		S& 	Doc 07-01 & Finite Product Measures and Fubini's  Theorem - A summary  & \pageref{doc07-01}\\
		D& 	Doc 07-02 & Finite Product Measures - Exercises & \pageref{doc07-02} \\
		SD& Doc 07-03  & Finite Product Measures - Exercises with Solutions & \pageref{doc07-03} \\
		%&  &   & \\
		%&  &   & \\
		%&  &   & \\
		\hline
		\end{tabular}
\end{table}

\newpage
\noindent \LARGE \textbf{Doc 07-01 : Finite Product Measures and Fubini's  Theorem - A summary}. \label{doc07-01}\\
\bigskip
\Large 

\bigskip \noindent \textbf{Important Advice}.  The learner is strongly advised to go back \textit{Doc 00-01, Doc 00-02 (Exercise 7), in Chapter \ref{00_sets}}, pages \pageref{exercise07_doc00-02} and \pageref{exercise07_sol_doc00-03}, for a reminder on Cartesian products on sets, and to \textbf{Point (01.08) in Doc 01-01 in Chapter \ref{01_setsmes}}, page \pageref{doc01-01} for revisiting the product sigma-algebra. As well, the properties of sections of subsets in a product space  and their properties are important to know before you proceed (\textit{See Exercise 3, Doc 01-06 in in Chapter \ref{01_setsmes}}, page \pageref{exercise03_doc01-06}).\\

\noindent \textbf{Remark 1} : Here, the important product measure will come out here as an application of the Monotone Convergence Theorem (MCT).\\

\noindent \textbf{Remark 2} : This document will focus on the simplest product $k=2$. The extension to the general case is a matter of induction.\newline

\noindent \textbf{Remark 3} : He we focus on the finite product measure of $\sigma$-finite measures. But after we treated all fundamental properties in details, we will open express the main results on a more general condition that $\sigma$-finiteness in Exercise 5 in Doc 07-03.\\

\bigskip \noindent Let be $k$ $\sigma$-finite measure spaces $(\Omega _{i},A_{i},m_{i})$, $i=1,...,k$. Consider the product space 
\begin{equation*}
\Omega =\prod_{1\leq i\leq k}\Omega _{i},
\end{equation*}

\bigskip \noindent endowed with the product $\sigma -$algebra 
\begin{equation*}
\mathcal{A}=\bigotimes_{1\leq i\leq k}\mathcal{A}_{i}.
\end{equation*}

\bigskip \noindent that is generated by the semi-algebra $\mathcal{S}$ of measurable rectangles or cylinders 
\begin{equation*}
\mathcal{S}=\{\prod_{1\leq i\leq k}A_{i},\text{ }A_{i}\in \mathcal{A}_{i}\}.
\end{equation*}

\bigskip \noindent In the first part, we focus on the simplest case where $k=2$.\\

\textbf{Part A - Product of two measure spaces}.

\bigskip \noindent Recall the definition of the sections of a subset $A$ of $\Omega $ by : for
$\omega _{1}\in \Omega _{1}$ fixed, the section of $A$ at  $\omega_{1}$ is defined by
\begin{equation*}
A_{\omega _{1}}=\{\omega _{2}\in \Omega _{2},(\omega _{1},\omega _{2})\in A\}
\end{equation*}

\noindent and for $\omega_{2}\in \Omega_{2}$ fixed, the section of $A$ at $\omega_{2}$ is defined by 
\begin{equation*}
A_{\omega _{2}}=\{\omega _{1}\in \Omega _{1},(\omega _{1},\omega _{2})\in A\}.
\end{equation*}

\bigskip \noindent \textbf{(07.01) Definition and construction of a product measure}.\\

\noindent Here, \textbf{one} product measure is constructed through the following steps.\\

\noindent \textbf{Step 1}.\\

\noindent \textbf{(a)}. For any $A\in \mathcal{A}$ and for any $\omega _{1}\in \Omega ,$ the function 
\begin{equation}
\omega _{2}\mapsto m_{2}(A_{\omega _{1}})
\end{equation}

\noindent is $\mathcal{A}_1$-measurable.\\

\noindent \textbf{(b)}. For any $\mathcal{A}_1 \otimes \mathcal{A}_2$-measurable function

\begin{equation*}
\begin{tabular}{lllll}
$f$ & $:$ & $\Omega =(\Omega _{1}\times \Omega _{2})$ & $\longmapsto $ & $%
\overline{\mathbb{R}}$
\end{tabular}
\end{equation*}

\bigskip \noindent the partial function, for $\omega _{1}\in \Omega ,$%

\begin{equation*}
\omega_{2}\longmapsto f(\omega _{1},\omega _{2})
\end{equation*}

\noindent is $\mathcal{A}_2$ measurable.\\

\noindent \textbf{(c)}. For any $\mathcal{A}_1 \otimes \mathcal{A}_2$-non-negative measurable application $f$, for any fixed $\omega _{1}\in \Omega$ the application%
\begin{equation*}
\omega _{1}\longmapsto \int_{\Omega _{2}}f(\omega _{1},\omega _{2})dm(\omega_{2})
\end{equation*}

\bigskip \noindent is $\mathcal{A}_1$-measurable.\newline

\bigskip \noindent \textbf{Step 2}. Based on the previous points, the application $m$ defined on $\mathcal{A}$ by%
\begin{equation}
m(A)=\int_{\Omega _{2}}m_{1}(A_{\omega _{2}})\text{ }dm(\omega _{2})
\label{mesprod}
\end{equation}

\noindent is well-defined and satisfies

\noindent \hskip 5cm  $m$ is a measure.\\

\noindent \hskip 5cm  $m$ satisfies the following formula on the class of measurable rectangles :\\

\begin{equation*}
m(A_{1}\times A_{2})=m_{1}(A_{1})\times m_{2}(A_{2})  \ \ (MESPROD)
\end{equation*}

\noindent for any $A=A_{1}\times A_{2}$, with $A_{i} \in \mathcal{A}_{i}$, $i=1,2$.

\bigskip \noindent \textbf{Step 3}. Since $m_{1}$ and $m_{2}$ are $\sigma $-finite, $m$ is $\sigma $-finite and is the unique measure on $(\Omega_1 \times \Omega_2, \mathcal{A}_{1} \times \mathcal{A}_{2}$ such that (MESPROD) is valid. The unique measure is called \textbf{the} product measure

$$
m_1 \otimes m_2.
$$

\bigskip \textbf{Part B - Fubini and Tonelli's Theorems the finite product measures}.

\bigskip \noindent \textbf{(07.02)}. Fubini's Theorem.\\

\bigskip \noindent \textbf{(07.02a)}. Non-negative functions. If 

\begin{equation*}
\begin{tabular}{lllll}
$f$ & $:$ & $\Omega =(\Omega _{1}\times \Omega _{2})$ & $\longmapsto $ & $%
\overline{\mathbb{R}}$ \\ 
&  & $(\omega _{1},\omega _{2})$ & $\hookrightarrow $ & $f(\omega
_{1},\omega _{2})$
\end{tabular}
\end{equation*}

\bigskip \noindent is a non-negative and measurable function $f$, we have  
\begin{equation*}
\int_{\Omega }f\text{ }dm=\int_{\Omega _{1}}\text{ \ }\left\{ \int_{\Omega
_{2}}f(\omega _{1},\omega _{2})\text{ }dm_{2}(\omega _{2})\right\} \text{ \ }%
dm_{1}(\omega _{1}). \ \ (FU1)
\end{equation*}

\bigskip \noindent \textbf{(07.02b)}. General case. If 
\begin{equation*}
\begin{tabular}{lllll}
$f$ & $:$ & $\Omega =(\Omega _{1}\times \Omega _{2})$ & $\longmapsto $ & $%
\overline{\mathbb{R}}$ \\ 
&  & $(\omega _{1},\omega _{2})$ & $\hookrightarrow $ & $f(\omega
_{1},\omega _{2})$
\end{tabular}
\end{equation*}

\bigskip \noindent is measurable and integrable. Then, we have

\begin{equation*}
\int_{\Omega }f\text{ }dm=\int_{\Omega _{1}}\text{ \ }\left\{ \int_{\Omega
_{2}}f(\omega _{1},\omega _{2})\text{ }dm_{2}(\omega _{2})\right\} \text{ \ }%
dm_{1}(\omega _{1}). \ \ (FU2)
\end{equation*}

\bigskip \noindent \textbf{(07.03)}. Tonelli's Theorem.\\

\noindent In the two forms of Fubini's Theorem, we have

\begin{equation*}
\int_{\Omega _{1}}\text{ \ }\left\{ \int_{\Omega
_{2}}f(\omega _{1},\omega _{2})\text{ }dm_{2}(\omega _{2})\right\} \text{ \ }%
dm_{1}(\omega _{1}). \ \  (FU3) 
\end{equation*}
 
\noindent An interesting question is : does Equality (FU3) with finite members implies that $f$ is integrable. We have the following facts. One of the member of Equality (FU3) may exists and not the other. Both may exist, be finite and equal and at the same time, $f$ is not integrable.\\

\noindent Fortunately we have the Tonelli Theorem concerning the non-negative functions. We have to remind that the (FU) holds for non-negative functions.\\

\noindent If $f$ is non-negative and one of the members of (FU3) is finite, then $f$ is integrable.\\

\textbf{Part C - General finite product measure}.\\

\noindent \textbf{(07.04)}. In the general case where $k$ $\sigma$-finite measure spaces $(\Omega_{i},A_{i},m_{i})$, $i=1,...,k$, and we have the product space  
\begin{equation*}
\Omega =\prod_{1\leq i\leq k}\Omega _{i}
\end{equation*}

\bigskip \noindent which is endowed with the product $\sigma -$algebra 
\begin{equation*}
\mathcal{A}=\bigotimes_{1\leq i\leq k}\mathcal{A}_{i}.
\end{equation*}

\bigskip \noindent and this latter is generated by the semi-algebra $\mathcal{S}$ of measurable cylinders
\begin{equation*}
\mathcal{S}=\{\prod_{1\leq i\leq k}A_{i},\text{ }A_{i}\in \mathcal{A}_{i}\}.
\end{equation*}

\bigskip \noindent We get the same conclusion we can easily get by induction. Since each measure $m_{i}$ is sigma-finite, there exists a unique measure $m$ on $\mathcal{A}$ such that for any element of $\mathcal{S}$
\begin{equation*}
\prod_{1\leq i\leq k}A_{i},\text{ }A_{i}\in \mathcal{A}_{i},
\end{equation*}%

\bigskip \noindent we have
\begin{equation*}
m(\prod_{1\leq i\leq k}A_{i})=\prod_{1\leq i\leq k}m_{i}(A_{i}). \ \ (MESPRODG)
\end{equation*}

\bigskip \noindent This unique measure is then $\sigma$-finite and is called the product measure :

$$
\bigotimes_{i=1}^{k} m_i.
$$

\bigskip \noindent \textbf{(07.04)}. Integration with respect to this measure use the Fubini's Theorem.\\

\noindent If $f$ is non-negative or $f$ is integrable, then we may proceed progressively integrate the partial functions as follows :

\begin{eqnarray*}
&&\int_{\Omega _{1}\times \Omega _{1}\times ...\times \Omega _{k-1}\times
\Omega _{k}}\text{ \ }f(\omega _{1},\omega _{2},...,\omega _{k-1},\omega_{k})\text{ }dm(\omega _{1},\omega _{2},...,\omega _{k-1},\omega _{k})\\
&=&\text{ \ }\int_{\Omega _{1}}dm_{1}(\omega _{1})\int_{\Omega_{2}}dm_{2}(\omega _{2})......\int_{\Omega _{k-1}}dm_{k-1}(\omega_{k-1})...\\
&...& \int_{\Omega _{k}}\text{ \ }f(\omega _{1},\omega _{2},...,\omega_{k-1},\omega _{k}) \ dm_{k}(\omega _{k}).
\end{eqnarray*}

\bigskip \noindent How to read Fubini's Theorem. First, fix all the variables $(\omega_{1},\omega _{2},...,\omega _{k-1})$ and consider the  function
$$
f(\omega _{1},\omega _{2},...,\omega _{k-1},\omega_{k})
$$ 

\noindent as depending only $\omega _{k}$ and integrate it with respect to $dm_{k}(\omega_{k})$. And the function

$$
\int_{\Omega _{k}}\text{ \ }f(\omega _{1},\omega _{2},...,\omega_{k-1},\omega _{k})\text{ }dm(\omega _{1},\omega _{2},...,\omega_{k-1},\omega _{k})dm_{k}(\omega _{k}),
$$

\noindent is a function of $(\omega _{1},\omega _{2},...,\omega _{k-1})$. We repeat the same procedure, by fixing $(\omega _{1},\omega _{2},...,\omega
_{k-2})$ and by integrating with respect to $dm_{k-1}(\omega _{k-1})$ so that the result is a function of $(\omega _{1},\omega _{2},...,\omega
_{k-3})$. Finally, after $k-1$ steps downward, we arrive at a function of
$\omega _{1}$ and we integrate it with respect to $dm_{1}(\omega _{1}).$\\

\bigskip \noindent Fubini's theorem allows downward sequential or successive int\'{e}gration.\\

\bigskip \noindent \textbf{Be careful}. If the whole domain of integration is a rectangle, that is for the form
\begin{equation*}
D=D_{1}\times D_{1}\times ...\times D_{k-1}\times D_{k},
\end{equation*}

\bigskip \noindent and we compute this integral
\begin{equation*}
\int_{D_{1}\times D_{1}\times ...\times D_{k-1}\times D_{k}}\text{ \ }f(\omega _{1},\omega _{2},...,\omega _{k-1},\omega _{k})\text{ }dm(\omega_{1},\omega _{2},...,\omega _{k-1},\omega _{k}),
\end{equation*}

\bigskip \noindent then, we still have Fubini's theorem by replacing 
$$
\Omega _{1}\times \Omega _{1}\times ...\times \Omega _{k-1}\times \Omega _{k}
$$ 

\noindent by 
$$
D_{1}\times D_{1}\times ...\times D_{k-1}\times D_{k}.
$$

\bigskip \noindent  But if the integration domain $D$ is not a rectangle with independent factors, we have be cautious. In the first step, we integrate on the section

\begin{equation*}
D(\omega _{1},\omega _{2},...,\omega _{k-1})=\{\omega _{k}\in \Omega_{k},(\omega _{1},\omega _{2},...,\omega _{k-1},\omega _{k})\in D\}.
\end{equation*}

\bigskip \noindent  Once this is done, the result is a function of  $(\omega _{1},\omega_{2},...,\omega _{k-1})$ defined on the projection of $D$ on $\Omega_{1}\times \Omega _{1}\times ...\times \Omega _{k-1}$ that is

\begin{equation*}
\pi _{k-1}(D)=\{(\omega _{1},\omega _{2},...,\omega _{k-1})\text{, }\exists
\omega _{k}\in \Omega _{k},(\omega _{1},\omega _{2},...,\omega _{k-1},\omega
_{k})\in D\}.
\end{equation*}

\bigskip \noindent And we proceed similarly downwards. Identifying the projections and the sections depends on particular cases and require experience. You certainly did a lot of that in undergraduate Calculus courses of on multiple integration. 

\bigskip \noindent \textbf{(07.05)}. It is important to notice that the products of $\sigma$-algebras and of measures (for $\sigma$-finite factors) are associative: \\

$$
\left(\mathcal{A}_1 \bigotimes \mathcal{A}_2 \right) \bigotimes \mathcal{A}_2 = \mathcal{A}_1 \bigotimes \left(\mathcal{A}_2 \bigotimes \mathcal{A}_3 \right)
$$

\bigskip \noindent and

$$
\left(m_1 \bigotimes m_2 \right) \bigotimes m_3 = 3_1 \bigotimes \left(m_2 \bigotimes m_3 \right).
$$

\noindent \LARGE \textbf{Doc 07-02 : Finite Product Measures - Exercises}. \label{doc07-02}\\
\bigskip
\Large

\bigskip \noindent In all this document, we work with $k$ $\sigma$-finite measure spaces $(\Omega _{i},A_{i},m_{i})$, $i=1,...,k$. The product space 
\begin{equation*}
\Omega =\prod_{1\leq i\leq k}\Omega _{i}
\end{equation*}

\bigskip \noindent endowed with the product $\sigma -$algebra 
\begin{equation*}
\mathcal{A}=\bigotimes_{1\leq i\leq k}\mathcal{A}_{i}.
\end{equation*}

\bigskip \noindent that is generated by the semi-algebra and $\pi$-system $\mathcal{S}$ of measurable rectangles or cylinders 
\begin{equation*}
\mathcal{S}=\{\prod_{1\leq i\leq k}A_{i},\text{ }A_{i}\in \mathcal{A}_{i}\}.
\end{equation*}

\bigskip \noindent The integer $k$ is taken as $k=2$. The sections of sets $A$ are denoted as $A_{\omega _{1}}$ and $A_{\omega _{2}}$ for $\omega _{1}\in \Omega_{1}$ and $\omega_{2}\in \Omega _{2}$.\\

\bigskip \noindent \textbf{Exercise 1}. \label{exercises01_doc07-02}\\

\noindent Question (a). Let $A=A_{1}\times A_{2} \in \mathcal{S}$, with $A_{i} \in \mathcal{A}_{i}$, $i=1,2$. By Question (b) Exercise 3 in Doc 01-06 if Chapter \ref{01_setsmes}, pages \pageref{exercise03_doc01-06} and \pageref{exercise03_doc01-06}, justify and show that for any $\omega_1 \in \Omega_1$,

\begin{equation*}
m_2\left(A_{\omega _{1}}\right)=\left\{ 
\begin{tabular}{lll}
$m_2(A_{2})$ & if & $\omega _{1}\in A_{1}$ \\ 
$0 $ & if & $\omega _{1}\notin A_{1}$%
\end{tabular}
\right. 
\end{equation*}

\bigskip \noindent and the application $\omega _{1}\mapsto m_{2}(A_{\omega _{1}})$ is $\mathcal{A}_1$-measurable.\\

\noindent Question (b). Suppose $m_1$ and $m_2$ are finite. Consider the collection of sets

$$
\mathcal{M}=\{A \in \mathcal{A}, \ the \ application \ \omega _{1}\mapsto m_{2}(A_{\omega _{1}}), \ is \ \mathcal{A}_1-measurable \}.
$$ 

\bigskip \noindent Show that $\mathcal{D}$ is Dynkin system containing the $\pi$-system $\mathcal{S}$. Use the $\lambda$-$\pi$ Lemma (\textit{See Point (01.04) in Doc 01-01 in Chapter \ref{01_setsmes}, 
page \pageref{doc01-01}}) to conclude that $\mathcal{D}=\mathcal{A}$ and conclude that  for any $A \in \mathcal{A}$, \ the \ application $\omega _{1}\mapsto m_{2}(A_{\omega _{1}})$  is  $\mathcal{A}_1$-measurable.\\

\noindent Question (c). Extend the conclusion of Question (a) to the case where $m_1$ and $m_2$ are $\sigma$-finite.\\

\noindent Question (d). Let $f$ be an $\mathcal{A}_1 \otimes \mathcal{A}_2$-non-negative measurable application $f$ and $\omega _{1}\in \Omega_1$ fixed. By using Point (03.22) in Doc 03.02 (Page
page \pageref{doc03-02}) in Chapter \ref{03_setsmes_applimes_cas_speciaux}, and by using the MCT, show that the partial function the application

$$
\omega_{2} \longmapsto f_{\omega_1}(\omega_2)=f(\omega_{1},\omega _{2})
$$

\bigskip \noindent is $\mathcal{A}_2$-measurable.\\

\noindent Question (e). By continuing the question (d), justify the definition of the function

$$
\omega_{1}\longmapsto \int_{\Omega _{2}}f(\omega _{1},\omega _{2})dm(\omega_{2})
$$

\bigskip \noindent and show it is $\mathcal{A}_1$-measurable.\\

\bigskip \noindent \textbf{Exercise 2}. \label{exercises02_doc07-02}.\\

\noindent Question (a). Based on Exercise 1, justify that  the application $m$ defined on $\mathcal{A}$ by%
\begin{equation}
m(A)=\int_{\Omega_{1}} m_{2}(A_{\omega_{1}}) \ dm_{1}(\omega_{1}) \label{mesprodL}
\end{equation}

\noindent is well-defined.\\

\noindent Question (b). Show that $m$ is a measure.\\

\noindent Question (c) Show that $m$ satisfies the following formula on the class of measurable rectangles :

\begin{equation*}
m(A_{1}\times A_{2})=m_{1}(A_{1})\times m_{2}(A_{2}),  \ \ (MESPROD)
\end{equation*}

\bigskip \noindent for any $A=A_{1}\times A_{2}$, with $A_{i} \in \mathcal{A}_{i}$, $i=1,2$.\\

\noindent Question (d) Show that there exists a unique measure $m$ on $(\Omega_1 \times \Omega_2, \mathcal{A}_{1} \times \mathcal{A}_{2})$, denoted 

$$
m=m_1 \otimes m_2.
$$

\bigskip \noindent such that (MESPROD) is valid. By exchanging the roles of $m_1$ and $m_2$ and by defining $m_2 \otimes m_1$ on $\mathcal{S}$ by

\begin{equation}
(m_2 \otimes m_1)(A)=\int_{\Omega_{2}} m_{1}(A_{\omega_{2}}) \ dm_{2}(\omega_{2}), \label{mesprodR}
\end{equation}

\bigskip \noindent conclude that $m_1 \otimes m_2=m_2 \otimes m_1$ and combine Formulas (\ref{mesprodL}) and (\ref{mesprodR}) to have for $f=1_A$, $A\in \mathcal{A}$

\begin{eqnarray*}
\int f \ d(m_1 \otimes m_2) &=& \int_{\Omega_1} \biggr( \int_{\Omega_2} f(\omega_1,\omega_2) \ dm_2(\omega_2) \biggr) \ dm_1(\omega_1)\\
&=&\int_{\Omega_2} \biggr( \int_{\Omega_1} f(\omega_1,\omega_2) \ dm_1(\omega_1) \biggr) \ dm_2(\omega_2). \ (FUB0)
\end{eqnarray*}

\bigskip \noindent \textbf{Exercise 3}. (Fubini's Theorem and Tonelli's fro non-negative functions) \label{exercises03_doc07-02}.\\

\noindent Here, we use the notation : $m=m_1 \otimes m_2 =m_2 \otimes m_1$.\\

\noindent Question (a). (One version of Fubini's Theorem) Let

\begin{equation*}
\begin{tabular}{lllll}
$f$ & $:$ & $\Omega =(\Omega _{1}\times \Omega _{2})$ & $\longmapsto $ & $%
\overline{\mathbb{R}}$ \\ 
&  & $(\omega _{1},\omega _{2})$ & $\hookrightarrow $ & $f(\omega
_{1},\omega _{2})$
\end{tabular}
\end{equation*}

\bigskip \noindent be a non-negative and measurable function, show that Formula (FUB0) of Exercise 2 still holds, by using the classical steps of the construction of the integral.\\

\noindent Question (b). (Tonelli's Theorem). Apply Question (a) to show that a non-negative and measurable function is integrable if and only if 
$$
\int_{\Omega _{1}} \biggr( \int_{\Omega_{2}}f(\omega _{1},\omega _{2}) \ dm_{2}(\omega _{2})\biggr) \ dm_{1}(\omega _{1}) <+\infty \  (TON1)
$$

\bigskip \noindent or

$$
\int_{\Omega _{2}} \biggr( \int_{\Omega_{1}}f(\omega _{1},\omega_{2}) \ dm_{1}(\omega_{1})\biggr) \ dm_{2}(\omega_{2}) <+\infty. \  (TON2).
$$

\bigskip \noindent \textbf{Exercise 4}. (Fubini's Theorem and Tonelli's fro non-negative functions) \label{exercises04_doc07-02}.\\

\noindent We still use the notation : $m=m_1 \otimes m_2 =m_2 \otimes m_1$.\\

\noindent Let 

\begin{equation*}
\begin{tabular}{lllll}
$f$ & $:$ & $\Omega =(\Omega _{1}\times \Omega _{2})$ & $\longmapsto $ & $%
\overline{\mathbb{R}}$ \\ 
&  & $(\omega _{1},\omega _{2})$ & $\hookrightarrow $ & $f(\omega
_{1},\omega _{2})$
\end{tabular}
\end{equation*}

\bigskip \noindent be a and integrable function with respect to $m$. Show that we have the Fubini's formula

\begin{eqnarray*}
\int f \ d(m_1 \otimes m_2) &=& \int_{\Omega_1} \biggr( \int_{\Omega_2} f(\omega_1,\omega_2) \ dm_2(\omega_2) \biggr) \ dm_1(\omega_1)\\
&=&\int_{\Omega_2} \biggr( \int_{\Omega_1} f(\omega_1,\omega_2) \ dm_1(\omega_1) \biggr) \ dm_2(\omega_2). \ (FUB).
\end{eqnarray*}

\noindent \textit{Hints}.\\

\noindent (a) Let $f$ be integrable. By using Tonelli's Formula, what what can say about the functions

$$
\omega_1 \mapsto \int_{\Omega_2} f^+(\omega_1,\omega_2) \ dm_2(\omega_2) \ and \ \omega_1 \mapsto \int_{\Omega_2} f^-(\omega_1,\omega_2) \ dm_2(\omega_2), \ \ (PNP) 
$$

\bigskip \noindent in terms of finiteness and measurability?\\

\noindent (b) Extend the conclusion to the function

$$
\omega_1 \mapsto \int_{\Omega_2} f(\omega_1,\omega_2) \ dm_2(\omega_2). 
$$

\bigskip \noindent (c) Try to conclude by using Formula (TON1) first and next Formula (TON2). Combine the two conclusion to get a final one.\\

\noindent \LARGE \textbf{Doc 07-03 : Finite Product Measures - Exercises with Solutions}. \label{doc07-03}\\
\bigskip
\Large

\bigskip \noindent In all this document, we work with $k$ $\sigma$-finite measure spaces $(\Omega _{i},A_{i},m_{i})$, $i=1,...,k$. The product space 
\begin{equation*}
\Omega =\prod_{1\leq i\leq k}\Omega _{i}
\end{equation*}

\bigskip \noindent endowed with the product $\sigma -$algebra 
\begin{equation*}
\mathcal{A}=\bigotimes_{1\leq i\leq k}\mathcal{A}_{i}.
\end{equation*}

\bigskip \noindent that is generated by the semi-algebra and $\pi$-system $\mathcal{S}$ of measurable rectangles or cylinders 
\begin{equation*}
\mathcal{S}=\{\prod_{1\leq i\leq k}A_{i},\text{ }A_{i}\in \mathcal{A}_{i}\}
\end{equation*}

\bigskip \noindent The integer $k$ is taken as $k=2$. The sections of sets $A$ are denoted as $A_{\omega _{1}}$ and $A_{\omega _{2}}$ for $\omega _{1}\in \Omega_{1}$ and $\omega_{2}\in \Omega _{2}$.\\

\bigskip \noindent \textbf{Exercise 1}. \label{exercises01_sol_doc07-03}.\\

\noindent Question (a). Let $A=A_{1}\times A_{2} \in \mathcal{S}$, with $A_{i} \in \mathcal{A}_{i}$, $i=1,2$. By Question (b) Exercise 3 in Doc 01-06 if Chapter \ref{01_setsmes}, pages \pageref{exercise03_doc01-06} and \pageref{exercise03_doc01-06}, justify and show that for any $\omega_1 \in \Omega_1$,

\begin{equation*}
m_2\left(A_{\omega _{1}}\right)=\left\{ 
\begin{tabular}{lll}
$m_2(A_{2})$ & if & $\omega _{1}\in A_{1}$ \\ 
$0 $ & if & $\omega _{1}\notin A_{1}$%
\end{tabular}
\right. 
\end{equation*}

\bigskip \noindent and the application $\omega _{1}\mapsto m_{2}(A_{\omega _{1}})$ is $\mathcal{A}_1$-measurable.\\

\noindent Question (b). Suppose $m_1$ and $m_2$ are finite. Consider the collection of sets

$$
\mathcal{M}=\{A \in \mathcal{A}, \ the \ application \ \omega _{1}\mapsto m_{2}(A_{\omega _{1}}), \ is \ \mathcal{A}_1-measurable \}.
$$ 

\bigskip \noindent Show that $\mathcal{D}$ is Dynkin system containing the $\pi$-system $\mathcal{S}$. Use the $\lambda$-$\pi$ Lemma (\textit{See Point (01.04) in Doc 01-01 in Chapter \ref{01_setsmes}, 
page \pageref{doc01-01}}) to conclude that $\mathcal{D}=\mathcal{A}$ and conclude that  for any $A \in \mathcal{A}$, \ the \ application $\omega _{1}\mapsto m_{2}(A_{\omega _{1}})$  is  $\mathcal{A}_1$-measurable.\\

\noindent Question (c). Extend the conclusion of Question (a) to the case where $m_1$ and $m_2$ are $\sigma$-finite.\\

\noindent Question (d). Let $f$ be an $\mathcal{A}_1 \otimes \mathcal{A}_2$-non-negative measurable application $f$ and $\omega _{1}\in \Omega_1$ fixed. By using Point (03.22) in Doc 03.02 (Page
page \pageref{doc03-02}) in Chapter \ref{03_setsmes_applimes_cas_speciaux}, and by using the MCT, show that the partial function the application

$$
\omega_{2} \longmapsto f_{\omega_1}(\omega_2)=f(\omega_{1},\omega _{2})
$$

\bigskip \noindent is $\mathcal{A}_2$-measurable.\\

\noindent Question (e). By continuing the question (d), justify the definition of the function

$$
\omega_{1}\longmapsto \int_{\Omega _{2}}f(\omega _{1},\omega _{2})dm(\omega_{2}),
$$

\bigskip \noindent and show it is $\mathcal{A}_1$-measurable.\\

\bigskip \noindent \textbf{SOLUTIONS}.\\

\noindent Question (a). We know from Exercise 3 in Doc 01-06 in Chapter \ref{01_setsmes} (page \pageref{exercise03_doc01-06}) that for any $\omega_1 \in \Omega_1$, for any measurable set 
$A \in \mathcal{A}$, the section $A_{\omega_{1}}$ is $\mathcal{A}_1$-measurable. In particular, if $A=A_{1}\times A_{2} \in \mathcal{S}$, with $A_{i} \in \mathcal{A}_{i}$, $i=1,2$, we have

\begin{equation*}
A_{\omega_{1}}=\left\{ 
\begin{tabular}{lll}
$A_{2}$ & if & $\omega _{1}\in A_{1}$ \\ 
$\emptyset $ & if & $\omega _{1}\notin A_{1}$%
\end{tabular}
\right. .
\end{equation*}

\noindent So by applying the measure $m_2$, we have 

\begin{equation*}
m_2\left(A_{\omega_{1}}\right)=\left\{ 
\begin{tabular}{lll}
$m_2(A_{2})$ & if & $\omega _{1}\in A_{1}$ \\ 
$0 $ & if & $\omega _{1}\notin A_{1}$%
\end{tabular}
\right. .
\end{equation*}

\bigskip \noindent This means that the application $\omega _{1}\mapsto m_{2}(A_{\omega _{1}})$ is an elementary function on $(\Omega, \mathcal{A}_1)$ and hence, is $\mathcal{A}_1$-measurable.\\

\noindent Question (b). The definition of $\mathcal{M}$ is meaningful since $A_{\omega _{1}}$ is $\mathcal{A}_2$-measurable. Let us begin to show that $\mathcal{M}$ is a Dynkin-System.\\

\noindent (i) We surely have that $\emptyset \in \mathcal{M}$, since $m_{2}(\emptyset_{\omega _{1}})=m_{2}(\emptyset)=0$ and the application $\omega _{1}\mapsto m_{2}(\emptyset_{\omega _{1}})=0$ is measurable. We similarly show that $\Omega \in \mathcal{M}$, since $m_{2}(\Omega_{\omega _{1}})=m_{2}(\Omega_2)$ which is constant.\\ 

\noindent (ii) Let $(A,B) \in \mathcal{A}^c$ with $A\subset B$. By the properties of the sections as seen in Exercise 3 in Doc 01-06 if Chapter \ref{01_setsmes} (page \pageref{exercise03_doc01-06}), we have for any $\omega_1 \in \Omega_1$,  $A_{\omega _{1}} \subset B_{\omega _{1}}$, and since $m_2$ is finite, we get 
$$
m_{2}( (B \setminus A)_{\omega_1})=m_2(B_{\omega_{1}})-m_2(A_{\omega_{1}}).
$$

\bigskip \noindent Hence the application $m_{2}( (B \setminus A)_{\omega _{1}})=m_2(B_{\omega _{1}})-m_2(A_{\omega _{1}})$ is well-defined and is $\mathcal{A}_1$-measurable as the difference of two $\mathcal{A}_1$-measurable applications.\\

\noindent (iii) Let $(A_n)_{n\geq 0} \subset \mathcal{A}$ pairwise disjoint. By the properties of the sections as seen in Exercise 3 in Doc 01-06 if Chapter \ref{01_setsmes} (page \pageref{exercise03_doc01-06}), we have for any $\omega_1 \in \Omega_1$

\begin{eqnarray*}
\left(\sum_{n \geq 0} A_{n} \right)_{\omega_1} &=&\sum_{n \geq 0}  \left( A_{n}\right)_{\omega_1}\\
&=&  \lim_{k\rightarrow +\infty} \sum_{n=0}^{k} \left( A_{n}\right)_{\omega_1}, \ (LIM01)
\end{eqnarray*}

\noindent which implies that for any for any $\omega_1 \in \Omega_1$, we have by finite additivity and the MCT,

\begin{eqnarray*}
m_2\left(\left(\sum_{n \geq 0} A_{n} \right)_{\omega_1}\right) &=&\sum_{n \geq 0}  m_2(\left( A_{n}\right)_{\omega1})\\
&=&  \lim_{k\rightarrow +\infty} \sum_{n=0}^{k} m_2(\left( A_{n}\right)_{\omega_1}), \ (LIM02)
\end{eqnarray*}

\noindent and thus, the application $\omega_1 \mapsto m_2\left(\left(\sum_{n \geq 0} A_{n} \right)_{\omega_1}\right)$ is $\mathcal{A}_1$-measurable as limit of sequence of $\mathcal{A}_1$-measurable applications.\\

\noindent  In conclusion, $\mathcal{D}$ is contained in $\mathcal{A}$, is a $\lambda$-system and contains the $\pi$-system $\mathcal{S}$ which generates $\mathcal{A}$. By the $\lambda$-$\pi$ Lemma, $\mathcal{A}$ is also the generated by $\mathcal{S}$, so that we will have $\mathcal{A} \subset \mathcal{D}$. We get $\mathcal{D}=\mathcal{S}$ and we may conclude.\\

\noindent (c) Since $m_2$ is $\sigma$-finite, there exists a measurable and countable partition of $\Omega_2$, that is

$$
\Omega_2=\sum_{j\geq 0} \Omega_{2,j}, \ (\forall j\geq 0, \ \Omega_{2,j} \in \mathcal{A}_2 \ and \ m_2(\Omega_{2,j} \cap B)<+\infty). 
$$

\bigskip \noindent As seen previously, $m_2$ is the countable sum of the finite measures 
$$
m_{2,j}(B)=m_2(\Omega_{2,j} \cap B), \ B\in \mathcal{A}_2,
$$

\bigskip \noindent  that is, for all $B\in \mathcal{A}_2$,

$$
m_2(B)= \sum_{j\geq 0} m_{2,j}(B). 
$$

\bigskip \noindent Hence, for any $\omega_1 \in \Omega_1$, for any $A \in \mathcal{A}$,

$$
m_2(A_{\omega_1})= \sum_{j\geq 0} m_{2,j}(A_{\omega_1})=\lim_{k \rightarrow +\infty} \sum_{j=1}^{k} m_{2,j}(A_{\omega_1}). 
$$
 
\bigskip \noindent Hence, by Question (b), the application $\omega_1 \mapsto m_2(A_{\omega_1})$ is $\mathcal{A}_1$-measurable as limit of sequence of $\mathcal{A}_1$-measurable applications.\\

\noindent \textbf{(d)}. Le  $f$ be an $\mathcal{A}_1 \otimes \mathcal{A}_2$-non-negative measurable application and let $\omega _{1}\in \Omega$ be fixed. Let us consider the partial application
\begin{equation*}
\omega _{1}\longmapsto f_{\omega_1}(\omega_2)=f(\omega _{1},\omega _{2}).
\end{equation*}

\bigskip \noindent By virtue of Point (03.22) in Doc 03.02 (Page page \pageref{doc03-02}) in Chapter \ref{03_setsmes_applimes_cas_speciaux}, we have a non-decreasing sequence of elementary functions

$$
f_n=\sum_{k=1}^{\ell_n} a_{n,k} 1_{A_{n,k}}, \ a_{n,k}\geq 0, \ A_{n,k} \in \mathcal{A})\ (\ 1\leq k \leq \ell_n), \ A_{n,1}+...+A_{n,\ell_n}=\Omega, \ \ell_n\geq 1,
$$

\bigskip \noindent such that $f_n \uparrow f$ as $n\uparrow +\infty$. It is also clear that for any $\omega _{1}\in \Omega$, the sequence of partial functions $(f_n)_{\omega_1}$ is a non-decreasing sequence to $f_{\omega_1}$.\\

\noindent Besides, we also know (See Doc 03-03, page \pageref{doc03-03}, Point b2) that for any $A \in \Omega_1 \times \Omega_2$, the partial function of $1_A$ satisfies for $\omega_1 \in \Omega_1$ fixed,

$$
\omega_2 \mapsto \biggr(1_A\biggr)_{\omega_1}(\omega_2)=1_A(\omega_1,\omega_2)=1_{A_{\omega_1}}(\omega_2).
$$

\bigskip \noindent So for any $\omega _{1}\in \Omega$, we have by the linearity of the partial function formation, 

$$
(f_n)_{\omega_1}(\omega_2)=\sum_{k=1}^{\ell_n} a_{n,k} 1_{(A_{n,k})_{\omega_1}}(\omega_2).
$$

\bigskip \noindent This is a sequence of elementary functions on $(\Omega_2, \mathcal{A}_2)$ non-decreasing to $f_{\omega_1}$. So $f_{\omega_1}$ is $\mathcal{A}_2$-measurable.\\

\noindent Question (e) We already know that, for each $n\geq 0$, the function

$$
\omega_1 \mapsto \int (f_n)_{\omega_1} \ dm_2=\sum_{k=1}^{\ell_n} a_{n,k} m_2((A_{n,k})_{\omega_1}).
$$

\bigskip \noindent is $\mathcal{A}_1$-measurable. By taking the limit and by using the MCT, or simply the definition of the integral for a non-negative function we get that the limit function

$$
\omega_1 \mapsto \int f_{\omega_1}(\omega_2) \ dm_2(\omega_2)=\int f(\omega_1, \omega_2) \ dm_2(\omega_2)
$$

\bigskip \noindent is $\mathcal{A}_1$-measurable.\\

\bigskip \noindent \textbf{Exercise 2}. \label{exercises02_sol_doc07-03}.\\

\noindent Question (a). Based on Exercise 1, justify that  the application $m$ defined on $\mathcal{A}$ by%
\begin{equation}
m(A)=\int_{\Omega_{1}} m_{2}(A_{\omega_{1}}) \ dm_{1}(\omega_{1}) \label{mesprodL}
\end{equation}

\noindent is well-defined.\\

\noindent Question (b). Show that $m$ is a measure.\\

\noindent Question (c) Show that $m$ satisfies the following formula on the class of measurable rectangles :

\begin{equation*}
m(A_{1}\times A_{2})=m_{1}(A_{1})\times m_{2}(A_{2}),  \ \ (MESPROD)
\end{equation*}

\bigskip \noindent for any $A=A_{1}\times A_{2}$, with $A_{i} \in \mathcal{A}_{i}$, $i=1,2$.\\

\noindent Question (d) Show that there exists a unique measure $m$ on $(\Omega_1 \times \Omega_2, \mathcal{A}_{1} \times \mathcal{A}_{2})$, denoted 

$$
m=m_1 \otimes m_2.
$$

\bigskip \noindent such that (MESPROD) is valid. By exchanging the roles of $m_1$ and $m_2$ and by defining $m_2 \otimes m_1$ on $\mathcal{S}$ by

\begin{equation}
(m_2 \otimes m_1)(A)=\int_{\Omega_{2}} m_{1}(A_{\omega_{2}}) \ dm_{2}(\omega_{2}), \label{mesprodR}
\end{equation}

\bigskip \noindent conclude that $m_1 \otimes m_2=m_2 \otimes m_1$ and combine Formulas (\ref{mesprodL}) and (\ref{mesprodR}) to have for $f=1_A$, $A\in \mathcal{A}$

\begin{eqnarray*}
\int f \ d(m_1 \otimes m_2) &=& \int_{\Omega_1} \biggr( \int_{\Omega_2} f(\omega_1,\omega_2) \ dm_2(\omega_2) \biggr) \ dm_1(\omega_1)\\
&=&\int_{\Omega_2} \biggr( \int_{\Omega_1} f(\omega_1,\omega_2) \ dm_1(\omega_1) \biggr) \ dm_2(\omega_2). \ (FUB0)
\end{eqnarray*}

\bigskip \noindent \textbf{SOLUTIONS}.\\

\noindent Question (a). By Exercise 1, for any $A\in \mathcal{A}$, the application $\omega_1 \mapsto m_{2}(A_{\omega_{1}})$ is non-negative and $\mathcal{A}_1$-measurable. Its integral with respect is 
$m_1$ is well-defined, so is that of $m(A)$.\\

\noindent It is clear that $m(\emptyset)=\int m_2(\emptyset_{\omega_1}) \ dm_1=\int 0 \ dm_1=0$. Next, let $(A_n)_{n\geq 0} \subset \mathcal{A}$ be pairwise disjoint. We have

\begin{eqnarray*}
m\left(\sum_{n \geq 0} A_{n} \right)&=& \int m_2\biggr(\left(\sum_{n \geq 0} A_{n} \right)_{\omega_1}\biggr) \ dm_1\\
&=&\int \sum_{n \geq 0} m_2(\left( A_{n} \right)_{\omega_1}) \ dm_1 \ (I1) \\
&=& \sum_{n \geq 0} \int m_2(\left( A_{n} \right)_{\omega_1}) \ dm_1 \ (I2) \\
&=& \sum_{n \geq 0} m(A_n).
\end{eqnarray*}

\noindent We have that (I1) is justified by the fact that the section of a sum of sets is the sum of their sections and (I2) is a consequence of the MCT. Hence, $m$ is a measure.\\

\bigskip \noindent Question (c).\\

\noindent Let Let $A=A_{1}\times A_{2} \in \mathcal{S}$, with $A_{i} \in \mathcal{A}_{i}$, $i=1,2$. By Question (a) in Exercise 1 above, we have for any $\omega_1 \in \Omega_1$

$$
m_2(A_{\omega_1})=m_2(A_2) 1_{A_{1}}(\omega_1),
$$

\bigskip \noindent and hence

\begin{eqnarray*}
m(A)&=& \int m_2(A_{\omega_1}) \ dm_1\\
&=& \int m_(A_2) 1_{A_{2}} \ dm_1\\
&=& m_2(A_2) \int  1_{A_{1}} \ dm_1\\
&=&m_1(A_1) m_2(A_2)
\end{eqnarray*}

\bigskip \noindent and hence, (MESPROD) holds.\\

\noindent Question (d). We already know that there exists a measure on $(\Omega_1 \times \Omega_2, \mathcal{A}_{1} \times \mathcal{A}_{2})$ which satisfies (MESPROD). Now let us define a application $m_0$ on $\mathcal{S}$ by (MESPROD), that is

$$
m_0(A_{1}\times A_{2})=m_{1}(A_{1})\times m_{2}(A_{2}),  \ \ (MESPROD)
$$

\bigskip \noindent for $A=A_{1}\times A_{2} \in \mathcal{S}$, with $A_{i} \in \mathcal{A}_{i}$, $i=1,2$. The application $m_0$ is non-negative and proper since $m_0(\emptyset \times \emptyset)=0$. It is additive since the only way to have that a sum of two elements of $\mathcal{S}$ as an element of $\mathcal{S}$, is to split one the factor of an elements of $\mathcal{S}$, that is, for example

$$
A_{1}\times A_{2}= A_{11}\times A_{2} + A_{12}\times A_{2},
$$

\bigskip \noindent with 

$$
A_{1}=A_{11}+A_{12}, \ A_{11} \in \mathcal{A}_1 \ and \ A_{11} \in \mathcal{A}_1.
$$

\bigskip \noindent In such a case, we have

\begin{eqnarray*}
m_0(A_{1}\times A_{2})&=&m_1(A_1) m_2(A_2)\\
&=&m_1(A_{11}+A_{12})=(m_1(A_{11})+m_1(A_{12})) m_2(A_2)\\
&=&=m_1(A_{11}) m_2(A_2)+ m_1(A_{12}) m_2(A_2)\\
&=&m_0(A_{11}\times A_{2})+m_0(A_{12}\times A_{2}).
\end{eqnarray*}

\noindent It is also clear that $m_0$ is $\sigma$-finite. By the Exercise 15 on extensions (See page \pageref{exercise15_doc04-02}) in Doc 04-02 in Chapter \ref{04_measures}, we know that $m_0$ is extensible to a proper and additive and $\sigma$-finite application on $a(\mathcal{S})$. It is clear that $m$ is an extension of $m_0$ which is $\sigma$-additive and $\sigma$-finite. By Caratheodory's Theorem (See Point (04-26) in Doc 03.04, page  \pageref{doc03-04}), or by the Exercise 10 on determining classes for $\sigma$-finite measures in Doc 04-05 (See page \pageref{exercise15_doc04-02}) in the same chapter, we know that the extension is unique.\\

\noindent We may conclude by saying :  there exists a unique measure $m$ on $(\Omega_1 \times \Omega_2, \mathcal{A}_{1} \times \mathcal{A}_{2})$, denoted 

$$
m=m_1 \otimes m_2,
$$

\bigskip \noindent such that (MESPROD) is valid.\\

\noindent The conclusions still hold if we exchange the roles of $m_1$ and $m_2$. We obtain a $\sigma$-finite measure denoted by $m_2 \otimes m_1$, defined on $\mathcal{S}$ by 
Formula (\ref{mesprodR}) and, of course ,satisfying (MESPROD).\\

\noindent The two measures $m_1 \otimes m_2$ and $m_2 \otimes m_1$ are $\sigma$-finite and are equal on the $\pi$-system $\mathcal{S}$. So they are equal 
(on $\mathcal{A}$). By combining Formulas (\ref{mesprod}) and (\ref{mesprodR}) , we see that Formula (FUB0) holds.

\bigskip \noindent \textbf{Exercise 3}. (Fubini's Theorem and Tonelli's fro non-negative functions) \label{exercises03_sol_doc07-03}.\\

\noindent Here, we use the notation : $m=m_1 \otimes m_2 =m_2 \otimes m_1$.\\

\noindent Question (a). (One version of Fubini's Theorem) Let

\begin{equation*}
\begin{tabular}{lllll}
$f$ & $:$ & $\Omega =(\Omega _{1}\times \Omega _{2})$ & $\longmapsto $ & $%
\overline{\mathbb{R}}$ \\ 
&  & $(\omega _{1},\omega _{2})$ & $\hookrightarrow $ & $f(\omega
_{1},\omega _{2})$
\end{tabular}
\end{equation*}

\bigskip \noindent be a non-negative and measurable function, show that Formula (FUB0) of Exercise 2 still holds, by using the classical steps of the construction of the integral.\\

\noindent Question (b). (Tonelli's Theorem). Apply Question (a) to show that a non-negative and measurable function is integrable if and only if 
$$
\int_{\Omega _{1}} \biggr( \int_{\Omega_{2}}f(\omega _{1},\omega _{2}) \ dm_{2}(\omega _{2})\biggr) \ dm_{1}(\omega _{1}) <+\infty \  (TON1)
$$

\bigskip \noindent or

$$
\int_{\Omega _{2}} \biggr( \int_{\Omega_{1}}f(\omega _{1},\omega_{2}) \ dm_{1}(\omega_{1})\biggr) \ dm_{2}(\omega_{2}) <+\infty. \  (TON2).
$$

\bigskip \noindent \textbf{SOLUTIONS}.\\

\noindent Question (a). Formula (FUB0) is readily extended to non-negative elementary functions by finite linearity of the integral on the class of non-negative functions. Next, it is extended on the class of non-negative and measurable functions by the MCT.\\

\noindent Question (b). The solution is obvious and is readily derived from Question (a).\\

\bigskip \noindent \textbf{Exercise 04}. (Fubini's Theorem and Tonelli's fro non-negative functions) \label{exercises04_sol_doc07-03}.\\

\noindent We still use the notation : $m=m_1 \otimes m_2 =m_2 \otimes m_1$.\\

\noindent Let 

\begin{equation*}
\begin{tabular}{lllll}
$f$ & $:$ & $\Omega =(\Omega _{1}\times \Omega _{2})$ & $\longmapsto $ & $%
\overline{\mathbb{R}}$ \\ 
&  & $(\omega _{1},\omega _{2})$ & $\hookrightarrow $ & $f(\omega
_{1},\omega _{2})$
\end{tabular}
\end{equation*}

\bigskip \noindent be a and integrable function with respect to $m$. Show that we have the Fubini's formula

\begin{eqnarray*}
\int f \ d(m_1 \otimes m_2) &=& \int_{\Omega_1} \biggr( \int_{\Omega_2} f(\omega_1,\omega_2) \ dm_2(\omega_2) \biggr) \ dm_1(\omega_1)\\
&=&\int_{\Omega_2} \biggr( \int_{\Omega_1} f(\omega_1,\omega_2) \ dm_1(\omega_1) \biggr) \ dm_2(\omega_2). \ (FUB)
\end{eqnarray*}

\noindent \textit{Hints}.\\

\noindent (a) Let $f$ be integrable. By using Tonelli's Formula, what what can say about the functions

$$
\omega_1 \mapsto \int_{\Omega_2} f^+(\omega_1,\omega_2) \ dm_2(\omega_2) \ and \ \omega_1 \mapsto \int_{\Omega_2} f^-(\omega_1,\omega_2) \ dm_2(\omega_2), \ \ (PNP) 
$$

\bigskip \noindent in terms of finiteness and measurability?\\

\noindent (b) Extend the conclusion to the function

$$
\omega_1 \mapsto \int_{\Omega_2} f(\omega_1,\omega_2) \ dm_2(\omega_2). 
$$

\bigskip \noindent (c) Try to conclude by using Formula (TON1) first and next Formula (TON2). Combine the two conclusion to get a final one.\\

\bigskip \noindent \textbf{SOLUTIONS}.\\

\noindent (a). Since $f$ is integrable, the integrals $\int f^+ \ dm$ and $\int f^+ \ dm$ are finite. By Exercise 2, we are going to apply (TON1) Formula and next Formula (TON2). From Formula (TON1), we see that the function

$$
\omega_1 \mapsto \int_{\Omega_2} f^+(\omega_1,\omega_2) \ dm_2(\omega_2), 
$$

\bigskip \noindent is $m_1$-\textbf{a.e.} finite (since it is $m_1$-integrable). We already knew from Exercise 1 that it is $\mathcal{A}_1$-measurable. By doing the same for the negative part, we get that the two functions given in Formula (PNP) are both $m_1$-\textbf{a.e.} finite and $\mathcal{A}_1$-measurable. By using the linearity of the integral on $\mathcal{L}^1(\Omega_2,\mathcal{A}_2,m_2)$, we may define 

$$
\omega_1 \mapsto \int_{\Omega_2} f^+(\omega_1,\omega_2) \ dm_2(\omega_2)-\int_{\Omega_2} f^-(\omega_1,\omega_2) \ dm_2(\omega_2)=\int_{\Omega_2} f(\omega_1,\omega_2) \ \ dm_2(\omega_2).
$$

\bigskip \noindent which becomes $m_1$-\textbf{a.e.} finite and $\mathcal{A}_1$-measurable. By integrating this function with respect to $m_1$ and by applying Fubuni's and Tonelli's Formulas for non-negative functions, we get (by beginning by the end), 

\begin{eqnarray*}
&&\int_{\Omega_1} \biggr( \int_{\Omega_2} f(\omega_1,\omega_2) \ dm_2(\omega_2)\biggr) dm_1(\omega_1)\\
&=& \int_{\Omega_1} \biggr( \int_{\Omega_2} f^+(\omega_1,\omega_2) \ dm_2(\omega_2) \biggr) dm_1(\omega_1)\\
&-& \int_{\Omega_1} \biggr( \int_{\Omega_2} f^-(\omega_1,\omega_2) \ dm_2(\omega_2) \biggr) dm_1(\omega_1)\\
&=& \int f^+ \ dm - \int f^- dm = \int f \ dm.
\end{eqnarray*}

\noindent By applying Formula (TON2), we proceed similarly and get

$$
\int_{\Omega_2}  \biggr( \int_{\Omega_1} f(\omega_1,\omega_2) \ dm_1(\omega_1)\biggr) dm_2(\omega_2)= \int f \ dm.
$$

\bigskip \noindent Combining the results of the two paths leads to Formula (FUB).\\

\newpage

\bigskip \noindent \textbf{General conclusion and extension}.\\

\noindent (a) The solutions of the exercise are summarized in Doc 07-01, which usually constitutes the main course in finite product $\sigma$-finite measures.\\

\noindent We did not need any topological properties on the spaces $\Omega_i$, $1\leq i \leq k$, to prove its existence. On topological spaces like $\mathbb{R}^d$, it is possible to derive the product measure of finite measures from the Caratheordory's Theorem by using the their topological properties, especially by using the nice fact that on such spaces, finite measures are determined by their values on compact sets. We will come back to such constructions in the Exercises book to come later.\\

\noindent (b) (\textbf{Extension}). We want to close this document replacing the $\sigma$-finiteness by a more general condition. If we summarize all the points discussed here, the $\sigma$-finiteness intervenes to ensure the uniqueness of the product measure and to show that : \\

\noindent (H1) For all $A \in \mathcal{A}$,  the  application $\omega _{1}\mapsto m_{2}(A_{\omega _{1}})$  is  $\mathcal{A}_1$-measurable.\\

\noindent From this, all the other properties readily follow. This leads to this extension.\\

\noindent Let $(\Omega_i, \mathcal{A}_i, m_i)$, $i \in \{1,2\}$ be two arbitrary measure spaces such that for all $A \in \mathcal{A}$,

$$
\omega_{2}\mapsto m_{1}(A_{\omega _{2}})
$$  

\bigskip \noindent is  $\mathcal{A}_1$-measurable. Then the following application, which is defined on $\mathcal{A}$ by 

\begin{equation}
(m_1\otimes m_2)(A)=\int_{\Omega_{1}} m_{2}(A_{\omega_{1}}) \ dm_{1}(\omega_{1}),
\end{equation}

\bigskip \noindent is a measure such that \\

\noindent (a) for any $A=A_{1}\times A_{2}$, with $A_{i} \in \mathcal{A}_{i}$, $i=1,2$, we have

\begin{equation*}
(m_1\otimes m_2)(A_{1}\times A_{2})=m_{1}(A_{1})\times m_{2}(A_{2})  \ \ (MESPROD)
\end{equation*}

\bigskip \noindent and \\

\noindent (b) for any real-valued function $f$ defined on $\Omega_1 \times \Omega_2$ and $\mathcal{A}_1\otimes \mathcal{A}_2$-measurable, non-negative or $(m_1\otimes m_2)$-integrable, we have

$$
\int f \ d(m_1 \otimes m_2) =\int_{\Omega_1} \biggr( \int_{\Omega_2} f(\omega_1,\omega_2) \ dm_2(\omega_2) \biggr) \ dm_1(\omega_1).\\
$$

\bigskip \noindent We also have that all $A \in \mathcal{A}$, the  application 
$$
\omega_{2}\mapsto m_{1}(A_{\omega _{2}})
$$  
 
\bigskip \noindent \bigskip \noindent is  $\mathcal{A}_2$-measurable, we may define $m_2 \otimes m_1$ similarly. And we cannot say that the have the equality $m_1 \otimes m_1=m_2 \otimes m_1$ without further conditions.

%\chapter{Radon-Nikodym Theorem}
%\newpage
\chapter{$\sigma$-additive sets applications and Radon-Nikodym's Theorem} \label{08_nikodym}

\noindent \textbf{Content of the Chapter}

\begin{table}[htbp]
	\centering
		\begin{tabular}{llll}
		\hline
		Type& Name & Title  & page\\
		S& Doc 08-01 &  Indefinite Integrals and Radon-Nikodym's Theorem & \pageref{doc08-01}\\
		D& Doc 08-02 &  $\sigma$-additive applications - Exercises & \pageref{doc08-02}\\
		D& Doc 08-03 &  Radon-Nikodym's Theorem - Exercises & \pageref{doc08-03}\\
		SD& Doc 08-04 &  $\sigma$-additive applications - Exercises/Solutions & \pageref{doc08-04}\\
		SD& Doc 08-05 &  Radon-Nikodym's Theorem - Exercises/Solutions & \pageref{doc08-05}\\
		D& Doc 08-06 &  Exercise on a continuous family of measures - Exercise & \pageref{doc08-06}\\
		SD& Doc 08-07 &  Exercise on a continuous family of measures/Solutions & \pageref{doc08-07}\\
		\hline
		\end{tabular}
\end{table}

\newpage
\noindent \LARGE \textbf{DOC 08-01 : Radon-Nikodym's Theorem - A summary} \label{doc08-01}\\
\bigskip
\Large

\bigskip \noindent Radon-Nikodym's Theorem is very important in modern analysis. As to probability theory, it is used for many purposes. But two of them are certainly the most important:\\

\bigskip \noindent (a) The probability density functions, with respect to the Lebesgue measure or to the counting measure.\\

\bigskip \noindent (b) The modern and general theory of conditional expectation which is founded on this theorem.\\

\bigskip \noindent  The proof of this theorem is lengthy. You will be provided of it. But the most important is to know where its comes and how it works. This is what we will do now.

\bigskip \noindent  This documents put it in the general context of $\sigma$-additive applications.\\

\bigskip \noindent  \textbf{Part I : $\sigma$-additive applications}.\\

\noindent Let us consider a measurable space $(\Omega, \mathcal{A})$. To be c, the definition of the $\sigma$-additivity  of a real-valued application $\phi$ on $\mathcal{A}$ requires that $\phi$ should not take both values $-\infty$ and $+\infty$. We have to avoid one them so that the relation $\phi(A+B)=\phi(A)+\phi(B)$ makes sense for any pair of disjoint elements of $\mathcal{A}$. By convention, we exclude the value $-\infty$ in our exposition. \textit{But, the situation where the value $+\infty$ us exclude is possible. In such a case, we need to adapt the results}.\\

\noindent So, by proper application on $\mathcal{A}$, we mean an application from $\mathcal{A}$ to $\mathbb{R}\cup{+\infty}$ which takes at least a finite value.\\

\noindent On top of $\sigma$-additive applications or measure (non-negative $\sigma$-additive applications) defined on $(\Omega, \mathcal{A})$, we have the following definitions. All the $\sigma$-additive applications considered here are supposed to be proper and defined on $(\Omega, \mathcal{A})$.\\

\noindent \textbf{(08.01a)}. A $\sigma$-additive application $\phi_2$ is continuous with respect to another A $\sigma$-additive application $\phi_1$, denoted $\phi_2 << \phi_1$, if and only if :

$$
\forall A\in \mathcal{A}, \  (\phi_1(A)=0) \Rightarrow (\phi_2(A)=0), 
$$

\noindent that is, any $\phi_1$-null set is a $\phi_2$-null set.\\

\noindent \textbf{(08.01b)}. A $\sigma$-additive application $\phi_2$ is singular with respect to another A $\sigma$-additive application $\phi_1$ if and only if :

$$
(\exists N\in \mathcal{A} : \ \phi_1(N)=0)  \ and \ (\forall A \in \mathcal{A}, \phi_2(A\cap N^c)=0).
$$

\noindent Remark that in this situation, for all $A \in \mathcal{A}$, $\phi_2(A)=\phi_2(A\cap N^c)$ and $\phi_1(A)=\phi_2(A\cap N)$.\\

\bigskip \noindent \textbf{(08.02)} Properties of $\sigma$-additive applications.\\
 .

\noindent \textbf{(08.02a)}  Let $\phi$ be $\sigma$-additive on $\mathcal{A}$. Let $(A_n)_{n\geq 0}$ be a sequence of mutually disjoint of elements of  such that $\phi\left( \sum_{n\geq 0} A_n \right)$ is finite. Then the series $\sum_{n\geq 0}  \phi\left( A_n \right)$ is absolutely convergent.\\

\noindent \textbf{(08.02b)} Let $\phi$ be $\sigma$-additive on $\mathcal{A}$. If $A$ and $B$ are two elements of $\mathcal{A}$ such that $A \subset B$. Then $\phi(A)$ is finite whenever $\phi(A)$. In particular, $\phi$ is finite if $\phi(\Omega)$ is.\\

\noindent If $\phi$ is a measure, it is non-decreasing and $\sigma$-sub-additive.\\

\noindent \textbf{(08.02c)} Let $\phi$ is $\sigma$-additive on $\mathcal{A}$ be a proper application on $\mathcal{A}$.\\

\noindent (A) If $\phi$ is $\sigma$-additive, it is continuous in the following sense :

 that $\phi$ is continuous in the following sense :

\noindent For any non-decreasing sequence $(A_n)_{n\geq 0}$, for any non-increasing sequence $(A_n)_{n\geq 0}$ for which there exists $n_0\geq 0$ such that $\phi(A_{n_0})$,

$$
\phi(\lim_{n\rightarrow +\infty} A_n)=\lim_{n\rightarrow +\infty} \phi(A_n). \ (LIM)
$$

\noindent \textbf{Warning}. The main difference with the continuity defined in \textit{Point (04.10) in Doc 04-01 (page \pageref{doc04-01}) of Chapter \ref{04_measures}} resides in the limit
$\lim_{n\rightarrow +\infty} \phi(A_n)$. For measures, this limit is monotone. In the general case, we know nothing about that.\\

\noindent \textbf{Terminology}. If Formula (LIM) holds for any  non-decreasing sequence $(A_n)_{n\geq 0}$, we say that $\phi$ is continuous below. If it holds for any non-increasing sequence $(A_n)_{n\geq 0}$ for which one of the $\phi(A_n)$ is finite, we say that $\phi$ is continuous above.\\

\noindent (B) Let $\phi$ be finitely additive and\\

\hskip 0.5cm either continuous from below,\\

\hskip 0.5cm or finite and continuous at $\emptyset$ (from above),\\

\noindent then $\phi$ is $\sigma$-additive.\\

\noindent \textbf{(08.03d)} Let $\phi$ be $\sigma$-additive on $\mathcal{A}$. There exists two measurable sets $C$ and $D$ such that $\phi(C)=\sup \phi$ and $\phi(D)=\inf \phi$.\\

\bigskip \noindent A proper $\sigma$-additive application may be called as a signed measure for the following Hahn-Jordan Theorem.\\

\noindent \textbf{(08.02)} (Hahn-Jordan Theorem). Let $\phi$ be a $\sigma$-additive application.\\

\noindent (A) There exists two measures $\phi_1$ and $\phi_2$ such 

$$
\phi=\phi_1-\phi_2.
$$

\bigskip \noindent (B) The measures $\phi_i$, $(i=1,2)$, are may taken as

$$
\phi_1(A)= \sup \{\phi(B), \ B \subset A, \ B\in \mathcal{A}\} \ and \ \phi_2(A)= -\inf \{\phi(B), \ B \subset A, \ B\in \mathcal{A}\}.
$$

\bigskip \noindent (C) Besides, there exists $D \in \mathcal{A}$ such that for all $A \in \mathcal{A}$

\begin{equation*}
\phi_{1}(A)=\phi(AD^c) \ and \ \phi_{2}(A)=- \phi (AD).
\end{equation*}

\newpage

\noindent \textbf{Part II : Indefinite integrals}.\\

\noindent \textbf{(08.04) Definition of indefinite integrals}. Let $(\Omega ,\mathcal{A},m)$ be a measure space. Let $X$ be a measurable $m$-quasi-integrable application $ X:(\Omega ,\mathcal{A},m)\mapsto \overline{\mathbb{R}}$, that is $\int X^{+} \ dm < +\infty$ or $\int X^{-} \ dm < +\infty$. For the needs of our exposition, we assume that $\int X^{-} \ dm < +\infty$.\\

\noindent The indefinite integral with respect $X$ and the measure $m$ is defined on $\mathcal{A}$ by 

$$
\mathcal{A} \ni A\mapsto \phi_{X} (A)=\int_{A}X \ dm. \ (I01).
$$

\bigskip \noindent \textbf{(08.05)} The Indefinite integral has the following properties :\\

\noindent \textbf{(08.05a)} $\phi_{X}$ is $\sigma$-additive.\\

\noindent \textbf{(08.05b)} $\phi_{X}$ is continuous with respect to $m$.\\

\noindent \textbf{(08.05c)} If $m$ is $\sigma$-finite and $X$ is \textit{a.e.} finite, $\phi_{X}$ is $\sigma$-finite.\\

\bigskip \noindent The reverses of such properties leads to the following Radon-Nikodym's theorem.
\newpage
\noindent \textbf{Part III : Radon-Nikodym's Theorem}.\\

\bigskip \noindent The Radon-Nikodym's Theorem is the reverse of Part II. The first result is referred to Lebesgue Decomposition Theorem.\\

\noindent \textbf{(08.06)} Let $\phi$ be a $\sigma$-additive application defined on $(\Omega, \mathbb{A})$ and let $m$ be a $\sigma$-finite measure on $(\Omega, \mathcal{A})$. There exist a $\sigma$-additive application $\phi_c$ defined on $\mathcal{A}$ and $m$-continuous, and a difference of finite $m$-singular and $\sigma$-additive application $\phi_s$ defined on $\mathcal{A}$  such that 
$$
\phi = \phi_c + \phi_s.
$$

\bigskip \noindent Besides, there exists an $m$-\textit{a.e.} finite and measurable application $X:(\Omega ,\mathcal{A})\mapsto \overline{\mathbb{R}}$ such that $\phi_c$ is the indefinite integral with respect $X$ and to the measure $m$, that is

$$
\forall A \in \mathcal{A}, \ \phi_{c} (A)=\int_{A} X \ dm, \ (I02).
$$

\bigskip \noindent \textbf{(08.07)} (Radon-Nikodym's Theorem).\\

\noindent \textbf{(08.07a)} ($\sigma$-finite version) Let $\phi$ be a $\sigma$-additive and $\sigma$-finite application defined on $(\Omega, \mathbb{A})$ and let $m$ be a $\sigma$-finite measure on
$(\Omega, \mathcal{A})$ such that $\phi$ is continuous with respect to $m$. Then there exists an $m$-\textit{a.e.} finite and measurable application $X:(\Omega ,\mathcal{A})\mapsto \overline{\mathbb{R}}$, $m$-\textit{a.e.} unique, such that 

$$
\forall A \in \mathcal{A}, \ \phi(A)=\int_{A}  X \ dm. \ (I03)
$$

\bigskip \noindent \textbf{Rephrasing}. Suppose  that the $\sigma$-additive and $\sigma$-finite application $\phi$ and the measure $m$ are defined on the same  $\sigma$-algebra. If $\phi$ is continuous with respect to measure $m$, then $\phi$ is an indefinite integral with respect to $m$. If $\phi$ is continuous with respect to measure $m$, then the singular part in Formula (I02) is the null application.

\bigskip \noindent  The function $X$ is called the Radon-Nikodym derivative of $\phi$ with respect to $m$ and denoted by  
\begin{equation*}
X=\frac{d\phi }{dm}
\end{equation*}

\bigskip \noindent Besides if $\phi$ is non-negative, then $X$ is $m$-\textit{a.e.} non-negative. \\

\noindent \textbf{(08.07c)} (Extended version) Let $\phi$ be a $\sigma$-additive application defined on $(\Omega, \mathbb{A})$, \textit{not necessarily $\sigma$-finite}, and let $m$ be $\sigma$-finite measure on $(\Omega, \mathcal{A})$ such that $\phi$ is continuous with respect to $m$. Then $\phi$ possesses a Radon-Nykodim derivative, $m$-\textit{a.e.} unique, not necessarily $m$-\textit{a.e.} finite.\\

\noindent \textbf{(08.07c)} (Additivity of the Radon-Nikodym derivative) Let $(\phi_i)_{i\in I}$, $I \subset \mathbb{N}$, be a countably family of $\sigma$-additive applications defined on $(\Omega, \mathbb{A})$, \textit{not necessarily $\sigma$-finite}, and let $m$ be $\sigma$-finite measure on $(\Omega, \mathcal{A})$ such that each $\phi_i$, $i\in I$, is continuous with respect to $m$. Suppose that one of the value $+\infty$ or $-\infty$ is simultaneously excluded for all the members to ensure that the application

$$
\phi=\sum_{i\in I} \phi_i
$$

\bigskip \noindent is well-defined.  The formula 
$$
\frac{d\phi}{dm}=\sum_{i \in I} \frac{d\phi_i}{dm}. \ \ (ARD)
$$  

\bigskip  \noindent is valid if $I$ is finite or the elements of the family are non-negative.\\

\noindent In the general, a condition for the validity of Formula (ARD) is

$$
\sum_{i\in I} \int_A 1_N |X_i| \ dm<+\infty,
$$

\bigskip  \noindent where

$$
N=\bigcap (X_i \text{ finite}).
$$

\bigskip \noindent Formula (ARD) still holds if $m(N)=0$ under the generated form of the equality with both sides infinite.\\

\bigskip \noindent We also have : for any $\sigma$-additive application $\phi$ defined on $(\Omega, \mathbb{A})$, \textit{not necessarily $\sigma$-finite}, and for any $\sigma$-finite measure on $(\Omega, \mathcal{A})$ $m$ be  such that $\phi$ is continuous with respect to $m$ and for any real and finite number $\lambda\neq 0$,

$$
\frac{\lambda \phi}{dm}=\lambda \frac{\lambda \phi}{dm}. (MRD)
$$

\noindent If $\phi$ is $\sigma$-finite, then Formula (MRD) is still true for $\lambda=0$.\\

\bigskip \noindent \textbf{(08.08) Terminology}. The expression of Radon-Nikodym derivative may be explained by the following formal writing. If $\phi$ is non-negative and is the indefinite integral with respect to $X\geq 0$, Formula (I03) may be written as
$$
\int h \ d\phi =\int h X dm  (RD) 
$$

\bigskip \noindent for $h=1_A$. By using the four step method in the construction of of the integral, we may see that Formula (RD) holds whenever both members are defined. Suppose that, by Borrowing the notation of differentiation, we could write 

$$
\int h \ d\phi = \equiv \int h \biggr(\frac{d\phi}{dm} \biggr) dm. \ (RD01)
$$

\bigskip \noindent By comparing Formulas (RD) and (RD02), it would be like we may make the replacement

$$
X=\frac{d\phi}{dm}
$$

\bigskip \noindent in Formula (RD). This explain why $X$ is called the Radon-Nikodym derivative of $phi$ with respect to $m$. We get the following rule.\\

\bigskip \noindent \textbf{Composition of derivatives}. Suppose that we have a finite number of measures $m_i$, $1\leq k \leq k$, $k\geq 20$, on $(\Omega, \mathcal{A})$ such that each $m_i$ is continuous with respect to $m_{i-1}$ for $i  \{2,...,k\}$. Then we have for any  $\mathcal{A}$-measurable and non-negative real-valued application $h$, we have

$$
\int h \ dm_k=\int h \biggr(\frac{dm_k}{dm_{k-1}}\biggr) \biggr(\frac{dm_k}{dm_{k-1}}\biggr)  \ldots \biggr(\frac{dm_2}{dm_{1}}\biggr) \ dm_1. \ (RDR)
$$

\noindent \textbf{(08.08)} (Extended Radon-Nikodym's Theorem).\\

\noindent Although Radon-Nikodym is frequently used with respect to a $\sigma$-finite measure, especially in Probability Theory, it holds without that assumption in the following statement :\\

\noindent (Extended Radon-Nikodym's Theorem).\\

\noindent Let $\phi$ be a $\sigma$-additive application defined on $(\Omega, \mathbb{A})$, \textit{not necessarily $\sigma$-finite}, and let $m$ be $\sigma$-finite measure on $(\Omega, \mathcal{A})$   such that $\phi$ is continuous with respect to $m$. Then there exists a measurable application $X:(\Omega ,\mathcal{A})\mapsto \overline{\mathbb{R}}$ \textit{not necessarily $m$-\textit{a.e.} finite} such that we have  

$$
\forall A \in \mathcal{A}, \ \phi(A)=\int_{A}  X \ dm.
$$

\noindent \LARGE \textbf{DOC 08-02 : $\sigma$-additive applications - Exercises} \label{doc08-02}\\
\bigskip
\Large

\noindent \textbf{Important remark}. In all these exercise, the function $\phi$ does never take the value $-\infty$, except in Exercise 5 where we allow excluding either $-\infty$ or $+\infty$.\\

\noindent \textbf{Exercise 1}. \label{exercise01_doc08-02}\\

\noindent Let $\phi$ be a proper application defined from $\mathcal{A}$ to $\mathcal{R}\cup {+\infty}$, which additive or $\sigma$-additive. Show that $\phi(\emptyset)=0$.\\

\noindent \textit{Hint}. Do nothing since the solution is already given in \textit{Exercise 1 in Doc 04-06 (page \pageref{exercise01_sol_doc04-06}) in Chapter \ref{04_measures}}.\\

\bigskip \noindent \textbf{Exercise 2}. \label{exercise02_doc08-02}\\

\noindent Let $\phi$ be $\sigma$-additive on $\mathcal{A}$. Let $(A_n)_{n\geq 0}$ be a sequence of mutually disjoint of elements of such that $\phi\left( \sum_{n\geq 0} A_n \right)$ is finite. Show that the series $\sum_{n\geq 0}  \phi\left( A_n \right)$ is absolutely convergent.\\

\noindent \textit{Hints}. Define

$$
I=\{n \in \mathbb{N} : \phi(A_n)\geq 0\} \ and \ J=\{n \in \mathbb{N} : \phi(A_n)<0\}
$$

\noindent and justify

\begin{eqnarray*}
\phi\left( \sum_{n\geq 0} A_n \right)&=&\phi\left( \sum_{n\in I} A_n + \sum_{n\in J} A_n \right) (L1)\\
&=&\phi\left( \sum_{n\in I} A_n \right) + \phi\left( \sum_{n\in J} A_n \right) (L2)\\
&=&\sum_{n\in I} \phi\left(  A_n \right) + \sum_{n\in J} \phi\left(  A_n \right) (L3)\\
&=:& S_1 + S_2.
\end{eqnarray*}

\noindent Why $S_2$ is finite from the fact that the value $-\infty$ is excluded for $\phi$? Combine this with the analysis of the signs of $S_1$ and $S_2$ to show that $S_1$ is also finite whenever $S_1+S_2$ is finite.\\

\noindent Express  $\sum_{n\geq 0} |\phi\left(  A_n \right)|$ as a function of $S_1$ and $S_2$ and conclude.\\

\bigskip \noindent \textbf{Exercise 3}. \label{exercise03_doc08-02}\\

\noindent Let $\phi$ be $\sigma$-additive on $\mathcal{A}$. If $A$ and $B$ are two elements of $\mathcal{A}$ such that $A \subset B$. Show that  $\phi(A)$ and
$\phi(B \setminus A)$ are finite whenever $\phi(B)$ and

$$
\phi(B \setminus A) =\phi(B) - \phi(A).
$$

\bigskip \noindent In particular, $\phi$ is finite if $\phi(\Omega)$ is.\\

\noindent \textbf{Hint}. Let $A$ and $B$ be two elements of $\mathcal{A}$ such that $A \subset B$ and $\phi(B)$ is finite. Use the relationship $\phi(B)=\phi(B\setminus A) + \phi(A)$. Consider two cases :
the case where one of them is non-positive and the case where both $\phi(B\setminus A)$ and $\phi(A)$ are both non-negative, and use again, when necessary, the fact that
the value $-\infty$ is excluded.\\

\noindent \textbf{Be careful}. $\phi(\Omega)$ is not necessarily a maximum value for $\phi$, while it is if $\phi$ is non-negative.\\

\bigskip \noindent \textbf{Exercise 4}. \label{exercise04_doc08-02}\\

\noindent Let $\phi$ is $\sigma$-additive on $\mathcal{A}$ be a proper application on $\mathcal{A}$.\\

\noindent Question (a) Let $\phi$ be $\sigma$-additive. Show that $\phi$ is continuous.\\

\noindent Question (b) Let $\phi$ be finitely additive and

\hskip 1cm (b1) either continuous from below

\hskip 1cm (b2) or finite and continuous at $\emptyset$ (from above).

\noindent Show that $\phi$ is $\sigma$-additive.\\

\noindent \textit{Hints}.\\

\noindent \textbf{Hints for Question (a)}. Adapt the solutions of Exercise 2 in Doc 04-06 (page \pageref{exercise02_sol_doc04-06}) in Chapter \ref{04_measures} in the following way.\\

\noindent Let $(A_n)_{n\geq 0}$ be a non-decreasing sequences in $\mathcal{A}$. Write 

$$
A=A_0 + \sum_{k\geq 1} (A_k \setminus A_{k-1}).
$$

\noindent and apply $\sigma$-additivity. Now discuss over two cases. \\

\noindent Case 1 : There exists $n_0\geq 0$ such that $\phi(A_n)=+\infty$ for all $n\geq n_0$. Use Exercise 3 to justify that Formula

$$
\phi(\lim_{n\rightarrow +\infty} A_n)=\lim_{n\rightarrow +\infty} \phi(A_n). \ (LIM)
$$

\noindent holds since both members are $+\infty$.\\

\noindent Case 2 : $\forall n\geq 0$, \ $\exists N\geq n$, \ $|\phi(A_n)|<+\infty$. Use Exercise 3 to see that all the $\phi(A_n)$'s are finite and apply the techniques of Chapter \ref{04_measures} :

\begin{eqnarray*}
&&\phi(A_0) + \sum_{k=1}^{+\infty} \phi(A_k \setminus A_{k-1})\\
&=&\lim_{n\rightarrow +\infty} \phi(A_0) + \sum_{k=1}^{n} \phi(A_k) - \phi(A_{k-1})\\
&=& \lim_{n\rightarrow +\infty} \phi(A_0) + (\phi(A_1)-\phi(A_0)) + (\phi(A_2)-\phi(A_1))+\cdots +(\phi(A_n)-\phi(A_{n-1}))\\
&=&\lim_{n\rightarrow +\infty} \phi(A_n).
\end{eqnarray*}

\noindent Next, let $(A_n)_{n\geq 0}$ be any non-increasing sequence for which there exists $n_0\geq 0$ such that $\phi(A_{n_0})$. Apply the results to $A_{n_0}\setminus A_n$ and conclude.\\

\noindent \textbf{Hints for Question (b1)}. Consider a sequence $(A_n)_{n\geq 0}$ of pairwise disjoint elements of $\mathcal{A}$. Use the formula

$$
\sum_{n=0}^{k} A_n \nearrow \sum_{n=0}^{+\infty} A_n 
$$ 

\noindent and the hypotheses to conclude.\

\noindent \textbf{Hints for Question (b2)}. Consider a sequence $(A_n)_{n\geq 0}$ of pairwise disjoint elements of $\mathcal{A}$. Use the formula

$$
\sum_{n=0}^{+\infty} A_n = \sum_{n=0}^{k} A_n  + \sum_{n=k+1}^{+\infty} A_n. 
$$ 

\noindent Combine the fact that $B_{k}=\sum_{n=k+1}^{+\infty} A_n \downarrow \emptyset$ with the hypotheses and conclude.\\

\bigskip \noindent \textbf{Exercise 5}. \label{exercise05_doc06-04}\\

\noindent Let $\phi$ be $\sigma$-additive on $\mathcal{A}$ which does not take either the value $-\infty$ or the value $+\infty$. Show that there exists two measurable sets $C$ and $D$ such that $\phi(C)=\sup \phi$ and $\phi(D)=\inf \phi$.\\

\noindent \textit{Hints}.\\

\noindent \textbf{Existence of C}.\\

\noindent Suppose that there is a measurable set $C$ such that $\phi(C)=+\infty$. Conclude.\\

\noindent Suppose that $\phi \leq 0$. Take $C=\emptyset$ and conclude.\\

\noindent Now suppose that $\phi(A)<+\infty$ for all finite for all $A \in \mathcal{A}$ and that there exists $A_0\in \mathcal{A}$ such that $\phi(A_0)>0$.\\

\noindent (1) Justify the existence of a sequence $(A_n)_{n\geq 1} \subset \mathcal{A}$ such that the $\phi(A_n)$'s are finite for large values on $n$ and  $\left(\phi(A_n)\right)_{n\geq 1}$ converges to $\alpha=\sup \phi$ in $\overline{\mathbb{R}}$.\\

\noindent (2) Put $A=\bigcup_{n\geq 1} A_n$. For each $n\geq 1$, consider the partitions $A=A_k + (A \setminus A_k)$ for $k=1,...,n$ and the the intersection $A$ $n$ times and conclude that $A$ may is partitioned into sets of the form $A_1^{\prime}A_2^{\prime}...A_n^{\prime}$ where $A_k^{\prime}=A_k$ or $A_k^{\prime}=A\setminus A_k$ for $k=1,...,n$.\\

\noindent What is the size $\ell_n$ of the partition of $A$ for $v\geq 1$ fixed? Denote the elements of that partition by $A_{n,j}$, $j=1,...,\ell_n$.\\

\noindent (3) For each $n\geq 1$, we have $A_n \subset A=\sum_{1\leq j \leq \ell_n} A_{n,j}$. By doing $A_n=A\cap A_n= \sum_{1\leq j \leq \ell_n} A\cap A_{n,j}$, can you say that $A_n$ is a sum of elements
of some elements $A_{n,j}$, $j=1,...,\ell_n$. What are these elements?\\

\noindent (4) As well, let $n <p$. By take the intersection of the partitions $A=A_k + (A \setminus A_k)$, $k=n+1,...,p$, $(p-n)$ times, make a conclusion similar to the one given in (2) relative a certain partition $\mathcal{P}_{n,p}$ to be precised.\\

\noindent Now Take any $A_{n,j}$, $k=1,...,\ell_n$. By remarking  that $A_{n,j}$ is a subset of $A$, and by forming the intersection $A_{n,j}=A\cap A_{n,j}$ using the partition $\mathcal{P}_{n,p}$, show that $A_{n,j}$ is some of elements $A_{p,h}$, $h=1,...,\ell_p$.\\

\noindent (5) Define
$$
I_n={j, \ 1\leq j \leq \ell_n, \ and \ \phi(A_{n,j})\geq 0 }, \ J_n={j, \ 1\leq j \leq \ell_n, \ \phi(A_{n,j})< 0 },  
$$
 
\noindent where $I_n+J_n=\{1,...,\ell_n\}$, and

$$
B_n=\emptyset \ if \ I_n=\emptyset, \ and  \ B_n=\sum_{j \in I_n} A_{n,j} \ if \ I_n\neq \emptyset.
$$

\noindent Define also, for $n\geq 1$ fixed, $\mathcal{S}_n$ as the collection of finite sums of sets $A_{n,j}$, $j\in I_n$ and $\mathcal{D}_n$ as the collection of finite sums of sets $A_{n,j}$, $j\in J_n$.

\noindent Remark that, by Exercise 14 (see page \pageref{exercise14_sol_doc01-03}) in Doc 01-03 of Chapter \ref{01_setsmes}, $\mathcal{S}_n$ and $\mathcal{D}_n$ are stable by finite unions. Try it for the collection of finite sums of the mutually disjoint sums of, say, $A$, $B$, $C$.\\

\noindent Show that for $n\geq 1$ fixed, $B_n \in \mathcal{S}_n$, $\phi$ is non-decreasing on $\mathcal{S}_n$ and non-increasing on $\mathcal{S}_n$, and that if $E\in \mathcal{S}_n$ and $F\in \mathcal{D}_n$ with $E \ cap EF=\emptyset$, use the non-positivity of $\phi(F)$ and the facts above to show that

$$
\phi(E+F) \leq \phi(E).
$$

\noindent (6a) From now, $n \geq $ is fixed. From (3), explain why $A_n$ is of the form  $A_n=E+F$ with $E\in \mathcal{S}_n$ and $F\in \mathcal{D}_n$. Use the non-decreasingness of $\phi$ on 
$\mathcal{S}_n$ to get

$$
\phi(A_n) \leq \phi(B_n).
$$

\bigskip \noindent (6b) You are going to prove the for any $p\geq n$, we have

$$
\phi\left( B_{n} \bigcup B_{n+1} \bigcup ... \bigcup B_{p} \right) \leq \phi\left( B_{n} \bigcup B_{n+1} \bigcup ... \bigcup B_{p} \bigcup \bigcup B_{p+1}\right). \ (AC)
$$

\noindent Fix $r \in\{n,...,p\}$. By (4), $B_r$ is a sum of elements $\mathcal{S}_{p+1}$ and we may write $B_r=E_{r,1}+E_{r,2}$ with $E_{r,1} \in \mathcal{S}_{n+1}$ and $E_{r,2}\in \mathcal{D}_{n+1}$. We have $E_{r,1} \subset B_{p+1}$ and $E_{r,2}$ is disjoint of $B_{p+1}$. Validate the computations : 

\begin{eqnarray*}
E&=:& B_{n} \bigcup B_{n+1} \bigcup ... \bigcup B_{p} \bigcup B_{p+1}\\
&=& \bigcup_{1\leq r \leq p} \left(B_r \bigcup B_{p+1}\right)\\
&=& \bigcup_{1\leq r \leq p} \left(E_{r,1} \bigcup E_{r,2}\right)\\
&=& \left(\bigcup_{1\leq r \leq p} E_{r,1} \right) \bigcup \left(\bigcup_{1\leq r \leq p} E_{r,2}\right).
\end{eqnarray*}

\noindent By exploiting that $\mathcal{S}_{p+1}$ and $\mathcal{D}_{p+1}$ are stable by finite unions, conclude that 

\begin{eqnarray*}
E&=:&B_{n} \bigcup B_{n+1} \bigcup ... \bigcup B_{p} \bigcup B_{p+1}\\
&=& \left(\bigcup_{1\leq r \leq p} E_{r,1} \right) \bigcup \left(\bigcup_{1\leq r \leq p} E_{r,2}\right)\\
&=:&E_1 + E_2,
\end{eqnarray*}

\noindent with $E_1 \in \mathcal{S}_{p+1}$ and $E_2 \in \mathcal{D}_{p+1}$, and then, $E_1 \subset B_{p+1}$ and $E_2$ is disjoint of $B_{p+1}$. Validate this,  

$$
E \cup B_{p+1}= E_1 \cup B_{p+1} \cup E_2=B_{^p+1} \cup E_2=B_{p+1} + E_2
$$

\noindent and, by using again the non-decreasingness of $\phi$ on $\mathcal{S}_{p+1}$, show that
 
$$
\phi(E)=\phi(E_1)+\phi(E_2) \leq \phi(B_{p+1})+\phi(E_2)=\phi(E \cup B_{p+1}).
$$

\noindent (7) Combine (5) and (6) and repeat Formula (A) as needed as possible to get, for all $n\geq n_0$, for all $p>n$;

$$
\phi(A_n) \leq \phi(B_n) \leq \phi\left( \bigcup_{n\leq r \leq p} B_r\right). 
$$

\noindent Conclude and find $C$ by letting first  $p\uparrow +\infty$ and next $n\uparrow +\infty$.\\

\noindent \textbf{Existence of D}. Consider $-\phi$ and apply the first part to get a set $D\in \mathcal{A}$ such that

$$
-\phi(D)=  \sup - \phi \ \ \Leftrightarrow \phi(D)=  \inf \phi. 
$$

\bigskip \noindent \textbf{Exercise 6}. \label{exercise06_doc08-02}\\

\noindent Let $\phi$ be a $\sigma$-additive application from $\mathcal{A}$ to $\mathbb{R}\cup \{+\infty\}$. Define for all $A \in \mathcal{A}$,

$$
\phi_1(A)= \sup \{\phi(B), \ B \subset A, \ B\in \mathcal{A}\} \ and \ \phi_2(A)= -\inf \{\phi(B), \ B \subset A, \ B\in \mathcal{A}\}.\\
$$

\noindent Consider the measurable space $D$, found in Exercise 5, such that $\phi(D)=\inf \phi$.\\

\noindent Question (a). If $\phi$ is a measure, check that : $\phi_1=\phi$, $\phi_2=0$,  $\phi_1(A)=\phi(AD^c)$ and $\phi_2=\phi(AD)$ for $D=\emptyset$.\\

\noindent Question (b). In the sequel, we suppose that there exists $A_0\in \mathcal{A}$ such that $\phi(A_0)<0$.\\

Show that for all $A\in \mathcal{A}$, $\phi(AD)\leq 0$ and $\phi(AD^)\geq 0$.\\

\noindent \textit{Hints}. Show that $\phi(D)$ is finite and by using Exercise 3, that $\phi(AD)$. By using meaning of the infimum, decompose $\phi(D \setminus AD)$ and show we have a contradiction if we don't have $\phi(AD)\leq 0$. Do the same for $\phi(D + AD^c)$ to show that $\phi(AD^)\geq 0$.\\

\noindent Question (c).\\

\noindent (c1) Justify each line of the formula below : For $A \in \mathcal{A}$, for all $B\in \mathcal{A}$ with $B \subset A$, we have

\begin{eqnarray*}
\phi(B)&=&\phi(BD)+\phi(BD^c)\\
&\leq& \phi(BD^c)\\
&\leq& \phi(BD^c) +\phi((A\setminus B)D^)\\
&=&\phi(AD^c).
\end{eqnarray*}

\noindent Deduce from this that : $\phi_1(A)=\phi(AD^c)$ for all $A\in \mathcal{A}$.\\

\noindent (c2) Justify each line of the formula below : For $A \in \mathcal{A}$, for all $B\in \mathcal{A}$ with $B \subset A$, we have

\begin{eqnarray*}
\phi(B)&=&\phi(BD)+\phi(BD^c)\\
&\geq& \phi(BD)\\
&\geq& \phi(BD) +\phi((A\setminus B)D) \\
&=&\phi(AD).
\end{eqnarray*}

\noindent Deduce from this that : $\phi_2(A)=-\phi(AD)$ for all $A\in \mathcal{A}$.\\

\noindent Question (d) Check (by two simple sentences) that $\phi_1$ and $\phi_2$ are measures. Is the application $\phi_1-\phi_2$ well-defined? Why?\\

\noindent Question (e) Conclude that $\phi=\phi_1-\phi_2$.\\

\noindent \emph{NB}. If the value $+\infty$ is excluded in place of $-\infty$, we establish the relation for $-\phi$ and adapt the results.\\

\noindent \LARGE \textbf{DOC 08-03 : Radon-Nikodym's Theorem - Exercises} \label{doc08-03}\\
\Large

%\bigskip \noindent \textbf{Exercise X}. \label{exercise0X_doc08-03}\\
%\bigskip \noindent \textbf{SOLUTIONS}.\\

\bigskip \noindent \textbf{NB}. Throughout the document, unless the contrary is specifically notified, $(\Omega,\mathcal{A},m)$ is a measure space $\phi$, $\phi_i$, $\psi$  and $\psi-i$ (i=1,2) are real-valued applications defined on $\mathcal{A}$ and $X : (\Omega,\mathcal{A}) \rightarrow \overline{\mathbb{R}}$ is a function.\\

\bigskip \noindent \textbf{Exercise 1}. \label{exercise01_doc08-03}\\

\noindent Let $(\Omega,\mathcal{A},m)$ be a measure space and $X : (\Omega,\mathcal{A}) \rightarrow \overline{\mathbb{R}}$ be a function.\\

\noindent Question (a) Let $X$ be non-negative ($X\geq 0$). Justify the definition

$$
\mathcal{A} \ni A \mapsto \phi_X(A)=\int_A X \ dm, \ (II)
$$ 

\bigskip 
\noindent and show, by using the \textit{MCT} (See Exercise 1, Doc 06-06, page \pageref{exercise01_doc06-06}) that $phi_X$ is $\sigma$-additive.\\ 

\noindent Question (b) Let $X^-$ be integrable or $X^+$ be integrable. Apply the the results of Question (a) to $X^-$ and $X^+$ and extend the definition \textit{(II)} by using
by first formula in Question (c), Exercise 6 in Doc 05-02 (page  \pageref{exercise06_doc05-02}) in Doc 05-02 in Chapter \ref{05_integration}) : for any real-valued function $f$, for any subset $A$ of $\Omega$, 
$$
(1_Af)^-=1_Af^- \ and \ (1_Af)^+=1_Af^+
$$

\bigskip \noindent and say why it is still $\sigma$-additive (Here, treat the case $X^-$ only).\\

\noindent Question (c) From now on, we suppose that $X^-$ is integrable or $X^+$ is integrable.\\

\noindent Show that $\phi_X$ is $m$-continuous.\\

\textit{Hints}. Use the first formula in Question (c), Exercise 6 in Doc 05-02 (page  \pageref{exercise06_doc05-02}) in Doc 05-02 in Chapter \ref{05_integration}) and Property (P3) for non-negative functions in Exercise 5 Doc 05-02 (page \pageref{doc05-02}) in Chapter \ref{05_integration}.\\

\noindent Question (d) Suppose that $m$ is $\sigma$-finite. Show that is $phi_X$ is $\sigma$-finite.\\

\noindent \textit{Hints}. Use a countable and measurable subdivision of $\Omega$ : $\Omega=\sum_{j\geq 0}$, $(\forall j\geq 0, \ \Omega_j \in \mathcal{A} \text{ and } m(\Omega_j)<+\infty)$ and the formula

$$
(f \ finite)=\sum_{n\in \mathcal{Z}} (n\leq X <n+1).
$$

\bigskip 
\noindent Show that $\phi_X$ is finite on each $A_{n,j}=(n\leq X <n+1)\bigcap \Omega_j$, $(n\in \mathcal{Z}, \ j\geq 0)$. Deduce from this that $\phi_X$ is $\sigma$-finite.

\noindent Question (e)  Readily check that $\phi_X$ is finite if $X$ is integrable.

\noindent Question (d)  Suppose that $X$ and $Y$ are both $m$-finite measurable and quasi-integrable. Show that $\phi_Y$ are  equal if and only if $X=Y$ $m$-a.e.\\

\bigskip \noindent \textbf{Exercise 2}. \label{exercise02_doc08-03}\\

\noindent Let $\phi$ and $\psi$ be two $\sigma$-additive applications which compatible in the sense that the value $-\infty$ is excluded for both of them or the value $+\infty$ is excluded for both of them, so that the application $\Phi=\phi+\psi$ is well-defined and is $\sigma$-finite.\\

\noindent Question (a) Show that if $\phi$ and $\psi$ are both $m$-continuous, then $\Phi$ is $m$-continuous.\\

\noindent Question (b) Show that if $\phi$ and $\psi$ are both $m$-singular \textbf{measures}, then $\Phi$ is $m$-continuous.\\

\noindent Question (c) Show that if $\phi$ is both $m$-continuous and  $m$-singular, then $\phi$ is the null application.\\

\bigskip \noindent \textbf{Exercise 3}. \label{exercise03_doc08-03}\\

\noindent Suppose that the measure $m$ is finite. Let $\phi$ another finite measure. Show that there exist a \textbf{finite} $m$-continuous measure $\phi_c$ and a \textbf{finite} $m$-singular measure $\phi_s$ such that $\phi=\phi_c+\phi_s$. Show that $\phi_c$ is of the form

$$
\phi_c(A)=\int_A X \ dm\leq \phi(A), \ A \in \mathcal{A},
$$
 
\bigskip 
\noindent and, hence, $X$ is $m$-\textit{a.e.}.\\

\noindent \textit{Hints}. Proceed as follows.\\

\noindent Question (1). Define $\Phi$ as the set of all Borel and non-negative applications $X : (\Omega,\mathcal{A}) \rightarrow \overline{\mathbb{R}}$ such that for all $A\in \mathcal{A}$, we have 
$$
\int_A X \ dm\leq \phi(A)
$$

\bigskip 
\noindent (a) Justify the existence and the finiteness of $\alpha=\sup \{\int_A X \ dm, \ X \in \Phi\} \in \mathbb{R}$, that find a least one element of the $\{\int_A X \ dm, \ X \in \Phi\}$ and exhibit a finite bound above.\\

\noindent (b) Consider a maximizing sequences $(X_n)_{n\geq 0} \subset \Phi$ such that $\int_A X \ dm \rightarrow \alpha$ as $n\rightarrow +\infty$. Justify that $Y_n=\max_{1\leq k\leq n} X_k$ has a limit $X$, a non-negative Borel application.\\

\noindent (c) Do you have $X \in \phi?$ Proceed as follows : For $n\geq 0$, denote $A_k=(Y_n=X_k)$ for $1\leq k\leq n$. Show that

$$
\Omega = \bigcap_{1\leq k \leq n} A_k.
$$

\bigskip 
\noindent Next, use the now well known technique by which one transforms the union of the $A_k$'s into a sum of sets (See Exercise 4, page \pageref{exercise03_sol_doc00-03}, in Doc 00-03 in Chapter \ref{00_sets} by taking

$$
B_1=A_1, \ B_2=A_1+...+A_{k-1}^c A_k
$$

\bigskip 
\noindent to get

$$
\Omega = \bigcap_{1\leq k \leq n} A_k = \sum_{1\leq k \leq n} B_k,
$$

\bigskip 
\noindent from which, you have for all $A\in \mathcal{A}$, is derived

$$
A=A\cap \Omega = \sum_{1\leq k \leq n} AB_k, \text{ and } Y_n=X_k \text{ on } B_k, \ 1\leq k \leq n. 
$$

\bigskip 
\noindent Question : Form the following formula which is valid for all $A\in \mathcal{A}$,

$$
\int_A Y_n \ dm = \int_{sum_{1\leq k \leq n} AB_k} Y_n \ dm.
$$

\bigskip 
\noindent Show that  $\int_A Y_n \ dm\leq \phi(A)$ for all $n\geq 0$ and deduce from this, that you have  $\int_A X \ dm\leq \phi(A)$. Conclude.

\noindent Conclude that

$$
\alpha = \int X \ dm.
$$

\bigskip 
\noindent Question (2). Take the indefinite integral with respect to $X$

$$
\mathcal{A} \ni A \mapsto \phi_c(A)=\int_A X\ dm
$$

\bigskip 
\noindent and define, based on the fact that $\phi$ is finite,

$$
\phi_s=\phi - \phi_c.
$$

\bigskip 
\noindent (a) Use the conclusion of Question (1) to show that $\phi_s$ is non-negative and thus is a measure.\\

\noindent (b) Show that $\phi_s$ is $m$-singular by using Exercise 1 and proceeding as follows.\\

\noindent For each $n\geq 1$, consider the $\sigma$-finite application, well-defined since $m$ is finite, 
$$
\psi_n=\phi_s - \frac{1}{n}m.
$$

\bigskip 
\noindent Denote, following Exercise 1, for each $n\geq 1$ the measurable set $D_n$ such $\psi_n(D_n)=\inf \phi_n$ which satisfies : for all $A\in \mathcal{A}$,

$$
\psi_n(AD_n)\leq 0 \ (C1) \text{  and  } \psi_n(AD_n^c)\geq 0 \ (C2).
$$

\bigskip 
\noindent Denote

$$
N=\bigcap_{n\geq 1} D_n^c \Leftrightarrow N^c=\bigcup_{n\geq 1} D_n.
$$

\bigskip 
\noindent Questions : \\

\noindent (A) From Part (C1) of the last equation and from the non-negativity of $\phi_s$, derive a bound of $\phi_s(AD_n)$ for each $n$ and from this, show that for all $A\in \mathcal{A}$,

$$
\phi_s(AN^c)=0.
$$

\bigskip 
\noindent (B) Justify each line of the following formula for $n\geq 1$ fixed and for all $A\in \mathcal{A}$ :

\begin{eqnarray*}
\int_A \biggr( X + \frac{1}{n}1_{AD_n^c}\biggr) \ dm &=&\phi_c(A)+\frac{1}{n}m(AD_n^c)\\
&=&\phi(A)-\phi_s(A)+\frac{1}{n}m(AD_n^c)\\
&=&\phi(A)-\biggr( \phi_s(AD_n^)+\frac{1}{n}m(AD_n^c) \biggr) +\biggr(\phi_s(AD_n^c)-\phi_s(A)\biggr)\\
&=&\phi(A)- \psi_n(AD_n^c) + \biggr(\phi_s(AD_n^c)-\phi_s(A)\biggr)\\
&\leq &\phi(A).
 \end{eqnarray*}

\noindent Deduce from this that 

$$
X + \frac{1}{n}1_{AD_n^c} \in \Phi,
$$

\bigskip 
\noindent and from this, deduce that for all $n\geq 1$, $m(D_n^c)=0$. Conclude that $m(N)=0$.\\

\noindent Conclude that $\phi_s$ is $m$-singular.\\

\noindent Make a general conclusion.\\

\bigskip \noindent \textbf{Exercise 4}. \label{exercise04_doc08-03}\\

\noindent Generalize the results of Exercise 3 to $\sigma$-finite measures $\phi$ and $m$ finite in the following way :  There exists $m$-\textit{a.e.} finite non-negative application $X$ and a \textbf{finite } $m$-singular measure $\phi_s$ such that for all $A \in \mathcal{A}$,

$$
\phi(A)=\int_A X \ dm + \phi_s(A).
$$

\bigskip 
\noindent \textit{Recommendations}. It is advised in a first reading, just to follow the solution of this extension. The ideas are the following : We consider a measurable partition of $\Omega$ : $\Omega=\sum_{j\geq 0}$, $(\forall j\geq 0, \ \Omega_j \in \mathcal{A})$ such that $\phi(\Omega_j)<+\infty$ and $m(\Omega_j)<+\infty$. Hence, we have for all $A \in \mathcal{A}$

$$
\phi = \sum_{j \geq 0} \phi_j \ (S1), \ \text{and} \ m = \sum_{j \geq 0} m_j, \ (S2)
$$

\bigskip 
\noindent where the applications

$$
\mathcal{B} \ni B \mapsto \phi_j(A)=\phi(A\cap \Omega_j) \text{ and } \mathcal{B} \ni B \mapsto m_j(A)=m(A\cap \Omega_j)
$$

\bigskip 
\noindent are measures which can be considered as defined on $\mathcal{A}$ (but with $\Omega_j$ as support) or, on the $(\Omega_j, \mathcal{A}_j)$, where $\mathcal{A}_j)$ is the induced 
$\sigma$-algebra of $\mathcal{A}$ on $\Omega_j$.\\

\noindent Hence on each $\Omega_j$, $j\geq 1$, $\phi_j$ and $m_j$ are finite and we get, by Exercise 3, a decomposition

$$
\phi_j(B)=\int_B X_j \ dm + \phi_{s,j}(B), \ B \in \mathcal{A}_j, \ (DE01)
$$

\bigskip 
\noindent From there, we sum members of Formula (DE01). But some technicalities are needed.\\

\bigskip \noindent \textbf{Exercise 5}. (General form of the Lebesgue Decomposition) \label{exercise06_doc08-03}\\

\noindent Generalize the results of Exercise 4 to a $\sigma$-finite and $\sigma$-additive application $\phi$ and a $\sigma$-finite measures $m$.\\

\noindent Question (a) Show that there exists an $m$-\textit{a.e.} finite application $X$ and a difference of \textbf{finite} $m$-singular measures $\phi_s$ such that for all $A \in \mathcal{A}$,

$$
\phi(A)=\int_A X \ dm + \phi_s(A).
$$

\bigskip 
\noindent Question (b) Show that the decomposition is unique if $\phi_s(A)$ is a difference of $m$-singular measures.\\

\noindent \textit{Hints} Combine Exercises 1  and 5.\\

\bigskip \noindent \textbf{Exercise 6}. (Theorem of Radon-Nikodym) \label{exercise07_doc08-03}\\

\noindent Let $\phi$ be a $\sigma$-finite and $\sigma$-additive application and a $\sigma$-finite measures $m$. Suppose that $\phi$ is $m$-continuous. Then there exists a real-valued measurable $X$ wich is $m$-finite such that for $A \in \mathcal{A}$, we have

$$
\phi(A)=\int_A X \ dm.
$$

\bigskip 
\noindent Show that $X$ is $m$-\textit{a.e.} non-negative if $\phi$ is a measure.\\

\bigskip \noindent \textbf{Exercise 7}. (Extended Theorem of Radon-Nikodym) \label{exercise07_doc08-03}\\

\noindent Let $\phi$ be a $\sigma$-additive application defined on $(\Omega, \mathbb{A})$, \textit{not necessarily $\sigma$-finite}, and let $m$ be $\sigma$-finite measure on $(\Omega, \mathcal{A})$   such that $\phi$ is continuous with respect to $m$. Show that there exists a measurable application $X:(\Omega ,\mathcal{A})\mapsto \overline{\mathbb{R}}$ \textit{not necessarily $m$-\textit{a.e.} finite} such that we have  

$$
\forall A \in \mathcal{A}, \ \phi(A)=\int_{A}  X \ dm.
$$

\bigskip 
\noindent \textit{Hints}. Show it only for $\phi$ is non-negative and $m$ finite. First for $\phi$ non-negative, $X_j$ will be found on a countable partition on $\Omega$ piece-wisely as in the solution of Exercise 6. For a general $\sigma$-additive application, the extension is done by using the difference of two measures, one of them being finite. Based on these remarks, proceed as follows.\\

\noindent So $\phi$ be a $m$-continuous measure, where $m$ is a finite measure.\\ 

\noindent Question (a). Consider the class of elements of $A \in \mathcal{A}_0$ such that $\phi$ is $\sigma$-finte of the induce $\sigma$-algebra $\mathcal{A}_A$ and a sequence $(B_n)_{n\geq 0} \subset \mathcal{A}_0$ such that $m(B_n)\rightarrow  s=\sup_{B\in \mathcal{A}_0} m(B)$. Define $\bigcup_{k=0}^{+\infty} B_k$ and define also $\mathcal{A}_1=\{AB, \ A\in \mathcal{A}\}$ and $\mathcal{A}_2=\{AB^c, \ A\in \mathcal{A}\}$.\\

\noindent Show that for all $n\geq 0$, $\bigcup_{k=0}^{n} B_k \in \mathcal{A}_0$ and deduce that $B \in \mathcal{A}_0$.\\

\noindent Apply the MCT continuity from above of $m$  to $\bigcup_{k=0}^{n} B_k \nearrow B$ as $n\nearrow +\infty$ to show that $s=m(B)$.\\

\noindent Question (b)  Show that if for some $C\in \mathcal{A}_2$ we have $0<\phi(C)<+\infty$, then $s$ would not be an infimum. Conclude that we have the equality 

$$
(C \in \mathcal{A}_2, \ \phi(C) \text{ finite }) \Rightarrow m(C)=0. (E1)
$$

\bigskip 
\noindent Deduce from this that $\phi$ is the indefinite integral of $X_2=+\infty$ of $\mathcal{A}_2$.\\

\noindent Use the convention of the integration on the class of elementary function that $\int_A \infty \ dm=m(A) \times \infty=0$ if $m(A)=0$.\\ 

\noindent Question (c) By remarking that $\phi$ is $\sigma$-finite on  $\mathcal{A}_1$, apply the first version of Radon-Nikodym with an  derivative $X_1$ on $B$.\\

\noindent Question (d) Form $X=X_1 1_B + X_2 1_{B^c}$ and conclude.\\

\noindent \LARGE \textbf{DOC 08-04 : $\sigma$-additive applications - Exercises/Solutions} \label{doc08-04}\\
\bigskip
\Large

\noindent \textbf{Important remark}. In all these exercise, the function $\phi$ does never take the value $-\infty$, except in Exercise 5 where we allow excluding either $-\infty$ or $+\infty$.\\

\noindent \textbf{Exercise 1}. \label{exercise01_sol_doc08-04}\\

\noindent Let $\phi$ be a proper application defined from $\mathcal{A}$ to $\mathcal{R}\cup {+\infty}$, which additive or $\sigma$-additive. Show that $\phi(\emptyset)=0$.\\

\noindent \textit{Hint}. Do nothing since the solution is already given in \textit{Exercise 1 in Doc 04-06 (page \pageref{exercise01_sol_doc04-06}) in Chapter \ref{04_measures}}.\\

\bigskip \noindent \textbf{SOLUTIONS}. The solution is already given in \textit{Exercise 1 in Doc 04-06 (page \pageref{exercise01_sol_doc04-06}) in Chapter \ref{04_measures}}, where the non-negativity of the proper $\sigma$-additive application $m$ is not used at all.\\

\bigskip \noindent \textbf{Exercise 2}. \label{exercise02_sol_doc08-04}\\

\noindent Let $\phi$ be $\sigma$-additive on $\mathcal{A}$. Let $(A_n)_{n\geq 0}$ be a sequence of mutually disjoint of elements of such that $\phi\left( \sum_{n\geq 0} A_n \right)$ is finite. Show that the series $\sum_{n\geq 0}  \phi\left( A_n \right)$ is absolutely convergent.\\

\noindent \textit{Hints}. Define

$$
I=\{n \in \mathbb{N} : \phi(A_n)\geq 0\} \ and \ J=\{n \in \mathbb{N} : \phi(A_n)<0\}
$$

\noindent and justify

\begin{eqnarray*}
\phi\left( \sum_{n\geq 0} A_n \right)&=&\phi\left( \sum_{n\in I} A_n + \sum_{n\in J} A_n \right) (L1)\\
&=&\phi\left( \sum_{n\in I} A_n \right) + \phi\left( \sum_{n\in J} A_n \right) (L2)\\
&=&\sum_{n\in I} \phi\left(  A_n \right) + \sum_{n\in J} \phi\left(  A_n \right) (L3)\\
&=:& S_1 + S_2.
\end{eqnarray*}

\noindent Why $S_2$ is finite from the fact that the value $-\infty$ is excluded for $\phi$? Combine this with the analysis of the signs of $S_1$ and $S_2$ to show that $S_1$ is also finite whenever $S_1+S_2$ is finite.\\

\noindent Express  $\sum_{n\geq 0} |\phi\left(  A_n \right)|$ as a function of $S_1$ and $S_2$ and conclude.\\

\bigskip \noindent \textbf{SOLUTIONS}.\\

\noindent Formula (L1) derives from the obvious relationship  $\sum_{n\geq 0} A_n =  \sum_{n\in I} A_n + \sum_{n\in J} A_n$ and the additivity of $\phi$. Formula (L2) derives from (L2) by
$\sigma$-additivity. Now, we see that $S_1$ is a series of non-negative terms and $S_2$ is a series of negative terms so that $S_1$ is non-negative and $S_2$ is negative.\\

\noindent By hypothesis $S_1 + S_2$ is finite. Besides, since the value $-\infty$ is excluded, the negative value $S_2=\phi\left( \sum_{n\in J} A_n \right)$ cannot be $-\infty$. Hence $S_2$ is finite. This implies that $S_1$ is finite, otherwise we would have $S_1=+\infty$ and $S_1+S2=+\infty$. In conclusion, $S_1$ and $S_2$ are finite.\\

\noindent But, we also have

\begin{eqnarray*}
\sum_{n\geq 0} |\phi\left(  A_n \right)|&=& \sum_{n\in I} |\phi\left(  A_n \right)| + \sum_{n\in J} |\phi\left(  A_n \right)| (L4)\\
&=:& S_1 - S_2,
\end{eqnarray*}

\noindent which is finite since both $S_1$ and $S_2$ are finite.\\

\bigskip \noindent \textbf{Exercise 3}. \label{exercise03_sol_doc08-04}\\

\noindent Let $\phi$ be $\sigma$-additive on $\mathcal{A}$. If $A$ and $B$ are two elements of $\mathcal{A}$ such that $A \subset B$. Show that  $\phi(A)$ and
$\phi(B \setminus A)$ are finite whenever $\phi(B)$ and

$$
\phi(B \setminus A) =\phi(B) - \phi(A).
$$

\bigskip \noindent In particular, $\phi$ is finite if $\phi(\Omega)$ is.\\

\noindent \textbf{Hint}. Let $A$ and $B$ be two elements of $\mathcal{A}$ such that $A \subset B$ and $\phi(B)$ is finite. Use the relationship $\phi(B)=\phi(B\setminus A) + \phi(A)$. Consider two cases :
the case where one of them is non-positive and the case where both $\phi(B\setminus A)$ and $\phi(A)$ are both non-negative, and use again, when necessary, the fact that
the value $-\infty$ is excluded.\\

\noindent \textbf{Be careful}. $\phi(\Omega)$ is not necessarily a maximum value for $\phi$, while it is if $\phi$ is non-negative.\\

\bigskip \noindent \textbf{SOLUTIONS}.\\

\noindent Let us start with $\phi(B)=\phi(B \setminus A) + \phi(A)$ with $\phi(B)$ finite and consider two cases :\\

\noindent Case 1 : One $\phi(B \setminus A)$ and $\phi(A)$ is non-positive, say $\phi(B \setminus A)$ is non-positive. Since $\phi$ does not take the value $-\infty$, it is finite and then the other term $\phi(A)$ is finite as the difference between two finite numbers.\\

\noindent Case 2 : Both $\phi(B\setminus A)$ and $\phi(A)$ are non-negative. Thus, they are both finite since their sum is finite.\\

\noindent Finally if $\phi(\Omega)$ is finite, $\phi$ is finite.\\

\bigskip \noindent \textbf{Exercise 4}. \label{exercise04_sol_doc08-04}\\

\noindent Let $\phi$ is $\sigma$-additive on $\mathcal{A}$ be a proper application on $\mathcal{A}$.\\

\noindent Question (a) Let $\phi$ be $\sigma$-additive. Show that $\phi$ is continuous.\\

\noindent Question (b) Let $\phi$ be finitely additive and

\hskip 5cm (b1) either continuous from below

\hskip 5cm (b2) or finite and continuous at $\emptyset$ (from above).

\noindent Show that $\phi$ is $\sigma$-additive.\\

\noindent \textit{Hints}.\\

\noindent \textbf{Hints for Question (a)}. Adapt the solutions of Exercise 2 in Doc 04-06 (page \pageref{exercise02_sol_doc04-06}) in Chapter \ref{04_measures} in the following way.\\

\noindent Let $(A_n)_{n\geq 0}$ a non-decreasing sequences in $\mathcal{A}$. Write 

$$
A=A_0 + \sum_{k\geq 1} (A_k \setminus A_{k-1}).
$$

\noindent and apply $\sigma$-additivity. Now discuss over two cases. \\

\noindent Case 1 : There exists $n_0\geq 0$ such that $\phi(A_n)=+\infty$ for all $n\geq n_0$. Use Exercise 3 to justify that Formula

$$
\phi(\lim_{n\rightarrow +\infty} A_n)=\lim_{n\rightarrow +\infty} \phi(A_n). \ (LIM)
$$

\noindent holds since both members are $+\infty$.\\

\noindent Case 2 : $\forall n\geq 0$, \ $\exists N\geq n$, \ $|\phi(A_n)|<+\infty$. Use Exercise 3 to see that all the $\phi(A_n)$'s are finite and apply the techniques of Chapter \ref{04_measures} :

\begin{eqnarray*}
&&\phi(A_0) + \sum_{k=1}^{+\infty} \phi(A_k \setminus A_{k-1})\\
&=&\lim_{n\rightarrow +\infty} \phi(A_0) + \sum_{k=1}^{n} \phi(A_k) - \phi(A_{k-1})\\
&=& \lim_{n\rightarrow +\infty} \phi(A_0) + (\phi(A_1)-\phi(A_0)) + (\phi(A_2)-\phi(A_1))+\cdots +(\phi(A_n)-\phi(A_{n-1}))\\
&=&\lim_{n\rightarrow +\infty} \phi(A_n).
\end{eqnarray*}

\noindent Next, let $(A_n)_{n\geq 0}$ be any non-increasing sequence for which there exists $n_0\geq 0$ such that $\phi(A_{n_0})$. Apply the results to $A_{n_0}\setminus A_n$ and conclude.\\

\noindent \textbf{Hints for Question (b1)}. Consider a sequence $(A_n)_{n\geq 0}$ of pairwise disjoint elements of $\mathcal{A}$. Use the formula

$$
\sum_{n=0}^{k} A_n \nearrow \sum_{n=0}^{+\infty} A_n 
$$ 

\noindent and the hypotheses to conclude.\

\noindent \textbf{Hints for Question (b2)}. Consider a sequence $(A_n)_{n\geq 0}$ of pairwise disjoint elements of $\mathcal{A}$. Use the formula

$$
\sum_{n=0}^{+\infty} A_n = \sum_{n=0}^{k} A_n  + \sum_{n=k+1}^{+\infty} A_n. 
$$ 

\noindent Combine the fact that $B_{k}=\sum_{n=k+1}^{+\infty} A_n \downarrow \emptyset$ with the hypotheses and conclude.\\

\bigskip \noindent \textbf{SOLUTIONS}.\\

\noindent Question (a). Let us consider a non-decreasing sequence $(A_n)_{n\geq 0}$ of elements in $\mathcal{A}$. We surely have 

$$
\bigcup_{n\geq 0}=A_0 + \sum_{k\geq 1} (A_k \setminus A_{k-1}).
$$

\noindent and by applying the $\sigma$-additivity, we have

$$
\phi\left(\bigcup_{n\geq 0} A_n\right)=\phi(A_0) + \sum_{k\geq 1} \phi(A_k \setminus A_{k-1}). \ (SU01)
$$

\noindent Now, let us consider two cases.\\

\noindent Case 1 : There exists $n_0\geq 0$ such that $\phi(A_n)=+\infty$ for all $n\geq n_0$. By Exercise 3, we have $\phi\left(\bigcup_{n\geq 0} A_n\right)=+\infty$, otherwise all the $\phi(A_n)$'s would be finite. Hence, it is clear that Formula (LIM) holds  with both members the equation being $+\infty$.\\

\noindent Case 2 : In the alternative of Case 1, for each $n\geq 0$, there exists $N\geq n$ such that $\phi(A_N)$ is finite. Since $A_n \subset A_N$, we have that $\phi(A_n)$ is finite by Exercise 3 and for all $k\geq 1$, we have

$$
\phi(A_k \setminus A_{k-1})=\phi(A_k) - \phi(A_{k-1}).
$$.

\noindent Thus Formula (SU01) becomes

\begin{eqnarray*}
\phi\left(\bigcup_{n\geq 0} A_n\right)&=&\phi(A_0) + \sum_{k=1}^{+\infty} \phi(A_k \setminus A_{k-1})\\
&=&\lim_{n\rightarrow +\infty} \phi(A_0) + \sum_{k=1}^{n} \phi(A_k) - \phi(A_{k-1})\\
&=& \lim_{n\rightarrow +\infty} \phi(A_0) + (\phi(A_1)-\phi(A_0)) +\cdots +(\phi(A_n)-\phi(A_{n-1}))\\
&=&\lim_{n\rightarrow +\infty} \phi(A_n).
\end{eqnarray*}

\noindent This establishes the continuity below. To handle the continuity above, consider a non-increasing sequence $(A_n)_{n\geq 0}$ of elements in $\mathcal{A}$ with $\phi(A_{n_0}$ finite. By Exercise 3,
the elements non-decreasing sequence $(A_{n_0}\setminus A_n)_{n\geq n_0}$ have finite values through $\phi$ and we have $\phi(A_{n_0} \setminus A_n)=\phi(A_{n_0}) - \phi(A_n)$ for $n\geq n_0$ and

\begin{eqnarray*}
\phi\left(\bigcup_{n\geq n_0} A_{n_0}\setminus A_n\right)&=&\phi\left(A_{n_0} \setminus \biggr(\bigcup_{n\geq n_0} \setminus A_n\biggr) \right)\\
&=&\phi\left(A_{n_0}\right) - \phi\biggr(\bigcup_{n\geq n_0} \setminus A_n\biggr).
\end{eqnarray*}

\bigskip \noindent We may and do apply the continuity below to get

\begin{eqnarray*}
\phi\left(\bigcup_{n\geq n_0} A_{n_0}\setminus A_n\right)&=&\lim_{n\rightarrow +\infty} \phi(A_{n_0} \setminus A_n)=\phi(A_{n_0})\\
&=&\lim_{n\rightarrow +\infty} \phi(A_{n_0}) - \phi(A_n)= \phi(A_{n_0}) - \phi(A_n) - \lim_{n\rightarrow +\infty} \phi(A_n),
\end{eqnarray*}

\noindent that is

$$
\phi\left(A_{n_0}\right) - \phi\biggr(\bigcup_{n\geq n_0} \setminus A_n\biggr)= \phi(A_{n_0}) - \phi(A_n) - \lim_{n\rightarrow +\infty} \phi(A_n).
$$

\noindent We conclude by dropping the finite number $\phi\left(A_{n_0}\right)$ from both sides.\\

\noindent \textbf{Question (b)}.\\

\noindent (b1) Suppose that $\phi$ is additive and continuous from below. Let $(A_n)_{n\geq 0}$ be a sequence of pairwise disjoint elements of $\mathcal{A}$. Ot is clear that

$$
\sum_{n=0}^{k} A_n \nearrow \sum_{n=0}^{+\infty} A_n.
$$ 

\noindent By using the below continuity and the finite additivity, we have

\begin{eqnarray*}
\phi\left(\sum_{n\geq 0} A_n\right)&=&\phi\left(\lim_{k\rightarrow +\infty} \sum_{n=0}^{k} A_n \right)\\
&=&\lim_{k\rightarrow +\infty} \phi\left( \sum_{n=0}^{k} A_n \right)\\
&=&\lim_{k\rightarrow +\infty} \sum_{n=0}^{k} \phi\left(  A_n \right)\\
&=&\sum_{n=0}^{+\infty} \phi\left( A_n \right).
\end{eqnarray*}

\noindent (b2) Suppose that $\phi$ is additive and continuous from above at $\emptyset$. Let $(A_n)_{n\geq 0}$ be a sequence of pairwise disjoint elements of $\mathcal{A}$. We have

$$
\sum_{n=0}^{+\infty} A_n = \sum_{n=0}^{k} A_n  + \sum_{n=k+1}^{+\infty} A_n. 
$$ 

\noindent It is clear that each $B_{k}=\sum_{n=k+1}^{+\infty} A_n$ is in $\mathcal{A}$ and the sequence $(B_{k})_{k\geq 1}$ is non-increasing to the empty set (See for example Exercise 14 in Doc 04-06, page \pageref{exercise14_sol_doc04-06}). Thus by our hypotheses, we have $\phi(B_k) \rightarrow \phi(\emptyset)=0$ as $k \rightarrow +\infty$. Now, by additivity, we have for all $k\geq 1$,

$$
\phi`\left(\sum_{n=0}^{+\infty} A_n\right) = \sum_{n=0}^{k} \phi(A_n)  + \phi(B_k).
$$

\noindent which leads to the $\sigma$-additivity with $k\rightarrow +\infty$.\\

\bigskip \noindent \textbf{Exercise 5}. \label{exercise05_doc06-04}\\

\noindent Let $\phi$ be $\sigma$-additive on $\mathcal{A}$ which does not take either the value $-\infty$ or the value $+\infty$. Show that there exists two measurable sets $C$ and $D$ such that $\phi(C)=\sup \phi$ and $\phi(D)=\inf \phi$.\\

\noindent \textit{Hints}.\\

\noindent \textbf{Existence of C}.\\

\noindent Suppose that there is a measurable set $C$ such that $\phi(C)=+\infty$. Conclude.\\

\noindent Suppose that $\phi \leq 0$. Take $C=\emptyset$ and conclude.\\

\noindent Now suppose that $\phi(A)<+\infty$ for all finite for all $A \in \mathcal{A}$ and that there exists $A_0\in \mathcal{A}$ such that $\phi(A_0)>0$.\\

\noindent (1) Justify the existence of a sequence $(A_n)_{n\geq 1} \subset \mathcal{A}$ such that the $\phi(A_n)$'s are finite for large values on $n$ and  $\left(\phi(A_n)\right)_{n\geq 1}$ converges to $\alpha=\sup \phi$ in $\overline{\mathbb{R}}$.\\

\noindent (2) Put $A=\bigcup_{n\geq 1} A_n$. For each $n\geq 1$, consider the partitions $A=A_k + (A \setminus A_k)$ for $k=1,...,n$ and the the intersection $A$ $n$ times and conclude that $A$ may is partitioned into sets of the form $A_1^{\prime}A_2^{\prime}...A_n^{\prime}$ where $A_k^{\prime}=A_k$ or $A_k^{\prime}=A\setminus A_k$ for $k=1,...,n$.\\

\noindent What is the size $\ell_n$ of the partition of $A$ for $v\geq 1$ fixed? Denote the elements of that partition by $A_{n,j}$, $j=1,...,\ell_n$.\\

\noindent For example, for three sets and $n=2$, we have
 
$$
A=A_1 \cup A_2 \cup A_3= A_1 \cap A_2 + (A\setminus A_1)\cap A_2 + A_1 \cap (A\setminus A_2) + (A\setminus A_1) \cap (A\setminus A_2).
$$ 

\noindent (3) Fr each $n\geq 1$, we have $A_n \subset A=\sum_{1\leq j \leq \ell_n} A_{n,j}$. By doing $A_n=A\cap A_n= \sum_{1\leq j \leq \ell_n} A\cap A_{n,j}$, can you say that $A_n$ is a sum of elements
of some elements $A_{n,j}$, $j=1,...,\ell_n$. What are these elements?\\

\noindent (4) As well, let $n <p$. By take the intersection of the partitions $A=A_k + (A \setminus A_k)$, $k=n+1,...,p$, $(p-n)$ times, make a conclusion similar to the one given in (2) relative a certain partition $\mathcal{P}_{n,p}$ to be precised.\\

\noindent Now Take any $A_{n,j}$, $k=1,...,\ell_n$. By remarking  that $A_{n,j}$ is a subset of $A$, and by forming the intersection $A_{n,j}=A\cap A_{n,j}$ using the partition $\mathcal{P}_{n,p}$, show that $A_{n,j}$ is some of elements $A_{p,h}$, $h=1,...,\ell_p$.\\

\noindent (5) Define
$$
I_n={j, \ 1\leq j \leq \ell_n, \ and \ \phi(A_{n,j})\geq 0 }, \ J_n={j, \ 1\leq j \leq \ell_n, \ \phi(A_{n,j})< 0 },  
$$
 
\noindent where $I_n+J_n=\{1,...,\ell_n\}$, and

$$
B_n=\emptyset \ if \ I_n=\emptyset, \ and  \ B_n=\sum_{j \in I_n} A_{n,j} \ if \ I_n\neq \emptyset.
$$

\noindent Define also, for $n\geq 1$ fixed, $\mathcal{S}_n$ as the collection of finite sums of sets $A_{n,j}$, $j\in I_n$ and $\mathcal{D}_n$ as the collection of finite sums of sets $A_{n,j}$, $j\in J_n$.

\noindent Remark that, by Exercise 14 (see page \pageref{exercise14_sol_doc01-03}) in Doc 01-03 of Chapter \ref{01_setsmes}, $\mathcal{S}_n$ and $\mathcal{D}_n$ are stable by finite unions. Try it for the collection of finite sums of the mutually disjoint sums of, say, $A$, $B$, $C$.\\

\noindent Show that for $n\geq 1$ fixed, $B_n \in \mathcal{S}_n$, $\phi$ is non-decreasing on $\mathcal{S}_n$ and non-increasing on $\mathcal{S}_n$, and that if $E\in \mathcal{S}_n$ and $F\in \mathcal{D}_n$ with $E \ cap EF=\emptyset$, use the non-positivity of $\phi(F)$ and the facts above to show that

$$
\phi(E+F) \leq \phi(E).
$$

\noindent (6a) From now, $n \geq $ is fixed. From (3), explain why $A_n$ is of the form  $A_n=E+F$ with $E\in \mathcal{S}_n$ and $F\in \mathcal{D}_n$. Use the non-decreasingness of $\phi$ on 
$\mathcal{S}_n$ to get

$$
\phi(A_n) \leq \phi(B_n).
$$

\bigskip \noindent (6b) You are going to prove the for any $p\geq n$, we have

$$
\phi\left( B_{n} \bigcup B_{n+1} \bigcup ... \bigcup B_{p} \right) \leq \phi\left( B_{n} \bigcup B_{n+1} \bigcup ... \bigcup B_{p} \bigcup \bigcup B_{p+1}\right). \ (AC)
$$

\noindent Fix $r \in\{n,...,p\}$. By (4), $B_r$ is a sum of elements $\mathcal{S}_{p+1}$ and we may write $B_r=E_{r,1}+E_{r,2}$ with $E_{r,1} \in \mathcal{S}_{n+1}$ and $E_{r,2}\in \mathcal{D}_{n+1}$. We have $E_{r,1} \subset B_{p+1}$ and $E_{r,2}$ is disjoint of $B_{p+1}$. Validate the computations : 

\begin{eqnarray*}
E&=:& B_{n} \bigcup B_{n+1} \bigcup ... \bigcup B_{p} \bigcup B_{p+1}\\
&=& \bigcup_{1\leq r \leq p} \left(B_r \bigcup B_{p+1}\right)\\
&=& \bigcup_{1\leq r \leq p} \left(E_{r,1} \bigcup E_{r,2}\right)\\
&=& \left(\bigcup_{1\leq r \leq p} E_{r,1} \right) \bigcup \left(\bigcup_{1\leq r \leq p} E_{r,2}\right).\\
\end{eqnarray*}

\noindent By exploiting that $\mathcal{S}_{p+1}$ and $\mathcal{D}_{p+1}$ are stable by finite unions, conclude that 

\begin{eqnarray*}
E&=:&B_{n} \bigcup B_{n+1} \bigcup ... \bigcup B_{p} \bigcup B_{p+1}\\
&=& \left(\bigcup_{1\leq r \leq p} E_{r,1} \right) \bigcup \left(\bigcup_{1\leq r \leq p} E_{r,2}\right)\\
&=:&E_1 + E_2,
\end{eqnarray*}

\noindent with $E_1 \in \mathcal{S}_{p+1}$ and $E_2 \in \mathcal{D}_{p+1}$, and then, $E_1 \subset B_{p+1}$ and $E_2$ is disjoint of $B_{p+1}$. Validate this,  

$$
E \cup B_{p+1}= E_1 \cup B_{p+1} \cup E_2=B_{^p+1} \cup E_2=B_{p+1} + E_2
$$

\noindent and, by using again the non-decreasingness of $\phi$ on $\mathcal{S}_{p+1}$, show that
 
$$
\phi(E)=\phi(E_1)+\phi(E_2) \leq \phi(B_{p+1})+\phi(E_2)=\phi(E \cup B_{p+1}).
$$

\noindent (7) Combine (5) and (6) and repeat Formula (A) as needed as possible to get, for all $n\geq n_0$, for all $p>n$;

$$
\phi(A_n) \leq \phi(B_n) \leq \phi\left( \bigcup_{n\leq r \leq p} B_r\right). 
$$

\noindent Conclude and find $C$ by letting first  $p\uparrow +\infty$ and next $n\uparrow +\infty$.\\

\noindent \textbf{Existence of D}. Consider $-\phi$ and apply the first part to get a set $D\in \mathcal{A}$ such that

$$
-\phi(D)=  \sup - \phi \ \ \Leftrightarrow \phi(D)=  \inf \phi. 
$$

\bigskip \noindent \textbf{SOLUTIONS}.\\

\noindent \textbf{Existence of C}.\\

\noindent If there is a measurable set $C$ such that $\phi(C)=+\infty$, we surely have $\phi(C)=\sup \phi=+\infty$. If $\phi \leq 0$, we also have $0=\phi(\emptyset)=\sup \phi$.\\

\noindent Now, if $\phi(A)<+\infty$ for all finite for all $A \in \mathcal{A}$ and there exists $A_0\in \mathcal{A}$ such that $\phi(A_0)>0$, we have that $\alpha=\sup \phi$ positive (at this step, ot might be $+\infty$).\\

\noindent (1) By the properties of the supremum in $\overline{\mathbb{R}}$ there exists a maximizing sequence, that is a sequence $(A_n)_{n\geq 1} \subset \mathcal{A}$ such that  $\left(\phi(A_n)\right)_{n\geq 1}$ converges to $\alpha=\sup \phi$ in $\overline{\mathbb{R}}$. Since $\alpha$ is positive, we surely have $\phi(A_n)>0$ for $n$ large enough, that is when $n\geq b_0$ for some $n_0\geq 1$. Hence 
the $\phi(A_n)$, $n\geq n_0$, are positive and strictly less that $+\infty$ : they are finite.\\

\noindent (2) Let us denote $A=\bigcup_{n\geq 1} A_n$. For each $n\geq 1$, we have the partitions $A=A_k + (A \setminus A_k)$ for $k=1,...,n$. We get

$$
\underset{n \ times}{\underbrace{A \cap ... \cap A}}=\bigcap_{1\leq k \leq n} A_k + (A \setminus A_k).
$$

\noindent By developing the intersection, we get a summation of sets of the form $A_1^{\prime}A_2^{\prime}...A_n^{\prime}$ where $A_k^{\prime}=A_k$ or $A_k^{\prime}=A\setminus A_k$ for $k=1,...,n$. By doing so, we have to choose for each $k=1,...,n$, between the two choices : putting $A_k$ or A $\setminus A_k$. In total, we have $\ell_n=2^n$ possibilities. So let us denote the elements of the partition by : $A_{n,j}$, $j=1,...,\ell_n$.\\

\noindent (3) For each $n\geq 1$, $A_n$ is a some of $A_{n,j}$. To see this, we remark that  $A_n \subset A=\sum_{1\leq j \leq \ell_n} A_{n,j}$ and we have

$$
A_n=A\cap A_n= \sum_{1\leq j \leq \ell_n} A_n\cap A_{n,j}.
$$

\noindent In each $A_{n,j}$, either $A_n$ is used or $(A\setminus A_n)$ is used. In the first case, $A_n\cap A_{n,j} =A_{n,j}$ and in the second, $A_n\cap A_{n,j}=\emptyset$. In total, $A_n$ is the sum of the set $A_{n,j}$ in which $A_n$ is used instead of $A\setminus A_n$.\\

\noindent (4) Le $n <p$. By take the intersection of the partitions $A=A_k + (A \setminus A_k)$, $k=n+1,...,p$, $(p-n)$ times, we also get that $A$ is the summation of sets of the form 
$A_n^{\prime}A_{n+1}{\prime}...A_p^{\prime}$ where $A_k^{\prime}=A_k$ or $A_k^{\prime}=A\setminus A_k$ for $k=n+1,...,p$. Let us denote the elements of that partition by $H_{n,p,j}$, 
$j=1, ... ,d(n,p)=2^{p-n}$.\\

\noindent Now for any $A_{n,j}$, $k=1,...,\ell_n$, we have $A_{n,j} \subset A=\sum_{1\leq j \leq d(n,p)} H_{n,p,j}$, and then

$$
A_{n,j}=A_{n,j} \cap A=\sum_{1\leq j \leq d(n,p)} A_{n,j} \cap H_{n,p,j}.
$$

\noindent But is clear that each $A_{n,j} \cap H_{n,p,j}$ is of the form 
$$
A_1^{\prime}A_2^{\prime}...A_n^{\prime} A_{n+1}^{\prime}...A_{p}^{\prime},
$$ 

\noindent where $A_k^{\prime}=A_k$ or $A_k^{\prime}=A\setminus A_k$ for $k=1,...,n,n+1,...,p$.\\

\noindent (5) Let us define
$$
I_n=\{j, \ 1\leq j \leq \ell_n, \ and \ \phi(A_{n,j})\geq 0 \}, \ J_n=\{j, \ 1\leq j \leq \ell_n, \ \phi(A_{n,j})< 0 \},  
$$
 
\noindent where we have $I_n+J_n=\{1,...,\ell_n\}$, and 

$$
B_n=\emptyset \ if \ I_n=\emptyset, \ and  \ B_n=\sum_{j \in I_n} A_{n,j} \ if \ I_n\neq \emptyset.
$$

\noindent Also, let us define, for $n\geq 1$ fixed, $\mathcal{S}_n$ as the collection of finite sums of sets $A_{n,j}$, $j\in I_n$ and $\mathcal{D}_n$ as the collection of finite sums of sets $A_{n,j}$, $j\in J_n$.\\

\noindent By \textit{Exercise 14 (see page \pageref{exercise14_sol_doc01-03}) in Doc 01-03 of Chapter \ref{01_setsmes}}, it is clear that $\mathcal{S}_n$ and $\mathcal{D}_n$ are stable by finite unions. Obviously, $\phi$ additive on  $\mathcal{S}_n$ and is non-negative. So it is non-decreasing on it. As well it is non-increasing on  $\mathcal{D}_n$  similar reasons. Besides,  $B_n \in \mathcal{S}_n$ obviously.\\

\noindent (6a) From now, $n \geq $ is fixed. From (3), we know that $A_n$ is a sum of elements $A_{n,j}$. By regrouping the $A_{n,j}$'s in $A_n$ and corresponding to $j\in I_n$ in a set $E$, and 
the $A_{n,j}$'s in $A_n$ and corresponding to $j\in J_n$ in a set $F$, we get $A_n=E+F$ with $E\in \mathcal{S}_n$ and $F\in \mathcal{D}_n$. Now, by using the non-positivity of $\phi(F)$ and the non-decreasingness of $\phi$ on $\mathcal{S}_n$, we get

$$
\phi(A_n)=\phi(E)+\phi(F) \leq \phi(E)\leq \phi(B_n)
$$

\bigskip \noindent (6b) We want to prove that for any $p> n\geq 1$, we have

$$
\phi\left( B_{n} \bigcup B_{n+1} \bigcup ... \bigcup B_{p} \right) \leq \phi\left( B_{n} \bigcup B_{n+1} \bigcup ... \bigcup B_{p} \bigcup \bigcup B_{p+1}\right). \ (AC)
$$

\bigskip \noindent Fix $r \in\{n,...,p\}$. By (4), $B_r$ is a sum of elements $\mathcal{S}_{p+1}$ and we may write $B_r=E_{r,1}+E_{r,2}$ with $E_{r,1} \in \mathcal{S}_{p+1}$ and $E_{r,2}\in \mathcal{D}_{p+1}$. We have $E_{r,1} \subset B_{p+1}$ and $E_{r,2}$ is disjoint of $B_{p+1}$. So 

\begin{eqnarray*}
E&=:& B_{n} \bigcup B_{n+1} \bigcup ... \bigcup B_{p} \bigcup B_{p+1}\\
&=& \bigcup_{1\leq r \leq p} \left(B_r \bigcup B_{p+1}\right)\\
&=& \bigcup_{1\leq r \leq p} \left(E_{r,1} \bigcup E_{r,2}\right)\\
&=& \left(\bigcup_{1\leq r \leq p} E_{r,1} \right) \bigcup \left(\bigcup_{1\leq r \leq p} E_{r,2}\right).
\end{eqnarray*}

\noindent Since $\mathcal{S}_{p+1}$ and $\mathcal{D}_{p+1}$ are stable by finite unions, we arrive at

\begin{eqnarray*}
E&=:&B_{n} \bigcup B_{n+1} \bigcup ... \bigcup B_{p} \bigcup B_{p+1}\\
&=& \left(\bigcup_{1\leq r \leq p} E_{r,1} \right) \bigcup \left(\bigcup_{1\leq r \leq p} E_{r,2}\right)\\
&=:&E_1 + E_2,
\end{eqnarray*}

\noindent with $E_1 \in \mathcal{S}_{p+1}$ and $E_2 \in \mathcal{D}_{p+1}$, and then, $E_1 \subset B_{p+1}$ and $E_2$ is disjoint of $B_{p+1}$. This implies  

$$
E \cup B_{p+1}= E_1 \cup B_{p+1} \cup E_2=B_{^p+1} \cup E_2=B_{p+1} + E_2
$$

\noindent and, by using again the non-decreasingness of $\phi$ on $\mathcal{S}_{p+1}$, 
 
$$
\phi(E)=\phi(E_1)+\phi(E_2) \leq \phi(B_{p+1})+\phi(E_2)=\phi(E \cup B_{p+1}).
$$

\noindent (7) By Combining (5) and (6), by iterating  Formula (AC), we get for all $n\geq n_0$, for all $p>n$,

$$
\phi(A_n) \leq \phi(B_n) \leq \phi\left( \bigcup_{n\leq r \leq p} B_r\right). 
$$

\noindent By letting $p\uparrow +\infty$, we have by continuity from below of $\phi$, for all $n\geq n_0$

$$
\phi(A_n) \leq \phi(B_n) \leq \phi\left( \bigcup_{r=n}^{+\infty} B_r\right). \ (FF)
$$

\noindent The sequence $C_n=\bigcup_{r=n}^{+\infty} B_r$, $n\geq 1$, is non-increasing to a set $C$, the values $\phi(C_n)$'s are finite for $n\geq n_0$. We may apply the continuity from below if $\phi$ in Formula (FF) to get, as $n\rightarrow +\infty$,

$$
\alpha=\sup \phi = \lim_{n \rightarrow +\infty} \phi(A_n) \leq \phi(C) \leq \sup \phi,
$$

\noindent that is $\phi(C) = \sup \phi$.\\

\noindent \textbf{Existence of D}. Proceed as suggested in the Hint, that amounts to doing almost nothing.\\

\bigskip \noindent \textbf{Exercise 6}. \label{exercise06_sol_doc08-04}\\

\noindent Let $\phi$ be a $\sigma$-additive application from $\mathcal{A}$ to $\mathbb{R}\cup \{+\infty\}$. Define for all $A \in \mathcal{A}$,

$$
\phi_1(A)= \sup \{\phi(B), \ B \subset A, \ B\in \mathcal{A}\} \ and \ \phi_2(A)= -\inf \{\phi(B), \ B \subset A, \ B\in \mathcal{A}\}.\\
$$

\noindent Consider the measurable space $D$, found in Exercise 5, such that $\phi(D)=\inf \phi$.\\

\noindent Question (a). If $\phi$ is a measure, check that : $\phi_1=\phi$, $\phi_2=0$,  $\phi_1(A)=\phi(AD^c)$ and $\phi_2=\phi(AD)$ for $D=\emptyset$.\\

\noindent Question (b). In the sequel, we suppose that there exists $A_0\in \mathcal{A}$ such that $\phi(A_0)<0$.\\

Show that for all $A\in \mathcal{A}$, $\phi(AD)\leq 0$ and $\phi(AD^)\geq 0$.\\

\noindent \textit{Hints}. Show that $\phi(D)$ is finite and by using Exercise 3, that $\phi(AD)$. By using meaning of the infimum, decompose $\phi(D \setminus AD)$ and show we have a contradiction if we don't have $\phi(AD)\leq 0$. Do the same for $\phi(D + AD^c)$ to show that $\phi(AD^)\geq 0$.\\

\noindent Question (c).\\

\noindent (c1) Justify each line of the formula below : For $A \in \mathcal{A}$, for all $B\in \mathcal{A}$ with $B \subset A$, we have

\begin{eqnarray*}
\phi(B)&=&\phi(BD)+\phi(BD^c)\\
&\leq& \phi(BD^c)\\
&\leq& \phi(BD^c) +\phi((A\setminus B)D^)\\
&=&\phi(AD^c).
\end{eqnarray*}

\noindent Deduce from this that : $\phi_1(A)=\phi(AD^c)$ for all $A\in \mathcal{A}$.\\

\noindent (c2) Justify each line of the formula below : For $A \in \mathcal{A}$, for all $B\in \mathcal{A}$ with $B \subset A$, we have

\begin{eqnarray*}
\phi(B)&=&\phi(BD)+\phi(BD^c)\\
&\geq& \phi(BD)\\
&\geq& \phi(BD) +\phi((A\setminus B)D) \\
&=&\phi(AD).
\end{eqnarray*}

\noindent Deduce from this that : $\phi_2(A)=-\phi(AD)$ for all $A\in \mathcal{A}$.\\

\noindent Question (d) Check (by two simple sentences) that $\phi_1$ and $\phi_2$ are measures. Is the application $\phi_1-\phi_2$ well-defined? Why?\\

\noindent Question (e) Conclude that $\phi=\phi_1-\phi_2$.\\

\noindent \emph{NB}. If the value $+\infty$ is excluded in place of $-\infty$, we establish the relation for $-\phi$ and adapt the results.\\

\bigskip \noindent \textbf{SOLUTIONS}.\\

\noindent The solution uses a $\sigma$-finite application which does not take the value $-\infty$. In the case where the value $+\infty$ is avoided, we apply the results below to $-\phi$ and adapt the results.\\

\noindent Question (a) If $\phi$ is itself a measure, we have the decomposition $\phi_1=\phi$, $\phi_2=0$ and $D=\emptyset$.\\

\noindent Question (b) Otherwise, there exists $A_0 \in \mathcal{A}$ such that $\phi(A_0)<0$. By Exercise 5 above, there exists $D \in \mathcal{A}$ such that

$$
-\infty < \inf \phi =\phi(D) \leq \phi(A_0)<0.
$$

\noindent Hence $\phi(D)$ is finite.\\

\noindent Now for any $A\in \mathcal{A}$, $AD \subset D$. And, since $\phi(D)$  is finite, we have by Exercise 3 :  

$$
\phi(A \setminus AD)=\phi(D) - \phi(AD).
$$

\noindent By the definition of the infimum, the above formula implies that $\phi(AD)\leq 0$. As well, we have by additivity

$$
\phi(D+AD^c)=\phi(D)+\phi(AD^c),
$$

\noindent which, for the same reason, requires that $\phi(AD^c)\geq 0$.\\

\noindent Question (c). Now, Let $A\in \mathcal{A}$.\\

\noindent (c1) We have
$$
\phi(AD^c) \leq \phi_1(A)= \sup \{\phi(A), \ B \subset A and B\in \mathcal{A}\},
$$

\noindent since $\phi(AD^c)$ is in the range over which the supremum is taken. Now for all $B\in \mathcal{A}$ with $B \subset A$, we have

\begin{eqnarray*}
\phi(B)&=&\phi(BD)+\phi(BD^c)\\
&\leq& \phi(BD^c) \ since \ \phi(BD)\leq 0\\
&\leq& \phi(BD^c) +\phi((A\setminus B)D^) \ since \ \phi((A\setminus B)D^c)\geq 0\\
&=&\phi(AD^c) \ by \ additivity.
\end{eqnarray*}

\noindent By taking the supremum over $B\subset A$, $B\in \mathcal{A}$, we get $\phi_1(A)\leq \phi(AD^c)$. We conclude that $\phi_1(A)=\phi(AD^c)$ for all $A\in \mathcal{A}$, and that 
$\phi_1$ is a measure.\\

\noindent (c2) We also that for $A\in \mathcal{A}$ :
$$
-\phi(AD) \leq \phi_2(A)= - \inf \{\phi(B), \ B \subset A and B\in \mathcal{A}\},
$$

\noindent since $\phi(AD)$ is in the range over which the infimum is taken. Now for all $B\in \mathcal{A}$ with $B \subset A$, we have

\begin{eqnarray*}
\phi(B)&=&\phi(BD)+\phi(BD^c)\\
&\geq& \phi(BD) \ since \ \phi(BD^c)\geq 0\\
&\geq& \phi(BD) +\phi((A\setminus B)D) \ since \ \phi((A\setminus B)D)\leq 0\\
&=&\phi(AD) \ by \ additivity.
\end{eqnarray*}

\noindent By taking the infimum over $B\subset A$, $B\in \mathcal{A}$, we get $-\phi_2(A) \geq \phi(AD)$, that is $-\phi(AD) \geq \phi_2(A) \geq$. We conclude that $\phi_2(A)=-\phi(AD)$ for all $A\in \mathcal{A}$, and that $\phi_2$ is a measure.\\

\noindent Question (d) We already saw that for all $A\in \mathcal{A}$, $\phi(AD)$ is finite. Hence the following relation is well justified :

$$
\phi(A)=\phi(AD^c)+\phi(AD)=\phi_1(A) - \phi_2(A),
$$

\noindent and the decomposition is established.\\

\noindent

\noindent \LARGE \textbf{DOC 08-05 : Radon-Nikodym's Theorem - Exercises/Solutions} \label{doc08-05}\\
\Large

%\bigskip \noindent \textbf{Exercise X}. \label{exercise0X_sol_doc08-05}\\
%\bigskip \noindent \textbf{SOLUTIONS}.\\

\bigskip \noindent \textbf{NB}. Throughout the document, unless the contrary is specifically notified, $(\Omega,\mathcal{A},m)$ is a measure space $\phi$, $\phi_i$, $\psi$  and $\psi-i$ (i=1,2) are real-valued applications defined on $\mathcal{A}$ and $X : (\Omega,\mathcal{A}) \rightarrow \overline{\mathbb{R}}$ is a function.\\

\bigskip \noindent \textbf{Exercise 1}. \label{exercise01_sol_doc08-05}\\

\noindent Let $(\Omega,\mathcal{A},m)$ be a measure space and $X : (\Omega,\mathcal{A}) \rightarrow \overline{\mathbb{R}}$ be a function.\\

\noindent Question (a) Let $X$ be non-negative ($X\geq 0$). Justify the definition

$$
\mathcal{A} \ni A \mapsto \phi_X(A)=\int_A X \ dm, \ (II)
$$ 

\bigskip 
\noindent and show, by using the \textit{MCT} (See Exercise 1, Doc 06-06, page \pageref{exercise01_doc06-06}) that $phi_X$ is $\sigma$-additive.\\ 

\noindent Question (b) Let $X^-$ be integrable or $X^+$ be integrable. Apply the the results of Question (a) to $X^-$ and $X^+$ and extend the definition \textit{(II)} by using
by first formula in Question (c), Exercise 6 in Doc 05-02 (page  \pageref{exercise06_doc05-02}) in Doc 05-02 in Chapter \ref{05_integration}) : for any real-valued function $f$, for any subset $A$ of $\Omega$, 
$$
(1_Af)^-=1_Af^- \ and \ (1_Af)^+=1_Af^+
$$

\bigskip \noindent and say why it is still $\sigma$-additive (Here, treat the case $X^-$ only).\\

\noindent Question (c) From now on, we suppose that $X^-$ is integrable or $X^+$ is integrable.\\

\noindent Show that $\phi_X$ is $m$-continuous.\\

\noindent \textit{Hints}. Use the first formula in Question (c), Exercise 6 in Doc 05-02 (page  \pageref{exercise06_doc05-02}) in Doc 05-02 in Chapter \ref{05_integration}) and Property (P3) for non-negative functions in Exercise 5 Doc 05-02 (page \pageref{doc05-02}) in Chapter \ref{05_integration}.\\

\noindent Question (d) Suppose that $m$ is $\sigma$-finite. Show that is $phi_X$ is $\sigma$-finite.\\

\noindent \textit{Hints}. Use a countable and measurable subdivision of $\Omega$ : $\Omega=\sum_{j\geq 0}$, $(\forall j\geq 0, \ \Omega_j \in \mathcal{A} \text{ and } m(\Omega_j)<+\infty)$ and the formula

$$
(f \ finite)=\sum_{n\in \mathbb{Z}} (n\leq X <n+1).
$$

\bigskip 
\noindent Show that $\phi_X$ is finite on each $A_{n,j}=(n\leq X <n+1)\bigcap \Omega_j$, $(n\in \mathcal{Z}, \ j\geq 0)$. Deduce from this that $\phi_X$ is $\sigma$-finite.\\

\noindent Question (e)  Readily check that $\phi_X$ is finite if $X$ is integrable.\\

\noindent Question (d)  Suppose that $X$ and $Y$ are both $m$-finite measurable and quasi-integrable. Show that $\phi_Y$ are  equal if and only if $X=Y$ $m$-a.e.\\

\bigskip \noindent \textbf{SOLUTIONS}.\\

\noindent Question (a) Let $X$ be non-negative ($X\geq 0$). Thus, the integral non-negative function $1_A X$ always exists and the definition

$$
\mathcal{A} \ni A \mapsto \phi_X(A)=\int_A X \ dm= \int 1_A X \ dm, \ (II01)
$$ 

\bigskip 
\noindent is well done. Now, let $(A_n)_{n\geq 0}$ be sequence of pairwise disjoint elements of $\mathcal{A}$. We apply \textit{MCT} to integral of a sum of non-negative functions 
(See Exercise 1, Doc 06-06, page \pageref{exercise01_doc06-06}) to have

\begin{eqnarray*}
\phi_X\left(\sum_{n \geq 0} A_n \right)& =& \int 1_{\left(\sum_{n \geq 0} A_n \right)} X \ dm\\
&=&\int 1_{\left(\sum_{n \geq 0} A_n \right)} X \ dm\\
&=& \int \sum_{n \geq 0}  \left(1_{A_n} X\right) \ dm\\
&=&  \sum_{n \geq 0}  \int \left( 1_{A_n} X\right) \ dm\\
&=& \sum_{n \geq 0}  \int_{A_n}  X \ dm\\
&=& \sum_{n \geq 0}  \phi_X(A_n),
\end{eqnarray*}

\noindent and thus $\phi_X$ is $\sigma$-additive.\\ 

\noindent Question (b) Suppose that Let $X^-$. By Question (a), we know that $\phi_{X^-}$ and $\phi_{X^+}$ are $\sigma$-additive. But $\phi_{X^-}$ is finite since

$$
\forall A \in \mathcal{A}, \  0 \leq \int 1_A X^- \ dm  \leq \int X^- \ dm <+\infty.
$$ 

\bigskip \noindent Thus the application $\phi_{X^+}-\phi_{X^-}$ is well-defined and is $\sigma$-finite. Moreover, by first formula in Question (c), Exercise 6 in Doc 05-02 (page  \pageref{exercise06_doc05-02}) in Doc 05-02 in Chapter \ref{05_integration}), we have  $(1_AX)^-= 1_AX^-$ and $(1_AX)^+= 1_AX^+$, and hence
$$
\int (1_AX)^- \ dm = \int 1_A X^- \ dm \leq \int X^- \ dm <+\infty.
$$

\bigskip 
\noindent Thus $(1_AX)$ is quasi-integrable and we have for all $A \in \mathcal{A}$,

\begin{eqnarray*}
\phi_X(A)&=&\int (1_AX)^- \ dm - \int (1_AX)^+  \ dm\\
&=& \int 1_A X^+ \ dm- \int X^- \ dm\\
&=&\phi_{X^+}(A)-\phi_{X^-}(A).
\end{eqnarray*}

\bigskip 
\noindent and finally is well-defined and $\phi_X(A)$ is $\sigma$-finite.\\

\noindent Question (c) Let $A \in \mathcal{A}$ such that $m(A)=0$. Property (P3) for non-negative functions in Exercise 5 Doc 05-02 (page \pageref{doc05-02}) in Chapter \ref{05_integration}, we have

$$
\phi_X^-(A)=\int_A X^- \ dm=0 \ and \ \phi_X^+(A)=\int_A X^+ \ dm=0
$$

\bigskip 
\noindent and next $\phi_X(A)=\phi_X^+(A)-\phi_X^-(A)=0$. $\square$\\

\noindent Question (d) Since $m$ is $\sigma$-finite, we have a measurable subdivision of $\Omega$ : $\Omega=\sum_{j\geq 0}$, $(\forall j\geq 0, \ \Omega_j \in \mathcal{A})$ such that
$m(\Omega_j)<+\infty$ for all $j\geq 0$. Next, since $X$ is \textit{a.e.} finite, we have

$$
m(X \ infinite)=0.
$$

\bigskip 
\noindent By denoting $N=(X \ infinite)$, we have $N^c=(f \ finite)=(f \in \mathbb{R})$. But $\mathbb{R}$ may be partitioned into

$$
\mathbb{R}^=\sum_{n\in \mathbb{Z}} [n, n+1[,
$$

\bigskip 
\noindent and this leads to

$$
N^c=\sum_{n\in \mathbb{Z}} (f \in [n, n+1[) =\sum_{n\in \mathbb{Z}} (n\leq f <n+1).
$$

\bigskip 
\noindent We get that subdivision :

$$
\Omega=\Omega \cap N + \Omega \cap N^c=\Omega \cap N + \sum_{j\geq 0} \sum_{n\in \mathbb{Z}} \Omega_j \cap  (n\leq f <n+1)=\Omega \cap N + \sum_{j\geq 0} \sum_{n\in \mathbb{Z}} A_{n,j}.
$$

\bigskip 
\noindent We have, since $m(\Omega \cap N)=0$, we have $\phi_X(\Omega \cap N)$ by Question (c), and for all $j\geq 0$, $n\in \mathbb{Z}$

$$
n \times m(\omega_j) \leq \phi_X(A_{n,j})=\int 1_{A_{n,j}} X \leq (n+1) \times m(\omega_j).
$$

\bigskip 
\noindent We conclude that $\phi_X$ is $sigma$-finite.

\noindent Question (e)  If $X$ is integrable, then for all $A \in \mathcal{A}$, $1_A X$ is integrable since it is dominated by $|X|$, and

$$
\left| \phi_X(A)=\int_A X \ dm \right| \leq \int \left|X\right| \ dm. \ \square
$$

\bigskip 
\noindent Question (d). Since $X$ and $Y$ are quasi-integrable, $\phi_X$ and $\phi_Y$ are defined. Since $X$ and $Y$ are $m$-finite, we have, for all $A\in \mathcal{A}$

$$
0=\phi_X(A)-\phi_Y(A)=\int_A (X-Y) \ dm.
$$

\bigskip 
\noindent Let $\varepsilon>0$. If we have $m(A)>0$ with $A=(X-Y>\varepsilon)$, we would get t

$$
0= \int_{A} (X-Y) \ dm \geq \varepsilon m(A) >0.
$$
 
\bigskip 
\noindent This contradictions leads to $m(X-Y>\varepsilon)=0$ for all $\varepsilon>0$. By the symmetry of the roles of $X$ and $Y$, we also have
$m(X-Y<\varepsilon)=0$ for all $\varepsilon>0$. Hence we have $m(|X-Y||>\varepsilon)=0$ for all $\varepsilon>0$. Hence $X=Y$ $m$.a.e.

\bigskip \noindent \textbf{Exercise 2}. \label{exercise02_sol_doc08-05}\\

\noindent Let $\phi$ and $\psi$ be two $\sigma$-additive applications which compatible in the sense that the value $-\infty$ is excluded for both of them or the value $+\infty$ is excluded for both of them, so that the application $\Phi=\phi+\psi$ is well-defined and is $\sigma$-finite.\\

\noindent Question (a) Show that if $\phi$ and $\psi$ are both $m$-continuous, then $\Phi$ is $m$-continuous.\\

\noindent Question (b) Show that if $\phi$ and $\psi$ are both $m$-singular \textbf{measures}, then $\Phi$ is $m$-continuous.\\

\noindent Question (c) Show that if $\phi$ is both $m$-continuous and  $m$-singular, then $\phi$ is the null application.\\

\bigskip \noindent \textbf{SOLUTIONS}.\\

\noindent Suppose that the value $-\infty$ is excluded for both of $\phi$ and $\psi$.\\

\noindent Question (a) Suppose that $\phi$ and $\psi$ are both $m$-continuous. Thus, for any $A \in \mathcal{A}$ such that $m(A)=0$, we have $\phi(A)=0$ and $\psi(A)=0$ and hence

$$
\Phi(A)=\phi(A)+\psi(A)=0.
$$

\bigskip 
\noindent Question (c) Suppose that $\phi$ and $\psi$ are both $m$-singular. Thus, there exist $N_i \in \mathcal{A}$ such that

$$
m(N_1)=m(N_2)=0 \text{ and } (\forall A \in \mathcal{A}, \ \phi(AN_1^c)=\psi(AN_2^c)=0).
$$

\bigskip 
\noindent Set $N=N_1 \cup N_2$. We have $N \in \mathcal{N}$ and $m(N)=0$. But $N_1 \subset N \Rightarrow N^c\subset N_1^c$. Hence, for any $A \in \mathcal{A}$,

$$
0\leq \phi(AN^c)\leq \phi(AN_1^c)=0 \text{ and } 0\leq \psi(AN^c)\leq \psi(AN_2^c)=0.
$$

\bigskip 
\noindent It follows that for any $A \in \mathcal{A}$, $\Phi(AN^c)=0$. $\square$.\\

\noindent Question (c) Let $\phi$ is both $m$-continuous and  $m$-singular. There exists $N \in \mathcal{A}$ such that for all $A \in \mathcal{A}$, $\phi(AN^c)=0$. But $\phi(N)=0$ since
$\phi$ is $m$-continuous and $m(N)=$. Hence, for all $A \in \mathcal{A}$,

$$
\phi(A)=\phi(AN)+\phi(AN^c)=0. \ \square.
$$

\bigskip \noindent \textbf{Exercise 3}. \label{exercise03_sol_doc08-05}\\

\noindent Suppose that the measure $m$ is finite. Let $\phi$ another finite measure. Show that there exist a \textbf{finite} $m$-continuous measure $\phi_c$ and a \textbf{finite} $m$-singular measure $\phi_s$ such that $\phi=\phi_c+\phi_s$. Show that $\phi_c$ is of the form

$$
\phi_c(A)=\int_A X \ dm\leq \phi(A), \ A \in \mathcal{A},
$$
 
\bigskip 
\noindent and, hence, $X$ is $m$-\textit{a.e.}.\\

\noindent \textit{Hints}. Proceed as follows.\\

\noindent Question (1). Define $\Phi$ as the set of all Borel and non-negative applications $X : (\Omega,\mathcal{A}) \rightarrow \overline{\mathbb{R}}$ such that for all $A\in \mathcal{A}$, we have 
$$
\int_A X \ dm\leq \phi(A).
$$

\bigskip 
\noindent (a) Justify the existence and the finiteness of $\alpha=\sup \{\int_A X \ dm, \ X \in \Phi\} \in \mathbb{R}$, that find a least one element of the $\{\int_A X \ dm, \ X \in \Phi\}$ and exhibit a finite bound above.\\

\noindent (b) Consider a maximizing sequences $(X_n)_{n\geq 0} \subset \Phi$ such that $\int_A X \ dm \rightarrow \alpha$ as $n\rightarrow +\infty$. Justify that $Y_n=\max_{1\leq k\leq n} X_k$ has a limit $X$, a non-negative Borel application.\\

\noindent (c) Do you have $X \in \phi?$ Proceed as follows : For $n\geq 0$, denote $A_k=(Y_n=X_k)$ for $1\leq k\leq n$. Show that

$$
\Omega = \bigcap_{1\leq k \leq n} A_k.
$$

\bigskip 
\noindent Next, use the now well known technique by which one transforms the union of the $A_k$'s into a sum of sets (See Exercise 4, page \pageref{exercise03_sol_doc00-03}, in Doc 00-03 in Chapter \ref{00_sets} by taking

$$
B_1=A_1, \ B_2=A_1+...+A_{k-1}^c A_k
$$

\bigskip 
\noindent to get

$$
\Omega = \bigcap_{1\leq k \leq n} A_k = \sum_{1\leq k \leq n} B_k,
$$

\bigskip 
\noindent from which, you have for all $A\in \mathcal{A}$, is derived

$$
A=A\cap \Omega = \sum_{1\leq k \leq n} AB_k, \text{ and } Y_n=X_k \text{ on } B_k, \ 1\leq k \leq n. 
$$

\bigskip 
\noindent Question : Form the following formula which is valid for all $A\in \mathcal{A}$,

$$
\int_A Y_n \ dm = \int_{sum_{1\leq k \leq n} AB_k} Y_n \ dm.
$$

\bigskip 
\noindent Show that  $\int_A Y_n \ dm\leq \phi(A)$ for all $n\geq 0$ and deduce from this, that you have  $\int_A X \ dm\leq \phi(A)$. Conclude.

\noindent Conclude that

$$
\alpha = \int X \ dm.
$$

\bigskip 
\noindent Question (2). Take the indefinite integral with respect to $X$

$$
\mathcal{A} \ni A \mapsto \phi_c(A)=\int_A X\ dm,
$$

\bigskip 
\noindent and define, based on the fact that $\phi$ is finite,

$$
\phi_s=\phi - \phi_c
$$

\bigskip 
\noindent (a) Use the conclusion of Question (1) to show that $\phi_s$ is non-negative and thus is a measure.\\

\noindent (b) Show that $\phi_s$ is $m$-singular by using Exercise 1 and proceeding as follows.\\

\noindent For each $n\geq 1$, consider the $\sigma$-finite application, well-defined since $m$ is finite, 
$$
\psi_n=\phi_s - \frac{1}{n}m
$$

\bigskip 
\noindent Denote, following Exercise 1, for each $n\geq 1$ the measurable set $D_n$ such $\psi_n(D_n)=\inf \phi_n$ which satisfies : for all $A\in \mathcal{A}$,

$$
\psi_n(AD_n)\leq 0 \ (C1) \text{  and  } \psi_n(AD_n^c)\geq 0 \ (C2).
$$

\bigskip 
\noindent Denote

$$
N=\bigcap_{n\geq 1} D_n^c \Leftrightarrow N^c=\bigcup_{n\geq 1} D_n.
$$

\bigskip 
\noindent Questions : \\

\noindent (A) From Part (C1) of the last equation and from the non-negativity of $\phi_s$, derive a bound of $\phi_s(AD_n)$ for each $n$ and from this, show that for all $A\in \mathcal{A}$,

$$
\phi_s(AN^c)=0.
$$

\bigskip 
\noindent (B) Justify each line of the following formula for $n\geq 1$ fixed and for all $A\in \mathcal{A}$ :

\begin{eqnarray*}
\int_A \biggr( X + \frac{1}{n}1_{AD_n^c}\biggr) \ dm &=&\phi_c(A)+\frac{1}{n}m(AD_n^c)\\
&=&\phi(A)-\phi_s(A)+\frac{1}{n}m(AD_n^c)\\
&=&\phi(A)-\biggr( \phi_s(AD_n^)+\frac{1}{n}m(AD_n^c) \biggr) +\biggr(\phi_s(AD_n^c)-\phi_s(A)\biggr)\\
&=&\phi(A)- \psi_n(AD_n^c) + \biggr(\phi_s(AD_n^c)-\phi_s(A)\biggr)\\
&\leq &\phi(A),
 \end{eqnarray*}

\noindent Deduce that 

$$
X + \frac{1}{n}1_{AD_n^c} \in \Phi,
$$

\bigskip 
\noindent and from this, deduce that for all $n\geq 1$, $m(D_n^c)=0$. Conclude that $m(N)=0$.\\

\noindent Conclude that $\phi_s$ is $m$-singular.\\

\noindent Make a general conclusion.\\

\bigskip \noindent \textbf{SOLUTIONS}.\\

\noindent Question (1).\\

\noindent (a). $X=0$ obviously belongs to $\{\int_A X \ dm, \ X \in \Phi\}$. All elements of $\{\int_A X \ dm, \ X \in \Phi\}$ are bounded $\phi(\Omega)<+\infty$. Thus the supremum $\alpha=\sup \{\int_A X \ dm, \ X \in \Phi\}$ exists in $\mathbb{R}$ and is bounded by $\phi(\Omega)$.\\

\noindent (b) Let $(X_n)_{n\geq 0} \subset \Phi$ be a  maximizing sequences, that $\int_A X \ dm \rightarrow \alpha$ as $n\rightarrow +\infty$.\\

\noindent The sequence $Y_n=\max_{1\leq k\leq n} X_k$, $n\geq 0$, is non-decreasing. Thus, it has a monotone limit $X$ with values in $\in \overline{\mathbb{R}}$, which of course remains non-negative as a limit of non-negative functions and measurable.\\

\noindent (c) Obviously, for all $\omega \in \Omega$, $Y_n(\omega)$ takes one of the values $X_k(\omega)$, $k\in \{1,...,k\}$, that is $\omega \in \bigcup_{1\leq k \leq n} (Y_n=X_k)$. Hence
$$
\Omega=\bigcap_{1\leq k \leq n} A_k.
$$

\bigskip 
\noindent By well-known technique by which one transforms the union of the $A_k$'s into a sum of sets (See Exercise 4, page \pageref{exercise03_sol_doc00-03}, in Doc 00-03 in Chapter \ref{00_sets}), we have
$$
\Omega = \bigcap_{1\leq k \leq n} A_k = \sum_{1\leq k \leq n} B_k,
$$

\bigskip 
\noindent with

$$
B_1=A_1, \ B_2=A_1+...+A_{k-1}^c A_k.
$$

\bigskip 
\noindent Hence, all $A\in \mathcal{A}$, we have

$$
A=A\cap \Omega = \sum_{1\leq k \leq n} AB_k.  
$$

\bigskip 
\noindent We remark that $Y_n=X_k$ on $B_k$, $\ 1\leq k \leq n$. We get for all $A\in \mathcal{A}$,

\begin{eqnarray*}
\int_A Y_n \ dm &=& \int_{\sum_{1\leq k \leq n} AB_k} Y_n \ dm\\
&=&  \sum_{1\leq k \leq n} \int_{ AB_k} Y_n \ dm \\
&=&  \sum_{1\leq k \leq n} \int_{ AB_k} X_k \ dm\\
&\leq&  \sum_{1\leq k \leq n} \phi(AB_k)\\
&=&   \phi\left(\sum_{1\leq k \leq n} AB_k\right)\\
&=&   \phi(A).
\end{eqnarray*}

\noindent By applying the MCT to the left-hand member of the inequality as $n\rightarrow$, we get $\int_A X \ dm\leq m(A)$ for all $A\in \mathcal{A}$. Hence $X \in \phi$. Thus, by the definition of the supremum, we have

$$
\int X \ dm \leq \alpha.
$$

\bigskip \noindent But for all $n\geq 0$,

$$
\int_A X_n \leq \int_A Y_n \ dm.
$$

\bigskip 
\noindent By letting $n\rightarrow +\infty$, we get $\alpha \leq \int X \ dm$. By combining the two last results, we have $\alpha=\int X \ dm$. $\square$

\noindent Question (2). Take the indefinite integral with respect to $X$

$$
\mathcal{A} \ni \mathcal{A} \mapsto \phi_c(A)=\int_A X\ dm.
$$

\noindent (a) By Exercise 2, $\phi_c$ is a $m$-continuous measure. Since $\phi$ is finite, the application $\phi_s=\phi - \phi_c$ is well-defined on $\mathcal{A}$ and is $\sigma$-finite. But for all 
$A \in \mathcal{A}$, by definition, $\phi_c(A)=\int_A X\ dm \leq \phi(A)$, that is $\phi(A)-\phi_c(A)\geq 0$. This $\phi_s$ is non-negative and is a measure.

\noindent (b) Let us show that $\phi_s$ is $m$-singular. By Exercise 1, let is for each $n\geq 1$, consider the $\sigma$-finite application 
$$
\psi_n=\phi_s - \frac{1}{n}m.
$$

\bigskip 
\noindent Since $m$ is finite, we will not face the the pathological situation $+\infty - \infty$. So $\psi_n$ is well-defined. By Exercise 1, for each $n\geq 1$, we have a measurable set $D_n$ such $\psi_n(D_n)=\inf \phi_n$ which satisfies : for all $A\in \mathcal{A}$, $\psi_n(AD_n)\leq 0$ and $\psi_n(AD_n^c)\geq 0$.\\

\noindent (A) From the fact that $\psi_n(AD_n)\leq 0$ for $A\in \mathcal{A}$, we derive that for all $n\geq 1$, for all $A\in \mathcal{A}$,

$$
\psi_n(AD_n)=\phi_s(AD_n) - \frac{1}{n}m(AD_n) \leq 0 \Rightarrow 0\leq \phi_s(AD_n) \leq \frac{1}{n}m(\Omega). 
$$

\bigskip 
\noindent By taking $N^c=\bigcap_{n\geq 1} D_n$, we have for all $n\geq 1$, for all $A\in \mathcal{A}$,

$$
\phi(AN^c) \leq \phi(AD_n) \leq \frac{1}{n}m(\Omega).
$$

\bigskip 
\noindent and by letting $n\rightarrow +\infty$, (remind that the measure $m$ us finite), we have $\phi(AN^c)=0$ for all $A \in \mathcal{A}$.\\

\noindent (B) From the fact that $\psi_n(AD_n^c)\geq 0$ for $A\in \mathcal{A}$, we derive that for all $n\geq 1$, for all $A\in \mathcal{A}$,

\begin{eqnarray*}
\int_A \biggr( X + \frac{1}{n}1_{AD_n^c}\biggr) \ dm &=&\phi_c(A)+\frac{1}{n}m(AD_n^c)\\
&=&\phi(A)-\phi_s(A)+\frac{1}{n}m(AD_n^c)\\
&=&\phi(A)-\biggr( \phi_s(AD_n^)+\frac{1}{n}m(AD_n^c) \biggr) +\biggr(\phi_s(AD_n^c)-\phi_s(A)\biggr)\\
&=&\phi(A)- \psi_n(AD_n^c) + \biggr(\phi_s(AD_n^c)-\phi_s(A)\biggr)\\
&\leq &\phi(A),
 \end{eqnarray*}

\noindent since $\psi_n(AD_n^c)\geq 0$ is non-negative and that $\biggr(\phi_s(AD_n^c)-\phi_s(A)\biggr)$ is non-positive because of the non-drcreasingness of $\phi_s$. It follows that

$$
X + \frac{1}{n}1_{AD_n^c} \in \Phi,
$$

\bigskip 
\noindent which implies that for all $n\geq 1$ (for $A=\Omega)$

$$
\int \biggr( X + \frac{1}{n}1_{AD_n^c}\biggr) \ dm=\alpha + \frac{1}{n}m(D_n^c)\leq \alpha.
$$

\bigskip 
\noindent Since $\alpha$ is finite and $m$ is non-negative, this ensures that $m(D_n^c)=0$ for all $n\geq 1$. This

$$
m\left(\bigcup_{n\geq 1} D_n^c \right)=:m(N)=0.
$$

\bigskip 
\noindent We conclude that $\phi_s$ is singular and we have the desired decomposition : $\phi=\phi_c+\phi_s$. We already know that $\phi_c$ is bounded by $\phi(\Omega)$ which is finite. Hence 
$\psi=\phi-\phi_c$ is also finite.\\

\bigskip \noindent \textbf{Exercise 4}. \label{exercise04_sol_doc08-05}\\

\noindent Generalize the results of Exercise 3 to $\sigma$-finite measures $\phi$ and $m$ finite in the following way :  There exists $m$-\textit{a.e.} finite non-negative application $X$ and a \textbf{finite } $m$-singular measure $\phi_s$ such that for all $A \in \mathcal{A}$,

$$
\phi(A)=\int_A X \ dm + \phi_s(A).
$$

\bigskip 
\noindent \textit{Recommendations}. It is advised in a first reading, just to follow the solution of this extension. The ideas are the following : We consider a measurable partition of $\Omega$ : $\Omega=\sum_{j\geq 0}$, $(\forall j\geq 0, \ \Omega_j \in \mathcal{A})$ such that $\phi(\Omega_j)<+\infty$ and $m(\Omega_j)<+\infty$. Hence, we have for all $A \in \mathcal{A}$

$$
\phi = \sum_{j \geq 0} \phi_j \ (S1), \ \text{and} \ m = \sum_{j \geq 0} m_j, \ (S2)
$$

\bigskip 
\noindent where the applications

$$
\mathcal{B} \ni B \mapsto \phi_j(A)=\phi(A\cap \Omega_j) \text{ and } \mathcal{B} \ni B \mapsto m_j(A)=m(A\cap \Omega_j),
$$

\bigskip 
\noindent are measures which can be considered as defined on $\mathcal{A}$ (but with $\Omega_j$ as support) or, on the $(\Omega_j, \mathcal{A}_j)$, where $\mathcal{A}_j)$ is the induced 
$\sigma$-algebra of $\mathcal{A}$ on $\Omega_j$.\\

\noindent Hence on each $\Omega_j$, $j\geq 1$, $\phi_j$ and $m_j$ are finite and we get, by Exercise 3, a decomposition

$$
\phi_j(B)=\int_B X_j \ dm + \phi_{s,j}(B), \ B \in \mathcal{A}_j, \ (DE01)
$$

\bigskip 
\noindent From there, we sum members of Formula (DE01). But some technicalities are needed.\\

\bigskip \noindent \textbf{SOLUTIONS}.\\

\noindent Let $\phi$ be a non-negative $\sigma$-additive application. Let us consider a measurable partition of $\Omega$ : $\Omega=\sum_{j\geq 0}$, $(\forall j\geq 0, \ \Omega_j \in \mathcal{A})$ such that
$\phi(\Omega_j)<+\infty$ and $m(\Omega_j)<+\infty$. Hence, we have for all $A \in \mathcal{A}$

$$
\phi = \sum_{j \geq 0} \phi_j \ (S1), \ \text{and} \ m = \sum_{j \geq 0} m_j, \ (S2)
$$

\bigskip 
\noindent where the applications

$$
\mathcal{B} \ni B \mapsto \phi_j(A)=\phi(A\cap \Omega_j) \text{ and } \mathcal{B} \ni B \mapsto m_j(A)=m(A\cap \Omega_j),
$$

\bigskip 
\noindent are measures which can be considered as defined on $\mathcal{A}$ (but with $\Omega_j$ as support) or, on the $(\Omega_j, \mathcal{A}_j)$, where $\mathcal{A}_j)$ is the induced 
$\sigma$-algebra of $\mathcal{A}$ on $\Omega_j$.\\

\noindent Hence on each $\Omega_j$, $j\geq 1$, $\phi_j$ and $m_j$ are finite and we get, by Exercise 3, a decomposition

$$
\phi_j(B)=\int_B X_j \ dm + \phi_{s,j}(B), \ B \in \mathcal{A}_j, \ (DE01)
$$

\bigskip 
\noindent where $X_j$ is $m_j$-\textit{a.e.} finite and measurable from $(\Omega_j, \mathcal{A}_j)$ to $\overline{\mathbb{R}}$ and $\phi_{s,j}$ is finite and $m_j$-singular on $\Omega_j$. Define piece-wisely,

$$
X = \sum_{j\geq 0} X_j 1_{\Omega_j}.
$$

\bigskip 
\noindent Formula (DE01) is also

$$
\phi_j(A)=\phi(A\cap \Omega_j)=\int_{(A\cap \Omega_j)} X_j \ dm + \\phi_{s,j}(A\cap \Omega_j), \ B \in \mathcal{A}_j, \ (DE02)
$$

\bigskip 
\noindent The application is measurable since for all $A\in \mathcal{A}$,

\begin{eqnarray*}
X^-1(A)&=&X^-1\left(\sum_{j\geq 0} A\cup \Omega_j \right)\\
&=&\sum_{j\geq 0} X^-1\left( A\cup \Omega_j \right)\\
&=&\sum_{j\geq 0} X_j^-1\left( A\cup \Omega_j \right).
\end{eqnarray*}

\noindent The last line is a sum of elements of $\mathcal{A}_j$'s which are sub-collections of $\mathcal{A}$.\\

\noindent $X$ is $m$-\textit{a.e.} finite since

$$
(X \text{ infinite })= \bigcup_{j\geq 0} (X_j \text{ infinite })
$$ 

\bigskip 
\noindent and hence

$$
m(X \text{ infinite })\leq \sum_{j\geq 0} m(X_j \text{ infinite })=\sum_{j\geq 0} m_j(X_j \text{ infinite })=0.
$$ 

\bigskip 
\noindent Further By Question (5) in Exercise 1 on Doc 06-10 (page \pageref{exercise01_sol_doc06-10}) of Chapter \ref{06_convergence},

\begin{eqnarray*}
\int_A X \ dm&=&\sum_{h\geq 0} \int 1_A X dm_j\\
&=&\sum_{h\geq 0} \int 1_A \biggr( \sum_{j\geq 0} 1_{\Omega_j} X_j\biggr) dm_h\\
&=&\sum_{h\geq 0} \int 1_A \biggr( \sum_{j\geq 0} 1_{\Omega_j} X_j\biggr) dm_h\\
&=&\sum_{h\geq 0} \sum_{j\geq 0} \int  1_{A\cap \Omega_j} X_j dm_h.
\end{eqnarray*}

\noindent For $h\neq j$, the integral $\int  1_{A\cap \Omega_j} X_j dm_h$ is zero since the measure is supported by $\Omega_h$ on which the function $1_{A\cap \Omega_j} X_j$ is zero. Hence

$$
\int_A X \ dm=\sum_{j\geq 0} \int_{A \cap \Omega_j} X_j dm_j. \ (DE03)
$$

\bigskip 
\noindent By combining (DE02) and (DE03), and by summing over $j\geq 1$, we get for all $A \in \mathcal{A}$

$$
\phi(A)=\int_A X \ dm + \phi_s(A),
$$

\bigskip 
\noindent where, for all $A \in \mathcal{A}$,

$$
\phi_s(A)=\sum_{j \geq 0}  \phi_{s,j}(A\cap \Omega_j).
$$

\bigskip 
\noindent By taking $N=\sum_{j\geq 0} N_j$, we get that $N$ is still a $m$-null set since

$$
m(N)=\sum_{j\geq 0}  m(N\cap \Omega_j)=\sum_{j\geq 0}  m_j(N\cap \Omega_j)=\sum_{j\geq 0}  m_j(N_j)=0,
$$

\bigskip 
\noindent and for $A \in \mathcal{A}$,

$$
\phi_s(AN^c)=\sum_{j \geq 0}  \phi_{s,j}(N^c\cap \Omega_j)=0,
$$

\bigskip 
\noindent where we use the fact that each $N^c\cap \Omega_j$ is the complement of $N_j$ in $\Omega_j$.\\
 
\noindent We get the searched decomposition with the attached conditions.\\

\bigskip \noindent \textbf{Exercise 5}. (General form of the Lebesgue Decomposition) \label{exercise06_sol_doc08-05}\\

\noindent Generalize the results of Exercise 4 to a $\sigma$-finite and $\sigma$-additive application $\phi$ and a $\sigma$-finite measures $m$.\\

\noindent Question (a) Show that there exists an $m$-\textit{a.e.} finite application $X$ and a difference of \textbf{finite} $m$-singular measures $\phi_s$ such that for all $A \in \mathcal{A}$,

$$
\phi(A)=\int_A X \ dm + \phi_s(A).
$$

\bigskip 
\noindent Question (b) Show that the decomposition is unique if $\phi_s(A)$ is a difference of $m$-singular measures.\\

\noindent \textit{Hints} Combine Exercises 1  and 6.\\

\bigskip \noindent \textbf{SOLUTIONS}.\\

\noindent Question (a). \noindent By Exercise 2, $\phi$ is difference of two measures $\phi_i$, $i=1,2$, which can be expressed as $\phi_1(A)=-\phi(AD)$ and $\phi_2(A)=\phi(AD^c)$ for $A\in \mathcal{A}$. Based on these expressions, we see that $\phi_1$ and $\phi_2$ are $\sigma$-finite of $\phi$ is.\\

\noindent One of $\phi_1$ and $\phi_1$ is finite at least. If $\phi$ does not take the value $-\infty$ for example, $\phi_1$ is finite. By Exercise 6, there exist  $m$-continuous measures $\phi_{i,c}$ and \textbf{finite} $m$-singular measures $\phi_{i,s}$, $i=1,2$, such that

$$
\phi_i=\phi_{i,c}+\phi_{i,s}, \ i=1,2.
$$  

\bigskip 
\noindent By the finiteness of $\phi_2$, it is clear that $\phi_{2,c}$ and $\phi_{2,s}$ are finite and we may get the difference

$$
\phi =(\phi_{1,c} -\phi_{2,c})+(\phi_{1,s} -\phi_{2,s}).
$$ 

\bigskip 
\noindent , which by Exercise 2 is a sum of an $m$-continuous application and an $m$-singular application.\\

\noindent Question (b) Consider the constructed decomposition $\phi=\phi_c + \phi_s$ with $\phi_s$ is a difference of finite $m$-singular measures. Since $\phi_c$ is an indefinite integral associated to a $m$-\textit{a.e.} application, it is $\sigma$-finite. Consider the measurable partition $\Omega=\sum_{j\geq 0} \Omega_j$, $(\forall j\geq, \ \Omega_j \in \mathcal{A} \text{ and } \phi_c(\Omega_j)<+\infty)$.\\

\noindent Now suppose we have another decomposition $\phi=\psi_c + \psi_s$ where $\psi_c$ is $m$-continuous and $\psi_s$ is a difference of $m$-singular measures and finite. By Exercise 2, $\psi_2-\phi_2$ is still a difference $m$-singular measure on $\mathcal{A}$. Thus, for any $A\in \mathcal{A}$, we have

$$
\psi_c(A\cap \Omega_j) + \psi_s(A\cap \Omega_j)=\phi_c(A\cap \Omega_j) + \psi_s(A\cap \Omega_j).
$$

\bigskip 
\noindent Hence, we may move the finite numbers from one side to another to get 

$$
\psi_c(A\cap \Omega_j)-\phi_c(A\cap \Omega_j) =\psi_s(A\cap \Omega_j)-\phi_s(A\cap \Omega_j).
$$

\bigskip \noindent Thus the application $\psi_c(A\cap \Omega_j)-\phi_c(A\cap \Omega_j)$ is both $m$-continuous and $m$-singular on $\Omega_j$and hence is the null application. This  $\psi_c$ and 
$\phi_c$ on one side and, $\psi_s$ and $\phi_c$ on another, are equal on each $\mathcal{A}_{\Omega_j}$. So they are equal. The decomposition is unique.\\

\bigskip \noindent \textbf{Exercise 6}. (Theorem of Radon-Nikodym) \label{exercise07_sol_doc08-05}\\

\noindent Question (a) Let $\phi$ be a $\sigma$-finite and $\sigma$-additive application and a $\sigma$-finite measures $m$. Suppose that $\phi$ is $m$-continuous. Then there exists a real-valued measurable $X$ wich is $m$-finite such that for $A \in \mathcal{A}$, we have

$$
\phi(A)=\int_A X \ dm.
$$

\bigskip 
\noindent Question (b) Show that $X$ is $m$-\textit{a.e.} non-negative if $\phi$ is a measure.\\

\noindent \textit{Hints} Suppose that $\phi$ is non-negative. Set for any $k\geq 1$, $A_k=(X <- 1/k)$. Show that $m(A_k)=0$ for any $k\geq 1$. Deduce that $m(X\leq 0)=0$.\\

\bigskip \noindent \textbf{SOLUTIONS}.\\

\noindent Question (a) By Exercise 5, we have a decomposition
$$
\phi(A)=\int_A X \ dm + \phi_s(A), \ A\in \mathcal{A}.
$$

\bigskip 
\noindent where $X$ is $m$-\textit{a.e.} finite and where $\phi_s$ is singular. Consider a measurable decomposition like in the solution of Question (b) of Exercise 5. We have $A\in \mathcal{A}$, we have

$$
\phi(A\cap \Omega_j)=\phi_c(A\cap \Omega_j) + \phi_s(A\cap \Omega_j), 
$$

\bigskip 
\noindent which implies

$$
\phi(A\cap \Omega_j)-\phi_c(A\cap \Omega_j)= \phi_s(A\cap \Omega_j).
$$

\bigskip 
\noindent Hence, on each $\Omega_j$, we have an application which is both $m$-continuous and $m$-singular. Hence $\phi_s$ is the null function on each $\mathcal{A}_{\Omega_j}$, $j\geq 0$. Hence
$\phi_s$ is the null function. The solution is complete. $\blacksquare$.\\

\bigskip Question (b) Suppose that $\phi$ is non-negative. Let us denote $A_k=(X <- 1/k)$, $k\geq 0$. We have for any $k\geq 1$,

$$
0\leq \phi(A_k)=\int_{A_k} X \ dm \leq -\frac{1}{k} m(A_k). 
$$

\noindent This is possible only if $m(A_k)=0$ for any $k\geq 1$. We conclude that $m(X \leq 0)=0$ since $(X <0)=A_k$.\\

\bigskip \noindent \textbf{Exercise 7}. (Extended Theorem of Radon-Nikodym) \label{exercise07_sol_doc08-05}\\

\noindent Let $\phi$ be a $\sigma$-additive application defined on $(\Omega, \mathbb{A})$, \textit{not necessarily $\sigma$-finite}, and let $m$ be $\sigma$-finite measure on $(\Omega, \mathcal{A})$   such that $\phi$ is continuous with respect to $m$. Show that there exists a measurable application $X:(\Omega ,\mathcal{A})\mapsto \overline{\mathbb{R}}$ \textit{not necessarily $m$-\textit{a.e.} finite} such that we have  

$$
\forall A \in \mathcal{A}, \ \phi(A)=\int_{A}  X \ dm.
$$

\bigskip 
\noindent \textit{Hints}. Show it only for $\phi$ is non-negative and $m$ finite. First for $\phi$ non-negative, $X_j$ will be found on a countable partition on $\Omega$ piece-wisely as in the solution of Exercise 6. For a general $\sigma$-additive application, the extension is done by using the difference of two measures, one of them being finite. Based on these remarks, proceed as follows.\\

\noindent So $\phi$ be a $m$-continuous measure, where $m$ is a finite measure.\\ 

\noindent Question (a). Consider the class of elements of $A \in \mathcal{A}_0$ such that $\phi$ is $\sigma$-finte of the induce $\sigma$-algebra $\mathcal{A}_A$ and a sequence $(B_n)_{n\geq 0} \subset \mathcal{A}_0$ such that $m(B_n)\rightarrow  s=\sup_{B\in \mathcal{A}_0} m(B)$. Define $\bigcup_{k=0}^{+\infty} B_k$ and define also $\mathcal{A}_1=\{AB, \ A\in \mathcal{A}\}$ and $\mathcal{A}_2=\{AB^c, \ A\in \mathcal{A}\}$.\\

\noindent Show that for all $n\geq 0$, $\bigcup_{k=0}^{n} B_k \in \mathcal{A}_0$ and deduce that $B \in \mathcal{A}_0$.\\

\noindent Apply the MCT continuity from above of $m$  to $\bigcup_{k=0}^{n} B_k \nearrow B$ as $n\nearrow +\infty$ to show that $s=m(B)$.\\

\noindent Question (b)  Show that if for some $C\in \mathcal{A}_2$ we have $0<\phi(C)<+\infty$, then $s$ would not be an infimum. Conclude that we have the equality 

$$
(C \in \mathcal{A}_2, \ \phi(C) \text{ finite }) \Rightarrow m(C)=0. (E1)
$$

\bigskip 
\noindent Deduce from this that $\phi$ is the indefinite integral of $X_2=+\infty$ of $\mathcal{A}_2$.\\

\noindent Use the convention of the integration on the class of elementary function that $\int_A \infty \ dm=m(A) \times \infty=0$ if $m(A)=0$.\\ 

\noindent Question (c) By remarking that $\phi$ is $\sigma$-finite on  $\mathcal{A}_1$, apply the first version of Radon-Nikodym with an $m$-\textit{a.e.} finite derivative $X_1$ on $B$.\\

\noindent Question (d) Form $X=X_1 1_B + X_2 1_{B^c}$ and conclude.\\

\bigskip \noindent \textbf{SOLUTIONS}.\\

\noindent Question (a) The set $\mathcal{A}_0$ contains $\emptyset$ and is bounded above by the finite number $m(\Omega)$. Hence the supremum $s=\sup_{B\in \mathcal{A}_0} m(B)$ exists and is finite. We may take a maximizing sequence $(B_n)_{n\geq 0} \subset \mathcal{A}_0$ such that $m(B_n)\rightarrow s=\sup_{B\in \mathcal{A}_0} m(B)$.\\

\noindent We have for all $n\geq 2$,

$$
\phi\left( \bigcup_{k=0}^{n} B_k \right)=\sum_{k=0}^{+\infty} \phi(C_k),
$$ 

\bigskip 
\noindent where $C_0=B_0 \subset B_0$, $C1=B_0^c B_1 \subset B_1$, $C_k=B_0^c \ldots B_{k-1}^c B_k$, $k=2,...,n$. Hence for $n\geq 0$, $\bigcup_{k=0}^{n} B_k \in \mathcal{A}_0$. This implied on one side that $\phi$ is $\sigma$-finite on $B$ and, on another side, 

$$
\forall n\geq 0, \ m(A_n) \leq m\left(\bigcup_{k=0}^{n} B_k\right) \leq s.
$$ 

\bigskip 
\noindent By letting $n\nearrow+\infty$ in the formula above and by using the continuity of $m$, we get $s=m(B)$.\\

\noindent Question (b) Suppose that for some $C=AB^c \in \mathcal{A}_1$, $A\in \mathcal{A}$, we have $\phi(C)<+\infty$. The inequality $\phi(C)<+\infty$ implies that $\phi$ is still in $\sigma$-finite on $B+C$, and thus we would have $m(B+C)=m(B)+m(C)\leq \alpha$, that is $\alpha+m(C)\leq \alpha$, which is impossible unless $m(C)=0$. So we get

$$
(\forall C \in \mathcal{A}_2),  \phi(C) \text{ finite } \Rightarrow m(C)=0,
$$

\bigskip 
\noindent that is

$$
\forall C \in \mathcal{A}_2, \ m(C)>0 \Rightarrow \phi(C)=+\infty. \ E2)
$$

\bigskip \noindent Combine this with $m$-continuity of $\phi$ with respect to $m$, to check the following, for all $C \in \mathcal{A}_2$\\

\noindent (i) If $m(C)>0$, then $\phi(C)=+\infty$ and $\int_C +\infty \ dm=m(C) \times +\infty$.\\

\noindent (ii) If $m(C)=0$, then $\phi(C)=0$ (by $m$-continuity of $\phi$ with respect to $m$) and $\int_C +\infty \ dm= m(C) \times +\infty=0$.\\

\noindent We conclude that for $C \in \mathcal{A}_2$ : \\

$$
\phi(C)=\int_C +\infty \ dm.
$$

\bigskip \noindent (c) On $\mathcal{A}_1$, apply the first version of Radon-Nikodym to get an $m$-finite and measurable application $X_1$ defined on $B$ such that for all $A\in \mathcal{A}$,

$$
\phi(AB)=\int_{AB} X_1 \ dm.
$$

\bigskip \noindent (c) By forming $X=X_1 1_B + X_2 1_{B^c}$, we get for all $A\in \mathcal{A}$,

$$
\phi(A)=\phi(AB)+\phi(AB^c)=\int_{AB} X_1 \ dm + \int_{AB^c} X_2 =\int_A X \ dm. \ \blacksquare
$$

\bigskip 
\bigskip \noindent \textbf{Exercise 8}. (Additivity of the Radon-Nikodym derivatives) \label{exercise08_sol_doc08-05}\\

\noindent Let $(\phi_i)_{i\in I}$, $I \subset \mathbb{N}$, be a countably family of $\sigma$-additive applications defined on $(\Omega, \mathbb{A})$, \textit{not necessarily $\sigma$-finite}, and let $m$ be $\sigma$-finite measure on $(\Omega, \mathcal{A})$ such that each $\phi_i$, $i\in I$, is continuous with respect to $m$. Suppose that one of thes value $+\infty$ or $-\infty$ is simultaneously excluded for all the members to ensure that the application

$$
\phi=\sum_{i\in I} \phi_i,
$$

\bigskip \noindent is well-defined. Of course, by the previous exercises, $\phi$ is $\sigma$-additive and is $m$-continuous. This exercise investigates the validity of the Formula

$$
\frac{d\phi}{dm}=\sum_{i \in I} \frac{d\phi_i}{dm}. (ARD)
$$  

\bigskip \noindent Question (a) Let $I$ be finite. Show that $X$ has Radom-Nikodym derivative with respect to $m$ and Formula (ARD) holds.\\

\noindent \textit{Hints}. Evoke Exercise 7 and consider the Radon-derivatives $X_i=\frac{d\phi_i}{dm}$, $i\in I$. Denote $N_i=(X_i=+\infty)$, $i\in I$. What do you know about the  $m$-\textit{a.e.} finiteness of any $X_i$ and that of 

$$
X=\left(\sum_{i\in I} X_i\right) 
$$

\bigskip \noindent on $N=\left(\bigcup N_i\right)^c$. What is the value of 

$$
1_{N^c} \left(\sum_{i\in I} X_i\right) \ ? 
$$

\bigskip \noindent Justify each line of the formula below. Actually, we only need to justify Line (S). 

\begin{eqnarray*}
\phi(A)&=&\sum_{i\in I} \phi_i\\
&=&\sum_{i\in I} \int_A 1_N X_i \ dm + \int_A 1_{N^c} X_i \ dm\\
&=& \int_A 1_N \left(\sum_{i\in I} X_i\right) \ dm + \int_A 1_{N^c} \left(\sum_{i\in I} X_i\right) \ dm. \text{ (S) }\\
\end{eqnarray*} 

\noindent Conclude.\\

\noindent Question (c) Suppose that the elements of $\mathcal{A}$ are non-negative and $I$ is infinite. Use the same method and justify Line (S) and conclude.\\ 

\noindent Question (d) Give a general condition for the validity of the Step (S) in the hints. Conclude by separating two cases : $m(N)=0$ and $m(N)>0$.\\

\bigskip \noindent \textbf{SOLUTIONS}.\\

\noindent Question (a).\\

\noindent By Exercise 7, there exist Radon-derivatives $X_i=\frac{d\phi_i}{dm}$, $i\in I$. Let us denote $N_i=(X_i=+\infty)$, $i\in I$. We have, by Exercise 7, that all the $X_i$ are $m$-\textit{a.e.} on
finite $N=\left(\bigcup N_i\right)^c$. Reminding that each $X_i$ takes $m$-\textit{a.e.} a finite value on $N_i^c$ and the value $+\infty$ on $N_i$ for any $i\in I$, we have that

$$
1_{N^c} \left(\sum_{i\in I} X_i\right)=+\infty. 
$$

\bigskip \noindent  and
$$
X=\left(\sum_{i\in I} X_i\right) 
$$

\bigskip \noindent is finite on $N$. \noindent Hence, for all $A\in \mathcal{A}$, we have

\begin{eqnarray*}
\phi(A)&=&\sum_{i\in I} \phi_i\\
&=&\sum_{i\in I} \int_A 1_N X_i \ dm + \int_A 1_{N^c} X_i \ dm\\
&=& \int_A 1_N \left(\sum_{i\in I} X_i\right) \ dm + \int_A 1_{N^c} \left(\sum_{i\in I} X_i\right) \ dm \text{ (S) }\\
&=& \int_A A_N X \ dm + \int_A 1_{N^c} (+\infty) \ dm.
\end{eqnarray*} 

\noindent Thus $X 1_N + (+\infty) 1_{N^c}$ is a Radon-derivative of $\phi$ and Formula (ARD) holds.\\

\noindent Question (b) By using the same method, we justify Line (S) by the MCT applied to a the integral of a sum non-negative measurable functions. We conclude similarly.

\noindent Question (c) In the general case, we use the DCT applied on the integral of a sum integrable functions. We say that if

$$
\sum_{i\in I} \int_A 1_N |X_i| \ dm<+\infty,
$$

\bigskip \noindent we have

$$
\sum_{i\in I} \int_A 1_N X_i \ dm=\int_A 1_N \left(\sum_{i\in I} X_i\right) \ dm.
$$

\bigskip \noindent To ensure that 

$$
\sum_{i\in I} \int_A 1_{N^c} X_i \ dm=\int_A 1_{N^c} (+\infty) \ dm.
$$

\bigskip \noindent we apply again the MCT. Hence, we are able to conclude.\\

\noindent \LARGE \textbf{DOC 08-06 : Exercise on a continuous family of measures} \label{doc08-06}\\
\bigskip
\Large

\noindent \textbf{NB}. This exercise is inspired by Lemma 2.1 and its proof in \cite{shao}, page 104. This result is important in Mathematical Statistics regarding the factorization theorem for sufficient statistics.\\

\noindent \textbf{Dominated family of $\sigma$-additives applications}. Let $\mathcal{M}$ be a family of real-valued $\sigma$-additive and non-null applications defined on $(\Omega, \mathcal{A})$ and let $m$ be a  measure defined on the same space. For a coherence purpose, let us suppose \\ 

\noindent The family $\mathcal{M}$ is dominated by $m$ if and only any element of $\mathcal{M}$ is $m$-continuous. The family $\mathcal{M}$ is dominated if and only if there is a measure on 
$(\Omega, \mathcal{A})$ which dominates it.\\ 

\noindent \textbf{NB}. In all this document, $\mathcal{M}$ is a family of real-valued $\sigma$-additive applications defined on $(\Omega, \mathcal{A})$ and $m$ a measure defined on the same space.\\

\noindent \textbf{Summary}. The objective of the exercise is to show this.\\

\noindent (A) If $m$ is finite, $\mathcal{A}$ is dominated by a countable linear combination of elements of its elements whose set of coefficient is bounded.\\

\noindent (B) if $m$ is $\sigma$-finite, $\mathcal{M}$ is dominated by a countable linear combination of elements of $\mathcal{M}$ supported by $m$-finite subspaces and is associated to a bounded set of coefficients.\\

\bigskip \noindent \textbf{Exercise 01}. \label{exercise01_doc08-06}\\

\noindent Consider a family $\mathcal{M}$ of $\sigma$-finite measures defined on $(\Omega, \mathcal{A},)$ and $m$ another $\sigma$-finite measures defined on the same space.\\

\noindent Question (a) Suppose that the family is countable, that is $\mathcal{M}=\{\phi_i, \ i\in I\}$, $I \subset \mathcal{A}$. For any $c=(c_i)_{i\in I} \in (]0,+\infty[)^I$, $\mathcal{M}$ is dominated by

$$
\phi_c=\sum_{i \in I} c_i \phi. \ (FC01)
$$ 

\noindent Question (c) Let us remain in the general case and let us suppose that \textbf{$m$ is finite}. Define for any $B \in \mathcal{A}$ with $m(B)>0$ :

$$
\mathcal{D}_B=\{C \in \mathcal{A}, \ C \in B, m(C)>0\}.
$$

\noindent Define (where $\neq$ is the negation logical operator) :

$$
\mathcal{A}_0=\{B \in \mathcal{A}, \ \exists \phi \in \mathcal{M} \phi(B)>0 \text{ and for any } C \in \mathcal{D}_B, \neg (1_B \frac{d\phi}{dm}=0 \ m-\textit{a.e})\}
$$

\noindent Show that class is not empty.\\

\noindent \textit{Hint}. Consider $B \notin \mathcal{A}_0$ and show that we cannot have $\phi(B)>0$. For this, suppose $\phi(B)>0$ and deduce that $m(B)=0$. For each $C \in \mathcal{D}_B$, use the formula 
$$
\phi(C)=\int 1_C \frac{d\phi}{dm} \ dm=0,
$$

\noindent show that assertion $\neg (1_C \frac{d\phi}{dm}=0 \ $m$-\textit{a.e})$ holds.\\

\noindent Question (b) Suppose that for all $(\phi_1, \phi_2) \in \mathcal{M}^2$, there exists a linear combination  $c_1 \phi_1 + c_2 \phi_2 \in \mathcal{M}$ with $c_i\neq 0$, $i\in \{1,2\}$. Show that $\mathcal{A}_0$ is stable under finite unions.\\

\noindent \textit{Hints}. Consider $(B_1,B_2)  \in \mathcal{A}_0^2$. For each $i \in \{1,2\}$, consider $\phi_i$ and $B_{i,0}\in \mathcal{D}_{B_i}$ such that $\phi_i(B_i)>0$ and $\neg (1_{B_{i}} \frac{d\phi}{dm}=0 \ $m$-\textit{a.e})$. Show that $\phi(B)>0$. To show that $\neg (1_C \frac{d\phi}{dm}=0 \ $m$-\textit{a.e})$ for all $C \in \mathcal{D}_B$, use Point (08.07c)

$$
\frac{d\phi}{dm}=c_1 \frac{d\phi_1}{dm}+c_2 \frac{d\phi_2}{dm},
$$ 
 
\noindent and by the non-negativity of the derivatives to show that the assertion : (There exists $C \in \mathcal{D}_B$, :$(1_C \frac{d\phi}{dm}=0 \ $m$-\textit{a.e})$), leads to contradiction related to 
$C\cap B_1$ or $C\cap B_2$ and conclude.\\

\noindent Question (c) Suppose that the family $\mathcal{M}$ is bounded by a finite number $M$, that is

$$
\sup_{A \in \mathcal{A}, \ \phi \in \mathcal{M}} \phi(B) \leq M <+\infty.
$$

\noindent Define for each $I \subset \mathcal{N}$, $I\neq \emptyset$,

$$
\mathcal{C}_\mathcal{N}=\{c=(c_i)_{i\in \mathcal{N}} \in ([0,+\infty[)^\mathcal{N} \ : \|c\|_{\infty}=\sup_{i\in \mathcal{N}} c_i<+\infty \text{ and } \|c\|_{1}=\sup_{i\in \mathcal{N}} c_i>0\}.
$$

\noindent and

$$
\mathcal{M}_{0}=\{ \sum_{i\in I} c_i \phi_i, \ c\in \mathcal{C}_\mathcal{N}, \ (\phi_i)_{i\in \mathcal{N}} \in \mathcal{M}^\mathcal{N}\}.
$$

\noindent \textbf{Important remark}. Beyond the notation, we have to understand that $\mathcal{M}_{0}$ is the class of all countable combinations of elements of $\mathcal{M}$ such that for each combination, the set of the coefficients contains at least a non-zero element and is bounded in $\mathcal{R}$.\\ 

\noindent (c1) Show that $\mathcal{M}_{0}$ is still dominated by $m$.\\

\noindent (c2) Define 

$$
\mathcal{A}_0^{\ast}=\{B \in \mathcal{A}, \ \exists \phi \in \mathcal{M}_0 \phi(B)>0 \text{ and } \neg (1_B \frac{d\phi}{dm}=0 \ m-\textit{a.e}) \}.
$$

\noindent We may easily see from Question (b) (with a small extra work) that $\mathcal{A}_0^{\ast}$ is not empty. Also, it is easy to see that for $(\phi_1, \phi_2) \in \mathcal{M}_0^2$,   $\phi_1 + \phi_2 \in \mathcal{M}_0$.\\

\noindent Consider, the characterization of the supremum allows to have a sequence $(B_n)_{n\geq 0}\subset \mathcal{A}_0^{\ast}$ such that

$$
m(C_n) \rightarrow \beta=\sup_{B\in \mathcal{A}_0^{\ast}} m(B).
$$

\noindent For $n\geq 0$, consider $\phi_n\in \mathcal{M}_0$  such that

$$
\phi_n(B)>0 \text{ and } \neg (1_{B_n} \frac{d\phi_n}{dm}=0 \ m-\textit{a.e.}).
$$

\noindent Question (c2-a) Fix $(\alpha_i)_{i\in \mathcal{N}} \in \mathcal{C}_{\mathcal{N}}$, \textbf{but} such that $\alpha_i\neq 0$ for all $i\in \mathcal{N}$. Put

$$
\hat{\phi}=\sum_{n\geq 0} \alpha_n \phi_n \in \mathcal{M}_1.
$$

\textit{Hints}. For each $n\geq 0$, consider $(\phi_{i,n})_{i\in \mathbb{N}} \in \mathcal{A}^{\mathbb{N}}$ and $(c_{i,n})_{i \in \mathcal{N}} \subset \mathcal{C}_{\mathcal{N}}$ such that 

$$
\phi_n=\sum_{i\geq 0}c_{i,n} \phi_{i,n}.
$$

\noindent Writing

$$
\hat{\phi}=\sum_{n\geq 0, i\geq 0} \alpha_n c_{i,n} \phi_{i,n},
$$

\noindent Show that $\hat{\phi} \in \mathcal{M}_0$.\\
 
\noindent (c2-3) Let $B=\bigcup B_n$. Use this application $\hat{\phi}$ to show that $B \in \mathcal{A}_0^{\ast}$. Use the method of Question (c) in the context of a countable sum.\\

\noindent (c2-4) Suppose that for some $A\in \mathcal{A}$, $\hat{\phi}(A)=0$. Take an arbitrary element $\phi in \mathcal{M}_0$. Denote $C=(\frac{d\phi}{dm} > 0)$. Use the following facts to be justified :

\noindent (i) Recall that : $B \in \mathcal{A}_0$.\\

\noindent (ii)  Show that $\hat{\phi}(A \cap B)=0 \rightarrow m(A \cap B) \rightarrow \phi(A \cap B)$.\\

\noindent (iii) Show that $\phi(A)=\phi(A\cap B^c)$ and next $\phi(A)=\phi(A\cap B^c \cap C)$.\\

\noindent (iv) Show that $\phi(A)>0 \rightarrow A\cap B \cap C^c \mathcal{A}_0 \rightarrow (B + A\cap B^c \cap C) \mathcal{A}_0$.\\

\noindent (v) From $m(B + A\cap B^c \cap C)=m(B) + m(A\cap B^c \cap C)$, show that $m(A\cap B^c \cap C^c)=0$.\\

\noindent (vi) Show that $m(A\cap B^c \cap C)=0 \rightarrow \phi(A\cap B \cap C^c)=\phi(A)$.\\

\noindent Conclude on the position of $\hat{\phi}$ with respect to $\mathcal{M}_0$?
 
\noindent (d) What is the position of $\mathcal{M}$ with respect to $\mathcal{M}_0$. Conclude.\\

\noindent (e) Le $m$ be $\sigma$-finite. Show $\mathcal{M}$ is dominated by a countable linear combination of elements of $\mathcal{A}$ supported by $m$-finite subspaces and is associated to a bounded sets of coefficients.\\

\noindent \textit{Hints}. Consider  a measurable partition of $\Omega$ : $\Omega=\sum_{j\geq 0}$, 
$$
\forall j\geq 0, \ \Omega_j \in \mathcal{A} \text{ and }m(\Omega_j)<+\infty.
$$

\noindent Take induced applications $\phi_j=\phi(\Omega_j\cap B)$ and $m_j(B)=m(\Omega_j\cap B)$ for $B\in \mathcal{A}_{\Omega_j}$, $j\geq j$.\\

\noindent Apply the results of Question (a)-(b) and conclude.\\

\noindent \LARGE \textbf{DOC 08-07 : Exercise on a continuous family of measures} \label{doc08-07}\\
\bigskip
\Large

\noindent \textbf{NB}. This exercise is inspired by Lemma 2.1 and its proof in \cite{shao}, page 104. This result is important in Mathematical Statistics regarding the factorization theorem for sufficient statistics.\\

\noindent \textbf{Dominated family of $\sigma$-additives applications}. Let $\mathcal{M}$ be a family of real-valued $\sigma$-additive and non-null applications defined on $(\Omega, \mathcal{A})$ and let $m$ be a  measure defined on the same space. For a coherence purpose, let us suppose \\ 

\noindent The family $\mathcal{M}$ is dominated by $m$ if and only any element of $\mathcal{M}$ is $m$-continuous. The family $\mathcal{M}$ is dominated if and only if there is a measure on 
$(\Omega, \mathcal{A})$ which dominates it.\\ 

\noindent \textbf{NB}. In all this document, $\mathcal{M}$ is a family of real-valued $\sigma$-additive applications defined on $(\Omega, \mathcal{A})$ and $m$ a measure defined on the same space.\\

\noindent \textbf{Summary}. The objective of the exercise is to show this.\\

\noindent (A) If $m$ is finite, $\mathcal{A}$ is dominated by a countable linear combination of elements of its elements whose set of coefficient is bounded.\\

\noindent (B) if $m$ is $\sigma$-finite, $\mathcal{M}$ is dominated by a countable linear combination of elements of $\mathcal{M}$ supported by $m$-finite subspaces and is associated to a bounded set of coefficients.\\

\bigskip \noindent \textbf{Exercise 01}. \label{exercise01_sol_doc08-07}\\

\noindent Consider a family $\mathcal{M}$ of $\sigma$-finite measures defined on $(\Omega, \mathcal{A},)$ and $m$ another $\sigma$-finite measures defined on the same space.\\

\noindent Question (a) Suppose that the family is countable, that is $\mathcal{M}=\{\phi_i, \ i\in I\}$, $I \subset \mathcal{A}$. For any $c=(c_i)_{i\in I} \in (]0,+\infty[)^I$, $\mathcal{M}$ is dominated by

$$
\phi_c=\sum_{i \in I} c_i \phi. \ (FC01)
$$ 

\noindent Question (c) Let us remain in the general case and let us suppose that \textbf{$m$ is finite}. Define for any $B \in \mathcal{A}$ with $m(B)>0$ :

$$
\mathcal{D}_B=\{C \in \mathcal{A}, \ C \in B, m(C)>0\}.
$$

\noindent Define (where $\neq$ is the negation logical operator) :

$$
\mathcal{A}_0=\{B \in \mathcal{A}, \ \exists \phi \in \mathcal{M} \phi(B)>0 \text{ and for any } C \in \mathcal{D}_B, \neg (1_B \frac{d\phi}{dm}=0 \ m-\textit{a.e})\}
$$

\noindent Show that class is not empty.\\

\noindent \textit{Hint}. Consider $B \notin \mathcal{A}_0$ and show that we cannot have $\phi(B)>0$. For this, suppose $\phi(B)>0$ and deduce that $m(B)=0$. For each $C \in \mathcal{D}_B$, use the formula 
$$
\phi(C)=\int 1_C \frac{d\phi}{dm} \ dm=0,
$$

\noindent show that assertion $\neg (1_C \frac{d\phi}{dm}=0 \ $m$-\textit{a.e})$ holds.\\

\noindent Question (b) Suppose that for all $(\phi_1, \phi_2) \in \mathcal{M}^2$, there exists a linear combination  $c_1 \phi_1 + c_2 \phi_2 \in \mathcal{M}$ with $c_i\neq 0$, $i\in \{1,2\}$. Show that $\mathcal{A}_0$ is stable under finite unions.\\

\noindent \textit{Hints}. Consider $(B_1,B_2)  \in \mathcal{A}_0^2$. For each $i \in \{1,2\}$, consider $\phi_i$ and $B_{i,0}\in \mathcal{D}_{B_i}$ such that $\phi_i(B_i)>0$ and $\neg (1_{B_{i}} \frac{d\phi}{dm}=0 \ $m$-\textit{a.e})$. Show that $\phi(B)>0$. To show that $\neg (1_C \frac{d\phi}{dm}=0 \ $m$-\textit{a.e})$ for all $C \in \mathcal{D}_B$, use Point (08.07c)

$$
\frac{d\phi}{dm}=c_1 \frac{d\phi_1}{dm}+c_2 \frac{d\phi_2}{dm},
$$ 
 
\noindent and by the non-negativity of the derivatives to show that the assertion : (There exists $C \in \mathcal{D}_B$, :$(1_C \frac{d\phi}{dm}=0 \ $m$-\textit{a.e})$), leads to contradiction related to 
$C\cap B_1$ or $C\cap B_2$ and conclude.\\

\noindent Question (c) Suppose that the family $\mathcal{M}$ is bounded by a finite number $M$, that is

$$
\sup_{A \in \mathcal{A}, \ \phi \in \mathcal{M}} \phi(B) \leq M <+\infty.
$$

\noindent Define for each $I \subset \mathcal{N}$, $I\neq \emptyset$,

$$
\mathcal{C}_\mathcal{N}=\{c=(c_i)_{i\in \mathcal{N}} \in ([0,+\infty[)^\mathcal{N} \ : \|c\|_{\infty}=\sup_{i\in \mathcal{N}} c_i<+\infty \text{ and } \|c\|_{1}=\sup_{i\in \mathcal{N}} c_i>0\}.
$$

\noindent and

$$
\mathcal{M}_{0}=\{ \sum_{i\in I} c_i \phi_i, \ c\in \mathcal{C}_\mathcal{N}, \ (\phi_i)_{i\in \mathcal{N}} \in \mathcal{M}^\mathcal{N}\}.
$$

\noindent \textbf{Important remark}. Beyond the notation, we have to understand that $\mathcal{M}_{0}$ is the class of all countable combinations of elements of $\mathcal{M}$ such that for each combination, the set of the coefficients contains at least a non-zero element and is bounded in $\mathcal{R}$.\\ 

\noindent (c1) Show that $\mathcal{M}_{0}$ is still dominated by $m$.\\

\noindent (c2) Define 

$$
\mathcal{A}_0^{\ast}=\{B \in \mathcal{A}, \ \exists \phi \in \mathcal{M}_0 \phi(B)>0 \text{ and } \neg (1_B \frac{d\phi}{dm}=0 \ m-\textit{a.e}) \}.
$$

\noindent We may easily see from Question (b) (with a small extra work) that $\mathcal{A}_0^{\ast}$ is not empty. Also, it is easy to see that for $(\phi_1, \phi_2) \in \mathcal{M}_0^2$,   $\phi_1 + \phi_2 \in \mathcal{M}_0$.\\

\noindent Consider, the characterization of the supremum allows to have a sequence $(B_n)_{n\geq 0}\subset \mathcal{A}_0^{\ast}$ such that

$$
m(C_n) \rightarrow \beta=\sup_{B\in \mathcal{A}_0^{\ast}} m(B).
$$

\noindent For $n\geq 0$, consider $\phi_n\in \mathcal{M}_0$  such that

$$
\phi_n(B)>0 \text{ and } \neg (1_{B_n} \frac{d\phi_n}{dm}=0 \ m-\textit{a.e.}).
$$

\noindent Question (c2-a) Fix $(\alpha_i)_{i\in \mathcal{N}} \in \mathcal{C}_{\mathcal{N}}$, \textbf{but} such that $\alpha_i\neq 0$ for all $i\in \mathcal{N}$. Put

$$
\hat{\phi}=\sum_{n\geq 0} \alpha_n \phi_n \in \mathcal{M}_1
$$

\textit{Hints}. For each $n\geq 0$, consider $(\phi_{i,n})_{i\in \mathbb{N}} \in \mathcal{A}^{\mathbb{N}}$ and $(c_{i,n})_{i \in \mathcal{N}} \subset \mathcal{C}_{\mathcal{N}}$ such that 

$$
\phi_n=\sum_{i\geq 0}c_{i,n} \phi_{i,n}
$$

\noindent Writing

$$
\hat{\phi}=\sum_{n\geq 0, i\geq 0} \alpha_n c_{i,n} \phi_{i,n},
$$

\noindent Show that $\hat{\phi} \in \mathcal{M}_0$.\\
 
\noindent (c2-3) Let $B=\bigcup B_n$. Use this application $\hat{\phi}$ to show that $B \in \mathcal{A}_0^{\ast}$. Use the method of Question (c) in the context of a countable sum.\\

\noindent (c2-4) Suppose that for some $A\in \mathcal{A}$, $\hat{\phi}(A)=0$. Take an arbitrary element $\phi in \mathcal{M}_0$. Denote $C=(\frac{d\phi}{dm} > 0)$. Use the following facts to be justified :

\noindent (i) Recall that : $B \in \mathcal{A}_0$.\\

\noindent (ii)  Show that $\hat{\phi}(A \cap B)=0 \rightarrow m(A \cap B) \rightarrow \phi(A \cap B)$.\\

\noindent (iii) Show that $\phi(A)=\phi(A\cap B^c)$ and next $\phi(A)=\phi(A\cap B^c \cap C)$.\\

\noindent (iv) Show that $\phi(A)>0 \rightarrow A\cap B \cap C^c \mathcal{A}_0 \rightarrow (B + A\cap B^c \cap C) \mathcal{A}_0$.\\

\noindent (v) From $m(B + A\cap B^c \cap C)=m(B) + m(A\cap B^c \cap C)$, show that $m(A\cap B^c \cap C^c)=0$.\\

\noindent (vi) Show that $m(A\cap B^c \cap C)=0 \rightarrow \phi(A\cap B \cap C^c)=\phi(A)$.\\

\noindent Conclude on the position of $\hat{\phi}$ with respect to $\mathcal{M}_0$?
 
\noindent (d) What is the position of $\mathcal{M}$ with respect to $\mathcal{M}_0$. Conclude.\\

\noindent (e) Le $m$ be $\sigma$-finite. Show $\mathcal{M}$ is dominated by a countable linear combination of elements of $\mathcal{A}$ supported by $m$-finite subspaces and is associated to a bounded sets of coefficients.\\

\noindent \textit{Hints}. Consider  a measurable partition of $\Omega$ : $\Omega=\sum_{j\geq 0}$, 
$$
\forall j\geq 0, \ \Omega_j \in \mathcal{A} \text{ and }m(\Omega_j)<+\infty.
$$

\noindent Take induced applications $\phi_j=\phi(\Omega_j\cap B)$ and $m_j(B)=m(\Omega_j\cap B)$ for $B\in \mathcal{A}_{\Omega_j}$, $j\geq j$.\\

\noindent Apply the results of Question (a)-(b) and conclude.\\

\bigskip \noindent \textbf{SOLUTIONS}.

\noindent Question (a) Let $c=(c_i)_{i\in I} \in (]0,+\infty[)^I$ and consider $\phi_c$ as given by Formula (FC01). For any $A\in \mathcal{A}$, it is clear that $\phi_c(A)=0$ implies that for each 
$i\in I$, $\phi_i(A)=0$, which establishes that $\phi_c$ dominates $\mathcal{A}$.\\

\noindent Question (b) Let $B$ any elements of $\in \mathcal{A}$. Suppose that $B \notin \mathcal{A}_0$. Let $\phi$ be any element of $\mathcal{M}$. We should have either $\phi(B)=0$ or $\phi(B)>0$. But if $\phi(B)>0$, wa have $m(B)>0$, otherwise $m(A)=0$ and the Radon-Nikodym property would imply

$$
\phi(A)=\int_A \frac{d\phi}{dm} \ dm=0.
$$

\noindent Hence $B \in \mathcal{D}_B$ and $\mathcal{D}_B$ is not empty. Now if, in addition, we have $B\notin \mathcal{A}_0$, we would get that assertion $(1_C \frac{d\phi}{dm}=0 \ m-\textit{a.e})$ for some $C \in \mathcal{D}_B$, which would implies
$$
\phi(C)=\int 1_C \frac{d\phi}{dm} \ dm=0=0,
$$

\noindent which in turn would imply, by the Radon-Nikodym Formula, $m(C)=0$. \noindent This contradiction ensures that $\phi(B)=0$. Hence if $\mathcal{A}_0$ is empty, we get that  $\phi(B)=0$ for all $B\in \mathcal{A}$ and for all $\phi \in \mathcal{M}$. This is in contradiction with our assumption that that all the elements of $\mathcal{M}$ are non-zero applications.\\

\noindent Question (b) Let $(B_1,B_2)  \in \mathcal{A}_0^2$. For each $i \in \{1,2\}$, consider $\phi_i$ such that $\phi_i(B_i)>0$ and $\neg (1_{C} \frac{d\phi_i}{dm}=0 \ $m$-\textit{a.e})$ for all 
$C \in \mathcal{D}_{B_i}$. Let us form $\phi=c_1 \phi_1+ c_2 \phi_1 \in \mathcal{M}$ and $B_1 \cup B_2$. We have
 
$$
\phi(B)=c_1 \phi_1(B)+ c_2 \phi_1(B)\geq c1\phi_1(B_1)+c_2\phi_2(B_2)>0.
$$

\bigskip \noindent By Point (08.07c), we have

$$
\frac{d\phi}{dm}=c_1 \frac{d\phi_1}{dm}+c_2 \frac{d\phi_2}{dm}. \ (RD1)
$$ 
 
\noindent Now for any $C\in \mathcal{D}_B$, $m(C\cap B_1)>0$ or $m(C\cap B_2)>0$, otherwise we would have $m(C)=m((C\cap B_1) \cup (C\cap B_2))=0$. Thanks to the non-negativity of the terms of the addition in Formula (RD1),  the assertion : (There exists $C\in \mathcal{D}_B$ : $(1_B \frac{d\phi}{dm}=0 \ $m$-\textit{a.e})$) would imply that, for each $i\in \{1,2\}$ :
$(1_{C\cup B_i} \frac{d\phi_i}{dm}=0 \ $m$-\textit{a.e})$ whenever $m(C\cup B_i)>0$. Since of the condition $m(C\cup B_i)>0$, $i\in \{1,2\}$, holds, we reach a contradiction. Hence we have 
$\neg (1_C \frac{d\phi}{dm}=0 \ m-\textit{a.e})$ for all $C \in \mathcal{D}_B$ and hence $B \in \mathcal{A}_0$.\\

\noindent Question (c) Define for each $I \subset \mathcal{N}$, $I\neq \emptyset$,

$$
\mathcal{C}_\mathcal{N}=\{c=(c_i)_{i\in \mathcal{N}} \in ([0,+\infty[)^\mathcal{N} \ : \|c\|_{\infty}=\sup_{i\in \mathcal{N}} c_i<+\infty \text{ and } \|c\|_{1}=\sup_{i\in \mathcal{N}} c_i>0\}.
$$

\noindent and

$$
\mathcal{M}_{0}=\{ \sum_{i\in I} c_i \phi_i, \ c\in \mathcal{C}_\mathcal{N}, \ (\phi_i)_{i\in \mathcal{N}} \in \mathcal{M}^\mathcal{N}\}.
$$

\noindent \textbf{Important remark}. Beyond the notation, we have to understand that $\mathcal{M}_{0}$ is the class of all countable combinations of elements of $\mathcal{M}$ such that for each combination, the set of the coefficients contains at least a non-zero element and is bounded in $\mathcal{R}$.\\ 

\noindent (c1) It is obvious that $\mathcal{M}_{0}$ is still dominated by $m$.\\

\noindent (c2) Define 

$$
\mathcal{A}_0^{\ast}=\{B \in \mathcal{A}, \ \exists \phi \in \mathcal{M}_0 \phi(B)>0 \text{ and } \neg (1_B \frac{d\phi}{dm}=0 \ m-\textit{a.e}) \}.
$$

\noindent We may easily see from Question (b) (with a small extra work) that $\mathcal{A}_0^{\ast}$ is not empty. Also, it is easy to see that for $(\phi_1, \phi_2) \in \mathcal{M}_0^2$,   $\phi_1 + \phi_2 \in \mathcal{M}_0$.\\

\noindent Consider, the characterization of the supremum allows to have a sequence $(B_n)_{n\geq 0}\subset \mathcal{A}_0^{\ast}$ such that

$$
m(C_n) \rightarrow \beta=\sup_{B\in \mathcal{A}_0^{\ast}} m(B).
$$

\noindent For $n\geq 0$, consider $\phi_n\in \mathcal{M}_0$  such that

$$
\phi_n(B)>0 \text{ and } \neg (1_{B_n} \frac{d\phi_n}{dm}=0 \ m-\textit{a.e.}).
$$

\noindent (c2-a) Let $(\alpha_i)_{i\in \mathcal{N}}$ be any element of $\mathcal{C}_{\mathcal{N}}$ such that $\alpha_i\neq 0$ for all $i\in \mathcal{N}$. Define

$$
\hat{\phi}=\sum_{n\geq 0} \alpha_n \phi_n \in \mathcal{M}_1.
$$

\bigskip \noindent By assumption there exist, for any $n\geq 0$, $(\phi_{i,n})_{i\in \mathbb{N}} \in \mathcal{A}^{\mathbb{N}}$ and $(c_{i,n})_{i \in \mathcal{N}} \subset \mathcal{C}_{\mathcal{N}}$ such that 

$$
\phi_n=\sum_{i\geq 0}c_{i,n} \phi_{i,n}.
$$

\bigskip \noindent we get

$$
\hat{\phi}=\sum_{n\geq 0, i\geq 0} \alpha_n c_{i,n} \phi_{i,n}=\sum_{n\geq 0, i\geq 0} d_{i,n} \phi_{i,n}.
$$

\noindent \noindent (c2-3) Let $B=\bigcup B_n$. We have, since all the $\alpha_n$'s are non-zero numbers,

$$
\hat{\phi}(B)=\phi=\sum_{n\geq 0} \alpha_n \phi_n(B) \geq \sum_{n\geq 0} \alpha_n \phi_n(B_n)>0. 
$$

\bigskip \noindent Next, by Point \textit{(08.07c)} in Doc 08-01, we have

$$
\frac{d\hat{\phi}}{dm}=\phi=\sum_{n\geq 0} \alpha_n \frac{d\phi_n}{dm}. \ (RD2) 
$$

\noindent \noindent Now for any $C\in \mathcal{D}_B$, there exists $N\geq 0$ such that $m(C\cap B_N)>0$, otherwise we would have $m(C)=m(\bigcup C\cap B_n)=0$. Thanks to the non-negativity of the terms of the sum in Formula (RD2),  the assertion : (There exists $C\in \mathcal{D}_B$ : $(1_C \frac{d\phi}{dm}=0 \ $m$-\textit{a.e})$) would imply that, $(1_{C\cup B_N} \frac{d\phi_i}{dm}=0 \ $m$-\textit{a.e})$. This contradiction ensures that we have $\neg (1_C \frac{d\phi}{dm}=0 \ $m$-\textit{a.e})$ for all $C \in \mathcal{D}_B$ and hence $B \in \mathcal{A}_0$. This leads to :

$$
m(B) \leq \beta.
$$

\noindent (c2-3) Let us show that the family is dominated by $\hat{\phi}$. Suppose that for some $A\in \mathcal{A}$, $\hat{\phi}(A)=0$. Let $\phi$ be ab arbitrary element of $\mathcal{M}_0$. Denote $C=(\frac{d\phi}{dm} > 0)$.\\

\noindent (i) We already have that $B \in \mathcal{A}_0$.\\

\noindent (ii)  Since $\hat{\phi}(A \cap B)\leq \hat{\phi}(A)=0$, we have

$$
\hat{\phi}(A \cap B)=\int 1_{A\cap B} \frac{d\hat{\phi}}{dm} \ dm =0.
$$

\noindent Hence $(1_{A\cap B} \frac{d\hat{\phi}}{dm}, \ m-\textit{a.e.})$. This would be in contradiction with the assumption if $m(A\cap B)>0$. Hence $0=m(A\cap B)$ which implies $\phi(A\cap B)=0$ since 
the class $\mathcal{M}_0$ is dominated by $m$.

\noindent (iii) Hence $\phi(A)=m(A\cap B)+\phi(A\cap B^c)=\phi(A\cap B^c)$ and next $\phi(A)=\phi(A\cap B^c \cap C)+\phi(A\cap B^c \cap C^c)$. But

$$
\phi(A\cap B^c \cap C^c)=\int_{A\cap B^c \cap C^c} \frac{d\phi}{dm} \ dm=\int_{A\cap B^c \cap C^c} 0 \ dm=0.
$$

\noindent Hence $\phi(A)=\phi(A\cap B^c \cap C)$.\\

\noindent (iv) Let $\phi(A)>0$. Let us prove that $B_0=\phi(A\cap B^c \cap C \in \mathcal{A}_0^{\ast}$. We first have $\phi(A)=\phi(B_0)$. Next, for all $D\in \mathcal{D}_{B_0}$, we have 
$m(A)>0$ and $\frac{d\phi}{dm} >0$ and hence $\neq (1_D \frac{d\phi}{dm}=0)$ holds.\\

\noindent (v) By Question (b), $B + (A\cap B^c \cap C) \in \mathcal{A}_0^{\ast}$ and

$$
m(B + (A\cap B^c \cap C))=m(B) + m(A\cap B^c \cap C)=\beta + m(A\cap B^c \cap C).
$$

\noindent (vi) By definition of the supremum, we necessarily have $m(A\cap B^c \cap C)=0$, which by continuity, leads to $0=\phi(A\cap B \cap C)=\phi(A)$. This contradiction ensures that we necessarily have $\phi(A)=0$.\\

\noindent We conclude that $\mathcal{M}_0$ is dominated by $\hat{\phi}$.\\
 
\noindent (d) Since $\mathcal{M}\subset \mathcal{M}_0$, $\mathcal{M}$ is dominated by $\hat{\phi}$. $\square$ \\

\noindent (e) Here $m$ is $\sigma$-finite. Let us consider a measurable partition of $\Omega$ : $\Omega=\sum_{j\geq 0}$, 
$$
\forall j\geq 0, \ \Omega_j \in \mathcal{A} \text{ and }m(\Omega_j)<+\infty.
$$

\noindent By taking the induced applications $\phi_j=\phi(\Omega_j\cap B)$ and $m_j(B)=m(\Omega_j\cap B)$ for $B\in \mathcal{A}_{\Omega_j}$, $j\geq j$, we still have that each family

$$
\mathcal{M}_j=\{\mathcal{A}_{\Omega_j} \ni B \mapsto \phi(\Omega_j\cap B) \},
$$

\noindent is a family measure dominated by the finite measure $m_j$. Hence by Questions (c)-(d), for each $j\geq 0$ there exists a countable linear combination

$$
\phi_j=\sum_{i\geq 0} c_{i,j} \phi_{i,j},
$$

\noindent where (for each $j\geq 0$,  $\phi_{i,j}\in \mathcal{M}_j$, $c_{i,j}\geq 0$ , $\sup_{i\geq 0} c_{i,j}<+\infty$) which dominates $\mathcal{M}_j$.\\

\noindent The choice of the $c_{i,j}$'s depends on us. We may choose them all below some $\gamma>0$. Define for $A\in \mathcal{A}$, by the Fubini's formula for non-negative series,

\begin{eqnarray*}
\hat{\phi}(A)&=&\sum_{j\geq 0} \phi_j(A)\\
&=&\sum_{i\geq 0} \sum_{j\geq 0} \phi_{i,j}(A).
\end{eqnarray*}

\noindent We readily get $\hat{\phi}$ is a dominating measure of $\mathcal{M}$ and $\hat{\phi}$ is a countable linear combination of elements of $\mathcal{M}$ supported by $m$-finite subspaces.

%\part{Part IV \\ Transitions to Lebesgue-Stieljes Probability Measures and to Functional Analysis}
%\newpage
%\chapter{Transition to Lebesgue-Stieljes Probability Measures}
%\include{09_transition1/mes_doc_09_01_ap}
%\include{09_transition1/mes_doc_09_02_ap}
%\include{09_transition1/mes_doc_09_03_ap}

%\newpage
%\chapter{Transition to Functional Analysis}

\part{Towards Probability Theory and Functional Analysis}
\chapter{Introduction to spaces $L^p$, $p\in[1,+\infty]$ and Daniel Integral} \label{09_lp}

\noindent \textbf{Content of the Chapter}

\begin{table}[htbp]
	\centering
		\begin{tabular}{llll}
		\hline
		Type& Name & Title  & page\\
		S& Doc 09-01 & Spaces $L^p$ - A summary & \pageref{doc09-01}\\
		D& Doc 09-02 & $L^p$ spaces by Exercises   & \pageref{doc09-02} \\
		SD& Doc 09-03  & $L^p$ spaces by Exercises with Exercises & \pageref{doc09-03} \\
		D& Doc 09-04  & Daniel Integral - Exercises & \pageref{doc09-04} \\
		SD& Doc 09-05  & Daniel Integral - Exercises with Solution& \pageref{doc09-05}\\
		\hline
		\end{tabular}
\end{table}

\newpage
\noindent \LARGE \textbf{DOC 09-01 : $L^p$ Spaces - A summary} \label{doc09-01}\\
\bigskip
\Large

\noindent \textbf{(09-01) Space $\mathcal{L}^{p}(\Omega, \mathcal{A},m)$, $p\in [1,+\infty[$}.\\

\noindent $\mathcal{L}^{p}(\Omega, \mathcal{A},m)$ is the space of all real-valued and measurable functions $f$ defined on $\Omega$ such that $|f|^p$ is integrable with respect to $m$, that is

$$
\mathcal{L}^{p}(\Omega, \mathcal{A},m)=\{f \in \mathcal{L}_0(\Omega, \mathcal{A}) \int |f|^p \ dm  < +\infty \}.
$$

\bigskip \noindent \textit{Recall}. $\mathcal{L}_0(\Omega, \mathcal{A})$ is the class of all real-valued and measurable functions defined on $\Omega$ as in \textit{Point (4.19a) in Doc 04-01 of Chapter \ref{04_measures}}, page \pageref{04_equivcal}.\\

\bigskip \noindent \textbf{(09-02) Space of $m$-a.e. bounded functions $\mathcal{L}^{\infty}(\Omega, \mathcal{A},m)$}.\\

\noindent $\mathcal{L}^{\infty}(\Omega, \mathcal{A},m)$ is the class of all real-valued and measurable functions which are bounded $m$-a.e., i.e,

$$
\mathcal{L}^{\infty}(\Omega, \mathcal{A},m)=\{f \in \mathcal{L}_0(\Omega, \mathcal{A}) 0, \ \exists C\in \mathbb{R}_{+}, \ |f|\leq C \ m-a.e \}.
$$

\bigskip \noindent \textbf{(05-03) Equivalence classes}.\\

\noindent  When we consider only the equivalence classes pertaining to the equivalence relation introduced in \textit{Point 4.19b in Doc 04-01 in Chapter \ref{04_measures}}, page \pageref{doc04-01}, we get the quotient spaces

$$
L^{p}(\Omega, \mathcal{A},m)=\mathcal{L}^{p}(\Omega, \mathcal{A},m)/ \mathcal{R}.
$$

\bigskip \noindent that is

$$
L^{p}(\Omega, \mathcal{A},m)=\{ \dot{f}, f \in  \mathcal{L}_0(\Omega, \mathcal{A}), \ \int |f|^p \ dm \text{ } < +\infty \}
$$

\bigskip \noindent and

$$
L^{\infty}(\Omega, \mathcal{A},m)=\mathcal{L}^{+\infty}(\Omega, \mathcal{A},m)/ \mathcal{R},
$$

\bigskip \noindent that is 

$$
L^{\infty}(\Omega, \mathcal{A},m)=\{ \dot{f}, f \mathcal{L}_0(\Omega, \mathcal{A}),\  |f| \ m-bounded \}.
$$

\bigskip \noindent \textbf{(09-04) Norms on $L^{p}$}.\\

\noindent The space $L^{p}(\Omega, \mathcal{A},m)$, $1\leq p <+\infty$, is equipped with the norm

$$
||f||_{p} = ( \int |f|^p dm )^{1/p}.
$$

\bigskip \noindent The space $L^{+\infty}(\Omega, \mathcal{A},m)$ is equipped with the norm

$$
||f||_{\infty} = \inf \{C \in \mathbb{R}_{+}, \ |f|\leq C \ m-.a.e  \}.
$$

\bigskip \noindent \textbf{(09-05) Remarkable inequalities on $L^{p}$}.\\

\noindent \textbf{09-05 Space $\mathcal{L}^{\infty}(\Omega, \mathcal{A},m)$}.\\

\noindent Define $\mathcal{L}^{^\infty}$ as the class of the equivalence classes of measurable functions $m$-\textit{a.e.} bounded. An other to clearly define it the following. Define

$$
\mathcal{L}^{\infty}=\{f \in \mathcal{L}_0, \ \exists M\in [0, +\infty[, \ |f|\leq M, \ m-a.e.\}.
$$

\bigskip \noindent Next, define  $L^{\infty}=\mathcal{L}/\mathcal{R}$ as the quotient space of by the equivalence relation $\mathcal{R}$ of the $m$-\textit{a.e.} equality. For $f \in L^{\infty}$, define

$$
\|f\|_{\infty}= \inf \{M\in [0, +\infty[, \ |f|\leq M, \ m-a.e.\}.
$$ 

\bigskip \noindent $\|.\|_{\infty}$ is a norm and $(L^{\infty},+, ., \|.\|_{\infty})$ is a real Banach space.\\

\bigskip \noindent \textbf{(09-06) Banach spaces  $L^{p}$, $p\in [1,+\infty]$}.\\

\noindent  For all $p \in [1, +\infty]$, $L^{p}(\Omega,\mathcal{A}, m) (+, ., ||f||_{p})$ is a Banach space.\\

\noindent $L^{2}(\Omega,\mathcal{A}, m) (+, ., ||f||_{2})$ is a Hilbert space.\\

\bigskip \noindent \textbf{(09-07) Convergence in $L^{p}$}.\\

\noindent \textbf({09-07a)} Let $(f_n)_{n\geq 0} \subset L^p$, $p\in [1,+\infty[$, converge in $L^p$ to $f$, denoted by 

$$
f_n \overset{L^p}{\longrightarrow} f.
$$

\bigskip \noindent We have the following implications. \\

\noindent (a) $(f_n)_{n\geq 0}$ is Cauchy in measure.\\

\noindent (b) There exists a sub-sequence $(f_{n_k})_{k\geq 0}$ of $(f_n)_{n\geq 0}$ converging $m$-\textit{a.e.} and in measure to $f$.\\
 
\noindent (c) $(f_n)_{n\geq 0}$ converges in measure to $f$.\\

\noindent (d) $\|f_n\|_p$ converges to $\|f\|_p$.\\

\noindent \textbf{(09-07b)} Let $(f_n)_{n\geq 0} \subset L^{\infty}$ converge in $L^{\infty}$ to $f$. Then $f_n \rightarrow f $ $m$-\textit{a.e.} 

\noindent \textbf{NB}. It is possible to reverse the implications, especially for $1\leq p <+\infty$, at the cost of extra conditions, among with equi-integrability ones, for finite measures. This is why we reserve that part to the probability theory monographs.

\noindent \LARGE \textbf{DOC 09-02 : $L^p$ Spaces by Exercises} \label{doc09-02}\\
\bigskip
\Large

%\bigskip \noindent \textbf{Exercise 1}. \label{exercise01_doc09-02}
%\bigskip \noindent \textbf{SOLUTIONS}.\\

\noindent \textbf{NB 1}. Throughout this document $(\Omega, \mathcal{A}, m)$ is a measure space. Functions $f$, $g$, $h$, $f_n$, $g_n$, $h_n$, $n\geq 0$, etc., are classes of measurable functions defined on
$(\Omega, \mathcal{A})$ with values in $\overline{\mathcal{A}}$.\\
 
\bigskip \noindent \textbf{Exercise 1}.  (Useful inequalities on real numbers)\label{exercise01_doc09-02}\\

\noindent Question (a) Let $a$ and $b$ be finite real numbers and let $p\geq 1$ be finite real number. Show that 

$$
|a+b|^p \leq 2^{p-1} (|a|^p+|b|^p). \ (I1)
$$

\bigskip \noindent \textit{Hints}. For $p\geq 1$, use the convexity of the function $]0,+\infty[ \ni x \mapsto x^p$.\\

\noindent Question (b) Two positive real numbers $p$ and $q$ are conjugate if and only if $1/p + 1/q=1$. Show that $p$ and $q$ are conjugate, then both of them are greater that the unity and that
$p=q/(q-1)$ and $q=p/(p-1)$.

\noindent Question (c) Let $a$ and $b$ be finite real numbers and let $p> 1$ and $p>1$ be $p\geq 1$ be conjugate real numbers. Show that 
Show that for all finite real numbers $a$ and $b$ and for all  , that

$$
|ab| \leq \frac{|a|^p}{p} + \frac{|b|^q}{q}. \ (I2)
$$

\bigskip \noindent \textit{Hints}. \\

\noindent (a) Consider $a\geq 0$ and $b\geq 0$. Since Formula $(I2)$ is obvious for $ab=0$, continue with $a> 0$ and $b>0$. Use the development of $(a-b)^2$ to solve the case $p=q$.\\

\noindent (b) Pass to the case where $p<q$ or $q<p$ and, because of the symmetry of the roles, consider $1<q<p$. Show that : $a^{p-1}>b$ and $b^{p-1}>a$ leads to $a<a^{(q-1)/(p-1)}$. Conclude that we have 
$b^{p-1}\leq a$ or $a^{p-1}\leq b$. In both cases we use the same method by making $a$ and $b$ playing symmetrical roles.\\

\noindent Suppose that $a^{p-1}\leq b$. The curve of the non-decreasing function $y=f(x)=x^{p-1}$ on $[0,a]$ join the points $A=(0,0)$ and $D=(a,a^{p-1}$. So it intersects with the line $y=b$ at a point $C=(c,b)$ with $c \in [0,a]$. Denote $E=(a,0)$ and $B=(0,b)$. Let $\Gamma_{1}=\{(x,x^{a-1}, \ 0\leq x\leq x\}$ and $\Gamma_{2}=\{(x,x^{a-1}, \ 0\leq x\leq c\}$ be the curve of $f$ respectively from $0$ and 
$a$ and from $0$ to $c$ respectively. With a simple drawing, we may see that the addition of the area $S_1$ between the segments $A-E$, $ED$ and the curve $\Gamma_{1}$, and the area $S_2$ between the segments $AB$, $BC$ and the curve $\Gamma_{2}$ exceed the area of the rectangle $A,B,(a,b),E$ which is $ab$. Show that $S_1+S_2$ is exactly the right-hand member of Formula (I2). Conclude.\\

\bigskip \noindent \textbf{Exercise 2}. (H\H{o}lder's Inequality) \label{exercise02_doc09-02}\\

\noindent Let us denote for any $p\geq 1$, for any $f \in L^p$,

$$
N_p(f)=\left(\int |f|^p\right)^{1/p}.
$$

\bigskip \noindent Question (a) Let $p>1$ and $q>1$ be two conjugate real numbers. Let $f \in L^p$ and $f \in L^q$. Show that $fg \in L^1$ and

$$
\left\vert \int fg \ dm \right\vert \leq N_p(f) N_q(g). (H) 
$$

\bigskip \noindent \textit{Hints}.\\

\noindent (1) Remark that $f$ and $g$ are $m$-\textit{a.e.} finite by assumption. Tell why Formula (H) holds if $N_p(f)=0$ or $N_q(g)=0$.\\

\noindent (2) Suppose that none of $N_p(f)$ or $N_q(g)$ is zero. Apply Question (c) of Exercise 1 to

$$
a=\frac{f}{N_p(f)} \text{    and   } b=\frac{g}{N_q(g)}
$$

\bigskip \noindent to conclude.\\

\noindent Question (b) Extend Formula (H) in the following way : For any measurable application such that the product $N_p(f)N_q(g)$ makes senses regardless of any integrability condition on $f$ or on $g$,

$$
 \int \left\vert fg \right\vert \ dm  \leq N_p(f) N_q(g). (HE) 
$$

\bigskip \noindent \textit{Hints} Just discuss on the finiteness or not of $N_p(f)$ and $N_q(g)$.\\

\bigskip \noindent \textbf{Exercise 3}. (Minkowski's Inequality) \label{exercise03_doc09-02}\\

\noindent Let $p\geq 1$ and let that $(f,g) \in (L^p)^2$.\\

\noindent Question (a) Show that $f+g$ is $m$-\textit{a.e.} defined and $f+g \in L^p$ and show that $f+g\in L^p$.\\

\noindent \textit{Hint}. Use Question (a) of Exercise 1.\\

\noindent Question (b) Show that

$$
N_p(f+g) \leq N_p(f)+ N_p(g).
$$

\bigskip \noindent \textit{Hints}.\\

\noindent Give a quick solution for $p=1$.\\

\noindent What happens if $N_p(f)=0$ or $N_p(g)=0$?\\

\noindent Remark that there is nothing to for  $p>1$ but  $N_p(f+g)=0$.\\

\noindent Now let $N_p(f)\neq 0$, $N_p(g)\neq 0$, $N_p(f+g)\neq 0$ and $p>1$. Use the following formula to be justify 

\begin{equation*}
\left\vert f+g\right\vert ^{p}=\left\vert f+g\right\vert ^{p-1}\left\vert
f+g\right\vert \leq \left\vert f+g\right\vert ^{p-1}\left\vert f\right\vert
+\left\vert f+g\right\vert ^{p-1}\left\vert g\right\vert 
\end{equation*}

\noindent to arrive at 

\begin{equation*}
\left\vert f+g\right\vert ^{p}\leq \left\vert f+g\right\vert
^{p-1}\left\vert f\right\vert +\left\vert f+g\right\vert ^{p-1}\left\vert
g\right\vert. \ (PM)
\end{equation*}
 
\noindent Apply the inequality of Exercise 2 with $q=p/(p-1)$, which will be shown to be conjugate of $p$, to both terms in the right-hand member of Inequality (PM) and try to conclude.\\

\bigskip \noindent \textbf{Exercise 4}. ($L^p$ is a Banach Space, $1\geq p <+\infty$) \label{exercise04_sol_doc09-02}\\

\noindent Let $p\geq 1$ be finite real number.\\

\noindent Question (a) Show that the application

$$
L^p \ni f \mapsto N_p(f)
$$

\bigskip \noindent is norm. From now one, put $N_p(f)=\|f\|_p$.\\

\noindent \textit{Hint}. Use the results of Exercises 2 and 3.\\

\noindent Question (b) Show (describe why?) that $(L^p, +, .)$ is a vector space, where the last operation is the external multiplication of function by real scalars.\\

\noindent Question (c) Show that when endowed with the norm $\|f\|_p$, $(L^p, +, .)$ is complete normed space, that is a Banach space.\\

\noindent \textit{Hints}. Consider a Cauchy sequence $(f_n)_{n\geq 0}$ in $(L^p, +, ., \|f\|_p)$.\\

\noindent (i) Use the Markov inequality (See Exercise 8 in Doc 05-04, page \pageref{exercise08_sol_doc05-05}, in Chapter \ref{05_integration})  and show that the sequence in Cauchy in measure.\\

\noindent (ii) By Exercise 6 in Doc 06-08, (see page \pageref{exercise06_sol_doc06-08}) in Chapter \ref{06_convergence}, consider a sub-sequence $(f_{n_k})_{k\geq 0}$ of $(f_n)_{n\geq 0}$ converging both in measure and $m$-\textit{a.e.} to a measurable function $f$.\\

\noindent Apply the (below version) Fatou-Lebesgue's lemma  (See Exercise 2 in Doc 06-09, page \pageref{exercise02_sol_doc06-09} in Chapter \ref{06_convergence}) to 
$$
\int |f_n-f_{n_k}|^p \ dm
$$ 

\bigskip \noindent as $k \rightarrow +\infty$ and use the \textit{Cauchyness} of the sequence to conclude.\\

\bigskip \noindent \textbf{Exercise 5}. \label{exercise05_doc09-02}\\

\noindent Define $\mathcal{L}^{^\infty}$ as the class of the equivalence classes of measurable functions $m$-\textit{a.e.} bounded. An other to clearly define it the following. Define

$$
\mathcal{L}^{\infty}=\{f \in \mathcal{L}_0, \ \exists M\in [0, +\infty[, \ |f|\leq M, \ m-a.e.\}.
$$

\bigskip \noindent Next, define  $L^{\infty}=\mathcal{L}/\mathcal{R}$ as the quotient space of by the equivalence relation $\mathcal{R}$ of the $m$-\textit{a.e.} equality. For $f \in L^{\infty}$, define

$$
N_{\infty}(f)= \inf \{M\in [0, +\infty[, \ |f|\leq M, \ m-a.e.\}.
$$ 

\bigskip \noindent Question (a) Show that an $m$-\textit{a.e.} bound $M$ of $f$ is also an $m$-\textit{a.e.} of any elements of the class $\overline{f}$. Apply this result to establish that $N_{\infty}$ is well-defined on $L^{\infty}$.\\
 
\noindent Question (b) Show that $(L^{\infty},+, .)$ is a linear space.\\

\noindent Question (c) Show that $N_{\infty}$ is a norm on $L^{\infty}$. From now on, we denote it $\|.\|_{\infty}$.\\

\noindent Question (d) Show that $(L^{\infty},+, ., \|.\|_{\infty})$ is a real Banach space.\\

\bigskip \noindent \textbf{Exercise 6}. ($L^2$ is a Hilbert Space) \label{exercise06_sol_doc09-02}\\

\noindent Question (a) Express the H\H{o}lder inequality for $p=q=2$, called the Cauchy-Schwarz inequality. Deduce that for $(f,g)\in (L^2)^2$, $fg \in L^1$.\\

\noindent Question (b) Define the application 

$$
(L^2)^2 \ni (f,g) \mapsto Q(f,g)=\int fg \ dm.
$$

\bigskip \noindent Show that $\Phi$ is a symmetrical bi-linear application and that $Q(f,f)=\|f\|_2$ is a scalar product. Conclude that $L^2$ is a Hilbert space.\\

\bigskip \noindent \textbf{Exercise 7}. (Convergence) \label{exercise07_sol_doc09-02}\\

\noindent Question (1) Suppose that a sequence $(f_n)_{n\geq 0} \subset L^p$, $p\in [1,+\infty[$, converges in $L^p$ to $f$, denoted by 

$$
f_n \overset{L^p}{\longrightarrow} f.
$$

\bigskip \noindent Establish the following implications. \\

\noindent (a) $(f_n)_{n\geq 0}$ is Cauchy in measure.\\

\noindent (b) There exists a sub-sequence $(f_{n_k})_{k\geq 0}$ of $(f_n)_{n\geq 0}$ converging $m$-\textit{a.e.} and in measure to $f$.\\
 
\noindent (c) $(f_n)_{n\geq 0}$ converges in measure to $f$.\\

\noindent (d) $\|f_n\|_p$ converges to $\|f\|_p$.\\

\noindent \textit{Hints}. For Parts (a), (b) and (c), just recall the results above. For Part (d), use the second triangle inequality for any norm : $\left\vert \|x\|-\|y\| \right\vert \leq \|x-y\|$.\\

\noindent Question 2) Suppose that a sequence $(f_n)_{n\geq 0} \subset L^{\infty}$, converges in $L^{\infty}$ to $f$. Show that $f_n \rightarrow f $ $m$-\textit{a.e.} 

\noindent \textit{Hint}. Just recall a result of Exercise 4.\

\noindent \LARGE \textbf{DOC 09-03 : $L^p$ Spaces by Exercises with Solutions} \label{doc09-03}\\
\bigskip
\Large

%\bigskip \noindent \textbf{Exercise 1}. \label{exercise01_sol_doc09-03}
%\bigskip \noindent \textbf{SOLUTIONS}.\\

\noindent \textbf{NB 1}. Throughout this document $(\Omega, \mathcal{A}, m)$ is a measure space. Functions $f$, $g$, $h$, $f_n$, $g_n$, $h_n$, $n\geq 0$, etc., are classes of measurable functions defined on
$(\Omega, \mathcal{A})$ with values in $\overline{\mathcal{A}}$.\\

\bigskip \noindent \textbf{Exercise 1}.  (Useful inequalities on real numbers)\label{exercise01_doc09-02}\\

\noindent Question (a) Let $a$ and $b$ be finite real numbers and let $p\geq 1$ be finite real number. Show that 

$$
|a+b|^p \leq 2^{p-1} (|a|^p+|b|^p). \ (I1)
$$

\bigskip \noindent \textit{Hints}. For $p\geq 1$, use the convexity of the function $]0,+\infty[ \ni x \mapsto x^p$.\\

\noindent Question (b) Two positive real numbers $p$ and $q$ are conjugate if and only if $1/p + 1/q=1$. Let $p$ and $q$ be conjugate, show that both of them are greater that the unity and that
$p=q/(q-1)$ and $q=p/(p-1)$. Besides, $p$ and $q$ are different unless $p=q=2$.\\

\noindent Question (c) Let $a$ and $b$ be finite real numbers and let $p> 1$ and $p>1$ be $p\geq 1$ be conjugate real numbers. Show that 
Show that for all finite real numbers $a$ and $b$ and for all  , that

$$
|ab| \leq \frac{|a|^p}{p} + \frac{|b|^q}{q}. \ (I2)
$$

\bigskip \noindent \textit{Hints}. \\

\noindent (a) Consider $a\geq 0$ and $b\geq 0$. Since Formula $(I2)$ is obvious for $ab=0$, continue with $a> 0$ and $b>0$. Use the development of $(a-b)^2$ to solve the case $p=q$.\\

\noindent (b) Pass to the case where $p<q$ or $q<p$ and, because of the symmetry of the roles, consider $1<q<p$. Show that : $a^{p-1}>b$ and $b^{p-1}>a$ leads to $a<a^{(q-1)/(p-1)}$. Conclude that we have 
$b^{p-1}\leq a$ or $a^{p-1}\leq b$. In both cases we use the same method by making $a$ and $b$ playing symmetrical roles.\\

\noindent Suppose that $a^{p-1}\leq b$. The curve of the non-decreasing function $y=f(x)=x^{p-1}$ on $[0,a]$ join the points $A=(0,0)$ and $D=(a,a^{p-1}$. So it intersects with the line $y=b$ at a point $C=(c,b)$ with $c \in [0,a]$. Denote $E=(a,0)$ and $B=(0,b)$. Let $\Gamma_{1}=\{(x,x^{a-1}, \ 0\leq x\leq x\}$ and $\Gamma_{2}=\{(x,x^{a-1}, \ 0\leq x\leq c\}$ be the curve of $f$ respectively from $0$ and 
$a$ and from $0$ to $c$ respectively. With a simple drawing, we may see that the addition of the area $S_1$ between the segments $A-E$, $ED$ and the curve $\Gamma_{1}$, and the area $S_2$ between the segments $AB$, $BC$ and the curve $\Gamma_{2}$ exceed the area of the rectangle $A,B,(a,b),E$ which is $ab$. Show that $S_1+S_2$ is exactly the right-hand member of Formula (I2). Conclude.\\

\bigskip \noindent \textbf{SOLUTIONS}.\\

\noindent Question (a). If $p=1$, Formula (I1) is the triangle inequality. Next, the formula is obvious if one of $|a|$ and $|b|$ is zero. For $p>1$, the second derivative of $f(x)=x^p$ is $^p(p-1)x^{p-2}$, which is positive on $]0,+\infty[$. Hence, for any $a>0$ and $b>0$,

$$
f\left(\frac{a+b}{2}\right)=\frac{(a+b)^p}{2^p} \leq \frac{1}{2}f(a)+\frac{1}{2}f(b)=\frac{1}{2}(a^p+b^p),
$$

\bigskip \noindent which leads to Formula (I1).\\

\noindent Question (b). This obvious question was asked to draw your attention to these simple facts.\\

\noindent Question (c). Formula $a=0$ or $b=0$, the formula is obviously true. For $p=2$ (and hence $q=2$), for $a\leq 0$ and $b\geq 0$, we have

$$
(a-b)^2=a^2+b^2-2ab\leq 0 \Rightarrow 0\leq ab \leq \frac{a^2}{2} + \frac{b^2}{2}.
$$

\bigskip \noindent For the rest of the solution, we suppose that $q<p$ and $a$ and $b$ positive and finite real number. The hints are very clear and we have

$$
0 < ab \leq S_1+S_2
$$

\noindent But

$$
S_1=\int_{0}{a} x^{p-1} \ dx=\frac{a^p}{p}
$$

\bigskip \noindent and

$$
S_2=1=cb - \int_{0}{c} x^{p-1} \ dx=ab-\frac{c^p}{p}.
$$

\bigskip \noindent By using that $b=c^{p-1}$, we get $cb=b^{1+1/(p-1)}=b^{p/(p-1)}=b^q$, $c^p=b^q$ and hence $S_2=b^q-b^q/p=b^q(1-1/p)=b^q/q$. This finishes the solution.\\

\bigskip \noindent \textbf{Exercise 2}. (H\H{o}lder's Inequality) \label{exercise02_sol_doc09-03}\\

\noindent Let us denote for any $p\geq 1$, for any $f \in L^p$,

$$
N_p(f)=\left(\int |f|^p\right)^{1/p}.
$$

\bigskip \noindent Question (a) Let $p>1$ and $q>1$ be two conjugate real numbers. Let $f \in L^p$ and $f \in L^q$. Show that $fg \in L^1$ and

$$
\left\vert \int fg \ dm \right\vert \leq N_p(f) N_q(g). (H) 
$$

\bigskip \noindent \textit{Hints}.\\

\noindent (1) Remark that $f$ and $g$ are $m$-\textit{a.e.} finite by assumption. Tell why Formula (H) holds if $N_p(f)=0$ or $N_q(g)=0$.\\

\noindent (2) Suppose that none of $N_p(f)$ or $N_q(g)$ is zero. Apply Question (c) of Exercise 1 to

$$
a=\frac{f}{N_p(f)} \text{    and   } b=\frac{g}{N_q(g)}
$$

\bigskip \noindent to conclude.\\

\noindent Question (b) Extend Formula (H) in the following way : For any measurable application such that the product $N_p(f)N_q(g)$ makes senses regardless of any integrability condition on $f$ or on $g$,

$$
 \int \left\vert fg \right\vert \ dm  \leq N_p(f) N_q(g). (HE) 
$$

\bigskip \noindent \textit{Hints} Just discuss on the finiteness or not of $N_p(f)$ and $N_q(g)$\\

\bigskip \noindent \textbf{SOLUTIONS}.\\

\noindent (1). Since $|f|^p$ is integrable, $|f|^p$ and hence $f$ is finite $m$-\textit{a.e.}. And as a class, $f$ is finite. As well, $g$ is finite. Now if $N_p(f)=0$, we get that
$f$ is the zero. Since $g$ is finite, we also have $fg=0$ $m$-\textit{a.e.}. Hence, both members of Formula (H) are 0. We get the same conclusion if $N_q(g)=0$.\\

\noindent (2) Let us suppose that none of $N_p(f)$ or $N_p(f)$ is zero. Apply Question (c) of Exercise 1 to

$$
a=\frac{f}{N_p(f)} \text{    and   } b=\frac{g}{N_p(f)}.
$$

\bigskip \noindent We get
$$
\frac{|fg|}{N_p(f)N_q(g)}\leq \frac{|f|^p}{pN_p(f)^p}+\frac{|g|^q}{qN_q(g)^q}.
$$

\bigskip \noindent By remarking that $N_p(f)^p=\int |f|^p \ dm$ and $N_q(g)^q=\int |g|^q \ dm$, and by integrating the formula, we have

$$
\frac{\int |fg| \ dm }{N_p(f)N_q(g)}\leq \frac{\int |f|^p \ dm }{pN_p(f)^p}+\frac{\int |g|^q \ dm}{qN_q(g)^q}=\frac{1}{p}+\frac{1}{q}=1.
$$

\bigskip \noindent Hence $fg \in L^1$ and we conclude that

$$
\left\vert \int fg \ dm \right\vert \leq \int \left\vert fg \right\vert \ dm  \leq N_p(f) N_q(g). (H) \ \square
$$

\bigskip \noindent Question (b). All the integrals involved in Formula (HE) exist in $\overline{\mathbb{R}}$. First, if If $N_p(f)$ and $N_q(g)$  are both finite, we go back to Question (a) and Formula (HE) holds. If they are both both infinite, (HE) holds since the right-hand member is infinite. If one of them is infinite and the other finite, Formula (HE) holds for the same reason.\\

\bigskip \noindent \textbf{Exercise 3}. (Minkowski's Inequality) \label{exercise03_sol_doc09-03}\\

\noindent Let $p\geq 1$ and let that $(f,g) \in (L^p)^2$.\\

\noindent Question (a) Show that $f+g$ is $m$-\textit{a.e.} defined and $f+g \in L^p$ and show that $f+g\in L^p$.\\

\noindent \textit{Hint}. Use Question (a) of Exercise 1.\\

\noindent Question (b) Show that

$$
N_p(f+g) \leq N_p(f)+ N_p(g).
$$

\bigskip \noindent \textit{Hints}.\\

\noindent Give a quick solution for $p=1$.\\

\noindent What happens if $N_p(f)=0$ or $N_p(g)=0$?\\

\noindent Remark that there is nothing to for  $p>1$ but  $N_p(f+g)=0$.\\

\noindent Now let $N_p(f)\neq 0$, $N_p(g)\neq 0$, $N_p(f+g)\neq 0$ and $p>1$. Use the following formula, to be justify, 

\begin{equation*}
\left\vert f+g\right\vert ^{p}=\left\vert f+g\right\vert ^{p-1}\left\vert
f+g\right\vert \leq \left\vert f+g\right\vert ^{p-1}\left\vert f\right\vert
+\left\vert f+g\right\vert ^{p-1}\left\vert g\right\vert 
\end{equation*}

\noindent to arrive at 

\begin{equation*}
\left\vert f+g\right\vert ^{p}\leq \left\vert f+g\right\vert
^{p-1}\left\vert f\right\vert +\left\vert f+g\right\vert ^{p-1}\left\vert
g\right\vert. \ (PM)
\end{equation*}
 
\noindent Apply the inequality of Exercise 2 with $q=p/(p-1)$, which will be shown to be conjugate of $p$, to both terms in the right-hand member of Inequality (PM) and try to conclude.\\

\bigskip \noindent \textbf{SOLUTIONS}.\\

\noindent Question (a) Since $f$ and $g$ are in $L^p$, $p\geq 1$, they are both $m$-\textit{a.e.} finite and $f+g$ is $m$-\textit{a.e.} defined.\\

\noindent By Question (a) of Exercise 1, we have

$$
\int |f+g|^p \ dm \leq 2^{p-1} \biggr( \int |f|^p \ dm + \int |g|^p \ dm \biggr) <+\infty. \ \square
$$ 

\bigskip \noindent Question (b)\\

\noindent If $N_p(f)=0$, hence $f=0$ $m$-\textit{a.e.} and Formula (M) becomes $N_p(g)=N_p(g)$. As well, Formula (M) is  $N_p(f)=N_p(f)$ if $N_p(f)=0$.\\

\noindent If $p=1$, Formula (M) is readily proved by using the triangle inequality\\

\noindent If $p>1$ and $N_p(f+g)=0$, (M) holds whatever be the right-hand member of the inequality.\\

\noindent If $p>1$ and $N_p(f+g)\neq 0$, (M) holds whatever be the right-hand member of the inequality.\\

\noindent Now let $N_p(f)\neq 0$, $N_p(g)\neq 0$, $N_p(f+g)\neq 0$ and $p>1$. Using the triangle inequality, leads do  

\begin{equation*}
\left\vert f+g\right\vert ^{p}=\left\vert f+g\right\vert ^{p-1}\left\vert
f+g\right\vert \leq \left\vert f+g\right\vert ^{p-1}\left\vert f\right\vert
+\left\vert f+g\right\vert ^{p-1}\left\vert g\right\vert 
\end{equation*}

\noindent which, after integration, becomes 

\begin{equation*}
N_p(f+g)^p \leq \int \left\vert f+g\right\vert ^{p-1}\left\vert f\right\vert \ dm + \int \left\vert f+g\right\vert ^{p-1}\left\vert \ dm.
g\right\vert. \ (PM)
\end{equation*}
 
\bigskip \noindent Set $h=|f+g|^(p-1)$. The positive number$q=p/(p-1)$ is conjugate of $p$ and $h^q=|f+g|^{q(p-1)}=|f+g|^{p}$ and $N_q(h)^q=N_p(f+g)^p$. Hence, any of $N_q(h)$, $N_p(f)$ and $N_p(g)=0$ is different of zero. We may apply Question (b) of Exercise 2 to both terms in the right-hand member of Formula (PM) to get :

$$
\int h \left\vert f\right\vert \ dm \in N_p(f)N_q(h)= N_p(f) (N_p(f+g))^{p/q}=N_p(f) (N_p(f+g))^{p-1}
$$

\noindent and

$$
\int h \left\vert g\right\vert \ dm \in N_p(g)N_q(h)= N_p(g) (N_p(f+g))^{p/q}=N_p(g) (N_p(f+g))^{p-1}.
$$

\bigskip \noindent We have

$$
N_p(f+g)^p \leq N_p(f) (N_p(f+g))^{p-1} +  N_p(g) (N_p(f+g))^{p-1}.
$$

\bigskip \noindent Since $N_p(f+g)$ is finite  by Question (a) and since it is positive, we may divide both members by it to conclude. $\blacksquare$\\

\bigskip \noindent \textbf{Exercise 4}. ($L^p$ is a Banach Space, $1\geq p <+\infty$) \label{exercise04_sol_doc09-03}\\

\noindent Let $p\geq 1$ be finite real number.\\

\noindent Question (a) Show that the application

$$
L^p \ni f \mapsto N_p(f)
$$

\bigskip \noindent is norm. From now one, put $N_p(f)=\|f\|_p$.\\

\noindent \textit{Hint}. Use the results of Exercises 2 and 3.\\

\noindent Question (b) Show (describe why?) that $(L^p, +, .)$ is a vector space, where the last operation is the external multiplication of function by real scalars.\\

\noindent Question (c) Show that when endowed with the norm $\|f\|_p$, $(L^p, +, .)$ is complete normed space, that is a Banach space.\\

\noindent \textit{Hints}. Consider a Cauchy sequence $(f_n)_{n\geq 0}$ in $(L^p, +, ., \|f\|_p)$.\\

\noindent (i) Use the Markov inequality (See Exercise 8 in Doc 05-05, page \pageref{exercise08_sol_doc05-05}, in Chapter \ref{05_integration})  and show that the sequence in Cauchy in measure.\\

\noindent (ii) By Exercise 6 in Doc 06-08, (see page \pageref{exercise06_sol_doc06-08}) in Chapter \ref{06_convergence}, consider a sub-sequence $(f_{n_k})_{k\geq 0}$ of $(f_n)_{n\geq 0}$ converging both in measure and $m$-\textit{a.e.} to a measurable function $f$.\\

\noindent Apply the (below version) Fatou-Lebesgue's lemma  (See Exercise 2 in Doc 06-09, page \pageref{exercise02_sol_doc06-09} in Chapter \ref{06_convergence}) to 
$$
\int |f_n-f_{n_k}|^p \ dm
$$ 

\bigskip \noindent as $k \rightarrow +\infty$ and use the \textit{Cauchyness} of the sequence to conclude.\\

\bigskip \noindent \textbf{SOLUTIONS}.\\

\noindent Question (a) Let us show the the application satisfies the three conditions of a norm.\\

\noindent (i) It is clear that  $N_p(f)=0$ iff $\int |f|^p \ dm=0$ iff $|f|^p=0$ $m$-\textit{a.e.} iff $f=0$ $m$-\textit{a.e.} iff $f=0$ as an equivalence class.\\

\noindent (i) For any $\lambda \in \mathcal{R}$, for any $f \in L^p$, we have

\begin{eqnarray*}
N_p(\lambda f)=\biggr( \int |\lambda f|^p \ dm \biggr)^{1/p}&=& \biggr( |\lambda |\int  |f|^p \ dm \biggr)^{1/p}\\
&=&|\lambda| \biggr( \int  |f|^p \ dm \biggr)^{1/p}=|\lambda| N_p(f),
\end{eqnarray*}

\bigskip \noindent and by the way, $\lambda f \in L^p$.\\

\bigskip \noindent (iii) Finally, we have the triangle inequality, that is the Minkowski's inequality : For all $f \in L^p$ and $g \in L^p$, we have $f+g \in L^p$ and

$$
N_p(f+f) \leq N_p-f) +N_p(g).
$$

\noindent Thus, $N_p$ is a norm on $L^p$.

\noindent Question (b) By Points (ii) and (iii) above, we see that $(L^p, +, .)$ is a linear space. Actually, it is a linear sub-space of all measurable and $m$-\textit{a.e.} finite functions endowed with the same operations.\\

\noindent Question (c) Let  $(f_n)_{n\geq 0}$ be a Cauchy sequence in $(L^p, +, ., \|f\|_p)$. We have to prove it converges in $(L^p, +, ., \|f\|_p)$.\\

\noindent (i) By Markov's inequality, we have for all $n\geq 0, \ m\geq 0$, for all $\varepsilon>0$,

$$
m(|f_n-f_s|>\varepsilon) \leq \frac{1}{\varepsilon} \|f_n-f_s\|_p^p \rightarrow 0 \text { as } (n,s)\rightarrow (+\infty,+\infty).
$$

\noindent (ii) Hence $(f_n)_{n\geq 0}$ is a Cauchy sequence in measure. By Exercise 6 in Doc 06-08, (see page \pageref{exercise06_sol_doc06-08}) in Chapter \ref{06_convergence}, consider a sub-sequence $(f_{n_k})_{k\geq 0}$ of $(f_n)_{n\geq 0}$ converging both in measure and $m$-\textit{a.e.} to a measurable function $f$.\\

\noindent Let us apply the Fatou-Lebesgue lemma's (See Exercise 2 in Doc 06-09, page \pageref{exercise02_sol_doc06-09} in Chapter \ref{06_convergence}), we have for all $n\geq 0$,

$$
\int \liminf_{k\rightarrow +\infty} |f_n-f_{n_k}|^p \ dm \leq \liminf_{k\rightarrow +\infty} \int |f_n-f_{n_k}|^p \ dm, 
$$

\bigskip \noindent that is, for all $n\geq 0$

$$
\int |f_n-f|^p \ dm \leq \liminf_{k\rightarrow +\infty} \|f_n-f_{n_k}\|_p^p, 
$$

\bigskip \noindent that is, by taking the $(1/p)$ powers,

$$
\|f_n-f\|^p \ dm \leq \liminf_{k\rightarrow +\infty} \|f_n-f_{n_k}\|_p, \ (PC1). 
$$

\bigskip \noindent Since $(f_n)_{n\geq 0}$ is a Cauchy sequence in $(L^p, +, ., \|.\|_p)$, there exists for $\varepsilon$, an integer $N_0\geq 0$, such that for all $n\geq N_0$ and $s\geq N_0$

$$
\|f_n-f_s\|^p <\varepsilon.
$$

\bigskip \noindent Since the sequence $(n_k)_{k\geq 0}$ is strictly increasing, there also exists $K_0$ such that $n_k\geq N_0$ for $k\geq K_0$. As a first implication of (PC1), we have 

$$
\|f\|_p \leq \|f_{N_0}\|_p + \|f_{N_0}-f\|_p\leq \|f_{N_0}\|_p + \varepsilon,
$$

\bigskip \noindent which ensures that $f\in L^p$. A second implication of (PC1) that for all $n\geq N_0$ and for $k\geq K_0$, 

$$
\|f_n-f_{n_k}\|_p \leq \varepsilon
$$

\bigskip \noindent and next, for all $n\geq N_0$,

$$
\|f_n-f\|^p \ dm \leq \varepsilon, 
$$

\bigskip \noindent which leads to

$$
\liminf_{n \rightarrow +\infty} \|f_n-f\|^p \ dm \leq \varepsilon,
$$

\bigskip \noindent for all $\varepsilon>0$. This ensure that $f_n \rightarrow f$ in $(L^p, +, ., \|.\|_p)$.

\bigskip \noindent \textbf{Exercise 5}. \label{exercise05_sol_doc09-03}\\

\noindent Define $\mathcal{L}^{^\infty}$ as the class of the equivalence classes of measurable functions $m$-\textit{a.e.} bounded. An other to clearly define it the following. Define

$$
\mathcal{L}^{\infty}=\{f \in \mathcal{L}_0, \ \exists M\in [0, +\infty[, \ |f|\leq M, \ m-a.e.\}.
$$

\bigskip \noindent Next, define  $L^{\infty}=\mathcal{L}/\mathcal{R}$ as the quotient space of by the equivalence relation $\mathcal{R}$ of the $m$-\textit{a.e.} equality. For $f \in L^{\infty}$, define

$$
N_{\infty}(f)= \inf \{M\in [0, +\infty[, \ |f|\leq M, \ m-a.e.\}.
$$ 

\bigskip \noindent Question (a) Show that an $m$-\textit{a.e.} bound $M$ of $f$ is also an $m$-\textit{a.e.} of any elements of the class $\overline{f}$. Apply this result to establish that $N_{\infty}$ is well-defined on $L^{\infty}$.\\
 
\noindent Question (b) Show that $(L^{\infty},+, .)$ is a linear space.\\

\noindent Question (c) Show that $N_{\infty}$ is a norm on $L^{\infty}$. From now on, we denote it $\|.\|_{\infty}$.\\

\noindent Question (d) Show that $(L^{\infty},+, ., \|.\|_{\infty})$ is a real Banach space.\\

\bigskip \noindent \textbf{SOLUTIONS}.\\

\noindent If $M$ is an $m$-\textit{a.e.} bound of $f$, there exists a measurable sub-space $\Omega_1$ such that $m(\Omega_1^c)=0$ and $|f|\leq M$ on $\Omega_1$. If $g \in \overline{f}$, there exists
a measurable sub-space $\Omega_2$ such that $m(\Omega_2^c)=0$ and $f=g$ on $\Omega_2$. Define $\Omega_0=\Omega_1 \cap \Omega_2$. We have $m(\Omega_0^c)=0$ and $|g|\leq M$ on $\Omega_0$. Thus
$M$ is an $m$-\textit{a.e.} bound $g$.\\

\noindent We conclude that two measurable functions $f$ and $g$ such that $f=g$ $m$-\textit{a.e.} have the same $m$-\textit{a.e.} bounds or not. Hence the definition

$$
N_{\infty}(\overline{f})= \inf \{M\in [0, +\infty[, \ |h|\leq M, \ m-a.e.\}
$$

\noindent for $h\in N_{\infty}(\overline{f})$ does not depend on the represent $h$ used in the formula above.\\

\noindent (b) Let us show the two points.\\

\noindent (b1) Let $f \in L^{\infty}$ and $\lambda \in \mathcal{R}$. Then $\lambda f  \in (L^{\infty}$.\\

\noindent But we  have : $f \in (L^{\infty}$ iff there exists $M \in [0, +\infty[$ such that $|f|\leq M$ $m$-\textit{a.e.}. Thus, 
$|\lambda f|\leq |\lambda| M<+\infty$. Thus $\lambda f \in L^{\infty}$.

\noindent (b2) Let $(f,g) \in L^{\infty})^2$. Then $f+g \in L^{\infty}$.\\

\noindent But we  have : $f \in (L^{\infty}$ iff there exists $M \in [0, +\infty[$ such that $|f|\leq M$ $m$-\textit{a.e.}. Thus, 
$|\lambda f|\leq |\lambda| M<+\infty$. Thus $\lambda f \in L^{\infty}$.\\

\noindent The assumptions implies there exist measurable sub-space $\Omega_i$ and $M_i\in [0, +\infty[$ such that $m(\Omega_i^c)=0$, $i\in \{1,2\}1$, and $|f|\leq M_1$ on $\Omega_1$ and 
$|g|\leq M_2$ on $\Omega_2$. Define $\Omega_0=\Omega_1 \cap \Omega_2$. We have $m(\Omega_0^c)=0$ and $|f+g|\leq M_1+M_2$ on $\Omega_0$. Thus $f+g \in L^{\infty}$.\\

\noindent In conclusion $L^{\infty}$ is a linear space.\\

\noindent Question (c). Let us begin to characterize the infimum by saying that for any $f \in L^{\infty}$,

$$
\forall \varepsilon >0, \ \exists M \in [0, +\infty[, \  0\leq M<\leq N_{\infty}(f)+\varepsilon \text{ and } |f|\leq M, \ m-a.e.,
$$

\bigskip \noindent or by taking a maximizing sequence : there exists a sequence $(M_n)_{n\geq 0} \subset [0, +\infty[$, there exists a sequence of measurable subsets $(\Omega_n)_{n\geq 0}$ such that :\\

\noindent (c1) $\forall n\geq 0$, $m(\Omega_n^c)=0$ and $|f|\leq M_n$ on $\Omega_n$.\\

\noindent (c2) $M_n \rightarrow N_{\infty}(f)$ as $n\rightarrow +\infty$.\\

\noindent (i) Let us prove that $N_{\infty}(f)=0$ iff $\overline{f}=0$. Is is clear that $N_{\infty}(f)=0$ if $\overline{f}=0$. Now suppose $N_{\infty}(f)=0$. Let us use Point (c2) above and take
$\overline{\Omega}=\bigcap_{n\geq 0} \Omega_n$. We have $m(\overline{\Omega})=0$ and 

$$
\forall n\geq 0, \ \ |f|\leq M_n \text{ on } \overline{\Omega}.
$$ 

\bigskip \noindent By letting $n\rightarrow+\infty$, we get $f=0$ on the co-null set $\overline{\Omega}$. Thus $f=0$ $m$-\textit{a.e.} and hence $\overline{f}=0$.\\

\noindent (ii) Let $f \in L^{\infty}$ and $\lambda \in \mathcal{R}$. Let us show that $N_{\infty}(\lambda f)=|\lambda| N_{\infty}(f)$.\\

\noindent If $\lambda=0$, it is obvious $\lambda f=0$ since $f$ is \textit{a.e.} bounded and hence $N_{\infty}(\lambda f)=0$ and $|\lambda| N_{\infty}(f)=0$.\\

\noindent If $\lambda \neq 0$, the formula comes from the following identity :

\begin{eqnarray*}
\inf \{M\in [0, +\infty[, \ |\lambda| |f|\leq M, \ m-a.e.\}&=&\inf \{M\in [0, +\infty[, \ |f|\leq M/|\lambda|, \ m-a.e.\}\\
&=& |\lambda| \inf \{C \in [0, +\infty[, \ |f|\leq C, \ m-a.e.\}.
\end{eqnarray*}

\bigskip \noindent (iii) Let $(f,g) \in L^{\infty})^2$. Let us show that $N_{\infty}(f+g) \leq N_{\infty}(f) + N_{\infty}(g)$.\\

\noindent Let us use Point (c1) above. Let $\varepsilon>0$. There exists $M_1$ and $M_2$ such that 

$$
0\leq M_1<\leq N_{\infty}(f)+\varepsilon/2 \text{ and } |f|\leq M_1, \ m-a.e.,
$$

\bigskip \noindent and

$$
0\leq M_2<\leq N_{\infty}(g)+\varepsilon/2 \text{ and } |g|\leq M_2, \ m-a.e.
$$

\bigskip \noindent Hence for all $\varepsilon>0$,

$$
0\leq M_1+M_2<\leq N_{\infty}(f)+N_{\infty}(g)+\varepsilon \text{ and } |f+g|\leq M1+M_2, \ m-a.e.
$$

\bigskip \noindent By $|f+g|\leq M1+M_2$  $m$-a.e., we deduce that $N_{\infty}(f+g)\leq M_1+M_2$ and thus for all $\varepsilon>0$,

$$
N_{\infty}(f+g) \leq N_{\infty}(f)+N_{\infty}(g)+\varepsilon.
$$
 
\bigskip \noindent We conclude by letting $\varepsilon \rightarrow 0$.\\

\noindent In conclusion $N_{\infty}(.)=\|.\|_{\infty}$ is a norm.\\

\noindent Question (d) Let $(f_n)_{n\geq 0}$ be a Cauchy sequence in $(L^{\infty},+, ., \|.\|_{\infty})$. Thus for all $k\geq 0$, there exists $n_k\geq 0$, for all $p\geq 0$, $q\geq 0$, $\|f_{n_k+p}-f_{n_k+q}\|_{\infty} < 1/(2(k+1))$, that is

$$
\inf \{M\in [0, +\infty[, \ |f_{n_k+p}-f_{n_k+q}|\leq M, \ m-a.e.\} < 1/(2(k+1)).
$$
 
\bigskip \noindent Hence for  all $p\geq 0$, $q\geq 0$, there exists $M(p,n)<1/(k+1)$ such that $|f_{n_k+p}-f_{n_k+q}|\leq M(p,q)\leq 1/(k+1)$, say on the measurable sub-space $\Omega_{k,p,q}$ such that
$m(\Omega_{p,q}^c)=0$.\\

\bigskip \noindent Set $\Omega_0=\bigcap_{(k,p,q)\in \mathbb{N}^3} \Omega_{k, p,q}$. We clearly have $m(\Omega_0^c)=0$.\\

\noindent Hence, for each $\omega \in \Omega$, $(f_n(\omega))_{n\geq 0}$ is a Cauchy sequence in $\mathbb{R}$. So it converges to a real number $f(\omega)$. We get that

$$
f_n 1_{\Omega_0} \rightarrow f_0 \text{ as } n\rightarrow +\infty.
$$

\bigskip \bigskip \noindent So $f$ is a measurable application which is defined on $(\Omega_0, \mathcal{A}_{\Omega_0})$. We  may extend it to
$$
f=f_0 1_{\Omega_0} + \infty 1_{\Omega_0^c}.
$$ 

\bigskip \noindent We clearly have $f_n \rightarrow f$, $m$-\textit{a.e.}. Besides, we have on  $\Omega_0$, on each $k\geq 0$, for $p\geq 0$, $q\geq 0$,

$$
|f_{n_k+p}-f_{n_k+q}|\leq 1/(k+1).
$$

\bigskip \noindent By letting $q \rightarrow +\infty$, we get 

$$
|f_{n_k+p}-f|\leq 1/(k+1) \text{ on } \Omega_0,
$$

\bigskip \noindent that is $\|f_{n_k+p}-f\|_{\infty}\leq 1/(k+1)$. This finally means that for all $\varepsilon>0$, there exists $k\geq 0$ such that $1/(k+) <\varepsilon$ and there exists $N=n_k$ such that for all $n\geq N$,

$$
\|f_{n_k+p}-f\|_{\infty} \leq \varepsilon,
$$

\bigskip \noindent which means that $f_n$ converges to $f$ in $(L^{\infty},+, ., \|.\|_{\infty})$.\\

\bigskip \noindent \textbf{Exercise 6}. ($L^2$ is a Hilbert Space) \label{exercise06_sol_doc09-03}\\

\noindent Question (a) Express the H\H{o}lder inequality for $p=q=2$, called the Cauchy-Schwarz inequality. Deduce that for $(f,g)\in (L^2)^2$, $fg \in L^1$.\\

\noindent Question (b) Define the application 

$$
(L^2)^2 \ni (f,g) \mapsto Q(f,g)=\int fg \ dm.
$$

\bigskip \noindent Show that $\Phi$ is a symmetrical bi-linear application and that $Q(f,f)=\|f\|_2$ is a scalar product. Conclude that $L^2$ is a Hilbert space.\\

\bigskip \noindent \textbf{SOLUTIONS}.\\

\noindent (a) The bi-linearity of $Q$ is a directed consequence of the linearity of the integral over the class of integrable function.\\

\noindent (b) The quadratic form $Q(f,f)=\|f\|_2^2$ is definite positive since $Q(f,f)>0$ for $f\neq 0$ (as an equivalence class), that $f\neq 0 \ m.a.e.$ as a simple function.\\

\noindent Now, the norm of complete space $(L^2, \|.\|_2)$ derives from the  quadratic form $Q(f,f)$. Hence  $(L^2, \|.\|_2)$, by definition, is Hilbert Space.\\

\bigskip \noindent \textbf{Exercise 7}. (Convergence) \label{exercise07_sol_doc09-03}\\

\noindent Question (1) Suppose that a sequence $(f_n)_{n\geq 0} \subset L^p$, $p\in [1,+\infty[$, converges in $L^p$ to $f$, denoted by 

$$
f_n \overset{L^p}{\longrightarrow} f.
$$

\noindent Establish the following implications. \\

\noindent (a) $(f_n)_{n\geq 0}$ is Cauchy in measure.\\

\noindent (b) There exists a sub-sequence $(f_{n_k})_{k\geq 0}$ of $(f_n)_{n\geq 0}$ converging $m$-\textit{a.e.} and in measure to $f$.\\
 
\noindent (c) $(f_n)_{n\geq 0}$ converges in measure to $f$.\\

\noindent (d) $\|f_n\|_p$ converges to $\|f\|_p$.\\

\noindent \textit{Hints}. For Parts (a), (b) and (c), just recall the results above. For Part (d), use the second triangle inequality for any norm : $\left\vert \|x\|-\|y\| \right\vert \leq \|x-y\|$.\\

\noindent Question 2) Suppose that a sequence $(f_n)_{n\geq 0} \subset L^{\infty}$, converges in $L^{\infty}$ to $f$. Show that $f_n \rightarrow f $ $m$-\textit{a.e.} 

\noindent \textit{Hint}. Just recall a result of Exercise 4.\

\bigskip \noindent \textbf{SOLUTIONS}.\\

\noindent Question (1).\\

\noindent (a) Since $(f_n)_{n\geq 0}$ converges in $L^p$, it is Cauchy in $L^p$. Thus, by Question (c) in Exercise 4, it is Cauchy in measure.\\

\noindent (b) See Exercise 6 in Doc 06-08, (see page \pageref{exercise06_sol_doc06-08}) in Chapter \ref{06_convergence}.\\

\noindent (c) See also Exercise 6 in Doc 06-08, (see page \pageref{exercise06_sol_doc06-08}) in Chapter \ref{06_convergence}.\\

\noindent (d) By using the second triangle inequality, we have

$$
\left\vert \|f_n\|_p-\|f\|_p \right\vert \leq \|f_n-f\|_p.
$$

\bigskip \noindent Hence $\|f_n\|_p$ converges to $\|f\|_p$ whenever $f_n$ converges to  $f$ in $L^p$.\\

\noindent Question (2). This is a by-product of the proof of Question (b) in Exercise 5. Since $(f_n)_{n\geq 0}$ is Cauchy in $L^{\infty}$ if it converges in $L^{\infty}$.\\

\noindent \LARGE \textbf{DOC 09-04 : Daniel Integral - Exercises} \label{doc09-04}\\
\bigskip
\Large

\noindent \textit{I - Introduction}.\\

\noindent The Daniel-Stone integral is a kind of reverse of the following results. Let $(\Omega, \mathcal{A},m)$ be a measure space. We have the following facts :\\

\noindent (L0) $\mathcal{F}=\mathcal{L}^1(\Omega, \mathcal{A},m)$ is a lattice linear space, meaning that : it is a linear space and it contains finite maxima and minima of its elements, where the maximum and the minimum are relative to the order relation on real-valued functions.\\

\noindent (L1) The functional

$$
f \ni \mathcal{F} \mapsto L(f)=\int f \ dm, \ (DS)
$$ 

\bigskip \noindent is linear and is non-negative, that if $L(f)\geq 0$ for $f\geq 0$.\\

\noindent (L2) For any sequence $(f_n)_{n\geq 0}\subset \mathcal{F}$ of non-decreasing to zero of non-negative functions such that $L(f_{n_0})<+\infty$ for some $n_0\geq 0$, we have

$$
\lim_{n\rightarrow +\infty} L(f_n)=0.
$$

\bigskip \noindent If $m$ is a probability measure we have in addition :\\

\noindent (L3) $1 \in \mathcal{F}$ and $L(1)=1$.\\

\noindent The principle of the Daniel integral is to begin with a functional $L$ on a class of functions $\mathcal{F}$ satisfying a set of assumptions and to construct a measure space
$(\Omega, \mathcal{A},m)$ such that $\mathcal{F} \subset \mathcal{L}^1(\Omega, \mathcal{A},m)$ and $L$ is of the form (DS).\\

\noindent Several statements, which are more or less general, are available (See for example \cite{bogachev2}). Here, we will restrict ourselves in the case where $m$ is a probability measure. Other forms are more or less generalizations of the results here given.\

\noindent This important construction of a measure plays an important role in the theory of weak convergence of probability measures.\\

\bigskip \noindent \textit{II - Simple statement}.\\

\noindent Let $\mathcal{F}$ be a lattice vector space of functions defined on $\Omega$ with values in $\mathbb{R}$ and $L$ be a linear form defined on $\mathcal{F}$ such that the assumptions (L1)-(L3) hold. Then there exists a unique probability measure $m$ on $(\Omega,  \mathcal{A}_{\mathcal{F})}$ such that :\\

\noindent (i)  $\mathcal{F} \subset \mathcal{L}^1(\Omega, \mathcal{A},m)$\\

\noindent and \\

\noindent (ii) for all $f \in \mathcal{F}$

$$
L(f)=\int f  \ dm,
$$

\bigskip \noindent where $\mathcal{A}=\mathcal{A}_{\mathcal{F})}=\sigma({\mathcal{F})})$ is the $\sigma$-algebra generated by the class $\mathcal{F}$, that is generated by the class of subsets of $\Omega$

$$
\mathcal{C}_{\mathcal{F}} =\{\ f^{-1}(B), \ f \in \mathcal{F}, \ B \in \mathcal{B}_{\infty}(\overline{\mathbb{R}}) \}.
$$

\bigskip \noindent \textit{III - Construction of the Daniel integral by the exercise}.\\

\noindent \textbf{NB}. These exercises are advanced ones. They greatly use the knowledge of fundamental measure theory which is constituted by the contents up to general $L^p$ spaces. In particular, a number of methods used in the study of elementary functions are re-conducted here.\\

\noindent In all the remaining part of the document, we suppose that the assumptions of the statement in Part II hold.\\

\noindent  Denote $\mathcal{F}^+=\{f \in\mathcal{F}, \ f\geq 0\}$. Define $\mathcal{L}^+$ the class of all non-negative and bounded real-valued functions defined on $\Omega$ which are limits of non-decreasing sequences of elements of $\mathcal{F}$, that is

$$
f \in \mathcal{L}^+ \Leftrightarrow f \text{ is bounded and } \exists (f_n)_{n\geq 0} \subset \mathcal{F}^+, \ f_n\nearrow f \ as \ n \nearrow +\infty. 
$$

\bigskip \noindent \textbf{Exercise 1}. \label{exercise01_doc09-04} \textbf{(Extension of $L$)}.\\

\noindent Question (a). Show $L$ is non-decreasing that is $L(f)\leq L(g)$, for any couple of elements of $\mathcal{F}$ such that $f\leq g$.\\

\noindent \textit{Hints}. Use Assumptions (L0) and (L1) to $h=f-g\geq 0$ for $f\leq g$.\\

\noindent Question (b). Show that if  $(f_n)_{n\geq 0}$ is non-decreasing sequences of elements of $\subset \mathcal{F}^+$ converging to  $f\in \mathcal{F}$, then

$$
\lim_{n\rightarrow +\infty} L(f_n)=L(f).
$$

\bigskip \noindent \textit{Hints}. Consider $g_n=f-f_n$, $n\geq 0$, and use (L1) and (L2).\\

\noindent Question (c). Let $(f_n)_{n\geq 0}$ and $(g_n)_{n\geq 0}$ be two non-decreasing sequences of elements of $\subset \mathcal{F}^+$. Show that

$$
\lim_{n\rightarrow +\infty} f_n \leq \lim_{n\rightarrow +\infty} g_n \Rightarrow \lim_{n\rightarrow +\infty} L(f_n) \leq \lim_{n\rightarrow +\infty} (g_n).
$$

\bigskip \noindent \textit{Hints}. Suppose that $\lim_{n\rightarrow +\infty} f_n \leq \lim_{n\rightarrow +\infty} g_n $ and apply Question (a) to the relation : 

$$
\forall n\geq 0, \ f_n=\lim_{k\rightarrow +\infty} min(f_n,g_k). \ (CP)
$$

\bigskip \noindent Question (d) Use Question (c) to justify the extension $L^{*}$ of $L$ on $\mathcal{L}^+$ by taking $L^{*}(f)$ as the limit of any sequence $L(f_n)_{n\geq 0}$ where  
$(f_n)_{n\geq 0} \subset \mathcal{F}^+$ is such that $f_n\nearrow f$  as $n \nearrow +\infty$, that is :\\

$$
L^{\ast}(f)=\lim_{n\rightarrow +\infty} L(f_n), \ (f_n)_{n\geq 0} \subset \mathcal{F}^+,  \ f_n\nearrow f \ as \ n \nearrow +\infty.
$$

\bigskip \noindent Question (e) Establish the following simple properties of $L$ of $\mathcal{L}^+$ :\\

\noindent (e1) $L$ is the restriction of $L^{\ast}$ on $\mathcal{F}_b^+$, the class of bounded and non-negative elements of $\mathcal{F}$. \textbf{In the sequel, $L^{\ast}$ is simply denoted by $L$}.\\

\noindent (e2) $L$ is non-decreasing on $\mathcal{L}^+$, that is $L(f)\leq L(g)$ for any couple of elements of $\mathcal{L}^+$ such that $f\leq g$ (use again Formula (CP) above).\\

\noindent (e3) $\mathcal{L}^+$ is lattice and is closed under finite sum of functions and under multiplication by non-negative scalars, that is : for any couple $(f,g)$ of elements of $\mathcal{L}^+$ and for any non-negative real scalars $c\geq 0$ and $d\geq 0$, we have

$$
L(cf+dg)=c L(f)+ d L(g).
$$

\bigskip \noindent (e4) Let $(f,g)$ be any any couple of elements of $\mathcal{L}^+$, then

$$
L(f)+L(g)=L(\max(f,g))+L(\min(f,g)),
$$

\bigskip \noindent from the identity $f+g=\max(f,g)+\min(f,g)$.\\

\bigskip \noindent Question (f) (A form the Monotone Convergence Theorem). Let $(f_n)_{n\geq 0}$ be a non-decreasing sequences of elements of $\subset \mathcal{L}^+$ which is uniformly bounded, that is, there exists $M\in \mathbb{R}$ such that

$$
\forall n\geq 0, \ for all \omega \in \Omega, \ 0\leq f(\omega) \leq M.
$$
 
\bigskip \noindent Show that

$$
\lim_{n\rightarrow +\infty} L(f_n)= L(\lim_{n\rightarrow +\infty} f_n).
$$

\bigskip \noindent \textit{Suggestions}. First, remark that $f=\lim_{n\rightarrow +\infty} f_n \in \mathcal{L}^+$ which justifies the writing of $L(\lim_{n\rightarrow +\infty} f_n)$. From there, the method used in 
Exercise 1 in Doc 06-09 of Chapter \ref{06_convergence} (see page \pageref{exercise01_sol_doc06-09}) may be re-conducted word by word. If you have done that exercise, you will gain nothing by doing it again. If not, simply re-conduct that reasoning. \\

\bigskip \noindent \textbf{Exercise 2}. \textbf{(Creation of $m$)}\label{exercise02_doc09-04}\\

\noindent Define 

$$
\mathcal{C}=\{A \subset \Omega \ : \ 1_A \in , \mathcal{L}^+\}.
$$

\bigskip \noindent Question (a) Show that $\{\emptyset, \Omega\}\subset \mathcal{C}$. Show that $\mathcal{C}$ is closed under finite intersections (that is, it is $\pi$-system) and countable unions.\\

\noindent \textit{Hints} To show that $\mathcal{C}$ is closed under finite intersections and unions, take only intersections or unions of two of its elements, which is enough thanks to the associativity of such operations. To show that To show that $\mathcal{C}$ is closed under countable unions, use partial unions and use Question (f) of Exercise 1.\\

\noindent Question (b) Define 

$$
\mathcal{C} \ni A \mapsto m^{\ast}(A)=L(1_A).
$$

\bigskip \noindent Show that $m^{\ast}$ is proper, non-negative, additive and continuous below, normed ($\mu^{\ast}(\Omega)=1$) and satisfies for all $G_i\in \mathcal{C}$, $i\in\{1,2\}$ 

$$
m^{\ast}(G_1\cap G_2) + m^{\ast}(G_1\cap G_2)=m^{\ast}(G_1)+m^{\ast}(G_2).
$$

\bigskip \noindent \textit{Hints}. All these facts are direct consequences of Exercise 1. To justify this, just indicate the exact formulas from which they derive.\\

\noindent Question (b) Show that the application

$$
\mathcal{P}(\Omega) \ni B \mapsto \mu(A)=\inf\{\mu^{\ast}(B), \ B \in \mathcal{C}, A \subset B\}
$$

\bigskip \noindent is an extension of $m$ and is an outer measure.\\

\noindent Question (c) By Part VII in Doc 04-01 (and its proof in Doc 04-10), $m$ is a measure on the $\sigma$-algebra 

$$
\mathcal{A}_0=\{A \in \mathcal{P}(\Omega) : \ \forall D\subset \Omega, \ \mu(D)=\mu(AD)+\mu(A^{c}D) \},
$$ 

\bigskip \noindent Tell why $m$ is the unique extension of $m^{\ast}$ to a measure, that $m$ is a probability measure. Use the additivity of $m^{\ast}$ of $\mathcal{C}$ to see that 
$\mathcal{C} \subset \mathcal{A}_0$.\\

\noindent Question (d) \noindent Let $\mathcal{A}=\mathcal{A}_{\mathcal{F}}$ the $\sigma$-algebra generated by $\mathcal{F}$, that is

$$
\mathcal{A}=\sigma(\{\ f^{-1}(B), \ f \in \mathcal{F}, \ B \in \mathcal{B}_{\infty}(\overline{\mathbb{R}}) \}).
$$

\bigskip \noindent Show that $\mathcal{A} \subset \mathcal{A}_0$.\\

\noindent \textit{Hints}.\\

\noindent For $f\in \mathcal{F}$, for any $c>0$, justify

$$
1_{(f>c)}=\lim_{n\rightarrow +\infty}  \min(1, \ n(f-c)^+). \ (AP01)
$$

\bigskip \noindent One one hand, deduce that $\mathcal{A}_{\mathcal{L}^+} \subset \mathcal{A}_{\mathcal{F}}$. Also easily check that, for $\mathcal{L}_0=\{1_A, \ A\in \mathcal{C} \}$, we have 
$\mathcal{A}_{\mathcal{L}_0} \subset \mathcal{A}_{\mathcal{L}^+}$.\\

\noindent On another hand show that $\mathcal{A}_{\mathcal{L}^+} \subset \mathcal{A}_{\mathcal{F}}$ from the definition $\mathcal{L}^+$. Conclude that $\mathcal{A}_{\mathcal{L}^+}=\mathcal{A}_{\mathcal{F}}$.

\noindent Finally, use the additivity of $m^{\ast}$ to show that $\mathcal{C}\in \mathcal{A}_0$. Make your final conclusion.\\

\noindent Question (e). Let $f\in \mathcal{F}$.\\

\noindent (e1) Check that the formula

$$
L(f) = \int f \ dm, \ (AP)
$$

\bigskip 
\noindent holds of $f\in \mathcal{L}_0$, that is $f$ is of the form $f=1_A$, $A \in \mathcal{A}$. Extend the validity of (AP) on the class of finite linear combinations of elements of $f\in \mathcal{L}_0$.\\

\noindent (e2) Let $f\in \mathcal{F}_b^+$, bounded by $M\geq 0$. Use Exercise 3 in Doc 03-08, page \pageref{exercise03_sol_doc03-08}, show how to drop the part relative to the segment 
$[2^n, +\infty[$, exchange the strict and broad inequalities to see that the sequence

\begin{eqnarray*}
f_n&=&\sum_{k=0}^{{k_M}-1} \frac{k-1}{2^n} \biggr(1_{(f>(k-1)/2^n)} - 1_{(f>k/2^n)} \biggr) + \biggr(1_{(f>(k_M-1)/2^n)} - 1_{(f>M)}.
\end{eqnarray*} 

\noindent broadly increases to $f$. Next use the remarkable identity

$$
\sum_{i=1}^{p} \alpha_{i-1} (\beta_{i-1}-\beta_i)=\alpha_0 \beta_0 + \biggr( \sum_{i=1}^{p-1} (\alpha_i - \alpha_{i-1} \biggr) +\alpha_{p-1} \beta_p,
$$

\bigskip \noindent where with $(\alpha_i)_{0\leq i \leq p_n}$ and $(\beta_i)_{0\leq i \leq p_n}$, $p\geq 1$, are finite real numbers and extend the formula (AP) to $f\in \mathcal{F}_b^+$. Be rigorous in justufying the limits.\\

\noindent (e3) For $f\in \mathcal{F}^+$, in general, use the approximation

$$
f=\lim_{p\rightarrow+\infty} \max(f,n)
$$

\bigskip \noindent and extend (AP) for $f\in \mathcal{F}^+$.  Give all the details.\\

\noindent (e4) For your general conclusion, use $f=f^+-f^-$.\\

\noindent (e5) Do we have $\mathcal{F} \subset \mathcal{L}^1(\omega, \mathcal{A},m)$?

\noindent \LARGE \textbf{DOC 09-05 : Daniel Integral - Exercises with solutions} \label{doc09-05}\\
\bigskip
\Large

\noindent In all the remaining part of the document, we suppose that the assumption of the statement part II (in Doc 09-04, page \pageref{doc09-04}) hold.\\

\noindent \textbf{Exercise 1}. \label{exercise01_sol_doc09-05}\\

\noindent Question (a). Show $L$ is non-decreasing that is $L(f)\leq L(g)$, for any couple of elements of $\mathcal{F}$ such that $f\leq g$.\\

\noindent \textit{Hints}. Use Assumptions (L0) and (L1) to $h=f-g\geq 0$ for $f\leq g$.\\

\noindent Question (b). Show that if  $(f_n)_{n\geq 0}$ is non-decreasing sequences of elements of $\subset \mathcal{F}^+$ converging to  $f\in \mathcal{F}$, then

$$
\lim_{n\rightarrow +\infty} L(f_n)=L(f).
$$

\bigskip \noindent \textit{Hints}. Consider $g_n=f-f_n$, $n\geq 0$, and use (L1) and (L2).\\

\noindent Question (c). Let $(f_n)_{n\geq 0}$ and $(g_n)_{n\geq 0}$ be two non-decreasing sequences of elements of $\subset \mathcal{F}^+$. Show that

$$
\lim_{n\rightarrow +\infty} f_n \leq \lim_{n\rightarrow +\infty} g_n \Rightarrow \lim_{n\rightarrow +\infty} L(f_n) \leq \lim_{n\rightarrow +\infty} (g_n).
$$

\bigskip \noindent \textit{Hints}. Suppose that $\lim_{n\rightarrow +\infty} f_n \leq \lim_{n\rightarrow +\infty} g_n $ and apply Question (a) to the relation : 

$$
\forall n\geq 0, \ f_n=\lim_{k\rightarrow +\infty} min(f_n,g_k). \ (CP)
$$

\bigskip \noindent Question (d) Use Question (c) to justify the extension $L^{*}$ of $L$ on $\mathcal{L}^+$ by taking $L^{*}(f)$ as the limit of any sequence $L(f_n)_{n\geq 0}$ where  
$(f_n)_{n\geq 0} \subset \mathcal{F}^+$ is such that $f_n\nearrow f$  as $n \nearrow +\infty$, that is :\\

$$
L^{\ast}(f)=\lim_{n\rightarrow +\infty} L(f_n), \ (f_n)_{n\geq 0} \subset \mathcal{F}^+,  \ f_n\nearrow f \ as \ n \nearrow +\infty.
$$

\bigskip \noindent Question (e) Establish the following simple properties of $L$ of $\mathcal{L}^+$ :\\

\noindent (e1) $L$ is the restriction of $L^{\ast}$ on $\mathcal{F}_b^+$, the class of bounded and non-negative elements of $\mathcal{F}$. \textbf{In the sequel, $L^{\ast}$ is simply denoted by $L$}.\\

\noindent (e2) $L$ is non-decreasing on $\mathcal{L}^+$, that is $L(f)\leq L(g)$ for any couple of elements of $\mathcal{L}^+$ such that $f\leq g$ (use again Formula (CP) above).\\

\noindent (e3) $\mathcal{L}^+$ is lattice and is closed under finite sum of functions and under multiplication by non-negative scalars, that is : for any couple $(f,g)$ of elements of $\mathcal{L}^+$ and for any non-negative real scalars $c\geq 0$ and $d\geq 0$, we have

$$
L(cf+dg)=c L(f)+ d L(g).
$$

\bigskip \noindent (e4) Let $(f,g)$ be any any couple of elements of $\mathcal{L}^+$, then

$$
L(f)+L(g)=L(\max(f,g))+L(\min(f,g)),
$$

\bigskip \noindent from the identity $f+g=\max(f,g)+\min(f,g)$.\\

\bigskip \noindent Question (f) (A form the Monotone Convergence Theorem). Let $(f_n)_{n\geq 0}$ be a non-decreasing sequences of elements of $\subset \mathcal{L}^+$ which is uniformly bounded, that is, there exists $M\in \mathbb{R}$ such that

$$
\forall n\geq 0, \ for all \omega \in \Omega, \ 0\leq f(\omega) \leq M.
$$
 
\bigskip \noindent Show that

$$
\lim_{n\rightarrow +\infty} L(f_n)= L(\lim_{n\rightarrow +\infty} f_n).
$$

\bigskip \noindent \textit{Suggestions}. First, remark that $f=\lim_{n\rightarrow +\infty} f_n \in \mathcal{L}^+$ which justifies the writing of $L(\lim_{n\rightarrow +\infty} f_n)$. From there, the method used in 
Exercise 1 in Doc 06-09 of Chapter \ref{06_convergence} (see page \pageref{exercise01_sol_doc06-09}) may be re-conducted word by word. If you have done that exercise, you will gain nothing by doing it again. If not, simply re-conduct that reasoning. \\

\bigskip \noindent \textbf{SOLUTIONS}.\\

\noindent Question (a). For any couple of elements of $\mathcal{F}$ such that $f\leq g$, we have $\mathcal{F} \ni g-f\geq 0$. By (L1) we have $0\leq L(g-f)=L(g)-L(f)$, that is $L(f)\leq L(g)$.\\

\noindent Question (b). Set $g_n=f-f_n$, $n\geq 0$. By the assumption $(g_n)_{n\geq 0}\subset \mathcal{F}^+$ and $g_n\downarrow 0$ as $n\uparrow +\infty$ and by (L3), we have $L(f)-L(f_n)\nearrow 0$,
that is  $L(f_n)\uparrow L(f)$.\\

\noindent Question (c). Suppose that $\lim_{n\rightarrow +\infty} f_n \leq \lim_{n\rightarrow +\infty} g_n $. We have 
$$
\forall n\geq 0, \ f_n=\lim_{k\rightarrow +\infty} \min(f_n,g_k),
$$

\bigskip \noindent By (L0) and (L1),  we have $\min(f_n,g_k) \in \mathcal{F}^+$ and $\min(f_n,g_k)\leq g_k$ for all $n\geq 0$ and $k \geq 0$. Thus by Question (c),

$$
L(f_n)=\lim_{k\rightarrow +\infty} L(\min(f_n,g_k)) \leq \lim_{k\rightarrow +\infty} L(g_k).
$$

\bigskip 
\noindent We conclude by letting $n\rightarrow +\infty$.\\

\noindent Question (d) Let $(f_n)_{n\geq 0}$ and $(g_n)_{n\geq 0}$ be two non-decreasing sequences of elements of $\subset \mathcal{F}^+$ converging to $f \in \mathcal{L}^+$. We have

$$
\lim_{n\rightarrow +\infty} f_n = \lim_{n\rightarrow +\infty} g_n.
$$

\bigskip 
\noindent Thus by Question (c), we have both
$$
\lim_{n\rightarrow +\infty} L(f_n) \leq \lim_{n\rightarrow +\infty} (g_n) \text{  and } \lim_{n\rightarrow +\infty} L(g_n) \leq \lim_{n\rightarrow +\infty} (f_n).
$$

\bigskip 
\noindent Hence $\lim_{n\rightarrow +\infty} L(f_n) =\lim_{n\rightarrow +\infty} (g_n)$. This show that the functional $L^{*}$ is well-defined on $\mathcal{L}^+$, since the value of $L^{*}(f)$ does not depend of the sequence that is used to compute it.\\

\noindent Question (e) These formulas are obvious and do not need solutions. But do them yourself.\\

\noindent Question (f). The suggestions given after the question are enough.\\

\bigskip \noindent \textbf{Exercise 2}. \label{exercise02_sol_doc09-05}\\

\noindent Define 

$$
\mathcal{C}=\{A \subset \Omega \ : \ 1_A \in , \mathcal{L}^+\}.
$$

\bigskip \noindent Question (a) Show that $\{\emptyset, \Omega\}\subset \mathcal{C}$. Show that $\mathcal{C}$ is closed under finite intersections (that is, it is $\pi$-system) and countable unions.\\

\noindent \textit{Hints} To show that $\mathcal{C}$ is closed under finite intersections and unions, take only intersections or unions of two of its elements, which is enough thanks to the associativity of such operations. To show that To show that $\mathcal{C}$ is closed under countable unions, use partial unions and use Question (f) of Exercise 1.\\

\noindent Question (b) Define 

$$
\mathcal{C} \ni A \mapsto m^{\ast}(A)=L(1_A).
$$

\bigskip \noindent Show that $m^{\ast}$ is proper, non-negative, additive and continuous below, normed ($\mu^{\ast}(\Omega)=1$) and satisfies for all $G_i\in \mathcal{C}$, $i\in\{1,2\}$ 

$$
m^{\ast}(G_1\cap G_2) + m^{\ast}(G_1\cap G_2)=m^{\ast}(G_1)+m^{\ast}(G_2).
$$

\bigskip \noindent \textit{Hints}. All these facts are direct consequences of Exercise 1. To justify this, just indicate the exact formulas from which they derive.\\

\noindent Question (b) Show that the application

$$
\mathcal{P}(\Omega) \ni B \mapsto \mu(A)=\inf\{\mu^{\ast}(B), \ B \in \mathcal{C}, A \subset B\}
$$

\bigskip \noindent is an extension of $m$ and is an outer measure.\\

\noindent Question (c) By Part VII in Doc 04-01 (and its proof in Doc 04-10), $m$ is a measure on the $\sigma$-algebra 

$$
\mathcal{A}_0=\{A \in \mathcal{P}(\Omega) : \ \forall D\subset \Omega, \ \mu(D)=\mu(AD)+\mu(A^{c}D) \},
$$ 

\bigskip \noindent Tell why $m$ is the unique extension of $m^{\ast}$ to a measure, that $m$ is a probability measure. Use the additivity of $m^{\ast}$ of $\mathcal{C}$ to see that 
$\mathcal{C} \subset \mathcal{A}_0$.\\

\noindent Question (d) \noindent Let $\mathcal{A}=\mathcal{A}_{\mathcal{F}}$ the $\sigma$-algebra generated by $\mathcal{F}$, that is

$$
\mathcal{A}=\sigma(\{\ f^{-1}(B), \ f \in \mathcal{F}, \ B \in \mathcal{B}_{\infty}(\overline{\mathbb{R}}) \}).
$$

\bigskip \noindent Show that $\mathcal{A} \subset \mathcal{A}_0$.\\

\noindent \textit{Hints}.\\

\noindent For $f\in \mathcal{F}$, for any $c>0$, justify

$$
1_{(f>c)}=\lim_{n\rightarrow +\infty}  \min(1, \ n(f-c)^+). \ (AP01)
$$

\bigskip \noindent One one hand, deduce that $\mathcal{A}_{\mathcal{L}^+} \subset \mathcal{A}_{\mathcal{F}}$. Also easily check that, for $\mathcal{L}_0=\{1_A, \ A\in \mathcal{C} \}$, we have 
$\mathcal{A}_{\mathcal{L}_0} \subset \mathcal{A}_{\mathcal{L}^+}$.\\

\noindent On another hand show that $\mathcal{A}_{\mathcal{L}^+} \subset \mathcal{A}_{\mathcal{F}}$ from the definition $\mathcal{L}^+$. Conclude that $\mathcal{A}_{\mathcal{L}^+}=\mathcal{A}_{\mathcal{F}}$.\\

\noindent Finally, use the additivity of $m^{\ast}$ to show that $\mathcal{C}\in \mathcal{A}_0$. Make your final conclusion.\\

\noindent Question (e). Let $f\in \mathcal{F}$.\\

\noindent (e1) Check that the formula

$$
L(f) = \int f \ dm, \ (AP)
$$

\bigskip \noindent holds of $f\in \mathcal{L}_0$, that is $f$ is of the form $f=1_A$, $A \in \mathcal{A}$. Extend the validity of (AP) on the class of finite linear combinations of elements of $f\in \mathcal{L}_0$.\\

\noindent (e2) Let $f\in \mathcal{F}_b^+$, bounded by $M\geq 0$. Use Exercise 3 in Doc 03-08, page \pageref{exercise03_sol_doc03-08}, show how to drop the part relative to the segment 
$[2^n, +\infty[$, exchange the strict and broad inequalities to see that the sequence

\begin{eqnarray*}
f_n&=&\sum_{k=0}^{{k_M}-1} \frac{k-1}{2^n} \biggr(1_{(f>(k-1)/2^n)} - 1_{(f>k/2^n)} \biggr) + \biggr(1_{(f>(k_M-1)/2^n)} - 1_{(f>M)}.
\end{eqnarray*} 

\noindent broadly increases to $f$. Next use the remarkable identity

$$
\sum_{i=1}^{p} \alpha_{i-1} (\beta_{i-1}-\beta_i)=\alpha_0 \beta_0 + \biggr( \sum_{i=1}^{p-1} (\alpha_i - \alpha_{i-1} \biggr) +\alpha_{p-1} \beta_p,
$$

\bigskip \noindent where with $(\alpha_i)_{0\leq i \leq p_n}$ and $(\beta_i)_{0\leq i \leq p_n}$, $p\geq 1$, are finite real numbers and extend the formula (AP) to $f\in \mathcal{F}_b^+$. Be rigorous in justufying the limits.\\

\noindent (e3) For $f\in \mathcal{F}^+$, in general, use the approximation

$$
f=\lim_{p\rightarrow+\infty} \max(f,n)
$$

\bigskip \noindent and extend (AP) for $f\in \mathcal{F}^+$.  Give all the details.\\

\noindent (e4) For your general conclusion, use $f=f^+-f^-$.\\

\noindent (e5) Do we have $\mathcal{F} \subset \mathcal{L}^1(\omega, \mathcal{A},m)$?\\

\bigskip \noindent \textbf{SOLUTIONS}.\\

\noindent Questions (a), (b), (c) : All these facts derive from properties of indicator functions and properties established in Exercise 1.\\

\noindent Question (d).\\

\noindent (d1) To show that Formula (AP01) holds for $f\in \mathcal{F}$, Is is enough to check it for $\omega \in (f>c)$ and $\omega \in (f\leq c)$. Remark that all the functions used here are non-negative and bounded. By the properties of $\mathcal{F}$ and $\mathcal{L}^+$ established in Exercise 1, we get that the indicator function $f(>c)$ is in $\mathcal{L}^+$ and thus $(f>c)$ lies in $\mathcal{C}$. We get that each $f$ in $\mathcal{F}$ is $\sigma(\mathcal{C})$-measurable. By definition, we have

\begin{eqnarray*}
\mathcal{A}_{\mathcal{F}} &=& \sigma(\{\ f^{-1}(B), \ f \in \mathcal{F}, \ B \in \mathcal{B}_{\infty}(\overline{\mathbb{R}}) \})\\
&\subset &  \sigma(\mathcal{C}).
\end{eqnarray*}

\noindent By denoting $\mathcal{L}_0=\{1_A, \ A\in \mathcal{C} \}$, it is clear that $\mathcal{A}\subset \mathcal{A}_{\mathcal{L}_0} \subset \mathcal{A}_{\mathcal{L}^+}$. Also we have
$\mathcal{A}_{\mathcal{L}^+} \subset \mathcal{A}_{\mathcal{F}}$, since the first class is generated by elements of the form $f^{-1}(B)$, $f \in \mathcal{F}$, $B \in \mathcal{B}_{\infty}$ and each of those sets $f^{-1}(B)$ satisfies

$$
f^{-1}(B)= \bigcup_{n\geq 0} f_n^{-1}(B),  \ (f_n)_{n\geq 0} \subset \mathcal{F}_b^+\subset \mathcal{L}^+.
$$

\bigskip 
\noindent By combining these facts, we get
$$
\mathcal{A}=\mathcal{A}_{\mathcal{F}}=\sigma(\mathcal{C}) \subset \mathcal{A}_0.
$$

\bigskip 
\noindent Hence $m$ is a probability measure on $(\Omega, \mathcal{A})$.

\noindent Question (e). Formula (SF) holds for $f=1_A$, $A \in \mathcal{C}$, that is $f\in \mathcal{L}_0$ and for $f$ beeing any finite linear combinations of elements of $\mathcal{L}_0$. Now for any function $f\geq 0$, we already have that $f$ is a non-decreasing limit of the sequence

$$
f_n=\sum_{k=1}^{2^{(2n)}} \frac{k-1}{2^n} 1_{(k-1)/2^n <f\leq k/2^n)} + 2^n 1_{(f>2^n)}.
$$
 
\bigskip 
\noindent (See Exercise 3 in Doc 03-08, page \pageref{exercise03_sol_doc03-08}, with a slight modification consisting in interchanging the strict and broad inequalities). Let $f\in \mathcal{F}_b^+$, where $M$ is the bound of $f$. If $n$ is large enough, we will have $M\leq 2^n$ and we can find a value $K_M$ such that ${(k_M-1)}/2^n <M\leq {k_M}/2^n)$ and $f_n$ becomes

\begin{eqnarray*}
f_n&=&\sum_{k=1}^{{k_M}-1} \frac{k-1}{2^n} 1_{(k-1)/2^n <f\leq k/2^n)} + \frac{{k_M}-1}{2^n} 1_{(k_M-1)/2^n <f\leq M}\\
&=&\sum_{k=0}^{{k_M}-1} \frac{k-1}{2^n} \biggr(1_{(f>(k-1)/2^n)} - 1_{(f>k/2^n)} \biggr) + \biggr(1_{(f>(k_M-1)/2^n)} - 1_{(f>M)}.
\end{eqnarray*}

\bigskip 
\noindent By the remarkable identity, 

$$
\sum_{i=1}^{p} \alpha_{i-1} (\beta_{i-1}-\beta_i)=\alpha_0 \beta_0 + \biggr( \sum_{i=1}^{p-1} (\alpha_i - \alpha_{i-1} \biggr) +\alpha_{p-1} \beta_p,
$$

\bigskip 
\noindent where with $(\alpha_i)_{0\leq i \leq p_n}$ and $(\beta_i)_{0\leq i \leq p_n}$, $p\geq 1$, are finite real numbers, we see that each $f_n$, for large enough values of $n\geq 0$, is if the form,

$$
f_n=\sum_{i=1}^{p_n} a_i 1_{(f>c_i)}, \ p_n\geq 1, \ c_i \in \mathbb{R}.
$$

\bigskip 
\noindent Hence for each $n$ large enough, is a finite linear combination of indicator functions of elements of $\mathcal{C}$. Hence, for $n \geq (\log M)/2$, we have
$$
L(f_n)=\int f_n \ dm. \ (AP02)
$$

\bigskip 
\noindent Since $(f_n)_{n\geq 0}$ converges to $f\in \mathcal{L}^+$, we may apply the MCT (for the measure $m$) and Question (a) of Exercise 1 (for the functional $L$) to get Formula (AP) by letting $n\rightarrow +\infty$.\\

\noindent Next for $f \in \mathcal{F}^+$, we have

$$
f=\lim_{p\rightarrow+\infty} \max(f,n)
$$

\bigskip 
\noindent and apply Question (a) of Exercise 1 and the MCT to get Formula (AP) again.\\

\noindent Finally for any $f\in \mathcal{F}$, we use the decomposition $f=f^+-f^-$ to conclude.

\chapter{Lebesgue-Stieljes measures on $\mathbb{R}^4$} \label{10_lsm}

\noindent \textbf{Content of the Chapter}

\begin{table}[htbp]
	\centering
		\begin{tabular}{llll}
		\hline
		Type& Name & Title  & page\\
		A& Doc 10-01 &  Lebesgue-Stieljes measure and Probability Theory  & \pageref{doc10-01}\\
		A& Doc 10-02 & Proof of the Existence of the Lebesgue-Stieljes measure on $\mathbb{R}^k$  & \pageref{doc10-02} \\
		A& Doc 10-03  & Lebesgue-Stieljes Integrals and Riemann-Stieljes Integrals & \pageref{doc10-03} \\
		%&  &   & \\
		%&  &   & \\
		%&  &   & \\
		\hline
		\end{tabular}
\end{table}

\newpage
\noindent \LARGE \textbf{DOC 10-01 : Lebesgue-Stieljes measure and Probability Theory} \label{doc10-01}\\
\bigskip

\Large
\noindent \textbf{(10.01) Introduction}. \label{doc10-01_intro}\\

\noindent A very large part of probability theory deals with real random variables, random vectors or stochastic processes with real states. Usually, formal
and general spaces are used in theory but almost always, applications come back to real data. The fundamental Theorem of Kolmogorov, as we will see it in the course of stochastic process, will establish that the probability law of an arbitrary family of real random variables is characterized by the laws of its finite subfamilies.\\

\noindent This means that we finally work with probability laws on $\mathbb{R}^{k}.$ The point is that probability laws in $\mathbb{R}^{k}$ are Lebesgue-Stieljes measures. This
allows us to say that mastering Lebesgue-Stieljes measures on \ $\mathbb{R}^{k}$ is the key that opens the castle of Probability Theory.\\

\noindent Of course the very foundation of probability theory is measure theory and its probability theory version. Normally, Lebesgue-Stieljes measures on \ $\mathbb{R}^{k}$
is part of the course of Measure theory. But the latter one is so heavy that there is no enough space to handle it properly. From a statistical point of view, the recent notion of copulas is now popular and studying some areas of Statistics without copulas is impossible. And it happens that copulas are special Lebesgue-Stieljes measures on $\mathbb{R}^{k}$.\\

\noindent Then it is of a great importance we provide a special cover of the Lebesgue-Stieljes measures on  $\mathbb{R}^{k}$ before we close this textbook.\\

\noindent \textbf{Warning}. We already have a first and complete handling in the univariate case in Doc 04-04, Doc 04-07 and Doc 04-10 (pages \pageref{doc04-04}, \pageref{doc04-07} and \pageref{doc04-10}). Although we will focus on the non-univariate case, we will represent, once again, the univariate case for the reader who wants to have all the materials in one document.\\

\bigskip \noindent \textbf{(10.02) Distribution functions and measures}. \label{doc10-02_intro}\\

\noindent We are going to discover the deep relations between distribution functions and measures on $\mathbb{R}^{k}$ assigning finite values to bounded above Borel sets. The general case of arbitrary dimension $k$ may make the reader feel that it is complicated. So we are going to guide him by beginning with the simplest case with $k=1$. Next, we treat the intermediate case $k=2$ and let him see principles that will guide the general case. While dealing with complex expressions in the general case, the meaning will be as simple as in the bivariate case.

\bigskip \noindent \textbf{(10.02a) Univariate case}. \label{doc10-02_intro}\\

\noindent  Let $m$ be a measure on the Borel set of $\mathbb{R}$\ such that for any  any $t\in \mathbb{R}$,\ 
\begin{equation*}
F_{m}(t)=m(]-\infty ,t])\text{ is finite.}  \ (RC01)
\end{equation*}

\noindent The function $F_{m}:\mathbb{R}\longmapsto $\ $\mathbb{R}$ defined by 
\begin{equation*}
t\mapsto F_{m}(t)  \ (FD01)
\end{equation*}

\noindent has the following properties.\\

\noindent \textbf{Claim 1}.  $F_{m}$\ assigns non-negative lengths to intervals, that is., 

$$
t\leq s\Longrightarrow \Delta F_{m}(]t,s])=F_{m}(s)-F_{m}(t)\geq 0. \ (PI01)
$$

\bigskip \noindent Here, we have to pay attention to the terminology. We say that $\Delta F_{m}(]t,s]=F_{m}(s)-F_{m}(t)$ is the length of the interval $]t,s]$ by $F_{m}.$ It is a pure coincidence that Formula (PI01) means that $F_{m}$ is non-decreasing in this particular case. In higher order spaces, this notion of non-decreasingness will disappear while the notion of positive lengths will be naturally extended to that positive areas and generally to positive volumes.\\

\noindent Now the proof of Formula (PI01) is immediately seen by remarking that
\begin{equation*}
(t\leq s)\Rightarrow (]-\infty ,t])\subseteq ]-\infty ,s])
\end{equation*}

\noindent and

$$
m_(]t,s])=F_m(s)-F_m(t)
$$

\bigskip \noindent and then Formula (PI01) holds based on the fact that $F_m(s)-F_m(t)$ is the measure of an interval. We say that $F_m(s)-F_m(t)$ is the length of $]t,s]$ by $F_m$ or that $F_m$ assigns to the interval $]t,s]$ the length $F_m(s)-F_m(t)$. With such a terminology, we are ready to move to area in dimension two and to volumes in higher dimensions.\\

\noindent \textbf{Claim 2}. $F_{m}$\ is right-continuous, that is,  for all $t\in \mathbb{R}$,
\begin{equation*}
F_{m}(t^{(n)})\downarrow F_{m}(t)
\end{equation*}

\noindent as $t^{(n)})\downarrow t$.\\

\bigskip \noindent  \textbf{Proof}.\\

\bigskip \noindent Let $t^{(n)})\downarrow t.$ We get\ 
\begin{equation*}
]-\infty ,t^{(n)})]=\text{ }\downarrow \text{ }]-\infty ,t].
\end{equation*}

\noindent Since $m(]-\infty ,t^{(n)})])$\ is finite for all $n$'s, we get by continuity of the measure $m,$ that 
\begin{equation*}
F_{m}(t^{n})=m(]-\infty ,t^{n}])\downarrow m(]-\infty ,t])=F_{m}(t).
\end{equation*}

\bigskip \noindent Suppose now that $m$ is a probability measure, that is $m(\mathbb{R})=1$, we
have two additional points.\\

\noindent \textbf{Claim 3}.\\

$F_{m}(t)\downarrow 0$ as $t\downarrow -\infty $ \ and $F_{m}(t)\uparrow 1$ as $t\uparrow +\infty$.\\

\bigskip \noindent These two points result from the continuity of the measure $m$. First 

\begin{equation*}
]-\infty ,t] \downarrow \emptyset \text{ as } t\downarrow -\infty \text{ and } m(]-\infty ,t])<+\infty \text{ for all } t\in \mathbb{R}
\end{equation*}

\bigskip \noindent implies

\begin{equation*}
m(]-\infty ,t]) \downarrow m(\emptyset)=0 \ \  as \ \ t \downarrow -\infty.
\end{equation*}

\noindent Next 

\begin{equation*}
]-\infty ,t]\uparrow \mathbb{R}\text{ as }t\uparrow +\infty 
\end{equation*}

\noindent implies

\begin{equation*}
m(]-\infty ,t])\uparrow m(\mathbb{R})=0\text{ as }t\uparrow -\infty.
\end{equation*}

\bigskip \noindent  In the sequel, we make these notation : $F(-\infty)=\lim_{t\downarrow -\infty }F(t)$ \ and $F(+\infty )=\lim_{t\uparrow+\infty }F(t)$. We summarize our study through theses two definitions.\\

\noindent \textbf{Definition}. A non-constant function $F:\mathbb{R}$ $\longmapsto$ $\mathbb{R}$ is a distribution function if and only if \\

\noindent (1) it assigns non-negative lengths to intervals

\noindent and \\

\noindent (2) it is right-continuous at any point $t\in \mathbb{R}$.

\bigskip \noindent This definition is very broad. $F$ is not necessarily non-negative. It is not required that $F(-\infty )=0.$ The second definition is more restrictive.\\

\noindent \textbf{Definition}. A \textbf{non-negative} function $F:\mathbb{R}\longmapsto $\ $\mathbb{R}$ is a probability
distribution function, or a cumulative distribution function (\textit{cdf}), if and only if the following assertions hold :\\

\noindent (1) it assigns non-negative lengths to intervals.\\

\noindent (2) it is right-continuous.\\

\noindent (3) $F(-\infty )=0$ \ and $F(+\infty )=1$.\\

\bigskip \noindent  We conclude by saying that a measure $m$ which assigns finite values to bounded above intervals generates the distribution functions $F_{m}$. A probability measure generates the probability distribution functions $F_{m}$.\\

\bigskip \noindent Conversely, let be given a  distribution function $F$. From the distribution function $F+c$, \ it is possible to create a measure on the semi-algebra 
\begin{equation*}
\mathcal{I}_{1}=\{]a,b], \ -\infty \leq a\leq b \leq +\infty \}
\end{equation*}

\bigskip \noindent that generates the usual $\sigma$-algebra $\mathcal{B}(\mathbb{R})$ by

\begin{equation}
m_{F}(]a,b])=\Delta F(]a,b])=F(b)-F(a).  \ (LS01)
\end{equation}

\bigskip \noindent Let us show that $m_{F}$ is additive on $\mathcal{I}_{1}$. The only possibility to have that an element $\mathcal{I}_{1}$ of is sum of two elements of is to split on element of into two others, that is

\begin{equation*}
]a,b]=]a,c]+]c,b]\text{, where }c\text{ is finite and }a<c<b
\end{equation*}

\noindent \noindent \bigskip and then

\begin{eqnarray*}
m_{F}(]a,c])+m_{F}(]c,b])&=&F(c)-F(a)+F(b)-F(c)\\
&=&F(b)-F(a)=m_{F}(]a,b]).
\end{eqnarray*}

\bigskip \noindent Extension of $m_{F}$ to an additive application on the algebra $\mathcal{C}$ generated by $\mathcal{I}_{1}$ which consists in the collection of finite sums of elements of $\mathcal{I}_{1}$, still denoted as $m_{F}$, is straightforward by simple arguments established in the course of Measure Theory. $m_F$ is also $\sigma$-finite since

$$
\mathbb{R}= \bigcup_{n\geq 1} ]-\infty, n] \ and, \ \forall n\geq 1, m_F(]-\infty, n])<+\infty.
$$

\noindent It remains to show that $m_{F}$ is sigma-sub-additive on $\mathcal{C}$. From there, the Caratheodory's Theorem ensures that $m_{F}$ is uniquely extensible to a measure, still denoted as 
$m_{F},$ on the usual sigma-algebra generated both by $\mathcal{I}_{1}$ and by $\mathcal{C}$.

\bigskip \noindent The distribution function $F_{m_{F}}$\ generated by \ $m_{F}$ is obtained by letting $a$ go to $-\infty$ in (\ref{def_m1}) that gives

\begin{equation*}
F_{m_{F}}(b)=m_{F}(]-\infty ,b])=F(b)-F(-\infty),
\end{equation*}

\noindent \noindent which is exactly $F$ if $F(-\infty )=0$. And, we have the following conclusion.\\

\noindent Let $\mathcal{M}_{0}$ the class of measures $m$ on $\mathbb{R}$, such that%

\begin{equation*}
\forall u\in \mathbb{R},m(]-\infty ,u])<+\infty 
\end{equation*}

\noindent and $\mathcal{F}_{0}$ the class of distribution functions $F$ on $\mathbb{R}$ such that
 
\begin{equation*}
F(-\infty)=0.
\end{equation*}

\bigskip \noindent There is a one-to-one mapping between $\mathcal{M}_{0}$ and $\mathcal{F}_{0}$ in the following way 
\begin{equation*}
\begin{tabular}{lll}
$F(u)=m(]-\infty ,u])$ & $\longleftrightarrow $ & $m(]a,b])=\Delta F(a,b)$.
\end{tabular}%
\end{equation*}

\bigskip \noindent In particular, there is a one-to-one mapping between the class $\mathcal{P}$ of probability measures on $\mathbb{R}$ and the class $\mathcal{F}_{P}$ of probability distribution functions on $\mathbb{R}$. Probability measures are characterized by their distribution functions.\\

\bigskip \noindent Before we consider the general case, let us do, in details, the case $k=2$ as an intermediate one.\\

\bigskip \noindent \textbf{(10.02b) Bivariate case}. \label{doc10-02b}\\

\noindent The usual $\sigma$-algebra on $\mathbb{R}^{2}$, the Borel sigma-algebra $\mathcal{B}(\mathbb{R}^{2}$), is generated by the semi-algebra%
\begin{equation*}
\mathcal{I}_{2}=\{]a,b],a\leq b\, \ (a,b)\in \mathbb{\overline{R}}^2 \},
\end{equation*}%

\noindent where $]a,b]=]a_{1},b_{1}]\times ]a_{2},b_{2}]$ and $a=(a_{1},a_{2})\leq b=(b_{1},b_{2})$ means : $\ a_{1}\leq a_{2}$ and $b_{1}\leq b_{2}$.\\ 

\noindent We begin to express $m(]a,b])$ by means of $F_{m}(u)=m(]-\infty,u])$ with $u=(u_1,u_2)$ and $]-\infty,u]=]-\infty,u_{1}] \times ]-\infty, u_2]$. We have

\begin{eqnarray*}
]a,b]&=&]a_{1},b_{1}]\times ]a_{2},b_{2}]=(]-\infty ,b_{1}]\setminus ]-\infty ,a_{1}])\times ]a_{2},b_{2}]\\
&=&(]-\infty ,b_{1}]\times ]a_{2},b_{2}])\setminus (]-\infty ,a_{1}]\times ]a_{2},b_{2}]).
\end{eqnarray*}

\bigskip \noindent Since the right-side term of the sets minus operator($\backslash$) is contained in the left-side term, we have
\begin{equation*}
m(]a,b])=m(]-\infty ,b_{1}]\times ]a_{2},b_{2}])-m(]-\infty ,a_{1}]\times ]a_{2},b_{2}]). \ (L1)
\end{equation*}

\bigskip \noindent We do the same thing with $A=m(]-\infty ,b_{1}]\times ]a_{2},b_{2}])$ by using
  
\begin{equation*}
]-\infty ,b_{1}]\times ]a_{2},b_{2}]=]-\infty ,b_{1}]\times (]-\infty
,b_{2}]\setminus ]-\infty ,a_{2}])
\end{equation*}

\bigskip \noindent to get
\begin{equation*}
]-\infty ,b_{1}]\times ]a_{2},b_{2}]=(]-\infty ,b_{1}]\times ]-\infty
,b_{2}])\setminus (]-\infty ,b_{1}]\times ]-\infty ,a_{2}]).
\end{equation*}

\bigskip \noindent Thus, we have  
\begin{eqnarray*}
A&=&m(]-\infty ,b_{1}]\times ]a_{2},b_{2}])=m(]-\infty ,b_{1}]\times ]-\infty,b_{2}])\\
&-&m(]-\infty ,b_{1}]\times ]-\infty ,a_{2}]) \\
&=&F_m(b_1,b_2)-F_m(b_1,a_2). \ (L2)
\end{eqnarray*}

\bigskip \noindent Let us repeat this for $B=m(]-\infty ,a_{1}]\times ]a_{2},b_{2}])$\ by using 
\begin{equation*}
]-\infty ,a_{1}]\times ]a_{2},b_{2}]=]-\infty ,a_{1}]\times (]-\infty
,b_{2}]\setminus ]-\infty ,a_{2}])
\end{equation*}

\bigskip \noindent and

\begin{equation*}
]-\infty ,a_{1}]\times ]a_{2},b_{2}]=(]-\infty ,a_{1}]\times ]-\infty
,b_{2}])\setminus (]-\infty ,a_{1}]\times ]-\infty ,a_{2}]).
\end{equation*}

\bigskip \noindent we get
\begin{eqnarray*}
&&B=m(]-\infty ,a_{1}]\times ]a_{2},b_{2}])  \\
&=&m(]-\infty ,a_{1}]\times ]-\infty,b_{2}])-m(]-\infty ,a_{1}]\times ]-\infty ,a_{2}])\\
&=&F_m(a_1,b_2)-F_m(a_1,a_2). \ (L3)
\end{eqnarray*}

\bigskip \noindent By putting together (L1), (L20) and (L3), we arrive at (L4) :

\begin{equation*}
m(]a,b])=F_{m}(b_{1},b_{2})-F_{m}(b_{1},a_{2})-F_{m}(a_{1},b_{2})+F_{m}(a_{1},a_{2}).
\label{lj04}
\end{equation*}

\bigskip \noindent We define
\begin{equation*}
\Delta
F_{m}(a,b)=F_{m}(b_{1},b_{2})-F_{m}(b_{1},a_{2})-F_{m}(a_{1},b_{2})+F_{m}(a_{1},a_{2})
\end{equation*}%

\noindent as the area of the rectangle $]a,b].$ To see why this is called an area, consider the case where

\begin{equation*}
F(u,v)=F_{1}(u)F_{2}(v),
\end{equation*}

\bigskip \noindent where each $F_{i}$ is a distribution function on $\mathbb{R}$. We have
\begin{equation*}
\Delta F(a,b)=(F_{1}(a_{2})-F_{1}(a_{1}))((F_{2}(b_{2})-F_{2}(b_{1})).
\end{equation*}

\bigskip \noindent For $F_{1}(u)=F_{2}(u)=u$, the identity function, we have the classical area of the rectangle
\begin{equation*}
\Delta F(a,b)=(a_{2}-a_{1})\times (b_{2}-b_{1}).
\end{equation*}

\bigskip \noindent From there, we may re-conduct the frame used for $\mathbb{R}$. Let $m$ be a measure on the Borel sets of $\mathbb{R}^{2}$\ such that for any $t\in \mathbb{R}^{2},$\ 
\begin{equation*}
F_{m}(t)=m(]-\infty ,t])\text{ is finite.} \ (RC02)
\end{equation*}

\bigskip \noindent The function $F_{m}:\mathbb{R}^{2}$ $\longmapsto$ $\mathbb{R}$ defined by 
\begin{equation*}
t\mapsto F_{m}(t) \ (FD02)
\end{equation*}

\bigskip \noindent has the following properties.\\

\noindent \textbf{Claim 4}. $F_{m}$\ assigns non-negative areas to rectangles, that is, 
\begin{equation*}
a=(a_{1},a_{2})\leq b=(b_{1},b_{2})\Longrightarrow \Delta F_{m}(a,b)\geq 0.
\end{equation*}

\bigskip \noindent \textbf{Warning!} Please, remark that $F_{m}$ is already component-wise non-decreasing and it is easy to see it by the same reasoning used in $\mathbb{R}$. Yet, this is not enough to get a distribution function! We go further by requiring that $F_{m}$ assigns non-negative areas to rectangles. When we go back to $\mathbb{R}$, rectangles become intervals and areas are lengths. It happens that assigning non-negative lengths is equivalent to be non-decreasing in $\mathbb{R}$. For higher dimensions, we will use cuboids and non-negative volumes.\\

\noindent \textbf{Claim 5}. $F_{m}$\ is right-continuous, that is, 
\begin{equation*}
F_{m}(t^{(n)})\downarrow F_{m}(t)
\end{equation*}

\bigskip \noindent as $t^{(n)}=(t_{1}^{(n)},t_{2}^{(n)})\downarrow t=(t_{1},t_{2})$ where 

\begin{equation*}
(t^{(n)}\downarrow t)\Leftrightarrow (\forall (i=1,2), \  t_{i}^{(n)}\downarrow t_{i}).
\end{equation*}

\bigskip \noindent If $m$ is a probability measure, that is $m(\mathbb{R}^{2})=1$,
we add these properties : \newline

\noindent \textbf{(3a)} $lim_{s\rightarrow -\infty \text{ \textbf{or} }t\rightarrow
-\infty }F_{m}(s,t)=0.$\newline

\noindent and \newline

\noindent \textbf{(3b)} $lim_{s\rightarrow +\infty \text{ \textbf{and} }t\rightarrow +\infty }F_{m}(t)=0$.\\

\noindent \textbf{Warning!} Please,pay attention to the connecting particles \textbf{and} and \textbf{or} in Formulas \textbf{(3a)} and \textbf{(3b)} . This leads us the same definitions.\\

\noindent \textbf{Definition} A non-constant function $F:\mathbb{R}^{2}\longmapsto $\ $\mathbb{R}$ is a distribution function if and only if\\

\noindent \textbf{(1)} it assigns non-negative areas to rectangles.\\

\noindent \textbf{(2)} it is right-continuous\\

\bigskip \noindent \textbf{Definition}. A \textbf{non-negative} function $F:\mathbb{R}^{2}\longmapsto $\ $\mathbb{R}$ is a probability distribution function, or a cumulative distribution function, if and only if \newline

\noindent \textbf{(1)} assigns non-negative areas to rectangles.\\

\noindent \textbf{(2)} it is right-continuous.\\

\noindent \textbf{(3)} We have
$$
\lim_{s\rightarrow -\infty \text{ or }t\rightarrow -\infty }\ F_{m}(s,t)=0
$$

\bigskip \noindent and

$$
\lim_{s\rightarrow +\infty \text{ and }t\rightarrow +\infty }F_{m}(t)=1.
$$

\bigskip \noindent Conversely we proceed similarly as follows.\\

\noindent Let $F$ be a distribution function on $\mathbb{R}^{2}$. Define the application $m_{F}$ on the semi-algebra $\mathcal{I}_{2}$

\begin{equation*}
\mathcal{I}_{2}=\{]a,b],a\leq b, \ \ (a,b)\in \mathbb{\overline{R}}^2 \}
\end{equation*}

\bigskip \noindent that generates the usual the sigma-algebra of $\mathcal{B}(\mathbb{R}^{2})$ by

\begin{equation*}
m_{F}(]a,b])=\Delta F(a,b). \ (LS02)
\end{equation*}

\bigskip \noindent Let us show that $m_{F}$ is additive on $\mathcal{I}_{2}$. The only possibility to have that an element of $\mathcal{I}_{2}$ is sum of two elements of $\mathcal{I}_{2}$
is to split on element of $\mathcal{I}_{2}$ into two others, by splitting one of its two factors, that is 
\begin{equation*}
]a,b]=]a_{1},b_{1}]\times ]a_{2},b_{2}]=(]a_{1},c]+]c,b_{1}])\times ]a_{2},b_{2}],
\end{equation*}

\bigskip \noindent where $c$ is finite and $a_{1}<c<b_{1}$,that is 
\begin{equation*}
A=]a,b]=]a_{1},c]\times ]a_{2},b_{2}]+]c,b_{1}]\times ]a_{2},b_{2}]=:B+C.
\end{equation*}

\bigskip \noindent Check that
\begin{eqnarray*}
m_{F}(B)+m_{F}(C)&=&F(c,b_{2})-F(a_{1},b_{2})-F(c,a_{2})+F(a_{1},a_{2})\\
&+&F(b_{1},b_{2})-F(c,b_{2})-F(b_{1},a_{2})+F(c,a_{2})\\
&=&F(b_{1},b_{2})-F(b_{1},a_{2})-F(a_{1},b_{2})+F(a_{1},a_{2})\\
&=&m_{F}(A).
\end{eqnarray*}

\bigskip \noindent Since $m_{F}$ is additive and non-negative, its extension of $m_{F}$ to an additive application on the algebra $\mathcal{C}$ generated by $\mathcal{I}%
_{1}$ which is the collection of finite sums of elements of $\mathcal{I}_{2}$, still denoted as $m_{F}$, is straightforward by a simple arguments established in the course of Measure Theory. $m_F$ is also $\sigma$-finite since

$$
\mathbb{R}= \bigcup_{n\geq 1} ]-\infty, (n,n)] \ and, \ \forall n\geq 1, m_F(]-\infty, (n,n)])<+\infty.
$$

\noindent It remains to show that $m_{F}$ is sigma-sub-additive on $\mathcal{C}$. From there, the Caratheodory's Theorem ensures that $m_{F}$ is uniquely extensible to a measure, still denoted
as $m_{F},$ on the usual sigma-algebra generated both by $\mathcal{I}_{1}$ and by $\mathcal{C}$.\\

\bigskip \noindent The distribution function $F_{m_{F}}$\ generated by \ $m_{F}$ is obtained by letting $a$ go to $-\infty $ in (\ref{def_m1}) that gives
\begin{eqnarray*}
F_{m_{F}}(b)&=&m_{F}(]-\infty ,b])\\
&=&F(b)+F(-\infty ,b_2)+F(b_1 ,-\infty)-F(-\infty ,-\infty).
\end{eqnarray*}

\bigskip \noindent which is exactly $F$ if 

$$
\lim_{s\rightarrow -\infty \text{ or } t\rightarrow -\infty }\ F(s,t)=0.
$$

\bigskip \noindent And, we have the following conclusion. Let $\mathcal{M}_{0}$ the class of measures $m$ on $\mathbb{R}^{2}$, such that
\begin{equation*}
\forall u\in \mathbb{R}^{2},m(]-\infty ,u])<+\infty 
\end{equation*}

\bigskip \noindent and $\mathcal{F}_{0}$ the class of distribution functions $F$ on $\mathbb{R}^{2}$ such that 
$$
\lim_{s\rightarrow -\infty \text{ or } t\rightarrow -\infty }\ F(s,t)=0.
$$

\bigskip \noindent There is a one-to-one mapping between $\mathcal{M}_{0}$ and $\mathcal{F}_{0}$
in the following way 
\begin{equation*}
\begin{tabular}{lll}
$F(u)=m(]-\infty ,u])$ & $\longleftrightarrow $ & $m(]a,b])=\Delta F(a,b)$.
\end{tabular}
\end{equation*}

\bigskip \noindent In particular, there is a one-to-one mapping between the $\mathcal{P}$ of probability measures on $\mathbb{R}$ and the class $\mathcal{F}_{P}$ of
probability distribution functions $\mathbb{R}$. Probability measures are characterized by their distribution functions.\\

\bigskip \noindent \textbf{(10.02c) General case : $k\geq 1$}. \label{doc10-02c}\\

\noindent We hope that the reader feels confident to go to the general case now. We use formulas that seem complicated but it is possible to use simple sentences to read them. For example, consider the expression of the areas with respect to the distribution function for $k=2$ :
\begin{equation*}
\Delta F(a,b)=F(b_{1},b_{2})-F(b_{1},a_{2})-F(a_{1},b_{2})+F(a_{1},a_{2}),
\end{equation*}

\bigskip \noindent with $a=(a_{1},a_{2})\leq b=(b_{1},b_{2}).$ We say that $\Delta F(a,b)$ is obtained according to the following rule :\\

\bigskip \noindent \textbf{Rule of forming $\Delta F(a,b)$}. First consider $F(b_{1},b_{2})$ the value of the distribution function at the right endpoint $b=(b_{1},b_{2})$ of the interval $]a,b]$. Next proceed to the replacements of each $b_{i}$ by $a_{i}$ by replacing exactly one of them, next two of them etc., and add each value of $F$ at the formed points, with a sign \textit{plus} $(+)$ if the number
of replacements is even and with a sign \textit{minus} $(-)$ if the number of replacements is odd.\\

\bigskip \noindent We also may use a formula. Let $\varepsilon =(\varepsilon _{1},\varepsilon_{2})\in \{0,1\}^{2}.$ We have four elements in $\{0,1\}^{2}:$ $(0,0),$ $%
(1,0),$ $(0,1),$ $(1,1).$ Consider a particular $\varepsilon _{i}=0$ or $1$, we have

\begin{equation*}
b_{i}+\varepsilon _{i}(a_{i}-b_{i})=\left\{ 
\begin{tabular}{lll}
$b_{i}$ & $if$  & $\varepsilon _{i}=0$ \\ 
$a_{i}$ & $if$ & $\varepsilon _{i}=1$%
\end{tabular}%
\right. .
\end{equation*}%

\noindent So, in
\begin{equation*}
F(b_{1}+\varepsilon _{1}(a_{1}-b_{1}),b_{1}+\varepsilon _{2}(a_{2}-b_{2})),
\end{equation*}

\bigskip \noindent the number of replacements of the $b_{i}$ by the corresponding $a_{i}$ is the number of the coordinates of $\varepsilon =(\varepsilon _{1},\varepsilon
_{2})$ which are equal to the unity (1). Clearly, the number of replacements is 
\begin{equation*}
s(\varepsilon )=\varepsilon _{1}+\varepsilon _{2}=\sum_{i=1}^{2}\varepsilon_{i}.
\end{equation*}

\bigskip \noindent We way rephrase the Rule of forming \ $\Delta F(a,b)$ into this formula
\begin{equation*}
\Delta F(a,b)=\sum_{\varepsilon =(\varepsilon _{1},\varepsilon _{2})\in
\{0,1\}}(-1)^{s(\varepsilon )}F(b_{1}+\varepsilon
_{1}(a_{1}-b_{1}),b_{1}+\varepsilon _{2}(a_{2}-b_{2})).
\end{equation*}

\bigskip \noindent We may be more compact by defining the product of vectors as the vector of the products of coordinates as

\begin{equation*}
(x,y)\ast (X,Y)=(x_1X_1, x_2X_2, ..., y_kY_k).
\end{equation*}

\bigskip \noindent The formula becomes 
\begin{equation*}
\Delta F(a,b)=\sum_{\varepsilon \in \{0,1\}}(-1)^{s(\varepsilon)}F(b+\varepsilon \ast (a-b)).
\end{equation*}

\bigskip \noindent This development has the merit to show us that the general case will be easy to get from the intermediate case $k=2.$ We may easily extend the previous rule like that.\\

\bigskip \noindent \textbf{General rule of forming $\Delta F(a,b)$}. Let $a=(a_{1},...,a_{k})\leq b=(b_{1},...,b_{k})$ two points of $\mathbb{R}^{k}$ and let $F$ an arbitrary function
from $\mathbb{R}^{k}$ to $\mathbb{R}$. We form $\Delta F(a,b)$ in this way. First consider $F(b_{1},b_{2},...,b_{k})$ the value of $F$ at right endpoint $b=(b_{1},b_{2},...,b_{k})
$ of the interval $]a,b].$ Next proceed to the replacement of each $b_{i}$ by $a_{i}$ by replacing exactly one of them, next two of them etc., and add the each value of $F$ at these points with a sign plus $(+)$ if the number of replacements is even and with a sign minus $(-)$ if the number of replacements is odd.\\

\bigskip \noindent The same reasoning gives the formula
\begin{equation*}
\Delta F(a,b)=\sum_{\varepsilon =(\varepsilon _{1},...,\varepsilon _{k})\in
\{0,1\}^{k}}(-1)^{s(\varepsilon )}F(b_{1}+\varepsilon
_{1}(a_{1}-b_{1}),...,b_{k}+\varepsilon _{k}(a_{k}-b_{k}))
\end{equation*}

\noindent and
\begin{equation*}
\Delta F(a,b)=\sum_{\varepsilon \in \{0,1\}^{k}}(-1)^{s(\varepsilon
)}F(b+\varepsilon \ast (a-b)).
\end{equation*}

\bigskip \noindent Now, we are ready to do same study already done for $k=1$ and $k=2.$ Let $m$ \ be a measure on $\mathbb{R}^{k}$\ such that for any $t=(t_{1},t_{2},...,t_{k})\in
\mathbb{R}^{k}$, 

\begin{equation*}
F_{m}(t)=m(]-\infty ,t])<\infty ,  \ (RC03)
\end{equation*}

\bigskip \noindent where 
\begin{equation*}
]-\infty ,t]=\prod_{i=1}^{k}]-\infty ,t_{i}].
\end{equation*}

\noindent The function $F_{m}:\mathbb{R}^{k}$\ $\longmapsto $ $\mathbb{R}$ defined by

\begin{equation*}
t\mapsto F_{m}(t) \ (DF03)
\end{equation*}

\bigskip \noindent has the following properties :\\

\noindent \textbf{Claim 5}. $F_{m}$\ is right-continuous at any point $t\in \mathbb{R}^{k}$, that is, 
\begin{equation*}
F_{m}(t^{(n)})\downarrow F_{m}(t)
\end{equation*}

\bigskip \noindent as 
\begin{equation*}
(t^{(n)}\downarrow t)\Leftrightarrow (\forall (1\leq i\leq
k),t_{i}^{(n)}\downarrow t_{i}).
\end{equation*}

\bigskip \noindent  \textbf{Proof}. To prove that, remark that, as $t^{(n)}\downarrow t,$\ 
\begin{equation*}
]-\infty ,t^{(n)}]=\prod_{i=1}^{k}]-\infty ,t_{i}^{(n)}]\text{ }\downarrow \text{
}]-\infty ,t]=\prod_{i=1}^{k}]-\infty ,t_{i}].
\end{equation*}

\bigskip \noindent Since $m(]-\infty ,t^{(n)}])$\ is finite for all $n>1$, we use the monotone convergence theorem to get 
\begin{equation*}
F_{m}(t^{(n)})=m(]-\infty ,t^{(n)}])\downarrow m(]-\infty ,t^{(n)}])=F_{m}(t).
\end{equation*}

\bigskip \noindent \textbf{Claim 6}. $F_{m}$ assigns non-negative volumes of cuboids, that is for  $a=(a_{1},...,a_{k})\leq b=(b_{1},...,b_{k}),$ we have  $\Delta F_{m}(a,b)\geq
0$.

\bigskip \noindent \textbf{Proof}. It is enough to establish that 
\begin{equation*}
\Delta F_{m}(a,b)=m(]a,b]).  \ (LS03)
\end{equation*}

\bigskip \noindent To do that, we proceed by induction. We know that (LS03) is true for $k=2.$ Let us make the induction assumption that (LS03) is true for any measure satisfying (RC03) and for $k$. Let us prove it for $k+1.$ Let $a=(a_{1},...,a_{k+1})\leq b=(b_{1},...,b_{k+1})$. We have 
\begin{equation*}
]a,b]=\prod\limits_{i=1}^{k+1}]a_{i},b_{i}]
\end{equation*}

\begin{equation*}
=\left\{ \left\{ \prod\limits_{i=1}^{k}]a_{i},b_{i}]\right\} \times ]-\infty
,b_{k+1}]\right\} \backslash \left\{ \left\{\prod\limits_{i=1}^{k}]a_{i},b_{i}]\right\} \times ]-\infty,a_{k+1}]\right\} .
\end{equation*}

\bigskip \noindent Since the right-side term of the sets minus operator ($\setminus$) is a subset of the left-side term, we have
\begin{eqnarray*}
m(]a,b])&=&m\left( \left\{ \prod\limits_{i=1}^{k}]a_{i},b_{i}]\right\} \times ]-\infty ,b_{k+1}]\right)\\
&-&m\left( \left\{\prod\limits_{i=1}^{k}]a_{i},b_{i}]\right\} \times ]-\infty ,a_{k+1}]\right) .
\end{eqnarray*}

\bigskip \noindent Set

\begin{equation*}
\overline{a}=(a_{1},...,a_{k})\leq \overline{b}=(b_{1},...,b_{k})
\end{equation*}%

\bigskip \noindent and denote
\begin{eqnarray*}
]\overline{a},\overline{b}]\longrightarrow m_{1}(]\overline{a},\overline{b} ])&=&m\left( ]\overline{a},\overline{b}]\times ]-\infty ,b_{k+1}]\right)\\
&=&m\left( \left\{ \prod\limits_{i=1}^{k}]a_{i},b_{i}]\right\} \times ]-\infty ,b_{k+1}]\right) 
\end{eqnarray*}

\bigskip \noindent and 

\bigskip 
\begin{eqnarray*}
]\overline{a},\overline{b}]\longrightarrow m_{2}(]\overline{a},\overline{b}])&=&m\left( ]\overline{a},\overline{b}]\times ]-\infty ,a_{k+1}]\right)\\
&=&m\left( \left\{ \prod\limits_{i=1}^{k}]a_{i},b_{i}]\right\} \times ]-\infty ,a_{k+1}]\right) .
\end{eqnarray*}

\bigskip \noindent When $a_{k+1}$ and $b_{k+1}$ are fixed, $m_{1}$ and $m_{2}$ are measures on $\mathbb{R}^{k}$ assigning finite values of bounded above intervals of $\mathbb{R}^{k}$. So we
may apply the induction assumption with
\begin{equation*}
F_{m_{1}}(x)=F_{m}(x,b_{k+1})\text{ and } F_{m_{2}}(x)=F_{m}(x,a_{k+1}),\text{ }x\in \mathbb{R}^{k}.
\end{equation*}

\bigskip \noindent and get

\begin{eqnarray*}
m(]a,b])&=&\sum_{\varepsilon \in \{0,1\}^{k}}(-1)^{s(\varepsilon )}F_{m}(\overline{b}+\varepsilon \ast (\overline{a}-\overline{b}),b_{k+1})\\
&-&\sum_{\varepsilon \in \{0,1\}^{k}}(-1)^{s(\varepsilon )}F(\overline{b}+\varepsilon \ast (\overline{a}-\overline{b}),a_{k+1}).
\end{eqnarray*}

\bigskip \noindent Let us write it in a more developed way, that is
\begin{eqnarray*}
&&m(]a,b])   \ \ \ \ (PI04) \\
&=&\sum_{\varepsilon \in \{0,1\}^{k}}(-1)^{s(\varepsilon)} F_{m}(b_{i}+\varepsilon _{i}(a_{i}-b_{i}),b_{k+1},1\leq i\leq k)\\
&-&\sum_{\varepsilon \in \{0,1\}^{k}}(-1)^{s(\varepsilon )}F_{m}(b_{i}+\varepsilon _{i}(a_{i}-b_{i}),a_{k+1},1\leq i\leq k). 
\end{eqnarray*}

\bigskip \noindent Now, a way of forming $\{0,1\}^{k+1}$ from $\{0,1\}^{k}$ consists of considering an arbitrary element $\varepsilon =(\varepsilon _{1},...,\varepsilon _{k})$ and
to add either $\varepsilon _{k+1}=0$ or $\varepsilon _{k+1}=0$ as the $(k+1)$-th coordinate. Thus, we have

\bigskip 
\begin{equation*}
\{0,1\}^{k+1}=\{(\varepsilon ,0),\varepsilon \in \{0,1\}^{k}\}\cup
\{(\varepsilon ,0),\in \{0,1\}^{k}\}.
\end{equation*}

\bigskip \noindent Hence, we get
\begin{eqnarray*}
\Delta F_{m}(a,b)&=&\sum_{\varepsilon \in \{0,1\}^{k+1}}(-1)^{s(\varepsilon )}F_{m}(b_{i}+\varepsilon _{i}(a_{i}-b_{i}),1\leq i\leq k+1)\\
&=&\sum_{(\varepsilon ,0),\varepsilon \in \{0,1\}^{k}}(-1)^{s((\varepsilon ,0))}F_{m}(b_{i}+\varepsilon _{i}(a_{i}-b_{i}),1\leq i\leq k+1)\\
&+&\sum_{(\varepsilon ,1),\varepsilon \in \{0,1\}^{k}}(-1)^{s((\varepsilon ,1))}F_{m}(b_{i}+\varepsilon _{i}(a_{i}-b_{i}),1\leq i\leq k+1).
\end{eqnarray*}

\bigskip \noindent Now $(-1)^{s((\varepsilon ,0))}=(-1)^{s(\varepsilon )}$ and $(-1)^{s((\varepsilon ,1))}=(-1)^{s(\varepsilon )+1}=-(-1)^{s(\varepsilon )}.$
Next in the first summation, we don't have replacement of $b_{k+1}$ since $\varepsilon _{k+1}=0,$\ while in the second $b_{k+1}$ is replaced by $%
a_{k+1}.$ We finally get that $\Delta F_{m}(a,b)$ is equal to : 

$$
\sum_{(\varepsilon ,0),\varepsilon \in \{0,1\}^{k}}(-1)^{s(\varepsilon
)}F_{m}(b_{i}+\varepsilon _{i}(a_{i}-b_{i}),b_{k+1},1\leq i\leq
k)
$$

$$
-\sum_{(\varepsilon ,0),\varepsilon \in \{0,1\}^{k}}(-1)^{s(\varepsilon
)}F_{m}(b_{i}+\varepsilon _{i}(a_{i}-b_{i}),b_{k+1},1\leq i\leq k),
$$

\bigskip \noindent which is the left-side member of (PI04). We conclude that 
\begin{equation*}
\Delta F_{m}(a,b)=m(]a,b])\geq 0.
\end{equation*}

\bigskip \noindent We also have this property. If $m$ is a probability measure, we have :\\

\noindent \textbf{Claim 7}.  $F_{m}$ satisfies (i) :
\begin{equation*}
\lim_{\exists i,1\leq i\leq k,t_{i}\rightarrow -\infty
}F_{m}(t_{1},...,t_{k})=0,
\end{equation*}%

\bigskip \noindent and (ii) :

\begin{equation*}
\lim_{\forall i,1\leq i\leq k,t_{i}\rightarrow +\infty
}F_{m}(t_{1},...,t_{k})=1.
\end{equation*}

\bigskip \noindent This is left as a simple exercise. We are in a position to give the following definitions.\\

\noindent \textbf{Definition}. A function $F$ 
\begin{equation*}
\begin{array}{ccc}
\mathbb{R}^{k} & \mapsto  & \mathbb{R} \\ 
t & \hookrightarrow  & F(t)%
\end{array}%
\end{equation*}

\bigskip \noindent is a distribution function on $\mathbb{R}^{k}$  if\\

\noindent (a) $F$ assigns non-negative volumes to cubes, that is $\Delta F(a,b)\geq 0$ for $a\leq b$\\

\noindent and \\

 \noindent (b) $F$ is right-continuous.\\

\bigskip \noindent \textbf{Definition}. A function $F$ 
\begin{equation*}
\begin{array}{ccc}
\mathbb{R}^{k} & \mapsto  & \mathbb{R} \\ 
t & \hookrightarrow  & F(t),
\end{array}%
\end{equation*}

\noindent is a probability distribution function on $\mathbb{R}^{k}$\ if\\

\noindent (a) $F$ assigns non-negative volumes to cubes, that is $\Delta F(a,b)\geq 0$, for $a\leq b$\\

\bigskip \noindent (b) $F$ is right-continuous;\\

\bigskip \noindent (c) $F_{m}$ satisfies (i)
\begin{equation*}
\lim_{\exists i,1\leq i\leq k,t_{i}\rightarrow -\infty
}F_{m}(t_{1},...,t_{k})=0
\end{equation*}

\bigskip \noindent and (ii)
\begin{equation*}
\lim_{\forall i,1\leq i\leq k,t_{i}\rightarrow +\infty
}F_{m}(t_{1},...,t_{k})=1.
\end{equation*}

\bigskip \noindent As we did twice, we may take the reverse way and establish for any distribution function $F$ this theorem.\\

\bigskip \noindent \textbf{Existence Theorem of Lebesgue-Stieljes Measures}.\\

\noindent Let $F$ be a distribution function on  $\mathbb{R}^{k},$\ then the application $m_{F}$
that is defined on the semi-algebra 
\begin{equation*}
\mathcal{S}=\{]a,b]=\prod_{i=1}^{k}]a_{i},b_{i}],\text{ }(a,b)\in (\overline{%
\mathbb{R}}^{k})^{2}\}
\end{equation*}

\bigskip \noindent by 
\begin{equation*}
\left\{ 
\begin{array}{cccc}
m_{F}: & \mathcal{S} & \mapsto  & \mathbb{R}_{+} \\ 
& ]a,b] & \hookrightarrow  & \Delta _{a,b}F%
\end{array}%
\right. 
\end{equation*}

\bigskip \noindent is additive and is uniquely extensible to a measure on $\mathbb{B}(\mathbb{R}^{k}),$\ named
as the Lebesgue-Stieljes associated with $F$. If $F$ satisfies%

$$
\lim_{\exists i,1\leq i\leq k,t_{i}\rightarrow -\infty }F(t_{1},...,t_{k})=0,
$$

\bigskip \noindent then $F$ is exactly the distribution function generated by the measure $m_{F}$ that is
\begin{equation*}
F(t)=m_{F}(]\infty ,t]),t\in \mathbb{R}^{k}.
\end{equation*}

\bigskip \noindent If $F$ is a probability distribution function, then $m_{F}$ is a probability measure.\\

\bigskip \noindent This theorem is announced in any of the three parts. It constitutes the foundation of the theory exposed here. But the reader does not necessarily need it to understand Lebesgue-Stieljes measures and its many usages. So, we postpone in the next document where fundamentalist mathematicians will find it.\\

\bigskip \noindent But to be coherent, let us prove the additivity of $m_{F}$ as we did for the two first cases. An element $]a,b]=\prod_{i=1}^{k}]a_{i},b_{i}]$\ of $S$%
\ is decomposable into two elements of $\mathcal{S}$ only by splitting of its $k$ components ;  $]a_{i},b_{i}]=]a_{i},c]+]c,b_{i}],$ with $a_{i}\leq
c\leq b_{i}.$\ For simplicity's sake, say the the first component is split. We have  
\begin{eqnarray*} 
A&=&\prod_{i=1}^{k}]a_{i},b_{i}]=]a_{1},c]\times \prod_{i=2}^{k}]a_{i},b_{i}]+]c,b_{1}]\times \prod_{i=2}^{k}]a_{i},b_{i}]\\
&\equiv& B+C.
\end{eqnarray*}

\bigskip \noindent We have
\begin{eqnarray*}
&&m_{F}(B)\\
&=&\sum_{\varepsilon =\{0,1\}^{k}}(-1)^{s(\varepsilon )}F_{m}(c+\varepsilon _{1}(a_{1}-c),b_{2}+\varepsilon _{k}(a_{k}-b_{k})+...,b_{k}+\varepsilon _{k}(a_{k}-b_{k}))
\end{eqnarray*}

\bigskip \noindent and 
\begin{eqnarray*}
&&m_{F}(C)\\
&=&\sum_{\varepsilon =\{0,1\}^{k}}(-1)^{s(\varepsilon )}F_{m}(b_{1}+\varepsilon _{1}(c-b_{1}),b_{2}+\varepsilon _{2}+(a_{2}-b_{2})+...,b_{k}+\varepsilon _{k}(a_{k}-b_{k})).
\end{eqnarray*}

\bigskip \noindent We say that a way of forming  $\ \{0,1\}^{k}$ is to consider an arbitrary element $\varepsilon ^{(k-1)}=(\varepsilon _{2},\varepsilon _{3},...,\varepsilon _{k})\in \{0,1\}^{k-1}$ and add it either $\varepsilon_{1}=0$ or $\varepsilon _{1}=1$ as first coordinate so that we have 

\begin{eqnarray*}
&&\{0,1\}^{k} \ \ (PI05)\\
&=&\{(0,\varepsilon ^{(k-1)}),\varepsilon ^{(k-1)}\in \{0,1\}^{k-1}\}\cup \{(1,\varepsilon ^{(k-1)}),\varepsilon ^{(k-1)}\in \{0,1\}^{k-1}\}. 
\end{eqnarray*}

\bigskip \noindent For both $m(B)$ and $m(C),$ we are going to develop the summation accordingly to (PI05) to get 
\begin{eqnarray*}
&&m_{F}(B)\\
&=&\sum_{\varepsilon =(0,\varepsilon ^{(k-1)})\in \{0,1\}^{k}}(-1)^{s(\varepsilon )}F_{m}(c,b_{2}+\varepsilon
_{2}(a_{2}-b_{2})+...,b_{k}+\varepsilon _{k}(a_{k}-b_{k}))\\
&+&\sum_{\varepsilon =(1,\varepsilon ^{(k-1)})\in \{0,1\}^{k}}(-1)^{s(\varepsilon )}F_{m}(a_{1},b_{2}+\varepsilon _{2}(a_{2}-b_{2})+...,b_{k}+\varepsilon _{k}(a_{k}-b_{k}))
\end{eqnarray*}

\bigskip \noindent and 
\begin{eqnarray*}
&&m_{F}(C)\\
&=&\sum_{\varepsilon =(0,\varepsilon ^{(k-1)})\in \{0,1\}^{k}}(-1)^{s(\varepsilon )}F_{m}(b_{1},b_{2}+\varepsilon _{2}(a_{2}-b_{2})+...,b_{k}+\varepsilon _{k}(a_{k}-b_{k}))\\
&+&\sum_{\varepsilon =(0,\varepsilon ^{(k-1)}) \in \{0,1\}^{k}}(-1)^{s(\varepsilon )}F_{m}(c,b_{2}+\varepsilon _{2}(a_{2}-b_{2})+...,b_{k}+\varepsilon _{k}(a_{k}-b_{k})).
\end{eqnarray*}

\bigskip \noindent The first line for the expression of $m(B)$ and the second line of the expression of $m(C)$ are opposite since for each $\varepsilon ^{(k-1)}$, the value of the associated $\varepsilon _{1}$ is zero for $m(B)$ and is one for $m(C)$ so that the signs given by the $(-1)^{s(\varepsilon )}$ are opposite while their absolute values are identical. It remains

\begin{eqnarray*}
&&m_{F}(B)+m_{F}(C)\\
&=&\sum_{\varepsilon =(0,\varepsilon ^{(k-1)}) \in \{0,1\}^{k}}(-1)^{s(\varepsilon )}F_{m}(b_{1},b_{2}+\varepsilon _{2}(a_{2}-b_{2})+...,b_{k}+\varepsilon _{k}(a_{k}-b_{k}))\\
&+&\sum_{\varepsilon =(1,\varepsilon ^{(k-1)}) \in \{0,1\}^{k}}(-1)^{s(\varepsilon )}F_{m}(a_{1},b_{2}+\varepsilon _{2}(a_{2}-b_{2})+...,b_{k}+\varepsilon _{k}(a_{k}-b_{k})),
\end{eqnarray*}

\bigskip \noindent which is the development of $m(A)$ accordingly with (PI05). So the additivity of $m_{F}$ is established.\\

\bigskip \noindent \textbf{(10.03) Factorized forms for a distribution functions}.\\

\noindent \textbf{(a) Totally Factorized form}.\\

\noindent Let $F_i$, $1\leq u \leq k\geq 2$, $k$ \textsl{df}'s (resp. probability distribution functions) on $\mathbb{R}$. We have

$$
F(x_1,...,x_k)=\prod_{1 \leq i \leq k} F_i(x_i), (x_1,...,x_k)^t \in \mathbb{R}^k, \ (FACT01)
$$

\bigskip \noindent defines  a \textsl{pf} (resp.  probability distribution function) on $\mathbb{R}^k$ and for any $a=(a_{1},...,a_{k})\leq b=(b_{1},...,b_{k})$,

$$
\Delta_{a,b}F= \prod_{1 \leq i \leq k} \delta_{a_i,b_i}F_i. \ (FACT02)
$$

\bigskip \noindent Let us give a quick proof based on an induction reasoning. Let $k=2$, we have 

\begin{eqnarray*}
\Delta_{a,b}&=&F_1(b_1)F_2(b_2)-F_1(a_1)F_2(b_2)-F_1(b_1)F_2(a_2)+F_1(a_1)F_2(a_2)\\
&=& F_2(b_2) \biggr(F_1(b_1)-F_1(a_1)\biggr) - F_2(a_2) \biggr(F_1(b_1)-F_1(a_1)\biggr)\\
&=& \biggr(F_1(b_1)-F_1(a_1)\biggr) \biggr(F_1(b_2)-F_1(a_2)\biggr).
\end{eqnarray*}

\bigskip \noindent Hence for $k=2$, the function $F$ is right-continuous and assigns non-negative areas of rectangles and Formula (FACT02) holds. Suppose that the results are true for some $k-1\geq 3$. For the value $k$, $F$ is still obviously right-continuous. Let us use the same techniques as above to have

\begin{eqnarray*}
&&\Delta_{a,b}F\\
&=&\sum_{(\varepsilon,0) \in \{0,1\}^{k}}(-1)^{s(\varepsilon )} F(b+(\varepsilon,0)*(b-a))\\
&-&\sum_{(\varepsilon,1) \in \{0,1\}^{k}} (-1)^{s(\varepsilon )}F(b+(\varepsilon,1)*(b-a))\\
&=&:A-B.
\end{eqnarray*}

\bigskip 
\noindent But, by denoting

$$
F^{\ast}(x_1,...,x_{k-1})=\prod_{1 \leq i \leq k-1} F_i(x_i), (x_1,...,x_{k-1})^t \in \mathbb{R}^{k-1},
$$

\bigskip 
\noindent and

$$
b^{\ast}=(b_1,...,b_{k-1}) \text{ and } a^{\ast}=(a_1,...,a_{k-1}),
$$

\bigskip 
\noindent we have

\begin{eqnarray*}
A&=&F_k(b_k) \biggr( \sum_{\varepsilon \in \{0,1\}^{k-1}}(-1)^{s(\varepsilon )} F^{\ast}(b^{\ast}+\varepsilon*(b^{\ast}-a^{\ast}))\biggr)\\
&=& F_k(b_k) \Delta_{b^{\ast},a^{\ast}}F^{\ast}
\end{eqnarray*}

\bigskip 
\noindent and 

$$
B=F_k(a_k) \Delta_{b^{\ast},a^{\ast}}F^{\ast}.
$$

\bigskip 
\noindent Thus, we have

$$
\Delta_{a,b}= \Delta_{b^{\ast},a^{\ast}}F^{\ast} \biggr(F_k(b_k)-F_k(a_k)\biggr). \ (FACT03)
$$

\bigskip 
\noindent By the induction hypothesis, we have

$$
\Delta_{b^{\ast},a^{\ast}}F^{\ast}=\prod_{1 \leq i \leq k-1} (F_i(b_i)-F_i(d_i)).\ (FACT04)
$$

\bigskip 
\noindent We conclude by combining Formulas (FACT02) and (FACT03) that $F$ assigns non-negative volumes to cuboids, and hence is a \textsl{df}, and Formula (FAC01) holds.\\

\noindent It is clear that if each $F_i$ is a the \textsl{df} of a probability measure, $F$ also is. $\square$\\

\noindent \textbf{Definition} The function $F$ is said to be in a factorized form if Formula (FACT01) holds.\\

\bigskip \noindent \textbf{(b) Another Factorized form}.\\

\noindent Let $G$ and $H$ be two \textsl{df} (resp. probability distribution functions) respectively defined on $\mathbb{R}^r$ and on $\mathbb{R}^s$, with $r\geq 1$ and $s\geq 1$. Set $k=r+s$. Thus the function 

$$
(y,z)\in \mathbb{R}^r \times \mathbb{R}^s \mapsto F(y,z)=F(y)H(z),
$$

\bigskip 
\noindent is a \textsl{df} (resp. probability distribution function) and for $a=(a_{1},...,a_{k})\leq b=(b_{1},...,b_{k})$ with
$$
a=(a^{(1)},a^{(2)}, \ b=(b^{(1)},b^{(2)}), \ \ (a^{(1)},b^{(1)}) \in \left(\mathbb{R}^r\right)^2, \ \ (a^{(2)},b^{(2)}) \in \left(\mathbb{R}^s\right)^2.  
$$

\bigskip 
\noindent we have

$$
\Delta_{a,b}F=\Delta_{a^{(1)},b^{(1)}}G \Delta_{a^{(2)},b^{(2)}}H. \ (FACT05)
$$

\bigskip 
\noindent \textbf{Proof}. Suppose that $F$ and $G$ are \textsl{df}'s. Let us consider $a=(a_{1},...,a_{k})\leq b=(b_{1},...,b_{k})$ and adopt the notations of the statements. We have

\begin{eqnarray*}
&&\Delta_{a,b}F\\
&=&\sum_{(\varepsilon^{(1)},\varepsilon^{(2)}) \in \{0,1\}^{r} \times \{0,1\}^{s}} (-1)^{s(\varepsilon^{(1)})+s(\varepsilon^{(2)}) )} F(b+((\varepsilon^{(1)},\varepsilon^{(2)})*(b-a))\\
&=&\sum_{\varepsilon^{(1)} \in \{0,1\}^{r}}\sum_{\varepsilon^{(2)} \in \{0,1\}^{s}} \biggr((-1)^{s(\varepsilon^{(1)})}H(b^{1}+\varepsilon^{(1)}*(b^{1}-a^{1}))\biggr)\\
&\times& \biggr( (-1)^{s(\varepsilon^{(2)})} H(b^{(2)}+\varepsilon^{(2)}*(b^{2}-a^{(2)}))\biggr)\\
&=&\biggr(\sum_{\varepsilon^{(1)} \in \{0,1\}^{r}} (-1)^{s(\varepsilon^{(1)})}H(b^{(1)}+\varepsilon^{(1)}*(b^{1}-a^{(1)}))\biggr)\\
&\times& \biggr( (-1)^{s(\varepsilon^{(2)})} H(b^{(2)}+\varepsilon^{(2)}*(b^{(2)}-a^{(2)}))\biggr)\\
&=&\Delta_{a^{(1)},b^{(1)}} \times G \Delta_{a^{(2)},b^{(2)}}H.
\end{eqnarray*}

\bigskip \noindent So, $F$ is right-continuous and assigns non-negative areas to cuboids and hence is a\textsl{df}. The conclusion for probability distribution functions follows easily.

\noindent \LARGE \textbf{DOC 10-02 : Proof of the Existence of the Lebesgue-Stieljes measure on $\mathbb{R}^k$} \label{doc10-02}\\
\bigskip

\Large
\noindent We have to prove the following :

\noindent \textbf{THEOREM}. 
Let $F$ be a distribution function on $\mathbb{R}^{k},$\ then the application $m_{F}$
that is defined on the semi-algebra 
\begin{equation*}
\mathcal{S}_{k}=\{]a,b]=\prod_{i=1}^{k}]a_{i},b_{i}],\text{ }(a,b)\in (%
\overline{\mathbb{R}}^{k})^{2}\}
\end{equation*}%
by 
\begin{equation*}
\left\{ 
\begin{array}{cccc}
m_{F}: & \mathcal{S}_{k} & \mapsto  & \mathbb{R}_{+} \\ 
& ]a,b] & \hookrightarrow  & \Delta _{a,b}F%
\end{array}%
\right. 
\end{equation*}%

\bigskip \noindent is additive and is uniquely extensible to a measure on $\mathcal{B}(\mathbb{R}^{k}).$

\bigskip \noindent \textbf{Proof of the theorem}.\\

\noindent Let us begin to see that $m$ is $\sigma$-finite since

$$
\mathbb{R^k}= \bigcup_{n\geq 1} ]-\infty, \{n\}^k] \ and, \ \forall n\geq 1, m_F(]-\infty, \{n\}^k])\geq 1.
$$

\bigskip \noindent \noindent We already know $m_{F}$ is proper, non-negative, additive and $\sigma$-finite on $\mathcal{S_k}$. Thus, by advanced version of Caratheodory's Theorem (see Point 4.28 of Doc 04-03) to uniquely extend $m_{F}$ to a measure on $\mathcal{R}^k$, it will be enough to show that $m_F$ is $\sigma$-sub-finite on $\mathcal{S_k}$. Since $m_F$ is already non-negative, additive, we use the characterization of $\sigma$-sub-additivity as established in \textit{Exercise 4 in Doc 04-02}.\\

\noindent We also get that $m_F$ extended to a non non-negative and additive application still denoted by $m_F$ so that $m_F$ is non-decreasing and sub-additive so $\mathcal{S_k}$ and on $\mathcal{C}$.\\

\noindent So we have to prove that : if the interval $]a,b]$, with $a\in \overline{\mathbb{R}}^k$, $b\in \overline{\mathbb{R}}^k$, decomposed like
\begin{equation*}
]a,b]=\sum_{n\geq 1}]a_n,\ b_n].
\end{equation*}

\noindent then we have

\begin{equation}
m_{F}(]a,b]) \leq\sum_{n\geq 1}m_{F}(]a_{n},b_{n}]).  \label{sousadditivite1}
\end{equation}

\noindent To prove that we proceed by two steps.\\

\noindent Step 1 : $a$ and $b$ are finite.\\

\noindent If  $a_{i}=b_{i}$ for some $i,1\leq i\leq k,$ we have that $]a,b]=\emptyset $
and (\ref{sousadditivite1}) holds. Next, let $a_{i}<b_{i}$\ for all  $1\leq
i\leq k.$\ We can find $\delta =(\delta _{1},...,\delta _{k})0$ such that $\delta _{i}>0$ for all $1\leq i\leq k$ so that   
\begin{equation}
\lbrack a+\delta ,b]\subseteq ]a,b].
\end{equation}

\bigskip \noindent The vectors $a^{n}$ and $b^{n}$\ have finite components since  $a\leq
a^{n}\leq b^{n}\leq b.$\ Since $F$ is right-continuous at $b_{n},$ we have%
\begin{equation*}
m(]a^{n},b^{n}+\varepsilon ])=\Delta F(a^{n},b^{n}+\varepsilon )\downarrow
\Delta F(a^{n},b^{n})=m(]a^{n},b^{n}])
\end{equation*}

\bigskip \noindent  as 
\begin{equation*}
\varepsilon =(\varepsilon _{1},\varepsilon _{2},...,\varepsilon
_{k})\downarrow 0.
\end{equation*}%

\bigskip \noindent So for any $\eta >0,$\ there exists $\varepsilon ^{n}=(\varepsilon
_{1}^{n},\varepsilon _{2}^{n},...,\varepsilon _{k}^{n})>0$\ such that  
\begin{equation}
m(]a^{n},b^{n}])\leq m(]a^{n},b^{n}+\varepsilon ^{n}])\leq
m(]a^{n},b^{n}])+\eta /2^{-n}.  \label{lj08}
\end{equation}

\bigskip \noindent Then 
\begin{equation}
\lbrack a+\delta ,b]\subseteq ]a,b]\subseteq \bigcup_{n\geq
1}]a^{n},b^{n}]\subseteq \bigcup_{n\geq 1}]a^{n},b^{n}+\varepsilon
^{n}/2[\subseteq \bigcup_{n\geq 1}]a^{n},b^{n}+\varepsilon ^{n}]
\label{lj09}
\end{equation}%

\bigskip \noindent and next 
\begin{equation}
\lbrack a+\delta ,b]\subseteq \bigcup_{n\geq 1}]a^{n},b^{n}+\varepsilon
^{n}/2[.  \label{lj10}
\end{equation}

\bigskip \noindent The compact set $[a+\delta ,b]$\ is included in the union of open sets $%
\bigcup_{1\leq j\leq p}]a^{n_{j}},b^{n_{j}}+\varepsilon ^{n_{j}}/2[.$ By the
Heine-Borel property, we may extract from this union a finite union $%
\bigcup_{1\leq j\leq p}]a^{n_{j}},b^{n_{j}}+\varepsilon ^{n_{j}}/2[$ such
that 
\begin{equation*}
\lbrack a+\delta ,b]\subseteq \bigcup_{1\leq j\leq
p}]a^{n_{j}},b^{n_{j}}+\varepsilon ^{n_{j}}/2[.
\end{equation*}

\bigskip \noindent Now, we get,  as $\delta \downarrow 0,$\ 
\begin{equation*}
]a,b]\subseteq \bigcup_{1\leq j\leq p}]a^{n_{j}},b^{n_{j}}+\varepsilon^{n_{j}}.
\end{equation*}

\bigskip \noindent Recall that any  finite union $A_{1}\cup A_{2}\cup ...\cup A_{p}$ may be
transformed into a sum like this%
\begin{eqnarray*}
A_{1}\cup A_{2}\cup ...\cup A_{p}&&=A_{1}+A_{1}^{c}A_{2}+A_{1}^{c}A_{2}^{c}A_{3}+...+A_{1}^{c}...A_{2}^{p-1}A_{3}^{p}\\
&=&B_{1}+...+B_{k},
\end{eqnarray*}

\bigskip \noindent where $B_{i}\subset A_{i},i=1,...,k$, from which we see that if the $A_{i}$ are in semi-algebra, the terms of the
sum are also in the semi-algebra. Here, put $A_{i}=]a^{n_{j}},b^{n_{j}}+%
\varepsilon ^{n_{j}}],$ we get%
\begin{eqnarray*}
m_{F}(]a,b])&\leq& m_{F}(\bigcup_{1\leq j\leq p}]a^{n_{j}},b^{n_{j}}+\varepsilon ^{n_{j}}])\\
&=&\sum_{i=1}^{k}m_{F}(B_{i})\\
&\leq&\sum_{i=1}^{k}m_{F}(]a^{n_{j}},b^{n_{j}}+\varepsilon ^{n_{j}}])\\
&\leq& \sum_{1\leq i\leq p}\left\{ m(]a^{n_{j}},b^{n_{j}}])+\eta 2^{n_{j}}\right\}\\
&\leq& \sum_{1\leq i\leq p}m(]a^{n},b^{n}])+\eta\\
&\leq& \sum_{n\geq 1}m(]a^{n},b^{n}])+\eta .
\end{eqnarray*}

\bigskip \noindent Since this is true for an arbitrary value of $\eta >0,$ we may let $\eta
\downarrow 0$ to have
\begin{equation*}
m(]a,b])\leq \sum_{n\geq 1}m(]a^{n},b^{n}]).
\end{equation*}

\bigskip \noindent \textbf{Step 2}. We consider a general element $]a,b]$ of $\mathcal{S}_{k}
$ with the possibility that some $a_{i}$ migth be $-\infty$ and that some $b_{i}$ might be $+\infty$. Let 
$I_{1}=\{i,1\leq i\leq k,a_{i}=-\infty \},I_{2}=\{i,1\leq i\leq
k,b_{i}=+\infty \}$%
\begin{equation*}
]a,b]=\sum_{n\geq 1}]a^{n},b^{n}].
\end{equation*}

\noindent If one of the $a_{i}$ is $+\infty $\ or one the $b_{i}$ is -$\infty ,$ $]a,b]
$ is empty and sub-additivity is automatic. We exclude these possibilities in
the sequel. Accordingly, let $\alpha \in \mathbb{R}^{k}$ with $\alpha \leq b$ and let us 
define it in the following way : $if$ $i\in I_{1}$ take $\alpha _{i}\leq b_{i}.
$ Take $\alpha _{i}=a_{i}$ if $a_{i}$ is finite. Define also $\beta
\in \mathbb{R}^{k}$ with $\beta \geq \alpha ,$ in the following way : $if$ $i\in I_{2}
$ take $\beta _{i}\leq \alpha _{i}$. Take $\beta_{i}=b_{i}$ if $b_{i}$
is finite. We get  
\begin{equation*}
]a,b]\cap ]\alpha ,\beta ]=]\alpha ,\beta ]=\sum_{n\geq 1}]a^{n},b^{n}]\cap
]\alpha ,\beta ].
\end{equation*}%

\noindent By the first step, we obtain%
\begin{equation*}
m_{F}(]\alpha ,\beta ])\leq \sum_{n\geq 1}m_{F}(]a^{n},b^{n}]\cap ]\alpha
,\beta ])\leq \sum_{n\geq 1}m_{F}(]a^{n},b^{n}]).
\end{equation*}

\bigskip \noindent Now by the definition of $F$ with possibly infinite arguments, we get that
for $\alpha _{i}\rightarrow -\infty $ for $i\in I_{1}$ and $b_{i}\rightarrow
+\infty$ for $i\in I_{2}$, and that $m_{F}(]\alpha ,\beta ])=\Delta
F(\alpha ,\beta )\rightarrow \Delta F(a,b)=m_{F}(]a,b])$ and we arrive at%

\begin{equation*}
m_{F}(]\alpha ,\beta ])\leq \sum_{n\geq 1}m_{F}(]a^{n},b^{n}]). \ \blacksquare
\end{equation*}

\noindent \LARGE \textbf{DOC 10-03 : Functions of bounded variation} \label{doc10-03}\\
\Large

\noindent Later and hopefully, the reader will study Stochastic Calculus. The origin of such a theory is related to the notion of bounded variation. We will shortly see that the Riemann-Stieljes related to a function $F$ are well-defined if $F$ is of bounded variation. Actually, we will show that of function of bounded variation is difference between two monotone functions and this makes possible the extension of the Lebesgue-Stieljes integral with respect to a bounded-variation function.\\

\noindent Since the Brownian motion $B$  is non-where of bounded variation, integrating with respect to $dB_t$, path-wisely, in the Stieljes sense does not make sense. This gave  birth to stochastic integration.\\

\noindent Let us begin by the concept of total variation.\\

\bigskip \noindent \textbf{A - Bounded variation functions}.\\

\noindent  Let $F:$ $]a,b]\longrightarrow \mathbb{R}$ be finite function. Consider a subdivision $\pi =(a=x_{0}<x_{1}<...<x_{\ell }<x_{\ell +1}=b)$ of $[a,b].$
The variation of $F$ on $[a,b]$ over $\pi $ is defined by%
\begin{equation*}
V(F,[a,b],\pi )=\sum_{i=0}^{\ell-1}\left\vert F(x_{i+1})-F(x_{i+1})\right\vert .
\end{equation*}

\bigskip \noindent Denote by $\mathcal{P}([a,b])$ the set of all subdivisions of $[a,b]$. We may define\\

\noindent \textbf{Definition}. The total variation $V(F,[a,b])$ of  $F$ on $[a,b]$ is the supremum of its variation over all the subdivisions of $[a,b]$, that is%
\begin{equation*}
V(F,[a,b])=\sup_{\pi \in \mathcal{P}([a,b])}V(F,[a,b],\pi ).
\end{equation*}

\bigskip \noindent The function $F$ is of bounded variation on $[a,b]$ if and only $V(F,[a,b])$ is finite.\\

\noindent \textbf{Examples}. Let $F$ be monotonic on $[a,b]$ , we easily see that any variation of $F$ on is exactly $\left\vert F(a)-F(b)\right\vert .$ To set that, let $F$ be non-decreasing. We have
\begin{equation*}
V(F,[a,b],\pi )=\sum_{i=0}^{m}F(x_{i+1})-F(x_{i})
\end{equation*}%
\begin{equation*}
=F(x_{1})-F(x_{0})+F(x_{2})-F(x_{1})+...+F(x_{\ell +1})-F(x_{\ell })
\end{equation*}%
\begin{equation*}
=F(x_{\ell +1})-F(x_{0})=F(a)-F(b).
\end{equation*}

\bigskip \noindent If $F$ is non-increasing, apply thus result to $-F$. We have this statement :\\

\noindent \textit{Monotonic functions $F$ on $[a,b]$ are of bounded variation}.\\

\noindent It is remarkable that non-decreasing functions characterize bounded variation functions in the following way.\\

\noindent \textbf{Proposition}. Any bounded variation $F$ on $[a,b]$ is a difference of two non-decreasing
functions.\\

\bigskip \noindent \textbf{Proof}. Let $a\leq x<y\leq b.$ Pick subdivision $\pi _{1}$ $\in \mathcal{P}([a,x])$ and a subdivision $\pi _{2}$ $\in \mathcal{P}([x,y]).$ It is clear that we
get a subdivision $\pi $ of $[a,y]$ by binding $\pi _{1}$ and $\pi _{2}$ at the point $x$ and 

\begin{equation*}
V(F,[a,y],\pi )=V(F,[a,x],\pi _{1})+V(F,[x,y],\pi _{2}).
\end{equation*}

\bigskip \noindent Since $V(F,[a,y],\pi )\leq V(F,[a,y])$, we get

\begin{equation*}
V(F,[a,x],\pi _{1})+V(F,[x,y],\pi _{2})\leq V(F,[a,y]).
\end{equation*}

\bigskip \noindent Since $\pi _{1}$ $\ $and $\pi _{2}$ are arbitrary elements of $\mathcal{P}([a,x])$ and $\mathcal{P}([x,y])$, we may get the suprema to get%

\begin{equation}
V(F,[a,x])+V(F,[x,y])\leq V(F,[a,y]).  \label{compvar}
\end{equation}

\bigskip \noindent Put

\begin{equation*}
F_{1}(x)=V(F,[a,x]),a\leq x\leq b.
\end{equation*}

\bigskip \noindent We see by \ref{compvar} that for $a\leq x<y\leq b$%
\begin{equation*}
F_{1}(y)\geq F_{1}(x)+V(F,[x,y])\geq F_{1}(x),
\end{equation*}

\bigskip \noindent since $V(F,[x,y])$. Then $F_{1}$ is non-decreasing. Let

\begin{equation*}
F_{2}(x)=F_{1}(x)-F(x).
\end{equation*}

\bigskip \noindent Let $a\leq x<y\leq b$. We have%
\begin{equation*}
F_{2}(y)-F_{2}(x)=(F_{1}(y)-F_{1}(x))-(F(y)-F(x)).
\end{equation*}

\bigskip \noindent But $F(y)-F(x)$ is a variation of $F$ over $[x,y]$ with the subdivision $\pi
=(x_{0}=x<x_{1}=y)$. Then $F(y)-F(x)\leq V(F,[x,y]).$ So

\begin{equation*}
F_{2}(y)-F_{2}(x)=(F_{1}(y)-F_{1}(x))-V(F,[x,y]),
\end{equation*}

\bigskip \noindent which is non-negative by\ref{compvar}. So $F_{2}$ is also non-decreasing and%
\begin{equation*}
F=F_{1}-F_{2}.
\end{equation*}

\bigskip \noindent We easily remark that if $F$ is right-continuous, we apply $h\downarrow 0$ to
\begin{equation*}
F(x+h)=F_{1}(x+h)-F_{2}(x+h),
\end{equation*}

\bigskip \noindent to get
\begin{equation*}
F(x)=F_{1}(x+)-F_{2}(x+).
\end{equation*}

\bigskip \noindent where $F_{i}(x+)$ is the right-limit of $F_{i}(x)$ that exists since each $F_{i}$, $i=1,2$, is monotonic. We also know that the functions $x\longrightarrow F_{i}(x+)$, $i=1,2$, are right-continuous. We get

\bigskip \noindent \textbf{Proposition}. Any right-continuous and bounded variation function on $[a,b]$ is difference
of two right-continuous and non-decreasing function.

\bigskip \noindent \textbf{B - Existence of the Riemann-Stieljes integrals}.\\

\bigskip \noindent Now we are going to prove this.\\\

\bigskip \noindent \textbf{Theorem}. Let $F$ be of bounded variation on $]a,b]$. Any continuous $f$  on $]a,b]$ is Riemann-Stieljes is integrable.\\

\bigskip \noindent The proof that follows cannot certainly be given in undergraduate courses. It uses the additivity properties of $\lambda_{F}$ on $\mathcal{I}_{1}$ and next on $\mathcal{C}=a(\mathcal{I}_{1})$. Remark that neither right-continuity nor monotonicity does play a role on finite additivity. The first condition only intervenes for $\sigma$-additivity and the second for the sign of $\lambda_{F}$ and its monotonicity. We do not deal with these concepts in this particular section.\\

\bigskip \noindent Let $\mathcal{F}_{e}$ the class of functions of the forms 
\begin{equation}
h=\sum_{i=0}^{\ell (n)-1}\alpha _{i,n}1_{]x_{i,n},x_{i+1,n}]},  \label{rep1}
\end{equation}

\bigskip \noindent where $\alpha _{i,n}$ is some real number and ($\alpha _{0,n},\alpha
_{1,n}...,\alpha _{m(n),n})=\alpha _{n}$ and $\pi _{n}=(a=x_{0,n}<x_{1,n}<...<x_{\ell(n),n}<x_{\ell(n)+1,n}=b)$ is a partition of $%
]a,b]$. Define for such a function

\begin{equation}
T_{n}(h,\pi _{n},\alpha _{n})=\sum_{i=0}^{\ell (n)}\alpha _{i,n}(F(x_{i+1,n})-F(x_{i,n})). \label{sumRiemann}
\end{equation}

\bigskip \noindent We see that
\begin{equation}
T_{n}(h,\pi _{n},\alpha _{n})=\sum_{i=0}^{\ell (n)}\alpha _{i,n}\lambda
_{F}(]x_{i,n},x_{i+1,n}]).  \label{integ}
\end{equation}

\bigskip \noindent Let $g$ be another element of $\mathcal{F}_{e}$
\begin{equation}
g=\sum_{i=0}^{\ell (r)}\beta _{i,r}1_{]y_{i,n},y_{i+1,n}]}, \label{rep2}
\end{equation}

\bigskip \noindent where $\pi _{r}=(a=y_{0,r}<y_{1,r}<...<y_{\ell (r),r}<y_{\ell (r)+1,r}=b)$ us a partition of ]a,b] and $\beta =(\beta _{0,r},...,\beta _{\ell(r),r})\in R^{\ell (r)+1}.$ We may see that $h$ and $g$ can be written in
\begin{equation*}
h=\sum_{(i,j)\in H}^{{}}\alpha
_{i,n}1_{]x_{i,n},x_{i+1,n}]}1_{]y_{i,n},y_{i+1,n}]},
\end{equation*}

\bigskip \noindent where $H=\{(i,j)\in \{0,...,\ell (n)\}\times \{0,...,\ell (r)\}$, $]x_{i,n},x_{i+1,n}]\cap ]y_{i,n},y_{i+1,n}]\neq \emptyset \}$ and the $]x_{i,n},x_{i+1,n}]\cap ]y_{i,n},y_{i+1,n}]$ (which are of the form $]z_{(i,j)},u_{(i,j)}])$ form a partition of $]a,b]$ for $(i,j)\in H$, and

\begin{equation*}
g=\sum_{(i,j)\in H}^{2}\beta_{i,r}1_{]x_{i,n},x_{i+1,n}]}1_{]y_{i,n},y_{i+1,n}]}.
\end{equation*}

\bigskip \noindent We have to prove that the definition (\ref{integ}) does not depend on the representation (\ref{rep1}) by showing that (\ref{rep2}) gives the same number. By using (\ref{rep2}) we get 

\begin{equation*}
T_{n}(\sum_{(i,j)\in H}^{{}}\alpha _{i,n}1_{]x_{i,n},x_{i+1,n}]\cap
]y_{i,n},y_{i+1,n}]})
\end{equation*}
\begin{equation*}
=\sum_{(i,j)\in H}^{{}}\alpha _{i,n}\lambda _{F}(]x_{i,n},x_{i+1,n}]\cap
]y_{i,n},y_{i+1,n}]).
\end{equation*}

\bigskip \noindent We may extend this the empty sets $]x_{i,n},x_{i+1,n}]\cap ]y_{i,n},y_{i+1,n}]$ whose values by $\lambda _{F}$ are zeros to get

\begin{eqnarray*}
&&\sum_{i=0}^{\ell (n)}\sum_{j=0}^{\ell (r)}\alpha _{i,n}\lambda _{F}(]x_{i,n},x_{i+1,n}]\cap ]y_{i,n},y_{i+1,n}])\\
&=&\sum_{i=0}^{\ell (n)}\alpha _{i,n}\left\{ \sum_{j=0}^{\ell (r)}\lambda_{F}(]x_{i,n},x_{i+1,n}]\cap ]y_{i,n},y_{i+1,n}])\right\}\\
&=&\sum_{i=0}^{m(n)}\alpha _{i,n}\left\{ \lambda _{F}(]x_{i,n},x_{i+1,n}]\cap \left\{ \sum_{j=0}^{m(r)}]y_{i,n},y_{i+1,n}]\right\} )\right\}\\
&=&\sum_{i=0}^{\ell (n)}\alpha _{i,n}\lambda _{F}(]x_{i,n},x_{i+1,n}]\cap]a,b]=T_{n}(h,\pi _{n},\alpha _{n}).
\end{eqnarray*}

\bigskip \noindent Before we turn back to the Riemann-Stieljes integration, just remark that Riemann-Stieljes sums satisfy

\begin{equation*}
S_{n}(f,F,a,b,\pi _{n},c_{n})=T_{n}(f_{n},\pi _{n},d_{n}(f)),
\end{equation*}

\bigskip \noindent where
\begin{equation*}
d_{n}(f)=(f(c_{0,n}),...,f(c_{\ell (n),n}).
\end{equation*}

\bigskip \noindent With this notations, we may show that sequence $S_{n}(f,F,a,b,\pi _{n},c_{n})
$ is Cauchy.\\

\bigskip Suppose that $F$ is of bounded variation. Consider two sequences of Riemann-Stieljes sums, $S_{n}(f,F,a,b,\pi _{n},c_{n})$ and $S_{r}(f,F,a,b,\delta _{r},e_{r}),$ where
we use the partitions $\pi _{n}$ and $\delta _{r}$ already defined and $e_{r}=(e_{0,r},...,e_{\ell (r),r})$. We have 
\begin{equation*}
S_{n}(f,F,a,b,\pi _{n},c_{n})=\sum_{(i,j)\in H}^{{}}f(c_{i,n})\lambda_{F}(]x_{i,n},x_{i+1,n}]\cap ]y_{i,n},y_{i+1,n}])
\end{equation*}

\bigskip \noindent and
\begin{equation*}
S_{r}(f,F,a,b,\delta _{r},e_{r})=\sum_{(i,j)\in H}^{{}}f(e_{i,n})\lambda_{F}(]x_{i,n},x_{i+1,n}]\cap ]y_{i,n},y_{i+1,n}]).
\end{equation*}

\bigskip \noindent We get

\begin{equation*}
\left\vert S_{n}(f,F,a,b,\pi _{n},c_{n})-S_{r}(f,F,a,b,\delta
_{r},e_{r})\right\vert 
\end{equation*}%
\begin{equation*}
\leq \sum_{(i,j)\in H}^{{}}\left\vert f(e_{i,n})-f(c_{i,n})\right\vert
\left\vert \lambda _{F}(]x_{i,n},x_{i+1,n}]\cap
]y_{i,n},y_{i+1,n}]\right\vert ).
\end{equation*}

\bigskip \noindent But both $f(e_{i,n})$ and $f(c_{i,n})$ lie in $]x_{i,n},x_{i+1,n}]\cap]y_{i,n},y_{i+1,n}]$ so that for any 
\begin{equation*}
\left\vert f(e_{i,n})-f(c_{i,n})\right\vert \leq \sup_{\left\vert
x-y\right\vert \leq \max (m(\pi _{n}),m(\delta _{r}))}\left\vert
f(x)-f(y)\right\vert =\gamma (f,\varepsilon (n,r)),
\end{equation*}

\bigskip \noindent with $\varepsilon (n,r)=\max (m(\pi _{n}),m(\delta _{r})$. Hence, we have 

\begin{eqnarray*}
&&\left\vert S_{n}(f,F,a,b,\pi _{n},c_{n})-S_{r}(f,F,a,b,\delta _{r},e_{r})\right\vert\\
&\leq&\gamma (f,\varepsilon (n,r))\sum_{(i,j)\in H}^{{}}\left\vert \lambda _{F}(]x_{i,n},x_{i+1,n}]\cap ]y_{i,n},y_{i+1,n}])\right\vert .
\end{eqnarray*}

\bigskip \noindent Since $\sum_{(i,j)\in H}^{{}}\left\vert \lambda _{F}(]x_{i,n},x_{i+1,n}]\cap]y_{i,n},y_{i+1,n}])\right\vert $ is a variation of $F$ on $]a,b]$, we get%
\begin{equation*}
=\gamma (f,\varepsilon (n,r))\times V(F,]a,b]).
\end{equation*}

\bigskip \noindent If $f$ is continuous on $[a,b]$, it is uniformly continuous on it so that $\gamma (f,\varepsilon (n,r))\rightarrow 0$ as $(n,r)\rightarrow (+\infty
,+\infty )$ since $\max (m(\pi _{n}),m(\delta _{r}))\rightarrow 0$ as $(n,r)\rightarrow (+\infty ,+\infty ).$ We conclude that the sequence $%
S_{n}(f,F,a,b,\pi _{n},c_{n})$ is Cauchy when $m(\pi _{n})\rightarrow 0.$ Then it has a unique limit in $\mathbb{R}$.\\

%\include{10_lsm/mes_doc_10_04_ap}
%\include{10_lsm/mes_doc_10_05_ap}

%\chapter{Appendix}
\part{Appendix, Conclusion and Bibliography}
\chapter{General Appendix} \label{11_appendix}

\noindent \textbf{Content of the Chapter}

\begin{table}[htbp]
	\centering
		\begin{tabular}{llll}
		\hline
		Type& Name & Title  & page\\
		SD& Doc 11-01 &  What should not be ignored on limits in $\overline{\mathbb{R}}$ & \pageref{doc11-01}\\
		\hline
		\end{tabular}
\end{table}

\noindent \LARGE \textbf{DOC 11-01 : What should not be ignored on limits in $\overline{\mathbb{R}}$ - Exercises with Solutions} \label{doc11-01}\\
\bigskip
\Large

\noindent \textbf{Definition} $\ell \in \overline{\mathbb{R}}$ is an accumulation point of a sequence 
 $(x_{n})_{n\geq 0}$ of real numbers finite or infinite, in $\overline{\mathbb{R}}$, if and only if there exists a sub-sequence $(x_{n(k)})_{k\geq 0}$ of
 $(x_{n})_{n\geq 0}$ such that $%
x_{n(k)}$ converges to $\ell $, as $k\rightarrow +\infty $.\newline

\noindent \textbf{Exercise 1}.\\

\noindent  Set $y_{n}=\inf_{p\geq n}x_{p}$ and $z_{n}=\sup_{p\geq n}x_{p} $ for all $n\geq 0$. Show that :\newline

\noindent \textbf{(1)} $\forall n\geq 0,y_{n}\leq x_{n}\leq z_{n}$.\newline

\noindent \textbf{(2)} Justify the existence of the limit of $y_{n}$ called limit inferior of the sequence $(x_{n})_{n\geq 0}$, denoted by $%
\liminf x_{n}$ or $\underline{\lim }$ $x_{n},$ and that it is equal to the following%
\begin{equation*}
\underline{\lim }\text{ }x_{n}=\lim \inf x_{n}=\sup_{n\geq 0}\inf_{p\geq
n}x_{p}.
\end{equation*}

\noindent \textbf{(3)} Justify the existence of the limit of $z_{n}$ called limit superior of the sequence $(x_{n})_{n\geq 0}$ denoted by $%
\lim \sup x_{n}$ or $\overline{\lim }$ $x_{n},$ and that it is equal%
\begin{equation*}
\overline{\lim }\text{ }x_{n}=\lim \sup x_{n}=\inf_{n\geq 0}\sup_{p\geq
n}x_{p}x_{p}.
\end{equation*}

\bigskip

\noindent \textbf{(4)} Establish that 
\begin{equation*}
-\liminf x_{n}=\limsup (-x_{n})\noindent \text{ \ \ and \ }-\limsup
x_{n}=\liminf (-x_{n}).
\end{equation*}

\newpage \noindent \textbf{(5)} Show that the limit superior is sub-additive and the limit inferior is super-additive, i.e. :  for two sequences
$(s_{n})_{n\geq 0}$ and $(t_{n})_{n\geq 0}$ 
\begin{equation*}
\limsup (s_{n}+t_{n})\leq \limsup s_{n}+\limsup t_{n}
\end{equation*}

\noindent and
\begin{equation*}
\lim \inf (s_{n}+t_{n})\geq \lim \inf s_{n}+\lim \inf t_{n}.
\end{equation*}

\noindent \textbf{(6)} Deduce from (1) that if%
\begin{equation*}
\lim \inf x_{n}=\lim \sup x_{n},
\end{equation*}%
then $(x_{n})_{n\geq 0}$ has a limit and 
\begin{equation*}
\lim x_{n}=\lim \inf x_{n}=\lim \sup x_{n}
\end{equation*}

\bigskip

\noindent \textbf{Exercise 2.} Accumulation points of $(x_{n})_{n\geq 0}$.\newline

\noindent \textbf{(a)} Show that if $\ell _{1}$=$\lim \inf x_{n}$ and $\ell
_{2}=\lim \sup x_{n}$ are accumulation points of $(x_{n})_{n\geq 0}.
$ Show one case and deduce the second one and by using Point (3) of Exercise 1.\newline

\noindent \textbf{(b)} Show that $\ell _{1}$ is the smallest accumulation point of $(x_{n})_{n\geq 0}$ and $\ell _{2}$ is the biggest.
(Similarly, show one case and deduce the second one and by using Point (3) of Exercise 1).\newline

\noindent \textbf{(c)} Deduce from (a) that if $(x_{n})_{n\geq 0}$ has a limit $\ell$,  then it is equal to the unique accumulation point and so,%
\begin{equation*}
\ell =\overline{\lim }\text{ }x_{n}=\lim \sup x_{n}=\inf_{n\geq
0}\sup_{p\geq n}x_{p}.
\end{equation*}

\noindent \textbf{(d)} Combine this result with Point \textbf{(6)} of Exercise 1 to show that a sequence $(x_{n})_{n\geq 0}$ of $\overline{\mathbb{R}}
$ has a limit $\ell $ in $\overline{\mathbb{R}}$ if and only if\ $\lim \inf
x_{n}=\lim \sup x_{n}$ and then%
\begin{equation*}
\ell =\lim x_{n}=\lim \inf x_{n}=\lim \sup x_{n}.
\end{equation*}

\newpage

\noindent \textbf{Exercise 3. } Let $(x_{n})_{n\geq 0}$ be a non-decreasing sequence
of $\overline{\mathbb{R}}$. Study its limit superior and its limit inferior and deduce that%
\begin{equation*}
\lim x_{n}=\sup_{n\geq 0}x_{n}.
\end{equation*}

\noindent Deduce that for a non-increasing sequence $(x_{n})_{n\geq 0}$
of $\overline{\mathbb{R}},$%
\begin{equation*}
\lim x_{n}=\inf_{n\geq 0}x_{n}.
\end{equation*}

\bigskip

\noindent \textbf{Exercise 4.} (Convergence criteria)\newline

\noindent \textbf{Prohorov Criterion} Let $(x_{n})_{n\geq 0}$ be a sequence of $\overline{%
\mathbb{R}}$ and a real number $\ell \in \overline{\mathbb{R}}$ such that: Every subsequence of $(x_{n})_{n\geq 0}$ 
also has a subsequence ( that is a subssubsequence of $(x_{n})_{n\geq 0}$ ) that converges to $\ell .$
Then, the limit of $(x_{n})_{n\geq 0}$ exists and is equal $\ell .$\newline

\noindent \textbf{Upcrossing or Downcrossing Criterion}. \newline

\noindent Let $(x_{n})_{n\geq 0}$ be a sequence in $\overline{\mathbb{R}}$ and two real numbers $a$ and $b$ such that $a<b.$
We define%
\begin{equation*}
\nu _{1}=\left\{ 
\begin{array}{cc}
\inf  & \{n\geq 0,x_{n}<a\} \\ 
+\infty  & \text{if (}\forall n\geq 0,x_{n}\geq a\text{)}%
\end{array}%
\right. .
\end{equation*}%
If $\nu _{1}$ is finite, let%
\begin{equation*}
\nu _{2}=\left\{ 
\begin{array}{cc}
\inf  & \{n>\nu _{1},x_{n}>b\} \\ 
+\infty  & \text{if (}n>\nu _{1},x_{n}\leq b\text{)}%
\end{array}%
\right. .
\end{equation*}%
.

\noindent As long as the $\nu _{j}'s$ are finite, we can define for $\nu
_{2k-2}(k\geq 2)$

\begin{equation*}
\nu _{2k-1}=\left\{ 
\begin{array}{cc}
\inf  & \{n>\nu _{2k-2},x_{n}<a\} \\ 
+\infty  & \text{if (}\forall n>\nu _{2k-2},x_{n}\geq a\text{)}%
\end{array}%
\right. .
\end{equation*}%
and for $\nu _{2k-1}$ finite, 
\begin{equation*}
\nu _{2k}=\left\{ 
\begin{array}{cc}
\inf  & \{n>\nu _{2k-1},x_{n}>b\} \\ 
+\infty  & \text{if (}n>\nu _{2k-1},x_{n}\leq b\text{)}%
\end{array}%
\right. .
\end{equation*}

\noindent We stop once one $\nu _{j}$ is $+\infty$. If $\nu
_{2j}$ is finite, then 
\begin{equation*}
x_{\nu _{2j}}-x_{\nu _{2j-1}}>b-a. 
\end{equation*}

\noindent We then say : by that moving from $x_{\nu _{2j-1}}$ to $x_{\nu
_{2j}},$ we have accomplished a crossing (toward the up) of the segment $[a,b]$
called \textit{up-crossings}. Similarly, if one $\nu _{2j+1}$
is finite, then the segment $[x_{\nu _{2j}},x_{\nu _{2j+1}}]$ is a crossing downward (down-crossing) of the segment $[a,b].$ Let%
\begin{equation*}
D(a,b)=\text{ number of up-crossings of the sequence of the segment }[a,b]\text{.}
\end{equation*}

\bigskip

\noindent \textbf{(a)} What is the value of $D(a,b)$ if \ $\nu _{2k}$ is finite and $\nu
_{2k+1}$ infinite.\newline

\noindent \textbf{(b)} What is the value of $D(a,b)$ if \ $\nu _{2k+1}$ is finite and $\nu
_{2k+2}$ infinite.\newline

\noindent \textbf{(c)} What is the value of $D(a,b)$ if \ all the $\nu _{j}'s$ are finite.%
\newline

\noindent \textbf{(d)} Show that $(x_{n})_{n\geq 0}$ has a limit iff
for all $a<b,$ $D(a,b)<\infty.$\newline

\noindent \textbf{(e)} Show that $(x_{n})_{n\geq 0}$ has a limit iff
for all $a<b,$ $(a,b)\in \mathbb{Q}^{2},D(a,b)<\infty .$\newline

\bigskip

\noindent \textbf{Exercise 5. } (Cauchy Criterion). Let $%
(x_{n})_{n\geq 0}$ $\mathbb{R}$ be a sequence of (\textbf{real numbers}).\newline

\noindent \textbf{(a)} Show that if $(x_{n})_{n\geq 0}$ is Cauchy,
then it has a unique accumulation point $\ell \in 
\mathbb{R}$ which is its limit.\newline

\noindent \textbf{(b)} Show that if a sequence $(x_{n})_{n\geq 0}\subset 
\mathbb{R}$ \ converges to $\ell \in \mathbb{R},$ then, it is Cauchy.%
\newline

\noindent \textbf{(c)} Deduce the Cauchy criterion for sequences of real numbers.

\newpage

\begin{center}
\textbf{SOLUTIONS}
\end{center}

\noindent \textbf{Exercise 1}.\newline

\noindent \textbf{Question (1)}. It is obvious that :%
\begin{equation*}
\underset{p\geq n}{\inf }x_{p}\leq x_{n}\leq \underset{p\geq n}{\sup }x_{p},
\end{equation*}

\noindent since $x_{n}$ is an element of $\left\{
x_{n},x_{n+1},...\right\} $ on which we take the supremum or the infimum.%
\newline

\noindent \textbf{Question (2)}. Let $y_{n}=\underset{p\geq 0}{\inf }%
x_{p}=\underset{p\geq n}{\inf }A_{n},$ where $A_{n}=\left\{
x_{n},x_{n+1},...\right\} $ is a non-increasing sequence of sets : $\forall n\geq 0$,
\begin{equation*}
A_{n+1}\subset A_{n}.
\end{equation*}

\noindent So the infimum on $A_{n}$ increases. If $y_{n}$ increases in $%
\overline{\mathbb{R}},$ its limit is its upper bound, finite or infinite. So%
\begin{equation*}
y_{n}\nearrow \underline{\lim }\text{ }x_{n},
\end{equation*}

\noindent is a finite or infinite number.\newline

\noindent \textbf{Question (3)}. We also show that $z_{n}=\sup A_{n}$ decreases and $z_{n}\downarrow \overline{\lim }$ $x_{n}$.\newline

\noindent \textbf{Question (4) \label{qst4}}. We recall that 
\begin{equation*}
-\sup \left\{ x,x\in A\right\} =\inf \left\{ -x,x\in A\right\}, 
\end{equation*}

\noindent which we write 
\begin{equation*}
-\sup A=\inf (-A).
\end{equation*}

\noindent Thus,

\begin{equation*}
-z_{n}=-\sup A_{n}=\inf (-A_{n}) = \inf \left\{-x_{p},p\geq n\right\}.
\end{equation*}

\noindent The right hand term tends to $-\overline{\lim}\ x_{n}$ and the left hand to $\underline{\lim} (-x_{n})$ and so 

\begin{equation*}
-\overline{\lim}\ x_{n}=\underline{\lim }\ (-x_{n}).
\end{equation*}

\bigskip \noindent Similarly, we show:
\begin{equation*}
-\underline{\lim } \ (x_{n})=\overline{\lim} \ (-x_{n}).
\end{equation*}

\noindent 

\noindent \textbf{Question (5)}. These properties come from the formulas, where $A\subseteq \mathbb{R},B\subseteq \mathbb{R}$ :%
\begin{equation*}
\sup \left\{ x+y,A\subseteq \mathbb{R},B\subseteq \mathbb{R}\right\} \leq
\sup A+\sup B.
\end{equation*}

\noindent In fact : 
\begin{equation*}
\forall x\in \mathbb{R},x\leq \sup A
\end{equation*}

\noindent and
\begin{equation*}
\forall y\in \mathbb{R},y\leq \sup B.
\end{equation*}

\noindent Thus 
\begin{equation*}
x+y\leq \sup A+\sup B,
\end{equation*}

\noindent where 
\begin{equation*}
\underset{x\in A,y\in B}{\sup }x+y\leq \sup A+\sup B.
\end{equation*}%
Similarly,%
\begin{equation*}
\inf (A+B\geq \inf A+\inf B.
\end{equation*}

\noindent In fact :

\begin{equation*}
\forall (x,y)\in A\times B,x\geq \inf A\text{ and }y\geq \inf B.
\end{equation*}

\noindent Thus 
\begin{equation*}
x+y\geq \inf A+\inf B,
\end{equation*}

\noindent and so
\begin{equation*}
\underset{x\in A,y\in B}{\inf }(x+y)\geq \inf A+\inf B
\end{equation*}

\noindent \textbf{Application}.\newline

\begin{equation*}
\underset{p\geq n}{\sup } \ (x_{p}+y_{p})\leq \underset{p\geq n}{\sup } \ x_{p}+\underset{p\geq n}{\sup } \ y_{p}.
\end{equation*}

\noindent All these sequences are non-increasing. By taking the infimum, we obtain the limits superior :

\begin{equation*}
\overline{\lim }\text{ }(x_{n}+y_{n})\leq \overline{\lim }\text{ }x_{n}+%
\overline{\lim }\text{ }x_{n}.
\end{equation*}

\bigskip

\noindent \textbf{Question (6)}. Set

\begin{equation*}
\underline{\lim } \ x_{n}=\overline{\lim } \ x_{n}.
\end{equation*}

\noindent Since : 
\begin{equation*}
\forall x\geq 1,\text{ }y_{n}\leq x_{n}\leq z_{n},
\end{equation*}%

\begin{equation*}
y_{n}\rightarrow \underline{\lim} \ x_{n}
\end{equation*}%

\noindent and 

\begin{equation*}
z_{n}\rightarrow \overline{\lim } \ x_{n},
\end{equation*}

\noindent we apply the Sandwich Theorem to conclude that the limit of $x_{n}$ exists and :

\begin{equation*}
\lim \text{ }x_{n}=\underline{\lim }\text{ }x_{n}=\overline{\lim }\text{ }%
x_{n}.
\end{equation*}

\bigskip 
\noindent \textbf{Exercice 2}.\newline

\noindent \textbf{Question (a).}\\

\noindent Thanks to Question (4) of Exercise 1, it suffices to show this property for one of the limits. Consider the limit superior and the three cases:\\

\noindent \textbf{The case of a finite limit superior} :

\begin{equation*}
\underline{\lim} x_{n}=\ell \text{ finite.}
\end{equation*}

\noindent By definition, 
\begin{equation*}
z_{n}=\underset{p\geq n}{\sup }x_{p}\downarrow \ell .
\end{equation*}

\noindent So: 
\begin{equation*}
\forall \varepsilon >0,\exists (N(\varepsilon )\geq 1),\forall p\geq
N(\varepsilon ),\ell -\varepsilon <x_{p}\leq \ell +\varepsilon .
\end{equation*}

\noindent Take less than that:

\begin{equation*}
\forall \varepsilon >0,\exists n_{\varepsilon }\geq 1:\ell -\varepsilon
<x_{n_{\varepsilon }}\leq \ell +\varepsilon.
\end{equation*}

\noindent We shall construct a sub-sequence converging to $\ell$.\\

\noindent Let $\varepsilon =1:$%
\begin{equation*}
\exists N_{1}:\ell -1<x_{N_{1}}=\underset{p\geq n}{\sup }x_{p}\leq \ell +1.
\end{equation*}

\noindent But if 
\begin{equation}
z_{N_{1}}=\underset{p\geq n}{\sup }x_{p}>\ell -1, \label{cc}
\end{equation}

\noindent there surely exists an $n_{1}\geq N_{1}$ such that%
\begin{equation*}
x_{n_{1}}>\ell -1.
\end{equation*}

\noindent If not, we would have 
\begin{equation*}
( \forall p\geq N_{1},x_{p}\leq \ell -1\ ) \Longrightarrow \sup \left\{
x_{p},p\geq N_{1}\right\} =z_{N_{1}}\geq \ell -1,
\end{equation*}

\noindent which is contradictory with (\ref{cc}). So, there exists $n_{1}\geq N_{1}$ such that
\begin{equation*}
\ell -1<x_{n_{1}}\leq \underset{p\geq N_{1}}{\sup }x_{p}\leq \ell -1.
\end{equation*}

\noindent i.e.

\begin{equation*}
\ell -1<x_{n_{1}}\leq \ell +1.
\end{equation*}

\noindent We move to step $\varepsilon =\frac{1}{2}$ and we consider the sequence%
 $(z_{n})_{n\geq n_{1}}$ whose limit remains $\ell$. So, there exists $N_{2}>n_{1}:$%
\begin{equation*}
\ell -\frac{1}{2}<z_{N_{2}}\leq \ell -\frac{1}{2}.
\end{equation*}

\noindent We deduce like previously that $n_{2}\geq N_{2}$ such that%
\begin{equation*}
\ell -\frac{1}{2}<x_{n_{2}}\leq \ell +\frac{1}{2}
\end{equation*}

\noindent with $n_{2}\geq N_{1}>n_{1}$.\\

\noindent Next, we set $\varepsilon =1/3,$ there will exist $N_{3}>n_{2}$ such that%
\begin{equation*}
\ell -\frac{1}{3}<z_{N_{3}}\leq \ell -\frac{1}{3}
\end{equation*}

\noindent and we could find an $n_{3}\geq N_{3}$ such that%

\begin{equation*}
\ell -\frac{1}{3}<x_{n_{3}}\leq \ell -\frac{1}{3}.
\end{equation*}

\noindent Step by step, we deduce the existence of $%
x_{n_{1}},x_{n_{2}},x_{n_{3}},...,x_{n_{k}},...$ with $n_{1}<n_{2}<n_{3}%
\,<...<n_{k}<n_{k+1}<...$ such that

$$
\forall k\geq 1, \ell -\frac{1}{k}<x_{n_{k}}\leq \ell -\frac{1}{k},
$$

\noindent i.e.

\begin{equation*}
\left\vert \ell -x_{n_{k}}\right\vert \leq \frac{1}{k},
\end{equation*}

\noindent which will imply: 
\begin{equation*}
x_{n_{k}}\rightarrow \ell 
\end{equation*}

\noindent Conclusion : $(x_{n_{k}})_{k\geq 1}$ is very well a subsequence since $n_{k}<n_{k+1}$ for all $k \geq 1$ 
and it converges to $\ell$, which is then an accumulation point.\\

\noindent \textbf{Case of the limit superior equal $+\infty$} : 
$$
\overline{\lim} \text{ } x_{n}=+\infty.
$$
\noindent Since $z_{n}\uparrow +\infty ,$ we have : $\forall k\geq 1,\exists
N_{k}\geq 1,$ 
\begin{equation*}
z_{N_{k}}\geq k+1.
\end{equation*}

\noindent For $k=1$, let $z_{N_{1}}=\underset{p\geq N_{1}}{\inf }%
x_{p}\geq 1+1=2.$ So there exists 
\begin{equation*}
n_{1}\geq N_{1}
\end{equation*}

\noindent such that :
\begin{equation*}
x_{n_{1}}\geq 1.
\end{equation*}

\noindent For $k=2$, consider the sequence $(z_{n})_{n\geq n_{1}+1}.$
We find in the same manner 
\begin{equation*}
n_2 \geq n_{1}+1
\end{equation*}%
\noindent and 
\begin{equation*}
x_{n_{2}}\geq 2.
\end{equation*}

\noindent Step by step, we find for all $k\geq 3$, an $n_{k}\geq n_{k-1}+1$ such that
\begin{equation*}
x_{n_{k}}\geq k,
\end{equation*}

\noindent which leads to $x_{n_{k}}\rightarrow +\infty $ as $k\rightarrow +\infty $.\\

\noindent \textbf{Case of the limit superior equal $-\infty$} : 

$$
\overline{\lim }x_{n}=-\infty.
$$

\noindent This implies : $\forall k\geq 1,\exists N_{k}\geq 1,$ such that%
\begin{equation*}
z_{n_{k}}\leq -k.
\end{equation*}

\noindent For $k=1$, there exists $n_{1}$ such that%
\begin{equation*}
z_{n_{1}}\leq -1.
\end{equation*}
But 
\begin{equation*}
x_{n_{1}}\leq z_{n_{1}}\leq -1.
\end{equation*}

\noindent Let $k=2$. Consider $\left( z_{n}\right) _{n\geq
n_{1}+1}\downarrow -\infty .$ There will exist $n_{2}\geq n_{1}+1:$%
\begin{equation*}
x_{n_{2}}\leq z_{n_{2}}\leq -2
\end{equation*}

\noindent Step by step, we find $n_{k1}<n_{k+1}$ in such a way that $x_{n_{k}}<-k$ for all $k$ bigger than $1$. So
\begin{equation*}
x_{n_{k}}\rightarrow +\infty 
\end{equation*}

\bigskip

\noindent \textbf{Question (b).}\\

\noindent Let $\ell$ be an accumulation point of $(x_n)_{n \geq 1}$, the limit of one of its sub-sequences $(x_{n_{k}})_{k \geq 1}$. We have

$$
y_{n_{k}}=\inf_{p\geq n_k} \ x_p \leq x_{n_{k}} \leq  \sup_{p\geq n_k} \ x_p=z_{n_{k}}.
$$

\noindent The left hand side term is a sub-sequence of $(y_n)$ tending to the limit inferior and the right hand side is a 
sub-sequence of $(z_n)$ tending to the limit superior. So we will have:

$$
\underline{\lim} \ x_{n} \leq \ell \leq \overline{\lim } \ x_{n},
$$

\noindent which shows that $\underline{\lim} \ x_{n}$ is the smallest accumulation point and $\overline{\lim } \ x_{n}$ is the largest.\\

\noindent \textbf{Question (c).} If the sequence $(x_n)_{n \geq 1}$ has a limit $\ell$, it is the limit of all its sub-sequences,
so subsequences tending to the limits superior and inferior. Which answers question (b).\\

\noindent \textbf{Question (d).} We answer this question by combining point (d) of this exercise and Point 6) of the Exercise 1.\\

\noindent \textbf{Exercise 3}. Let $(x_{n})_{n\geq 0}$ be a non-decreasing sequence, we have:%
\begin{equation*}
z_{n}=\underset{p\geq n}{\sup} \ x_{p}=\underset{p\geq 0}{\sup} \ x_{p},\forall
n\geq 0.
\end{equation*}

\noindent Why? Because by increasingness,%
\begin{equation*}
\left\{ x_{p},p\geq 0\right\} =\left\{ x_{p},0\leq p\leq n-1\right\} \cup
\left\{ x_{p},p\geq n\right\}.
\end{equation*}

\bigskip

\noindent Since all the elements of $\left\{ x_{p},0\leq p\leq
n-1\right\} $ are smaller than than those of $\left\{ x_{p},p\geq n\right\} ,$
the supremum is achieved on $\left\{ x_{p},p\geq n\right\} $ and so 
\begin{equation*}
\ell =\underset{p\geq 0}{\sup } \ x_{p}=\underset{p\geq n}{\sup }x_{p}=z_{n}.
\end{equation*}

\noindent Thus
\begin{equation*}
z_{n}=\ell \rightarrow \ell .
\end{equation*}

\noindent We also have $y_n=\inf \left\{ x_{p},0\leq p\leq n\right\}=x_n$, which is a non-decreasing sequence and so converges to
$\ell =\underset{p\geq 0}{\sup } \ x_{p}$. \\

\bigskip

\noindent \textbf{Exercise 4}.\\

\noindent Let $\ell \in \overline{\mathbb{R}}$ having the indicated property. Let $\ell ^{\prime }$ be a given accumulation point.%
\begin{equation*}
 \left( x_{n_{k}}\right)_{k\geq 1} \subseteq \left( x_{n}\right) _{n\geq 0}%
\text{ such that }x_{n_{K}}\rightarrow \ell ^{\prime}.
\end{equation*}

\noindent By hypothesis this sub-sequence $\left( x_{n_{K}}\right) $
has in turn a sub-sub-sequence $\left( x_{n_{\left( k(p)\right) }}\right)_{p\geq 1} $ such that $x_{n_{\left( k(p)\right) }}\rightarrow
\ell $ as $p\rightarrow +\infty $.\newline

\noindent But as a sub-sequence of $\left( x_{n_{\left( k\right)
}}\right) ,$ 
\begin{equation*}
x_{n_{\left( k(\ell )\right) }}\rightarrow \ell ^{\prime }.
\end{equation*}%
Thus
\begin{equation*}
\ell =\ell ^{\prime}.
\end{equation*}

\noindent Applying that to the limit superior and limit inferior, we have:%
\begin{equation*}
\overline{\lim} \ x_{n}=\underline{\lim}\ x_{n}=\ell.
\end{equation*}

\noindent And so $\lim x_{n}$ exists and equals $\ell$.\\

\noindent \textbf{Exercise 5}.\\

\noindent \textbf{Question (a)}. If $\nu _{2k}$ is finite and if $\nu _{2k+1}$ is infinite, then there ate exactly $k$ up-crossings : 
$[x_{\nu_{2j-1}},x_{\nu _{2j}}]$, $j=1,...,k$, that is, we have $D(a,b)=k$.\\

\noindent \textbf{Question (b)}. If $\nu _{2k+1}$ is finite and $\nu _{2k+2}$ is infinite, then there are exactly $k$ up-crossings:
$[x_{\nu_{2j-1}},x_{\nu_{2j}}]$, $j=1,...,k$, that is we have $D(a,b)=k$.\\

\noindent \textbf{Question (c)}. If all the $\nu_{j}'s$ are finite, then there are an infinite number of up-crossings : 
$[x_{\nu_{2j-1}},x_{\nu_{2j}}]$, $j\geq 1k$ : $D(a,b)=+\infty$.\\

\noindent \textbf{Question (d)}. Suppose that there exist $a < b$ rationals such that $D(a,b)=+\infty$. 
Then all the $\nu _{j}'s$ are finite. The subsequence $x_{\nu_{2j-1}}$ is strictly below $a$. 
So its limit inferior is below $a$. This limit inferior is an accumulation point of the sequence $(x_n)_{n\geq 1}$, 
so is more than $\underline{\lim}\ x_{n}$, which is below $a$.\\

\noindent Similarly, the subsequence $x_{\nu_{2j}}$ is strictly below $b$. So the limit superior is above $a$. 
This limit superior is an accumulation point of the sequence $(x_n)_{n\geq 1}$, so it is below $\overline{\lim}\ x_{n}$, 
which is directly above $b$. This leads to :
$$
\underline{\lim}\ x_{n} \leq a < b \leq \overline{\lim}\ x_{n}. 
$$

\noindent That implies that the limit of $(x_n)$ does not exist. In contrary, we just proved that the limit of $(x_n)$ exists, 
meanwhile for all the real numbers $a$ and $b$ such that $a<b$, $D(a,b)$ is finite.\\

\noindent Now, suppose that the limit of $(x_n)$ does not exist. Then,

$$
\underline{\lim}\ x_{n} < \overline{\lim}\ x_{n}. 
$$

\noindent We can then find two rationals $a$ and $b$ such that $a<b$ and a number $\epsilon$ such that $0<\epsilon$, such that 

$$
\underline{\lim}\ x_{n} < a-\epsilon < a < b < b+\epsilon <  \overline{\lim}\ x_{n}. 
$$

\noindent If $\underline{\lim}\ x_{n} < a-\epsilon$, we can return to Question \textbf{(a)} of Exercise 2 and construct a sub-sequence of $(x_n)$
which tends to $\underline{\lim}\ x_{n}$ while remaining below $a-\epsilon$. Similarly, if $b+\epsilon < \overline{\lim}\ x_{n}$, 
we can create a sub-sequence of $(x_n)$ which tends to $\overline{\lim}\ x_{n}$ while staying above $b+\epsilon$. 
It is evident with these two sequences that we could define with these two sequences all $\nu_{j}$ finite and so $D(a,b)=+\infty$.\\

\noindent We have just shown by contradiction that if all the $D(a,b)$ are finite for all rationals $a$ and $b$ such that $a<b$, 
then, the limit of $(x_n)_{n\geq 0}$ exists.\\

\noindent \textbf{Exercise 5}. Cauchy criterion in $\mathbb{R}$.\\

\noindent Suppose that the sequence is Cauchy, $i.e.$,
$$
\lim_{(p,q)\rightarrow (+\infty,+\infty)} \ (x_p-x_q)=0.
$$

\noindent Then let $x_{n_{k,1}}$ and $x_{n_{k,2}}$ be two sub-sequences converging respectively to $\ell_1=\underline{\lim}\ x_{n}$ and $\ell_2=\overline{\lim}\ x_{n}$. So

$$
\lim_{(p,q)\rightarrow (+\infty,+\infty)} \ (x_{n_{p,1}}-x_{n_{q,2}})=0.
$$

\noindent, By first letting $p\rightarrow +\infty$, we have

$$
\lim_{q\rightarrow +\infty} \ \ell_1-x_{n_{q,2}}=0,
$$

\noindent which shows that $\ell_1$ is finite, else $\ell_1-x_{n_{q,2}}$ would remain infinite and would not tend to $0$. 
By interchanging the roles of $p$ and $q$, we also have that $\ell_2$ is finite.\\

\noindent Finally, by letting $q\rightarrow +\infty$, in the last equation, we obtain
$$
\ell_1=\underline{\lim}\ x_{n}=\overline{\lim}\ x_{n}=\ell_2.
$$

\noindent which proves the existence of the finite limit of the sequence $(x_n)$.\\

\noindent Now suppose that the finite limit $\ell$ of $(x_n)$ exists. Then

$$
\lim_{(p,q)\rightarrow (+\infty,+\infty)} \ (x_p-x_q)=\ell-\ell=0,
$$0
\noindent which shows that the sequence is Cauchy.\\

\chapter{Conclusion and Bibliography} \label{12_conbib}

\noindent \textbf{Content of the Chapter}

\begin{table}[htbp]
	\centering
		\begin{tabular}{llll}
		\hline
		Type& Name & Title  & page\\
		S& Doc 12-01 & Conclusion & \pageref{doc12-01}\\
		R&  & Bibliography  & \pageref{doc12-02} \\
		\hline
		\end{tabular}
\end{table}

\newpage
\noindent \LARGE \textbf{DOC 12-01 : Conclusion} \label{doc12-01}\\
\bigskip
\Large

\noindent The learner who went through all the material is expected to have covered a serious and solid part of measure theory. In the lines below, we will see where this takes him.\\

\noindent (1) Without express it properly, a large part of fundamental Probability theory is already acquired with the proper terminology : Normed measures are probability measures, measurable functions are random variables, measurable sets are events, almost everywhere (\textit{a.e.}) properties are  almost-sure (\textit{a.s.}) ones. Here are some important facts.\\

\noindent (2) The probability law of a random element is the image-measure of the underlying probability measure by this random element.\\

\noindent (3) The mathematical expectation with respect to a probability measure of a real-valued random variable is nothing else than its integral if it exists.\\

\noindent (4) The independence of two random elements means that their joint probability laws is the product measure of their probability laws.\\

\noindent (5) The Radon-Nikodym's Theorem and the results of $L^p$ spaces, as well as the handling of absolutely convergent series as integral with respect to counting measures are very useful for the learner in the course of Functional Analysis. If this course has been taken before Measure Theory, a re-reading of the first course ensure him a deeper understanding.\\

\noindent (6) Also, the The Radon-Nikodym's Theorem is the more general tool for dealing with the conditional expectation of a real-valued $X$ random variable with respect to a sub-$\sigma$-algebra$\mathcal{B}$. It is simply the Radon-Nikodym derivative of the indefinite integral with respect to $X$ defined on $\mathcal{B}$.\\

\noindent The theory of mathematical expectation is the general frame to deal with the main dependent sequence or stochastic processes as martingales, Markov chains or Markov stochastic processes, etc.\\

\noindent (7) Another approach of handling the conditional mathematical expectation is to considered as a projection in Hilbert spaces. This will be considered in Functional Analysis.\\
   
\noindent (8) The fundamental law of Kolmogorov is the generalization of the finite product measure to the infinite (countably or not) case. This generalization needs topological properties to hold.\\

\noindent (9) Specially, Documents Doc 08-06 and Doc 08-07 have been introduced as a key tool for Mathematical Statistics in relation with the Factorization Theorem for sufficient statistics.\\   

\noindent (10) Mathematical Statistics is largely based on Measure Theory, especially on probability laws.\\ 

\noindent (11) The Daniel integral studied in Doc 09-04 is the master key to establish the convergence of stochastic processes.

\bigskip \noindent These are only a sample of remarkable facts we provided just to enlighten the power of measure theory and its impact of a a considerable number of Mathematics fields.\\

\noindent I finish by advising the learner to have as many readings of the book as possible to the point he is able to repeat the whole content himself.\\

\end{document}